\documentclass{amsart}
\usepackage{enumitem}
\usepackage[figuresleft]{rotating}
\usepackage{wrapfig}
\usepackage[letterpaper,margin=1in,centering,marginparwidth=.95in,marginparsep=.01in]{geometry}
\usepackage{amsthm,amsmath,amssymb,amsfonts,mathtools,microtype,cancel}
\usepackage[utf8]{inputenc} 
\usepackage[T1]{fontenc}    
\usepackage[hidelinks,pagebackref]{hyperref}       
\usepackage{booktabs}       
\usepackage{nicefrac}       
\usepackage[widespace]{fourier} 
\usepackage[dvipsnames]{xcolor}
\usepackage{pgf}

\allowdisplaybreaks[1]
\newcommand\hidesubsection[2][]{{\def\addcontentsline##1##2##3{}\expandafter\ifx\expandafter\relax\detokenize{#1}\relax\subsection*{#2}\else\subsection*{#2\textnormal{\ (#1)}}\fi\strut}}
\DeclareMathOperator{\sinc}{sinc}
\DeclareSymbolFontAlphabet{\mathcal}   {symbols}
\DeclareMathAlphabet{\mathcal}{OMS}{Zplm}{m}{n}
\SetMathAlphabet{\mathcal}{bold}{OMS}{zplm}{b}{n}

\usepackage{graphicx}
\graphicspath{ {figs/} }
\usepackage{doi}
\usepackage{placeins}
\usepackage{caption}
\usepackage{subcaption}
\usepackage{tikz}
\usepackage{tkz-euclide}
\usepackage{yfonts}
\usepackage[capitalize,nameinlink,noabbrev]{cleveref}
\usetikzlibrary{spy}
\usetikzlibrary {angles,quotes}
\newtheorem{theorem}{Theorem}[section]
\newtheorem{ltheorem}{Theorem} 
\newtheorem{corollary}[theorem]{Corollary}
\newtheorem{proposition}[theorem]{Proposition}
\newtheorem{lemma}[theorem]{Lemma}
\theoremstyle{definition}
\newtheorem{definition}[theorem]{Definition}
\newtheorem{notation}[theorem]{Notation}
\crefname{notation}{Notation}{Notations}
\crefname{ltheorem}{Theorem}{Theorem}
\newtheorem{remark}[theorem]{Remark}
\numberwithin{equation}{section}
\usetikzlibrary{intersections, calc, angles}
\def\bigstrut{\rule[-0.4\baselineskip]{0pt}{1.2\baselineskip}}
\def\Bigstrut{\rule[-0.5\baselineskip]{0pt}{1.3\baselineskip}}
\newif\ifcomments
\newcounter{commentlabel}
  \newlength{\commentlift}
\makeatletter\@mparswitchfalse\def\commentfl@g{%
    \vbox to\z@{%
      \setlength{\fboxrule}{0.75pt}%
      \setlength{\fboxsep}{0.75pt}%
      \setlength{\commentlift}{1ex}%
      \addtolength{\commentlift}{\fboxrule}%
      \addtolength{\commentlift}{\fboxsep}%
\vss\color{blue}\rlap{\rlap{\vrule\@height\commentlift\@width\fboxrule}\raise \commentlift%
      \hbox{\fcolorbox{blue}{blue!3}{\normalfont\tiny\bfseries\thecommentlabel}}}}}
  \DeclareRobustCommand{\COMMENT}[1]{\commentstrue\stepcounter{commentlabel}%
\hypertarget{\the\value{section}-\the\value{commentlabel}}{}%
\edef\WRITECOM##1{\noexpand\write\COM{\thecommentlabel, S.##1}}\WRITECOM{\the\value{section}, \noexpand\hyperlink{\the\value{section}-\the\value{commentlabel}}{p.\the\value{page}}: {#1}}\@bsphack%
    \commentfl@g%
    \marginpar{\noindent\lineskip=1pt\lineskiplimit=\maxdimen\raggedright%
    \tiny\def\Hrule{\strut\hrule height 1pt\strut}\textbf{\color{blue}\llap{\textbullet}\thecommentlabel:}{\color{blue}\thinspace#1\endgraf}}%
  \@esphack}%
  \DeclareRobustCommand{\COMMENTINLINE}[1]{\commentstrue\stepcounter{commentlabel}%
    \write\COM{\thecommentlabel, S.\the\value{section}, p.\the\value{page}: {#1}}\@bsphack%
      \noindent\setlength{\fboxsep}{0pt}\newline\fcolorbox{blue}{lightgray!20}{\begin{minipage}{\textwidth}\tiny\textbf{\color{red}\thecommentlabel}:\thinspace\Tiny\bfseries#1\par\end{minipage}}\allowbreak%
  \@esphack}\makeatother
\newwrite\COM
\immediate\openout\COM=comments
\AtEndDocument{\ifcomments\section*{Marginal Comments}\label{commentstex}{\let\tiny\relax\let\Tiny\relax\let\Hrule\newline\def\pullout#1, S.{\strut\llap{\textbf{#1: }}\S\ }\obeylines\parindent0pt\everypar{\normalsize\normalfont\pullout}\immediate\closeout\COM\def\MT_extended_eqref:n {\eqref}\def\math@bb{\mathbb }\input comments }\fi}
\makeatother
\let\NEWPAGE\newpage
\def\COMMENT#1{}\def\COMMENTINLINE#1{}\let\NEWPAGE\relax
\newcommand{\arc}[1]{\wideparen{#1}}
\newcommand{\Int}{\operatorname{int}}
\newcommand{\Fc}{\mathcal F}

\let\Hrule\relax
\def\Y{\overset{\textcolor{red}{\text{\Tiny\upshape\bfseries WHY?}}}}
\newcommand{\p}{\ensuremath{\text{\upshape p}}}
\newcommand{\phistar}{\phi_*}
\newcommand{\Phistar}{\Phi_*}
\newcommand{\Bcal}{\ensuremath{\mathcal{B}}}
\newcommand{\Fcal}{\ensuremath{\mathcal{F}}}
\newcommand{\Scal}{\ensuremath{\mathcal{S}}}
\newcommand{\constant}{\ensuremath{\Psi}}
\newcommand{\const}{\ensuremath{\text{const.}}}
\newcommand{\dfn}{\ensuremath{\mathbin{{:}{=}}}}
\newcommand{\nfd}{\ensuremath{\mathbin{{=}{:}}}}
\newcommand{\aeq}{\ensuremath{\mathbin{\overset{\text{\upshape\ae}}{=}}}}
\newcommand{\Mout}{{\color{red}M}^{\text{{\color{red}\upshape out}}}}
\newcommand{\Min}{{\color{blue}M}^{\text{{\color{blue}\upshape in}}}}
\newcommand{\Mrout}{\ensuremath{{\color{red}M^{\text{\upshape out}}_r}}}
\newcommand{\Mrin}{\ensuremath{{\color{blue}M^{\text{\upshape in}}_r}}}
\newcommand{\MRout}{\ensuremath{{\color{red}M^{\text{\upshape out}}_R}}}
\newcommand{\MRin}{\ensuremath{{\color{blue}M^{\text{\upshape in}}_R}}}
\newcommand{\Nin}{\ensuremath{{\color{blue}N^{\text{\upshape in}}}}}
\newcommand{\Nout}{\ensuremath{{\color{red}N^{\text{\upshape out}}}}}
\def\PET{\mathbb{Pet}}
\def\LEM{\mathbb{Lem}}
\newcommand{\boris}{\strut\textcolor{red}{\emph{\bfseries--Boris says: }}}
\newcommand{\wentao}{\strut\textcolor{green!70!black}{\emph{\bfseries--Wentao says: }}}
\title{Ergodicity of asymmetric lemon billiards}
\author{Wentao Fan}
\address{Department of Mathematics\\
        Tufts University\\
        Medford, MA 02155}
\email{wentao.fan@tufts.edu}
\author{Boris Hasselblatt}
\address{Department of Mathematics\\
        Tufts University\\
        Medford, MA 02155}
\email{boris.hasselblatt@tufts.edu}

\hypersetup{
pdftitle={Uniform Hyperbolicity and ergodicity of asymmetric Lemon Billiard},
pdfkeywords={Expansion, Uniform Hyperbolic, Wojtkowski Billiards, Cone, Wave Fronts Curvature, singularity curve, local ergodic theorem, global ergodic},
}
\begin{document}
\begin{abstract}
We show ergodicity of (asymmetric) lemon billiards, billiard tables that are the intersection of two circles of which one contains the centers of both. These do not satisfy the Wojtkowski criteria for hyperbolicity, but we establish \emph{uniform} expansion of vectors in an invariant cone family and alignment of singularity curves. Both of these are difficult, and the approach to the latter seems new. Together, these properties imply ergodicity.
\end{abstract}
\maketitle
\setcounter{tocdepth}{1}\tableofcontents
\section{Introduction}
\label{sec:IntroMainResult}
\label{sec:imr}
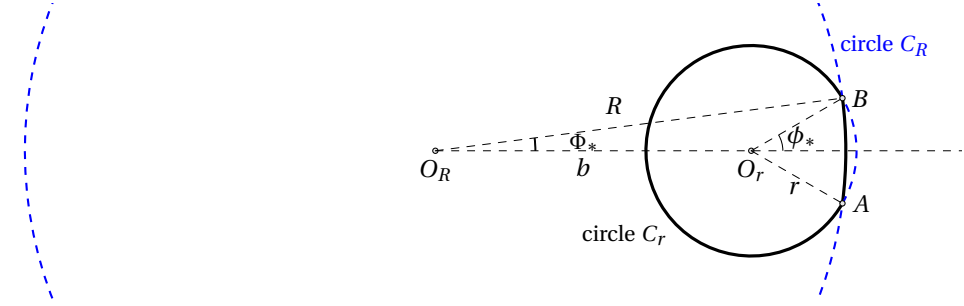
\begin{figure}[ht]
    \begin{center}
    \begin{tikzpicture}[scale=0.7]
    \tkzDefPoint(-6,0){OR};
    \tkzDefPoint(0,0){Or};
    \tkzDefPoint(2,0){C};
    \pgfmathsetmacro{\XValueArc}{2*sin(60)}
    \tkzDefPoint(\XValueArc,-1){A};
    \tkzDefPoint(\XValueArc,1){B};
    \clip(-14,-2.8) rectangle (4, 2.8);
    \node[] at ([shift={(0.8,1.0)}]B) {\textcolor{blue}{\small circle $C_R$}};
    \node[] at ([shift={(-2.4,-1.6)}]Or) {\small circle $C_r$};
    \tkzDrawArc[very thick](OR,A)(B);
    \tkzDrawArc[very thick](Or,B)(A);
    \tkzDrawArc[thick,dashed, blue](OR,B)(A);
    \tkzDrawArc[thick,dashed, blue](Or,A)(B);
    \tkzDrawPoints(OR,Or,A,B);
    \tkzLabelPoint[right](A){$A$};
    \tkzLabelPoint[right](B){$B$};
    \tkzLabelPoint[below](OR){$O_R$};
    \tkzLabelPoint[below](Or){$O_r$};
    \draw [very thin, dashed] (OR) --(5,0);
    \coordinate (C) at (5,0);
    \tkzDrawSegments[very thin,dashed](OR,B);
    \tkzDrawSegments[very thin,dashed](Or,B);
    \tkzDrawSegments[very thin,dashed](Or,A);
    \tkzLabelSegment[above left=0.5pt and 3pt](OR,B){$R$};
    \tkzLabelSegment[above left=-10pt and -5pt](Or,A){$r$};
    \draw (0.6,0) arc (0:30:0.6);
    \node[] at (15:1)  {$\phi_{*}$};
    \node[] at (-3.2,-0.3){$b$};
    \tkzMarkAngle[arc=lll,ultra thin,size=1.9](C,OR,B);
    \tkzLabelPoint[right]([shift={(2.4,0.15)}]OR){\small $\Phistar$};
    \end{tikzpicture}
    \end{center}
    \caption{Lemon billiard Table} 
    \label{fig:BilliardTable}
    \end{figure}

We\COMMENT{\huge\bfseries\color{red}For arXiv posting define \textbackslash COMMENT and \textbackslash newpage to do nothing.} prove ergodicity of the (asymmetric) \emph{lemon billiards} in \cref{fig:BilliardTable}. This involves establishing \emph{uniform} hyperbolicity of the return map to a suitable section (not only nonzero Lyapunov exponents everywhere, see \cref{thm:MroutreturnUniformExpansion} on page \pageref{thm:MroutreturnUniformExpansion}), control of singularity curves (which here look different than elsewhere in the literature---near vertical with sign change of slope), and the use of Hopf chains.

Mathematical billiard systems describe a point particle moving freely in a domain with piecewise smooth boundary and reflecting off the boundary the way a photon does in a mirror, with equal angles before and after collision (specular reflection). The introduction of these to dynamical systems is usually associated with Birkhoff \cite[p.~288f]{HF}, and dispersing (or, rather, hyperbolic) billiards are associated in no small part with the Sinai school. Among the highlights of those developments are Wojtkowski's criteria for hyperbolicity in the presence of focusing boundary circle arcs (the circle to which such an arc belongs must lie inside the billiard table \cite{MR848647,ICLE1985}) and the book by Chernov and Markarian \cite{cb}, which carefully presents the needed technology and is a central reference for the present work. In more recent years, significant efforts have been directed at the study of billiards which are expected to be hyperbolic but do not satisfy the Wojtkowski criteria. Our main result is uniform hyperbolicity of some such billiards---the \emph{lemon billiards} shown in \cref{fig:BilliardTable} for suitable parameters.

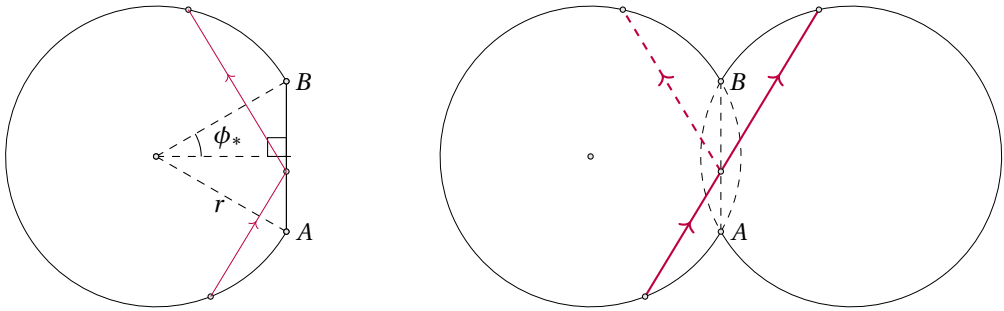
\begin{figure}[hb]
\strut\hfill
    \begin{tikzpicture}
    \tkzDefPoint(0,0){Or};
    \tkzDefPoint(2,0){C};
    \pgfmathsetmacro{\XValueArc}{2*sin(60)}
    \tkzDefPoint(\XValueArc,-1){A};
    \tkzDefPoint(\XValueArc,1){B};
    \tkzDrawArc[name path=Cr](Or,B)(A);
    \tkzDrawPoints(Or,A,B);
    \draw[name path=lineAB,thin](A)--(B);
    \tkzLabelPoint[right](A){$A$};
    \tkzLabelPoint[right](B){$B$};
    \draw [name path=xaxis, thin, dashed] (Or) --(1.83,0);
    \path[name intersections={of=xaxis and lineAB,by={C}}];
    \coordinate (Or_mirror) at ($(Or)!2!(C)$);
    \draw[thin](A)--(B);
    \tkzMarkRightAngle(B,C,Or);
    \tkzDrawSegments[thin,dashed](Or,B);
    \tkzDrawSegments[thin,dashed](Or,A);
    \tkzLabelSegment[above left=-10pt and -5pt](Or,A){$r$};
    \draw (0.6,0) arc (0:30:0.6);
    \node[] at (15:1)  {$\phi_{*}$};
    
    \coordinate (P) at (\XValueArc,-0.2);
    \tkzDefPoint(2*cos(-1.2),2*sin(-1.2)){T};
    \tkzDefPointBy[projection=onto B--A](T)\tkzGetPoint{H}; 
    \coordinate (TReflect) at ($(T)!2.0!(H)$);
    \coordinate (QQ) at ($(TReflect)!2.4!(P)$);
    \path[name path=linePQQ](P)--(QQ);
    \path[name intersections={of=Cr and linePQQ,by={Q}}];
    \tkzDrawPoints(Or,A,B,P,T,Q);
    \begin{scope}[very thin,decoration={
    markings,
    mark=at position 0.6 with {\arrow{>}}}
    ] 
    \draw[purple,postaction={decorate}](T)--(P);
    \draw[purple,postaction={decorate}](P)--(Q);
    \end{scope}

    \end{tikzpicture}\hfill
    \begin{tikzpicture}
    \tkzDefPoint(0,0){Or};
    \tkzDefPoint(2,0){C};
    \pgfmathsetmacro{\XValueArc}{2*sin(60)}
    \tkzDefPoint(\XValueArc,-1){A};
    \tkzDefPoint(\XValueArc,1){B};
    \tkzDrawArc[name path=Cr](Or,B)(A);
    \draw[name path=lineAB,thin,dashed](A)--(B);
    \tkzLabelPoint[right](A){$A$};
    \tkzLabelPoint[right](B){$B$};
    \path[name path=xaxis, thin, dashed] (Or) --(1.83,0);
    \path[name intersections={of=xaxis and lineAB,by={C}}];
    \coordinate (Or_mirror) at ($(Or)!2!(C)$);
    \tkzDrawArc[name path=mirrorCr](Or_mirror,A)(B);
    \tkzDrawArc[dashed](Or,A)(B);
    \tkzDrawArc[dashed](Or_mirror,B)(A);
        \coordinate (P) at (\XValueArc,-0.2);
    \tkzDefPoint(2*cos(-1.2),2*sin(-1.2)){T};
    \tkzDefPointBy[projection=onto B--A](T)\tkzGetPoint{H};

    \coordinate (TReflect) at ($(T)!2.0!(H)$);
    \coordinate (QQ) at ($(TReflect)!2.4!(P)$);
    \path[name path=linePQQ](P)--(QQ);
    \path[name intersections={of=Cr and linePQQ,by={Q}}];
    \tkzDefPointBy[projection=onto B--A](Q)\tkzGetPoint{G};
    \coordinate (QReflect) at ($(Q)!2.0!(G)$);
    \begin{scope}[thick,decoration={
    markings,
    mark=at position 0.6 with {\arrow{>}}}
    ] 
    \draw[purple,postaction={decorate}](T)--(P);
    \draw[purple,postaction={decorate}](P)--(QReflect);
    \draw[purple,postaction={decorate},dashed](P)--(Q);
    \end{scope}
    \tkzDrawPoints(Or,A,B,P,T,Q,QReflect);
    \end{tikzpicture}
    \hfill\strut
    \caption{The 1-petal billiard and its unfolding to the 2-petal billiard}
\label{fig:WojtkowskiBilliardTable}
\end{figure}

These may have first appeared in \cite{Tomsovic} in the form of the intersection of two circles of equal radius. \cite{oaedhalb16} proved their hyperbolicity. We here consider \emph{asymmetric lemon billiards} (see \cref{fig:BilliardTable}), which consist of the intersection of two disks of which one contains the centers of both and were introduced in the article \cite{MR3389766}, which contains the conjecture that they are hyperbolic when the radii differ greatly \cite{oaedhalb16}. The extreme cases of this are on one hand a circular billiard, which is not hyperbolic (indeed, integrable), and on the other hand a disk with a small part of the boundary flattened to a straight line (\cref{fig:WojtkowskiBilliardTable}). The latter is significantly more tractable, and we call it the 1-petal billiard \cite{MR357736,MR530154}. When unfolded along the line segment, it becomes the two-petal billiard shown in \cref{fig:WojtkowskiBilliardTable}, and any focusing circle arc of the billiard table is such that the disk in whose boundary it lies is contained in the billiard table. This produces a defocusing (or overfocusing) mechanism which traces back to Bunimovich \cite{MR357736} and was greatly refined by Wojtkowski  \cite{MR848647} and then Markarian \cite{MR954676}, Donnay \cite{MR1133266}, and Bunimovich \cite{MR1179172}. By design, lemon billiards do not satisfy these criteria, so ergodicity must be established by arguments closely tied to the specifics of the billiard table. Consequently, this work adapts techniques from the theory of dispersing billiards without being able to invoke ready-made theorems that lead to ergodicity. We note that the tractability of the 1-petal billiard is being used in our arguments.

From this perspective, the subject of the present article can be viewed as a perturbation in which the straight boundary segment of the 1-petal billiard in \cref{fig:WojtkowskiBilliardTable} is replaced by    
an arc of a circle with large radius; these systems do not satisfy the Wojtkowski criteria and cannot be unfolded to do so in the way shown in \cref{fig:WojtkowskiBilliardTable} for the 1-petal billiard. Nonetheless, it was shown relatively recently, that they possess an invariant cone family \cite{oaedhalb16,hoalb}, which implies (nonuniform) hyperbolicity and is a preliminary step towards ergodicity.

The present article begins with the next step, which is to establish uniform expansion 
of vectors in cones of a suitable family. This constitutes \emph{uniform} hyperbolicity (\cref{Corollary:EuclideanExpansionMhat,corollary:MroutreturnUniformExpansion,Col:EuclideanExpansionMhat1}). We next control the singularity curves, an essential step towards local ergodicity, and one where significant novelty resides because in this context, the slope of some of these changes sign (\cref{fig:ArchSingularities})! Ergodicity is finally proved in \cref{corollary:LETonSection,Proposition:globalErgodic}.
\begin{ltheorem}\label{RMASUH}
For almost all $\phistar\in \big(0,\tan^{-1}{({1}/{3})}\big)$, there exists an $R_{\text{\upshape HF}}(r,\phistar)$
such that if $R>R_{\text{\upshape HF}}(r,\phistar)$, 
then:\COMMENT{\Large This now includes everything, so anything like "proof of Theorem A" must be edited!} \begin{enumerate}
    \item The ``half-quadrant'' cone family $C_x\dfn\big\{(d\phi,d\theta)\bigm|\frac{d\theta}{d\phi}\in [0,1]\big\}$ is strictly invariant under the return map $\hat F$ to the global section $\hat{M}$ from \eqref{def:Mhat} of next to last or sole collisions with $\Gamma_r$.
    \item This return map is uniformly hyperbolic: there exist $c>0$ and $\Lambda>1$ such that $\dfrac{\|D\hat{F}^n(dx)\|}{\|dx\|}>c\Lambda^n$ for nonsingular $x\in\hat{M}$, $dx\in C_x$ and $n\in\mathbb{N}$ (\cref{Corollary:EuclideanExpansionMhat,corollary:MroutreturnUniformExpansion,Col:EuclideanExpansionMhat1}).
    \item The singularity curves are properly aligned.
    \item The billiard map is ergodic (\cref{corollary:LETonSection,Proposition:globalErgodic}).
\end{enumerate}
\end{ltheorem}
In fact, we expect the K-property and higher mixing rates (\cref{REMFutureWork}). 

We consider a different section from those in the prior works, though several choices of sections yield hyperbolicity. The cone family is a little narrower as well (see \cref{thm:MroutreturnUniformExpansion}). An important feature of our choice of section is that it adresses hitherto insurmountable issues with singularity curves.

One price we pay for these results is that $R_{\text{\upshape HF}}$ in \cref{RMASUH} has to be quite large; part of our arguments are in effect of a perturbative nature with respect to the 1-petal billiard, where ``$R=\infty$.''
The main result of \cite{hoalb} is that for $\phistar\in(0,\frac{\pi}{2})$ and
\begin{equation}\label{eqJZHypCond0}
R>R_{\text{JZ}}(r,\phistar)\dfn r\cdot\max\big\{\frac{14.6}{\min{\{\phistar,\frac{\pi}{2}-\phistar\}}},\frac{147}{\sin^2{\phistar}},1773.7\big\},
\end{equation} the return map to the global section \(\hat M\) from \cref{def:SectionSetsDefs} preserves the cones defined by the positive quadrant in \cref{corollary:InvQuadrant}
\cite[Equations (3.1) and (3.2)]{hoalb}.\COMMENT{Refer to here when "1700" arises.\newline\wentao{Agree. Work In Progress}}

This present work establishes that for $a.e.\, \phistar\in\big(0,\tan^{-1}(1/3)\big)$, there exists an
\begin{equation}\label{eqJZHypCond}
    R_{\text{\upshape HF}}(r,\phistar)
    \ge r\cdot\max{\Big\{\frac{30000}{\sin^2\phistar},\frac{324}{\sin^2{(\phistar/2)\sin^2(\phistar)}}\Big\}}\ge R_{\text{JZ}}(r,\phistar)\text{ such that }R>R_{\text{\upshape HF}}(r,\phistar)\Rightarrow\text{ hyperbolicity,}
\end{equation} 
that is, the return maps on certain global sections are uniformly hyperbolic.\COMMENT{It will be good to keep thinking about what references to \cref{eqJZHypCond} we need.\newline\wentao{Agree. Work in Progress}}

Unlike in \cite{hoalb}, we are not able to provide a closed form formula for $R_{\text{\upshape HF}}(r,\phistar)$ because its dependence on $\phistar$ is significantly more subtle and related to the need to exclude some values of $\phistar$. These are not numerous and can be described as follows: The chord $\overline{AB}$\COMMENT{Put the exceptional $\phistar$s definitions in section 8 which comes from 1-petal. Forward reference to that from here.\wentao{done}} is part of a periodic orbit in the 1-petal billiard bounded by $C_r$ that does not hit (the interior of) $\arc{AB}$, so it subtends an angle $2\pi/n$ of $C_r$ for some integer $n$. It is the \emph{generalized $\phistar$} defined in \cref{def:ExceptionPhistar}.\COMMENT{How does this help here?} Note that the exceptional $\phi_*$ set is defined by the 1-Petal billiard so that it is independent of $R$. We provide an algorithm in \cref{def:AcomputationForRHF} to determine $R_{\text{\upshape HF}}$.

\begin{proof}[Proof of uniform hyperbolicity in \cref{RMASUH}]
As is common in the literature, we prove expansion of the so-called p-metric\COMMENT{ref}, and expansion with respect to the Eucliden metric then follows\COMMENT{ref.}

For the section $\hat{M}$ in \cref{def:SectionSetsDefs}\eqref{def:Mhat} with its return map $\hat{F}$, for $x\in\hat{M}$ the return orbit segment $x, \Fc(x),\cdots,\Fc^{\sigma(x)}(x)\in\hat{M}$ in \cref{def:MhatReturnOrbitSegment} with $x_1=(\Phi_1,\theta_1)\in M_R$ in \eqref{eqPtsFromeqMhatOrbSeg} and for all $dx\in C_x$, there are the following cases.\COMMENT{Comment on transition from p-metric to euclidean}
    
In case (a0) of  \eqref{eqMainCases}, $\sin\theta_1>\sqrt{4r/R}$, which corresponds to large collisions angles. Here, $\frac{\|D\hat{F}_x(dx)\|_{\p}}{\|dx\|_{\p}}>1+\lambda_c$ with constant $\lambda_c>0$ from \cref{TEcase_a0}. Thus we get expansion in this situation.

Not so in cases (a1)(b)(c) of \eqref{eqMainCases}, \cref{TEcase_a1,TEcase_b,TEcase_c} give that $\frac{\|D\hat{F}_x(dx)\|_{\p}}{\|dx\|_{\p}}>0.05$. Hence, there is an uniform $\overline{l}>0$  given in \cref{def:ApettrajectoryExpansion} determined by the 1-petal billiard in \cref{lemma:PETExpExpansion}, thus independent of $R$ such that 
    \[\frac{\|D\hat{F}^{\overline{l}+1}_x(dx)\|_{\p}}{\|dx\|_{\p}}>0.05\times22.5=1.125.\]

    Let $\overline{\Lambda}=\min\{1+\lambda_c,1.125\}$, where $\lambda_c$ is from \cref{TEcase_a0}. \cref{def:ApettrajectoryExpansion,def:AcomputationForRHF,lemma:Trajectory_approximation_1}\eqref{eq:xorbitA0Condition} ensure that among every consecutive $\overline{l}+1$ times $\hat{M}$ return orbit segment defined in \eqref{eqMhatOrbSeg} there can exist at most $1$ time return orbit segment not in case (a0) of \eqref{eqMainCases}. Therefore, We can conclude that for all nonsingular $x$ and $\forall n\ge1$, \[\frac{\|D\hat{F}^{n}_x(dx)\|_{\p}}{\|dx\|_{\p}}>(\overline{\Lambda})^{\left\lfloor\frac{n}{\overline{l}+1}\right\rfloor}(0.05).\] 

    Therefore, $\frac{\|D\hat{F}^{n}_x(dx)\|_{\p}}{\|dx\|_{\p}}>(\overline{\Lambda})^{\frac{n}{\overline{l}+1}-1}(0.05)=c\Lambda^n$, where $c=(0.05)/\overline{\Lambda}$, $\Lambda=(\overline{\Lambda})^{\frac{1}{1+\overline{l}}}$.
\end{proof}
\subsection{Structure of this paper}

This proof invoked results that will be established in subsequent parts of this work, all the way through \cref{sec:SACTWFELBWB}. These and subsequent arguments are organized as follows.

\cref{sec:ECE} states lower bounds for the differential of the billiard map in a case-by-case fashion (\cref{eqMainCases,TEcase_a0,TEcase_a1,TEcase_b,TEcase_c}) and establishes that those cases in which these lower bounds do not produce expansion are confined to a small portion of the section (\cref{contraction_region}). This is sufficient information to show strict invariance of the positive quadrant (\cref{corollary:InvQuadrant}; see also \cref{corollary:halfquadrantcone}) and strict expansion in it upon iteration (\cref{def:ApettrajectoryExpansion}). The proofs of these items in \cref{sec:ECE} are carried out in later sections as follows. 

\cref{sec:FPCLA} prepares for delicate estimates across various parameter ranges by establishing monotonicity results which reduce the needed estimates to considerations at the endpoints of the various ranges. These are lengthy arguments, right up to page~\pageref{sec: PPDSR}. 

\cref{sec: PPDSR} provides the proof that the aforementioned ``bad region’’ is indeed small (\cref{contraction_region}) and then gathers and adapts ingredients from \cite{hoalb,cb}. 

\cref{sec:ProofTECase_a0_a1} then establishes expansion in cases (a0) and (a1). These may be the easiest cases, but this takes up pages \pageref{subsec:AFELB}--\pageref{sec:ProofTECase_b}. 

\cref{sec:ProofTECase_b} covers Case (b) (\cref{TEcase_b,subsec:ProofTECase_b}). 

\cref{sec:ProofTECase_c} finally tackles the remaining case by proving \cref{TEcase_c}. 

One important deferred item remains for \cref{sec:SACTWFELBWB}, which is to control the dynamics in the small ``bad region’’ by approximating those trajectories by their counterpart for the 1-petal billiard (\cref{fig:WojtkowskiBilliardTable}). This necessitates the assumption that the larger of the two radii is so much larger than the other.\COMMENT{I am sure that this could be better, but it is a start. Remove this note if this is good enough.\wentao{Mention the three approximation theorems? The radii ratio is so large so that the three approximation gives a way to see the expansion is good?}}

As the preceding narrative and the page count of the efforts it describes indicate, establishing uniform hyperbolicity is the most laborious undertaking in this work. The next steps towards ergodicity only begin in \cref{sec:LET}. \cite{cb} provides the blueprint for working towards ergodicity from uniform hyperbolicity. \cref{sec:LET} starts work towards local ergodicity. It establishes local ergodicity subject to the assumption of proper alignment of singularity curves. This, and then global ergodicity, are finally established in \cref{sec:SCLB}. This is a nontrivial coda in that it deals with singularity curves whose geometry had not previously been handled.
\subsection*{Acknowledgments}
We owe much gratitude to Pengfei Zhang and Xin Jin. This work builds on their seminal article \cite{hoalb}, and their support and encouragement were invaluable, notably conversations with Pengfei Zhang at various stages of this project. Conversations with Marco Lenci, Domokos Szasz, and Maciej Wojtkowski were also quite helpful in directing us to pertinent works which proved useful \cite{MR1636897,OKPPHB}.\COMMENT{Insert all pertinent references here.} Finally, the support and perspective of Hongkun Zhang provided a significant boost as this article neared completion.
\section{Preliminaries}\label{subsec:NotationsdefinitionsAndpropositions}
We introduce basic notations and observations.
\begin{notation}[Configuration parameters, phase space, singularities, and sections]\label{def:BasicNotations}\label{eqMroutMrin}\label{def:LemonTableConfiguration}\label{def:BasicNotations0}\label{def:StandardCoordinateTable}
     The (asymmetric) lemon billiard table is the intersection of two disks with radii $r$ and $R$, respectively,\footnote{The notations would be slightly lighter if we were to choose $r=1$ throughout, but we retain both radii to conform with the notation in other articles on asymmetric lemon billiards.} and centers $O_r$ and $O_R$ at distance $b$ as in \Cref{fig:BilliardTable}.
     \begin{itemize}
\item When we use cartesian coordinates, they will be chosen such that $O_r=(0,0)$ and $O_R=(0,b)$ (see, for example, \cref{fig:CaseAMonotone,fig:R1_R2_d_0,fig:lengthfraction,fig:twoTangentLines,fig:QonUpperHalfCircle}).
\item The pertinent notations are as in \Cref{fig:BilliardTable}: $A$, $B$ are the intersection points, $2\phistar=\measuredangle A O_r B$ and $2\Phistar=\measuredangle A O_R B$, so $R\sin\Phistar=r\sin\phistar$.\COMMENT{I suggest that using distinctive greek letters for these two angles is nicer than always needing to attach a *. $\Sigma,\vartheta,\varrho,\zeta,\varsigma$ are not currently used elsewhere. \small I defined macros \textbackslash phistar and \textbackslash Phistar for these angles \Tiny and  globally replaced "hi\_*" by "histar." That could defer a final decision on which letters to use. One criterion would be that $R\sin\Phistar=r\sin\phistar$ looks reasonable with that choice. Also, it seems sensible to use an uppercase-lowercase pair as is currently the case.\newline But maybe this is all a bad idea\dots}\COMMENT{Is the notation $\theta_*$ in \cref{rmk:extensionDefs} ideal?\wentao{ehr...,let's discuss. using * meaning on the boundary by consistency.}}
\item We denote by $\Gamma_r$ and $\Gamma_R$  the two smooth pieces of the table boundary with radii $r$, $R$, respectively, and by $C_r$ and $C_R$ the circles to which they belong.
    \item  The phase space $M\dfn M_r\cup M_R$ consists of inward-pointing vectors on \(\Gamma_r\cup\Gamma_R\), where 
    $M_r=[\phistar,2\pi-\phistar]\times[0,\pi]$ and $M_R=[-\Phistar,\Phistar]\times[0,\pi]$ for the collisions on $\Gamma_r$, $\Gamma_R$, respectively. Here, the position coordinates are the angle parameters $\phi\in[\phistar,2\pi-\phistar]$ on $C_r$ and  $\Phi\in [-\Phistar,\Phistar]$ on $C_R$. The direction coordinates $\theta\in [0,\pi]$ are the collision angle with the tangent direction. We can also write $\theta=\pi/2-\varphi$, where $\varphi$ is the collision angle with the inward normal vector. Then $d\phi=\nicefrac{ds}{r}$, $d\Phi=\nicefrac{ds}{R}$, where $ds$ is the arc-length differential on $\Gamma_r$ and $\Gamma_R$, respectively. The tangential direction on $\Gamma_r$ and $\Gamma_R$ is given by the counter-clockwise orientation on those boundary.
    \item A piece of a billiard table boundary is a focusing or dispersive arc if it has negative or positive curvatures given the counterclockwise orientation of the arc and a flat arc if it has 0 curvature \cite[Section 2.1, Equation (2.4)]{cb}. The lemon billiard has two focusing boundary pieces.
    \item For $x\in M$, we denote by $p(x)$ its position on the lemon billiard table.
    \item Denote by $\mathcal{F}$ the billiard map which sends a vector in \(M\) to the outward vector after the next collision---provided that this next collision is not at a corner \cite[Equation (2.18)]{cb}. (This is the return map of the billiard flow to the boundary \cite[Sections 2.9--2.12]{cb}.
    \item The boundary of the phase space is $\partial M\dfn \Scal_0\dfn\partial M_r\cup\partial M_R$

    $\Scal_1\dfn\big\{x\text{ is an interior point of }M\text{ with }\mathcal{F}(x)\in \Scal_0\big\}\cup \Scal_0$, 
    
    $\Scal_{-1}\dfn\big\{x\text{ is an interior point of } M \text{ with }\mathcal{F}^{-1}(x)\in \Scal_0\big\}\cup \Scal_0$. 
    \item Recursively define $\Scal_k\dfn \Scal_{k-1}\cup\big\{x\bigm|x\in M\setminus \Scal_{k-1},\mathcal{F}(x)\in \Scal_{k-1}\big\}$, $\Scal_{-k}\dfn \Scal_{-k+1}\cup\big\{x\bigm|x\in M\setminus \Scal_{-k+1},\mathcal{F}(x)\in \Scal_{-k+1}\big\}$ for $k\ge 1$ \cite[Equation (2.27)]{cb}. 
    \item Singular points. The forward and backward singularity sets are $\Scal_{\infty}\dfn \bigcup_{k=0}^{\infty}\Scal_k$ and $\Scal_{-\infty}\dfn \bigcup_{k=0}^{\infty}\Scal_{-k}$ respectively \cite[Equation (2.28)]{cb}. Those singular points are either ``grazing'' vectors on $\Gamma_r\cup\Gamma_R$  (with $\theta=0\text{ or }\pi$) tangential to the billiard boundary (this does happen in the lemon billiard) or those vectors coming from a corner or going to a corner. The singularity sets have Lebesgue measure 0 \cite[Section 2.11 and equation (2.28)]{cb}, i.e., a.e. $x\in M$ is a nonsingular point.
    \item The billiard map $\mathcal{F}$ is a diffeomorphism from $M\setminus \Scal_{1}$ to $M\setminus \Scal_{-1}$ \cite[Theorem 2.33 and equation (2.26)]{cb}. 
    \item We consider the sections \[\Mrout\dfn\Int M_r\cap \mathcal{F}^{-1}(\Int M_R)\text{ and }\Mrin\dfn\Int M_r\cap \mathcal{F}(\Int M_R) \]
shown in \cref{fig:Mr}\COMMENT{Explain the dark disks in that figure\wentao{I think it is explained in the figure caption?}} and  
\[\MRout\dfn\Int M_R\cap \mathcal{F}^{-1}(\Int M_r)\text{ and }\MRin\dfn M_R\cap \mathcal{F}(\Int M_r)\]
shown in \cref{fig:MR}. Here  ``interior'' merely omits corners. These sections consist of the last and first collisions with $\Gamma_r$ and $\Gamma_R$, respectively. $\Mrin,\Mrout$ are global sections of $M$. That is, $\bigcup_{n\ge0}\mathcal{F}^n(\Mrin)\aeq M\aeq\bigcup_{n\ge0}\mathcal{F}^n(\Mrout)$ since every $x=(\phi,\theta)\in M_r$ (except for finitely many $\theta$ such that $\theta/\pi$ is rational) has an $\mathcal F$-orbit across $\Mrout$ and $\Mrin$.

\end{itemize}
\end{notation}
\begin{remark}
We noted that the billiard map $\mathcal{F}$ is not always well defined on the boundary of the phase space, i.e., for the corners (and in other cases, grazing collisions) in the table. Since the angle at corners of the lemon billiard does not divide $\pi$, there are two ways to continuously extend the definition of the billiard map to the phase-space boundary by considering the collision as one with one or the other of the adjacent boundary pieces  \cite[Section 2.8]{cb}.\COMMENT{Forward reference to where we do this.\wentao{Actually no need...Let's discuss.}}
\end{remark}
\begin{figure}[ht]
\begin{minipage}{.6\textwidth}
\begin{center}
    \begin{tikzpicture}[xscale=.6,yscale=.6]
        \pgfmathsetmacro{\PHISTAR}{0.8}
        \pgfmathsetmacro{\WIDTH}{7}
        \pgfmathsetmacro{\LENGTH}{7}
        \tkzDefPoint(-\LENGTH,0.5*\WIDTH){LLM};
        \tkzDefPoint(-\LENGTH,0){LLL};
        \tkzDefPoint(-\LENGTH+\PHISTAR,0){LL};
        \tkzDefPoint(-\LENGTH+\PHISTAR,\WIDTH){LU};
        \tkzDefPoint(-\LENGTH,\WIDTH){LLU};
        \tkzDefPoint(\LENGTH-\PHISTAR,0){RL};
        \tkzDefPoint(\LENGTH-\PHISTAR,\WIDTH){RU};
        \tkzDefPoint(\LENGTH,0){RRL};
        \tkzDefPoint(\LENGTH,\WIDTH){RRU};
        \tkzDefPoint(-\LENGTH+\PHISTAR,\PHISTAR){X};
        \tkzDefPoint(-\LENGTH+\PHISTAR,\WIDTH-\PHISTAR){IX};
        \tkzDefPoint(\LENGTH-\PHISTAR,\PHISTAR){IY};
        \tkzDefPoint(\LENGTH-\PHISTAR,\WIDTH-\PHISTAR){Y};
        \tkzLabelPoint[below](LL){$\phistar$};
        \tkzLabelPoint[below](RL){$2\pi-\phistar\qquad$};
        \tkzLabelPoint[left](LLL){0};
        \tkzLabelPoint[left](LLM){$\frac{\pi}{2}$};
        \tkzLabelPoint[left](LLU){$\pi$};
        \tkzLabelPoint[below](LLL){0};
        \tkzLabelPoint[below](RRL){\ $2\pi$};
        \tkzLabelPoint[right](IY){$x_*$};
        \tkzLabelPoint[left](IX){$y_*$};
        \tkzLabelPoint[right](Y){$I x_*$};
        \tkzLabelPoint[left](X){$I y_*$};

        \draw [thick] (LL) --(LU);
        \draw [thick] (LU) --(RU);
        \draw [thick] (RL) --(RU);
        \draw [thick] (LL) --(RL);
        \draw [thin,dashed](LLL)--(LL);
        \draw [thin,dashed](LLU)--(LU);
        \draw [thin,dashed](LLU)--(LLL);
        \draw [thin,dashed](RRL)--(RL);
        \draw [thin,dashed](RRU)--(RU);
        \draw [thin,dashed](RRU)--(RRL);
        \fill [blue, opacity=10/30](LL) -- (Y) -- (RU) -- (X) -- cycle;
        \fill [red, opacity=10/30](LU) -- (IY) -- (RL) -- (IX) -- cycle;
        \begin{scope}
            \clip (LL) -- (X) -- (RU) -- cycle;
            \draw[fill=blue] circle[at=(X),radius=0.35];
        \end{scope}
        \begin{scope}
            \clip (RU) -- (Y) -- (LL) -- cycle;
            \draw[fill=blue] circle[at=(Y),radius=0.35];
        \end{scope}
        \begin{scope}
            \clip (RL) -- (IX) -- (LU) -- cycle;
            \draw[fill=red] circle[at=(IX),radius=0.35];
        \end{scope}
        \begin{scope}
            \clip (LU) -- (IY) -- (RL) -- cycle;
        \draw[fill=red] circle[at=(IY),radius=0.35];
        \end{scope}
    \end{tikzpicture}
    \end{center}
    \caption{The region $M_r$. The \textcolor{blue}{blue} parallelogram region is \Mrin, and the \textcolor{red}{red} one is \Mrout. The small sectors are `bad' regions (\cref{REMBadRegions}).} 
    \label{fig:Mr}    
\end{minipage}\hfill
\begin{minipage}{.35\textwidth}
\begin{center}
    \begin{tikzpicture}[scale=.28]
        \pgfmathsetmacro{\PHISTAR}{2.5};
        \pgfmathsetmacro{\WIDTH}{15};
        \pgfmathsetmacro{\LENGTH}{7.5};
        \tkzDefPoint(-\LENGTH,0.5*\WIDTH){LLM};
        \tkzDefPoint(-\LENGTH,0){LLL};
        \tkzDefPoint(-\PHISTAR,0){LL};
        \tkzDefPoint(-\PHISTAR,\WIDTH){LU};
        \tkzDefPoint(-\LENGTH,\WIDTH){LLU};
        \tkzDefPoint(\PHISTAR,0){RL};
        \tkzDefPoint(\PHISTAR,\WIDTH){RU};
        \tkzDefPoint(\LENGTH,0){RRL};
        \tkzDefPoint(\LENGTH,\WIDTH){RRU};
        \tkzDefPoint(-\PHISTAR,\PHISTAR){X};
        \tkzDefPoint(-\PHISTAR,\WIDTH-\PHISTAR){IX};
        \tkzDefPoint(\PHISTAR,\WIDTH-\PHISTAR){Y};
        \tkzDefPoint(\PHISTAR,\PHISTAR){IY};
        \draw [thick] (LL) --(LU);
        \draw [thick] (LU) --(RU);
        \draw [thick] (RL) --(RU);
        \draw [thick] (LL) --(RL);
        \draw [thin,dashed](LLL)--(LL);
        \draw [thin,dashed](LLU)--(LU);
        \draw [thin,dashed](LLU)--(LLL);
        \draw [thin,dashed](RRL)--(RL);
        \draw [thin,dashed](RRU)--(RU);
        \draw [thin,dashed](RRU)--(RRL);
         \fill [blue, opacity=10/30](LL) -- (IY) -- (RU) -- (IX) -- cycle;
        \fill [red, opacity=10/30](LU) -- (X) -- (RL) -- (Y) -- cycle;
        \tkzLabelPoint[below](LL){$-\Phistar$};
        \tkzLabelPoint[left](LLL){$0$};
        \tkzLabelPoint[below](RL){$+\Phistar$};
        \tkzLabelPoint[left](LLU){$\pi$};
        \tkzLabelPoint[left](LLM){${\pi}/{2}$};
    \end{tikzpicture}
\end{center}
  \captionsetup{width=\linewidth}
    \caption{The region $M_R$. The \textcolor{red}{red} region is \MRout, the \textcolor{blue}{blue} region is \MRin.}
    \label{fig:MR}    
\end{minipage}
\end{figure}

\begin{definition}[Symmetries]\label{DefSymmetries}
    \cref{fig:Mr,fig:MR,fig:MrN} exhibit some symmetry, and there are indeed two symmetries in the lemon billiard system that give rise to this,\COMMENT{Is this the right place for this remark, or is there a better earlier place?\wentao{I now moved this remark to be right after \cref{fig:Mr,fig:MR,fig:MrN}}} One of these is \emph{reversibility} of the billiard \emph{flow} (which holds for all billiard systems): changing the sign of the velocity vector in any billiard reverses the motion and the collision sequence. Equivalently, the flip map $I$ from \cref{def:SectionSetsDefs} conjugates the billiard map to its inverse: $I\Fc=\Fc^{-1}I$.
    \begin{equation}\label{eq:InversionSymmetryDef}
        \begin{aligned}
        I\colon M_r\sqcup  M_R&\rightarrow M_r\sqcup M_R,\\ 
        (\phi,\theta)&\mapsto(\phi,\pi-\theta)\in M_r,\text{ if }(\phi,\theta)\in M_r.\\
        (\Phi,\theta)&\mapsto(\Phi,\pi-\theta)\in M_R,\text{ if }(\Phi,\theta)\in M_R.
        \end{aligned}
    \end{equation}
    In the phase space (\cref{fig:Mr,fig:MR}), this symmetry is manifested as an up-down symmetry, which interchanges blue and red in \cref{fig:Mr,fig:MR}.\COMMENT{Possibly improve after discussing} 
    
    The other symmetry corresponds to the mirror\COMMENT{\color{red}Replace "fold-position" and (worse!) "fold-position map" {\Large THROUGHOUT} by "mirror" or similar; "fold-position" and "fold" will confuse readers (it confused me!) because nothing is being folded.} symmetry of the lemon billiard table with respect to the line through the centers of the defining circles. In the phase space, this fold-position map symmetry is manifested as follows.\COMMENT{Insert lucid text here.\wentao{edited}}
    \begin{equation}\label{eq:FoldSymmetryDef}
    \begin{aligned}
        J\colon M_r\sqcup M_R&\rightarrow M_r\sqcup M_R,\\ 
        (\phi,\theta)&\mapsto(2\pi-\phi,\theta)\in M_r,\text{ if }(\phi,\theta)\in M_r.\\
        (\Phi,\theta)&\mapsto(-\Phi,\theta)\in M_R,\text{ if }(\Phi,\theta)\in M_R,
    \end{aligned}
    \end{equation}which is a left-right symmetry and also interchanges blue and red in \cref{fig:Mr,fig:MR}. It also conjugates the billiard map to its inverse: $J\Fc=\Fc^{-1}J$.\COMMENT{This sounds wrong, coming from a symmetry of the configuration space.} And $I $, $J$ are symmetries/diffeomorphisms between $\Mrin$ and $\Mrout$, between $\MRin$ and $\MRout$.\COMMENT{\color{red} Not sure this is true; we must definitely discuss $J$!}
        
        In \cref{fig:Mr,fig:MrN} $I \circ J=J\circ I$ is the mirror symmetry that folds\COMMENT{What does this mean???} the position and flips\COMMENT{Flips???\wentao{right}} the collision direction. When restricted on $M_r$, it is the symmetry with respect to the center $(\pi,\pi/2)$ of the phase space $M_r$.  When restricted on $M_R$, it is the symmetry with respect to the center $(0,\pi/2)$ of the phase space $M_R$. $I \circ J\circ\Fc=\Fc\circ I\circ J$.\COMMENT{{\Large\color{red}Important: }Neither $M_r$ nor $M_R$ is a phase space; both are \emph{components} of the (one and only) phase space of the lemon billiard. Find {\bfseries every} occurrence of "phase space" to make sure that each instance is fixed.}
\end{definition}

\begin{notation}[\emph{p-metric}, grading of sections, ``bad'' regions]\label{def:SectionSetsDefs}\hfill

\begin{itemize}
\item Define the \emph{p-metric} by $\|dx\|_\p=\sin{\theta}|ds|$ for $dx\in T_x(M)$, where $ds$ is the differential of the arc-length parameter \cite[page 58]{cb}.\footnote{Here ``p'' stands for ``\textbf perpendicular,'' per Leonid Bunimovich---and this is a \textbf pseudometric.}
\item For $(\phi,\theta)=x\in M$, define the flip (or reversal) operation $I x=(\phi,\pi-\theta)$. This interchanges pre- and post-collision directions. In particular, for the corner points $y_*=(\phistar,\pi-\phistar)$ and $x_*=(2\pi-\phistar,\phistar)$ in \cref{fig:MrN}, we have $I y_*=(\phistar,\phistar)$ and $I x_*=(2\pi-\phistar,\pi-\phistar)$. (See \cref{DefSymmetries} below.) $y_*$ and $x_*$ in \cref{fig:Mr,fig:MrN} correspond to the two trajectories $\overrightarrow{AB}$ and $\overrightarrow{BA}$ between the corners $A$, $B$ of the billiard table in \cref{fig:BilliardTable}.
\item$B(Ix_*,\delta_*)$, $B(Iy_*,\delta_*)$, $B(x_*,\delta_*)$, $B(y_*,\delta_*)$ are the neighborhoods of radius $\delta_*\dfn{17}\sqrt{r/R}/\sin{(\phistar/2)}$.\COMMENT{Is $\delta_*$ ideal?}
\item ``Bad'' regions as in \cref{fig:MrN}:
\[\Nin\dfn\big(B(Ix_*,\delta_*)\cup B(Iy_*,\delta_*)\big)\cap \Mrin\text{
and }\Nout\dfn\big(B(x_*,\delta_*)\cup B(y_*,\delta_*)\big)\cap \Mrout.\]
\item For $x\in \Mrin$, let 
\[
m(x)\dfn\max\big\{m\ge 0\mid p(\mathcal{F}^i(x))\in\operatorname{int}\Gamma_r\text{ for }0\le i\le m\big\} 
\]
be the number of subsequent collisions with $\Gamma_r$.
\item For $x\in \MRin$, let 
\[
m(x)\dfn\max\big\{m\ge 0\mid p(\mathcal{F}^i(x))\in\operatorname{int}\Gamma_R\text{ for }0\le i\le m\big\} 
\]
be the number of subsequent collisions with $\Gamma_R$.
\item Let $\Min_{r,n}\dfn\big\{x\in \Mrin\big| m(x)=n\big\}$ and $\Mout_{r,n}\dfn I(\Min_{r,n})$ (see \Cref{fig:MrN}).\label{def-Mri}
We note that
\begin{equation}\label{EQMrRinoutgrading}\begin{aligned}\Mrin&=(\bigsqcup_{n=0}^{\infty}\Min_{r,n})\cup\{\text{singular points in }\Mrin\}\aeq\bigsqcup_{n=0}^{\infty}\Min_{r,n},\\ 
\Mrout&=(\bigsqcup_{n=0}^{\infty}\Mout_{r,n})\cup\{\text{singular points in }\Mrout\}\aeq\bigsqcup_{n=0}^{\infty}\Mout_{r,n},\\ 
\MRin&=(\bigsqcup_{n=0}^{\infty}\Min_{R,n})\cup\{\text{singular points in }\MRin\}\aeq\bigsqcup_{n=0}^{\infty}\Min_{R,n},\\
\MRout&=(\bigsqcup_{n=0}^{\infty}\Mout_{R,n})\cup\{\text{singular points in }\Mrout\}\aeq\bigsqcup_{n=0}^{\infty}\Mout_{R,n}.\end{aligned}\end{equation}
\item The domain \begin{equation}\label{def:Mhat}\hat{M}_{-}\dfn\hat M\dfn\big(\bigsqcup_{n\ge1}\mathcal{F}^{n-1}\Min_{r,n}\big)\sqcup\Min_{r,0}\aeq\mathcal{F}^{-1}\big(\Mrout\smallsetminus{\Mrin}\big)\sqcup\big(\Mrout\cap\Mrin\big)\end{equation} of next to last or sole collisions with $\Gamma_r$ is shown in \cref{fig:MrN},\COMMENT{Might it be possible to have a picture of $\hat M$ in the \emph{\Large configuration space}?} and $\hat{F}$ is the return map to $\hat{M}$.\label{def:Mhatdefinition} (The implicit null set consists of those points in $\Mrin$ whose backward orbit never hits $\Gamma_R$.)\COMMENT{\huge This does not seem to make sense. That set should be empty.}
\item 
As shown in \cref{fig:hatUplusUminus},
\begin{equation}\label{def:symetricsectionMhatPlusMinus}
I\hat M=J\hat M=\hat{M}_{+}\dfn\hat{M}_1\dfn\big(\bigsqcup_{n\ge1}\Fc\Min_{r,n}\big)\sqcup\Min_{r,0}\aeq(\Mrin\cap\Mrout)\sqcup\Fc(\Mrin\smallsetminus\Mrout).
\end{equation}
\end{itemize}
\end{notation}
\begin{remark}\label{remark:M1hatMhattwoglobalsection}
    $\hat{M}_{\pm}$ are global sections, i.e., $M\aeq\bigcup_{k\in\mathbb{Z}}\Fc^{k}(\hat{M}_{\pm})$.
\end{remark}
\cref{EQMrRinoutgrading} implies that a nonsingular $x\in\hat{M}$ has a return orbit segment as follows, with $x_0\in\Mrout$, $x_1\in\MRin$, and $x_2\in\Mrin$.
%
%
%
%
\begin{figure}[ht]
\begin{center}
    \begin{tikzpicture}[xscale=1.0,yscale=1.0]
        \pgfmathsetmacro{\PHISTAR}{0.68}
        \pgfmathsetmacro{\WIDTH}{7}
        \pgfmathsetmacro{\LENGTH}{7}
        \tkzDefPoint(-\LENGTH,0.5*\WIDTH){LLM};
        \tkzDefPoint(-\LENGTH,0){LLL};
        \tkzDefPoint(-\LENGTH+\PHISTAR,0){LL};
        \tkzDefPoint(-\LENGTH+\PHISTAR,\WIDTH){LU};
        \tkzDefPoint(-\LENGTH,\WIDTH){LLU};
        \tkzDefPoint(\LENGTH-\PHISTAR,0){RL};
        \tkzDefPoint(\LENGTH-\PHISTAR,\WIDTH){RU};
        \tkzDefPoint(\LENGTH,0){RRL};
        \tkzDefPoint(\LENGTH,\WIDTH){RRU};
        \tkzDefPoint(-\LENGTH+\PHISTAR,\PHISTAR){X};
        \tkzDefPoint(-\LENGTH+\PHISTAR,\WIDTH-\PHISTAR){IX};
        \tkzDefPoint(\LENGTH-\PHISTAR,\PHISTAR){IY};
        \tkzDefPoint(\LENGTH-\PHISTAR,\WIDTH-\PHISTAR){Y};
        \tkzLabelPoint[below](LL){$\phistar$};
        \tkzLabelPoint[below]([shift={(-0.5,0)}]RL){$2\pi-\phistar$};
        \tkzLabelPoint[left](LLL){0};
        \tkzLabelPoint[left](LLM){$\frac{\pi}{2}$};
        \tkzLabelPoint[left](LLU){$\pi$};
        \tkzLabelPoint[below](LLL){0};
        \tkzLabelPoint[below](RRL){$2\pi$};
        \tkzLabelPoint[right](IY){$x_*$};
        \tkzLabelPoint[left](IX){$y_*$};
        \tkzLabelPoint[right](Y){$I x_*$};
        \tkzLabelPoint[left](X){$I y_*$};
        
        \draw [thick] (LL) --(LU);
        \draw [thick] (LU) --(RU);
        \draw [thick] (RL) --(RU);
        \draw [thick] (LL) --(RL);
        \draw [thin,dashed](LLL)--(LL);
        \draw [thin,dashed](LLU)--(LU);
        \draw [thin,dashed](LLU)--(LLL);
        \draw [thin,dashed](RRL)--(RL);
        \draw [thin,dashed](RRU)--(RU);
        \draw [thin,dashed](RRU)--(RRL);
        \fill [blue, opacity=17/30](LL) -- (Y) -- (RU) -- (X) -- cycle;
        \fill [red, opacity=17/30](LU) -- (IY) -- (RL) -- (IX) -- cycle;
        \tkzDefPoint(-\LENGTH+\PHISTAR,\WIDTH/2-\PHISTAR/4-\PHISTAR/4){M1LL};
        \tkzDefPoint(-\LENGTH+\PHISTAR,\WIDTH/2-\PHISTAR/4+\PHISTAR/4){M1LU};
        \tkzDefLine[parallel=through M1LU](M1LL,RL)\tkzGetPoint{M11RU};
        \path[name path=lineM1U] (M11RU) -- (M1LU);
        \path[name path=lineM0L] (IX) -- (RL);
        \path[name intersections={of=lineM1U and lineM0L,by={M11RU}}];
        \fill [red, opacity=10/30](RL) -- (M1LL) -- (M1LU) -- (M11RU) -- cycle;
        \path[name path=RBoundry](RL) -- (RU);
        \tkzDefLine[parallel=through LU](M1LL,RL)\tkzGetPoint{UM11RU};
        \path[name path=lineUM1U](LU)--(UM11RU);
        \path[name intersections={of=lineUM1U and RBoundry,by={UM1RU}}];
        \tkzDefShiftPoint[UM1RU](0,-\PHISTAR/2){UM1RL};
        \tkzDefLine[parallel=through UM1RL](UM1RU,LU)\tkzGetPoint{UM11LL};
        \path[name path=lineLUIY](LU) -- (IY);
        \path[name path=UM11LLUM1RL](UM11LL)--(UM1RL);
        \path[name intersections={of=lineLUIY and UM11LLUM1RL,by={UM1LL}}];
        \fill [red, opacity=10/30](LU) -- (UM1RU) -- (UM1RL) -- (UM1LL) -- cycle;
        \tkzDefPoint(-\LENGTH+\PHISTAR,\WIDTH/3){L_LUM2};
        \tkzDefShiftPoint[L_LUM2](0,-\PHISTAR/3){L_LLM2};
        \tkzDefLine[parallel=through L_LUM2](L_LLM2,RL)\tkzGetPoint{Far_L_RUM2};
        \path[name path=line_L_LUM2_Far_L_RUM2] (L_LUM2)--(Far_L_RUM2);
        \path[name path=line_RL_M1LL] (RL)--(M1LL);
        \path[name intersections={of=line_L_LUM2_Far_L_RUM2 and line_RL_M1LL,by={L_RUM2}}];
        \fill [green, opacity=15/30](RL)--(L_LLM2)--(L_LUM2)--(L_RUM2) -- cycle;
        \tkzDefPoint(\LENGTH-\PHISTAR,2*\WIDTH/3){R_RLM2};
        \tkzDefShiftPoint[R_RLM2](0,\PHISTAR/3){R_RUM2};
        \tkzDefLine[parallel=through R_RLM2](R_RUM2,LU)\tkzGetPoint{Far_R_LLM2};
        \path[name path=line_R_RLM2_Far_R_LLM2] (R_RLM2)--(Far_R_LLM2);
        \path[name path=line_LU_UM1RU] (LU)--(UM1RU);
        \path[name intersections={of=line_R_RLM2_Far_R_LLM2 and line_LU_UM1RU,by={R_LLM2}}];
        \fill [green, opacity=15/30](LU)--(R_RUM2)--(R_RLM2)--(R_LLM2)-- cycle;

        \tkzDefPoint(\LENGTH-\PHISTAR,1*\WIDTH/3){MR_RLM2};
        \tkzDefShiftPoint[MR_RLM2](0,\PHISTAR/3){MR_RUM2};
        \path[name path=line_IY_LU] (IY)--(LU);
        \tkzDefLine[parallel=through MR_RLM2](R_RLM2,R_LLM2)\tkzGetPoint{far_MR_LLM2};
        \tkzDefLine[parallel=through MR_RUM2](R_RLM2,R_LLM2)\tkzGetPoint{far_MR_LUM2};
        \path[name path=line_MR_RLM2_far_MR_LLM2] (far_MR_LLM2)--(MR_RLM2);
        \path[name path=line_MR_RUM2_far_MR_LUM2] (far_MR_LUM2)--(MR_RUM2);
        \path[name path=line_LU_IY] (LU)--(IY);
        \path[name intersections={of=line_LU_IY and line_MR_RLM2_far_MR_LLM2,by={MR_LLM2}}];
        \path[name intersections={of=line_LU_IY and line_MR_RUM2_far_MR_LUM2,by={MR_LUM2}}];
        \fill [green, opacity=15/30](MR_RUM2)--(MR_RLM2)--(MR_LLM2)--(MR_LUM2)-- cycle;

        \tkzDefPoint(-\LENGTH+\PHISTAR,2*\WIDTH/3){ML_LUM2};
        \tkzDefShiftPoint[ML_LUM2](0,-\PHISTAR/3){ML_LLM2};
        \tkzDefLine[parallel=through ML_LUM2](R_LLM2,R_RLM2)\tkzGetPoint{far_ML_RUM2};
        \tkzDefLine[parallel=through ML_LLM2](R_LLM2,R_RLM2)\tkzGetPoint{far_ML_RLM2};
        \path[name path=line_ML_LUM2_far_ML_RUM2] (ML_LUM2)--(far_ML_RUM2);
        \path[name path=line_ML_LLM2_far_ML_RLM2] (ML_LLM2)--(far_ML_RLM2);
        \path[name path=line_RL_IX] (RL)--(IX);
        \path[name intersections={of=line_RL_IX and line_ML_LUM2_far_ML_RUM2,by={ML_RUM2}}];
        \path[name intersections={of=line_RL_IX and line_ML_LLM2_far_ML_RLM2,by={ML_RLM2}}];
        \fill [green, opacity=15/30](ML_RUM2)--(ML_RLM2)--(ML_LLM2)--(ML_LUM2)-- cycle;
    \begin{scope}
      \clip (LL) -- (X) -- (RU) -- cycle;
      \draw[fill=blue] circle[at=(X),radius=0.35];
    \end{scope}
    \begin{scope}
      \clip (RU) -- (Y) -- (LL) -- cycle;
      \draw[fill=blue] circle[at=(Y),radius=0.35];
    \end{scope}
    \begin{scope}
      \clip (RL) -- (IX) -- (LU) -- cycle;
      \draw[fill=red] circle[at=(IX),radius=0.35];
    \end{scope}
    \begin{scope}
      \clip (LU) -- (IY) -- (RL) -- cycle;
      \draw[fill=red] circle[at=(IY),radius=0.35];
    \end{scope}
    
    \node[] at (0,0.5*\WIDTH)  {$M^{\text{\upshape in}}_{r,0}$};
    
    \path[name path=line_RU_X] (RU)--(X);
    \path[name intersections={of=line_RU_X and lineM1U,by={Mr1Label}}];
    \node[] at ([shift={(0.07,-0.18)}]Mr1Label) {\small $M^{\text{\upshape in}}_{r,1}$};
    \node[] at ([shift={(2.1*\LENGTH/3,\WIDTH/3.27)}]Mr1Label) {\small $M^{\text{\upshape in}}_{r,1}$};

    \path[name intersections={of=line_RU_X and line_L_LUM2_Far_L_RUM2,by={Mr2Label}}];
    \node[] at ([shift={(0.14,-0.16)}]Mr2Label) {\small $M^{\text{\upshape in}}_{r,2}$};
    \node[] at ([shift={(3.2*\LENGTH/3,\WIDTH*0.48)}]Mr2Label) {\small $M^{\text{\upshape in}}_{r,2}$};
    \node[] at ([shift={(0.2,-0.05)}]X) {\tiny $N^{\text{\upshape in}}$};
    \node[] at ([shift={(-0.14,0.13)}]Y) {\tiny $N^{\text{\upshape in}}$};
    \node[] at ([shift={(0.25,0.1)}]IX) {\tiny $N^{\text{\upshape out}}$};
    \node[] at ([shift={(-0.08,-0.1)}]IY) {\tiny $N^{\text{\upshape out}}$};
    \end{tikzpicture}
    \caption{\cref{fig:Mr} with $\phistar<\pi/6$ and more detail. The two pink strips are $\mathcal{F}^{-1}(\Mrout\setminus\Mrin)$ and have slope $-1/4$. 
    The green strips are $\mathcal{F}^{-2}(\Mrout\setminus(\Mout_{r,0}\cup\Mout_{r,1}\cup\Mout_{r,2}))$ and have slope $-1/6$.\newline The intersections $\Min_{r,0}$, $\textcolor{blue}{M^{\text{\upshape in}}_{r,1}}$ and $\textcolor{blue}{M^{\text{\upshape in}}_{r,2}}$ are marked.  \newline The leftmost point of the closure of \textcolor{blue}{$M^{\text{\upshape in}}_{r,0}$} is $(\pi-\phistar,\frac{\pi}{2})$. The leftmost point of the closure of \textcolor{blue}{$M^{\text{\upshape in}}_{r,1}$} is $(\frac{2\pi}{3}-\phistar,\frac{\pi}{3})$. The leftmost point of the closure of \textcolor{blue}{$M^{\text{\upshape in}}_{r,2}$} is $(\frac{\pi}{2}-\phistar,\frac{\pi}{4})$\newline$\hat{M}$ is the union of the two pink strips and the center rhombus $\Mout_{r,0}=\Min_{r,0}=\Mrin\cap\Mrout$.\newline\Nin\ and \Nout\ are neighborhoods of $x_*,y_*,Ix_*,Iy_*$ respectively.}\label{fig:MrN}.
\end{center}
\end{figure}
\begin{definition}[Returns to $\hat M$]\label{def:MhatReturnOrbitSegment}
Write the orbit segment of $x\in\hat{M}$ to its first return to $\hat M$ as 
\begin{equation}\label{eqMhatOrbSeg}
        \underbracket{x}_{\in \hat{M}},\underbracket{\mathcal{F}(x),\cdots}_{\notin\hat M},\underbracket{\mathcal{F}^{\sigma(x)}(x)=\hat{F}(x)}_{\in \hat{M}},\text{ where }\sigma(x)\dfn\inf\{k>0\mid \mathcal{F}^{k}(x)\in\hat{M}\}.  
\end{equation}
By\COMMENT{It looks like none of what follows is a definition. Why not end the definition environment before "By"?} the definition of $\hat{M}$ (\cref{def:SectionSetsDefs}), we have either $x\in\mathcal{F}^{-1}(\Mrout\setminus\Mrin)$, in which case we set $x_0\dfn\mathcal{F}(x)\in\Mrout\cap\Mrin$, or $x\in\Mrout\cap\Mrin$, in which case we set $x_0\dfn x$. Hence, this orbit segment includes an $x_0\in\Mrout$ among its first two points. Then $x_1\dfn \mathcal{F}(x_0)\in \MRin$ and, taking $n_1\dfn m(x_1)$ from \cref{def:SectionSetsDefs}, $x_2\dfn\mathcal{F}^{n_1+1}(x_1)\in \Mrin$ (because  $\hat M\subset M_r$).

Then $m(x_2)$ determines whether $\hat{F}(x)=x_2$ or $\hat{F}(x)=\Fcal^{m(x_2)-1}(x_2)$ as follows.
    
Case i): $m(x_2)\le1$, that is $x_2\in(\textcolor{blue}{M^{\text{\upshape{in}}}_{r,0}}\cup\textcolor{blue}{M^{\text{\upshape in}}_{r,1}})\subset\hat{M}$. Then $x_2=\hat{F}(x)$.
    
Case ii): $m(x_2)\ge2$, that is $x_2\in{\Min_{r,n}}$ for some $n=m(x_2)\ge2$. Then $\mathcal{F}^{m(x_2)}(x_2)\in\Mrout\setminus\Mrin$ and $\mathcal{F}^{m(x_2)-1}(x_2)=\hat{F}(x)\in\mathcal{F}^{-1}(\Mrout\setminus\Mrin)\subset\hat{M}$ since $\mathcal{F}^i(x_2)\notin\mathcal{F}^{-1}(\Mrout\setminus\Mrin)$ for $i=0,\cdots,m(x_2)-2$.\COMMENT{This is not a definition. Maybe a remark? If so, combine with the next one?}
\end{definition}

\begin{remark}[``Bad regions'']\label{REMBadRegions}
    One inconspicuous feature of \cref{fig:Mr,fig:MrN} is worth attention. It turns out that the small regions \Nin\ and \Nout\ in \cref{def:SectionSetsDefs,fig:MrN} contain those points which are problematic with respect to the expansion we wish to establish (\cref{contraction_region}). Therefore we now prove, in effect, that they are as small as they look in \cref{fig:MrN} in the sense that they do not meet the other regions shown. Studying the dynamics on these regions entails significant effort and constrains the parameters for which we can prove uniform hyperbolicity. In effect, that work entails treating the lemon billiard as a perturbation of a 1-petal billiard in \cref{sec:SACTWFELBWB}.
\end{remark}\begin{proposition}[Location of bad regions]\label{Prop:Udisjoint}
    If $\phistar\in\big(0,\tan^{-1}(1/3)\big)$ and $R$ satisfies \eqref{eqJZHypCond}, then $\Nin\cap\big(\textcolor{blue}{M^{\text{\upshape in}}_{r,0}}\cup{\textcolor{blue}{M^{\text{\upshape in}}_{r,1}}}\cup{\textcolor{blue}{M^{\text{\upshape in}}_{r,2}}}\big)=\emptyset$ and $\Nout\cap\big(\textcolor{red}{M^{\text{\upshape in}}_{r,0}}\cup\textcolor{red}{M^{\text{\upshape in}}_{r,1}}\cup\textcolor{red}{M^{\text{\upshape in}}_{r,2}}\big)=\emptyset$.
\end{proposition}
\begin{proof}
$\Nin$ has two components shown in \cref{fig:MrN}, a left component near $(\phistar,\phistar)$, and a right component near $(2\pi-\phistar,\pi-\phistar)$. They are symmetric with respect to the center of $M_r$. ($I \circ J$ in \cref{DefSymmetries}).
    
The green strips in \cref{fig:MrN} have slope $-1/6$ and the blue strip  \Mrin\  has slope $1/2$. \textcolor{blue}{$M^{\text{\upshape in}}_{r,2}$} has two components. The left component is the interior of the parallelogram with vertices $(\pi/2-\phistar,\pi/4)$, $(\pi/2-\phistar+\phistar/2,\pi/4+\phistar/4)$, $(\pi/2+\phistar,\pi/4)$ and $(\pi/2+\phistar-\phistar/2,\pi/4-\phistar/4)$ among which $(\pi/2-\phistar,\pi/4)$ is the leftmost. The left and right components are symmetric with respect to the center of $M_r$ (i.e., under $I \circ J$, see \cref{DefSymmetries}). It suffices to show that the left components of \textcolor{blue}{$M^{\text{\upshape in}}_{r,2}$} and $\Nin$ are disjoint.\COMMENT{State somewhere the symmetries in \cref{fig:MrN} \newline\wentao{More details are added}} \COMMENT{How?\newline\wentao{More details are added.}\boris More to do. Find a place for a remark about symmetry of the billiard table and reversibility and what symmetries this produces in the phase space. Make sure that it is written carefully and clearly; it makes sense to connect this (in a clear way!) with the symmetries in \cref{fig:MrN}. Then refer to it wherever either symmetry is invoked\wentao{More details referred and added in \cref{DefSymmetries}}.\newline Here, the second issue is the proof that the left component of \textcolor{blue}{$M^{\text{\upshape in}}_{r,2}$} is disjoint with the left component of $\Nin$.\wentao{proof much changed now}}\NEWPAGE

Note that in \cref{fig:MrN} by elementary geometry, the Euclidean distance\COMMENT{Euclidean distance in $\mathbb R^2$?\wentao{Yes}} between $(\pi/2-\phistar,\pi/4)$ and $(\phistar,\phistar)$ is the smallest Euclidean distance from any point in the closure of the left component of $\textcolor{blue}{M^{\text{\upshape in}}_{r,2}}$ to $(\phistar,\phistar)=Iy_*$, and this smallest Euclidean distance is $\sqrt{5}(\nicefrac{\pi}{4}-\phistar)$ since $0<\phistar<\tan^{-1}(1/3)<\pi/4$.\COMMENT{What does this mean????? Do you want to say that $\phistar\in\big(0,\tan^{-1}(1/3)\big)$ and $\tan^{-1}(1/3)<\pi/4$? Maybe say $\phistar\in\big(0,\tan^{-1}(1/3)\big)$ and $\pi/4-\tan^{-1}1/3=\tan^{-1}1/2>.46$} We get the following.\COMMENT{Explain the estimate of $\cos^2\phistar/2$ (why 500??), maybe mention $\pi/4-\tan^{-1}1/3=\tan^{-1}1/2>.46$\wentao{thanks, more estimate is explained now.}} 
        \[
R\overset{\eqref{eqJZHypCond}}>\frac{30000r}{\sin^2\phistar}=\frac{30000r}{4\sin^2(\phistar/2)\cos^2(\phistar/2)}>\frac{500r}{\sin^2(\phistar/2)}\overbracket{>}^{\mathclap{\substack{0<\phistar<\tan^{-1}(1/3)\Rightarrow \pi/4-\tan^{-1}1/3=\tan^{-1}1/2>0.46\Rightarrow \frac{289}{5(\pi/4-\phistar)^2}<\frac{289}{5(\pi/4-\tan^{-1}(1/3))^2}<500\\\vphantom.}}}\frac{289r}{5(\pi/4-\phistar)^2\sin^2{(\phistar/2)}},
        \]
    which implies $\sqrt{5}(\nicefrac{\pi}{4}-\phistar)>\delta_*=17\sqrt{r/R}/\sin{(\phistar/2)}$. Hence, the left components of $\textcolor{blue}{M^{\text{\upshape in}}_{r,2}}$ and \Nin are disjoint. Thus $\Nin\cap\textcolor{blue}{M^{\text{\upshape in}}_{r,2}}=\emptyset$. 
    
Similar arguments about the distance between the leftmost points of the closure of $\textcolor{blue}{M^{\text{\upshape in}}_{r,0}}$ and $I y_*$ and the distance between the leftmost points of the closure of $\textcolor{blue}{M^{\text{\upshape in}}_{r,1}}$ and $I y_*$ yield $\Nin\cap\textcolor{blue}{M^{\text{\upshape in}}_{r,0}}=\emptyset$ and $\Nin\cap\textcolor{blue}{M^{\text{\upshape in}}_{r,1}}=\emptyset$. Thus, we have $\Nin\cap\big(\textcolor{blue}{M^{\text{\upshape in}}_{r,0}}\cup\textcolor{blue}{M^{\text{\upshape in}}_{r,1}}\cup\textcolor{blue}{M^{\text{\upshape in}}_{r,2}}\big)=\emptyset$, and, by symmetry $I \circ J$ in \cref{DefSymmetries},  $\Nout\cap\big(\textcolor{red}{M^{\text{\upshape out}}_{r,0}}\cup\textcolor{red}{M^{\text{\upshape out}}_{r,1}}\cup\textcolor{red}{M^{\text{\upshape out}}_{r,2}}\big)=\emptyset$.
\end{proof}
\begin{remark}\label{remark:MoreUdisjoint}
    In fact, using the same proof we see that if $0<\phistar<\tan^{-1}\big(1/3\big)$ with $R>\frac{30000r}{\sin^2\phistar}$, then $\Nin\cap{\textcolor{blue}{M^{\text{\upshape in}}_{r,3}}}$, $\Nin\cap{\textcolor{blue}{M^{\text{\upshape in}}_{r,4}}}$, $\Nout\cap\textcolor{red}{M^{\text{\upshape out}}_{r,3}}$, $\Nout\cap\textcolor{red}{M^{\text{\upshape out}}_{r,4}}$ are also empty.
\end{remark}
\COMMENT{I am trying to explain to the reader what $\hat M$ is, rather than just defining it. We should discuss two things soon.\newline 0) Is there a better way of describing $\hat M$?\newline 1) Why is it preferable to prove \cref{RMASUH} and derive \cref{thm:MroutreturnUniformExpansion} from it, rather than proving \cref{thm:MroutreturnUniformExpansion} directly?\newline 2) How does \cref{RMASUH} imply \cref{thm:MroutreturnUniformExpansion}?}

\begin{remark}
Here are a few places where we use that $R$ is large. In the proof of \cref{Proposition:DF4ConeUpperBound} (around page \pageref{proposition:dx3dx2expansion}), we use $R>\frac{1758r}{\sin^2\phistar}$.\COMMENT{Add more as we think of it\newline\wentao{Agree. Work In Progress}}
\end{remark}

Unless\COMMENT{Rethink this subsection title\wentao{Main theorems now.}} stated otherwise, we henceforth assume $0<\phistar<\tan^{-1}{(1/3)}$ and \eqref{eqJZHypCond}.
The statement of \cref{thm:MroutreturnUniformExpansion} uses notions from \cref{def:SectionSetsDefs} below, which includes further notation needed for the companion result promised in \cref{rmk:othersectionUH}, \cref{RMASUH} below.

[Here we provide the proof sketch. More details are given/reiterated in \cref{subsec:UniformExpansionOnSection,RMASUHwithProof}]

\begin{remark}[Uniform hyperbolicity of return map on other sections]\label{rmk:othersectionUH}
\cref{thm:UniformExpansionOnMhat1,Col:EuclideanExpansionMhat1} give the uniform hyperbolicity of the return map $\hat{F}_1$ on $\hat{M}_1$. \cref{thm:MroutreturnUniformExpansion,corollary:MroutreturnUniformExpansion} give the uniform hyperbolicity of the return map $\tilde{F}$ on $\Mrout$.\COMMENT{Think where this belongs.}
\end{remark}
\begin{remark}
\strut\COMMENT{\textcolor{red}{Expand here on \cref{rmk:othersectionUH} with an explanation from the comments about technical issues about applying them.} Here is the part of the comments that might be adaptable (at this point it is hard to understand):\newline The issue of \cref{thm:MroutreturnUniformExpansion} UH is the unstable cone is half quadrant cone slopes in [0,1] but its stable cone is not in the flipped slopes in [-1,0]. We can check the singularity curves, $\Scal_{-1}$ is in the unstable cone. But $\Scal_{1}$ is not. That's why JZ21 needs a shear to have the cone. And we have good unstable/stable cones in the many middle sections so that we have alignment for both positive/negative orders singularity curves on middle sections but only negative orders singularity curves alignment for \Mrout. We need both singularity alignments to prove local ergodicity SC4 double singularity condition and SC5 transversal condition.\wentao{Let's discuss. It may not be needed. I added some explanation in the remark.}}
Although sections $\hat{M}$ in \cref{def:SectionSetsDefs}\eqref{def:Mhat} and $\hat{M}_1$ in \cref{def:symetricsectionMhatPlusMinus} are fragmental in \cref{fig:MrN,fig:hatUplusUminus}, in \cref{thm:UniformExpansionOnMhat,thm:UniformExpansionOnMhat1} the unstable invariant cone on those sections is inside the positive quadrant cone with associated quadratic norm being nonnegative (see \cref{def:ExpansionInCone}). This is convenient for our local ergodicity proof.

Choosing reference sections $\Mrout$ also works fine for local ergodicity proof and we also prove \cref{thm:MroutreturnUniformExpansion,corollary:MroutreturnUniformExpansion} the uniform hyperbolicity for the return map on $\Mrout$. In this case, with the stable cone on $\Nout$ is larger than the quadrant cone, we have to exclude those ``bad'' regions inside $\Nout$ in the verification for L4, L3 local ergodicity conditions in \cite{MR3082535}. We alternatively pick the open reference sets $\mathfrak{U}_{-}$, $\mathfrak{U}_{+}$ defined in \cref{def:UIntRegionofD} with $\Mrout\setminus\Nout\subset\mathfrak{U}_{-}\subset\Mrout$ and $\Mrin\setminus\Nin\subset\mathfrak{U}_{+}\subset\Mrin$. 

Because they come with nicer cones, we still use the more fragmented sections. This choice will not cause difficulties in the proof of global ergodicity so long as the set of double-order singularities is countable \cite[Proposition 6.19] {cb} (Also See \cite[Theorem 6.20]{cb}). 
\end{remark}
\NEWPAGE

\section[Expansion estimates]{Expansion estimates}\label{sec:ECE}
\label{sec:FPCFICBMD}
This section states lower bounds for the differential of the billiard map in a case-by-case fashion (\cref{eqMainCases,TEcase_a0,TEcase_a1,TEcase_b,TEcase_c}) and establishes that those cases in which these lower bounds do not produce expansion are confined to a small portion of the section (\cref{contraction_region}). This is sufficient information to show strict invariance of the positive quadrant (\cref{corollary:InvQuadrant}; see also \cref{corollary:halfquadrantcone}) and strict expansion in it upon iteration (\cref{def:ApettrajectoryExpansion}).\COMMENT{This may need editing since our emphasis has now shifted ot half-quadrants\wentao{Let's discuss, Actually on different section, the chosen cone is different.We may say the cones in \cref{thm:UniformExpansionOnMhat,thm:UniformExpansionOnMhat1,thm:MroutreturnUniformExpansion} etc.}} The proofs of these items in \cref{sec:ECE} are carried out in later sections.

Toward establishing uniform expansion in invariant cones, a core step is to estimate the expansion obtained in the transitions from one of the boundary arcs to the other. This turns out to require a case-by-case approach, which already occurred when the existence of an invariant cone family was established in \cite{hoalb}. We begin with a presentation of the main cases and the notation and prior results about invariant cones.
\subsection{Invariant cones}
\label{subsec:pclf}
For nonsingular\COMMENT{Needed?\wentao{When the section we state/conclude arbitrarily return map differential/expansion, we should assume $x$ nonsingular}\boris Can you explain this sentence? Does it answer the question?\wentao{\small\color{red} We should a discussion here for what is nonsingular first order or higher order, for billiard map $\Fc$ or for the return map $\hat{F}$. We need different tailored/meanings in different propositions for different section return map.\Large This seems important}} $x\in\hat{M}$, in the $\phi\theta$-coordinates from \cref{def:BasicNotations0}, consider an orbit segment \eqref{eqMhatOrbSeg} (notations as in \cref{def:BasicNotations,def:SectionSetsDefs,def:MhatReturnOrbitSegment}) and the points in it found right after \eqref{eqMhatOrbSeg}:
\begin{equation}\label{eqPtsFromeqMhatOrbSeg}
    x_{0}=(\phi_0,\theta_0)\in\Mrout,\ x_1=(\Phi_1,\theta_1)\dfn \mathcal{F}(x_0)\in \MRin
    \text{ and }x_{2}=(\phi_2,\theta_2)=\mathcal{F}^{n_1+1}(x_1)\in\Mrin\text{ with }n_1\ge0.
\end{equation}
\begin{figure}[h]
\begin{center}
\begin{sideways}
\begin{tikzpicture}[xscale=0.53, yscale=0.53]
    \tkzDefPoint(0,0){Or};
    \tkzDefPoint(0,-4.8){Y};
    \pgfmathsetmacro{\rradius}{4.5};
    \pgfmathsetmacro{\bdist}{8.7};
    \tkzDefPoint(0,\bdist){OR};
    \pgfmathsetmacro{\phistardeg}{40};
    \pgfmathsetmacro{\XValueArc}{\rradius*sin(\phistardeg)};
    \pgfmathsetmacro{\YValueArc}{\rradius*cos(\phistardeg)};
    \clip(-4.9,-13) rectangle (6.6, 11.5);
    \tkzDefPoint(\XValueArc,-1.0*\YValueArc){B};
    \tkzDefPoint(-1.0*\XValueArc,-1.0*\YValueArc){A};
    \pgfmathsetmacro{\Rradius}{veclen(\XValueArc,\YValueArc+\bdist)};
    \draw[name path=Cr,dashed,ultra thin] (Or) circle (\rradius);
    \tkzDrawArc[thick](Or,B)(A);
    \tkzDrawArc[name path=CR, thick](OR,A)(B);
    \tkzDrawArc[blue, dashed, ultra thin](OR,B)(A);
   
    \pgfmathsetmacro{\PXValue}{\Rradius*sin(7)};
    \pgfmathsetmacro{\PYValue}{\bdist-\Rradius*cos(7)};
    \tkzDefPoint(\PXValue,\PYValue){P};
    \draw[name path=PQ_Far,blue,ultra thin](P)--+(70:7)coordinate(Q_Far);
    \draw[name path=PT_Far,blue,ultra thin](P)--+(124:9);
    \draw[name path=PT0_Far,blue,dashed,ultra thin](P)--+(-56:7);
    \draw[name path=PQ0_Far,blue,dashed,ultra thin](P)--+(-110:9);
    \draw[red,dashed,->](P)--+(7:4)coordinate(PposX);
    \draw[name path=PPnegX,red,dashed](P)--+(-173:6)coordinate(PnegX);
    \draw[name path=linePOR,red,dashed,->](P)--(OR);
     \path[name path= linePOr] (P)--(Or);
    \tkzDefPoint(0,-5.5){Y};
    \tkzDefPoint(0,\bdist-\Rradius){C};
    \draw[red,dashed](OR)--(C);
    \tkzMarkAngle[arc=lll,ultra thin,size=1.7](Or,OR,P);
    \path[name intersections={of=Cr and PQ_Far,by={Q}}];
    \draw[,blue,ultra thin,dashed](Q)--+(70:7);
    \path[name intersections={of=Cr and PQ0_Far,by={Q0}}];
    \path[name intersections={of=Cr and PT0_Far,by={T0}}];
    \path[name intersections={of=Cr and PT_Far,by={T}}];
    \tkzMarkAngle[arc=lll,very thin,size=3.2](PposX,P,Q);
    \tkzMarkRightAngle(OR,P,PnegX);
    \node[rotate=-90] at ($(OR)+(-85:2.3)$)  {$\Phi_1$};
    \node[rotate=-90] at ($(P)+(30:3.5)$)  {\Large {$\theta_1$}};

    \node[rotate=-90] at ([shift={(0.4,2.2)}]PposX) {\textcolor{blue}{\small circle $C_R$}};
    
    \tkzLabelPoint[left,rotate=-90](A){$A$};
    \tkzLabelPoint[left,rotate=-90](B){$B$};
    \tkzLabelPoint[left,rotate=-90](OR){$O_R$};
    \tkzLabelPoint[below,rotate=-90](Or){$O_r$};
    \tkzLabelPoint[above,rotate=-90]([shift={(0,0.35)}]P){$P$};
    \tkzLabelPoint[above,rotate=-90](Q){$Q$};
    \tkzLabelPoint[left,rotate=-90](T){$T$};
    \tkzLabelPoint[below right,rotate=-90](Q0){$Q_0$};
    \tkzLabelPoint[right,rotate=-90](T0){$T_0$};
    \tkzLabelPoint[right,rotate=-90](PposX){$TP$};
    \tkzDefLine[perpendicular=through Or](P,PnegX)\tkzGetPoint{M_far};
    \path[name path= lineOrMfar] (Or)--(M_far);
    \path[name path= PPnegX] (P)--(PnegX);
    \path[name path= PPposX] (P)--(PposX);
    \path[name intersections={of=lineOrMfar and PPnegX,by={M}}];
    \path[name intersections={of=Cr and PPnegX,by={A0}}];
    \path[name intersections={of=Cr and PPposX,by={B0}}];
    \coordinate (Or_mirror) at ($(Or)!2!(M)$);
    \tkzDrawArc[name path=mirrorCr, ultra thin, dashed](Or_mirror,A0)(B0);
    \tkzLabelPoint[below,rotate=-90](Or_mirror){$O'_r$};

\path[name intersections={of=mirrorCr and PQ0_Far,by={Q1}}];
    \path[name intersections={of=mirrorCr and PT0_Far,by={T1}}];
    \tkzLabelPoint[above,rotate=-90](Q1){$Q_1$};
    \tkzLabelPoint[above,rotate=-90](T1){$T_1$};
    \tkzDrawPoints(A,Or,OR,P,B,T,Q,Or_mirror,A0,T0,Q0,T1,Q1);
    \begin{scope}[ultra thin,decoration={markings, mark=at position 0.75 with {\arrow{>}}}
    ] 
    \draw[blue,postaction={decorate}](P)--(Q);
    \draw[blue,postaction={decorate}](T)--(P);
    \end{scope}
    \tkzLabelPoint[above,rotate=-90]([shift={(0.15*\rradius,1.1*\rradius)}]Or){\Large $\Gamma_r$};
    \tkzLabelPoint[above,rotate=-90]([shift={(0.35*\rradius,-1.05*\rradius)}]Or_mirror){\Large $\Gamma_r'$};
\end{tikzpicture}
\end{sideways}
\captionsetup{type=figure}
\caption{$n_1=0$ case: the trajectory hitting $P=p(x_1)=p(\Fcal^{-1}(x_2))$ with $T=p(\Fcal^{-1}(x_1))$ and $Q=p(x_2)$. $\tau_0=|TP|$, $\tau_1=|PQ|$. With $O'_r$ is symmetric to $O_r$ and arc $\Gamma'_r$ is symmetric to arc $\Gamma_r$ with respect to the line $P-TP$ which is the tangent line to $\Gamma_R$ at $P$.\newline 
Extend line segment $TP$ to meet the circle $C_r$ at $T_0$ and the arc $\Gamma'_r$ at $T_1$. 
\newline Extend line segment $QP$ to meet the circle $C_r$ at $Q_0$ and the arc $\Gamma'_r$ at $Q_1$.
}\label{fig:lengthfunctonFigure}
\end{center}
\end{figure}
\NEWPAGE\begin{notation}\label{def:fpclf}
For $x_0\in\Mrout$, $x_1\in\MRin$ and $x_2\in\Mrin$, the free-path and chord-length functions (see \cref{fig:lengthfunctonFigure}) are:\COMMENT{The caption text of \cref{fig:heuristicproof} needs editing; it is hard to understand (what is ``mirror table''? ``mirrored billiard''?). It does not show any of these functions.\wentao{edited}\newline Also, in that figure I replaced green by green!50!black for visibility. This might help in other figures as well.\wentao{edited to use black.}\Hrule\textcolor{red}{It would really help to have a figure near here which shows ONLY these quantities in order to illustrate this definition. \cref{fig:heuristicproof} does not help.}\wentao{added \cref{fig:lengthfunctonFigure} with more explanations}}
\begin{itemize}
    \item $\tau_0$ is the distance between $p(x_0)\in\Gamma_r$ and $p(\mathcal{F}(x_0))\in\Gamma_R$.
    \item $\tau_1$ is the distance between $p(x_2)\in\Gamma_r$ and $p(\mathcal{F}^{-1}(x_2))\in\Gamma_R$, or $L_1\dfn\tau_1$ to denote length.
    \item $d_0\dfn r\sin{\theta_0}$ is the half-length of the chord that connects the two intercepts of $C_r$ with the line that goes through $p(x_0)$ and $p(\mathcal{F}(x_0))$, and in \cref{fig:lengthfunctonFigure} $|TT_0|$ is such a chord\COMMENT{Incomprehensible\wentao{More illustrations added. Let's discuss.}}.
    \item $L_0\dfn2d_0-\tau_0$ and in \cref{fig:lengthfunctonFigure}, $|PT_0|=L_0$.
    \item $d_1\dfn R\sin{\theta_1}$ is the half-length of the chord that connects the two intercepts of $C_R$ with the line that goes through $p(x_1)$ to $p(\mathcal{F}(x_1))$.\COMMENT{Incomprehensible\wentao{More illustrations added. Let's discuss.}} In \cref{fig:lengthfunctonFigure}, the line $PQ$ has two intersections with $C_R$. $P=p(x_1)$ is one of the two intersections.
    \item $d_2\dfn r\sin{\theta_2}$ is the half-length of the chord that connects the two intercepts of $C_r$ with the line that goes through $p(x_2)$ and $p(\mathcal{F}^{-1}(x_2))$, and in \cref{fig:lengthfunctonFigure} $|QQ_0|$ is such a chord.\COMMENT{Incomprehensible\wentao{More illustrations added. Let's discuss.}}
\end{itemize}
\end{notation}
\begin{remark}\label{remark:deffpclf}
$d_0,d_1,d_2,\tau_0,\tau_1$ are continuous functions of nonsingular $x_1\in\MRin\subset\Int M_R$. In \cref{subsec:continuousextension,subsec:A1contractionCase,sec:ProofTECase_c}, they are continuously extended to the boundary of $M_R$. Because any free path ( i.e., a segment between collisions) is contained in the disk surrounded by $C_r$ and also inside the disk surrounded by $C_R$, we have the following.
\[2d_2>\tau_1=L_1,\quad 2d_0-\tau_0=L_0>0,\quad2d_1>\tau_0, \quad2d_1>\tau_1.
\]
\end{remark}
\begin{remark}[\cref{fig:lengthfunctonFigure}]\label{remark:symmetrywiththetangentline}
If $n_1=0$, let $P=p(x_1)=p(\Fcal^{-1}(x_2))$ as in \cref{fig:lengthfunctonFigure}, then with $|TP|=\tau_0$\COMMENT{Does not match definition}, $|PQ|=\tau_1$, $|TT_0|=2d_0$ since $|PQ|$ is symmetric to $|PT_1|$ with respect to the line $P-TP$\COMMENT{What does this notation mean?}, we see $\tau_0+\tau_1-2d_0=|TP|+|PQ|-|TT_0|=|TP|+|PT_1|-|TT_0|=|TT_1|-|T_0T_1|=|T_0T_1|>0$. Similarly, by symmetry $|Q_1P|=|TP|=\tau_0$ and $|QQ_0|=2d_2$, we have $\tau_0+\tau_1-2d_2=|TP|+|PQ|-|QQ_0|=|Q_1P|+|PQ|-|QQ_0|=|QQ_1|-|Q_0Q_1|=|Q_0Q_1|>0$.\COMMENT{This remark needs checking} Therefore, we have the following inequalities.
\[
    \tau_0+\tau_1-2d_0>0,\quad \tau_0+\tau_1-2d_2>0,\quad \tau_0+\tau_1-d_0-d_2>0.
\]
\end{remark}

Our case-by-case approach is related to different features of the derivative matrix of the return map---the good scenarios are those when all entries are positive or all are negative; in both cases, the positive quadrant is preserved. 
We introduce these cases now.
\begin{equation}\label{eqMainCases}
    \begin{cases}\text{Case (a) in \eqref{eqPtsFromeqMhatOrbSeg}:}&d_1=R\sin{\theta_1}\ge 2r,\quad\begin{cases}
        \text{Case (a0):}&\sin{(\theta_1)}\ge\sqrt{{4r}/{R}},\\\text{Case (a1):}&\sin{(\theta_1)}\in [{2r}/{R},\sqrt{{4r}/{R}}).
    \end{cases}\\
    \text{Case (b) in \eqref{eqPtsFromeqMhatOrbSeg}:}&d_1=R\sin{\theta_1}< 2r\text{ and }n_1=0\text{, i.e., the orbit segment has exactly one collision on }\Gamma_R,\\
    \text{Case (c) in \eqref{eqPtsFromeqMhatOrbSeg}:}&d_1=R\sin{\theta_1}< 2r\text{ and }n_1\ge1\text{, i.e., the orbit segment has $n_1+1\ge 2$ collisions on }\Gamma_R.
    \end{cases}
\end{equation}
\cref{JZProp3.6a,JZProp3.6b,JZProp3.6c} below (contained in \cite[Proposition 3.6]{hoalb}) provide invariance of the first-quadrant cone family. They assume that $R$ satisfies \eqref{eqJZHypCond} and apply to the points in \eqref{eqPtsFromeqMhatOrbSeg} from \eqref{eqMhatOrbSeg}.
\begin{proposition}[{\cite[Proposition 3.6]{hoalb}}]\label{JZProp3.6a} In Case {(a)}\label{caseA}, $x_2=\mathcal{F}(x_1)\in\Mrin$, $D\mathcal{F}^2_{x_0}$ is a 2 by 2 matrix with entries depending on $\tau_0,\tau_1,d_0,d_1,d_2$, and all 4 entries are negative \cite[equations (4.2) through (4.5)]{hoalb}.
\end{proposition}
\begin{proposition}[{\cite[Proposition 3.6]{hoalb}}]\label{JZProp3.6b} In Case {(b)}\label{caseB}, $x_2=\mathcal{F}(x_1)\in\Mrin$, $x_{-1}:=\mathcal{F}^{-1}(x_0)\in M_r$, $x_3:=\mathcal{F}(x_2)\in M_r$, $D\mathcal{F}^4_{x_{-1}}$ is a 2 by 2 matrix with entries depending on $\tau_0,\tau_1,d_0,d_1,d_2$, and all 4 entries are negative \cite[Equations (4.10) through (4.19)]{hoalb}.
\end{proposition}
\begin{proposition}[{\cite[Proposition 3.6]{hoalb}}]\label{JZProp3.6c} In Case {(c)}\label{caseC},  $x_3\dfn \mathcal{F}(x_2)\in M_r$, $x_{-1}\dfn \mathcal{F}^{-1}(x_0)\in M_r$, $D\mathcal{F}_{x_{-1}}^{n_1+4}$ is a 2 by 2 matrix with entries depending on $\tau_0,\tau_1,d_0,d_1,d_2,n_1$, and  all 4 entries are positive in this case \cite[Section 5]{hoalb}.
\end{proposition}
\begin{remark}
The calculations of the entries of the derivative matrices in these three cases can be found in \cite[Sections 3--5]{hoalb} except that in subcase (a0) \cref{TEcase_a0}, we need to extend the definition of $D\mathcal{F}_{x_0}^2$ to the boundary of $M_R$ (\cref{subsec:continuousextension}).
\end{remark}
Unlike the invariance of cones, expansion of vectors in cones (the main novelty of this paper) is quite difficult outside Case~(a0) in \eqref{eqMainCases} (see \cref{TEcase_a0,TEcase_a1,TEcase_b,TEcase_c} below). Let us note here that those harder cases are confined to a small part of the phase space:
\begin{proposition}[Contraction region]\label{contraction_region}
Cases {(a1) (\Cref{TEcase_a1}), (b) (\Cref{TEcase_b}), and  (c) (\Cref{TEcase_c}}) in \eqref{eqMainCases} arise only when $x_0\in\Nout$ and $x_2\in\Nin$ (see \cref{def:SectionSetsDefs}).
\end{proposition}
The proof is in \Cref{subsec:ProofContractionRegion}.

\subsection{Control of \texorpdfstring{$D\hat F$}{Control of return map differential}}
We now proceed to the central object of this article, control of $D\hat F$---which is clearly crucial for \cref{RMASUH}. We will eventually establish expansion of the p-metric, beginning with \cref{TEcase_a0,TEcase_a1,TEcase_b,TEcase_c}, one proposition for each of the four cases in \eqref{eqMainCases}.

{\bfseries Standing assumptions for \cref{TEcase_a0,TEcase_a1,TEcase_b,TEcase_c}.} Fix $r,R$ satisfying the hyperbolicity condition \eqref{eqJZHypCond}, fix $\phistar\in(0,\tan^{-1}{(1/3)})$, $\hat F$ from \cref{def:SectionSetsDefs}, and the cone family $C_x$ from \cref{RMASUH}.

\begin{proposition}[Case (a0)]\label{TEcase_a0} 
$\exists$ a constant $\lambda_c>0$ such that 
if $\sin{\theta_1}\ge\sqrt{{4r}/{R}}$ for the points in \eqref{eqPtsFromeqMhatOrbSeg} with nonsingular $x$, then
\begin{enumerate}
\item\label{itemTEcase_a0-1} the cone family $C_x$ is strictly invariant under $D\hat{F}$, i.e., $D\hat{F}_{x}(C_x)\subset\text{\big\{interior of }C_{\hat{F}(x)}\text{\big\}}$,
\item\label{itemTEcase_a0-2} $\frac{\|D\hat{F}_{x}(dx)\|_\p}{\|dx\|_\p}>1+\lambda_c$ for all $dx\in C_x$ and $\lambda_c>0$ is a constant determined by $r$, $R$, $\phistar$ configuration.
\end{enumerate}
\end{proposition}
\begin{proof}
    \eqref{itemTEcase_a0-1},\eqref{itemTEcase_a0-2} are the conclusions of \cref{theorem:UniformExpansionA0} below (see page \pageref{theorem:UniformExpansionA0} and the proof in \cref{sec:ProofTECase_a0_a1}).
\end{proof}
It is the kind of statement one would wish to prove in all subsequent cases as well, but instead, the corresponding propositions have lower bounds less than 1, necessitating subsequent work. Here they are. 
\begin{proposition}[Case (a1)]\label{TEcase_a1} 
If $\sin{\theta_1}\in [{2r}/{R},\sqrt{{4r}/{R}})$ for the points in \eqref{eqPtsFromeqMhatOrbSeg} with nonsingular $x$, then
\begin{enumerate}
    \item\label{itemTEcase_a1-1} the cone family $C_x$ is strictly invariant under $D\hat{F}$, i.e., $D\hat{F}_{x}(C_x)\subset\text{\big\{interior of }C_{\hat{F}(x)}\text{\big\}}$,
    \item\label{itemTEcase_a1-2} $\frac{\|D\hat{F}_{x}(dx)\|_\p}{\|dx\|_\p}>0.26$, for all $dx\in C_x$.
\end{enumerate}
\end{proposition}
\begin{proof}
\eqref{itemTEcase_a1-1}\eqref{itemTEcase_a1-2} are the conclusions of \cref{theorem:contractioncontrolA1} below (see page \pageref{theorem:contractioncontrolA1} and its proof in \cref{sec:ProofTECase_a0_a1}).
\end{proof} Note that this does not establish expansion; it controls how far from expansion the return map is in this case. Subsequent iterates are needed to convert this control into actual expansion. The same goes for the subsequent cases; the next case is the most daunting because the lower bound is quite small. For ease of exposition, we will take the liberty to refer to these lower bounds as lower bounds for expansion even when they are less than 1.\COMMENT{\color{red}With this having been said, almost all instances of "expansion/contraction" throughout the paper can now be replaced by "expansion." Even some occurrences of "contraction" could be changed to "expansion"---with great care!}
\begin{proposition}[Case (b)]\label{TEcase_b}
If 
$\sin{\theta_1}<{2r}/{R}$ and $n_1=0$ for the points in \eqref{eqPtsFromeqMhatOrbSeg}, then\COMMENT{Each item below is \textbackslash labeled. \textcolor{red}{At the start of  \cref{sec:ProofTECase_b} there should be a list of explicit statements about which of these five items is proved where.}}
\begin{enumerate}
    \item\label{itemTEcase_b1} the cone family $C_x$ is strictly invariant under $D\hat{F}$, i.e., $D\hat{F}_{x}(C_x)\subset\text{\big\{interior of }C_{\hat{F}(x)}\text{\big\}}$,
    \item\label{itemTEcase_b2} $\frac{\|D\hat{F}_{x}(dx)\|_\p}{\|dx\|_\p}>0.05$ for all $dx\in C_x$.
\end{enumerate}
\end{proposition}
\begin{proof}
\eqref{itemTEcase_b1} is the conclusion of \cref{corollary:halfquadrantcone}. \eqref{itemTEcase_b2} is the conclusion of \cref{theorem:proofofTEcase_b} (see its proof in \cref{sec:ProofTECase_b}).
\end{proof}
\begin{proposition}[Case (c)]\label{TEcase_c}
If $\sin{\theta_1}<{2r}/{R}$ and $n_1\ge1$ for the points in \eqref{eqPtsFromeqMhatOrbSeg}, then\COMMENT{Make sure \cref{sec:ProofTECase_c} EXPLICITLY proves each of the items below!}
\begin{enumerate}
    \item\label{itemTEcase_c1} the cone family $C_x$ is strictly invariant under $D\hat{F}$, i.e., $D\hat{F}(C_x)\subset\text{\big\{interior of }C_{\hat{F}(x)}\text{\big\}}$.
    \item\label{itemTEcase_c2} $\frac{\|D\hat{F}_{x}(dx)\|_\p}{\|dx\|_\p}>1-1743\big(r/R\big)-15450\big(r/R\big)^2>0.9$, for all $dx\in C_x$.
    \item\label{itemTEcase_c3} when $n_1\ge4$, $\frac{\|D\hat{F}(dx)\|_\p}{\|dx\|_\p}>\frac{0.96}{1.01}n_1^2+5.655n_1-8.628>29.1$, for all $dx\in C_x$. See \cref{subsec:expansionLarge}.
\end{enumerate}
\end{proposition}
The proof is in \cref{sec:ProofTECase_c}.
\subsection{Expansion for an iterate}
The purpose of this subsection is to obtain from the preceding results a strictly invariant cone family on $\hat M$ (\cref{corollary:InvQuadrant}) and then state the uniform expansion implied by  \Cref{TEcase_a0,TEcase_a1,TEcase_b,TEcase_c,contraction_region} and three approximation theorems in  \Cref{sec:SACTWFELBWB}.

As noted above, collisions between transitions are often needed to obtain expansion.\COMMENT{We need to discuss and edit this text.} We now present this mechanism. The first ingredient is \cref{contraction_region}, which narrows down the possibilities for the problematic cases, i.e., those other than Case~(a0).
From \cref{contraction_region,Prop:Udisjoint,caseA,caseB,caseC} we get a strictly invariant cone family on $\hat M$ as follows.
\begin{corollary}\label{corollary:InvQuadrant}
If $0<\phistar<\tan^{-1}{\big(1/3\big)}$ and $R$ satisfies \eqref{eqJZHypCond}, then the quadrant cone $\mathcal{Q}_x(\mathrm{I,III})\dfn\big\{\frac{d\theta}{d\phi}\in[0,+\infty]\big\}$ in the tangent space at $x\in\hat{M}$ is strictly invariant under the differential $D\hat{F}$ of the return map from \cref{def:SectionSetsDefs}.\COMMENT{Two small things: 1) The proof is a little long for calling this a corollary. Maybe proposition?\newline2) The wording needs a little fix because a cone at a given point can only be invariant if the point is fixed. We should refine the statement to say what is actually being proved.}
\end{corollary}
\begin{proof}
As noted before \cref{def:MhatReturnOrbitSegment}, nonsingular $x\in\hat{M}$ have a return orbit segment containing $x_0\in\Mrout$, $x_1=(\Phi_1,\theta_1)\in\MRin$, $x_2\in\Mrin$. By \cref{def:MhatReturnOrbitSegment} we have:
\[
x_0\in\Mrout\xrightarrow{\mathcal{F}}x_1\in\MRin\xrightarrow{\mathcal{F}}\cdots \xrightarrow{\mathcal{F}}\mathcal{F}^{n_1}(x_1)\in\MRout\xrightarrow{\mathcal{F}}x_2\in\Mrin.
\]
In cases (a1), (b) and (c) of \eqref{eqMainCases}, \cref{contraction_region} gives $x_0\in\Nout$ and $x_2\in\Nin$, then $x\notin\Mrin\cap\Mrout$ because otherwise $x=x_0\in \Nout\cap \textcolor{red}{M^{\text{\upshape{out}}}_{r,0}}=\emptyset$ by \cref{Prop:Udisjoint}. Thus, by \cref{def:SectionSetsDefs}, we have $x\in \mathcal{F}^{-1}(\Mrout\cap\Mrin)$ and $x_0=\mathcal{F}(x)\in\Mrout\setminus\Mrin$ .

On the other hand, $x_2\in\Nin$ and \cref{Prop:Udisjoint} imply $x_2\notin\big(\textcolor{blue}{M_{r,0}^{\text{\upshape in}}}\cup\textcolor{blue}{M_{r,1}^{\text{\upshape in}}}\cup\textcolor{blue}{M_{r,2}^{\text{\upshape in}}}\big)$, so $m(x_2)\ge3$, hence $x_2\notin\Mrout$ and $x_3\dfn\mathcal{F}(x_2)\notin\Mrout$, and $\mathcal{F}^2(x_2)=\mathcal{F}(x_3)\notin\Mrout$. Hence, the \cref{def:MhatReturnOrbitSegment} return orbit segment in cases (a1), (b) and (c) is 
\begin{equation}\label{eq:a1bccasesorbitsegment}
\underbracket{x\in\hat{M}\xrightarrow{\mathcal{F}}x_0\in\Mrout}_{\text{collisions on }\Gamma_r}\xrightarrow{\mathcal{F}}\underbracket{x_1\in\MRin\xrightarrow{\mathcal{F}}\cdots \xrightarrow{\mathcal{F}}\mathcal{F}^{n_1}(x_1)\in\MRout}_{\text{collisions on }\Gamma_R}\xrightarrow{\mathcal{F}}\underbracket{x_2\in\Mrin\xrightarrow{\mathcal{F}}x_3\cdots\xrightarrow{\substack{m(x_2)-2\\\text{ steps}}}\mathcal{F}^{m(x_2)-2}(x_3)=\hat{F}(x)}_{\text{collisions on }\Gamma_r},
\end{equation} with 
$D\mathcal{F}_{x_3}^{m(x_2)-2}=\begin{pmatrix}
1 & 2(m(x_2)-2) \\
0 & 1 
\end{pmatrix}$ in $\phi\theta$-coordinates. Hence, $m(x_2)\ge3$ indicates that
\begin{equation}\label{eq:mx2ge3halfquadrantinvariant}
D\mathcal{F}^{m(x_2)-2}_{x_3}\big(\Int\mathcal{Q}_{x_3}(\mathrm{I,III})\big)\subset\Int\big\{(d\phi,d\theta)\in \mathcal{Q}_{\hat{F}(x)}(\mathrm{I,III})\bigm|\frac{d\theta}{d\phi}\in[0,1]\big\}.
\end{equation}\newline In case (a1) of \eqref{eqMainCases}, \cref{caseA} says that $D\mathcal{F}^2_{x_0}$ is a negative matrix with $n_1=0$. Since $x$, $x_0$ are the collisions on the same arc $\Gamma_r$, $D\Fc_x=\begin{pmatrix}
1 & 2 \\
0 & 1 
\end{pmatrix}$. Since $x_2$, $x_3$ are the collisions on the same arc $\Gamma_r$, $D\Fc_{x_2}=\begin{pmatrix}
1 & 2 \\
0 & 1 
\end{pmatrix}$. Therefore, $D\mathcal{F}^{n_1+4}_x=D\mathcal{F}^4_x\overbracket{=}^{\mathrlap{\text{chain rule}}}D\mathcal{F}_{x_2}\cdot D\mathcal{F}^2_{x_0}\cdot D\mathcal{F}_x$ is a negative matrix. So $D\Fc^{n_1+4}_x\big(\mathcal{Q}_x(\mathrm{I,III})\big)=D\mathcal{F}^4_x\big(\mathcal{Q}_x(\mathrm{I,III})\big)\subset \Int\mathcal{Q}_{x_3}(\mathrm{I,III})$.
\newline In case (b) of \eqref{eqMainCases}, \cref{caseB} says that $D\mathcal{F}^4_x$ is a negative matrix with $n_1=0$, so $D\Fc^{n_1+4}_x\big(\mathcal{Q}_x(\mathrm{I,III})\big)=D\mathcal{F}^4_x\big(\mathcal{Q}_x(\mathrm{I,III})\big)\subset \Int\mathcal{Q}_{x_3}(\mathrm{I,III})$.\newline In case (c) of \eqref{eqMainCases}, \cref{caseC} says that $D\mathcal{F}^{4+n_1}_x$ is a positive matrix with $n_1\ge1$, so $D\mathcal{F}^{4+n_1}_x\big(\mathcal{Q}_x(\mathrm{I,III})\big)\subset \Int\mathcal{Q}_{x_3}(\mathrm{I,III})$.

Therefore, in cases (a1), (b) and (c) of \eqref{eqMainCases}, by the chain rule, the quadrant is strictly invariant. In fact, since $m(x_2)\ge3$, we have the image of the quadrant under $D\hat{F}_x$ is 
\begin{equation}\label{eq:strictlyinvcone_casea1bc}
    \begin{aligned}
        D\hat{F}_x(\mathcal{Q}_x(\mathrm{I,III}))&= D\mathcal{F}^{n_1+4+m(x_2)-2}_x(\mathcal{Q}_x(\mathrm{I,III}))\overbracket{\subset}^{\mathclap{\substack{D\mathcal{F}^{4+n_1}_x\big(\mathcal{Q}_x(\mathrm{I,III})\big)\subset \Int\mathcal{Q}_{x_3}(\mathrm{I,III})\text{ and chain rule}}}} D\mathcal{F}_{x_3}^{m(x_2)-2}(\Int\mathcal{Q}_{x_3}(\mathrm{I,III}))\\&\overbracket{\subset}^{\mathclap{\text{\eqref{eq:mx2ge3halfquadrantinvariant}}}}\Big\{\text{interior of }\big\{(d\phi,d\theta)\in \mathcal{Q}_{\hat{F}(x)}(\mathrm{I,III})\bigm|\frac{d\theta}{d\phi}\in[0,1]\big\}\Big\}\nfd\mathcal{HQ}_{\hat{F}(x)}(\mathrm{I,III})
    \end{aligned}
\end{equation}

Now we check case (a0) of \eqref{eqMainCases}. There are 4 subcases.

Subcase 1) $x\in\mathcal{F}^{-1}(\Mrout\setminus\Mrin)$ and $x_2\in\Mrout\cap\Mrin$, that is, $m(x_2)=0$.

In this subcase, $x_2=\hat{F}(x)$. Since $x_0=\mathcal{F}(x)$, and $p(x)$, $p(\mathcal{F}(x))$ are on the same boundary $\Gamma_r$, we have $D\mathcal{F}_{x}=\begin{pmatrix}
1 & 2 \\
0 & 1 
\end{pmatrix}$ in $\phi\theta$-coordinates. Thus, the quadrant is invariant under $D\mathcal{F}$. i.e. $D\mathcal{F}_x(\mathcal{Q}_x(\mathrm{I,III}))\subset\mathcal{Q}_{x_0}(\mathrm{I,III})$. By \cref{caseA}, $D\mathcal{F}^2_{x_0}$ is a negative matrix, $D\mathcal{F}_{x_0}^2(\mathcal{Q}_{x_0}(\mathrm{I,III}))\subset\Int\mathcal{Q}_{x_2}(\mathrm{I,III})=\Int\mathcal{Q}_{\hat{F}(x)}(\mathrm{I,III})$. Hence, by the chain rule, $D\mathcal{F}^3_x(\mathcal{Q}_x(\mathrm{I,III}))\subset D\mathcal{F}_{x_0}^2(\mathcal{Q}_{x_0}(\mathrm{I,III}))\subset\Int\mathcal{Q}_{x_2}(\mathrm{I,III})=\Int\mathcal{Q}_{\hat{F}(x)}(\mathrm{I,III})$.

Subcase 2) $x\in\mathcal{F}^{-1}(\Mrout\setminus\Mrin)$ and $x_2\notin\Mrout\cap\Mrin$, that is, $m(x_2)\ge1$. 

From $x_2$ to $\hat{F}(x)$ are the $m(x_2)$ collisions on the same arc, the differential $D\mathcal{F}_{x_2}^{m(x_2)-1}=\begin{pmatrix}
1 & 2(m(x_2)-2) \\
0 & 1 
\end{pmatrix}$ in the coordinates
$(\phi,\theta)$. Then for the same reasons as in subcase 1) and by \cref{caseA}. The quadrant is strictly invariant.

Subcase 3) $x\in\Mrout\cap\Mrin$ and $x_2\in\Mrout\cap\Mrin$, that is, $m(x_2)=0$. The \cref{caseA} immediately indicates that the quadrant is strictly invariant.

Subcase 4) $x\in\Mrout\cap\Mrin$ and $x_2\notin\Mrout\cap\Mrin$, that is, $m(x_2)\ge1$. The same reasons as in subcase 2) and by \cref{caseA} give that the quadrant is strictly invariant.
\end{proof}
\begin{remark}
    We will see that even the half-quadrant family $\big\{\frac{d\theta}{d\phi}\in[0,1]\big\}$ is strictly invariant (\cref{corollary:halfquadrantcone}).
\end{remark}

\section{Monotonicity Properties}\label{sec:FPCLA}\hfill

The proofs of the results in  \Cref{sec:ECE} (\Cref{RMASUH,TEcase_a0,TEcase_a1,TEcase_b,TEcase_c,contraction_region} and \Cref{def:ApettrajectoryExpansion,RMASUHwithProof} ) require delicate estimates across various parameter ranges.\COMMENT{Later sections also seem to use monotonicity. Search through the code to find references to facts from this section and here make forward references.} In order to make those feasible, we now establish monotonicity of various of the quantities involved, so the inequalities will later only have to be established at end-points of parameter ranges.\COMMENT{Is there a way to refer forward to the exact places where the monotonicity results are being used? Wherever they are used, there should be explicit references to the corresponding statement here and in EXACTLY the right places.} Readers will probably be inclined to first skip this section and refer back to it as needed later.

\begin{figure}[ht]
\begin{center}
    \begin{tikzpicture}[xscale=0.7,yscale=0.7]
        \pgfmathsetmacro{\Dlt}{0.5}
        \pgfmathsetmacro{\PHISTAR}{3};
        \pgfmathsetmacro{\WIDTH}{12};
        \pgfmathsetmacro{\LENGTH}{5};
        \tkzDefPoint(-\PHISTAR,0){LL};
        \tkzDefPoint(-\PHISTAR,\WIDTH){LU};
        \tkzDefPoint(-\PHISTAR,0.5*\WIDTH){LM};
        \tkzDefPoint(\PHISTAR,0){RL};
        \tkzDefPoint(\PHISTAR,\WIDTH){RU};
        \tkzDefPoint(\PHISTAR,0.5*\WIDTH){RM};
        \tkzDefPoint(-\PHISTAR,\PHISTAR){X};
        \tkzDefPoint(-\PHISTAR,\WIDTH-\PHISTAR){IX};
        \tkzDefPoint(\PHISTAR,\WIDTH-\PHISTAR){Y};
        \tkzDefPoint(\PHISTAR,\PHISTAR){IY};
        \tkzDefPoint(-\PHISTAR,\PHISTAR+\Dlt){LLD};
        \tkzDefPoint(-\PHISTAR,\WIDTH-\PHISTAR-\Dlt){LUD};
        \tkzDefPoint(\PHISTAR,\PHISTAR+\Dlt){RLD};
        \tkzDefPoint(\PHISTAR,\WIDTH-\PHISTAR-\Dlt){RUD};
        \tkzDefPoint(0,0.5*\WIDTH){CP};
        \tkzDefPoint(0,\PHISTAR+\Dlt){DLM};
        \tkzDefPoint(0,\WIDTH-\PHISTAR-\Dlt){DUM};
        \draw [thick] (DLM) --(DUM);
        \draw [thin] (LL) --(LU);
        \draw [thin] (LU) --(RU);
        \draw [thin] (RL) --(RU);
        \draw [thin] (LL) --(RL);
        \draw [blue,thin] (LL) --(IX);
        \draw [blue,thin,dashed] (LL) --(IY);
        \draw [blue,thick,dashed] (IX) --(RU);
        \draw [blue,thick] (IY) --(RU);
        \draw [red,thin,dashed]  (X)--(RL);
        \draw [red,thin,dashed]  (Y)--(LU);
        \draw [purple,thick] (LUD)--(RUD);
        \draw [purple,thick] (LLD)--(RLD);
        \draw [purple,thick] (LM)--(RM);
        \coordinate (x1W) at ($(X)!0.7!(RL)$);
        \coordinate (x1E) at ($(IY)!0.7!(RL)$);
        \draw [green,thick] (x1W)--(x1E);
        \tkzLabelPoint[right](x1E){\small $ (\Phi,\theta=\constant_3<{\pi}/{2})\in\MRout$}
        \coordinate (x2W) at ($(X)!0.3!(LL)$);
        \coordinate (x2E) at ($(IY)!0.3!(LL)$);
        \draw [green,thick] (x2W)--(x2E);
        \tkzLabelPoint[left](x2W){\small $ (\Phi,\theta=\constant_4)\in\MRin$}
        \node[] at (0.5*\PHISTAR, 0.5*\PHISTAR+0.5*\Dlt+0.25*\WIDTH){\Large $R_1$};
        \node[] at (0.5*\PHISTAR, \WIDTH-0.5*\PHISTAR-0.5*\Dlt-0.25*\WIDTH){\Large $R_4$};
        \node[] at (-0.5*\PHISTAR,         0.5*\PHISTAR+0.5*\Dlt+0.25*\WIDTH){\Large $R_2$};
        \node[] at (-0.5*\PHISTAR,         \WIDTH-0.5*\PHISTAR-0.5*\Dlt-0.25*\WIDTH){\Large $R_3$};
        \fill [orange, opacity=18/30](DLM) -- (RLD) -- (RM) -- (CP) -- cycle;
        \tkzDefPoint(0.7*\PHISTAR,\PHISTAR+\Dlt){X1L};
        \tkzDefPoint(0.7*\PHISTAR,0.47*\WIDTH){X1U};
        \tkzDefPoint(\PHISTAR,0.47*\WIDTH){X1UR};
        \draw [green,thick] (X1L)--(X1U);
        \draw [red,ultra thin,dashed] (X1U)--(X1UR);
        \tkzLabelPoint[right](X1UR){$\theta=\frac{\pi}{2}-(\sin^{-1}{(\frac{R\sin{\Phi}}{|O_rP|})}-\Phi)=\frac{\pi}{2}-\measuredangle{O_rPO_R}>\frac{\pi}{2}-\phistar$};
        \node[] at (0.5*\PHISTAR, 0.52*\WIDTH){ \small$ (\Phi=\constant_1,\theta)$};
        \tkzLabelPoint[below](LL){$-\Phistar$};
        \tkzLabelPoint[left](LL){$0$};
        \tkzLabelPoint[below](RL){$+\Phistar$};
        \tkzLabelPoint[left](LU){$\pi$};
        \tkzLabelPoint[left](LM){${\pi}/{2}$};
        \tkzLabelPoint[right](RLD){$\theta=\sin^{-1}{({2r}/{R})}$};
        \tkzLabelPoint[right](RUD){$\theta=\pi-\sin^{-1}{({2r}/{R})}$};
        \tkzDefPoint(0,0.35*\WIDTH){HL};
        \tkzDefPoint(\PHISTAR,0.35*\WIDTH){HR};
        \draw[green,thick](HL)--(HR);
        \tkzLabelPoint[right](HR){\small $(\Phi\ge0,\theta=\constant_2<{\pi}/{2})\in R_1\subset\MRout\cap\MRin$}
    \end{tikzpicture}
\end{center}
    \caption{$\mathfrak D^\circ$ is the union of $R_1\dfn\big\{(\Phi,\theta)\bigm|0\le\Phi<\Phistar,\,\sin^{-1}{({2r}/{R})}\le\theta\le{\pi}/{2}\big\}\subset\MRin\cap\Mrout$,\newline
    $R_2\dfn\big\{(\Phi,\theta)\bigm|-\Phistar<\Phi<0,\,\sin^{-1}{({2r}/{R})}\le\theta\le{\pi}/{2}\big\}$,\newline
    $R_3\dfn\big\{(\Phi,\theta)\bigm|-\Phistar<\Phi\le0,\,\pi-\sin^{-1}{({2r}/{R})}\ge\theta\ge{\pi}/{2}\big\}$,\newline
    $R_4\dfn\big\{(\Phi,\theta)\bigm|\Phistar>\Phi>0,\,\pi-\sin^{-1}{({2r}/{R})}\ge\theta\ge{\pi}/{2}\big\}$.}\label{fig:R1}
\end{figure}

For $x_0\in \Mrout$, $\mathcal{F}(x_0)=x_1=(\Phi,\theta)\in \MRin$
 (see \cref{fig:MR}). If $\sin^{-1}{(2r/R)}\le\theta\le\frac{\pi}{2}$, then $x_1\in R_1\cup R_2$.
 The following propositions establish monotonicity when we shift a line with a fixed angle crossing the flatter boundary or rotate a line around a fixed point on the flatter boundary. i.e., We move $x_1\in\MRin$ horizontally or vertically in \cref{fig:R1}. Note that by symmetry, the following conclusions on $R_1$ and $R_2$ of the phase space shown in \Cref{fig:R1} will apply to regions $R_3$ and $R_4$. The conclusion of the case on $\theta \in[0,\pi/2]$ will also imply the same conclusion of the case on $\theta \in[\frac{\pi}{2},\pi]$. The subsequent discussions are limited to $\theta\le\pi/2$.
\begin{figure}[ht]
\begin{center}
\begin{tikzpicture}
    \tkzDefPoint(0,0){Or};
    \pgfmathsetmacro{\rradius}{2.3};
    \pgfmathsetmacro{\bdist}{6};
    \tkzDefPoint(0,\bdist){OR};
    \pgfmathsetmacro{\phistardeg}{40};
    \pgfmathsetmacro{\XValueArc}{\rradius*sin(\phistardeg)};
    \pgfmathsetmacro{\YValueArc}{\rradius*cos(\phistardeg)};
    \tkzDefPoint(\XValueArc,-1.0*\YValueArc){B};
    \tkzDefPoint(-1.0*\XValueArc,-1.0*\YValueArc){A};
    \pgfmathsetmacro{\Rradius}{veclen(\XValueArc,\YValueArc+\bdist)};
\clip (-\Rradius,-3.2) rectangle (\Rradius,6.5);
    \draw[name path=Cr,thin] (Or) circle (\rradius);
    \draw[name path=CR,ultra thin,blue!50] (OR) circle (\Rradius);
    \pgfmathsetmacro{\PXValue}{\Rradius*sin(5)};
    \pgfmathsetmacro{\PYValue}{\bdist-\Rradius*cos(5)};
    \tkzDefPoint(\PXValue,\PYValue){P};
    \path[name path=lineQ] (P) -- +(30:7);
    \draw [blue,ultra thin](P) -- +(30:7);
    \path[name path=lineT] (P) -- +(-20:3);
    \draw [blue,ultra thin](P) -- +(-20:3);
    \path[name path=lineT0] (P) -- +(160:7);
    \draw [blue,ultra thin](P) -- +(160:7);
    \path[name intersections={of=Cr and lineQ,by={Q}}];
    \path[name intersections={of=Cr and lineT,by={T}}];
    \path[name intersections={of=Cr and lineT0,by={T0}}];
    \path[name intersections={of=CR and lineT0,by={T1}}];
    
    \tkzDefLine[orthogonal=through P](P,OR)\tkzGetPoint{LF};
    \tkzDefLine[perpendicular=through P](OR,P)\tkzGetPoint{RF};
    \draw[-latex,draw=red,ultra thin,dashed] (OR)--(P);
    \draw[-latex,draw=red,ultra thin,dashed] (P)--(LF);
    \draw[-latex,draw=red,ultra thin,dashed] (P)--(RF);
    
    \tkzLabelPoint[above left,xshift=-3,yshift=3](P){\large $P$};
    \tkzLabelPoint[right,xshift=1](Q){$Q$};
    \tkzLabelPoint[below,xshift=2,yshift=-2](T){\large $T$};
    \tkzLabelPoint[right,above,xshift=8.5](T0){\large $T_0$};
    \tkzLabelPoint[left,xshift=-3.5](Or){$O_r$};
    \tkzLabelPoint[above](OR){$O_R$};
    \tkzDefPointBy[projection=onto P--Q](OR)\tkzGetPoint{M1};
    \tkzDefPointBy[projection=onto P--T0](OR)\tkzGetPoint{M2};
    \tkzDefPointBy[projection=onto P--Q](Or)\tkzGetPoint{M3};
    \tkzDefPointBy[projection=onto P--T0](Or)\tkzGetPoint{M4};
    \tkzDefPointBy[projection=onto OR--M1](Or)\tkzGetPoint{M5};
     \tkzDefPointBy[projection=onto OR--M2](Or)\tkzGetPoint{M6};
    \draw[red,dashed,ultra thin] (OR)--(M1);
    \draw[red,dashed,ultra thin] (OR)--(M2);
    \draw[red,dashed,ultra thin] (Or)--(M3);
    \draw[red,dashed,ultra thin] (Or)--(M4);
    \draw[red,dashed,ultra thin] (Or)--(M5);
    \draw[red,dashed,ultra thin] (Or)--(M6);
    \tkzDrawPoints(OR,Or,P,Q,T,T0,M1,M2,M3,M4,M5,M6);
    \draw[purple](Or)--(OR);
    \tkzLabelPoint[below,right,xshift=4](M1){\large E};
    \tkzLabelPoint[below,left,xshift=-3,yshift=-3](M2){\large F};
    \tkzLabelPoint[above,right,xshift=0.5,yshift=12](M3){\large G};
    \tkzLabelPoint[below](M4){\large H};
    \tkzLabelPoint[above,right,xshift=0.9,yshift=1.8](M5){\large I};
    \tkzLabelPoint[left,xshift=-8.5](M6){\large J};
    \tkzLabelPoint[right,below,xshift=50](LF){Tangent line at P};
    \tkzLabelSegment[left,xshift=-4](OR,Or){$b$}
    \tkzMarkRightAngle(OR,M1,P);
    \tkzMarkRightAngle(OR,M2,T1)
    \tkzMarkRightAngle(Or,M3,Q);
    \tkzMarkRightAngle(Or,M4,T1);
    \tkzMarkRightAngle(Or,M5,M1);
    \tkzMarkRightAngle(Or,M6,M2);
    \tkzMarkAngle[arc=lll,ultra thin,size=1](P,OR,M1);
    \tkzMarkAngle[arc=lll,blue,ultra thin,size=2.5](Or,OR,P);
    \tkzMarkAngle[arc=lll,ultra thin,size=1.5](M6,OR,P);
    
    \node at ($(OR) + (-75:1.3)$)  {\Large $\theta$};
    \node at ($(OR) + (-88:3.2)$)  {
    \Large \textcolor{blue}{$\Phi$}};
    \node at ($(OR) + (-100:2.1)$)  {\Large $\theta$};   
    \coordinate (Y) at (0,-2.8);
    \draw[red,ultra thin,dashed](Or)--(Y);
    \tkzLabelPoint[below](Y){$Y$}
    \draw[very thin,dashed](Or)--(P);
\end{tikzpicture}
    \caption{
    $T_0$, $P$, $Q$ are the positions $p(x_0)$, $p(x_1)$, $p(x_2)$ in billiard table.\\
    $E$ is the foot of $PQ$'s perpendicular through $O_R$. $F$ is the foot of $PT_0$'s perpendicular through $O_R$.\\
    $G$ is the foot of $PQ$'s perpendicular through $O_r$. $H$ is the foot of $PT_0$'s perpendicular through $O_r$.\\
    $I $ is the foot of $O_RE$'s perpendicular through $O_r$.$J$ is the foot of $O_RF$'s perpendicular through $O_r$.}\label{fig:CaseAMonotone}
\end{center}
\end{figure}

\begin{proposition}[$\theta\mapsto L_1$ increases]\label{monotone_p1}
With $\constant_1\in[0,\Phistar)$ fixed, for $x_1$ on the vertical green line\COMMENT{\textcolor{red}{Note that \Cref{fig:CaseAMonotone} has two different $\theta$.}\newline\wentao{I should mean they are equal to the $\theta$ parameter of $x_1$ } Which may be better put somewhere else. Still Thinking. Lower priority.\newline\boris They look quite different, and there is no $x_1$ in the figure.\wentao{In caption, it's mentioned $P=p(x_1)$. $x_1$ is in the phase space not in the billiard table.}} \[\Big\{x_1=(\Phi=\constant_1,\theta)\in R_1=\big\{(\Phi,\theta)\in \MRin \bigm| 0\le\Phi\le\Phistar,\ \sin^{-1}{({2r}/{R})}\le\theta\le{\pi}/{2}\big\}\Big\}\]
in  \cref{fig:R1}, $P=p(x_1)$ is independent of $\theta$ and $L_1=\tau_1$ is a strictly increasing function of $\theta\in[\sin^{-1}{({2r}/{R})},{\pi}/{2}]$ (notations as in \cref{def:fpclf}).
\end{proposition}
\begin{proof}
We show that $\frac{d L_1}{d\theta}>0$. Use cartesian coordinates with $O_r=(0,0)$ and $O_R=(0,b)$. Let $Q=p(\mathcal{F}(x_1))$ then $\tau_1=L_1=|PQ|$  (\cref{def:fpclf}). Since $x_1=(\Phi,\theta)\in\MRin$,  $P=p(x_1)$ has cartesian coordinates  
\begin{equation} \label{xP}
P=(x_P,y_P)=(R\sin{\Phi},b-R\cos(\Phi)),
\end{equation}
and likewise $Q$ has coordinates
\begin{equation}\label{xQ}
Q=(x_Q,y_Q)=(x_P+L_1\cos{(\Phi+\theta)},y_P+L_1\sin{(\Phi+\theta)}.
\end{equation}
Thus,
\begin{equation}\label{eqLsquare}
r^2=x^2_Q+y^2_Q=L^2_1+[2x_P\cos{(\Phi+\theta)}+2y_P\sin{(\Phi+\theta)}]L_1+x^2_P+y^2_P.
\end{equation}
Since $\Phi$ is held constant, and hence so are $x_P$, and $y_P$, we have $\frac{dx_P}{d\theta}=0$ and $\frac{dy_P}{d\theta}=0$. Differentiating this last equation with respect to $\theta$ and using the chain rule therefore yields $$2L_1\frac{d L_1}{d\theta}+L_1[-2x_P\sin{(\Phi+\theta)}+2y_P\sin{(\Phi+\theta)}]+\frac{d L_1}{d\theta}[2x_P\cos{(\Phi+\theta)}+2y_P\sin{(\Phi+\theta)}]=0,$$
or
\begin{equation} \label{dL1Dtheta}
\frac{d L_1}{d\theta}=\frac{x_P\sin{(\Phi+\theta)}-y_P\cos{(\Phi+\theta)}}{L_1+x_P\cos{(\Phi+\theta)}+y_P\sin{(\Phi+\theta)}}L_1\stackrel{\eqref{xP}}{=}\frac{\overbracket{R\cos{(\theta)}-b\cos{(\Phi+\theta)}}^{\mathclap{=\text{distance}(O_r,\text{line }PQ)\text{ (proof below)}}}}{\underbracket{L_1+b\sin{(\Phi+\theta)}-R\sin{\theta}}_{=d_2>0\text{ (proof below)}}}L_1>0,
\end{equation}
so  if $x_1=(\Phi=\text{constant},\theta)$, then $\tau_1=L_1$ is an increasing function of $\theta\in[\sin^{-1}{(2r/R)},\pi/2]$.

It remains to compute the numerator and denominator in \eqref{dL1Dtheta}.

To prove that the denominator satisfies $L_1+b\sin{(\Phi+\theta)}-R\sin{\theta}=d_2$, let $x_2=(\Phi_2,\theta_2)$, and observe from \cref{fig:CaseAMonotone}  that $p(x_2)=Q$, 
$|PQ|=L_1=\tau_1$, $|PE|=R\sin\theta=d_1$, $|GQ|=d_2=r\sin{\theta_2}$, $|GE|=b\sin{(\theta+\Phi)}$.

There are six possible configurations of the points $P,G,Q,E$ on the line $PQ$, depending on the lengths $\tau_1,d_1,d_2$.

Case 1) $\tau_1\le d_2\le d_1$: the ordering is $\overrightarrow{GPQE}$, and $|PQ|+|GE|-|PE|=|GQ|$.

Case 2) $\tau_1\le d_1\le d_2$: the ordering is $\overrightarrow{GPQE}$, and $|PQ|+|GE|-|PE|=|GQ|$. This case does not occur when $d_1\ge2r$.

Case 3) $d_1\le\tau_1\le d_2$: the ordering is $\overrightarrow{GPEQ}$, and $|PQ|+|GE|-|PE|=|GQ|$. This case does not occur when $d_1\ge2r$.

Case 4) $d_2\le\tau_1\le d_1$: the ordering is $\overrightarrow{PGQE}$, and $|PQ|+|GE|-|PE|=|GQ|$.

Case 5) $d_2\le d_1<\tau_1$: the ordering is $\overrightarrow{PGEQ}$ or $\overrightarrow{PEGQ}$, but since $\tau_1\le2d_1,\tau_1\le2d_2,\tau_1\le d_1+d_2$, the ordering   $\overrightarrow{PEGQ}$ cannot happen. So, for the ordering $\overrightarrow{PGEQ}$ observe that $|PQ|+|GE|-|PE|=|GE|+|EQ|=|GQ|$, even though this case also does not occur when $d_1\ge2r$.

Case 6) $d_1\le d_2<\tau_1$: reasoning as in case 5) we see that the ordering is $\overrightarrow{PGEQ}$, hence $|PQ|+|GE|-|PE|=|GE|+|EQ|=|GQ|$, even though this case does not occur when $d_1\ge2r$.

Hence in all six cases, $L_1+b\sin{(\theta+\Phi)}-R\sin{\theta}=|PQ|+|GE|-|PE|=|GQ|=d_2=r\sin{\theta_2}$.



To prove that the numerator $R\cos{(\theta)}-b\cos{(\Phi+\theta)}=\text{distance}(O_r,\text{line }PQ)$. In \cref{fig:CaseAMonotone} suppose that a line $l_1$ passing through $O_r$ is parallel to the line $PQ$. For $\Phi\in[0,\Phistar)$ fixed, we observe that: if $\pi/2-\theta<\Phi$, then $|O_RI|=b\cos{(\theta+\Phi)}$, $|O_RE|=R\cos\theta$, distance$(O_r,\text{line }PQ)=|IE|=|O_RE|-|O_RI|=R\cos\theta-b\cos{(\theta+\Phi)}$. And if $\pi/2-\theta\ge\Phi$, then $|O_RI|=-b\cos{(\theta+\Phi)}$, $|O_RE|=R\cos\theta$, distance$(O_r,\text{line }PQ)=|IE|=|O_RE|+|O_RI|=R\cos\theta-b\cos{(\theta+\Phi)}$.\COMMENT{Is it clear what this means?\wentao{Much rewritten}}Therefore, \[R\cos\theta-b\cos{(\theta+\Phi)}=\text{distance}(O_r,\text{line }PQ)>0.\qedhere\]
\end{proof}
\begin{proposition}[$\theta\mapsto L_0$ decreases]\label{monotone_p2}
With $\constant_1\in[0,\Phistar)$ fixed,For $x_1=(\Phi_1=\constant_1,\theta)$ on the vertical green line in the region $R_1=\big\{(\Phi,\theta)\in \MRin \bigm| 0\le\Phi<\Phistar,\;\sin^{-1}{(2r/R)}\le\theta\le\pi/2\big\}$ in \Cref{fig:R1} (notation as in \Cref{def:fpclf}), $P=p(x_1)$ is held a constant, thus independent of $\theta$, and 
$L_0=2d_0-\tau_0$ is a strictly decreasing function of $\theta\in\big[\sin^{-1}{(2r/R)},\ \pi/2-(\sin^{-1}{(R\sin{\Phi_1}/\rho)}-\Phi_1))\big]$.\COMMENT{Eliminate \textbackslash frac in most inline cases elsewhere.\wentao{Work In Progress}}
\end{proposition}
\begin{remark}\label{rmk:monotone_p2}\hfill

\begin{enumerate}
    \item\label{item01rmkmonotone_p2} In this context, $P=p(x_1)$ is independent of $\theta$, hence so is $|O_rP|=\rho(\Phi_1)=:\rho$. And $\sin^{-1}((R/\rho)\sin{\Phi_1})-\Phi_1<\phistar$ for all constant $\Phi_1=\constant_1\in[0,\Phistar)$.\COMMENT{\textcolor{red}{\boris Does this (proposition and remark) say what you mean? Then we can keep it and adapt the proof and picture notations---as well as the wording of like propositions, starting with \Cref{monotone_p1}!\Hrule} Isn't $p(\cdot)$ always independent of $\theta$?\wentao{Yes, because $\Phi_1$ is held as a constant $\constant_1$.}}
    \item\label{item02rmkmonotone_p2} For $\Phi_1\in[0,\Phistar)$, in \cref{fig:CaseAMonotone} we see $\Phi=\Phi_1=\measuredangle{O_rO_RP}$, $\measuredangle{PO_rY}=\sin^{-1}{(\frac{R\sin\Phi_1}{\rho})}$. Therefore, $\measuredangle O_rPO_R=\measuredangle{PO_rY}-\measuredangle{O_rO_RP}=\sin^{-1}{(\frac{R\sin\Phi_1}{\rho})}-\Phi_1$. Since $0\le\measuredangle{O_rPO_R}\le\measuredangle{YO_rP}<\phistar$, $0\le\sin^{-1}{(\frac{R\sin\Phi_1}{\rho})}-\Phi_1<\phistar$,
    \[
        \pi/2-(\sin^{-1}{(R\sin{\Phi_1}/\rho)}-\Phi_1))>\pi/2-\phistar.
    \] This fact is referenced in \cref{prop:casea0smallthetaR1}.
    \item\label{item03rmkmonotone_p2} In \cref{fig:CaseAMonotone} suppose that a line $l$ passing through $O_r$ is parallel to the line $PT_0$. Therefore, distance($l,O_R$)$=|O_RJ|=b\cos(\theta-\Phi)$, distance($\text{line } PT_0,O_R$)$=|O_RF|=R\cos\theta$. For $\Phi\in[0,\Phistar)$ fixed, we observe that in \cref{fig:CaseAMonotone}, the order of $O_R,F,J$ depends on the comparing of $\measuredangle{T_0PO_R}$ with $\measuredangle{O_rO_RP}$. If $\measuredangle{T_0PO_R}=\measuredangle{O_rPO_R}$, then $PT_0$ passes through $O_r$, so $J$ and $F$ are the same point, $|O_RJ|=|O_RF|$.
    
    if $\pi/2-\theta=\measuredangle{T_0PO_R}>\measuredangle{O_rPO_R}=\measuredangle{YO_rP}-\measuredangle{O_rO_RP}=\sin^{-1}{(\frac{R\sin\Phi}{\rho})}-\Phi$, then $|O_RJ|<|O_RF|$.
    
    And if $\pi/2-\theta=\measuredangle{T_0PO_R}\le\measuredangle{O_rPO_R}=\measuredangle{YO_rP}-\measuredangle{O_rO_RP}=\sin^{-1}{(\frac{R\sin\Phi}{\rho})}-\Phi$, then $|O_RJ|\ge |O_RF|$.

\end{enumerate}
\end{remark}\NEWPAGE
\begin{proof}
With fixed $\Phi=\constant_1$, in cartesian coordinates (\cref{def:StandardCoordinateTable}) for  \cref{fig:CaseAMonotone}, 
\begin{equation} \label{xPP}
        P=(x_P,y_P)=(R\sin{\Phi},b-R\cos(\Phi)),
\end{equation}
is independent of $\theta$. 
\COMMENT{Where is ``$\sin^{-1}{(\frac{R\sin{\Phi}}{\rho})}=\measuredangle{PO_rY}$ and $\sin^{-1}{(\frac{R\sin{\Phi}}{\rho})}-\Phi=\measuredangle{O_rPO_R}\in[0,\phistar)$'' used? \textcolor{red}{Why is this stated right here?? It looks like it belongs elsewhere. Since there is no label attached, I suspect that if it is ever used later, the reader has no idea where it came from. I see no purpose.} The equations seem backwards. Is this meant to prove what is in the remark?:\textcolor{red}{F: It will be used in proposition 2.5 proof.[\textcolor{blue}{But that's already completed!!}] There are two subsection more transversal case and fewer transversal case. This inequality [\textcolor{blue}{I am not asking about an inequality! I am asking about "$\sin^{-1}{(\frac{R\sin{\Phi}}{\rho})}=\measuredangle{PO_rY}$ and $\sin^{-1}{(\frac{R\sin{\Phi}}{\rho})}-\Phi=\measuredangle{O_rPO_R}\in[0,\phistar)$."}] will make sure that the monotone vertical line is long enough to cover the fewer transversal case region.}\newline\boris There have to be explanations IN THE TEXT of what is used (and proved!) exactly where. Specific forward references are best. To \cref{prop:casea0smallthetaR1}?\wentao{I moved this claim to the \cref{rmk:monotone_p2}.}}

By \Cref{def:fpclf}, we have $L_0=2d_0-
\tau_0=|TP|$ in \cref{fig:CaseAMonotone}.
And by elementary geometry, we can observe from \cref{fig:CaseAMonotone} that if $x_1=(\Phi,\theta)\in \MRin$, then $\theta>\Phi$ and\COMMENT{\textcolor{red}{Add {\bfseries correct} punctuation in ALL displays!}}
\begin{equation}\label{xT}
        r^2=x^2_T+y^2_T\overbracket=^{\mathrlap{T=(x_T,y_T)=(x_P+L_0\cos{(\theta-\Phi)},y_P-L_0\sin{(\theta-\Phi)}}\quad}L^2_0+[2x_P\cos{(\theta-\Phi)}-2y_P\sin{(\theta-\Phi)}]L_0+x^2_P+y^2_P.
\end{equation}
Since $\Phi$, $x_P$, $y_P$ are constants, $\frac{dx_P}{d\theta}=0$, $\frac{dy_P}{d\theta}=0$. Differentiating \eqref{xT} with respect to $\theta$ and using the chain rule gives
\[
2L_0\frac{dL_0}{d\theta}-[2x_{P}\sin{(\theta-\Phi)}+2y_{P}\cos(\theta-\Phi)]L_0+[2x_{P}\cos{(\theta-\Phi)-2y_{P}\sin{(\theta-\Phi)}}]\frac{dL_0}{d\theta}=0,
\]
so
\begin{equation} \label{dL2Dtheta}
\begin{aligned}
\frac{d L_0}{d\theta}  &=\frac{x_{P}\sin{(\theta-\Phi)}+y_{P}\cos{(\theta-\Phi)}}{L_0+x_{P}\cos{(\theta-\Phi)}-y_{P}\sin{(\theta-\Phi)}}L_0\\
&\overbracket{=}^{\mathllap{\eqref{xPP}}}\frac{R\sin{(\Phi)}\sin{(\theta-\Phi)}-R\cos{\Phi}\cos{(\theta-\Phi)}+b\cos{(\theta-\Phi)}}{L_0+R\sin{\Phi}\cos{(\theta-\Phi)}+R\cos{\Phi}\sin{(\theta-\Phi)}}L_0=\frac{\overbracket{b\cos(\theta-\Phi)-R\cos{\theta}}^{\mathclap{<0\text{ if } 0<\theta<\frac{\pi}{2}-\sin^{-1}{(\frac{R\sin\Phi}{\rho})}\text{\cref{rmk:monotone_p2}\eqref{item03rmkmonotone_p2}}}}}{\underbracket{L_0-b\sin{(\theta-\Phi)+R\sin{\theta}}}_{=d_0>0\text{ (proof below)}}}L_0.
 \end{aligned}
\end{equation}

Now it remains to compute the numerator and denominator in the rightmost term of \eqref{dL2Dtheta}.

For the numerator, \cref{rmk:monotone_p2}\eqref{item03rmkmonotone_p2}: if $0<\theta<\frac{\pi}{2}-\sin^{-1}{(\frac{R\sin{\Phi}}{\rho})}+\Phi$,\COMMENT{Why are the other cases needed?}, then $b\cos{(\Phi-\theta)}-R\cos\theta=|O_RJ|-|O_RF|<0$.

For the denominator, 
\cref{fig:CaseAMonotone} gives $L_0=|TP|$,  $b\sin{(\theta-\Phi)}=|HF|$, $|PF|=R\sin\theta=d_1$, and $|TH|=d_0$.
The ordering of $T,P,H,F$  on the line $TP$ is determined by the lengths $L_0=2\tau_0-d_0,d_1,d_0$. Since\COMMENT{"When" here says that we add an assumption. I don't think you mean to do that. In that case, replace "When" by "Since"!\newline If this is an additional assumption, then there is a problem.\wentao{Done}} $d_1\ge2r$, there are two cases according to the sign of $d_0-\tau_0$:
\newline$d_1\ge2r> d_0\ge2d_0-\tau_0$: the ordering is $\overrightarrow{TPHF}$, and $L_0-b\sin{(\theta-\Phi)}+R\sin{\theta}=|TP|-|HF|+|PF|=|TH|=d_0$.
\newline$d_1\ge2r>2d_0-\tau_0>d_0$: the ordering is $\overrightarrow{THPF}$, and $|TP|-|HF|+|PF|=|TH|=d_0$ (this case will not occur when $x_1\in R_1$).

Hence, $\frac{dL_0}{d\theta}=\frac{|O_RJ|-|O_RF|}{d_0}L_0<0$, so $L_0=2d_0-\tau_0$ is a decreasing function of $\theta$.\COMMENT{This is a hopefully helpful example of how to insert just a few words to make text comprehensible. Please apply this elsewhere.\wentao{Done}}
\end{proof}

\begin{proposition}\label{monotone_p3}
 In cartesian coordinates (\cref{def:StandardCoordinateTable}):\COMMENT{This proposion needs editing. The commentary distracts and confuses. The enumerate environment is questionable given the size of each item. I have to understand this more clearly to decide what exactly to do.}
 \begin{enumerate}
     \item\label{itemmonotone_p3-1} In \cref{fig:R1}, suppose $x_1=(\Phi,\theta)$ is on a horizontal line segment in $\MRout$ with $\theta<\pi/2$ constant. Let $P=p(x_1)=(x_P,y_P)$, $Q=p(\mathcal{F}(x_1))=(x_Q,y_Q)$. Then $\tau_1$ from \cref{def:fpclf} is a function of $\Phi$ with
\begin{equation}\label{eq1monotone_p3}
    \frac{d\tau_1}{d\Phi}=-\frac{bx_Q}{d_2}.
\end{equation}
Two such horizontal line segments are shown in green in \cref{fig:R1}: $\big\{(\Phi,\theta=\constant_4<\frac{\pi}{2})\in \MRout,\Phi\in [-\Phistar,\Phistar]\big\}$ and $\big\{(\Phi,\theta=\constant_3<\Phistar)\in \MRout\big\}$.
\newline Note that if $\Fc(x_1)=(\phi_2,\theta_2)\in\Mrin$, then by \cref{def:fpclf}, $d_2=r\sin{\theta_2}$, and $x_Q=r\sin{\phi_2}$.
\newline If $x_Q>0$, then $\Phi\mapsto\tau_1$ decreases with respect to $\Phi$.

\item\label{itemmonotone_p3-2} Suppose $x_1$ is on a horizontal line segment $\big\{(\Phi,\theta=\constant_2<\frac{\pi}{2})\subseteq R_1,\ \Phi\in [0,\Phistar]\big\}\subset\MRin$ as shown in \cref{fig:R1},\COMMENT{This wording is the best among these items; the first and the last item should be edited into like form.\Hrule However, maybe there could be a global reference to \cref{fig:R1} in the first line of the proposition so we don't have to keep repeating that reference.} where $\sin^{-1}{(\frac{2r}{R})}\le \constant_2\le\frac{\pi}{2}$. Let $T_0=p(\mathcal{F}^{-1}(x_1))$ and $T$ be the intersection of $\Gamma_R$ with the line through $T_0$ and $P$. $T$ has cartesian coordinates $(x_T,y_T)$. Then $L_0=2d_0-\tau_0$ is a function of $\Phi$ with 
\begin{equation}\label{eq2monotone_p3}
    \frac{dL_0}{d\Phi}=-\frac{bx_T}{d_0}.
\end{equation} If $x_T>0$, then $\Phi\mapsto L_0$ decreases with respect to $\Phi$.
\item\label{itemmonotone_p3-3} In \cref{fig:R1}, suppose $x_1=(\Phi,\theta)$ is on a horizontal line segment in $\MRin$ with $\theta<\pi/2$ constant. Let $P=p(x_1)=(x_P,y_P)$, $T_0=p(\mathcal{F}^{-1}(x_1))=(x_{T_0},y_{T_0})$.  Then $\tau_0$ from \cref{def:fpclf} is a function of $\Phi$ with
\begin{equation}\label{eq3monotone_p3}
    \frac{d\tau_0}{d\Phi}=-\frac{bx_{T_0}}{d_0}.
\end{equation}
Two examples of such horizontal line segments are shown in green in \cref{fig:R1}: $\big\{(\Phi,\theta=\constant_4<\frac{\pi}{2})\in \MRin,\Phi\in [-\Phistar,\Phistar]\big\}$ and $\big\{(\Phi,\theta=\constant_4<\Phistar)\in \MRin\big\}$.
\newline Note that if $\Fc^{-1}(x_1)=(\phi_0,\theta_0)\in\Mrout$, then by \cref{def:fpclf}, $d_0=r\sin{\theta_0}$, and $x_{T_0}=r\sin{\phi_0}$.\newline If $x_{T_0}<0$, then $\Phi\mapsto\tau_0$ increases with respect to $\Phi$.
\end{enumerate}
\end{proposition}
\begin{proof}
\textbf{Proof of \eqref{itemmonotone_p3-1}}: In \cref{fig:CaseAMonotone}, $L_1=\tau_1=|PQ|,L_0=2d_0-\tau_0=|PT|$. The coordinates of $P,Q,T$ are the same as in \cref{monotone_p1,monotone_p2}. The equations in\eqref{xP}\eqref{eqLsquare} are still true.

Differentiating \eqref{xP} with respect to $\Phi$ gives

\begin{equation}
\begin{split}
    & \frac{dx_P}{d\Phi}=R\cos{\Phi}\\
    & \frac{dy_P}{d\Phi}=R\sin{\Phi}
\end{split}
\end{equation}

Differentiating \eqref{eqLsquare} with respect to $\Phi$ and using the chain rule yields
\begin{equation} \label{eq:dphiL1}
\begin{multlined}[.85\textwidth]
 2L_1\frac{dL_1}{d\Phi}+2\frac{dL_1}{d\Phi}\big[x_P\cos{(\theta+\Phi)}+y_P\sin{(\theta+\Phi)}\big]\\+2L_1\big[R\cos{\theta}-x_P\sin{(\theta+\Phi)+y_P\cos{(\theta+\Phi)}}\big]
 +2x_PR\cos{\Phi}+2y_PR\sin{\Phi}=0,
\end{multlined}
\end{equation}
that is,
\begin{equation}\label{eq:dphiL1_1}
\begin{aligned}
\big(L_1+b\sin{(\theta+\Phi)}-R\sin{\theta}\big)\frac{dL_1}{d\Phi} &=
{-L_1(R\cos{(\theta)}-\overbracket{x_P}^{\mathclap{=R\sin{\Phi}}}\sin{(\theta+\Phi)}+\overbracket{y_P}^{\mathclap{=b-R\cos{\Phi}}}\cos{(\theta+\Phi)})-bR\sin{\Phi}}
\\
& =
{-L(R\cos{\theta}-R\cos{\theta}+b\cos{(\theta+\Phi)})-bR\sin{\Phi}}
\\
& =
{-bL_1\cos{(\theta+\Phi)}-bR\sin{\Phi}}
.\end{aligned}
\end{equation}
Note that if $x_1\in\MRout$, then $\theta+\Phi>0$. In \cref{fig:CaseAMonotone}, $|PQ|=L_1=\tau_1$, $|PE|=d_1$, $|GQ|=d_2$, and $|GE|=b\sin{(\theta+\Phi)}$.

Using the same proof as in \cref{monotone_p1}, note that in all six cases of that proposition, we have$$L_1+b\sin{(\theta+\Phi)}-R\sin{\theta}=d_2.$$
Then by \eqref{xP} and from \cref{fig:CaseAMonotone}, $$L_1\cos{(\theta+\Phi)}+R\sin{\Phi}=(x_Q-x_P)+x_P=x_Q.$$
Therefore, we have proved item \eqref{itemmonotone_p3-1}, i.e.,
$$\frac{d\tau_1}{d\Phi}=\frac{dL_1}{d\Phi}=-b\frac{x_Q}{d_2}.$$
\newline \textbf{Proof of \eqref{itemmonotone_p3-2}}:
Differentiating \eqref{xT} with respect to $\Phi$ yields $\frac{dL_0}{d\Phi}=-b\frac{R\sin\Phi+L_0\cos{(\theta-\Phi)}}{L_0+R\sin\theta-b\sin{(\theta-\Phi)}}$.

From \cref{fig:CaseAMonotone} and equation \eqref{xP}, we get $R\sin\Phi+L_0\cos{(\theta-\Phi)}=x_P+(x_T-x_P)=x_T$, and from \eqref{dL2Dtheta} $d_0=L_0+R\sin\theta-b\sin{(\theta-\Phi)}$,
therefore $\frac{d(2d_0-\tau_0)}{d\Phi}=\frac{dL_0}{d\Phi}=-b\frac{x_T}{d_0}.$ This proves \eqref{itemmonotone_p3-2}.\COMMENT{Immediately after proof of (1) explain how it implies (3). That is currently missing. Refer to remark on "reversal symmetry" as needed\wentao{So far (3) is not a must in anywhere else. I will make it a later decision it's needed or not. The symmetry seems a worthwhile mentioning.}}
\newline \textbf{Proof of \eqref{itemmonotone_p3-3}}: The fold\COMMENT{"fold" seems to make no sense.} symmetry map $J:(\Phi,\theta=\Psi_4)\mapsto (-\Phi,\theta=\Psi_4)$ in \eqref{eq:FoldSymmetryDef} diffeomorphically maps the line segment in $\mathfrak{l}\subset\MRin$ with $\theta$ being a constant $\Psi_4$ to a line segment $J(\mathfrak{l})\subset\MRout$ with $\theta$ being the same constant $\Psi_4$. Since concerning $\tau_0$ as a function of $x_1\in\MRin$ and $\tau_1$ as a function on $\MRout$, $\tau_0$ and $\tau_1$ satisfy $\tau_0\restriction_{\MRin}=\tau_1\restriction_{\MRout}\circ J$.\COMMENT{How about replacing this whole sentence by $\tau_0\restriction_{\MRin}=\tau_1\restriction_{\MRout}\circ J$?\wentao{agree, done}} When $\tau_0$ is restricted to the line segment $\mathfrak{l}$ with $\Phi$ as parameter, $\frac{dJ}{d\Phi}=-1$ and $\frac{d\tau_0}{d\Phi}$ satisfies the following.
\begin{align*}
    \frac{d\tau_0}{d\Phi}\Bigm|_{x_1=(\Phi,\theta=\Psi_4)}=\frac{d(\tau_1\circ J)}{d\Phi}\overbracket{=}^{\text{chain rule}}\frac{d\tau_1}{d\Phi}\Bigm|_{J(x_1)}\overbracket{\frac{d J}{d\Phi}}^{=-1}=-\frac{d\tau_1}{d\Phi}\Bigm|_{J(x_1)}.
\end{align*}
Let $\Fc J(x_1)\nfd x_{2,J}=(\phi_{2,J},\theta_{2,J})\in\Mrin$, then, when restricted to the line segment $J(\mathfrak{l})$, \eqref{itemmonotone_p3-1} gives $$\frac{d\tau_1}{d\Phi}\Bigm|_{J(x_1)}=-b\frac{r\sin{\phi_{2,J}}}{r\sin{\theta_{2,J}}}.$$ Hence $\frac{d\tau_0}{d\Phi}\Bigm|_{x_1=(\Phi,\theta=\Psi_4)}=b\frac{r\sin{\phi_{2,J}}}{r\sin{\theta_{2,J}}}$. Since $\Fc^{-1}(x_1)=x_0=(\phi_0,\theta_0)$ and $J\Fc^{-1}=\Fc J$, it gives $$x_{2,J}=\Fc J(x_1)=J\Fc^{-1}(x_1)=J(x_0),$$, by \eqref{eq:FoldSymmetryDef} which means that $\phi_0=2\pi-\phi_{2,J}$ and $\theta_0=\theta_{2,J}$. Thus, \[\frac{d\tau_0}{d\Phi}\Bigm|_{x_1=(\Phi,\theta=\Psi_4)}=b\frac{r\sin{\phi_{2,J}}}{r\sin{\theta_{2,J}}}=-b\frac{r\sin{\phi_0}}{r\sin{\theta_0}}=-b\frac{x_{T_0}}{d_0}\qedhere\]
\end{proof}

\begin{proposition}[$\Phi\mapsto d_0$ increases]\label{monotone_p4}
Given $0<\phistar<\tan^{-1}(1/3)$ and $R$ which meets the hyperbolicity condition \eqref{eqJZHypCond}, $x\in\hat{M}$ has an orbit segment defined in \cref{def:fpclf} with $x_1\in \MRin$. If $x_1$ moves along a horizontal line segment: $\{\,x_1=(\Phi,\,\theta=\constant_4\,)\,\}\subset\MRin$ in \cref{fig:R1} for some $\constant_4<\sin^{-1}{(\sqrt{4r/R})}$,\COMMENT{It would be good to explain the origin of this upper bound when it first appears.} then $d_0$ is an increasing function of $\Phi$.
\end{proposition}
\begin{proof}
In the \cref{fig:R1_R2_d_0} and in the \cref{def:StandardCoordinateTable} coordinate system, 
$P=p(x_1)=(x_P,y_P)=(R\sin(\Phi),b-R\cos{\Phi})$ for $x_1=(\Phi,\theta)\in\MRin$. And $x_0=\mathcal{F}^{-1}(x_1)=(\phi_0,\theta_0)\in\Mrout, Q=p(x_0)=(x_Q,y_Q)=(r\sin\phi_0,-r\cos\phi_0)$. $T$ is the another intercept of line $PQ$ with circle $C_r$. Thus, the collision angle \(\theta\) is the angle between $\overrightarrow{TP}$ and the tangential direction of $\Gamma_R$ at $P$ given the counterclockwise orientation on $\Gamma_R$. 

By \cref{contraction_region} and its proof in \cref{subsec:ProofContractionRegion} with \eqref{eq:47}, $x_0$ and $x_*$ satisfy $\|x_0-x_*\|<\frac{16.46}{\sin{(\phistar/2)}}\sqrt{r/R}$. Therefore,
\[
|\phi_0-(2\pi-\phistar)|<\|x_0-x_*\|<\frac{16.46}{\sin{(\phistar/2)}}\sqrt{r/R}\overset{\eqref{eqJZHypCond}}<\frac{16.46}{\sin{(\phistar/2)}}\sqrt{\sin^2{\phistar}/30000}\overbracket{=}^{\mathrlap{\sin^2\phistar=4\sin^2(\phistar/2)\cos^2(\phistar/2)}\quad}2\cdot16.46\cdot30000^{-1/2}\cos(\phistar/2)<0.2.
\]
Hence, $\phi_0-2\pi+\phistar>-0.2$, so $ \phi_0>2\pi-\phistar-0.2>2\pi-0.2-\tan^{-1}(1/3)>5.7>3\pi/2$. This means $Q$ is on the lower left half of the circle $C_r$ in \cref{fig:R1_R2_d_0}. Hence $x_Q<0$, $y_Q<0$.\COMMENT{It seems like here we really need 30000. This should be highlighted.\wentao{3000 also works}\boris Because of an estimate for the cos?}


When $\phistar<\frac{\pi}{6}$ and $PT$ moves with constant collision angle $\theta<\sin^{-1}{\sqrt{4r/R}}$ with $\Gamma_R$.\COMMENT{Nonsense} In \Cref{fig:R1_R2_d_0}, as $\Phi$ increases, $P$ moves from $A$ and stops the rightmost where $P$ hits $B$ or $T$ hits $A$.\COMMENT{Fix this sentence.} Equivalently, the range of $\Phi$ is determined by the two constraints $P\in\Gamma_R$ and $T\notin\Gamma_R$.\COMMENT{Better to just omit the previous broken sentence?}
\begin{figure}[ht]
\begin{center}
\begin{tikzpicture}
    \tkzDefPoint(0,0){Or};
    \pgfmathsetmacro{\rradius}{1.7};
    \pgfmathsetmacro{\bdist}{3.5};
    \tkzDefPoint(0,\bdist){OR};
    \pgfmathsetmacro{\phistardeg}{40};
    \pgfmathsetmacro{\XValueArc}{\rradius*sin(\phistardeg)};
    \pgfmathsetmacro{\YValueArc}{\rradius*cos(\phistardeg)};
    \tkzDefPoint(\XValueArc,-1.0*\YValueArc){B};
    \tkzDefPoint(-1.0*\XValueArc,-1.0*\YValueArc){A};
    \pgfmathsetmacro{\Rradius}{veclen(\XValueArc,\YValueArc+\bdist)};
\clip(-\Rradius,-1.8) rectangle (\Rradius,4);
    \draw[name path=Cr,thin] (Or) circle (\rradius);
    \draw[name path=CR,ultra thin,blue!50] (OR) circle (\Rradius);
    \pgfmathsetmacro{\PXValue}{\Rradius*sin(5)};
    \pgfmathsetmacro{\PYValue}{\bdist-\Rradius*cos(5)};
    \pgfmathsetmacro{\LPXValue}{-\Rradius*sin(70)};
    \pgfmathsetmacro{\LPYValue}{\bdist-\Rradius*cos(70)};
    \tkzDefPoint(\PXValue,\PYValue){P};
    \tkzDefPoint(0,-2){Y};
    \draw[red,ultra thin,dashed](OR)--(Y);
    \path[name path=lineQ] (P) -- +(152:10);
    \path[name path=lineY] (Or) -- (Y);
    \path[name intersections={of=CR and lineQ,by={P0}}];
    \path[name intersections={of=Cr and lineQ,by={Q}}];
    \path[name intersections={of=lineQ and lineY,by={S}}];
    \draw[blue,thick](P)--(P0);
    \tkzDefPointBy[projection=onto P--Q](OR)\tkzGetPoint{G};
    \tkzDefPointBy[projection=onto P--Q](Or)\tkzGetPoint{H};
    \draw[red,ultra thin,dashed](Or)--(H);
    \draw[red,ultra thin,dashed](OR)--(G);
    \draw[red,ultra thin,dashed](OR)--(P);
    
    \tkzMarkAngle[arc=lll,ultra thin,size=0.8](G,OR,P);
    \tkzMarkAngle[blue,arc=lll,ultra thin,size=2](Y,OR,P);
    \node at ($(OR) + (-105:1.2)$)  {\Large $\theta$};
    \node at ($(OR) + (-88:2.6)$)  {
    \Large \textcolor{blue}{$\Phi$}};
    \tkzLabelPoint[above,yshift=3](OR){$O_R$};
    \tkzLabelPoint[left,xshift=-3](Or){$O_r$};
    \tkzLabelPoint[above,right,yshift=10](P){$P$};
    \tkzLabelPoint[left](S){$S$};
    \tkzLabelPoint[below left](G){$G$};
    \tkzLabelPoint[below left](H){$H$};
    \tkzLabelPoint[above right](P0){$T$};
    \tkzLabelPoint[above right](Q){$Q$};
    \tkzLabelPoint[below left](A){$A$};
    \tkzLabelPoint[below right](B){$B$};
    \tkzDrawPoints(OR,Or,P,Q,S,G,H,P0);
    \tkzMarkRightAngle(Or,H,P0);
    \tkzMarkRightAngle(OR,G,P0);
\end{tikzpicture}
\end{center}
    \caption{$PT$ moves along $\Gamma_R$ with constant collision angle $\theta$ small. $G$ is the perpendicular foot of $O_R$ on $PT$. $H$ is the perpendicular foot of $O_r$ on $PT$.}\label{fig:R1_R2_d_0}
\end{figure}

Let $S=(0,y)$ be the intersection of line $PT$ and the $y$-axis (in \cref{fig:R1_R2_d_0}). From \eqref{xP} and by checking the slopes of $\overrightarrow{SP}$, we see that $\frac{y-b+R\cos{\Phi}}{-R\sin{\Phi}}=-\tan{(\theta-\Phi)}$. This gives $y=b-R\frac{\cos{\theta}}{\cos{(\theta-\Phi)}}$.  
Now,
$x_1=(\Phi,\,\theta=\constant_4\,)\in\MRin$ in \cref{fig:R1}, so
\[
0<\theta-\Phi<\theta+\Phistar\overbracket{=}^{\strut\mathclap{\text{\cref{def:BasicNotations0}}}}\theta+\sin^{-1}(r\sin\phistar/R)\overbracket{<}^{\mathclap{\theta=\constant_4<\sin^{-1}\sqrt{4r/R}}}\sin^{-1}{\sqrt{4r/R}}+\sin^{-1}(r\sin\phistar/R)\overbracket{<}^{\mathclap{\strut\text{\eqref{eqJZHypCond}: }R>1700r}}\pi/2.
\]
Also, $R$, $b$, $\theta$ are constants, so $y$ increases with $\Phi$, that is, if $\Phi\big\uparrow$, then $y\big\uparrow$. Then we will show that $y$ is always negative.

The slopes of lines $PQ$ and $PS$ are the same: $\frac{y-y_Q}{0-x_Q}=\frac{y_P-y_Q}{x_P-x_Q}$, which gives $y=y_Q-x_Q\frac{y_P-y_Q}{x_P-x_Q}$. The positions of $P,Q$ in the table imply $x_P-x_Q>0$, $y_P-y_Q<0$ (so the fraction is negative). Since $Q$ is in the lower left quarter circle,  $x_Q<0>y_Q$, hence 
$y=y_Q-x_Q\frac{y_P-y_Q}{x_P-x_Q}<y_Q<0$, that is, $y$ is always negative as $x_1$ moves along the horizontal line: $\{x_1=(\Phi,\,\theta=\constant_4\,)\,\}\subset\MRin$ in \cref{fig:R1} with some $\constant_4<\sin^{-1}{(\sqrt{4r/R})}$. Hence, if $\Phi\Big\uparrow$, then $|y|\Big\downarrow$.\COMMENT{I am not a fan of these vertical arrows. See the edits in the last paragraph.}

Note that we also have $\frac{|O_rH|}{|O_RG|}=\frac{|O_rS|}{|O_RS|}=\frac{|y|}{|y|+b}=1-\frac{b}{b+|y|}\overbracket{>}^{\mathclap{|y|>0}}0$, and $|O_RG|=R\cos\theta$ remains constant since $\theta$ is held as constant. If $|y|\big\downarrow$, then $|O_rH|\big\downarrow$.

As $P$ moves on $\Gamma_R$ with $P=(R\sin\Phi,b-R\cos\Phi)\in\Gamma_R$ and $T\in\Gamma_r$ such that $\overrightarrow{TP}$ maintains a constant $\theta=\constant_4$ as the angle with tangential direction of $\Gamma_R$ at $P$, if $\Phi$ increases, then $|y|$ decreases and $0<|O_rH|$ decreases, so $d_0=\sqrt{r^2-|O_rH|^2}$ increases. Hence $d_0$ is an increasing function of $\Phi$.
\end{proof}
\begin{proposition}[$\Phi\mapsto L_0/\tau_1$ decreases]\label{monotone_p5}
Given $\phistar\in(0,\tan^{-1}{(1/3)})$ and $R$ satisfying \eqref{eqJZHypCond}, if $x_1=(\Phi, \theta)$ moves on the horizontal line segment $\big\{(\Phi, \theta=\sin^{-1}{\sqrt{4r/R}}) \bigm| \Phi\in[0,\Phistar] \big\}$ in the region $R_1$ of \cref{fig:R1}, then the length ratio $\frac{2d_0-\tau_0}{\tau_1}$ is a decreasing function of $\Phi$.
\end{proposition}

\begin{proof}
Consider \cref{fig:lengthfraction} with cartesian coordinates such that $O_r=(0,0)$, $O_R=(0,b)$. Let $P=p(x_1)$ with $x_1=(\Phi,\theta=\sin^{-1}{\sqrt{4r/R}})\in R_1\subset\MRin$. As shown in \Cref{fig:lengthfraction}, $T_0$ is the position $p(\mathcal{F}^{-1}(x_1))$, and  $Q$ is the position $p(\mathcal{F}(x_1))$. The extension of $T_0P$ intersects with $\Gamma_r$ at $T$. $Q_1$ is the point of the horizontal line through $Q$ intersecting $\Gamma_r$. $T_1$ is the point of the horizontal line through $T$ intersecting $\Gamma_r$.
Suppose 
\begin{equation}\label{eq:15}
\begin{split}
    & Q=(x_Q,y_Q)=(r\sin{\alpha},-\cos{\alpha}),\\
    & T=(x_T,y_T)=(r\sin{\beta},-\cos{\beta}).
\end{split}
\end{equation}

It is clear that $\pi/2\ge\alpha>\phistar>\beta>0$ from \cref{fig:lengthfraction}. Because $0\le\Phi\le\Phistar=\sin^{-1}{(r\sin\phistar/R)}<\sin^{-1}{\sqrt{4r/R}}=\theta$, the slope of the line $T_0P$ is $\tan{(\Phi-\theta)}<0$, the slope of the line $Q_0P$ is $\tan{(\Phi+\theta)}>0$ which is larger than $\tan{\Big(\Phi+(\Phistar-\Phi)/2\Big)}>0$, and both slopes are bounded. If $\alpha>\pi$, then $\theta\ge\pi/2$. If $3\pi/2>\alpha>\pi/2$, then by \cref{lemma:returnorbittheta1Large} $\theta\in(\frac{\pi}{6},\frac{5\pi}{6})$. With \eqref{eqJZHypCond}: $R>1700r$, it is impossible for $\theta=\sin^{-1}{\sqrt{4r/R}}$ to satisfy $\theta\ge \pi/2$ or $\theta\in(\frac{\pi}{6},\frac{5\pi}{6})$.

Next, as shown in \Cref{fig:lengthfraction}, we pick $P_1$ on segment $QQ_1$ and $P_2$ on segment $TT_1$ such that $\measuredangle{PP_1Q_1}=\measuredangle{YO_rQ}=\alpha,\measuredangle{PP_2T}=\measuredangle{YO_rT}=\beta$. By elementary geometry,  $\measuredangle{T_1T_0T}=\measuredangle{PP_2T}=\beta$ and $\measuredangle{Q_1Q_0Q}=\measuredangle{PP_1Q}=\pi-\alpha$, so $\Delta TP_2P$ is similar to $\Delta TT_0T_1$ and $\Delta QP_1P$ is similar to $\Delta QQ_0Q_1$. Therefore, \begin{equation}\label{eq:16}
\begin{split}
& \frac{1}{|P_1Q|}=\frac{|Q_1Q|}{|Q_0Q|\times|PQ|}=\frac{x_Q}{d_2|PQ|},\\
& \frac{1}{|P_2T|}=\frac{|T_1T|}{|T_0T|\times|PT|}=\frac{x_T}{d_0|PT|}.
\end{split}
\end{equation}
Since $|PQ|=\tau_1=L_1$ and $|PT|=2d_0-\tau_0=L_0$, the quotient rule gives
\begin{equation}
\label{eq:18}
\frac{d}{d\Phi}\big(\frac{2d_0-\tau_0}{\tau_1}\big)=\frac{1}{L^2_1}\big(\overbracket{\underbracket{\frac{dL_0}{d\Phi}}_{\qquad\mathllap{\overset{\eqref{eq3monotone_p3}}=-b{x_T}/{d_0}}}L_1-\underbracket{\frac{d
L_1}{d\Phi}}_{\mathrlap{\overset{\eqref{eq1monotone_p3}}=-b{x_Q}/{d_2}}\qquad} L_0}^{\mathclap{=b\big(-\frac{x_T}{d_0}L_1+\frac{x_Q}{d_2}L_0\big)
=bL_1L_0\big(-\frac{x_T}{d_0L_0}+\frac{x_Q}{d_2L_1}\big)\qquad\qquad}}\big)
=\frac{bL_0}{L_1}\big(\overbracket{-\underbracket{\frac{x_T}{d_0|PT|}}_{\mathclap{\overset{\eqref{eq:16}}=1/|P_2T|\quad}}+\underbracket{\frac{x_Q}{d_2|PQ|}}_{\mathclap{\quad\overset{\eqref{eq:16}}=1/|P_1Q|}}}^{\mathclap{\qquad\qquad\overset{\eqref{eq:16}}=\frac{|P_2T|-|P_1Q|}{|P_1Q||P_2T|}=\frac{x_Qx_T}{d_0d_2L_0L_1}\big(|P_2T|-|P_1Q|\big)}}\big)=\frac{bx_Qx_T}{d_0d_2L^2_1}\big(|P_2T|-|P_1Q|\big).
\end{equation}
\begin{figure}[ht]
\begin{tikzpicture}
    \tkzDefPoint(0,0){Or};
    \tkzDefPoint(0,-4.5){Y};
    \pgfmathsetmacro{\rradius}{4.5};
    \pgfmathsetmacro{\bdist}{19};
    \tkzDefPoint(0,\bdist){OR};
    \pgfmathsetmacro{\phistardeg}{40};
    \pgfmathsetmacro{\XValueArc}{\rradius*sin(\phistardeg)};
    \pgfmathsetmacro{\YValueArc}{\rradius*cos(\phistardeg)};
    \tkzDefPoint(\XValueArc,-1.0*\YValueArc){B};
    \tkzDefPoint(-1.0*\XValueArc,-1.0*\YValueArc){A};
    \pgfmathsetmacro{\Rradius}{veclen(\XValueArc,\YValueArc+\bdist)};
    \clip(-\rradius,-\rradius-0.3) rectangle (\rradius,0.3);
    \draw[name path=Cr,ultra thin,dashed] (Or) circle (\rradius);
    \tkzDrawArc[name path=CR](OR,A)(B);
    \tkzDrawArc[](Or,B)(A);
    \pgfmathsetmacro{\PXValue}{\Rradius*sin(3)};
    \pgfmathsetmacro{\PYValue}{\bdist-\Rradius*cos(3)};
    \tkzDefPoint(\PXValue,\PYValue){P};
    \path[name path=lineQ] (P) -- +(18:3.5);
    \path[name path=lineQL] (P) -- +(198:4);
    \path[name path=lineT] (P) -- +(-12:4);
    \path[name path=lineTL] (P) -- +(-192:5);
    \path[name intersections={of=Cr and lineQ,by={Q}}];
    \path[name intersections={of=Cr and lineQL,by={QL}}];
    \path[name intersections={of=Cr and lineT,by={T}}];
    \path[name intersections={of=Cr and lineTL,by={TL}}];
    \path[name path=lineQLH] (Q) -- +(-8,0);
    \path[name path=lineTLH] (T) -- +(-8,0);
    \path[name intersections={of=Cr and lineQLH,by={QLH}}];
    \path[name intersections={of=Cr and lineTLH,by={TLH}}];
    \draw[red,ultra thin](Q)--(QL);
    \draw[blue,ultra thin](T)--(TL);
    \draw[name path=QQLH,red,ultra thin,dashed](Q)--(QLH);
    \draw[name path=TTLH,blue,ultra thin,dashed](T)--(TLH);
    \draw[purple,ultra thin,dashed](Or)--(Y);
    \coordinate (PQ) at ($(P)-(Q)$);
    \coordinate (QQL)at ($(Q)-(QL)$);
    \coordinate (QQLH) at ($(Q)-(QLH)$);

    \coordinate (PT) at ($(P)-(T)$);
    \coordinate (TTL)at ($(T)-(TL)$);
    \coordinate (TTLH) at ($(T)-(TLH)$);
    

    \draw[name path=lineP1,red,dashed] let \p1=($(PQ)$),
    \p2=($(QQL)$),
    \p3=($(QQLH)$),
    \n1={veclen(\x1,\y1)},
    \n2={veclen(\x2,\y2)},
    \n3={veclen(\x3,\y3)}
    in (P)--([shift={(-\n1*\n2/\n3,0)}]Q);
    
    \path[name intersections={of=QQLH and lineP1,by={P1}}];

    \draw[name path=lineP2,blue,dashed] let \p1=($(PT)$),
    \p2=($(TTL)$),
    \p3=($(TTLH)$),
    \n1={veclen(\x1,\y1)},
    \n2={veclen(\x2,\y2)},
    \n3={veclen(\x3,\y3)}
    in (P)--([shift={(-\n1*\n2/\n3,0)}]T);

     \path[name intersections={of=TTLH and lineP2,by={P2}}];

    \tkzDrawPoints(Or,P,Q,QL,QLH,A,B,P1,P2,T,TL,TLH,Y);
    \draw[red,ultra thin,dashed](Or)--(T);
    \draw[red,ultra thin,dashed](Or)--(Q);
    
    
    \tkzLabelPoint[left,xshift=-3](Or){$O_r$};
    \tkzLabelPoint[above left,yshift=0.1,xshift=-0.3](P){$P$};
    \tkzLabelPoint[right,xshift=1](Q){$Q$};
    \tkzLabelPoint[right,xshift=1](T){$T$};
    \tkzLabelPoint[above,yshift=1](TL){$T_0$};
    \tkzLabelPoint[left,xshift=-1](QL){$Q_0$};
    \tkzLabelPoint[left,xshift=-1](QLH){$Q_1$};
    \tkzLabelPoint[left,xshift=-1](TLH){$T_1$};
    \tkzLabelPoint[above right](P1){$P_1$};
    \tkzLabelPoint[below](P2){$P_2$};
    \tkzLabelPoint[left](Y){$Y$};
    \tkzLabelPoint[left](A){$A$};
    \tkzLabelPoint[right](B){$B$};
    \tkzMarkAngle[arc=lll,ultra thin,size=1.7](Y,Or,Q);
    \tkzMarkAngle[blue,arc=lll,ultra thin,size=0.7](Y,Or,T);
    \tkzMarkAngle[arc=lll,ultra thin,size=0.3](QLH,P1,P);
    \tkzMarkAngle[blue,arc=lll,ultra thin,size=0.3](T,P2,P);
    \node at ($(Or) + (-65:1.9)$)  {\Large $\alpha$};
    \node at ($(Or) + (-78:1.0)$)  {
    \Large \textcolor{blue}{$\beta$}};
    \node at ($(P1) + (-160:0.5)$)  {$\alpha$};
    \node at ($(P2) + (-20:0.45)$)  {
    \textcolor{blue}{$ \beta$}};
    \draw[name path=lineE,blue,dashed] (T)-- +(1.5,0);
    \tkzDefLine[perpendicular=through Q](T,TLH)\tkzGetPoint{QY};
    \path[name path=LineQY](Q)--(QY);
    \path[name intersections={of=LineQY and lineE,by={E}}];
    \draw[red,dashed,ultra thin](Q) -- (E);    
    \path[name path=LinePY](P)--+(90:1.5);
    \path[name path=LineP2Y](P2)--+(90:1.5);
    \path[name intersections={of=LinePY and QQLH,by={F}}];
    \path[name intersections={of=LineP2Y and QQLH,by={G}}];
    \draw[red,dashed,ultra thin](P)--(F);
    \draw[red,dashed,ultra thin](P2)--(G);
    \tkzMarkRightAngle(T,E,Q);
    \tkzMarkRightAngle(P,F,QLH);
    \tkzMarkRightAngle(P2,G,QLH);
    \tkzDrawPoints(E,F,G);
    \tkzLabelPoint[below](E){$E$};
    \tkzLabelPoint[above](F){$H$};
    \tkzLabelPoint[above](G){$G$};
    \path[name path=LineOrY](Or)--(Y);
    \path[name intersections={of=LineOrY and CR,by={M}}];
    \tkzLabelPoint[below](M){$M$};
    \tkzDrawPoints(M);
\end{tikzpicture}
\caption{$\tau_1$ and $2d_0-\tau_0$ ratio analysis with $\beta<\alpha$}\label{fig:lengthfraction}
\end{figure}
We now compute this last expression (see \eqref{eq:23}) and then show that it is negative, as required.

In \cref{fig:lengthfraction}, $QQ_1,TT_1$ are the horizontal lines, and $QE\perp TE$, $PH\perp HQ$, $P_2G\perp GQ$.
The $Q,T$ coordinates satisfy:
\begin{equation}\label{eq:hthqyqytyp}
    \begin{aligned}
    -r\cos{\alpha}&=y_Q=y_P+h_Q,\\
    -r\cos{\beta}&=y_T=y_P-h_T.
    \end{aligned}
\end{equation}
When $(\Phi,\theta)=x_1\in\Mrin$, $P=p(x_1)$\COMMENT{$(\Phi,\theta)=x_1\in\Mrin,P=p(x_1)$ must be replaced by $(\Phi,\theta)=x_1\in\Mrin$, $P=p(x_1)$, and likewise everywhere else throughout the text, and also with periods. Commas and periods need to be followed by space.\newline\wentao{Agree, Work in progress}}, then $\theta>\Phi$ from either \cref{fig:MR} or \cref{fig:R1}. Hence for $P=p(x_1)$, $\overrightarrow{PT}$ always has negative slope, and $h_Q>0$, $h_T>0$. Then \cref{fig:lengthfraction} further shows that
\begin{equation}\label{eq:anglelengthobservationforalphabeta}
\begin{aligned}
(i)&\quad |P_2G|=|QE|=h_T+h_Q,\\
(ii)&\quad |PH|=h_Q,|HP_1|=h_Q\cot{\alpha},\\
(iii)&\quad |GH|=h_T\cot{\beta},\\
(iv)&\quad \measuredangle{QTE}=\pi-\measuredangle{T_1TO_r}-\measuredangle{O_rTQ}=\pi-(\frac{\pi}{2}-\beta)-(\frac{\pi}{2}-\frac{\alpha-\beta}{2})=\frac{\alpha+\beta}{2},\\
(v)&\quad |TE|\overbracket{=}^{\mathllap{\text{ \cref{fig:lengthfraction}: }TE\perp QE}}|QE|\cot\angle{QTE}\overset{(iv)}=|QE|\cot{(\frac{\alpha+\beta}{2})}\overset{(i)}{=}(h_T+h_Q)\cot{(\frac{\alpha+\beta}{2})},\\
(vi)&\quad |GH|+|HP_1|+|P_1Q|=|GQ|=|P_2E|=|P_2T|+|TE|.
\end{aligned}
\end{equation}
From equations \eqref{eq:anglelengthobservationforalphabeta} and\eqref{eq:hthqyqytyp} and some trigonometry, we further obtain:
\begin{equation}\label{eq:21}
\begin{aligned}
|P_2T|-|P_1Q| &\overbracket{=}^{\mathclap{\text{\eqref{eq:anglelengthobservationforalphabeta} (vi)}}} |GH|+|HP_1|-|TE|\overbracket{=}^{\mathclap{\text{\eqref{eq:anglelengthobservationforalphabeta} (ii)(iii)(v)}}}h_T\cot{\beta}+h_Q\cot{\alpha}-(h_T+h_Q)\cot{(\frac{\alpha+\beta}{2})}\\
&\overbracket{=}^{\mathclap{\text{\eqref{eq:hthqyqytyp}}}} (y_P+r\cos{\beta})\cot{\beta}-(r\cos{\alpha}+y_P)\cot{\alpha}-(r\cos{\beta}-r\cos{\alpha})\cot{(\frac{\alpha+\beta}{2})}\\
&=(y_P+r\cos{\beta})\cot{\beta}-(r\cos{\alpha}+y_P)\cot{\alpha}+2r\sin{(\frac{\beta-\alpha}{2})}\sin{(\frac{\beta+\alpha}{2})}\cot{(\frac{\alpha+\beta}{2})}\\
&=(y_P+r\cos{\beta})\cot{\beta}-(r\cos{\alpha}+y_P)\cot{\alpha}+2r\sin{(\frac{\beta-\alpha}{2})}\cos{(\frac{\alpha+\beta}{2})}\\
&=(y_P+r\cos{\beta})\cot{\beta}-(r\cos{\alpha}+y_P)\cot{\alpha}-r\sin\alpha+r\sin\beta\\
&=y_P\cot{\beta}-y_P\cot{\alpha}+r\frac{\cos^2\beta}{\sin\beta}-r\frac{\cos^2\alpha}{\sin\alpha}-r\sin\alpha+r\sin\beta\\
&=y_P\cot{\beta}-y_P\cot{\alpha}+r\frac{\cos^2\beta+\sin^2\beta}{\sin\beta}-r\frac{\cos^2\alpha+\sin^2\alpha}{\sin\alpha}\\
&= y_P\cot{\beta}-y_P\cot{\alpha}-r\frac{1}{\sin{\alpha}}+r\frac{1}{\sin{\beta}}=y_P\frac{\cos\beta}{\sin\beta}-y_P\frac{\cos\alpha}{\sin\alpha}+r\frac{\sin{\alpha}-\sin{\beta}}{\sin{\alpha}\sin{\beta}}\\
&= y_P\frac{\sin{\alpha}\cos{\beta}-\cos{\alpha}\sin{\beta}}{\sin{\alpha}\sin{\beta}}+r\frac{\sin{\alpha}-\sin{\beta}}{\sin{\alpha}\sin{\beta}}= y_P\frac{\sin{(\alpha-\beta)}}{\sin{\alpha}\sin{\beta}}+r\frac{\sin{\alpha}-\sin{\beta}}{\sin{\alpha}\sin{\beta}}\\
&=y_P\frac{\sin{(\alpha-\beta)}}{\sin{\alpha}\sin{\beta}}+r\frac{2\sin{(\frac{\alpha-\beta}{2})}\cos{(\frac{\alpha+\beta}{2})}}{\sin{\alpha}\sin{\beta}}=y_P\frac{\sin{(\alpha-\beta)}}{\sin{\alpha}\sin{\beta}}+r\frac{2\sin{(\frac{\alpha-\beta}{2})}\cos{(\frac{\alpha+\beta}{2})\cos{(\frac{\alpha-\beta}{2})}}}{\sin{\alpha}\sin{\beta}\cos{(\frac{\alpha-\beta}{2})}}\\
&=y_P\frac{\sin{(\alpha-\beta)}}{\sin{\alpha}\sin{\beta}}+r\frac{\sin{(\alpha-\beta)}\cos{(\frac{\alpha+\beta}{2})}}{\sin{\alpha}\sin{\beta}\cos{(\frac{\alpha-\beta}{2})}}=\frac{\sin{(\alpha-\beta)}}{\sin{\alpha}\sin{\beta}}\big(y_P+r\frac{\cos{(\frac{\alpha+\beta}{2})}}{\cos{(\frac{\alpha-\beta}{2})}}\big).
\end{aligned}
\end{equation}
In \cref{fig:twoTangentLines} below,  the tangent lines at $T,Q$ intersect at $S$. Then by $|O_rS|=\frac{|O_rQ|}{\cos{\measuredangle{QO_rS}}}=\frac{r}{\cos{\frac{\alpha-\beta}{2}}}$ and $\measuredangle{MO_rS}=\frac{\alpha+\beta}{2}$
\begin{equation}\label{eq:xsys}
S=(x_S,y_S)=(|O_rS|\sin\angle{MO_rS},-|O_rS|\cos\angle{MO_rS})=(r\frac{\sin{(\frac{\alpha+\beta}{2})}}{\cos{(\frac{\alpha-\beta}{2})}},-r\frac{\cos{(\frac{\alpha+\beta}{2})}}{\cos{(\frac{\alpha-\beta}{2})}}),
\end{equation}
hence $|P_2T|-|P_1Q|\overset{\eqref{eq:21}}=\frac{\sin{(\alpha-\beta)}}{\sin{\alpha}\sin{\beta}}(y_P-y_S)$. Combine equations ~\eqref{eq:15}, \eqref{eq:18}, \eqref{eq:21}, and \eqref{eq:xsys} to obtain
\begin{equation}\label{eq:23}
\frac{d}{d\Phi}(\frac{2d_0-\tau_0}{\tau_1})=\frac{bx_Qx_T}{d_0d_2\tau^2_1}\frac{\sin{(\alpha-\beta)}}{\sin{\alpha}\sin{\beta}}(y_P-y_S)\overset{\eqref{eq:15}}=\frac{br^2\sin{(\alpha-\beta)}}{d_0d_2\tau^2_1}(y_P-y_S).
\end{equation}
\Cref{lemma:fractionLemma1,lemma:fractionLemma2} below and the Intermediate-Value Theorem imply the following. If $p(x_1)=P$ moves from the midpoint $M$ of $\arc{AB}$ to the corner $B$ with constant small crossing angle with the arc $AB$, then $y_P\le y_S$ with equality  only when $P$ is at the mid-point $M$ or the corner $B$. Since $\alpha>\beta$, the right-hand side of \eqref{eq:23} is negative (except for two points), so $\frac{2d_0-\tau_0}{\tau_1}$ is a decreasing function of $\Phi\in[0,\Phistar]$.
\end{proof}

\begin{figure}[h]
\begin{tikzpicture}
    \tkzDefPoint(0,0){Or};
    \tkzDefPoint(0,-4.3){Y};
    \pgfmathsetmacro{\rradius}{4.5};
    \pgfmathsetmacro{\bdist}{19};
    \tkzDefPoint(0,\bdist){OR};
    \pgfmathsetmacro{\phistardeg}{40};
    \pgfmathsetmacro{\XValueArc}{\rradius*sin(\phistardeg)};
    \pgfmathsetmacro{\YValueArc}{\rradius*cos(\phistardeg)};
    \tkzDefPoint(\XValueArc,-1.0*\YValueArc){B};
    \tkzDefPoint(-1.0*\XValueArc,-1.0*\YValueArc){A};
    \pgfmathsetmacro{\Rradius}{veclen(\XValueArc,\YValueArc+\bdist)};
     \clip(-\rradius,-\rradius-0.3) rectangle (\rradius,0.3);
    \draw[name path=Cr,ultra thin,dashed] (Or) circle (\rradius);
    \tkzDrawArc[name path=CR](OR,A)(B);
    \tkzDrawArc[](Or,B)(A);
    \pgfmathsetmacro{\PXValue}{\Rradius*sin(3)};
    \pgfmathsetmacro{\PYValue}{\bdist-\Rradius*cos(3)};
    \tkzDefPoint(\PXValue,\PYValue){P};
    \path[name path=lineQ] (P) -- +(18:3.5);
    \path[name path=lineQL] (P) -- +(198:4);
    \path[name path=lineT] (P) -- +(-12:4);
    \path[name path=lineTL] (P) -- +(-192:5);
    \path[name intersections={of=Cr and lineQ,by={Q}}];
    \path[name intersections={of=Cr and lineQL,by={QL}}];
    \path[name intersections={of=Cr and lineT,by={T}}];
    \path[name intersections={of=Cr and lineTL,by={TL}}];
    \draw[red,ultra thin,dashed] (P) --(Q);
    \draw[blue,ultra thin,dashed] (P) --(T);
    \draw[red,ultra thin,dashed] (Or) --(Q);
    \draw[blue,ultra thin,dashed] (Or) --(T);
    \tkzDefLine[orthogonal=through Q](Q,Or)\tkzGetPoint{QI};
    \tkzDefLine[orthogonal=through T](Or,T)\tkzGetPoint{TI};
    \path[name path=ST] (T)--(TI);
    \path[name path=SQ] (Q)--(QI);
    \path[name intersections={of=ST and SQ,by={S}}];
    \draw[blue,ultra thin,dashed] (S) --(T);
    \draw[red,ultra thin,dashed] (S) --(Q);
    \draw[red,ultra thin,dashed] (Or)--(Y);
    \draw[blue,ultra thin,dashed](TL)--(P);
    \tkzMarkRightAngle(Or,T,S);
    \tkzMarkRightAngle(Or,Q,S);
    \tkzLabelPoint[below](P){$P$};
    \tkzLabelPoint[below](T){$T$};
    \tkzLabelPoint[right](Q){$Q$};
    \tkzLabelPoint[right](S){$S$};
    \tkzLabelPoint[left](Or){$O_r$};
    \tkzLabelPoint[below](Y){$Y$};
    \tkzLabelPoint[right](TL){$T_0$};
    \path[name path=LineOrY](Or)--(Y);
    \path[name intersections={of=LineOrY and CR,by={M}}];
    \path[name path=LineHS](S)--+(-5,0);
    \path[name intersections={of=LineHS and LineOrY,by={M0}}];
    \draw[red,ultra thin,dashed](M0)--(S);
    \tkzLabelPoint[below](M){$M$};
    \tkzLabelPoint[below left](A){$A$};
    \tkzLabelPoint[below](B){$B$};
    \tkzLabelPoint[above left](M0){$M_0$};
    \tkzDrawPoints(Or,P,Q,T,S,M,A,B,TL,M0);
    \tkzMarkAngle[arc=lll,ultra thin,size=1.7](Y,Or,Q);
    \tkzMarkAngle[blue,arc=lll,ultra thin,size=0.7](Y,Or,T);
    \node at ($(Or) + (-65:1.9)$)  {\Large $\alpha$};
    \node at ($(Or) + (-78:1.0)$)  {
    \Large \textcolor{blue}{$\beta$}};
\end{tikzpicture}
\caption{Tangent lines at $T$ and $Q$ intersecting at $S$}\label{fig:twoTangentLines}
\end{figure}

\begin{lemma}[$y_P\ne y_S$ in \cref{fig:twoTangentLines}]\label{lemma:fractionLemma1}With the notations in \cref{def:LemonTableConfiguration}, suppose that our lemon billiard satisfies the condition \eqref{eqJZHypCond} with the coordinate system in \cref{def:StandardCoordinateTable} such that $O_r=(0,0)$ and $O_R=(0,b)$. In \cref{fig:twoTangentLines}, let $P$, $Q$, $T_0$ be the positions in the table of $x_1=(\Phi,\theta)\in \Mrin$, $x_2=\mathcal{F}(x_1)\in \Mrin$, $x_0=\mathcal{F}^{-1}(x_1)\in \Mrout$ and $M$ be the midpoint of $\Gamma_R$. Suppose $T$ is the intersection of the extension of the line $\overrightarrow{T_0P}$ with the circle $C_r$. Let $S=(x_S,y_S)$ be the intersection of the tangent line of circle $C_r$ at $T$ with the tangent line of circle $C_r$ at $Q$. If $\overrightarrow{PQ}$ moves on $\arc{MB}$ with a constant crossing angle $\theta\in[\sin^{-1}{(2r/R)},\sin^{-1}{(\sqrt{4r/R})}]$ with the tangential direction of $\Gamma_R$ at $P=(x_P,y_P)$ and $P$ is the interior point of the arc segment $\arc{MB}$, then $y_P\ne y_S$.
\end{lemma}
\begin{proof}
In \cref{fig:twoTangentLines}, suppose that the horizontal line through $S$ intersects the $y$-axis at $M_0$. From the coordinates of $Q$,$T$,$S$ in \cref{fig:twoTangentLines}, the equations~\eqref{eq:15} and \eqref{eq:xsys}, we have $M_0=(0,y_S)=(0,-r\frac{\cos{(\frac{\alpha+\beta}{2})}}{\cos{(\frac{\alpha-\beta}{2})}})$.

Note that since $\measuredangle{O_rM_0S}=\measuredangle{O_rTS}=\pi/2$, $O_r,M_0,T,S$ are all on the circle with $O_rS$ as its diameter. Hence $\measuredangle{TM_0S}=\measuredangle{TO_rS}=\frac{\alpha-\beta}{2}$, because the angles at the circumference subtended by the same arc are equal and $S$ is on the bisector of $\angle{TO_rQ}$.

Similarly, since $\measuredangle{O_rM_0S}=\measuredangle{O_rQS}=\pi/2$, $O_r,M_0,S,Q$ are all on the circle with $O_rS$ as its diameter. Hence $\measuredangle{QM_0S}=\measuredangle{QO_rS}=\frac{\alpha-\beta}{2}$. Therefore, 
\[\measuredangle{QM_0S}=\measuredangle{TM_0S}=\frac{\alpha-\beta}{2}\] in \cref{fig:twoTangentLines} and \cref{fig:M0PS}.



\COMMENT{Can this be made less complicated and more clear? Also, is the sign correct?\newline\wentao{I find There is another easier way to show. That is to show $O_r$, $S$, $T$, $M_0$ are on same circle. This only requires to show $\measuredangle{TSO_r}=\beta$ which is quite obvious(\boris In \cref{fig:twoTangentLines} that angle looks VERY different from $\beta$). Similarly $Q$, $S$, $M_0$, $O_r$ being on same circle is obvious, hence $\measuredangle{QM_0S}=\measuredangle{QO_rS}=\frac{\alpha-\beta}{2}$}\boris I will let you rewrite this eventually. It sounds like a major improvement.\wentao{You can check that the current computations using trigomnometry is correct.}\small\bfseries\boris I'd like to understand the alternative argument.\wentao{Now the trigonometry are avoided/not used. The proof uses the Thales theorem}}

To prove that $y_P\ne y_S$ in \Cref{fig:twoTangentLines} we now assume $y_S=y_P$ in \cref{fig:M0PS} for purposes of contradiction.\COMMENT{Look for further shortcuts} This means that $P$ is on the line segment $M_0S$.
\begin{figure}[ht]
\begin{center}
\begin{tikzpicture}
\tkzDefPoint(-3,0){M0};
\tkzDefPoint(-0.6,0){P};
\tkzDefPoint(1.8,0){S};
\tkzDefPoint(2.25,1.5){Q};
\tkzDefPoint(0.5,-1){T};
\tkzDefPoint(4.25,0.4){U};
\draw[blue,ultra thin,dashed](M0)--(T);
\draw[blue,ultra thin,dashed](M0)--(Q);
\draw[red,ultra thin](M0)--(S)--+(3.0,0);
\draw[red,ultra thin,dashed](P)--(U);
\draw[blue,ultra thin](P)--(Q);
\draw[blue,ultra thin](P)--(T);
\tkzDrawPoints(M0,P,S,Q,T);
\tkzLabelPoint[left](M0){$M_0$};
\tkzLabelPoint[below](P){$P$};
\tkzLabelPoint[below](S){$S$};
\tkzLabelPoint[right](Q){$Q$};
\tkzLabelPoint[right](T){$T$};
\tkzLabelPoint[right](U){$U$};
\tkzMarkAngle[arc=lll,ultra thin,size=1.1,mark=|](P,M0,Q);
\tkzMarkAngle[arc=lll,ultra thin,size=1.1,mark=|](T,M0,P);
\end{tikzpicture}
\end{center}
    \caption{Assume $P$ were on the line  $M_0S$ which is a horizontal line}\label{fig:M0PS}
\end{figure}
In \cref{fig:M0PS} denote by $PU$ the tangent line at P of the billiard boundary $\Gamma_R=\arc{AB}$. $PU$ has positive slope and bisects $\angle{QPT}$.

Since $\theta>\Phi$, $PT$ has negative slope, while $PQ$ has positive slope. And since $\theta,\Phistar$ are both small, $\angle{TPM_0}$, $\angle{QPM_0}$ are always obtuse angles when $P$ is on $\arc{MB}$.

By the sine law, we have \[\frac{\sin{\measuredangle{PTM_0}}}{|M_0P|}=\frac{\sin{\measuredangle{PM_0T}}}{|PT|}\text{ and } \frac{\sin{\measuredangle{PQM_0}}}{|M_0P|}=\frac{\sin{\measuredangle{PM_0Q}}}{|PQ|}\]

Since $\angle{PTM_0}$ and $\angle{PQM_0}$ are acute, $\measuredangle{PM_0Q}=\measuredangle{PM_0T}$, and by \cite[Equation 4.1]{hoalb}, $|PT|=2d_0-\tau_0<\tau_1=|PQ|$,  we get $\measuredangle{PTM_0}>\measuredangle{PM_0Q}$, and thus
\begin{equation}\label{eq:27}
    \measuredangle{TPS}>\measuredangle{QPS}.
\end{equation}
But PU is the bisector of $\angle{QPT}$ with positive slope, i.e., $\measuredangle{SPU}>0$, hence
\begin{equation}\label{eq:28}
\begin{aligned}
\measuredangle{QPS}&=\measuredangle{QPU}+\measuredangle{SPU}>\measuredangle{QPU}=\measuredangle{TPU}>\measuredangle{TPU}-\measuredangle{SPU}=\measuredangle{TPS},
\end{aligned}
\end{equation}
contrary to \eqref{eq:27}, so $P$, $S$ cannot be on the same horizontal line, hence $y_P\ne y_S$.
\end{proof}
\begin{lemma}[$y_S-y_P$ is increasing at $\Phi=0$ (for \cref{monotone_p5})]\label{lemma:fractionLemma2}
In the context of \cref{lemma:fractionLemma1}, that is, $S=(x_S,y_S)$ as the intersection of tangent lines of $C_r$ at $Q$ and $T$, suppose that $\overrightarrow{PQ}$ moves in $\arc{MB}$ with a constant crossing angle $\theta\in[\sin^{-1}{(2r/R)},\sin^{-1}{(\sqrt{4r/R})}]$ with the tangential direction of $\Gamma_R$ in $P$. Then $P=(x_P,y_P)=(R\sin{\Phi},b-R\sin{\Phi})$ with $\Phi\in[0,\Phistar)$, and $y_P,y_S$ from \eqref{eq:xsys} are continuous differentiable functions of $\Phi$ such that $\frac{dy_P}{d\Phi}\bigr|_{
\Phi=0}<\frac{dy_S}{d\Phi}\bigr|_{\Phi=0}$, that is, $y_S-y_P$ is increasing at $M$.
\COMMENT{The short version would seem to be that $y_S-y_P$ is increasing at $M$? \wentao{agree.}}
\end{lemma}
\begin{proof}
When $\Phi=0$, i.e., $P=M$ is at the midpoint of the arc $\Gamma_R$, then $x_P=0$, so \[\frac{dy_P}{d\Phi}\bigr|_{\Phi=0}=\frac{d(b-R\cos{\Phi})}{d\Phi}=R\sin{\Phi}=x_P=0.\] 
It remains to show that $\frac{dy_S}{d\Phi}\bigr|_{\Phi=0}>0$, and this is a brute-force computation.
First, differentiate \eqref{eq:xsys} with respect to $\Phi$, using the chain rule and trigonometry:
\begin{equation}\label{eq:30}
\begin{aligned}
\frac{dy_S}{d\Phi}&=\frac{d}{d\Phi}\Big[\frac{-r\cos{(\frac{\alpha+\beta}{2})}}{\cos{(\frac{\alpha-\beta}{2})}}\Big]\\&=\frac{-r}{\cos^{2}{(\frac{\alpha-\beta}{2})}}\Big[-\cos{(\frac{\alpha-\beta}{2})}\sin{(\frac{\alpha+\beta}{2})}\frac{d}{d\Phi}(\frac{\alpha+\beta}{2})+\cos{(\frac{\alpha+\beta}{2})}\sin{(\frac{\alpha-\beta}{2})}\frac{d}{d\Phi}(\frac{\alpha-\beta}{2})\Big]\\&=\frac{r}{2\cos^{2}{(\frac{\alpha-\beta}{2})}}\Big[\sin{(\frac{\alpha+\beta}{2})}\cos{(\frac{\alpha-\beta}{2})}-\sin{(\frac{\alpha-\beta}{2})}\cos{(\frac{\alpha+\beta}{2})}\Big]\frac{d\alpha}{d\Phi}\\
& \qquad+\frac{r}{2\cos^{2}{(\frac{\alpha-\beta}{2})}}\Big[\sin{(\frac{\alpha+\beta}{2})}\cos{(\frac{\alpha-\beta}{2})}+\sin{(\frac{\alpha-\beta}{2})}\cos{(\frac{\alpha+\beta}{2})}\Big]\frac{d\beta}{d\Phi}\\
& =\frac{r}{\cos^{2}{(\frac{\alpha-\beta}{2})}}\Big(\frac{1}{2}\sin{\beta}\frac{d\alpha}{d\Phi}+\frac{1}{2}\sin{\alpha}\frac{d\beta}{d\Phi}\Big)
\end{aligned}
\end{equation}
We next find the derivatives of \(\alpha\) and \(\beta\), with $L_1=\tau_1,L_0=2d_0-\tau_0$ from \cref{def:fpclf}.
In \cref{fig:twoTangentLines},
\begin{equation}\label{eq:29}
\begin{split}
& Q=(x_Q,y_Q) =(r\sin{\alpha},-r\cos{\alpha}),\\
& T=(x_T,y_T) =(r\sin{\beta},-r\cos{\beta}),
\end{split}
\end{equation}
where $\pi/2\ge\alpha>\beta>0$ are smooth functions of $\Phi$.
Also in \cref{fig:twoTangentLines},  $\alpha-\theta-\Phi$ is the angle between $\overrightarrow{PQ}$ and the tangent $\overrightarrow{SQ}$, so 
\begin{equation}\label{eqd2}
    d_2=r\sin{(\alpha-\theta-\Phi)}.
\end{equation} 
Observing in \cref{fig:CaseAMonotone,fig:twoTangentLines}, we have the following coordinate equations for the points $P,Q$. 

\begin{gather}\label{eq:gatheredcoordinatesEquationsPQT_1}
r\sin{\alpha}=x_Q=x_P+L_1\cos{(\Phi+\theta)}=R\sin\Phi+L_1\cos{(\Phi+\theta)}\\\label{eq:gatheredcoordinatesEquationsPQT_2}
\overbracket{b-y_Q}^{\mathllap{y\text{-coordinate difference between } O_R\text{ and }Q.}}=\overbracket{-L_1\sin{(\Phi+\theta)}}^{\mathrlap{y\text{-coordinate difference between } P\text{ and }Q.}}+\underbracket{R\cos\Phi}_{\mathclap{y\text{-coordinate difference between } O_R\text{ and }P.}}
\end{gather}
Since $\theta$ is held as a constant, differentiating \eqref{eq:gatheredcoordinatesEquationsPQT_1} with respect to $\Phi$ in leftmost and rightmost sides, using \Cref{monotone_p3}\eqref{itemmonotone_p3-1}, \eqref{eq:29}, chain rule and trigonometry will further give
\COMMENT{Instead of listing a ton of possibly useful things, label each relation with its justification.\wentao{More reference added in each step.}} 
\begin{align*}
r\cos{\alpha}\frac{d\alpha}{d\Phi}\quad\qquad\qquad&\overbracket{=}^{\mathclap{\substack{\text{In \eqref{eq:gatheredcoordinatesEquationsPQT_1} taking derivative with respect}\\\text{to }\Phi\text{ and using chain rule}}}}\qquad\quad\qquad\frac{d}{d\Phi}(R\sin{\Phi})-L_1\sin{(\Phi+\theta)}+\frac{dL_1}{d\Phi}\cos{(\Phi+\theta)}\\
&\overbracket{=}^{\mathclap{\text{\Cref{monotone_p3}\eqref{itemmonotone_p3-1}}}}R\cos{\Phi}-L_1\sin{(\Phi+\theta)}-\frac{bx_Q}{d_2}\cos{(\Phi+\theta)}\overbracket{=}^{\mathclap{\eqref{eq:gatheredcoordinatesEquationsPQT_2}}}-y_Q+b-\frac{b\sin\alpha}{\sin{(\alpha-\theta-\Phi)}}\cos{(\Phi+\theta)}\\
&\overbracket{=}^{\mathllap{\text{\eqref{eq:29}}}}r\cos{\alpha}-b\frac{\sin\alpha\cos{(\Phi+\theta)}-\sin{(\alpha-\theta-\Phi)}}{\sin{(\alpha-\theta-\Phi)}}=r\cos{\alpha}-b\frac{\sin{(\Phi+\theta)}\cos{\alpha}}{\underbracket{\sin{(\alpha-\theta-\Phi)}}_{\overset{\eqref{eqd2}}{=}d_2/r}}.
\end{align*}
So we have found that 
\begin{equation}\label{eq:32}
\frac{d\alpha}{d\Phi}=
1-\frac{b\sin{(\Phi+\theta)}}{d_2}.
\end{equation}
In \cref{fig:twoTangentLines}, $\beta+\theta-\Phi$ is the angle between  $\overrightarrow{TS}$, and $\overrightarrow{T_0T}$, so
\begin{equation}\label{eqd0}
    d_0=r\sin{(\beta+\theta-\Phi)}.
\end{equation}

Observing in \cref{fig:CaseAMonotone,fig:twoTangentLines}, we have the following coordinate equations for the points $P,Q$. 

\begin{gather}\label{eq:gatheredcoordinatesEquationsPQT_3}
r\sin{\beta}=x_T=x_P+L_0\cos{(\theta-\Phi)}=R\sin\Phi+L_0\cos{(\theta-\Phi)}\\\label{eq:gatheredcoordinatesEquationsPQT_4}
\underbracket{b-y_T}_{\mathllap{y-\text{coordinate difference between }O_R\text{ and }T}}=\overbracket{R\cos\Phi}^{\mathllap{y-\text{coordinate difference between }O_R\text{ and }P}}+\underbracket{L_0\sin{(\theta-\Phi)}}_{\mathclap{y-\text{coordinate difference between }P\text{ and }T}}
\end{gather}
Since $\theta$ is held as a constant, differentiating \eqref{eq:gatheredcoordinatesEquationsPQT_3} with respect to $\Phi$ in leftmost and rightmost sides, using \Cref{monotone_p3}\eqref{itemmonotone_p3-2}, \eqref{eq:29}, chain rule and trigonometry will further give
\COMMENT{Once again, attach each reference to the step where it is used.\wentao{More reference added in each step.}}
\[
\begin{aligned}
    r\cos\beta\frac{d\beta}{d\Phi}\quad\qquad\qquad&\overbracket{=}^{\mathclap{\substack{\text{In \eqref{eq:gatheredcoordinatesEquationsPQT_3} taking derivative with respect}\\\text{to }\Phi\text{ and using chain rule}}}}\qquad\quad\qquad\frac{d}{d\Phi}(R\sin\Phi)+L_0\sin{(\theta-\Phi)}+\cos{(\theta-\Phi)}\frac{dL_0}{d\Phi}\\
    &\overbracket{=}^{\mathllap{\text{\Cref{monotone_p3}\eqref{itemmonotone_p3-2}}}}R\cos\Phi+L_0\sin{(\theta-\Phi)}-\frac{bx_T}{d_0}\cos{(\theta-\Phi)}\overbracket{=}^{\mathclap{\eqref{eq:gatheredcoordinatesEquationsPQT_4},\eqref{eqd0}}}b-y_T-\frac{bx_T\cos{(\theta-\Phi)}}{r\sin{(\beta+\theta-\Phi)}}\\
    &\overbracket{=}^{\mathllap{\eqref{eq:29}}}r\cos\beta+b-\frac{b\sin\beta\cos{(\theta-\Phi)}}{\sin{(\beta+\theta-\Phi)}}=r\cos\beta+b\frac{\sin{(\beta+\theta-\Phi)}-\sin\beta\cos{(\theta-\Phi)}}{\sin{(\beta+\theta-\Phi)}}\\
    &=r\cos\beta+b\frac{\cos\beta\sin{(\theta-\Phi)}}{\underbracket{\sin{(\beta+\theta-\Phi)}}_{\overset{\eqref{eqd0}}=d_0/r}}
\end{aligned}
\]
We have thus found that 
\begin{equation}\label{eq:33}
\frac{d\beta}{d\Phi}
{=}1+b\frac{\sin{(\theta-\Phi)}}{d_0}.
\end{equation}
Inserting \eqref{eq:32}, \eqref{eq:33} into \eqref{eq:30}, we get
\begin{equation}\label{eq:34}
\frac{dy_S}{d\Phi}=\frac{r}{2\cos^2{(\frac{\alpha-\beta}{2})}}\big[\sin{\beta}\frac{d_2-b\sin{(\theta+\Phi)}}{d_2}+\sin{\alpha}\frac{d_0+b\sin{(\theta-\Phi)}}{d_0}\big].
\end{equation}
If $\Phi=0$, then $d_0=d_2$, and \eqref{eq:34} becomes
\[
\frac{dy_S}{d\Phi}\bigr|_{\Phi=0}=\frac{r}{2\cos^2{(\frac{\alpha-\beta}{2})}}\Big[\sin\alpha+\sin\beta+\frac{b\sin\theta}{d_0}(\sin\alpha-\sin\beta)\Big]\overset{\pi/2\ge\alpha>\beta>0}>0,
\]
so 
$\displaystyle\big(\frac{dy_S}{d\Phi}-\frac{dy_P}{d\Phi}\big)\bigr|_{\Phi=0}=\frac{dy_S}{d\Phi}\bigr|_{\Phi=0}>0$, proving \Cref{lemma:fractionLemma2}.
\end{proof}
\section{Contraction region, return map, cones, shears}\label{sec: PPDSR}
The purpose of this section is twofold. We first prove  \cref{contraction_region}.\\dots \COMMENT{Outline what this section does and where. Notably, we need to \color{red}figure out the purpose of the second subsection.}\subsection[Proof of \texorpdfstring{\cref{contraction_region}}{contraction region} (contraction region)]{Proof of \texorpdfstring{\cref{contraction_region}}{contraction region} (Region of Contraction)}\label{subsec:ProofContractionRegion}
We here prove \cref{contraction_region}: Cases {(a1), (b), and  (c) in \eqref{eqMainCases} arise only when $x_0\in\Nout$ and $x_2\in\Nin$.
\begin{proof}
We use the notation from \cref{contraction_region,fig:Mr}: $x_*=(\phistar,\pi-\phistar)$, $y_*=(2\pi-\phistar,\phistar)$ and the inversions $I x_*=(\phistar,\phistar),Iy_*=(2\pi-\phistar,\pi-\phistar)$.\COMMENT{From here to \eqref{eq:caseCNinNout2}, this is an endless string of computations that look random. Explain often what comes next and why, and when a step has been completed. Are there steps? Is there a strategy?\small\bfseries\Hrule We should probably discuss this whole section carefully.}

In Case (a1) of  \eqref{eqMainCases},\COMMENT{It should always be made clear when we go from one case to another.} suppose $x_1=(\Phi_1,\theta_1)\in\Mrin$ with $\Phi_1\in [-\Phistar,+\Phistar], \theta_1\in[\sin^{-1}{(2r/R)},\sin^{-1}{(\sqrt{4r/R})}]$,\newline then $(\phi_2,\theta_2)=x_2=\mathcal{F}(x_1)\in\Mrin,x_0=\mathcal{F}^{-1}(x_1)\in \Mrout$. As shown in \cref{fig:CaseAMonotone}, in Cartesian coordinates with $O_r=(0,0)$, $O_R=(0,b)$ (\cref{def:StandardCoordinateTable}) we have $T_0=p(x_0)$, 
\begin{equation}\label{eq:XqXt}
\begin{split}
Q=p(x_2) & =(x_Q,y_Q)=(r\sin{\phi_2},-r\cos{\phi_2}),\\
P=p(x_1) & =(x_P,y_P)=(R\sin{\Phi_1},b-R\cos{\Phi_1}).
\end{split}
\end{equation}
In \cref{fig:QonUpperHalfCircle} it is straightforward to see that $Q$ cannot be on the arc $\arc{AC}$ since otherwise 
$\pi/2<\measuredangle{QP_1B}=\Phi_1+\theta_1<\Phistar+\sin^{-1}{\sqrt{4r/R}}=\sin^{-1}{(r\sin\phistar/R)}+\sqrt{4r/R}$, which is impossible when $R>1700r$. Also, $Q$ cannot be on the arc $\arc{CD}$ which is the upper half of the circle $C_r$ since otherwise $\theta_1\in\big(\frac{\pi}{6},\frac{5\pi}{6}\big)$ by \cref{lemma:returnorbittheta1Large}, which is also impossible since $R>1700r$ and $\sin\theta_1\le\sqrt{4r/R}$. Thus, $Q$ must be on the arc $\arc{BD}$ i.e. $\phi_2\in\big(\phistar,\pi/2\big]$.\COMMENT{Will this ne used later? Where? How is the reader supposed to remember? And for how long?}

$\overrightarrow{PQ}$ has angle $\theta_1+\Phi_1$ with respect to positive x-axis, so by slope calculation we get:\COMMENT{Fix this display. The alignment and labeling are random, and this is a haphazard jumble. I suspect only the last equation needs a label.\Hrule Equations should be aligned if there is a sensible reason for it, such as repeated steps or defensible aesthetics. Otherwise each should be centered. (Sometimes multline but definitely not here.)}
\begin{equation}\label{eq:38}
\begin{split}
& \frac{b-R\cos{\Phi_1}+r\cos{\phi_2}}{R\sin{\Phi_1}-r\sin{\phi_2}}=\tan{(\theta_1+\Phi_1)}\\
& r[\cos{\phi_2}+\sin{\phi_2}\tan{(\theta_1+\Phi_1)}]=R\sin{\Phi_1}\tan{(\theta_1+\Phi_1)}+R\cos{\Phi_1}-b\\
&
r\cos{(\phi_2-\theta_1-\Phi_1)}=R\cos{\theta_1}-b\cos{(\theta_1+\Phi_1)}.
\end{split}
\end{equation}

We observe that in the last equation in \eqref{eq:38}, as $R\rightarrow\infty,\theta_1\rightarrow 0,\Phi_1\rightarrow 0$ $LHS\longrightarrow r\cos{\phi_2},RHS\longrightarrow r\cos{\phistar}$.\COMMENT{It sounds like what follows is supposed to prove this. Where does the proof end? And where is this used? Where does the next thing start, and what is it about?}

Indeed, subtracting $R-b$ from both sides of the last equation in \eqref{eq:38} gives:
\begin{equation}\label{eq:39}
\begin{split}
r\cos{(\phi_2-\theta_1-\Phi_1)}-(R-b)& =R(\cos{\theta_1}-1)+[1-\cos{(\theta_1+\Phi_1)}]b= -2R\sin^{2}{(\frac{\theta_1}{2})}+2b\sin^{2}{(\frac{\theta_1}{2}+\frac{\Phi_1}{2})}\\
& =-2(R-b)\sin^{2}{(\frac{\theta_1}{2}+\frac{\Phi_1}{2})}+2R\sin^{2}{(\frac{\theta_1}{2}+\frac{\Phi_1}{2})}-2R\sin^{2}{\frac{\theta_1}{2}}.
\end{split}
\end{equation}

Using further trigonometry and triangle inequality,\COMMENT{Not helpful, fix} we obtain\COMMENT{Again, chaos-display. Fix carefully.}

\begin{equation}\label{eq:40}
\begin{aligned}
r\cos{(\phi_2-\theta_1-\Phi_1)} & = (R-b)[1-2\sin^{2}{(\frac{\theta_1}{2}+\frac{\Phi_1}{2})}]+2R[\sin{(\frac{\theta_1}{2}+\frac{\Phi_1}{2})}-\sin{(\frac{\theta_1}{2})}][\sin{(\frac{\theta_1}{2}+\frac{\Phi_1}{2})}+\sin{(\frac{\theta_1}{2})}]\\
& = (R-b)[1-2\sin^2{(\frac{\theta_1}{2}+\frac{\Phi_1}{2})}]+2R\times 2\sin{(\frac{\Phi_1}{4})}\cos{(\frac{\theta_1}{2}+\frac{\Phi_1}{4})}2\sin{(\frac{\theta_1}{2}+\frac{\Phi_1}{4})}\cos{(\frac{\Phi_1}{4})}\\
& = R-b-2(R-b)\sin^{2}{(\frac{\theta_1}{2}+\frac{\Phi_1}{2})}+2R\sin{(\frac{\Phi_1}{2})}\sin{(\theta_1+\frac{\Phi_1}{2})}\\
r\cos{(\phi_2-\theta_1-\Phi_1)}& - r\cos\phistar  = R-b-r\cos\phistar-2(R-b)\sin^{2}{(\frac{\theta_1}{2}+\frac{\Phi_1}{2})}+2R\sin{(\frac{\Phi_1}{2})}\sin{(\theta_1+\frac{\Phi_1}{2})}\\
\end{aligned}
\end{equation}.

We also have the Cosine Law\COMMENT{Never capitalize this; fix throughout. Remove this comment when done.} $R^2=r^2+b^2-2rb\cos{(\pi-\phistar)}$ to give\COMMENT{Sentence needs fixing}

\begin{equation}\label{eq:41}
0<R-b-r\cos{\phistar}\overset{\text{Cosine Law}}=\frac{r^2\sin^2{\phistar}}{R+b+r\cos{\phistar}}<\frac{r^2\sin^2{\phistar}}{R}.
\end{equation}

Because $
\Phi_1\in [-\Phistar,+\Phistar],\Phistar=\sin^{-1}{(\frac{r\sin\phistar}{R})},\theta_1\in[\sin^{-1}{(\frac{2r}{R})},\sin^{-1}{(\sqrt{\frac{4r}{R}})}]$, we have $-\Phi_1\le\Phistar=\sin^{-1}{(\frac{r\sin\phistar}{R})}<\sin^{-1}{(\frac{2r}{R})}\le\theta_1$ and $\Phi_1\le\Phistar=\sin^{-1}{(\frac{r\sin\phistar}{R})}<\sin^{-1}{(\frac{2r}{R})}\le\theta_1$ that is $0<\Phi_1+\theta_1<2\theta_1$.\COMMENT{What are we trying to do right now?}

It is obvious that $0<R-b<r$. Thus, these combined with \eqref{eq:40} and \eqref{eq:41} yield

\begin{equation}\label{eq:42}
\begin{aligned}
|r\cos{(\phi_2-\theta_1-\Phi_1)}-r\cos\phistar|  &\overset{\eqref{eq:40}} \le |R-b-r\cos{\phistar}|+2(R-b)\sin^{2}{(\frac{\theta_1}{2}+\frac{\Phi_1}{2})}+2R|\sin(\frac{\Phi_1}{2})||\sin{(\theta_1+\frac{\Phi_1}{2})}|\\
&\overset{\eqref{eq:41}} <\frac{r^2\sin{\phistar}}{R}+2(R-b)\sin^{2}{(\frac{\theta_1}{2}+\frac{\Phi_1}{2})}+2R|\sin(\frac{\Phi_1}{2})||\sin{(\theta_1+\frac{\Phi_1}{2})}|\\
&\overbracket{<}^{\quad\mathllap{0<R-b<r,0<\theta_1+\Phi_1<\theta_1}} \frac{r^2\sin{\phistar}}{R}+2r\sin^2{(\theta_1)}+2R|\sin(\frac{\Phi_1}{2})||\sin{(\theta_1+\frac{\Phi_1}{2})}|\\
&\overbracket{<}^{\quad\mathllap{|\Phi_1|<\Phistar \text{ and trigonometry}}} \frac{r^2\sin{\phistar}}{R}+2r\sin^2{(\theta_1)}+2R\sin{(\Phistar)}(|\sin{(\theta_1)}|+|\sin{(\frac{\Phi_1}{2})}|)
\end{aligned}
\end{equation}

Therefore,
\begin{equation}\label{eq:43}
\begin{aligned}
|\cos{(\phi_2-\theta_1-\Phi_1)}-\cos{\phistar}| &\overset{\eqref{eq:42}}<\frac{r\sin{\phistar}}{R}+2\sin^2{\theta_1}+2\sin{\phistar}(|\sin{(\theta_1)}|+|\sin{(\frac{\Phi_1}{2})}|)\\
&<(r/R)+2\cdot(4r/R)+2(\sqrt{4r/R}+r\sin{\phistar}/R)<(10r/R)+4\sqrt{r/R}\overbracket{<}^{\mathclap{R>100r}}5\sqrt{r/R}\\
\end{aligned}
\end{equation}


The hyperbolicity condition \eqref{eqJZHypCond} implies $R>\frac{147r}{\sin^2\phistar}\overbracket{>}^{\qquad\mathllap{\substack{1>\cos(\phistar/4),\ 1>\cos{(\phistar/2)},\\147>16,\ \sin^2\phistar<\sin\phistar}}}\frac{16r\cos{(\phistar/4)}\cos{(\phistar/2)}}{\sin{\phistar}}\overbracket{=}^{\mathrlap{\substack{4\sin{(\phistar/4)}\cos{(\phistar/4)}\cos{(\phistar/2)}\\=2\sin{(\phistar/2)}\cos{(\phistar/2)}=\sin\phistar}}}\frac{4r}{\sin{(\phistar/4)}}$.

Since $\phi_2\in (\phistar,\frac{\pi}{2})$, we have ${\pi}/{2}\overset{\mathclap{0<\theta_1+\Phi_1}}>\phi_2-\theta_1-\Phi_1\overset{\mathclap{\phi_2>\phistar,\  \theta_1>\Phi_1}}>\phistar-2\theta_1\ge\phistar-2\sin^{-1}{\sqrt{4r/R}}\overbracket{>}^{\mathclap{R>\frac{4r}{\sin{(\phistar/4)}}}}\frac{\phistar}{2}$.\COMMENT{Inconsistent about using \textbackslash overbracket versus \textbackslash overset.}


Then we apply the Mean-Value Theorem to the left-hand side of \eqref{eq:43}: there exists $\tilde{\phi}_2\in (\frac{\phistar}{2},\frac{\pi}{2})$ such that $|\sin{\tilde{\phi}_2}||\phi_2-\theta_1-\Phi_1-\phistar|=|\cos{(\phi_2-\theta_1-\Phi_1)}-\cos{\phistar}|<5\sqrt{r/R}$. Hence
\begin{equation}\label{eq:44}
\begin{split}
&|\phi_2-\theta_1-\Phi_1-\phistar|<\frac{5}{\sin{\tilde{\phi}_2}}\sqrt{r/R}<\frac{5}{\sin{(\phistar/2)}}\sqrt{r/R}.
\end{split}
\end{equation}
If $0<z<\frac{1}{10}$, then 
$\frac{\sin z}{z}>\frac{25}{26}$. 
If $R>1700r$, then $0<z\dfn\sin^{-1}{(\sqrt{4r/R})}<\sin^{-1}(\sqrt{4/1700})<\frac{1}{10}$, hence $\frac{\sqrt{4r/R}}{\sin^{-1}{(\sqrt{4r/R})}}=\frac{\sin{(\sin^{-1}{(\sqrt{4r/R})})}}{\sin^{-1}{(\sqrt{4r/R})}}>\frac{25}{26}$.

Therefore, from  \cref{eq:44} and using the triangle inequality we obtain\COMMENT{In the next equations, mathllap is actually a good choice---because there is nothing to the left of the "=" and a lot of space to the left.}
\begin{equation}\label{eq:45}
\begin{aligned}
|\phi_2-\phistar|& \le|\theta_1|+|\Phi_1|+|\phi_2-\theta_1-\Phi_1-\phistar|\\
&\overbracket{<}^{\quad\mathllap{0<\theta_1<\sin^{-1}{(\sqrt{4r/R})},|\Phi_1|<\Phistar, \eqref{eq:44}}} \sin^{-1}{(\sqrt{4r/R})}+\sin^{-1}{(\frac{r\sin\phistar}{R})}+\frac{5}{\sin{(\phistar/2)}}\sqrt{r/R}\\
& \overbracket{<}^{\quad\mathllap{0<\sin{(\frac{\phistar}{2})}<1,\frac{r\sin\phistar}{R}<\sqrt{4r/R}}}\frac{1}{\sin{(\phistar/2)}}\sin^{-1}{(\sqrt{4r/R})}+\frac{1}{\sin{(\phistar/2)}}\sin^{-1}{(\sqrt{4r/R})}+\frac{2.5}{\sin{(\phistar/2)}}\sqrt{4r/R}\\
& \overbracket{<}^{\quad\mathllap{\sqrt{4r/R}<\sin^{-1}{(\sqrt{4r/R})}}}\frac{4.5}{\sin{(\phistar/2)}}\sin^{-1}{(\sqrt{4r/R})}\overbracket{<}^{\frac{26}{25}\sqrt{4r/R}>\sin^{-1}{(\sqrt{4r/R})}}\frac{4.68}{\sin{(\phistar/2)}}\sqrt{4r/R}.
\end{aligned}
\end{equation}


On the other hand, for in \cref{def:fpclf}\COMMENT{Fix words} $x_1\in\MRin$ with $n_1=0$ that is, $x_1\in\MRin\cap \MRout$, by some elementary geometry\COMMENT{Not helpful} we have $\theta_2=\phi_2-\theta_1-\Phi_1$ (also see the reasoning for \cref{eqd2,fig:twoTangentLines}). Therefore, we obtain\COMMENT{The second step is mysterious. Explain. Not clear why the things above "<" are relevant.\small\Hrule Discuss?}
\begin{equation}\label{eq:46}
\begin{split}
|\theta_2-\phistar| =|\phi_2-\phistar-\theta_1-\Phi_1| \le |\phi_2-\phistar|+|\theta_1|+|\Phi_1| &< \frac{4.68}{\sin{(\phistar/2)}}\sqrt{4r/R}+\sin^{-1}{(\sqrt{4r/R})}+\sin^{-1}{(r\sin{\phistar}/R)}\\
&\overbracket{<}^{\quad\mathllap{\eqref{eq:45},0<\theta_1<\sin^{-1}{\sqrt{(4r/R)}},|\Phi_1|<\Phistar,0<\sin{(\phistar/2)}<1}} \frac{4.68}{\sin{(\phistar/2)}}\sqrt{\frac{4r}{R}}+\frac{2}{\sin{(\phistar/2)}}\sin^{-1}{(\sqrt{4r/R})}\\
&\overbracket{<}^{\quad\mathllap{\frac{26}{25}\sqrt{4r/R}>\sin^{-1}{(\sqrt{4r/R})}}}\frac{6.76}{\sin{(\phistar/2)}}\sqrt{4r/R}
\end{split}
\end{equation}
Inequalities \eqref{eq:45}, \eqref{eq:46} imply\COMMENT{Bad display. Also, what is the symmetry? $\|x_2-Iy_*\|=\|x_0-x_*\|$?} 
\begin{equation}\label{eq:47}
    \begin{aligned}
        \|x_2-Iy_*\|=\sqrt{(\theta_2-\phistar)^2+(\phi_2-\phistar)^2}<&\frac{16.46}{\sin{(\phistar/2)}}\sqrt{r/R},\\
        \text{so by symmetry, } \|x_0-x_*\|<&\frac{16.46}{\sin{(\phistar/2)}}\sqrt{r/R}.
    \end{aligned}
\end{equation}
Here we have assumed $x_1=(\Phi_1,\theta_1)\in\Mrin$ with $\Phi_1\in [-\Phistar,+\Phistar], \theta_1\in[\sin^{-1}{(2r/R)},\sin^{-1}{(\sqrt{4r/R})}]$. And if we alternatively assume $\theta_1\in [\pi-\sin^{-1}{(\sqrt{4r/R})},\pi-\sin^{-1}{(2r/R)}]$, then by the symmetry argument we have $\|x_0-y_*\|<\frac{16.46}{\sin{(\phistar/2)}}\sqrt{r/R}$ and $
\|x_2-Ix_*\|<\frac{16.46}{\sin{(\phistar/2)}}\sqrt{r/R}$.\COMMENT{\color{red}Are we still in \Cref{TEcase_a1} Case (a1)? Is this the end of the proof for that case? It is completely unclear why it would.}

Then\COMMENT{This sounds like we are continuing in the same case as before. Look up the meaning of "then." Something here is wrong or at least confusing.} in case (b) and case (c) for \cref{TEcase_b,TEcase_c}, by symmetry and \cite[lemma 3.2, Equations (3.14), (3.21)]{hoalb}\COMMENT{Another useless pile of references instead of explanations in what follows. This makes the work unverifiable.} we have the following results.

We have in case (b), if $\theta_1<\sin^{-1}(2r/R)$ then\COMMENT{Below, mathclap aould be better.}

\begin{equation}\label{eq:casebcontractionregionx2_1}
\begin{aligned}
    \|x_2-Iy_*\|&\overbracket{<}^{\mathllap{\text{\cite[equation (3.14)]{hoalb}}}}\frac{14.6r}{R\sin{\phistar}}\overset{R>r}<\frac{16.46}{\sin{(\phistar/2)}}\sqrt{r/R},\\
    \text{ and }\|x_0-x_*\|&\overbracket{<}^{\mathllap{\text{\cite[equation (3.14)]{hoalb}}}}\frac{14.6r}{R\sin{\phistar}}\overset{R>r}<\frac{16.46}{\sin{(\phistar/2)}}\sqrt{r/R}.
\end{aligned}
\end{equation}
In case (b), if $\theta_1>\pi-\sin^{-1}(2r/R)$, then 
\begin{equation}\label{eq:casebcontractionregionx2_2}
\begin{aligned}
\|x_2-Ix_*\|&\overbracket{<}^{\mathllap{\text{\cite[equation (3.14)]{hoalb}}}}\frac{14.6r}{R\sin{\phistar}}\overset{R>r}<\frac{16.46}{\sin{(\phistar/2)}}\sqrt{r/R},\\ 
\text{ and }\|x_0-y_*\|&\overbracket{<}^{\mathllap{\text{\cite[equation (3.14)]{hoalb}}}}\frac{14.6r}{R\sin{\phistar}}\overset{R>r}<\frac{16.46}{\sin{(\phistar/2)}}\sqrt{r/R}.
\end{aligned}
\end{equation}
In case (c), if $\theta_1<\sin^{-1}(2r/R)$, then 
\begin{equation}\label{eq:caseCNinNout1}
\begin{aligned}
\|x_2-Iy_*\|&\overbracket{<}^{\qquad\mathllap{\text{\cite[equation (3.21)]{hoalb}}}}\frac{5.84r\sin{\phistar}}{R}<\frac{14.6r}{R\sin{\phistar}}\overset{R>r}<\frac{16.46}{\sin{(\phistar/2)}}\sqrt{r/R},\\
\text{ and }\|x_0-x_*\|&\overbracket{<}^{\qquad\mathllap{\text{\cite[equation (3.21)]{hoalb}}}}\frac{5.84r\sin{\phistar}}{R}<\frac{14.6r}{R\sin{\phistar}}\overset{R>r}<\frac{16.46}{\sin{(\phistar/2)}}\sqrt{r/R}.
\end{aligned}
\end{equation}
In case (c), if $\theta_1>\pi-\sin^{-1}(2r/R)$, then 
\begin{equation}\label{eq:caseCNinNout2}
    \begin{aligned}
    \|x_2-Ix_*\|&\overbracket{<}^{\qquad\mathllap{\text{\cite[equation (3.21)]{hoalb}}}}\frac{5.84r\sin{\phistar}}{R}<\frac{14.6r}{R\sin{\phistar}}\overset{R>r}<\frac{16.46}{\sin{(\phistar/2)}}\sqrt{r/R},\\
    \text{ and }\|x_0-y_*\|&\overbracket{<}^{\qquad\mathllap{\text{\cite[equation (3.21)]{hoalb}}}}\frac{5.84r\sin{\phistar}}{R}<\frac{14.6r}{R\sin{\phistar}}\overset{R>r}<\frac{16.46}{\sin{(\phistar/2)}}\sqrt{r/R}.
    \end{aligned}
\end{equation}

Thus, for all three cases (a1),(b) and (c)  in  \eqref{eqMainCases}, we can choose the desired neighborhood radius to be $\frac{16.46}{\sin{\phistar/2)}}\sqrt{r/R}$. Thus,\COMMENT{????} we have established \cref{contraction_region}: Cases (a1), (b) and (c) 
arise only when $x_0\in\Nout$ and $x_2\in\Nin$.\COMMENT{This is far from clear. Explicitly connect ``$x_0\in\Nout$ and $x_2\in\Nin$'' with the inequalities. Put a preview of this final argument at the start of the proof to explain what will be shown. Edit the proof so it is always clear which case is being considered.}
\end{proof}

\subsection{Coordinate change and expansion in the invariant cone by shears}This subsection discusses the definitions and calculations from \cite{hoalb} and \cite{cb}.\COMMENT{\color{red}Explain why these are being discussed here. Where and how is this eventually put to use? In short, what is the purpose of this subsection?\Hrule\small Many comments here!}

Suppose $x=(\phi,\theta)\in M_r,dx=(d\phi,d\theta)\in T_x(M)$, In the arc-length parameter $s$ and angle with inward normal vector parameter $\varphi$, the coordinate change from $(\phi,\theta)$ to $(s,\varphi)$ is given by

\begin{equation}\label{eq:coordinateschange}
\begin{aligned}
& \theta = \frac{\pi}{2}-\varphi\\
& d\varphi = -d\theta\\
& ds = rd\phi
\end{aligned}
\end{equation}

Then \cref{def:SectionSetsDefs} introduced \cite[page 58]{cb} p-metric $\|dx\|_{p}=\cos{\varphi}|ds|=\sin{\theta}|ds|$.

\cite[Equations (3.39), (3.40)]{cb} give the mirror equation of before/after collision infinitesimal wave front curvatures
$\mathcal{B}^{\pm}$ and the expansion of tangent vector in p-metric.\COMMENT{Fix this sentence}\quad\COMMENT{The next display is a random jumble of symbols with no context at all. Is the previous sentence meant to introduce this? Then the period at its end must be moved to the end of the last of these equations. \textcolor{red}{Every sentence must end with a period, and every period ends a sentence. Every formula must in some way be part of a sentence, so there must almost always be punctuation. You basically never have any, and you must learn how to use it in order to be able to communicate mathematics.} In the case ofthese three equations, do the preceding words introduce them? They seem to be part of the sentence which continues after them, but commas are still required.}
\begin{equation}\label{eq:MirrorEquationExpansionFormula}
\begin{aligned}
-\frac{d\theta}{rd\phi}\overset{\eqref{eq:coordinateschange}}=&\frac{d\varphi}{ds}=\mathcal{V}=\mathcal{B}^{-}\cos{\varphi}+\mathcal{K}=\mathcal{B}^{+}\cos{\varphi}-\mathcal{K}\\
\mathcal{B}^{+}=&\mathcal{B}^{-}+\frac{2\mathcal{K}}{\cos{\varphi}}\\
\frac{\|D\mathcal{F}_x(dx)\|_\p}{\|dx\|_\p}=&|1+\tau\mathcal{B^+}|
\end{aligned}
\end{equation}
where $\mathcal{K}$ is the billiard boundary's curvature\COMMENT{Search the source code for "'s". Almost all instances should be edited out. In the present case, "the billiard boundary's curvature" should become "the curvature of the billiard boundary". Again, this is not an isolated comment. If you find any occurrence of "'s" which you don't know how to fix, flag it so we can discuss it.} at the collision point $x$. $\tau$ is the distance between $p(x)$ and $p(\mathcal{F}(x))$ in billiard table.\COMMENT{\textcolor{red}{You sometimes forget to use math mode; a thing to fix.}} In the case if $x\in M_r$, then $\mathcal{K}=-\frac{1}{r}$.

Next we state three propositions and two corollaries for expansions of \COMMENT{What are expansions of 1?} a series ($\ge1$ times) collisions on focus circle\COMMENT{What is a focus circle? Maybe after you explain this sentence we can replace it by one that explains what is going on.}. They are presented in \cite[Chapter 8]{cb} for the Bunimovich Stadium Billiard. We now prove them for lemon billiards.\COMMENT{It does not look like we do. In \cite{cb} they are proved for the stadium, so we must exercise great care in establishing them for the present context. Each must come with proof. \bfseries UNLESS it is clear that they are established in \cite{cb} in greater generality even though stadia are the context.}

\begin{proposition}\label{proposition:ShearOneStepExpansion}
For nonsingular $(\phi,\theta)=x\in M_r$, $\mathcal{F}(x)\in M_r$ with $dx$ a tangent vector at $x$, if $dx=(d\phi,d\theta)\in \big\{(d\phi,d\theta)\bigm|\frac{d\theta}{d\phi}\in\big[0,+\infty\big]\big\}$, then $\frac{\|D\mathcal{F}_x(dx)\|_\p}{\|dx\|_\p}\ge1$
\end{proposition}
\begin{proof}
    See \cite[Sections 8.2, 8.3]{cb}.\COMMENT{We must say something about why this holds for lemon billiards. After all, we just sais "We now prove them for lemon billiards"!}
\end{proof}

\begin{corollary}
For $(\phi,\theta)=x\in M_r,\mathcal{F}(x)\in M_r,dx=(d\phi,d\theta)$ such that  $\frac{d\theta}{d\phi}\in[c,+\infty]$\COMMENT{Fix} for some $c>0$, then $\frac{\|D\mathcal{F}_{x}(dx)\|_\p}{\|dx\|_\p}>1+\delta$ for some $\delta>0$
\end{corollary}
\begin{proof}
For $\tau=2r\cos{\varphi}$,\COMMENT{Is this a new definition ($\tau$ appeared previously!!) or a reminder (then tell people where this first appeared) or a new fact (then tell people why it is true). The way it is worded, it is a case distinction or assumption.}
\[\frac{\|D\mathcal{F}_x(dx)\|_\p}{\|dx\|_\p}=|1+\tau \Bcal^+|=|1+2r\cos{\varphi}\frac{\frac{-d\theta}{rd\phi}-\frac{1}{r}}{\cos{\varphi}}|=|-1-2\frac{d\theta}{d\phi}|.\qedhere\]
\end{proof}
\NEWPAGE\begin{proposition}\label{proposition:ImaginaryExpansion}
{
Now let\COMMENT{This is not how a proposition starts} an infinitesimal wave front move along the trajectory experiencing k+1 successive collisions with focusing arc $\Gamma_r$.\COMMENT{Edit} Denote by $\mathcal{B}^{-}_{0}$ for before the very first and $\mathcal{B}^{-}_{i}$ before the $i$-th collision of its curvatures.\COMMENT{Edit} Denote by $\tau$ the consecutive focusing arc collisions distance in table.\COMMENT{Edit} Then the expansion of the front in p-metric, in the course of the whole series, i.e., the expansion between the very first to the very last\COMMENT{last what? collision?} is\COMMENT{Is $\mathcal{J}$ used anywhere else? Why not $\frac{\|d\mathcal{F}(dx)\|_\p}{\|dx\|_\p}$?}
\begin{equation}
    \mathcal{J}=\prod_{i=0}^{k-1}|1+\tau(-\frac{4}{\tau}+\mathcal{B}_{i}^{-})|
\end{equation}
And
\begin{equation}\label{eqJfromCM}
    \mathcal{J}=\mathcal{J}_{im}=\big|1-k\tau(-\frac{2}{\tau}+\mathcal{B}^{-}_0)\big|
\end{equation}
}\end{proposition}
\begin{proof}
The proof is in \cite[equation (8.12),  Excercise 8.29]{cb} using the idea of an imaginary billiard.\COMMENT{We must say something about why this holds for lemon billiards.\newline Why mention imaginary billiards?}
\end{proof}
\begin{proposition}[{\cite[Excercise 8.30]{cb}}]\label{proposition:circleserieslength}There\COMMENT{We must say something about why this holds for lemon billiards.} is a positive constant $c_r$ (depending on the length of $\Gamma_r$, hence on $r$, $\phistar$) such that in any series of $k+1\ge2$ successive collisions with $\Gamma_r$ the interval ($k\tau$)\COMMENT{Why the parentheses? $k\tau$ is a number, not an interval.} between the very first and very last collision always exceeds $c_r$.\COMMENT{Incomprehensible, edit! Does the \emph{length} of the interval exceed\dots?}
\end{proposition}
\begin{corollary}\label{corollary:shearlowerboundexpansion}
For $(\phi,\theta)=x\in {\Min_{r,k}}$, $k\ge2$, $(d\phi,d\theta)=dx \in T_xM$ such that $\frac{d\theta}{d\phi}\in[c,+\infty]$ for some $c>0$, then $\mathcal{F}^{k}(x)\in \Mrout$, $\mathcal{F}^{k-1}(x)\in \mathcal{F}^{-1}(\Mrout\smallsetminus\Mrin)$, then $\frac{\|D\mathcal{F}_{x}^{k-1}(dx)\|_\p}{\|dx\|_\p}\ge1+\frac{c_rc}{2}$\COMMENT{Incomprehensible, edit!} with $c_r$ from \cref{proposition:circleserieslength}.\COMMENT{Rewrite as if-then statement to be clear.}
\end{corollary}
\begin{proof}
By \cref{proposition:ImaginaryExpansion}, mirror equation in \eqref{eq:MirrorEquationExpansionFormula},\COMMENT{\textbackslash eqref throughout, not just here} and $\tau=2r\cos{\varphi}=2r\sin{\theta}$,\COMMENT{Connect $\mathcal{J}$ and $\frac{\|d\mathcal{F}^{k-1}(dx)\|_\p}{\|dx\|_\p}$ first!}
\begin{align*}
\frac{\|D\mathcal{F}^{k-1}(dx)\|_\p}{\|dx\|_\p}\overset{\eqref{eqJfromCM}}=&|1-(k-1)\tau(-2/\tau+B^{-})|
\overset{\eqref{eq:MirrorEquationExpansionFormula}}=\Big|1-(k-1)\tau\Big(-\frac{2}{\tau}+\frac{\frac{-d\theta}{rd\phi}+\frac{1}{r}}{\cos{\varphi}}\Big)\Big|\\
\overset{\mathclap{\tau=2r\cos{\varphi}}}{\scalebox{6}[1]{$=$}}&\Big|1-(k-1)\tau\Big(-\frac{2}{2r\cos{\varphi}}+\frac{\frac{-d\theta}{rd\phi}+\frac{1}{r}}{\cos{\varphi}}\Big)\Big|=\Big|1+(k-1)\tau\frac{{d\theta}/{d\phi}}{r}\Big|.
\end{align*}
The reasoning in \cref{proposition:circleserieslength} also shows that there exists a constant $\overline{c}_r=\frac{1}{2}c_r$ such that for all $k\ge2$, the interval $(k-1)\tau$\COMMENT{This is a number, not an interval} between the first and next to last collisions always exceeds $\overline{c}_r$.\COMMENT{What does "exceeds" mean here?} Thus the conclusion follows.
\end{proof}
\section
{Expansion in Case (a)}\label{sec:ProofTECase_a0_a1}
\strut\COMMENT{The title is quite descriptive, but it would still help to explain what happens where in this section}
\subsection[\texorpdfstring{$\mathcal{I}$}{I} A Formula for expansion Lower Bound in \texorpdfstring{\cref{TEcase_a0,TEcase_a1}}{case a0,a1}]{A Formula for expansion Lower Bound in Case (a)}\label{subsec:AFELB}\strut\COMMENT{Say what we do in this subsection and where the formula is.}
\begin{notation}\label{def:AformularForLowerBoundOfExpansion}
For $x_0\in\Mrout$, $x_1\in\MRin$, $x_2\in\Mrin$, $n_1$ from \cref{def:MhatReturnOrbitSegment} and the free path and chord length functions  $\tau_0$, $\tau_1$, $d_0$, $d_1$, $d_2$\COMMENT{Is there a picture as well? Might not hurt to make one if not.\newline\wentao{still need to refer \cref{fig:heuristicproof}}} from \cref{def:fpclf}, considered as functions of $x_1\in\big\{(\Phi_1,\theta_1)\in\MRin\bigm|\sin\theta_1\ge2r/R\big\}$ (see \cref{fig:R1,fig:A_compactRegion}), let
\[\mathcal{I}\dfn-1+\frac{\tau_1}{d_0}\Big[\frac{2(\tau_0-d_0)}{d_1}-\frac{\tau_0+\tau_1-2d_0}{\tau_1}\Big].\] 
Note that if $d_1\ge2r$, then $n_1=0$, and if $n_1=0$, then $\tau_0+\tau_1>2d_i$, $i=0,2$. See \cref{fig:lengthfunctonFigure,remark:symmetrywiththetangentline}.\COMMENT{This sentence makes no sense. Is the last part trying to describe \cref{fig:A_compactRegion}?\wentao{OK, there was some wrong reference. Now it's correct referencing the already mentioned remark and figure.}}
\end{notation}
\NEWPAGE
\begin{lemma}[Orbit configuration and expansion estimate in Case (a) of \eqref{eqMainCases}]\label{lemma:ExpansionLowerBoundA0A1}
Suppose a nonsingular $x\in \hat{M}$  (\cref{def:SectionSetsDefs}) has a return orbit segment as in \cref{def:MhatReturnOrbitSegment}, which is $x,\mathcal{F}(x),\cdots,\mathcal{F}^{\sigma(x)}(x)=\hat{F}(x)\in \hat{M}$ with $x_{0}\in\Mrout$, $x_1=\mathcal{F}(x_0)=(\phi_1,\theta_1)\in \MRin$, $x_2\in\Mrin$\COMMENT{Isn't this so far simply recalling notation instead of something we suppose.} with $\sin\theta_1\ge 2r/R$\COMMENT{Isn't this the defining inequality for Case (a)?} in cases (a0) and (a1) of \eqref{eqMainCases} and context of \cref{TEcase_a0,TEcase_a1}.\COMMENT{What does the "context" part add? \color{red}Isn't it just "in Case (a) of \eqref{eqMainCases}"?\Hrule\small Also, is there any actual assumption so far or simply a long statement about what Case (a) is???}

With $x\in\hat{M}=(\Mrin\cap\Mrout)\sqcup\mathcal{F}^{-1}(\Mrout\setminus\Mrin)$  (\cref{def:SectionSetsDefs}),\COMMENT{Too late to recall definition of $\hat M$; see first line.} suppose that 
\[dx=(d\phi,d\theta)\in\Big\{(d\phi,d\theta)\bigm|\frac{d\theta}{d\phi}\in\big[0,+\infty\big]\Big\}\]
is a tangent vector at $x$, i.e., $dx$ is in the first and\COMMENT{``or''??} third quadrants in $\phi$, $\theta$ coordinates.\COMMENT{Change all occurrences of "$\phi$, $\theta$ coordinates" to "$\phi\theta$-coordinates"!} Then: 
\begin{itemize} 
    \item If $x\in\Mrin\cap\Mrout$, then $dx_0\dfn dx$ and $dx_2\dfn D\mathcal{F}^2_x(dx)$. Otherwise, $x\in\mathcal{F}^{-1}(\Mrout\setminus\Mrin)$, then\COMMENT{??? Is there an "if" to match?} $dx_0\dfn D\mathcal{F}_x(dx)$ and $dx_2\dfn D\mathcal{F}_{x}^3(dx)$. In either case, $dx_0$ is a tangent vector at $x_0$ and $dx_2$ is a tangent vector at $x_2$. $dx_1\dfn D\mathcal{F}_{x_0}(dx_0)$ is a tangent vector at $x_1\in\MRin$.\COMMENT{This sentence sounds largely vacuous. What is the point?}
    \item $\frac{\|dx_0\|_\p}{\|dx\|_\p}\ge 1$ and $\frac{\|D\hat{F}(dx)\|_\p}{\|dx_2\|_\p}\ge 1$.
    \item If $\mathcal{I}<0$ (\cref{def:AformularForLowerBoundOfExpansion}), then $\frac{\|dx_2\|_\p}{\|dx_0\|_\p}\ge|\mathcal{I}|$.
\end{itemize}
\end{lemma}
\begin{proof}
Since $dx$ is in the first and\COMMENT{And?? Only the origin is in both!} third quadrant, by the definition of $\hat{M}$  (\cref{def:SectionSetsDefs}) and the return orbit segment (\cref{def:MhatReturnOrbitSegment}), if $x\in\Mrin\cap\Mrout$, then $x_0=x,dx_0\dfn dx$ and otherwise $x\in\mathcal{F}^{-1}(\Mrout\setminus\Mrin)$, we have\COMMENT{Something wrong here} $x_0=\mathcal{F}(x)$, $D\mathcal{F}_{x}(dx)=dx_0=(d\phi_0,d\theta_0)$ and $D\mathcal{F}_{x}=\begin{pmatrix}
1 & 2\\ 0 & 1
\end{pmatrix}$ in $\phi,\theta$ coordinate\COMMENT{Plural (always) because there are two variables}. Thus for either $x\in\Mrin\cap\Mrout$ or $x\in\mathcal{F}^{-1}(\Mrout\setminus\Mrin)$, $dx_0=(d\phi_0,d\theta_0)$ is in the invariant quadrant $\{(d\phi,d\theta)|\frac{d\theta}{d\phi}\in[0,+\infty]\}$ of tangent space at $x_0$. Therefore, by \cref{proposition:ShearOneStepExpansion} we have $\frac{\|dx_0\|_\p}{\|dx\|_\p}\ge 1$.

As described in \eqref{eqPtsFromeqMhatOrbSeg}, case (a) of \eqref{eqMainCases} and \cref{caseA},\COMMENT{Does the reader really need to look at all three places?} we have $x_1=\mathcal{F}(x_0)$, $x_2=\mathcal{F}^2(x_0)$ and $D\mathcal{F}^2_{x_0}$ is a negative matrix in coordinates $\phi,\theta$. Hence $(d\phi_2,d\theta_2)=dx_2=D\mathcal{F}^2_{x_0}(dx_0)$ has $\frac{d\theta_2}{d\phi_2}>0$, that is $dx_2$ is also in the quadrant.

From \cref{def:MhatReturnOrbitSegment}, there are two cases for $m(x_2)$ from \cref{def:SectionSetsDefs}. Either $m(x_2)\le1$ or $m(x_2)\ge2$. If $m(x_2)\le1$, then $\hat{F}(x)=x_2\in\hat{M}$, $dx_2=D\hat{F}_x(dx)$ and $\frac{\|D\hat{F}_x(dx)\|_\p}{\|dx_2\|_\p}=1$. If $m(x_2)\ge2$, then for $k=m(x_2)$, $(\phi_2,\theta_2)=x_2\in \textcolor{blue}{M^{\text{\upshape{in}}}_{r,k}}$,  applying \cref{corollary:shearlowerboundexpansion}, we have $\hat{F}(x)=\mathcal{F}^{k-1}(x_2)\in\hat{M}$ and $\frac{\|D\hat{F}_x(dx)\|_\p}{\|dx_2\|_\p}=\frac{\|D\mathcal{F}^{k-1}_{x_2}(dx_2)\|_\p}{\|dx_2\|_\p}>1$.

To analyze $\frac{\|dx_2\|_\p}{\|dx_0\|_\p}$ we use the same notation as in \eqref{eq:MirrorEquationExpansionFormula} with subscripts. We denote by $\mathcal{B}^{\pm}_0,\mathcal{B}^{\pm}_1$ the infinitesimal wave front curvatures after/before collisions at $x_0$, $x_1$, respectively.
Then
\begin{equation}\label{eqB+d0}
\mathcal{B}^{+}_0\overset{\eqref{eq:MirrorEquationExpansionFormula}}=\frac{\frac{-d\theta_0}{rd\phi_0}-\frac{1}{r}}{\cos{\varphi_0}}\le-\frac{1}{r\cos{\varphi_0}}\overset{\eqref{eq:coordinateschange}}=-\frac{1}{r\sin{\theta_0}}\overset{\text{\tiny\cref{def:fpclf}}}{\scalebox{4}[1]{$=$}}-\frac1{d_0}<0.
\end{equation}
\begin{remark}[Was commented out. Why? \eqref{eqBBr} is needed.]
\eqref{eq:MirrorEquationExpansionFormula} and \cref{def:fpclf} give
\begin{equation}\label{eqBBd}
\mathcal{B}^{+}_1\overset{\text{\eqref{eq:MirrorEquationExpansionFormula}}}=\mathcal{B}^{-}_1-\frac{2}{d_1}.
\end{equation}
\cite[equation (3.31)]{cb} gives
\begin{equation}\label{eqBBr}
\frac{1}{\mathcal{B}_{1}^{-}}\overset{\text{\cite[equation (3.31)]{cb}}}=\frac{1}{\mathcal{B}_{0}^{+}}+\tau_0
\end{equation}
\end{remark}
\(
\frac{\|dx_1\|_\p}{\|dx_0\|_\p}=|1+\tau_0\mathcal{B}^{+}_0|\), which in turn implies by reindexing
\[\frac{\|dx_2\|_\p}{\|dx_1\|_\p}=|1+\tau_1\mathcal{B}^{+}_1|\overbracket{=}^{\mathclap{\mathcal{B}^{+}_1\overset{\text{\eqref{eq:MirrorEquationExpansionFormula}}}=\mathcal{B}^{-}_1-\frac{2}{d_1}}}\big|1+\tau_1(\mathcal{B}^{-}_1-\frac2{d_1})\big|\overbracket{=}^{\frac{1}{\mathcal{B}_{1}^{-}}\overset{\mathclap{\bigstrut\text{\cite[equation (3.31)]{cb}}}}=\frac{1}{\mathcal{B}_{0}^{+}}+\tau_0}\Big|1+\tau_1(\frac{1}{\frac{1}{\mathcal{B}^{+}_0}+\tau_0}-\frac{2}{d_1})\Big|=\Big|1+\frac{\mathcal{B}^{+}_0\tau_1}{1+\mathcal{B}^{+}_0\tau_0}-\frac{2\tau_1}{d_1}\Big|.
\]
Therefore
\begin{equation}\label{eq:dx2overdx0}
\frac{\|dx_2\|_\p}{\|dx_0\|_\p}=\frac{\|dx_2\|_\p}{\|dx_1\|_\p}\frac{\|dx_1\|_\p}{\|dx_0\|_\p} =|1+\tau_0\mathcal{B}^{+}_0|\Big|1+\frac{\mathcal{B}^{+}_0\tau_1}{1+\mathcal{B}^{+}_0\tau_0}-\frac{2\tau_1}{d_1}\Big| =\Big|\underbracket{(\tau_0+\tau_1-\frac{2\tau_0\tau_1}{d_1})\mathcal{B}^{+}_0+1-\frac{2\tau_1}{d_1}}_{\nfd\mathcal{I}(\mathcal{B}^{+}_0)}\Big|
\end{equation}
In Case (a) of \eqref{eqMainCases} 
we have $d_1\ge2r$, so $\frac{2}{d_1}\le\frac{1}{r}\le\frac{1}{2r}+\frac{1}{2r}\le\frac{1}{\tau_0}+\frac{1}{\tau_1}$\COMMENT{Where do we first state that $\tau_i\le2r$?} (also see \cite[Proposition 3.6]{hoalb})\COMMENT{Why?}, hence $\tau_0+\tau_1-\frac{2\tau_0\tau_1}{d_1}\ge0$. This combined with \eqref{eqB+d0}\COMMENT{How used?} gives: 
\begin{equation}\label{eq:56}
\begin{split}
\mathcal{I}(\mathcal{B}^{+}_0) & \le(\tau_0+\tau_1-\frac{2\tau_0\tau_1}{d_1})(-\frac{1}{d_0})+1-\frac{2\tau_1}{d_1}\\
& =1-\frac{\tau_0+\tau_1}{d_0}+\frac{2\tau_0\tau_1}{d_0d_1}-\frac{2\tau_1}{d_1} = -1-\frac{\tau_0+\tau_1-2d_0}{d_0}+\frac{2\tau_1}{d_1}(\frac{\tau_0}{d_0}-1)\\
& =-1-\frac{\tau_0+\tau_1-2d_0}{d_0}+\frac{2\tau_1}{d_1}(\frac{\tau_0-d_0}{d_0}) = -1+\frac{\tau_1}{d_0}\Big[\frac{2(\tau_0-d_0)}{d_1}-\frac{\tau_0+\tau_1-2d_0}{\tau_1}\Big]=\mathcal{I}
\end{split}
\end{equation}
Hence \eqref{eq:56} combined with \eqref{eq:dx2overdx0} gives that if $\mathcal{I}<0$, then $\frac{\|dx_2\|_\p}{\|dx_0\|_\p}\ge|\mathcal{I}|$.
\end{proof}

\begin{proposition}\label{Prop:Ilowerbound}
For $\mathcal{I}$\COMMENT{Since this all depends on $x$, should this be $\mathcal I(x)$?} as in \Cref{def:AformularForLowerBoundOfExpansion}, if the orbit segment for $x\in\hat{M}$ defined in \cref{def:MhatReturnOrbitSegment,def:fpclf} with $x_0\in\Mrout$, $\mathcal{F}(x_0)=x_1=(\Phi_1,\theta_1)\in\MRin,\,\sin{\theta_1}\ge\sqrt{4r/R}$, $x_2\in\Mrin$ in the context of \cref{TEcase_a0}, i.e. case (a0) in \eqref{eqMainCases}, then $\mathcal{I}<-1$, and hence $\frac{\|dx_2\|_\p}{\|dx_0\|_\p}>1$ by \cref{lemma:ExpansionLowerBoundA0A1}.
\COMMENT{Rewrite as if-then statement for clarity\wentao{Done}\color{red}\boris BUT IT IS NOT CLEAR WHAT THE IF IS. Maybe we have to discuss this.\Hrule\small Needs a lot of work}
\end{proposition}
\begin{proof}
We invoke  \cref{prop:casea0largetheta,prop:casea0smallthetaR2,prop:casea0smallthetaR1} as follows (see also \cref{fig:R1_fraction_extreme}).

If $\sin\theta_1>1/2$, then $\mathcal{I}<-1.37$ (\cref{prop:casea0largetheta} below). If $\sin\theta_1\ge\sqrt{4r/R}$ and $x_1\in R_2$ in \cref{fig:R1}, then $\mathcal{I}<-1$ (\cref{prop:casea0smallthetaR2} below). If $1/2\ge\sin\theta_1\ge\sqrt{4r/R}$ and $x_1\in R_1$ in \cref{fig:R1}, then $\mathcal{I}<-1$ (\cref{prop:casea0smallthetaR1} below). Then, by symmetry on the entire region, $\mathcal{I}<-1$ for $x_1\in\big\{(\Phi_1,\theta_1)\in\Mrin\bigm|\sin\theta_1\ge\sqrt{4r/R}\big\}$. 
\end{proof}

\subsection{Proof of \texorpdfstring{\cref{Prop:Ilowerbound}}{Proposition of Lower bound for I} part I}\label{SBSProp53PartI} We now consider the case when $x$ return orbit segment\COMMENT{This sounds weird. I recommend to replace "$x$ return orbit segment" by "returning orbit segment of $x$" \color{red}throughout.} has a more transverse reflection on $\Gamma_R$, specifically, with $x_1=(\Phi_1,\theta_1), \frac{\pi}{6}<\theta_1<\frac{5\pi}{6}$.\COMMENT{"$x_1=(\Phi_1,\theta_1), \frac{\pi}{6}<\theta_1<\frac{5\pi}{6}$" should be "$x_1=(\Phi_1,\theta_1)$, $\frac{\pi}{6}<\theta_1<\frac{5\pi}{6}$"---look at the tex file to see the difference; \color{red}this is a lesson that should be applied throughout.} The goal is to prove that $\mathcal{I}<-1.37$ (\cref{prop:casea0largetheta}). Less transverse reflections will be considered starting on page \pageref{SBSProp53PartII}.
\begin{lemma}\label{lemma:returnorbittheta1Large}
 Given that $\phistar\in\big(0,\tan^{-1}(1/3)\big)$ and the lemon billiard configuration satisfies \eqref{eqJZHypCond}. Suppose $x\in \hat{M}$ the return orbit segment as defined in \cref{def:MhatReturnOrbitSegment,def:fpclf} has the from \eqref{eqPtsFromeqMhatOrbSeg} $x_0\in\Mrout,x_1=(\Phi_1,\theta_1)\in\MRin,x_2=(\phi_2,\theta_2)\in\Mrin$, we have the following two conclusions.\COMMENT{What's the difference between them?}
 
 If $\frac{\pi}{2}<\phi_2<\frac{3\pi}{2}$, then $\theta_1\in \big(\frac{\pi}{6},\frac{5\pi}{6}\big)$ and $\sin\theta_1>\sqrt{4r/R}$, $n_1=0$ in \eqref{eqPtsFromeqMhatOrbSeg}.

 If $\frac{\pi}{2}<\phi_0<\frac{3\pi}{2}$, then $\theta_1\in \big(\frac{\pi}{6},\frac{5\pi}{6}\big)$ and $\sin\theta_1>\sqrt{4r/R}$, $n_1=0$ in \eqref{eqPtsFromeqMhatOrbSeg}.
\end{lemma}
\begin{proof}
We will prove the first conclusion since the second conclusion proof will be just a symmetric argument of the first conclusion.

By symmetry we can assume $0<\theta_1\le\frac{\pi}{2}$ and we look into the \cref{fig:QonUpperHalfCircle} and use the Cartesian coordinate system defined in \cref{def:StandardCoordinateTable} so that in  \cref{def:BasicNotations} two centers $O_r=(0,0)$ and $O_R=(0,b)$. Then the lemon billiard has two corners coordinates $A=(-r\sin\phistar,-r\cos\phistar),B=(r\sin\phistar,-r\cos\phistar)$. Denote $C=(-r,0),D=(r,0)$ both on $C_r$. In \cref{def:fpclf}. 

For $x\in \hat{M}$ the return orbit segment as defined in \cref{def:fpclf} has $x_0\in\Mrout,x_1=(\Phi_1,\theta_1)\in\MRin,x_2=(\phi_2,\theta_2)\in\Mrin$, in \cref{fig:QonUpperHalfCircle}, denote by $P=p(\mathcal{F}^{-1}(x_2)),Q=p(x_2)$ the positions in billiard table. Then $\mathcal{F}^{-1}(x_2)=(\Phi,\theta_1)\in\MRout$ for some $\Phi\in(-\Phistar,\Phistar)$.

On the other hand $x_2=(\phi_2,\theta_2)\in\Mrin,\phi_2\in(\frac{\pi}{2},\frac{3\pi}{2})$ implies $Q$ is on the upper half circle of boundary $\Gamma_r$ in \cref{fig:QonUpperHalfCircle}. i.e. $Q$ has a positive $y-$coordinate. 

Since $P$ is an interior point of arc $\arc{AB}$, $Q$ is the interior point of the upper half-part of the circle $C_r$ in \Cref{fig:QonUpperHalfCircle}, line segment $\overline{PQ}$ must intercept the line segment $\overline{AB}$ at some $P_1$ and line segment $\overline{PQ}$ must intercept the line segment $\overline{CD}$ at some $P_2$, where $P_1$ is some interior point of line segment $\overline{AB}$ and $P_2$ is some interior point of line segment $\overline{CD}$. Note that the line $\overline{CD}$ going through center $O_r$ is parallel to line $\overline{AB}$. By some elementary geometry, $\measuredangle{BAD}=\frac{\pi}{4}-\frac{\phistar}{2}$ and $\theta_1+\Phi=\measuredangle{BP_1Q}$. 

In this coordinate system (\cref{def:StandardCoordinateTable}), write the points coordinates as $Q=(x_Q,y_Q),P=(x_P,y_P),P_1=(x_{P_1},y_{P_1}),P_2=(x_{P_2},y_{P_2}),A=(x_A,y_A),B=(x_B,y_B),D=(x_D,y_D)$. Note that $y_Q-y_P>y_{P_2}-y_{P_1}=y_D-y_A>0$ and if $x_Q>x_P$ then $0<x_Q-x_P<x_D-x_A$. Hence if $0\ge x_Q-x_P$ then the line $PQ$ has a negative slope or infinite slope. Otherwise if $x_Q>x_P$ then the line $PQ$ has a positive slope larger than the slope of line $AD$. Also see \cref{fig:QonUpperHalfCircle}. In either case, $\measuredangle{BP_1Q}>\measuredangle{BAD}$.
Therefore,\COMMENT{The alignment below is pointless; the equations are unrelated, so they should both be centered instead. See \url{https://www.overleaf.com/learn/latex/Aligning_equations_with_amsmath}}
\begin{equation}\label{eq:57}
\begin{aligned}
 \theta_1+\Phi &=\measuredangle{BP_1Q}>\measuredangle{BAD}=\frac{\pi}{4}-\frac{\phistar}{2}\\
 \theta_1 &>\frac{\pi}{4}-\frac{\phistar}{2}-\Phi\overbracket{>}^{\mathclap{\strut\Phi<\Phistar}}\frac{\pi}{4}-\frac{\phistar}{2}-\Phistar=\frac{\pi}{4}-\frac{\phistar}{2}-\sin^{-1}{(\frac{r\sin\phistar}{R})}\overbracket{>}^{\mathclap{\bigstrut\phistar\in\big(0,\tan^{-1}(1/3)\big), R>1700r}}\frac{\pi}{6}
\end{aligned}
\end{equation}
Here we have assumed $\theta_1\in(0,\frac{\pi}{2}]$. 

When we alternatively assume $\theta_1\in[\frac{\pi}{2},\pi)$, by symmetry we can show $\pi-\theta_1-\Phi=\measuredangle{AP_1Q}>\measuredangle{ABC}=\frac{\pi}{4}-\frac{\phistar}{2}$. Similar calculations as in \eqref{eq:57} yield $\theta_1<\frac{5\pi}{6}$.

Since $\sin\theta_1>0.5\overbracket{>}^{\mathclap{R>1700r}}\sqrt{4r/R}>\frac{r\sin\phistar}{R}$, we have $n_1=0$ in the return orbit segment of $x$ from \cref{def:fpclf}, that is, $x_1=\mathcal{F}^{-1}(x_2)$.
\end{proof}

\begin{figure}[h]
\begin{minipage}{.5\textwidth}
\begin{center}
\begin{tikzpicture}[xscale=0.6, yscale=0.6]
    \tkzDefPoint(0,0){Or};
    \tkzDefPoint(0,-4.8){Y};
    \pgfmathsetmacro{\rradius}{4.5};
    \pgfmathsetmacro{\bdist}{4.7};
    \tkzDefPoint(0,\bdist){OR};
    \pgfmathsetmacro{\phistardeg}{40};
    \pgfmathsetmacro{\XValueArc}{\rradius*sin(\phistardeg)};
    \pgfmathsetmacro{\YValueArc}{\rradius*cos(\phistardeg)};
    \tkzDefPoint(\XValueArc,-1.0*\YValueArc){B};
    \tkzDefPoint(-1.0*\XValueArc,-1.0*\YValueArc){A};
    \pgfmathsetmacro{\Rradius}{veclen(\XValueArc,\YValueArc+\bdist)};
    \tkzDrawArc[name path=Cr,very thick](Or,B)(A);
    \tkzDrawArc[name path=CR,very thick](OR,A)(B);
    \pgfmathsetmacro{\PXValue}{\Rradius*sin(10)};
    \pgfmathsetmacro{\PYValue}{\bdist-\Rradius*cos(10)};
    \tkzDefPoint(\PXValue,\PYValue){P};
    \pgfmathsetmacro{\QXValue}{\rradius*cos(10)};
    \pgfmathsetmacro{\QYValue}{\rradius*sin(10)};
    \tkzDefPoint(\PXValue,\PYValue){P};
    \tkzDefPoint(\QXValue,\QYValue){Q};
    \tkzDefPoint(-\rradius,0){C};
    \tkzDefPoint(\rradius,0){D};
    \path[name path=linePQ](P)--(Q);
    \begin{scope}[ultra thin,decoration={
    markings,
    mark=at position 0.55 with {\arrow{>}}}
    ] 
    \draw[blue,postaction={decorate}](P)--(Q);
    \end{scope}
    \draw[name path=lineCD,red,dashed,ultra thin](C)--(D);
    \draw[name path=lineAB,red,dashed,ultra thin](A)--(B);
    \draw[name path=lineORP,red,dashed,ultra thin](OR)--(P);
    \draw[name path=lineOrA,red,dashed,ultra thin](Or)--(A);
    \draw[name path=lineAB,red,dashed,ultra thin](B)--(A);
    \draw[name path=lineOROr,red,dashed,ultra thin](OR)--(Or);
    \draw[name path=lineAD,red,dashed,ultra thin](A)--(D);
    \path[name intersections={of=lineAB and linePQ,by={P1}}];
    \path[name intersections={of=lineCD and linePQ,by={P2}}];
    \tkzLabelPoint[above](OR){$O_R$};
    \tkzLabelPoint[above right](Q){$Q$};
    \tkzLabelPoint[below](P){$P$};
    \tkzLabelPoint[above](Or){$O_r$};
    \tkzLabelPoint[left](A){$A$};
    \tkzLabelPoint[right](B){$B$};
    \tkzLabelPoint[left](C){$C$};
    \tkzLabelPoint[right](D){$D$};
    \tkzLabelPoint[above left](P2){$P_2$};
    \tkzLabelPoint[above, yshift=1.1](P1){$P_1$};
    \tkzDrawPoints(A,B,C,D,Or,OR,P,Q,P1,P2);
    \tkzMarkAngle[arc=lll,ultra thin,size=0.6](C,Or,A);
    \tkzMarkAngle[arc=lll,ultra thin,size=0.4](B,P1,Q);
    \tkzMarkAngle[arc=lll,ultra thin,size=0.8](Or,OR,P);
    \tkzMarkAngle[arc=lll,ultra thin,size=0.8](B,A,D);
    \node at ($(OR) + (-83:1.6)$)  { $\Phi_1=\Phi$};
    \node at ($(Or) + (-160:1.4)$)  {\Large $\frac{\pi}{2}-\phistar$};
    \node at ($(A) + (10:1.8)$)  {$\frac{\pi}{4}-\frac{\phistar}{2}$};
    \node at ($(P1) + (18:1.4)$)  {\small $\theta_1+\Phi$};
\end{tikzpicture}
\end{center}
\caption{$p(x_2)$ is on the upper half circle. $O_r=(0,0)$, $O_R=(0,b)$,
\newline $A=(-r\sin\phistar,-r\cos\phistar)$,
\newline $B=(r\sin\phistar,-r\cos\phistar)$.
}\label{fig:QonUpperHalfCircle}
\end{minipage}\hfill
\begin{minipage}{.5\textwidth}
\begin{center}
\begin{tikzpicture}[xscale=0.6, yscale=0.6]
    \tkzDefPoint(0,0){Or};
    \tkzDefPoint(0,-4.8){Y};
    \pgfmathsetmacro{\rradius}{4.5};
    \pgfmathsetmacro{\bdist}{8.7};
    \tkzDefPoint(0,\bdist){OR};
    \pgfmathsetmacro{\phistardeg}{40};
    \pgfmathsetmacro{\XValueArc}{\rradius*sin(\phistardeg)};
    \pgfmathsetmacro{\YValueArc}{\rradius*cos(\phistardeg)};
    \tkzDefPoint(\XValueArc,-1.0*\YValueArc){B};
    \tkzDefPoint(-1.0*\XValueArc,-1.0*\YValueArc){A};
    \pgfmathsetmacro{\Rradius}{veclen(\XValueArc,\YValueArc+\bdist)};
    \clip(-\rradius-1.5,-\rradius-0.6) rectangle (\rradius+1.5,\bdist+0.4);
    \draw[name path=Cr,dashed,ultra thin] (Or) circle (\rradius);
    \tkzDrawArc[very thick](Or,B)(A);
    \tkzDrawArc[name path=CR,very thick](OR,A)(B);
    \pgfmathsetmacro{\PXValue}{\Rradius*sin(7)};
    \pgfmathsetmacro{\PYValue}{\bdist-\Rradius*cos(7)};
    \tkzDefPoint(\PXValue,\PYValue){P};
    \tkzDefPoint(0,\bdist-\Rradius){C};
    \draw[name path=PQ_Far,blue,ultra thin,dashed](P)--+(70:7)coordinate(PposX);
    \draw[name path=PT_Far,blue,ultra thin,dashed](P)--+(124:9);
    \draw[name path=PT0_Far,blue,ultra thin,dashed](P)--+(-56:0.7);
    \draw[name path=PQ0_Far,blue,ultra thin,dashed](P)--+(-110:1);
    \draw[red,dashed,->](P)--+(7:4)coordinate(PposX);
    \draw[red,dashed](P)--+(-173:6)coordinate(PnegX);
    \draw[name path=linePOR,red,dashed,->](P)--(OR);
    \draw[name path=linePOr,red,dashed](P)--(Or);
    \tkzDefPoint(0,-\rradius){Y};
    \draw[red,dashed](OR)--(Y);
    \tkzMarkAngle[arc=lll,ultra thin,size=1.7](Or,OR,P);
    \tkzMarkAngle[arc=lll,ultra thin,size=0.7](Y,Or,P);
    \tkzMarkAngle[arc=lll,ultra thin,size=0.8](OR,P,Or);
    
    \path[name intersections={of=Cr and PQ_Far,by={Q}}];
    \path[name intersections={of=Cr and PQ0_Far,by={Q0}}];
    \path[name intersections={of=Cr and PT0_Far,by={T0}}];
    \path[name intersections={of=Cr and PT_Far,by={T}}];
    \tkzMarkAngle[arc=lll,ultra thin,size=1.8](PposX,P,Q);
    \tkzMarkRightAngle(OR,P,PnegX);
    \node[] at ($(Or)+(-78:1)$)  {\Large $\alpha$};
    \node[] at ($(OR)+(-85:2.3)$)  {\Large $\Phi_1$};
    \node[] at ($(P)+(105:1.3)$)  {\Large {$\delta$}};
    \node[] at ($(P)+(35:2.3)$)  {\Large {$\theta_1$}};
    \tkzLabelPoint[left](A){$A$};
    \tkzLabelPoint[below](B){$B$};
    \tkzLabelPoint[above](C){$M$};
    \tkzLabelPoint[left](OR){$O_R$};
    \tkzLabelPoint[left](Or){$O_r$};
    \tkzLabelPoint[below left](P){$P$};
    \tkzLabelPoint[right](Q){$Q$};
    \tkzLabelPoint[left](T){$T$};
    \tkzLabelPoint[below](Q0){$Q_0$};
    \tkzLabelPoint[below](T0){$T_0$};
    \tkzLabelPoint[below](Y){$Y$};
    \tkzLabelPoint[below](PposX){$TP$};
    \tkzDrawPoints(A,Or,OR,P,B,Y,T,Q,T0,Q0,C);
    \begin{scope}[ultra thin,decoration={markings, mark=at position 0.55 with {\arrow{>}}}
    ] 
    \draw[blue,postaction={decorate}](P)--(Q);
    \draw[blue,postaction={decorate}](T)--(P);
    \end{scope}
\end{tikzpicture}
\caption{A coordinate system $TP-P-O_R$ originated at $P$}\label{fig:fractionEstimate2}
\end{center}
\end{minipage}
\end{figure}

Unlike \cref{def:StandardCoordinateTable}, we now define a new coordinate system for the billiard table with origin at the collision position $P\in\Gamma_R$.
\begin{definition}[Coordinate System with origin at $P$]\label{def:TPPORsystem}
    Based on the billiard table notation from \cref{def:LemonTableConfiguration}, suppose that $M$ is the midpoint of $\Gamma_R$. Let $x_1=(\Phi_1,\theta_1)\in\MRin$ with $\theta_1\in[\sin^{-1}(2r/R),\pi-\sin^{-1}(2r/R)]$ in case (a) of \eqref{eqMainCases} and $P=p(x_1)\in\Gamma_R$ of the billiard boundary. As shown in \cref{fig:fractionEstimate2}, suppose that $\overrightarrow{TTP}$ is the counterclockwise tangential direction on $P$ of $\Gamma_R$. We have $TP-P-O_R$ as the coordinate system of the billiard table. Then we define several length and angle variables marked in \cref{fig:fractionEstimate2}.
    \begin{itemize}
        \item $\rho\dfn|O_rP|\overset{\text{Cosine Law}}=\sqrt{R^2+b^2-2bR\cos\Phi_1}$
        \item $\alpha\dfn\tan^{-1}{\big(\frac{R\sin\Phi_1}{R\cos\Phi_1-b}\big)}$
        \item $\delta\dfn\alpha-\Phi_1$
    \end{itemize}
Hence\COMMENT{?} $\rho$, $\alpha$, $\delta$ are smooth functions of $x_1$. Note that 
$P=(R\sin\Phi_1,b-R\cos\Phi_1)$, $O_r=(0,0)$, $O_R=(0,b)$, therefore, in \cref{fig:fractionEstimate2}, $\alpha$ is the angle between $\overrightarrow{O_rM}$ and $\overrightarrow{O_rP}$ given the counterclockwise orientation on the billiard table. Then $\alpha$, $\delta$ and $\Phi_1$ have the same sign, and $|\alpha|=\measuredangle{MO_rP}\ge\measuredangle{MO_RP}=|\Phi_1|$,  $|\delta|=\measuredangle{O_rPO_R}$.\COMMENT{\color{red}If it can be done with reasonable effort, then it might be good to rotate \cref{fig:fractionEstimate2} by 90 degrees to save space. Same with \cref{fig:heuristicproof,fig:mirrorR2,fig:Fraction_tau1_d0_min,fig:n1_equal3_leftendpoint_case} and maybe others. Additional clipping may help further. Maybe rotating is easy: \url{https://tex.stackexchange.com/questions/438466/how-may-i-rotate-a-tikzpicture} might help.\newline\wentao{Many figures are rotated and several figures are combined into one.}}
\end{definition}



\begin{lemma}\label{lemma:newcoordinateslength}
    For the length functions $\tau_0$, $\tau_1$, $d_0$, $d_1$, $d_2$ from \cref{def:fpclf,def:AformularForLowerBoundOfExpansion} and $\rho$, $\alpha$, $\delta$ from \cref{def:TPPORsystem} as functions of $x_1=(\Phi_1,\theta_1)\in\Mrin$ with $\theta_1\in[\sin^{-1}(2r/R),\pi-\sin^{-1}(2r/R)]$, we have the following.
    \begin{equation}\label{eq:NewCoordinateLegnth}
        \begin{aligned}
        \tau_1&=\rho\sin{(\theta_1-\delta)}+\sqrt{r^2-\rho^2\cos^2{(\theta_1-\delta)}}\\
        d_2&=\sqrt{r^2-\rho^2\cos^2{(\theta_1-\delta)}}\\
        \tau_0&=\rho\sin{(\theta_1+\delta)}+\sqrt{r^2-\rho^2\cos^2{(\theta_1+\delta)}}\\
        d_0&=\sqrt{r^2-\rho^2\cos^2{(\theta_1+\delta)}}
        \end{aligned}
    \end{equation}
And these formulas imply that $\tau_0,\tau_1,d_0,d_1,d_2$ are continuous functions of $x_1$ in $\mathfrak{D}^\circ\dfn R_1\cup R_2\cup R_3\cup R_4$ in \cref{fig:R1}.
\end{lemma}
\begin{proof}
    We first assume $\theta_1\in[\sin^{-1}(2r/R),\pi/2]$. Using the \cref{def:TPPORsystem} $TP-P-O_R$ coordinate system and its defined variables, we have $P=(0,0)$, $O_r=(x_r,y_r)=(-\rho\sin\delta,\rho\cos\delta)$ and suppose $Q=p(x_2)=(x_Q,y_Q)$ with $x_Q>0$, $y_Q>0$. In the \cref{fig:fractionEstimate2} $TP-P-O_R$ coordinate system, $Q_0=(x_{Q_0},y_{Q_0})$ with $x_{Q_0}<0$, $y_{Q_0}<0$ is another intersection point of the line $PQ$ with $C_r$.
    
    We can observe that $(x_Q,y_Q)$, $(x_{Q_0},y_{Q_0})$ satisfy the following equations. That is, $x=x_Q$, $y=y_Q$ or $x=x_{Q_0}$, $y=y_{Q_0}$ together with $x_r,y_r$ will achieve the equality of the following equations.\COMMENT{Which following equations? The first three? The display is a jumble that could be improved.}
    \begin{align*}
         &\left.
        \begin{aligned}
            &(x-x_r)^2+(y-y_r)^2=r^2\\
            &y=(\tan{\theta_1})x \\
            &x^2_r+y^2_r=\rho^2 \\
        \end{aligned}
        \right\} \Longrightarrow (\sec^2\theta_1)x^2-2(x_r+y_r\tan\theta_1)x+\rho^2-r^2=0\\
        &(x_r,y_r)=(-\rho\sin{\delta},\rho\cos{\delta}) \Longrightarrow 2x_r+2y_r\tan\theta_1=2(-\rho\sin\delta+\rho\cos\delta\tan\theta_1)=2\rho\sec\theta_1\sin{(\theta_1-\delta)}
    \end{align*}
    Hence, $x=x_Q$ and $x=x_{Q_0}$ both satisfy
    \begin{equation}\label{eq:60}
        \begin{aligned}
            &(\sec^2\theta_1)x^2-2\rho\sin{(\theta_1-\delta)}(\sec\theta_1)x+\rho^2-r^2=0\\
        \end{aligned}
    \end{equation}
    If we write $L=x\sec\theta_1$, then \eqref{eq:60} becomes the quadratic equation $L^2-2\rho\sin{(\theta_1-\delta)}L+\rho^2-r^2=0$. Since $\rho^2-r^2<0$, $\sec\theta_1>0$, $x_Q>0$, $x_{Q_0}<0$, this quadratic equation has a positive root $\rho\sin{(\theta_1-\delta)}+\sqrt{r^2-\rho^2\cos^2{(\theta_1-\delta)}}=L_+=x_Q\sec\theta_1$ and a negative root $\rho\sin{(\theta_1-\delta)}-\sqrt{r^2-\rho^2\cos^2{(\theta_1-\delta)}}=L_-=x_{Q_0}\sec\theta_1$. Note that in \cref{fig:fractionEstimate2}, $x_Q\sec\theta_1=|PQ|=\tau_1$, $-x_{Q_0}\sec\theta_1=|PQ_0|$ and $x_Q\sec\theta_1-x_{Q_0}\sec\theta_1=|PQ|+|PQ_0|=|QQ_0|=2d_2$ (see \cref{def:fpclf}). Thus, $\tau_1=\rho\sin{(\theta_1-\delta)}+\sqrt{r^2-\rho^2\cos^2{(\theta_1-\delta)}}$, $d_2=\sqrt{r^2-\rho^2\cos^2{(\theta_1-\delta)}}$.

    Similarly in the coordinate system $TP-P-OR$, for $(x_T,y_T)=T=p(x_0)=p(\mathcal{F}^{-1}(x_1))$ on the billiard table \cref{fig:fractionEstimate2} and $T_0$ to be another intersection of line $PT$ with $C_r$, we have $x_T<0$, $y_T>0$, $x_{T_0}>0$, $y_{T_0}<0$. We can observe that $(x_T,y_T)$ and $(x_{T_0},y_{T_0})$ satisfy the following equations. That is, $x=x_T$, $y=y_T$ or $x=x_{T_0}$, $y=y_{T_0}$ together with $x_r,y_r$ will achieve the equality of the following equations.
    \begin{align*}
         &\left.
        \begin{aligned}
            &(x-x_r)^2+(y-y_r)^2=r^2\\
            &y=-(\tan{\theta_1})x \\
            &x^2_r+y^2_r=\rho^2 \\
        \end{aligned}
        \right\} \Longrightarrow (\sec^2\theta_1)x^2-2(x_r-y_r\tan\theta_1)x+\rho^2-r^2=0\\
        &(x_r,y_r)=(-\rho\sin{\delta},\rho\cos{\delta}) \Longrightarrow 2x_r-2y_r\tan\theta_1=2(-\rho\sin\delta-\rho\cos\delta\tan\theta_1)=-2\rho\sec\theta_1\sin{(\theta_1+\delta)}
    \end{align*}
Hence, $x=x_T$ and $x=x_{T_0}$ both satisfy
    \begin{equation}\label{eq:61}
        \begin{aligned}
            &(\sec^2\theta_1)x^2+2\rho\sin{(\theta_1+\delta)}(\sec\theta_1)x+\rho^2-r^2=0\\
        \end{aligned}
    \end{equation}
    If we write $L=-x\sec\theta_1$, then \eqref{eq:61} becomes a quadratic equation $L^2-2\rho\sin{(\theta_1+\delta)}L+\rho^2-r^2=0$. Since $\rho^2-r^2<0$, $\sec\theta_1>0$, $x_T<0$, $x_{T_0}>0$, this quadratic equation has a positive root $\rho\sin{(\theta_1+\delta)}+\sqrt{r^2-\rho^2\cos^2{(\theta_1+\delta)}}=L_+=-x_T\sec\theta_1$ and a negative root $\rho\sin{(\theta_1+\delta)}-\sqrt{r^2-\rho^2\cos^2{(\theta_1+\delta)}}=L_-=-x_{T_0}\sec\theta_1$. Note that in \cref{fig:fractionEstimate2}, $-x_T\sec\theta_1=|PQ|=\tau_0$, $x_{T_0}\sec\theta_1=|PT_0|$ and $-x_T\sec\theta_1+x_{T_0}\sec\theta_1=|PT|+|PT_0|=|TT_0|=2d_0$ (see \cref{def:fpclf}). Thus, $\tau_0=\rho\sin{(\theta_1+\delta)}+\sqrt{r^2-\rho^2\cos^2{(\theta_1+\delta)}}$, $d_0=\sqrt{r^2-\rho^2\cos^2{(\theta_1+\delta)}}$.
    
    Hence, we proved \eqref{eq:NewCoordinateLegnth} for $\theta_1\in[\sin^{-1}(2r/R),\pi/2]$. 
    
    Then for $\theta_1\in[\pi/2,\pi-\sin^{-1}(2r/R)]$, by symmetry $\tau_1$ with $\theta_1\in[\pi-\sin^{-1}(2r/R),\pi/2]$ is $\tau_0$ with $(\pi-\theta_1)\in[\sin^{-1}(2r/R),\pi/2]]$. $\tau_0$ with $\theta_1\in[\pi-\sin^{-1}(2r/R),\pi/2]$ is $\tau_1$ with $(\pi-\theta_1)\in[\sin^{-1}(2r/R),\pi/2]]$. $d_0$ with $\theta_1\in[\pi-\sin^{-1}(2r/R),\pi/2]$ is $d_2$ with $(\pi-\theta_1)\in[\sin^{-1}(2r/R),\pi/2]$. $d_2$ with $\theta_1\in[\pi-\sin^{-1}(2r/R),\pi/2]$ is $d_0$ with $(\pi-\theta_1)\in[\sin^{-1}(2r/R),\pi/2]$.
    That is, in \eqref{eq:NewCoordinateLegnth} replacing $\theta_1$ with $\pi-\theta_1$ and interchanging the subindex between $\tau_0$ and $\tau_1$ and between $d_0$ and $d_2$. This will give the same equations as in \eqref{eq:NewCoordinateLegnth}.
\end{proof}
\begin{proposition}\label{prop:casea0largetheta}
    Suppose $\phistar\in\big(0,\tan^{-1}(1/3)\big)$,  \eqref{eqJZHypCond} holds,  the orbit segment \eqref{eqPtsFromeqMhatOrbSeg} of $x$ from \cref{def:MhatReturnOrbitSegment,def:fpclf} has $x_0\in \Mrout$, $(\Phi_1,\theta_1)=\mathcal{F}(x_0)=x_1\in \MRin, \mathcal{F}(x_1)=x_2\in \Mrin$. If $\theta_1\in\big(\frac{\pi}{6},\frac{5\pi}{6}\big)$, then $\mathcal{I}$ from \cref{def:AformularForLowerBoundOfExpansion} satisfies $\mathcal{I}<-1.37$.
\end{proposition}
\begin{proof}
 In \cref{fig:fractionEstimate2}, $P=p(x_1)$, $Q=p(x_2)$, $T=p(x_0)$ in $\Gamma_R$. By symmetry, we assume $\theta_1\in (\frac{\pi}{6},\frac{\pi}{2}]$. $\Phi_1\in(-\Phistar,\Phistar)$. Then in \cref{fig:fractionEstimate2}, it has $\measuredangle{P O_R Or}=|\Phi_1|$, $\measuredangle{P O_r Y}=|\alpha|$, $\measuredangle{O_R P O_r}=|\delta|$. 
 Using the coordinates and notation of \cref{def:TPPORsystem}, we have $\rho=|O_rP|$ and the range of $\delta\in(-\phistar+\Phistar,+\phistar-\Phistar)\subset (-\phistar,+\phistar)$.

From \cref{def:AformularForLowerBoundOfExpansion} we note that 
\[
\mathcal{I} = -1+\frac{\tau_1}{d_0}\big[\frac{2(\tau_0-d_0)}{d_1}-\frac{\tau_0+\tau_1-2d_0}{\tau_1}\big]=-1-\frac{\tau_0+\tau_1-2d_0}{d_0}+\frac{\tau_1}{d_0}\frac{2(\tau_0-d_0)}{d_1}.
\]
To estimate the two fractions, we begin with \cref{lemma:newcoordinateslength} and some trigonometry:\COMMENT{How about omitting $d_0$ from what follows and \eqref{eq:casea0largethetalengthsfraction}?\newline\wentao{It will save two lines but add more sentence for conclusion.}}
\begin{equation}\label{eq:casea0largetheta_trigonometrynegativedelta}
\begin{aligned}
\frac{\tau_0+\tau_1-2d_0}{d_0}& \overset{\eqref{eq:NewCoordinateLegnth}}= \frac{\rho\sin{(\theta_1-\delta)}+\rho\sin{(\theta_1+\delta)}-\sqrt{r^2-\rho^2\cos^2{(\theta_1+\delta)}}+ \sqrt{r^2-\rho^2\cos^2{(\theta_1-\delta)}}}{d_0}\\
& = \frac{1}{d_0}\big[2\rho\sin\theta_1\cos\delta+\frac{r^2-\rho^2\cos^2{(\theta_1-\delta)}-(r^2-\rho^2\cos^2{(\theta_1+\delta)})}{\sqrt{r^2-\rho^2\cos^2{(\theta_1+\delta)}}+ \sqrt{r^2-\rho^2\cos^2{(\theta_1-\delta)}}}\big]\\
& \overset{\quad\mathllap{\cos^2x=\frac12\cos2x+\frac12}}= \frac{\rho}{d_0}\big[2\sin{\theta_1}\cos{\delta}+\frac{1}{2}\frac{\rho[\cos{(2\theta_1+2\delta)}-\cos{(2\theta_1-2\delta)}]}{\sqrt{r^2-\rho^2\cos^2{(\theta_1+\delta)}}+ \sqrt{r^2-\rho^2\cos^2{(\theta_1-\delta)}}}\big]\\
& \overset{\quad\mathllap{2\sin x\sin y=\cos(x-y)-\cos(x+y)}}=\frac{\rho}{d_0}\big[2\sin{\theta_1}\cos{\delta}-\frac{\rho\sin{(2\delta)}\sin{(2\theta_1)}}{\sqrt{r^2-\rho^2\cos^2{(\theta_1+\delta)}}+ \sqrt{r^2-\rho^2\cos^2{(\theta_1-\delta)}}}\big]\\
&\overset{\quad\mathllap{\sin2x=2\sin x\cos x}}=\frac{2\rho\sin{\theta_1}\cos{\delta}}{d_0}\big[1-\frac{2\rho\cos{\theta_1}\sin{\delta}}{\sqrt{r^2-\rho^2\cos^2{(\theta_1+\delta)}}+ \sqrt{r^2-\rho^2\cos^2{(\theta_1-\delta)}}}\big]
\end{aligned}
\end{equation}
Now, if $\delta\le0$, then $\displaystyle1-\frac{2\rho\cos{\theta_1}\sin{\delta}}{\sqrt{r^2-\rho^2\cos^2{(\theta_1+\delta)}}+ \sqrt{r^2-\rho^2\cos^2{(\theta_1-\delta)}}}\ge1$.
If $\delta>0$, then $0<\delta<\phistar<\frac{\pi}{6}$, $\frac{\phi}{6}<\theta_1<\frac{\pi}{2}$ and trigonometry give $\tan\delta<1/3$ and  $\cot\theta_1<\sqrt{3}$, so
\begin{equation}\label{eq:casea0largetheta_trigonometrypositivedelta}
\begin{aligned}
1-\frac{2\rho\cos{\theta_1}\sin{\delta}}{\sqrt{r^2-\rho^2\cos^2{(\theta_1+\delta)}}+ \sqrt{r^2-\rho^2\cos^2{(\theta_1-\delta)}}}&
\overbracket{>}^{\mathclap{\rho<r}}  1-\frac{2\rho\cos{\theta_1}\sin{\delta}}{\sqrt{\rho^2-\rho^2\cos^2{(\theta_1+\delta)}}+ \sqrt{\rho^2-\rho^2\cos^2{(\theta_1-\delta)}}}\\
&=  1-\frac{2\cos{\theta_1}\sin{\delta}}{\sin{(\theta_1+\delta)}+\sin{(\theta_1-\delta)}}=1-\cot{\theta_1}\tan{\delta}\\
&>  1-\sqrt{3}\cdot\frac{1}{3}=\frac{\sqrt3-1}{\sqrt3}>0.4.
\end{aligned}
\end{equation}
$0<\phistar<\tan^{-1}{(\frac{1}{3})}$, $\theta_1\in(\frac{\pi}{6},\frac{\pi}{2}]$ together with $\rho>r\cos\phistar$ imply that $|\delta|<\phistar<\tan^{-1}{(\frac{1}{3})},\,\cos{\delta}>\cos{\phistar}>\cos\tan^{-1}\frac13=\sqrt8/3>0.94$,\COMMENT{Do we need 0.94? Or 0.4 above?\newline\wentao{no need now. we have 0.375 now in use. Do you agree no harm to mention 0.4,0.94?}} $\sin\theta_1>1/2$. Therefore, for either $\delta>0$ or $\delta\le 0$, we have
\begin{equation}\label{eq:casea0largethetalengthsfraction}
\frac{\tau_0+\tau_1-2d_0}{d_0}\overbracket{>}^{{\textstyle\strut}\mathclap{(\ref{eq:casea0largetheta_trigonometrynegativedelta},\ref{eq:casea0largetheta_trigonometrypositivedelta})}}
\frac{\sqrt3-1}{\sqrt3}\frac{\overbracket{\rho}^{\mathclap{>r\cos\phistar}}}{\underbracket{d_0}_{\mathclap{\le r}}}\overbracket{2\sin{\theta_1}}^{>1}\overbracket{\cos{\delta}}^{>\cos\phistar}>\frac{\sqrt3-1}{\sqrt3}\cos^2\phistar>\frac89\frac{\sqrt3-1}{\sqrt3}>0.375.
\end{equation}
On the other hand, $\theta_1\in(\frac{\pi}{6},\frac{\pi}{2}]$ therefore $\sin{\theta_1}>\frac{1}{2}$\COMMENT{Where is this used? We seem to only need $\sin\theta_1>1/2$, so might as well say that instead.\wentao{agree. said now}}, $-d_0\le\tau_0-d_0\le d_0$,\COMMENT{Are the r's needed here?\wentao{agree removed}} $0<\tau_1\le 2r$, $R>1700r$ and $d_1=R\sin\theta_1>\frac12R>850r$ imply\COMMENT{Should the 4 below be a 2? In that case also edit the last line of the proof and the statement of \cref{prop:casea0largetheta}.\newline\wentao{Done}} 
\begin{equation}\label{eq:casea0largethetalengthsfractionII}
    \begin{aligned}
        \frac{\overbracket{\tau_1}^{\le 2r}}{d_0}\frac{2(\tau_0-d_0)}{\underbracket{d_1}_{\mathclap{=R\sin\theta_1>\frac12R>850r}}}\le\frac{2}{850}\cdot2\big|\frac{\tau_0-d_0}{d_0}\big|<\frac{1}{425}\cdot2=\frac{4}{850}.
    \end{aligned}
\end{equation}
Hence, from \cref{def:AformularForLowerBoundOfExpansion}, \[
\mathcal{I} = -1+\frac{\tau_1}{d_0}\big[\frac{2(\tau_0-d_0)}{d_1}-\frac{\tau_0+\tau_1-2d_0}{\tau_1}\big]=-1-\frac{\tau_0+\tau_1-2d_0}{d_0}+\frac{\tau_1}{d_0}\frac{2(\tau_0-d_0)}{d_1}\overbracket{<}^{\mathrlap{\eqref{eq:casea0largethetalengthsfraction},\eqref{eq:casea0largethetalengthsfractionII}}\quad}
-1-0.375+\frac{4}{850}<-1.37.\qedhere\]
\end{proof}
For later use in \cref{theorem:UniformExpansionA0,corollary:halfquadrantcone}, we note:\begin{corollary}\label{Corollary:CollisionOnMrin0Mrin1}
  Suppose $\phistar\in\big(0,\tan^{-1}(1/3)\big)$ and \eqref{eqJZHypCond} holds. If the return orbit segment  \eqref{eqPtsFromeqMhatOrbSeg} of $x$ from \cref{def:MhatReturnOrbitSegment,def:fpclf} has $x_0\in \Mrout$, and $\mathcal{F}(x_0)=x_1=(\Phi_1,\theta_1)\in \MRin$, $\mathcal{F}(x_1)=x_2=(\phi_2,\theta_2)\in \textcolor{blue}{M^{\text{\upshape{in}}}_{r,0}}\cup \textcolor{blue}{M^{\text{\upshape{in}}}_{r,1}}$, then
  \begin{itemize}
      \item $\phi_2\in(\pi/2,3\pi/2)$, i.e. $p(x_2)$ is on the upper half of circle $C_r$ in \cref{fig:fractionEstimate2}.
      \item $\theta_1\in(\pi/6,5\pi/6)$.
      \item  $\mathcal{I}$ from \cref{def:AformularForLowerBoundOfExpansion} satisfies $\mathcal{I}<-1.37$.
  \end{itemize}
\end{corollary}
\begin{proof}
    In \cref{fig:MrN}, the leftmost point of the closure of \textcolor{blue}{$M^{\text{\upshape in}}_{r,1}$} has $\phi$-coordinate $\frac{2\pi}{3}-\phistar>\frac{\pi}{2}$, and by symmetry the rightmost point of the closure of \textcolor{blue}{$M^{\text{\upshape in}}_{r,1}$} has $\phi$-coordinate $\frac{4\pi}{3}+\phistar<\frac{3\pi}{2}$. The leftmost point of the closure of \textcolor{blue}{$M^{\text{\upshape in}}_{r,0}$} has $\phi$-coordinate $\pi-\phistar>\frac{\pi}{2}$ and the rightmost point of the closure of \textcolor{blue}{$M^{\text{\upshape in}}_{r,0}$} has $\phi$-coordinate $\pi+\phistar<\frac{3\pi}{2}$. Therefore,  $\theta_1\in (\frac{\pi}{6},\frac{5\pi}{6})$ by \cref{lemma:returnorbittheta1Large}. Furthermore,\COMMENT{"Furthermore" is the wrong word; please look up what it means. "Therefore" is the right word here. Please look it up as well. It would be good to search the file for both and make sure that each is always used correctly.} $\mathcal{I}<-1.37$ by \cref{prop:casea0largetheta}.
\end{proof}
\begin{remark}
In the computation \eqref{eq:casea0largetheta_trigonometrypositivedelta} and \cref{lemma:returnorbittheta1Large}, we see that we can relax the condition in \eqref{eqJZHypCond} to be $\phistar<\pi/6$ and still obtain $\mathcal{I}<-1$. However, when $\phistar$ is close to $\pi/6$, in \eqref{eq:casea0largetheta_trigonometrypositivedelta}, $R$ may need to be arbitrarily large to yield $\mathcal{I}<-1$. We avoid computing $R_{\text{\upshape HF}}$ in such a complicated situation.






\end{remark}

\begin{figure}[h]
\begin{center}
\begin{sideways}
\begin{tikzpicture}[xscale=0.7, yscale=0.7]
    \tkzDefPoint(0,0){Or};
    \tkzDefPoint(0,-4.8){Y};
    \pgfmathsetmacro{\rradius}{4.5};
    \pgfmathsetmacro{\bdist}{8.7};
    \tkzDefPoint(0,\bdist){OR};
    \pgfmathsetmacro{\phistardeg}{40};
    \pgfmathsetmacro{\XValueArc}{\rradius*sin(\phistardeg)};
    \pgfmathsetmacro{\YValueArc}{\rradius*cos(\phistardeg)};
    \tkzDefPoint(\XValueArc,-1.0*\YValueArc){B};
    \tkzDefPoint(-1.0*\XValueArc,-1.0*\YValueArc){A};
    \pgfmathsetmacro{\Rradius}{veclen(\XValueArc,\YValueArc+\bdist)};
    \draw[name path=Cr,dashed,ultra thin] (Or) circle (\rradius);
    \tkzDrawArc[very thick](Or,B)(A);
    \tkzDrawArc[name path=CR, very thick](OR,A)(B);
    \pgfmathsetmacro{\PXValue}{\Rradius*sin(7)};
    \pgfmathsetmacro{\PYValue}{\bdist-\Rradius*cos(7)};
    \tkzDefPoint(\PXValue,\PYValue){P};
    \draw[name path=PQ_Far,blue,ultra thin](P)--+(70:7)coordinate(PposX);
    \draw[name path=PT_Far,blue,ultra thin](P)--+(124:9);
    \draw[name path=PT0_Far,blue,dashed,ultra thin](P)--+(-56:7);
    \draw[name path=PQ0_Far,blue,dashed,ultra thin](P)--+(-110:9);
    \draw[red,dashed,->](P)--+(7:4)coordinate(PposX);
    \draw[name path=PPnegX,red,dashed](P)--+(-173:6)coordinate(PnegX);
    \draw[name path=linePOR,red,dashed,->](P)--(OR);
    \draw[name path=linePOr,red,dashed](P)--(Or);
    \tkzDefPoint(0,-5.5){Y};
    \tkzDefPoint(0,\bdist-\Rradius){C};
    \draw[red,dashed](OR)--(Y);
    \tkzMarkAngle[arc=lll,ultra thin,size=1.7](Or,OR,P);
    \tkzMarkAngle[arc=lll,ultra thin,size=0.7](Y,Or,P);
    \tkzMarkAngle[arc=lll,ultra thin,size=0.8](OR,P,Or);
    \path[name intersections={of=Cr and PQ_Far,by={Q}}];
    \path[name intersections={of=Cr and PQ0_Far,by={Q0}}];
    \path[name intersections={of=Cr and PT0_Far,by={T0}}];
    \path[name intersections={of=Cr and PT_Far,by={T}}];
    \tkzMarkAngle[arc=lll,very thin,size=1.8](PposX,P,Q);
    \tkzMarkRightAngle(OR,P,PnegX);
    \node[rotate=-90] at ($(Or)+(-78:1)$)  {\Large $\alpha$};
    \node[rotate=-90] at ($(OR)+(-85:2.3)$)  {\Large \textcolor{blue}{$\Phi_1$}};
    \node[rotate=-90] at ($(P)+(105:1.3)$)  {\Large {$\delta$}};
    \node[rotate=-90] at ($(P)+(45:2.3)$)  {\Large {$\theta_1$}};
    
    \tkzLabelPoint[left,rotate=-90](A){$A$};
    \tkzLabelPoint[left,rotate=-90](B){$B$};
    \tkzLabelPoint[left,rotate=-90](OR){$O_R$};
    \tkzLabelPoint[below,rotate=-90](Or){$O_r$};
    \tkzLabelPoint[above,rotate=-90]([shift={(0,0.2)}]P){$P$};
    \tkzLabelPoint[above,rotate=-90](Q){$Q$};
    \tkzLabelPoint[left,rotate=-90](T){$T$};
    \tkzLabelPoint[below,rotate=-90](Q0){$Q_0$};
    \tkzLabelPoint[above,rotate=-90](T0){$T_0$};
    \tkzLabelPoint[above,rotate=-90](PposX){$TP$};
    \tkzDefLine[perpendicular=through Or](P,PnegX)\tkzGetPoint{M_far};
    \path[name path= lineOrMfar] (Or)--(M_far);
    \path[name path= PPnegX] (P)--(PnegX);
    \path[name path= PPposX] (P)--(PposX);
    \path[name intersections={of=lineOrMfar and PPnegX,by={M}}];
    \path[name intersections={of=Cr and PPnegX,by={A0}}];
    \path[name intersections={of=Cr and PPposX,by={B0}}];
    \coordinate (Or_mirror) at ($(Or)!2!(M)$);
    \tkzDrawArc[name path=mirrorCr, ultra thin, dashed](Or_mirror,A0)(B0);
    \tkzLabelPoint[below,rotate=-90](Or_mirror){$O'_r$};

\path[name intersections={of=mirrorCr and PQ0_Far,by={Q1}}];
    \path[name intersections={of=mirrorCr and PT0_Far,by={T1}}];
    \tkzLabelPoint[above,rotate=-90](Q1){$Q_1$};
    \tkzLabelPoint[above,rotate=-90](T1){$T_1$};
    \tkzLabelPoint[left,rotate=-90](C){$M$};    \tkzDrawPoints(A,Or,OR,P,B,T,Q,Or_mirror,A0,T0,Q0,T1,Q1,C);
    \begin{scope}[ultra thin,decoration={markings, mark=at position 0.55 with {\arrow{>}}}
    ] 
    \draw[blue,postaction={decorate}](P)--(Q);
    \draw[blue,postaction={decorate}](T)--(P);
    \end{scope}
    \tkzLabelPoint[above,rotate=-90]([shift={(0.15*\rradius,1.1*\rradius)}]Or){\Large $\Gamma_r$};
    \tkzLabelPoint[above,rotate=-90]([shift={(0.35*\rradius,-1.05*\rradius)}]Or_mirror){\Large $\Gamma_r'$};
\end{tikzpicture}
\end{sideways}
\captionsetup{type=figure}
\caption{Reflecting the billiard over the tangent line of $\Gamma_R$ at $P\dfn p(x_1)$. 
Here, $|TT_1|=\tau_0+\tau_1$, $|TT_0|=2d_0$, $|T_0T_1|=\tau_0+\tau_1-2d_0>0$, $|QQ_0|=2d_2$, $|Q_1Q_0|=\tau_0+\tau_1$,$|Q_0Q_1|=\tau_0+\tau_1-2d_2>0$. (See also \cite[Figure 7, equation (4.1)]{hoalb}.)}\label{fig:heuristicproof}
\end{center}
\end{figure}
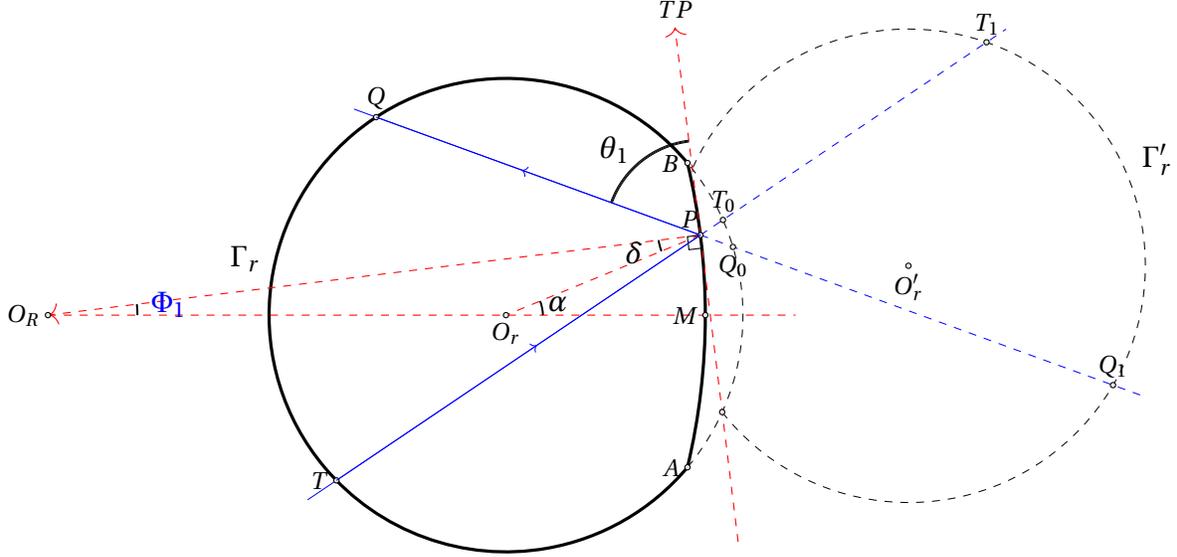

\subsection{Proof of \texorpdfstring{\cref{Prop:Ilowerbound}}{Proposition of Lower bound for I} part II}\label{SBSProp53PartII}
We now consider the situation when the orbit segment \eqref{eqPtsFromeqMhatOrbSeg} of $x$ in \cref{def:MhatReturnOrbitSegment,def:fpclf} has a less transverse reflection on $\Gamma_R$, that is, with $x_1=(\Phi_1,\theta_1)\in R_1\cup R_2\cup R_3\cup R_4\subset\MRin\cap\MRout$ and $\sqrt{4r/R}\le\sin\theta_1
\le\frac{1}{2}$. We know that $\mathcal{I}<-1.37$ if $\sin\theta_1>1/2$ (\cref{prop:casea0largetheta}). 
Now we prove that $\mathcal{I}<-1$ if $x_1=(\Phi_1,\theta_1)\in R_1\cup R_2\cup R_3\cup R_4$ with $\sin\theta_1\ge\sqrt{4r/R}$. By symmetry, we can restrict our analysis to $\theta_1\in [\sin^{-1}(2r/R),\frac{\pi}{2}]$. i.e., the region $x_1\in R_1\cup R_2$  in \cref{fig:R1}. We deal with $x_1\in R_2$ in \cref{prop:casea0smallthetaR2} and with $x_1\in R_1$ in \cref{prop:casea0smallthetaR1}.

\begin{proposition}[$x_1\in R_2$]\label{prop:casea0smallthetaR2}
Suppose $\phistar\in\big(0,\tan^{-1}(1/3)\big)$,  \eqref{eqJZHypCond} holds, and the return orbit segment \eqref{eqPtsFromeqMhatOrbSeg} of $x$ from \cref{def:MhatReturnOrbitSegment,def:fpclf} has $x_0\in\Mrout$, $\mathcal{F}(x_0)=x_1=(\Phi_1,\theta_1)\in\MRin$, $x_2=\mathcal{F}(x_1)\in \Mrin$ with $\theta_1\in\big[\sin^{-1}{(2r/R)},\pi/2\big]$, and $\Phi_1\in(-\Phistar,0)$, that is, $x_1\in R_2$ in \cref{fig:R1} with $\sin\theta_1\ge2r/R$. Then $\mathcal{I}\le-1$  (\cref{def:AformularForLowerBoundOfExpansion}). In the special case where $\sin{\theta_1}\ge\sqrt{4r/R}$, we moreover obtain $\mathcal{I}<-1$.
\end{proposition}
\begin{proof}
For $x_1\in R_2$ in \cref{fig:R1}, thus, with $d_1=R\sin\theta_1\ge2r$ and $\tau_1>0$ in \cref{def:AformularForLowerBoundOfExpansion}, we can analyze $\mathcal{I}=-1+\frac{\tau_1}{d_0}\Big[\frac{2(\tau_0-d_0)}{d_1}-\frac{\tau_0+\tau_1-2d_0}{\tau_1}\Big]$ in the following three cases.
\newline
\textbf{Case i} $\tau_0\le d_0$, then since $\tau_0+\tau_1-2d_0>0$ (see \cref{remark:deffpclf}), we have $\frac{2(\tau_0-d_0)}{d_1}-\frac{\tau_0+\tau_1-2d_0}{\tau_1}<0$. Thus $\mathcal{I}<-1$ .
\newline
\textbf{Case ii} $\tau_0>d_0$ and $\tau_1\ge d_0$, then we have
\begin{equation}\label{eq:casea0smallthetaR2case2}
\begin{aligned}
& \frac{\tau_0+\tau_1-2d_0}{\tau_1}=1-\frac{2d_0-\tau_0}{\tau_1}\overbracket{\ge}^{2d_0-\tau_0>0,\tau_1\ge d_0} 1-\frac{2d_0-\tau_0}{d_0}=\frac{2(\tau_0-d_0)}{2d_0}\overbracket{\ge}^{2d_0 \le 2r\le d_1,\tau_0>d_0}\frac{2(\tau_0-d_0)}{d_1}
\end{aligned}
\end{equation}
Thus, \eqref{eq:casea0smallthetaR2case2} implies $\mathcal{I}\le-1$. In particular, if $\sin\theta_1\ge\sqrt{4r/R}>2r/R$, then $d_1=R\sin\theta_1>2r\ge2d_0$ and the last inequality in \eqref{eq:casea0smallthetaR2case2} becomes $\frac{2(\tau_0-d_0)}{2d_0}\overbracket{>}^{2d_0 \le 2r< d_1,\tau_0>d_0}\frac{2(\tau_0-d_0)}{d_1}$. Thus, if $\sin\theta_1\ge\sqrt{4r/R}$, then $\mathcal{I}<-1$.
\newline
\textbf{Case iii} $\tau_0>d_0>\tau_1$. We show that this cannot occur if $x_1\in R_2=\{(\phi_1,\theta_1)|\theta_1\in[\sin^{-1}{(2r/R)},\pi/2],\Phi_1\in(-\Phistar,0)\}$.

In \cref{fig:mirrorR2}, let $P=p(x_1)\in \Gamma_R$ be on the table boundary. $Q=p(x_2)=p(\mathcal{F}(x_1))$, $T=p(x_0)=p(\mathcal{F}^{-1}(x_1))$. In the coordinate system of \cref{def:StandardCoordinateTable}, $P$ has a negative $x$ coordinate. And since $O_R=(0,b)$, $O_r=(0,0)$, suppose $C_1$ and $C_2$ are the two intersections  of the line $O_RP$ with the circle $C_r$. We have $|C_1C_2|<2r$ and the arc $\arc{C_1TC_2}$ with a length less than half the perimeter of $C_r$. $Q$ and $T$ are on two sides of line $O_RP$. The arc $\arc{C_1QC_2}$ is having a length more than the half of the perimeter of $C_r$. 

Then the arc $\arc{C_1TC_2}$'s mirror symmetry of the line $O_RP$ is contained inside the $C_r$ surrounding disk. $T_1$ which is the mirror symmetry $T$ of the line $O_RP$ must be inside the $C_r$ surrounding disk \big(see \cref{fig:mirrorR2} dashed arc $\arc{C_1T_1C_2}$\big). 

Hence, in \cref{fig:mirrorR2}, $\tau_1=|PQ|\ge|PT_1|=|TP|=\tau_0$,\COMMENT{Fix caption of \Cref{fig:mirrorR2} and of \Cref{fig:R1_fraction_extreme}. Note that often \textbackslash frac is NOT the preferred solution. (But it must be kept if it would otherwise be unclear what is in the numerator and what is in the denominator.) Look at my edits in \Cref{fig:R1_fraction_extreme} and apply that elsewhere (this would be an improvement in MANY places, possibly the majority of the instances of using frac, definitely in \Cref{fig:A_compactRegion,fig:A1_region} and also SELECTIVELY in \Cref{fig:A1contractionCase}). We can talk about this.\wentao{In most cases it's fixed. We may talk more about this.}} so \textbf{Case iii} $\tau_0>d_0>\tau_1$ cannot occur.
\end{proof}
\begin{figure}[h]
\begin{minipage}{.5\textwidth}
\begin{center}
\begin{tikzpicture}[xscale=0.5,yscale=0.5]
    \tkzDefPoint(0,0){Or};
    \tkzDefPoint(0,-4.8){Y};
    \pgfmathsetmacro{\rradius}{4.5};
    \pgfmathsetmacro{\bdist}{8.7};
    \tkzDefPoint(0,\bdist){OR};
    \pgfmathsetmacro{\phistardeg}{40};
    \pgfmathsetmacro{\XValueArc}{\rradius*sin(\phistardeg)};
    \pgfmathsetmacro{\YValueArc}{\rradius*cos(\phistardeg)};
    \tkzDefPoint(\XValueArc,-1.0*\YValueArc){B};
    \tkzDefPoint(-1.0*\XValueArc,-1.0*\YValueArc){A};
    \pgfmathsetmacro{\Rradius}{veclen(\XValueArc,\YValueArc+\bdist)};
    \draw[name path=Cr,ultra thin,dashed] (Or) circle (\rradius);
    \tkzDrawArc[name path=CR,very thick](OR,A)(B);
    \tkzDrawArc[very thick](Or,B)(A);
    \pgfmathsetmacro{\PXValue}{-\Rradius*sin(7)};
    \pgfmathsetmacro{\PYValue}{\bdist-\Rradius*cos(7)};
    \tkzDefPoint(\PXValue,\PYValue){P};
    \path[](OR) -- (P) -- ([turn]90:6)coordinate (PosX);
    \path[](OR) -- (P) -- ([turn]-90:3)coordinate (NegX);
    \path[](OR) -- (P) -- ([turn]0:1)coordinate (NegY);
    \path[](OR) -- (P) -- ([turn]155:9)coordinate (Q_far);
    \path[](OR) -- (P) -- ([turn]-155:9)coordinate (T_far);
    \draw[red,ultra thin,dashed,->](P)--(PosX);
    \draw[red,ultra thin,dashed](P)--(NegX);
    \draw[red,ultra thin,dashed](OR)--(NegY);
    \path[name path=Xpos](P)--(PosX);
    \path[name path=Xneg](P)--(NegX);
    \path[name path=Yneg](P)--(NegY);
    \path[name path=linePQ_far](P)--(Q_far);
    \path[name path=linePT_far](P)--(T_far);
    \path[name path=lineOR_P,ultra thin, dashed,->](P)--(OR);
    \path[name intersections={of=lineOR_P and Cr,by={D}}];
    \path[name intersections={of=linePQ_far and Cr,by={Q}}];
    \path[name intersections={of=linePT_far and Cr,by={T}}];
    \path[name intersections={of=Yneg and Cr,by={C}}];
    \begin{scope}[decoration={
    markings,
    mark=at position 0.5 with {\arrow{>}}}
    ] 
    \draw[blue,postaction={decorate}](T)--(P);
    \draw[blue,postaction={decorate}](P)--(Q);
    \end{scope}
    \coordinate (M) at ($(C)!0.5!(D)$);
    \coordinate (Or_mirror) at ($(Or)!2!(M)$);  
    \tkzDrawArc[name path= Crmirror,red,dashed](Or_mirror,C)(D);
    \tkzDefPointBy[projection=onto P--D](T)\tkzGetPoint{H};
    \coordinate (T1) at ($(T)!2!(H)$);
    \tkzDrawPoints(Or,OR,P,Q,T,T1,C,D);
    \tkzLabelPoint[above](OR){$O_R$};
    \tkzLabelPoint[above](Or){$O_r$};
    \tkzLabelPoint[above left](P){$P$};
    \tkzLabelPoint[right](Q){$Q$};
    \tkzLabelPoint[left](T){$T$};
    \tkzLabelPoint[left](T1){$T_1$};
    \tkzLabelPoint[below left](C){$C_1$};
    \tkzLabelPoint[above left](D){$C_2$};
\end{tikzpicture}
    \caption{Comparing $\tau_0,\tau_1$ for $x_1=(\Phi_1,\theta_1)\in R_2$ with $\Phi_1\in(-\Phistar,0)$, $\theta_1\in [\sin^{-1}{(2r/R)},\pi/2]$}\label{fig:mirrorR2}
\end{center}
\end{minipage}\hfill
\begin{minipage}{.5\textwidth}
\begin{center}
    \begin{tikzpicture}[xscale=0.55,yscale=.55]
        \pgfmathsetmacro{\PHISTAR}{1.5};
        \pgfmathsetmacro{\Dlt}{0.4};
        \pgfmathsetmacro{\Dltplus}{\Dlt+0.6}
        \pgfmathsetmacro{\WIDTH}{12};
        \tkzDefPoint(-\PHISTAR,0){LL};
        \tkzDefPoint(-\PHISTAR,\WIDTH){LU};
        \tkzDefPoint(\PHISTAR,0){RL};
        \tkzDefPoint(\PHISTAR,\WIDTH){RU};
        \tkzDefPoint(-\PHISTAR,\PHISTAR){X};
        \tkzDefPoint(\PHISTAR,\PHISTAR){IY};
        \tkzDefPoint(-\PHISTAR,\PHISTAR+\Dlt){LLD};
        \tkzDefPoint(\PHISTAR,\PHISTAR+\Dlt){RLD};
        \tkzDefPoint(0,\PHISTAR+\Dltplus){LLDp};
        \tkzDefPoint(\PHISTAR,\PHISTAR+\Dltplus){RLDp};
        \tkzDefPoint (0.3*\PHISTAR,\PHISTAR+\Dltplus){montonestart};
        \tkzDefPoint (0.3*\PHISTAR,0.333*\WIDTH){montoneend};
        \draw [thin] (LL) --(LU);
        \draw [thin] (LU) --(RU);
        \draw [thin] (RL) --(RU);
        \draw [thin] (LL) --(RL);
        \draw [blue,ultra thin,dashed] (LL) --(IY);
        \draw [red,ultra thin,dashed]  (X)--(RL);
        \draw [purple,ultra thin,dashed]  (LLD)--(RLD);
        \draw [green,ultra thin]  (LLDp)--(RLDp);
        \draw [green,ultra thin]  (montonestart)--(montoneend);
        \tkzLabelPoint[right](RU){$\pi/2$};
        \tkzLabelPoint[right](RLD){$\sin{\theta_1}=2r/R$};
        \tkzLabelPoint[right](RLDp){$\sin{\theta_1}=\sqrt{4r/R}$};
        \coordinate (RM) at ($(RL)!0.333!(RU)$);
        \coordinate (UpperM) at ($(LU)!0.5!(RU)$);
        \coordinate (LowerM) at ($(LL)!0.5!(RL)$);
        \coordinate (LM) at ($(LowerM)!0.333!(UpperM)$);
        \coordinate (x1) at ($(montonestart)!0.5!(montoneend)$);
        \tkzLabelPoint[right](RM){$\pi/6$};
        \draw [purple,ultra thin]  (LM)--(RM);
        \draw [dashed,ultra thin]  (LowerM)--(UpperM);
        \node[] at (0.5*\PHISTAR, 0.8*\PHISTAR+0.8*\Dlt+0.6*\WIDTH){\Large $R_1$};
        \node[] at (-0.5*\PHISTAR, 0.8*\PHISTAR+0.8*\Dlt+0.6*\WIDTH){\Large $R_2$};
        \tkzLabelPoint[below](LL){$-\Phistar$};
        \tkzLabelPoint[below](RL){$+\Phistar$};
        \tkzLabelPoint[below](LowerM){$0$};
        \tkzLabelPoint[left](LLDp){\tiny $x_M$};
        \tkzLabelPoint[below](montonestart){\tiny $x_L$};
        \tkzLabelPoint[above](montoneend){\tiny $x_U$};
        \tkzLabelPoint[right](x1){\tiny $x_1$};
        \tkzDrawPoints(LLDp,montonestart,montoneend,x1);
    \end{tikzpicture}
    \caption{For $x_1\in R_1$ with \newline $\frac{1}{2}\ge\sin{(\theta_1)}\ge\sqrt{4r/R}$, $\frac{2d_0-\tau_0}{\tau_1}$ \newline takes maximum value at $x_M$}\label{fig:R1_fraction_extreme}
\end{center}
\end{minipage}
\end{figure}

\begin{proposition}[$x_1\in R_1$]\label{prop:casea0smallthetaR1}
Suppose $\phistar\in\big(0,\tan^{-1}(1/3)\big)$, \eqref{eqJZHypCond} holds, and the return orbit segment \eqref{eqPtsFromeqMhatOrbSeg} of $x$ from \cref{def:MhatReturnOrbitSegment,def:fpclf} has $x_0\in\Mrout$, $\mathcal{F}(x_0)=x_1=(\Phi_1,\theta_1)\in\MRin$ and $x_2=\mathcal{F}(x_1)\in \Mrin$ with $\theta_1\in[\sin^{-1}{(\sqrt{4r/R})},\pi/6],\Phi_1\in[0,\Phistar)$, that is, $x_1\in R_1$ in \cref{fig:R1_fraction_extreme} with $1/2\ge\sin\theta_1\ge\sqrt{4r/R}$. Then $\mathcal{I}<-1$ in \cref{def:AformularForLowerBoundOfExpansion}.
\end{proposition}


\begin{proof}
$x_1=(\Phi_1,\theta_1)\in R_1$ with $\theta_1\in[\sin^{-1}{(\sqrt{4r/R})},\frac{\pi}{2}]$ and $\Phi_1\in[0,\Phistar)$ is on a vertical line segment inside $R_1$ connected to its upper and lower end point $x_U=(\Phi_1,\pi/6)$ and $x_L=(\Phi_1,\sin^{-1}(\sqrt{4r/R}))$ as shown in \cref{fig:R1_fraction_extreme}. And $x_M=(0,\sin^{-1}(\sqrt{4r/R}))$ and $x_L$ are on the same green horizontal line in \cref{fig:R1_fraction_extreme}.

We consider $\tau_1(x_1)$, $\tau_0(x_1)$, $d_2(x_1)$, $d_0(x_1)$, $d_1(x_1)$ in \cref{def:fpclf} as functions of $x_1=(\Phi_1,\theta_1)$ on $R_1$. 

For $x_1$ being on the line segment $[x_L,x_U]$ (vertical \textcolor{green}{green} line in \cref{fig:R1_fraction_extreme}), in the context of \cref{monotone_p1,monotone_p2} since $\pi/6<\pi/2-\phistar$, \cref{monotone_p2,rmk:monotone_p2} imply that $(2d_0-\tau_0)\bigm|_{x_1=x_L}\ge (2d_0-\tau_0)(x_1)>0$ and \cref{monotone_p1} implies that $0<\tau_1\bigm|_{x_1=x_L}\le \tau_1(x_1)$. Therefore, $\frac{2d_0-\tau_0}{\tau_1}(x_1)\le \frac{2d_0-\tau_0}{\tau_1}\bigm|_{x_1=x_L}$.

On the other hand, for $x_1$ being on the line segment $[x_M,x_L]$ (horizontal \textcolor{green}{green} line in \cref{fig:R1_fraction_extreme}), \cref{monotone_p5} implies that $\frac{2d_0-\tau_0}{\tau_1}\bigm|_{x_1=x_L}\le \frac{2d_0-\tau_0}{\tau_1}\bigm|_{x_1=x_M}$.

Hence, for all $x_1\in R_1$ with $\sqrt{4r/R}\le\sin\theta_1\le1/2$, we have $\frac{2d_0-\tau_0}{\tau_1}(x_1)\le\frac{2d_0-\tau_0}{\tau_1}\bigm|_{x_1=x_M}$. That is $\frac{\tau_1+\tau_0-2d_0}{\tau_1}(x_1)=1-\frac{2d_0-\tau_0}{\tau_1}(x_1)\ge 1-\frac{2d_0-\tau_0}{\tau_1}\bigm|_{x_1=x_M}=\frac{\tau_1+\tau_0-2d_0}{\tau_1}\bigm|_{x_1=x_M}$. Then we will go to compute this lower bound $\frac{\tau_1+\tau_0-2d_0}{\tau_1}\bigm|_{x_1=x_M}$.

Note that if $(\Phi_M,\theta_M)=x_1=x_M=(0,\sin^{-1}(\sqrt{4r/R}))$, that is $p(x_1)=P=M$ the midpoint of the arc $\Gamma_R$ in \cref{fig:fractionEstimate2} with $\delta=0$, then in \eqref{eq:NewCoordinateLegnth} we have $\rho=|O_rP|=|O_rM|=r\cos\phistar+2R\sin^2{(\frac{\Phistar}{2})}>r\cos\phistar\overset{\phistar<\pi/3}>\frac{r}{2}$. We further get

\begin{equation}\label{eq:new73}
\big(\tau_0+\tau_1-2d_0\big)\bigm|_{x_1=x_M} \overset{\eqref{eq:NewCoordinateLegnth}\text{ with }\delta=0}{\scalebox{6}[1]{$=$}}2\rho\sin{\theta_M}.
\end{equation}
Therefore,
\begin{equation}\label{eq:new74}
    \begin{aligned}
        \big(\frac{\tau_1+\tau_0-2d_0}{\tau_1}\big)(x_1)&\ge\big(\frac{\tau_0+\tau_1-2d_0}{\tau_1}\big)\bigm|_{x_1=x_M}\overset{\eqref{eq:new73}}=\frac{2\rho\sin{\theta_M}}{\tau_1\bigm|_{x_1=x_M}}\overset{\rho>r/2}>r\frac{\sin{\theta_M}}{\tau_1\bigm|_{x_1=x_M}}\\
        &\overbracket{\ge}^{\mathllap{\tau_1\le2r}}\frac{\sin{\theta_M}}{2}\overset{\sin\theta_M=\sqrt{4r/R}}{\scalebox{8}[1]{$=$}}\frac{2r}{R\sin{\theta_M}}\overset{\sin\theta_1\ge\sin{\theta_M}}\ge\frac{2r}{d_1(x_1)}\overset{r\ge d_0}\ge\frac{2d_0}{d_1}(x_1)\overset{2d_0>\tau_0}>\frac{2(\tau_0-d_0)}{d_1}(x_1)
    \end{aligned}
\end{equation}
for all $x_1\in R_1$ with $\sqrt{4r/R}\le\sin\theta_1\le1/2$. Then $\mathcal{I}<-1$ from \cref{def:AformularForLowerBoundOfExpansion}.
\end{proof}

\subsection{Continuous extension of the second derivative to corners}\label{subsec:continuousextension}\strut\COMMENT{The second derivative or the derivative of the second iterate? These are VERY different!}
\begin{remark}[$D$ is negative on $\mathfrak D^\circ$]\label{rmk:extensionDefs}
\label{def:DF2matrix}
    In the orbit segment $x_0\in \Mrout, \mathcal{F}(x_0)=x_1=(\Phi_1,\theta_1)\in \MRin\cap\MRout$, $x_2=\mathcal{F}(x_1)\in\Mrin$ corresponding to cases (a) and (b) in \eqref{eqMainCases},  \cite[Equations (4.4), (4.5)]{hoalb} imply
\begin{equation}\label{eq:DF2matrix}
\begin{aligned}
D=&
\begin{bmatrix}
D_{11} & D_{12}\\
D_{21} & D_{22}
\end{bmatrix}
\dfn d_1 d_2D\mathcal{F}^2_{x_0}\\
=&
\begin{bmatrix}
2\tau_1(\tau_0-d_0)-d_1(\tau_0-d_0+\tau_1) & 2\tau_1\tau_0-d_1(\tau_0+\tau_1) \\
2(\tau_1-d_2)(\tau_0-d_0)-d_1(\tau_0-d_0+\tau_1-d_2) & 2(\tau_1-d_2)\tau_0-d_1(\tau_0+\tau_1-d_2)
\end{bmatrix}
\end{aligned}
\end{equation}with $d_0$, $d_1$, $\tau_0$, $\tau_1$, $\tau_2$ as in \cref{def:fpclf}.
If $x_1=(\Phi_1,\theta_1)\in\mathfrak{D}^\circ= R_1\cup R_2\cup R_3\cup R_4=\big\{(\Phi_1,\theta_1)\in\MRin\bigm|\sin\theta_1\ge2r/R\big\}\subset\MRin\cap \MRout$  (\cref{fig:R1}), then $D$ is a negative matrix by \cref{caseA}. 
We extend $\tau_0,\tau_1,d_0,d_1,d_2$ to $x_1=(\pm\Phistar,\theta_*)$\COMMENT{Is $\theta_*$ ideal notation? (Also: where is it defined?)} on the boundary of $\mathfrak{D}^\circ$ with $p(x_1)=A$ or $B$ at a corner of the table.\COMMENT{This sounds like it has no purpose. Presumably, before this sentence you want to say that the aim is to extend negativity of $D$ to $\mathfrak D\dfn\overline{\mathfrak{D}}^\circ$. Is that the purpose?}
\end{remark}

\begin{definition}[Continuous extensions]\label{def:extensionDefs}
The extensions of $\tau_0,\tau_1,d_0,d_1,d_2$ for $p(x_1)$ to the corners of the billiard table are given by the following (with $\theta_*\in [\sin^{-1}{(2r/R)},\pi-\sin^{-1}{(2r/R)}]$).
\par\smallskip\noindent\begin{minipage}{.495\textwidth}
For $p(x_1)=A$ in Fig.~\ref{fig:heuristicproof}, $x_1=(-\Phistar,\theta_*)\in\partial\MRin$  in Fig.~\ref{fig:R1}, 
\begin{equation}\label{eq:lengthextensionA}
\begin{aligned}
& \tau_0=r\sin{(\theta_*-\phistar+\Phistar)}+r|\sin{(\theta_*-\phistar+\Phistar)}|\text{,}\\
& d_0=r|\sin{(\theta_*-\phistar+\Phistar)}|\text{,}\\
& \tau_1=r\sin{(\theta_*+\phistar-\Phistar)}+r|\sin{(\theta_*+\phistar-\Phistar)}|\text{,}\\
& d_2=r|\sin{(\theta_*+\phistar-\Phistar)}|\text{,}\\
& d_1=R\sin{\theta_*}.
\end{aligned}
\end{equation}    
\end{minipage}
%
\hfill\vrule\hfill\begin{minipage}{.495\textwidth}
For $p(x_1)=B$ in Fig.~\ref{fig:heuristicproof}, $x_1=(\Phistar,\theta_*)\in\partial\MRin$ in Fig.~\ref{fig:R1}, 
\begin{equation}\label{eq:lengthextensionB}
\begin{aligned}
& \tau_0=r\sin{(\theta_*+\phistar-\Phistar)}+r|\sin{(\theta_*+\phistar-\Phistar)}|\text{,}\\
& d_0=r|\sin{(\theta_*+\phistar-\Phistar)}|\text{,}\\
& \tau_1=r\sin{(\theta_*-\phistar+\Phistar)}+r|\sin{(\theta_*-\phistar+\Phistar)}|\text{,}\\
& d_2=r|\sin{(\theta_*-\phistar+\Phistar)}|\text{,}\\
& d_1=R\sin{\theta_*}.
\end{aligned}
\end{equation}
\end{minipage}
\end{definition}


\NEWPAGE\begin{lemma}\label{lemma:extendedDomainDF2}
Suppose $\phistar\in\big(0,\tan^{-1}(1/3)\big)$ and \eqref{eqJZHypCond} holds.
\begin{itemize}
    \item The extensions in \cref{def:extensionDefs} of $\tau_0$, $\tau_1$, $d_0$, $d_1$, $d_2$ from \cref{def:fpclf} are continuous on the extended domain $\mathfrak{D}\dfn\overline{\mathfrak{D}}^\circ=\,\overline{\!R_1\cup R_2\cup R_3 \cup R_4}=\big\{(\Phi_1,\theta_1)\bigm|-\Phistar\le\Phi_1\le\Phistar,\,\sin^{-1}{(2r/R)}\le\theta_1\le\pi-\sin^{-1}{(2r/R)}\big\}$\COMMENT{Correct?} in \cref{fig:A_compactRegion}.\COMMENT{\cref{fig:A_compactRegion} seems to represent $\mathfrak{D}$ incorrectly because it suggests that $\mathfrak{D}$ is the shaded region. The figure should be altered so $\mathfrak{D}$ is next to a vertical brace that spans the right vertical range.\Hrule Also, why is that figure so far from here?}
    \item The matrix elements $D_{11}$, $D_{21}$, $D_{12}$, $D_{22}$ from 
    \eqref{eq:DF2matrix} are negative on $\mathfrak{D}$.\COMMENT{\cref{fig:A_compactRegion} only mentions that two of them are negative.}
\end{itemize}
\end{lemma}
\begin{proof} We first check that \eqref{eq:lengthextensionA} and \eqref{eq:lengthextensionB} indeed provide continuous extensions of $\tau_0,\tau_1,d_0,d_1,d_2$.  \cref{lemma:newcoordinateslength} and \eqref{eq:NewCoordinateLegnth} show that 
    $\tau_0$, $\tau_1$, $d_0$, $d_2$ are continuous functions of $\rho,\alpha,\delta$ from \cref{def:TPPORsystem}. And $\rho,\alpha,\delta$ are also continuous functions of $x_1$ on $\Int{\mathfrak{D}}\subset\MRin\cap\MRout$.
    
    As $\Int{\mathfrak{D}}\ni x_1=(\Phi_1,\theta_1)\rightarrow(-\Phistar,\theta_*)\in\partial\MRin$ with $\theta_1\in [\sin^{-1}{(2r/R)},\pi-\sin^{-1}{(2r/R)}]$, by 
\cref{def:TPPORsystem} $\rho\rightarrow r$,  $\alpha\rightarrow -\phistar$, and $\delta\rightarrow -\phistar+\Phistar$. Then by \eqref{eq:NewCoordinateLegnth} and with $d_1$ from \cref{def:fpclf}, we get
\begin{equation}\label{eq:continuousExtension}
    \begin{aligned}
        \tau_1&=\rho\sin{(\theta_1-\delta)}+\sqrt{r^2-\rho^2\cos^2{(\theta_1-\delta)}}\longrightarrow r\sin{(\theta_*+\phistar-\Phistar)}+r|\sin{(\theta_*+\phistar-\Phistar)}|\text{,}\\
        d_2&=\sqrt{r^2-\rho^2\cos^2{(\theta_1-\delta)}}\longrightarrow r|\sin{(\theta_*+\phistar-\Phistar)}|\text{,}\\
        \tau_0&=\rho\sin{(\theta_1+\delta)}+\sqrt{r^2-\rho^2\cos^2{(\theta_1+\delta)}}\longrightarrow r\sin{(\theta_*-\phistar+\Phistar)}+r|\sin{(\theta_*-\phistar+\Phistar)}|\text{,}\\
        d_0&=\sqrt{r^2-\rho^2\cos^2{(\theta_1+\delta)}}\longrightarrow r|\sin{(\theta_*-\phistar+\Phistar)}|\text{,}\\
        d_1&=R\sin\theta_1\longrightarrow R\sin\theta_*
    \end{aligned}
\end{equation}
By the same reasoning, the limits of $\tau_0$, $\tau_1$, $d_0$, $d_0$, $d_1$ are given in \eqref{eq:lengthextensionB} as $\Int{\mathfrak{D}}\ni x_1=(\Phi_1,\theta_1)\rightarrow(\Phistar,\theta_*)\in\partial\MRin$. Therefore, with the extended definitions in \eqref{eq:lengthextensionA} and \eqref{eq:lengthextensionB} $\tau_0$, $\tau_1$, $d_0$, $d_1$, $d_2$ are continuous on $\mathfrak{D}$.

The remainder of the proof is dedicated to showing that $D$ is (entrywise) negative.
By \cref{caseA} and \cite[Lemma 4.1]{hoalb},\COMMENT{Where is this (as opposed to \cref{caseA}) being used?} we have $D_{11}<0$, $D_{21}<0$ for $x_1\in\Big\{x_1=(\Phi_1,\theta_1)\bigm|-\Phistar<\Phi_1<\Phistar,\,\sin^{-1}{(2r/R)}\le\theta_1\le\pi-\sin^{-1}{(2r/R)}\Big\}$. It suffices to verify $D_{11}<0$, $D_{21}<0$ $D_{12}<0$, $D_{22}<0$ for $x_1=(\pm\Phistar,\theta_*)$ with $\sin\theta_*\ge2r/R$.\COMMENT{Why?} We have the following cases.
\newline
\textbf{Case 1:} 
    $x_1=(-\Phistar,\theta_*)$, $p(x_1)=A$, $d_1=R\sin\theta_*\ge2r$ with $\tau_0,\tau_1,d_0,d_1,d_2$ as in \eqref{eq:lengthextensionA}.
\newline
\textbf{Subcase 1.1:}  $\theta_*\in\big(\phistar-\Phistar,\,\pi-\phistar+\Phistar\big)$.
By \eqref{eq:lengthextensionA} we have $\tau_0=2d_0=2r\sin{(\theta_*-\phistar+\Phistar)}>0$ and $\tau_1=2d_2=2r\sin{(\theta_*+\phistar-\Phistar)}>0$, hence
\begin{align*}\frac{2}{d_1}&\overset{d_1\ge2r}\le\frac{1}{r}\overset{0<\tau_0-d_0=d_0\le r}\le\frac{1}{\tau_0-d_0}\overset{\tau_1>0}<\frac{1}{\tau_0-d_0}+\frac{1}{\tau_1},\\
\frac{2}{d_1}&\overset{d_1\ge2r}\le\frac{1}{r}<\frac{2}{r}\overset{d_0\le r, d_2\le r}\le\frac{1}{d_0}+\frac{1}{d_2}\overset{\tau_0=2d_0,\tau_1=2d_2}{\scalebox{7}[1]{$=$}}\frac{1}{\tau_0-d_0}+\frac{1}{\tau_1-d_2},\\
\frac{2}{d_1}&\overset{d_1\ge2r}\le\frac{1}{r}\overset{0<\tau_1-d_2=d_2\le r}\le\frac{1}{\tau_1-d_2}\overset{\tau_0>0}<\frac{1}{\tau_0}+\frac{1}{\tau_1-d_2},\\
\frac{2}{d_1}&\overset{d_1\ge2r}\le\frac{1}{2r}+\frac{1}{2r}<\frac{1}{\tau_0}+\frac{1}{\tau_1},\text{ since }\big(\sin{(\theta_*-\phistar+\Phistar)},\sin{(\theta_*+\phistar-\Phistar)}\big)\neq(1,1).
\end{align*}
This and $d_0,d_1,d_2,\tau_0,\tau_1>0$ give
\begin{equation}\label{eq:DextensionCase1-1}
\begin{gathered}
\begin{multlined}[.85\textwidth]
D_{11}=2\tau_1(\tau_0-d_0)-d_1(\tau_0-d_0+\tau_1)\\
=(\tau_0-d_0)d_1\tau_1\big(\frac{2}{d_1}-\frac{1}{\tau_1}-\frac{1}{\tau_0-d_0}\big)
\overbracket{=}^{\mathclap{\tau_0=2d_0}}\underbracket{d_0d_1\tau_1}_{>0}\big(\overbracket{\frac{2}{d_1}-\frac{1}{\tau_1}-\frac{1}{\tau_0-d_0}}^{<0}\big)<0,
\end{multlined}\\
\begin{multlined}[.85\textwidth]
D_{21}=2(\tau_1-d_2)(\tau_0-d_0)-d_1(\tau_0-d_0+\tau_1-d_2)\\
=(\tau_1-d_2)(\tau_0-d_0)d_1(\frac{2}{d_1}-\frac{1}{\tau_1-d_2}-\frac{1}{\tau_0-d_0})
\overbracket{=}^{\mathclap{\tau_0=2d_0,\tau_1=2d_2}}\underbracket{d_2d_0d_1}_{>0}\big(\overbracket{\frac{2}{d_1}-\frac{1}{\tau_1-d_2}-\frac{1}{\tau_0-d_0}}^{<0}\big)<0,    
\end{multlined}\\
\begin{multlined}[.85\textwidth]
D_{22}=2(\tau_1-d_2)\tau_0-d_1(\tau_0+\tau_1-d_2)\\
=(\tau_1-d_2)d_1\tau_0(\frac{2}{d_1}-\frac{1}{\tau_0}-\frac{1}{\tau_1-d_2})=\overbracket{d_2d_1\tau_0}^{>0}\big(\overbracket{\frac{2}{d_1}-\frac{1}{\tau_0}-\frac{1}{\tau_1-d_2}}^{<0})<0,\\
\end{multlined}\\
\begin{multlined}[.85\textwidth]
D_{12}=2\tau_1\tau_0-d_1(\tau_0+\tau_1)=\overbracket{\tau_0\tau_1d_1}^{>0}\big(\overbracket{\frac{2}{d_1}-\frac{1}{\tau_0}-\frac{1}{\tau_1}}^{<0}\big)<0.
\end{multlined}
\end{gathered}
\end{equation}
\textbf{Subcase 1.2:} $\theta_*\in \big[\sin^{-1}(2r/R),\,\phistar-\Phistar\big]$. By \eqref{eq:lengthextensionA}, we have $\tau_0=0$,
\[
        d_0=r|\sin{(\theta_*-\phistar+\Phistar)}|=r\sin{(-\theta_*+\phistar-\Phistar)}\ge0
        \text{, and  }\tau_1=2d_2=2r\sin{(\theta_*+\phistar-\Phistar)}>0.
\]
Since $0\le -\theta_*+\phistar-\Phistar< \theta_*+\phistar-\Phistar\le2\phistar-2\Phistar<\pi/2$, we have \[d_2=r\sin{(\theta_*+\phistar-\Phistar)}>r\sin{(-\theta_*+\phistar-\Phistar)}=d_0\ge0\]
and
\begin{equation}\label{eq:DextensionCase1-2}
\begin{aligned}
& D_{11}=2\tau_1(\tau_0-d_0)-d_1(\tau_0-d_0+\tau_1)\overset{\tau_0=0,\tau_1=2d_2}{\scalebox{7}[1]{$=$}}-2\tau_1d_0-d_1(-d_0+2d_2)\overbracket{<}^{\mathclap{d_2>d_0\ge0,\tau_1>0,d_1\ge2r}}0,\\
& D_{21}=2(\tau_1-d_2)(\tau_0-d_0)-d_1(\tau_0-d_0+\tau_1-d_2)\overset{\tau_0=0,\tau_1=2d_2}{\scalebox{7}[1]{$=$}}-2d_2d_0-d_1(d_2-d_0)\overbracket{<}^{\mathclap{d_2>d_0\ge0,d_2>0,d_1\ge2r}}0.\\
& D_{22}=2\tau_0(\tau_1-d_2)-d_1(\tau_1-d_2+\tau_0)\overset{\tau_0=0,\tau_1=2d_2}{\scalebox{7}[1]{$=$}}-d_1d_2\overbracket{<}^{\mathclap{d_2>d_0\ge0,d_1\ge2r}}0,\\
&D_{12}=2\tau_1\tau_0-d_1(\tau_0+\tau_1)\overset{\tau_0=0,\tau_1=2d_2}{\scalebox{7}[1]{$=$}}-2d_1d_2\overbracket{<}^{\mathclap{d_2>d_0\ge0,d_1\ge2r}}0.
\end{aligned}
\end{equation}
\textbf{Subcase 1.3:} $\theta_*\in \big[\pi-\phistar+\Phistar,\,\pi-\sin^{-1}(2r/R)\big]$. By \eqref{eq:lengthextensionA}, we have 
\begin{align*}
    &\tau_1=0\text{ , }d_2=r|\sin{(\theta_*+\phistar-\Phistar-\pi)}| = r\sin{(\theta_*+\phistar-\Phistar-\pi)}\ge0
    \text{ and }\\&\tau_0=2d_0=2r|\sin{(\theta_*-\phistar+\Phistar)}|=2r\sin{(\pi-\theta_*+\phistar-\Phistar)}>0.
\end{align*}
Since $0\le\theta_*+\phistar-\Phistar-\pi<\pi-\theta_*+\phistar-\Phistar\le2(\phistar-\Phistar)<\pi/2$, we have \[d_0=r\sin{(\pi-\theta_*+\phistar-\Phistar)}>r\sin{(\theta_*+\phistar-\Phistar-\pi)}=d_2\ge0\] 
and
\begin{equation}\label{eq:DextensionCase1-3}
\begin{aligned}
&D_{11}=2\tau_1(\tau_0-d_0)-d_1(\tau_0-d_0+\tau_1)\overset{\tau_1=0,\tau_0=2d_0}{\scalebox{5}[1]{$=$}}-d_1d_0\overbracket{<}^{\mathclap{d_0>0,d_1\ge2r}}0,\\
&D_{21}=2(\tau_1-d_2)(\tau_0-d_0)-d_1(\tau_0-d_0+\tau_1-d_2)\overset{\tau_1=0,\tau_0=2d_0}{\scalebox{5}[1]{$=$}}-2d_2d_0-d_1(d_0-d_2)\overbracket{<}^{\mathclap{d_0>d_2\ge0,d_1\ge2r}}0.\\
&D_{22}=2(\tau_1-d_2)\tau_0-d_1(\tau_0+\tau_1-d_2)\overset{\tau_1=0,\tau_0=2d_0}{\scalebox{5}[1]{$=$}}-4d_2d_0-d_1(2d_0-d_2)\overbracket{<}^{\mathclap{d_0>d_2\ge0,d_1\ge2r}}0\\
&D_{12}=2\tau_1\tau_0-d_1(\tau_0+\tau_1)\overset{\tau_1=0,\tau_0=2d_0}{\scalebox{5}[1]{$=$}}-2d_0d_1\overbracket{<}^{\mathclap{d_0>d_2\ge0,d_1\ge2r}}0
\end{aligned}
\end{equation}
\textbf{Case 2} (symmetric to Case 1): $x_1=(\Phistar,\theta_*),p(x_1)=B,d_1=R\sin\theta_*\ge2r$, $\tau_0,\tau_1,d_0,d_1,d_2$ are given in \eqref{eq:lengthextensionB}.
\newline
\textbf{Subcase 2.1:} $\theta_*\in \big(\phistar-\Phistar,\,\pi-\phistar+\Phistar\big)$. By \eqref{eq:lengthextensionB}, we have $\tau_0=2d_0=2r\sin{(\theta_1-\phistar+\Phistar)}>0$ and
$\tau_1=2d_2=2r\sin{(\theta_1+\phistar-\Phistar)}>0$. Using the same reasoning as in \eqref{eq:DextensionCase1-1}, we see that $D_{11}<0$, $D_{21}<0$, $D_{12}<0$, $D_{22}<0$ in this subcase.
\newline
\textbf{Subcase 2.2:} $\theta_*\in [\sin^{-1}(2r/R),\,\phistar-\Phistar]$. By \eqref{eq:lengthextensionB}, we have \[\tau_1=0\text{ , }d_2=r|\sin{(\theta_*-\phistar+\Phistar)}|=r\sin{(-\theta_*+\phistar-\Phistar)}\text{ and }\tau_0=2d_0=2r\sin{(\theta_*+\phistar-\Phistar)}>0\]

Since $0\le -\theta_*+\phistar-\Phistar< \theta_*+\phistar-\Phistar\le2\phistar-2\Phistar<\pi/2$, we have \[d_0=r\sin{(\theta_*+\phistar-\Phistar)}>r\sin{(-\theta_*+\phistar-\Phistar)}\ge0\]

Using the same reasoning as in \eqref{eq:DextensionCase1-3}, we see $D_{11}<0$, $D_{21}<0$, $D_{12}<0$, $D_{22}<0$ in this subcase.
\newline
\textbf{Subcase 2.3:} $\theta_*\in \big[\pi-\phistar+\Phistar,\,\pi-\sin^{-1}(2r/R)\big]$. By \eqref{eq:lengthextensionB}, we have 
\begin{align*}
    &\tau_0=0\text{ , }d_0=r|\sin{(\theta_*+\phistar-\Phistar)}|=r\sin{(\theta_*+\phistar-\Phistar-\pi)}\ge0\\
    \text{ and }&\tau_1=2d_2=2r|\sin{(\theta_*-\phistar+\Phistar)}|=2r\sin{(\pi-\theta_*+\phistar-\Phistar)}>0
\end{align*}
Since $0\le\theta_*+\phistar-\Phistar-\pi<\pi-\theta_*+\phistar-\Phistar\le2(\phistar-\Phistar)<\pi/2$, we have \[d_2=r\sin{(\theta_*+\phistar-\Phistar-\pi)}>d_0=r\sin{(\pi-\theta_*+\phistar-\Phistar)}\ge0\] 
Using the same reasoning as in \eqref{eq:DextensionCase1-2}, we see $D_{11}<0$, $D_{21}<0$, $D_{12}<0$, $D_{22}<0$ in this subcase.

Hence, $D_{11}$, $D_{21}$, $D_{12}$, $D_{22}$ are negative on $\mathfrak{D}$.
\end{proof}

\begin{figure}[h]
\begin{center}
    \begin{tikzpicture}[xscale=.5, yscale=.5]
        \pgfmathsetmacro{\PHISTAR}{1.5};
        \pgfmathsetmacro{\Dlt}{0.4};
        \pgfmathsetmacro{\Dltplus}{\Dlt+0.6}
        \pgfmathsetmacro{\WIDTH}{12};
        \tkzDefPoint(-\PHISTAR,0){LL};
        \tkzDefPoint(-\PHISTAR,\WIDTH){LU};
        \tkzDefPoint(\PHISTAR,0){RL};
        \tkzDefPoint(\PHISTAR,\WIDTH){RU};
        \tkzDefPoint(-\PHISTAR,\PHISTAR){X};
        \tkzDefPoint(\PHISTAR,\WIDTH-\PHISTAR){Y};
        \tkzDefPoint(-\PHISTAR,\WIDTH-\PHISTAR){IX};
        \tkzDefPoint(\PHISTAR,\PHISTAR){IY};
        \tkzDefPoint(-\PHISTAR,\PHISTAR+\Dlt){LLD};
        \tkzDefPoint(\PHISTAR,\PHISTAR+\Dlt){RLD};
        \tkzDefPoint(-\PHISTAR,\PHISTAR+\Dlt){LLD};
        \tkzDefPoint(\PHISTAR,\PHISTAR+\Dlt){RLD};
        \tkzDefPoint(-\PHISTAR,\WIDTH-\PHISTAR-\Dlt){LUD};
        \tkzDefPoint(\PHISTAR,\WIDTH-\PHISTAR-\Dlt){RUD};
        \tkzDefPoint(-\PHISTAR,\PHISTAR+\Dltplus){LLDp};
        \tkzDefPoint(\PHISTAR,\PHISTAR+\Dltplus){RLDp};
        \tkzDefPoint(-\PHISTAR,\WIDTH-\PHISTAR-\Dltplus){LUDp};
        \tkzDefPoint(\PHISTAR,\WIDTH-\PHISTAR-\Dltplus){RUDp};
        \tkzDefPoint(-\PHISTAR,\PHISTAR+\Dltplus+0.8){LLDdp};
        \tkzDefPoint(\PHISTAR,\PHISTAR+\Dltplus+0.8){RLDdp};
        \tkzDefPoint(-\PHISTAR,\WIDTH-\PHISTAR-\Dltplus-0.8){LUDdp};
        \tkzDefPoint(\PHISTAR,\WIDTH-\PHISTAR-\Dltplus-0.8){RUDdp};
        \tkzDefPoint (0.3*\PHISTAR,\PHISTAR+\Dltplus){montonestart};
        \tkzDefPoint (0.3*\PHISTAR,0.333*\WIDTH){montoneend};
        \draw [thin] (LL) --(LU);
        \draw [thin] (LU) --(RU);
        \draw [thin] (RL) --(RU);
        \draw [thin] (LL) --(RL);
        \draw [blue,ultra thin,dashed] (LL) --(IY);
        \draw [red,ultra thin,dashed]  (X)--(RL);
        \draw [blue,ultra thin,dashed] (RU) --(IX);
        \draw [red,ultra thin,dashed]  (LU)--(Y);
        \draw [purple,thick]  (LLD)--(RLD);
        \draw [purple,thick]  (LUD)--(RUD);
        \draw [green,ultra thin]  (LLDp)--(RLDp);
        \draw [green,ultra thin]  (LUDp)--(RUDp);
        \draw [green,ultra thin]  (LLDdp)--(RLDdp);
        \draw [green,ultra thin]  (LUDdp)--(RUDdp);
        \tkzLabelPoint[right](RU){$\pi$};
        \tkzLabelPoint[right](RL){$0$};
        \tkzLabelPoint[right](RLD){$\sin{\theta_1}=2r/R$};
        \tkzLabelPoint[right](RUD){$\sin{\theta_1}=2r/R$};
        \tkzLabelPoint[right](RLDp){$\sin{\theta_1}=\sqrt{4r/R}$};
        \tkzLabelPoint[right](RUDp){$\sin{\theta_1}=\sqrt{4r/R}$};
        \tkzLabelPoint[right](RLDdp){$\theta_1=\phistar-\Phistar$};
        \tkzLabelPoint[right](RUDdp){$\theta_1=\pi-\phistar+\Phistar$};
        \coordinate (RM) at ($(RL)!0.333!(RU)$);
        \coordinate (UpperM) at ($(LU)!0.5!(RU)$);
        \coordinate (LowerM) at ($(LL)!0.5!(RL)$);
        \coordinate (LM) at ($(LowerM)!0.333!(UpperM)$);
        \draw [dashed,ultra thin]  (LowerM)--(UpperM);
        \tkzLabelPoint[below](LL){$-\Phistar$};
        \tkzLabelPoint[below](RL){$+\Phistar$};
        \tkzLabelPoint[below](LowerM){$0$};
        \fill [orange, opacity=15/30](LLDp) -- (RLDp) -- (RUDp) -- (LUDp) -- cycle;
        \node[] at (0,0.5*\WIDTH){\Large $\mathfrak{D}$};
    \end{tikzpicture}
    \caption{Compact region $\mathfrak{D}=\overline{R_1\cup R_2\cup R_3\cup R_4}=\big\{(\Phi_1,\theta_1)\bigm|-\Phistar\le\Phi_1\le\Phistar,\sin{\theta_1}\ge2r/R\big\}$ where $D_{11}<0,D_{21}<0$}\label{fig:A_compactRegion}
\end{center}
\end{figure}
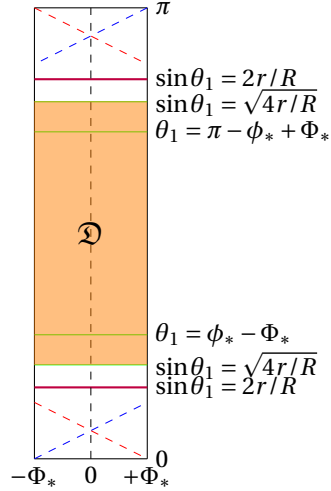


\begin{corollary}\label{corollary:cone-lowerbound}
If $\phistar\in\big(0,\tan^{-1}(1/3)\big)$ and \eqref{eqJZHypCond} holds, then for $D_{11}$, $D_{12}$, $D_{21}$, $D_{12}$ from \cref{def:DF2matrix,def:extensionDefs}, there exist  $\lambda_0(r,R,\phistar)>0$, $\lambda_1(r,R,\phistar)>0$ such that $\frac{D_{21}}{D_{11}}\ge \lambda_0(r,R,\phistar)$, $\frac{D_{22}}{D_{12}}\le \lambda_1(r,R,\phistar)$ for all $x_1$ in the compact region $\mathfrak{D}$ in \cref{lemma:extendedDomainDF2}.
\COMMENT{\cref{fig:A_compactRegion} suggests that $\mathfrak D$ is the shaded region. It is not. I suggest a brace on the left to mark the right height.\wentao{$\mathfrak{D}$ now is marked and mentioned in caption. And now the highlighted color for the region $\mathfrak{D}$ is correct, I think. I don't need to highlight case (a1) region I think in this figure now.}}
\end{corollary}

\begin{proof}
    By \cref{lemma:newcoordinateslength,lemma:extendedDomainDF2}, $\frac{D_{21}}{D_{11}}>0$ and is a continuous function of $x_1$ on the compact rectangle $\big\{x_1=(\Phi_1,\theta_1)\bigm|-\Phistar\le\Phi_1\le\Phistar,\sin^{-1}{(2r/R)}\le\theta_1\le(\pi-\sin^{-1}{(2r/R)})\big\}$ in \cref{fig:A_compactRegion}. Hence, on the compact rectangle $\frac{D_{21}}{D_{11}}$ attains its minimum  $\lambda_0(r,R,\phistar)>0$  determined by the billiard configuration. Similarly $\frac{D_{22}}{D_{12}}>0$ is a continuous function on $\mathfrak{D}$, so it attains a maximum $\lambda_1(r,R,\phistar)>0$ determined by the billiard configuration.
\end{proof}

\subsection[Narrow cone on \texorpdfstring{$\hat{M}$}{Mhat} and Proof of  \texorpdfstring{\cref{TEcase_a0}}{case a0 expansion}]{Invariance of half-quadrants and expansion in Case (a0)}\label{subsec:narrowcone}

\begin{corollary}[Strictly invariant half-quadrant]\label{corollary:halfquadrantcone}If $0<\phistar<\tan^{-1}{\big(1/3\big)}$ and \eqref{eqJZHypCond} holds, then the half-quadrant $\mathcal{HQ}_x(\mathrm{I,III})\dfn\big\{(d\phi,d\theta)\in T_x M\bigm|\frac{d\theta}{d\phi}\in[0,1]\big\}=C_x\subset T_x\hat{M}$ is strictly invariant under the differential $D\hat{F}$ of the return map (\cref{def:SectionSetsDefs}), i.e., $D\hat{F}_{x}(C_x)\subset\text{\big\{interior of }C_{\hat{F}(x)}\text{\big\}}$.\COMMENT{Of what is this a corollary? This is not obvious, so maybe "corollary" is the wrong word.}
\end{corollary}

\begin{proof}
Suppose a nonsingular $x\in\hat{M}$ has a return orbit segment containing nonsingular $x_0\in\Mrout$, $x_1=(\Phi_1,\theta_1)\in\MRin$, $x_2\in\Mrin$. By \cref{def:MhatReturnOrbitSegment} we have 
\[
x_0\in\Mrout\xrightarrow{\mathcal{F}}x_1=(\Phi_1,\theta_1)\in\MRin\xrightarrow{\mathcal{F}}\cdots \xrightarrow{\mathcal{F}}\mathcal{F}^{n_1}(x_1)\in\MRout\xrightarrow{\mathcal{F}}x_2\in\Mrin.
\]
Write 
$\mathcal{Q}_x(\mathrm{I,III})=\big\{\frac{d\theta}{d\phi}\in[0,+\infty]\big\}$.
In cases (a1), (b) and (c) of \eqref{eqMainCases}, \eqref{eq:strictlyinvcone_casea1bc} in \cref{corollary:InvQuadrant} gives 
\[D\hat{F}_x(\mathcal{HQ}_x(\mathrm{I,III}))\subset D\hat{F}_x(\mathcal{Q}_x(\mathrm{I,III}))\overbracket{\subset}^{\text{\eqref{eq:strictlyinvcone_casea1bc}}}\big\{\text{interior of }\mathcal{HQ}_{\hat{F}(x)}(\mathrm{I,III})\big\}.\]

In case (a0) of \eqref{eqMainCases}, since $x\in\hat{M}=(\Mrin\cap\Mrout)\sqcup\mathcal{F}^{-1}(\Mrout\setminus\Mrin)$, either $x=x_0$ or $\Fc(x)=x_0$ with $x$, $x_0\in\Mrout$ representing consecutive collisions on the same arc $\Gamma_r$ and $D\mathcal{F}_{x}=\begin{pmatrix}
1 & 2 \\
0 & 1 
\end{pmatrix}$ in $\phi\theta$-coordinates.\COMMENT{Unclear sentence---not clear where "or" ends.} So we have two possibilities:
\begin{equation}\label{eq:twoalternativefacts}
\begin{aligned}
    \text{If } x=x_0, &\text{ then } \mathcal{HQ}_x(\mathrm{I,III})=\mathcal{HQ}_{x_0}(\mathrm{I,III}).\\
    \text{If } x_0=\Fc(x), &\text{ then } D\Fc_x(\mathcal{HQ}_x(\mathrm{I,III}))\subset\mathcal{HQ}_{x_0}(\mathrm{I,III}).
\end{aligned}
\end{equation}
Then also in (a0) there are two subcases based on nonsingular $x_2\in\Mrin$.
\newline
\textbf{Subcase (i):} $m(x_2)\le1$, i.e., $x_2\in \textcolor{blue}{M^{\text{\upshape{in}}}_{r,0}}\cup\textcolor{blue}{M^{\text{\upshape{in}}}_{r,1}}\subset\hat{M}$. In this subcase, $\hat{F}(x)=x_2\in\hat{M}$. From \cref{def:DF2matrix} and \eqref{eq:DF2matrix}, 
\begin{equation}\label{eq:narrowconereason1}
    \begin{aligned}
        D_{22}-D_{12}=d_2d_1-2\tau_0d_2=d_2(d_1-2\tau_0)&\underbracket{>}_{\mathclap{\tau_0<2r}}d_2(\overbracket{d_1}^{\mathclap{\substack{\text{By \cref{lemma:returnorbittheta1Large,Corollary:CollisionOnMrin0Mrin1}}\\\sin\theta_1>1/2 \text{ and }d_1=R\sin\theta_1>R/2>850r}}}-4r)>0. 
    \end{aligned}
\end{equation} Together with $D_{22}<0$, $D_{12}<0$, \eqref{eq:narrowconereason1} indicates that $\frac{D_{22}}{D_{12}}\in(0,1)$. And the billiard map (and its iteration $D\Fc^2_{x_0}$) maintains orientation (\cite[equation (2.29)]{cb}) which indicates that $D_{22}D_{11}-D_{12}D_{21}>0$. Since by \cref{caseA}, $D_{22}$, $D_{11}$, $D_{12}$, $D_{21}$ are all negative, we get the following.
\begin{equation}\label{eq:narrowconereason2}
    \frac{D_{22}}{D_{12}}-\frac{D_{21}}{D_{11}}=\frac{D_{22}D_{11}-D_{12}D_{21}}{D_{12}D_{11}}>0.
\end{equation}
Hence $0<\frac{D_{21}}{D_{11}}<\frac{D_{22}}{D_{12}}<1$.\newline For $dx_0=(d\phi_0,d\theta_0)\in\mathcal{Q}_{x_0}(\mathrm{I,III})$, that is, $\frac{d\theta_0}{d\phi_0}=k$ for some $k\in[0,\infty]$, $D\Fc_{x_0}^2(dx_0)=dx_2=(d\phi_2,d\theta_2)$ satisfies 
\[
\frac{D_{21}+kD_{22}}{D_{11}+kD_{12}}-\frac{D_{21}}{D_{11}}=\frac{\overbracket{k}^{\mathclap{\ge0}}(\overbracket{D_{11}D_{22}-D_{21}D_{12}}^{\mathrlap{>0}})}{\underbracket{D_{11}(D_{11}+kD_{12})}_{>0}}\ge0>\frac{(\overbracket{D_{11}D_{22}-D_{21}D_{12}}^{\mathrlap{>0}})}{\underbracket{D_{12}(D_{11}+kD_{12})}_{>0}}=
    \frac{D_{22}}{D_{12}}-\frac{D_{21}+kD_{22}}{D_{11}+kD_{12}}\rlap{, so}
\]
\begin{equation}\label{eq:narrowconereason3}
    0<\frac{D_{21}}{D_{11}}\le\frac{d\theta_2}{d\phi_2}=\frac{D_{21}+kD_{22}}{D_{11}+kD_{12}}<\frac{D_{22}}{D_{12}}<1,
\end{equation}
and $D\mathcal{F}^2_{x_0}(\mathcal{Q}_{x_0}(\mathrm{I,III}))\subset\big\{\text{interior of }\mathcal{HQ}_{x_2}(\mathrm{I,III})\big\}$. We get two conclusions:\COMMENT{Use overbracket and mathclap on the last equal signs in each}\newline
If $x=x_0$, then 
\begin{align*}
D\hat{F}_x(\mathcal{HQ}_x(\mathrm{I,III}))&\overbracket{=}^{\mathclap{x=x_0,\Fc^2(x)=x_2=\hat{F}(x)}}D\Fc^2_x(\mathcal{HQ}_x(\mathrm{I,III}))\subset D\mathcal{F}^2_{x_0}(\mathcal{Q}_{x_0}(\mathrm{I,III}))\\
&\subset\big\{\text{interior of }\mathcal{HQ}_{x_2}(\mathrm{I,III})\big\}\overset{\hat{F}(x)=x_2}{=}\big\{\text{interior of }\mathcal{HQ}_{\hat{F}(x)}(\mathrm{I,III})\big\}.\\
\end{align*} If $x_0=\Fc(x)$, then
\begin{align*}
D\hat{F}_x(\mathcal{HQ}_x(\mathrm{I,III}))&\overbracket{=}^{\mathclap{\Fc(x)=x_0,\Fc^2(x_0)=x_2=\hat{F}(x)\quad}}D\Fc^3_x(\mathcal{HQ}_x(\mathrm{I,III}))\overbracket{=}^{\mathclap{\quad\text{chain rule and }x_0=\Fc(x)}}D\Fc^2_{x_0}(D\Fc_x(\mathcal{HQ}_x(\mathrm{I,III})))\overbracket{\subset}^{\mathclap{\text{\eqref{eq:twoalternativefacts}}}} D\Fc^2_{x_0}(\mathcal{HQ}_{x_0}(\mathrm{I,III}))\\
&\subset\big\{\text{interior of }\mathcal{HQ}_{x_2}(\mathrm{I,III})\big\}\overset{\hat{F}(x)=x_2}{=}\big\{\text{interior of }\mathcal{HQ}_{\hat{F}(x)}(\mathrm{I,III})\big\}.
\end{align*}\newline
\textbf{Subcase (ii):} $m(x_2)\ge2$, i.e., $\mathcal{F}^{m(x_2)-1}(x_2)=\hat{F}(x)$. Here, $D\mathcal{F}^{m(x_2)-1}_{x_2}=\begin{pmatrix}
1 & 2(m(x_2)-1) \\
0 & 1 
\end{pmatrix}$ satisfies \[D\mathcal{F}^{m(x_2)-1}_{x_2}\big(\big\{\text{interior of }\mathcal{Q}_{x_2}(\mathrm{I,III})\big\}\big)\subset\big\{\text{interior of }\mathcal{HQ}_{\hat{F}(x)}(\mathrm{I,III})\big\}.\]
If $x=x_0$, then $x_2=\Fc^2(x)$, $\Fc^{m(x_2)+1}=\hat{F}(x)$ and 
\begin{align*}
D\hat{F}_x(\mathcal{HQ}_x(\mathrm{I,III}))&\overbracket{=}^{\mathclap{\Fc^{m(x_2)+1}=\hat{F}(x)}}D\Fc^{m(x_2)+1}_x(\mathcal{HQ}_x(\mathrm{I,III}))\overbracket{\subset}^{\mathclap{x=x_0,\Fc^2(x)=x_2}} D\Fc^{m(x_2)-1}_{x_2}(D\mathcal{F}^2_{x_0}(\mathcal{Q}_{x_0}(\mathrm{I,III})))\\
&\overbracket{\subset}^{\mathclap{D\mathcal{F}^2_{x_0}\text{ is a matrix with negative entries}}}D\mathcal{F}^{m(x_2)-1}_{x_2}\Big(\big\{\text{interior of }\mathcal{Q}_{x_2}(\mathrm{I,III})\big\}\Big)\subset\big\{\text{interior of }\mathcal{HQ}_{\hat{F}(x)}(\mathrm{I,III})\big\}.\\
\end{align*}
If $x=\Fc(x_0)$, then $x_2=\Fc^3(x)$,  $\Fc^{m(x_2)+2}=\hat{F}(x)$ and
\begin{align*}
D\hat{F}_x(\mathcal{HQ}_x(\mathrm{I,III}))&\overbracket{=}^{\mathclap{\Fc^{m(x_2)+2}=\hat{F}(x)}}D\Fc^{m(x_2)+2}_x(\mathcal{HQ}_x(\mathrm{I,III}))\overbracket{\subset}^{\mathclap{\Fc(x)=x_0,\Fc^2(x_0)=x_2\text{ and chain rule}}} D\Fc^{m(x_2)-1}_{x_2}(D\mathcal{F}^2_{x_0}(D\Fc_x(\mathcal{Q}_{x}(\mathrm{I,III}))))\overbracket{\subset}^{\mathclap{\text{\eqref{eq:twoalternativefacts}}}} D\Fc^{m(x_2)-1}_{x_2}(D\mathcal{F}^2_{x_0}(\mathcal{Q}_{x_0}(\mathrm{I,III})))\\
&\overbracket{\subset}^{\mathclap{D\mathcal{F}^2_{x_0}\text{ is a matrix with negative entries}}}D\mathcal{F}^{m(x_2)-1}_{x_2}\Big(\big\{\text{interior of }\mathcal{Q}_{x_2}(\mathrm{I,III})\big\}\Big)\subset\big\{\text{interior of }\mathcal{HQ}_{\hat{F}(x)}(\mathrm{I,III})\big\}.\\
\end{align*}In both \textbf{subcases (i) and (ii)}, $D\hat{F}_x(\mathcal{HQ}_x(\mathrm{I,III}))\subset\big\{\text{interior of }\mathcal{HQ}_{\hat{F}(x)}(\mathrm{I,III})\big\}$.
\end{proof}

\begin{theorem}[Case (a0) \cref{TEcase_a0}\eqref{itemTEcase_a0-1}, \eqref{itemTEcase_a0-2}]\label{theorem:UniformExpansionA0}
Suppose a nonsingular $x\in \hat{M}$ (\cref{def:SectionSetsDefs}) has a return orbit segment $x,\mathcal{F}(x),\cdots,\mathcal{F}^{\sigma(x)}(x)=\hat{F}(x)\in \hat{M}$ as in \cref{def:MhatReturnOrbitSegment} with $x_{0}\in\Mrout$, $x_1=\mathcal{F}(x_0)=(\phi_1,\theta_1)\in \MRin$, $x_2\in\Mrin$ and $\sin\theta_1\ge \sqrt{4r/R}$ (case (a0) of \eqref{eqMainCases}; \cref{TEcase_a0}).

Then there exists a $\lambda_c(r,R,\phistar)>0$ such that for all $x\in\hat{M}=(\Mrin\cap\Mrout)\sqcup\mathcal{F}^{-1}(\Mrout\setminus\Mrin)$  (\cref{def:SectionSetsDefs}) and any tangent vector $dx=(d\phi,d\theta)\in C_x=\mathcal{HQ}_{x}=\Big\{(d\phi,d\theta)\bigm|\frac{d\theta}{d\phi}\in\big[0,1\big]\Big\}$ at $x$ (the half-quadrant from \cref{corollary:halfquadrantcone}) we have 
\begin{itemize} 
    \item $D\hat{F}(dx)\in\big\{\text{interior of }\mathcal{HQ}_{\hat{F}(x)}\big\}$,
    \item $\frac{\|D\hat{F}(dx)\|}{\|dx\|_\p}>1+\lambda_c$.
\end{itemize}
\end{theorem}These two conclusions are exactly the statements of \cref{TEcase_a0}\eqref{itemTEcase_a0-1}, \eqref{itemTEcase_a0-2}.

\begin{proof}
By \cref{corollary:halfquadrantcone}, $D\hat{F}(dx)\in\big\{\text{interior of }\mathcal{HQ}_{\hat{F}(x)}\big\}$. And by \cref{lemma:ExpansionLowerBoundA0A1}, we have
\begin{itemize} 
    \item If $x\in\Mrin\cap\Mrout$, then $dx_0\dfn dx$ and $(d\phi_2,d\theta_2)=dx_2\dfn D\mathcal{F}_{x}^2(dx)$. Otherwise, $x\in\mathcal{F}^{-1}(\Mrout\setminus\Mrin)$, then $dx_0\dfn D\mathcal{F}_x(dx)$ and $(d\phi_2,d\theta_2)=dx_2\dfn D\mathcal{F}_x^3(dx)$. In either case, $dx_0$ is a tangent vector at $x_0$ and $dx_2$ is a tangent vector at $x_2$.
    \item $\frac{\|dx_0\|_\p}{\|dx\|_\p}\ge 1$ and $\frac{\|D\hat{F}(dx)\|_\p}{\|dx_2\|_\p}\ge 1$.
    \item If $\mathcal{I}<0$ (see \cref{def:AformularForLowerBoundOfExpansion}),\COMMENT{What does this mean? "If  $\mathcal{I}<0$ (see \cref{def:AformularForLowerBoundOfExpansion})"?\wentao{OK, replaced with this statement.}} then $\frac{\|dx_2\|_\p}{\|dx_0\|_\p}\ge|\mathcal{I}|$.
\end{itemize}

Then we have two subcases on $x_2$.
\newline
\textbf{Case i:} $m(x_2)\le1$, i.e., $x_2\in\textcolor{blue}{M^{\text{\upshape{in}}}_{r,0}}\sqcup\textcolor{blue}{M^{\text{\upshape{in}}}_{r,1}}$. \cref{Corollary:CollisionOnMrin0Mrin1} gives $\mathcal{I}<-1.37$. So, $\frac{\|D\hat{F}(dx)\|}{\|dx\|_\p}=\frac{\|dx_0\|_\p}{\|dx\|_\p}\frac{\|dx_2\|_\p}{\|dx_0\|_\p}\frac{\|D\hat{F}(dx)\|_\p}{\|dx_2\|_\p}>1.37$.
\newline
\textbf{Case ii:}  $m(x_2)\ge2$. \cref{Prop:Ilowerbound} gives $\mathcal{I}<-1$. Since the billiard map (and its iteration $D\Fc^2_{x_0}$) maintains orientation (\cite[equation (2.29)]{cb}), by \cref{corollary:cone-lowerbound}, $\frac{d\theta_2}{d\phi_2}\in [\frac{D_{21}}{D_{11}}(x_1),\frac{D_{22}}{D_{12}}(x_1)]\subset[\lambda_0(r,R,\phistar),\infty)$ for all $x_1=(\Phi_1,\theta_1)$ with $\sin\theta_1\ge\sqrt{4r/R}>2r/R$. Therefore, by \cref{corollary:shearlowerboundexpansion}, we have $\frac{\|D\hat{F}(dx)\|_\p}{\|dx_2\|_\p}\ge 1+\frac{c_r\lambda_0(r,R,\phistar)}{2}$. We also get $\frac{\|D\hat{F}(dx)\|_\p}{\|dx\|_\p}=\frac{\|dx_0\|_\p}{\|dx\|_\p}\frac{\|dx_2\|_\p}{\|dx_0\|_\p}\frac{\|D\hat{F}(dx)\|_\p}{\|dx_2\|_\p}>1+\frac{c_r\lambda_0(r,R,\phistar)}{2}$. 

Hence, we can choose $\lambda_c=\min{\big\{0.37,\frac{c_r\lambda_0(r,R,\phistar)}{2}\big\}}$, which depends on $r$, $R$, $\phistar$.
\end{proof}
\subsection[Proof of \texorpdfstring{\cref{TEcase_a1}}{case (a1) contraction control}]{Contraction control in Case (a1)}\label{subsec:ProofTECase_a1}
\strut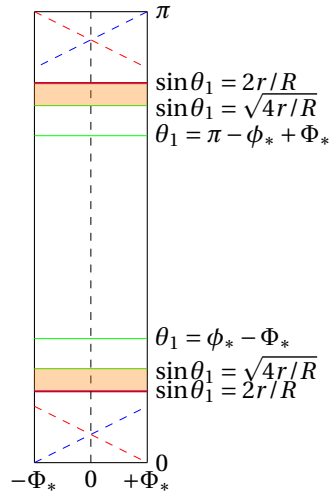
\begin{figure}[hb]
\begin{center}
    \begin{tikzpicture}[xscale=.5,yscale=.5]
        \pgfmathsetmacro{\PHISTAR}{1.5};
        \pgfmathsetmacro{\Dlt}{0.4};
        \pgfmathsetmacro{\Dltplus}{\Dlt+0.6}
        \pgfmathsetmacro{\WIDTH}{12};
        \tkzDefPoint(-\PHISTAR,0){LL};
        \tkzDefPoint(-\PHISTAR,\WIDTH){LU};
        \tkzDefPoint(\PHISTAR,0){RL};
        \tkzDefPoint(\PHISTAR,\WIDTH){RU};
        \tkzDefPoint(-\PHISTAR,\PHISTAR){X};
        \tkzDefPoint(\PHISTAR,\WIDTH-\PHISTAR){Y};
        \tkzDefPoint(-\PHISTAR,\WIDTH-\PHISTAR){IX};
        \tkzDefPoint(\PHISTAR,\PHISTAR){IY};
        \tkzDefPoint(-\PHISTAR,\PHISTAR+\Dlt){LLD};
        \tkzDefPoint(\PHISTAR,\PHISTAR+\Dlt){RLD};
        \tkzDefPoint(-\PHISTAR,\PHISTAR+\Dlt){LLD};
        \tkzDefPoint(\PHISTAR,\PHISTAR+\Dlt){RLD};
        \tkzDefPoint(-\PHISTAR,\WIDTH-\PHISTAR-\Dlt){LUD};
        \tkzDefPoint(\PHISTAR,\WIDTH-\PHISTAR-\Dlt){RUD};
        \tkzDefPoint(-\PHISTAR,\PHISTAR+\Dltplus){LLDp};
        \tkzDefPoint(\PHISTAR,\PHISTAR+\Dltplus){RLDp};
        \tkzDefPoint(-\PHISTAR,\WIDTH-\PHISTAR-\Dltplus){LUDp};
        \tkzDefPoint(\PHISTAR,\WIDTH-\PHISTAR-\Dltplus){RUDp};
        \tkzDefPoint(-\PHISTAR,\PHISTAR+\Dltplus+0.8){LLDdp};
        \tkzDefPoint(\PHISTAR,\PHISTAR+\Dltplus+0.8){RLDdp};
        \tkzDefPoint(-\PHISTAR,\WIDTH-\PHISTAR-\Dltplus-0.8){LUDdp};
        \tkzDefPoint(\PHISTAR,\WIDTH-\PHISTAR-\Dltplus-0.8){RUDdp};
        \tkzDefPoint (0.3*\PHISTAR,\PHISTAR+\Dltplus){montonestart};
        \tkzDefPoint (0.3*\PHISTAR,0.333*\WIDTH){montoneend};
        \draw [thin] (LL) --(LU);
        \draw [thin] (LU) --(RU);
        \draw [thin] (RL) --(RU);
        \draw [thin] (LL) --(RL);
        \draw [blue,ultra thin,dashed] (LL) --(IY);
        \draw [red,ultra thin,dashed]  (X)--(RL);
        \draw [blue,ultra thin,dashed] (RU) --(IX);
        \draw [red,ultra thin,dashed]  (LU)--(Y);
        \draw [purple,thick]  (LLD)--(RLD);
        \draw [purple,thick]  (LUD)--(RUD);
        \draw [green,ultra thin]  (LLDp)--(RLDp);
        \draw [green,ultra thin]  (LUDp)--(RUDp);
        \draw [green,ultra thin]  (LLDdp)--(RLDdp);
        \draw [green,ultra thin]  (LUDdp)--(RUDdp);
        \tkzLabelPoint[right](RU){$\pi$};
        \tkzLabelPoint[right](RL){$0$};
        \tkzLabelPoint[right](RLD){$\sin{\theta_1}=2r/R$};
        \tkzLabelPoint[right](RUD){$\sin{\theta_1}=2r/R$};
        \tkzLabelPoint[right](RLDp){$\sin{\theta_1}=\sqrt{4r/R}$};
        \tkzLabelPoint[right](RUDp){$\sin{\theta_1}=\sqrt{4r/R}$};
        \tkzLabelPoint[right](RLDdp){$\theta_1=\phistar-\Phistar$};
        \tkzLabelPoint[right](RUDdp){$\theta_1=\pi-\phistar+\Phistar$};
        \coordinate (RM) at ($(RL)!0.333!(RU)$);
        \coordinate (UpperM) at ($(LU)!0.5!(RU)$);
        \coordinate (LowerM) at ($(LL)!0.5!(RL)$);
        \coordinate (LM) at ($(LowerM)!0.333!(UpperM)$);
        \draw [dashed,ultra thin]  (LowerM)--(UpperM);
        \tkzLabelPoint[below](LL){$-\Phistar$};
        \tkzLabelPoint[below](RL){$+\Phistar$};
        \tkzLabelPoint[below](LowerM){$0$};
        \fill [orange, opacity=10/30](LLD) -- (RLD) -- (RLDp) -- (LLDp) -- cycle;
        \fill [orange, opacity=10/30](LUD) -- (RUD) -- (RUDp) -- (LUDp) -- cycle;
    \end{tikzpicture}
    \caption{Orange Region ($a1$) case,\, $2r/R\le\sin{\theta_1}< \sqrt{4r/R}$}\label{fig:A1_region}
\end{center}
\end{figure}
\begin{theorem}[In Case (a1) of \eqref{eqMainCases} contraction control]\label{theorem:contractioncontrolA1}
Suppose\COMMENT{Let us discuss the caption of \cref{fig:A1_region}} a nonsingular $x\in \hat{M}$ (\cref{def:SectionSetsDefs}) has a return orbit segment as in \cref{def:MhatReturnOrbitSegment}: $x,\mathcal{F}(x),\cdots,\mathcal{F}^{\sigma(x)}(x)=\hat{F}(x)\in \hat{M}$ with $x_{0}\in\Mrout$, $x_1=\mathcal{F}(x_0)=(\phi_1,\theta_1)\in \MRin$, $x_2\in\Mrin$ with $2r/R\ge\sin\theta_1< \sqrt{4r/R}$ in cases (a1) of \eqref{eqMainCases} and context of \cref{TEcase_a1}.\COMMENT{The last part of this sentence is incoherent.}

For $\forall x\in\hat{M}=(\Mrin\cap\Mrout)\sqcup\mathcal{F}^{-1}(\Mrout\setminus\Mrin)$ by \cref{def:SectionSetsDefs}, suppose a tangent vector at $x$, $dx=(d\phi,d\theta)\in C_x=\mathcal{HQ}_{x}=\Big\{(d\phi,d\theta)\bigm|\frac{d\theta}{d\phi}\in\big[0,1\big]\Big\}$ half-quadrant as defined in \cref{corollary:halfquadrantcone}.
\begin{itemize} 
    \item $D\hat{F}(dx)\in\big\{\text{interior of }\mathcal{HQ}_{\hat{F}(x)}\big\}$.
    \item $\frac{\|D\hat{F}(dx)\|}{\|dx\|_\p}>0.26$.
\end{itemize}
\end{theorem} These two conclusions are exactly the statements of \cref{TEcase_a1}\eqref{itemTEcase_a1-1}, \eqref{itemTEcase_a1-2}.
\begin{proof}
By \cref{corollary:halfquadrantcone}, $D\hat{F}(dx)\in\big\{\text{interior of }\mathcal{HQ}_{\hat{F}(x)}\big\}$.

By \cref{contraction_region}, we have $x_0\in\Nout$ and $x_2\in\Nin$. Then by \cref{def:SectionSetsDefs}, we have $x\in\mathcal{F}^{-1}(\Mrin\setminus\Mrout)$, $\mathcal{F}(x)=x_0$. Otherwise, $x=x_0\in\Mrin\cap\Mrout$ and $x_0\in\textcolor{red}{M^{\text{\upshape{out}}}_{r,0}}$, which contradicts \cref{Prop:Udisjoint} that $\Nout\cap\textcolor{red}{M^{\text{\upshape{out}}}_{r,0}}=\emptyset$. By \cref{lemma:ExpansionLowerBoundA0A1}, we have
\begin{itemize} 
    \item $x\in\mathcal{F}^{-1}(\Mrout\setminus\Mrin)$, $dx_0\dfn D\mathcal{F}_x(dx)$ and $(d\phi_2,d\theta_2)=dx_2\dfn D\mathcal{F}_x^3(dx)$. $dx_0$ is a tangent vector at $x_0$ and $dx_2$ is a tangent vector at $x_2$.
    \item $\frac{\|dx_0\|_\p}{\|dx\|_\p}\ge 1$ and $\frac{\|D\hat{F}(dx)\|_\p}{\|dx_2\|_\p}\ge 1$.
    \item If in \cref{def:AformularForLowerBoundOfExpansion} it has $\mathcal{I}<0$, then $\frac{\|dx_2\|_\p}{\|dx_0\|_\p}\ge|\mathcal{I}|$.
\end{itemize}
We assume $\theta_1\in(0,\pi/2]$, then $\theta_1\in\big[\sin^{-1}{(2r/R)},\sin^{-1}{\sqrt{4r/R}}\big)$. We have two subcases for $\Phi_1$.

\noindent\textbf{Subcase i}, $\Phi_1\in(-\Phistar,0)$, then $x_1\in R_2$ region of \cref{fig:R1}. By \cref{prop:casea0smallthetaR2}, we have $\mathcal{I}\le -1$.

\noindent\textbf{Subcase ii}, $\Phi_1\in[0,\Phistar)$, then $x_1\in R_1$ region of \cref{fig:R1}. Note that from \cref{def:AformularForLowerBoundOfExpansion}, we have $\tau_0+\tau_1>2d_0$, $\tau_1>0$, $d_0>0$ thus $\mathcal{I}<-1+\frac{\tau_1}{d_0}\cdot\frac{2(\tau_0-d_0)}{d_1}$.

For a fixed $\theta_1\in\Big[\sin^{-1}{(2r/R)},\sin^{-1}{\big(\sqrt{4r/R}\big)}\Big)$, as $\Phi_1$ varies from $0$ to $\Phistar$, $x_1=(\Phi_1,\theta_1)$ is variable depending on $\Phi_1\in[0,\Phistar)$, $\tau_1$, $d_0$, $d_2$ are functions of $x_1$ thus of $\Phi_1$. 

By \cref{monotone_p4}, $d_0$ as functions of $x_1$ is a monotone increasing function of $\Phi_1$. 

In \cref{subsec:ProofContractionRegion}, in the coordinate system defined in \cref{def:StandardCoordinateTable}, we have $(r\sin{\phi_2},-r\cos\phi_2)=(x_Q,y_Q)=Q=p(x_2)$ by \eqref{eq:XqXt}. By \eqref{eq:45}, in case (a1) we further have \[|\phi_2-\phistar|\overset{\eqref{eq:45}}<\frac{4.68}{\sin{(\phistar/2)}}\sqrt{r/R}\overset{\eqref{eqJZHypCond}}<\frac{4.68}{\sin{(\phistar/2)}}\sqrt{\sin^2{\phistar}/30000}<0.0148\cos{\big(\phistar/2\big)}<0.0270.\] Since $\phi_2\in(\phistar,2\pi-\phistar)$, $0<\phistar<\tan^{-1}{\big(1/3\big)}$, we have $\phi_2<\phistar+0.0270<\pi/2$. Thus, $x_Q>0$, by \cref{monotone_p3}(1), have $\frac{d\tau_1}{d\Phi_1}=-\frac{bx_Q}{d_2}<0$, so $\tau_1$ is monotonically decreasing with respect to $\Phi_1$. 

By symmetry, $d_2\bigm|_{x_1=(0,\theta_1)}=d_0|_{x_1=(0,\theta_1)}$. Hence we have 
\begin{equation}\label{eq:a1Ibound1}
0<\frac{\tau_1}{d_0}\le \frac{\tau_1}{d_0}\bigm|_{x_1=(\Phi_1,\theta_1)}\overbracket{\le}^{\mathclap{\text{\cref{monotone_p3,monotone_p4}}}}\frac{\tau_1}{d_0}\bigm|_{x_1=(0,\theta_1)}=\frac{\tau_1}{d_2}\bigm|_{x_1=(0,\theta_1)}\overbracket{<}^{\mathclap{\text{\cref{remark:deffpclf}}}}2
\end{equation}

In \cref{fig:heuristicproof}, $P=p(x_1)$ moves\COMMENT{This "moves" terminology is not helpful, but rampant. Over time, these occurrences should be rewritten. Low priority.} from the midpoint $M$  of $\arc{AB}$ to the corner $B$ with fixed angle $\theta_1$ between $\overrightarrow{PQ}$ and the tangential direction of $\Gamma_R$ at $P$. Using the $TP-P-O_R$ coordinate system,\COMMENT{"$TP-P-O_R$ coordinate system" is bizarre terminology. What does it mean?} $\rho$, $\delta$ the length and angle variables of \cref{def:TPPORsystem,lemma:newcoordinateslength,fig:fractionEstimate2}, in this subcase $\Phi_1\in[0,\Phistar)$,\COMMENT{This sentence is incomprehensible.} we have $\Phi_1\in[0,\Phistar)$, $0\le\delta<\phistar$ and\COMMENT{The next display needs fixing}
\begin{equation}\label{eq:a1Ibound2}
\begin{aligned}
     \frac{2(\tau_0-d_0)}{d_1}&\overset{\eqref{eq:NewCoordinateLegnth}}=\frac{2\rho\sin{(\theta_1+\delta)}}{d_1}\overset{0<\rho\le r, 2r\le d_1}\le\sin{(\theta_1+\delta)}\\
    &<|\sin{\theta_1}|+|\sin{\delta}|<\sqrt{4r/R}+\sin{\phistar}\overset{R>1700r,0<\phistar<\tan^{-1}{(1/3)}}<0.37\\
    \mathcal{I}&=-1+\frac{\tau_1}{d_0}\frac{2(\tau_1-d_0)}{d_1}-\overbracket{\frac{\tau_0+\tau_1-2d_0}{d_0}}^{\mathllap{\text{\cref{remark:symmetrywiththetangentline}:}>0}}<-1+\frac{\tau_1}{d_0}\frac{2(\tau_1-d_0)}{d_1}\\
    &\overbracket{<}^{\mathrlap{\text{\eqref{eq:a1Ibound1},\eqref{eq:a1Ibound2}}}}-1+2\times0.37=-0.26
\end{aligned}
\end{equation}    

Therefore, $\mathcal{I}<-0.26$ in either subcase. By \cref{lemma:ExpansionLowerBoundA0A1} we have $\frac{\|D\hat{F}(dx)\|_\p}{\|dx\|_\p}=\frac{\|D\hat{F}(dx)\|_\p}{\|dx_2\|_\p}\cdot\frac{\|dx_2\|_\p}{\|dx_0\|_\p}\cdot\frac{\|dx_0\|_\p}{\|dx\|_\p}>0.26$.\end{proof}

\section[Proof of \texorpdfstring{\cref{TEcase_b}}{case (b)}]{Contraction control in case (b)}\label{sec:ProofTECase_b}


In this section, we prove \cref{TEcase_b}\eqref{itemTEcase_b2}, 
 which covers Case (b) of \eqref{eqMainCases}. (This is \cref{subsec:ProofTECase_b};  \cref{TEcase_b}\eqref{itemTEcase_b1} is the conclusion of \cref{corollary:halfquadrantcone}.)
 Accordingly, all statements in this section are made in the context of \cref{TEcase_b}, which we now recall for the convenience of the reader. 

\hidesubsection{Standing assumptions for this section}We consider a nonsingular $x\in \hat{M}=(\Mrin\cap\Mrout)\sqcup\mathcal{F}^{-1}(\Mrout\setminus\Mrin)$ (\cref{def:SectionSetsDefs}) with return orbit segment $x,\mathcal{F}(x),\cdots,\mathcal{F}^{\sigma(x)}(x)=\hat{F}(x)\in \hat{M}$ as in \cref{def:MhatReturnOrbitSegment}, where $x_{0}\in\Mrout$, $x_1=\mathcal{F}(x_0)=(\phi_1,\theta_1)\in \MRin$, $x_2\in\Mrin$, $\sin\theta_1< 2r/R$, $n_1=0$ (i.e., Case (b) of \eqref{eqMainCases}),
and the tangent vector $dx=(d\phi,d\theta)\in\Big\{(d\phi,d\theta)\bigm|\frac{d\theta}{d\phi}\in\big[0,1\big]\Big\}$ at $x$ is in the half-quadrant cone in $\phi\theta$-coordinates.\COMMENT{It would be good to write ``$\phi\theta$-coordinates'' consistently!} \cref{lemma:ExpansionLowerBoundB} provides more information about this orbit segment.\NEWPAGE 
\begin{theorem}[\cref{TEcase_b}\eqref{itemTEcase_b2}]\label{theorem:proofofTEcase_b}\label{subsec:ProofTECase_b}
    Under the standing assumptions for this section, $\dfrac{\|D\hat{F}_x(dx)\|_\p}{\|dx\|_\p}>0.05$.
\end{theorem}
\begin{proof}
From facts below, 
    $\displaystyle\frac{\|D\hat{F}(dx)\|_\p}{\|dx\|_\p}=\underbracket{\frac{\|D\hat{F}(dx)\|_\p}{\|dx_3\|_\p}}_{\mathclap{>1\text{ (\cref{lemma:ExpansionLowerBoundB}\eqref{itemExpansionLowerBoundB2})}}}\overbracket{\frac{\|dx_3\|_\p}{\|dx_2\|_\p}}^{\mathclap{\ge1 \text{ (\cref{proposition:dx3dx2expansion})}}}\underbracket{\frac{\|dx_2\|_\p}{\|dx_0\|_\p}}_{\mathclap{=|\mathcal{I}(k)|\text{ (\cref{lemma:ExpansionLowerBoundB}\eqref{itemExpansionLowerBoundB3})}}}\overbracket{\frac{\|dx_0\|_\p}{\|dx\|_\p}}^{\mathclap{\ge1\text{ (\cref{lemma:ExpansionLowerBoundB}\eqref{itemExpansionLowerBoundB1})}}}>|\mathcal{I}(k)|\underbracket>_{\mathclap{\text{\cref{proposition:IktransionInCaseb}}}}0.05$.
\end{proof}
This provides the estimate in  \cref{TEcase_b}. The rest of this section  proves the ingredients just invoked. Readers wishing to skip the proofs of these facts can advance to page \pageref{eq:casebcase2minorterm4}.
\subsection{Lemmas for \texorpdfstring{\cref{TEcase_b}}{Case b}}\label{subsec:AFELBincaseb}
\begin{lemma}[Orbit configuration and contraction bound]\label{lemma:ExpansionLowerBoundB}
Under the standing assumptions for this section
\begin{enumerate}
    \item\label{itemExpansionLowerBoundB1} $x\in\mathcal{F}^{-1}(\Mrout\setminus\Mrin)\subset M_r$, $x_0=\mathcal{F}(x)\in\Mrout$, $dx_0\dfn D\mathcal{F}_x(dx)$ and $dx_2\dfn D\mathcal{F}_{x}^3(dx)$, where $dx_0$ is a tangent vector at $x_0$ and $dx_2$ is a tangent vector at $x_2$. We have $\frac{\|dx_0\|_\p}{\|dx\|_\p}\ge1$.
    \item\label{itemExpansionLowerBoundB2} $m(x_2)\ge3$, $x_3\dfn\mathcal{F}(x_2)\in M_r$, $dx_3\dfn D\mathcal{F}^4_x(dx)$ and  $\mathcal{F}^{m(x_2)-2}(x_3)=\mathcal{F}^{m(x_2)-1}(x_2)\overset{\text{\cref{def:SectionSetsDefs}}}{\scalebox{7}[1]{$=$}}\hat{F}(x)\in\hat{M}$, i.e. $x_3=\mathcal{F}(x_2)$, $\cdots$, $\mathcal{F}^{m(x_2)-1}(x_2)$ are all in $M_r\setminus\Mrout$. We have $\frac{\|D\hat{F}_x(dx)\|_\p}{\|dx_3\|_\p}>1$.
    \item\label{itemExpansionLowerBoundB3} We use the same notation as in \eqref{eq:MirrorEquationExpansionFormula} with index subscript denoting by $\mathcal{B}^{\pm}_0,\mathcal{B}^{\pm}_1$ the infinitesimal wave front curvatures after/before collisions at $x_0$, $x_1$, respectively.
    We note that $\Bcal^+_0\in\big[-\frac{-4}{3d_0},\frac{-1}{d_0}\big]$.
    
    Then $\frac{\|dx_2\|_\p}{\|dx_0\|_\p}=\big|\mathcal{I}(k)\big|$ where $\mathcal{I}(k)\dfn(\tau_0+\tau_1-\frac{2\tau_0\tau_1}{d_1})(\frac{-k}{d_0})+1-\frac{2\tau_1}{d_1}$ and $k=1+\frac{1}{2+\frac{1}{{d\theta}/{d\phi}}}\in[1,\frac{4}{3}]$ (since $\frac{d\theta}{d\phi}\in[0,1]$), and $\tau_0$, $\tau_1$, $d_0$, $d_1$, $d_2$ are the length functions (\cref{def:fpclf}) of $x_1\in\big\{(\Phi_1,\theta_1)\in\MRin\cap\MRout\bigm|\sin\theta_1<2r/R\big\}=\mathcal{L}\sqcup\mathcal{U}$ in \cref{fig:A1contractionCase}. 
    Hence the expansion is determined by the cone vector $dx$ and length functions at $x_1$.
\end{enumerate}
\end{lemma}
\begin{figure}[ht]
\begin{center}
    \begin{tikzpicture}[xscale=0.6, yscale=0.6]
        \pgfmathsetmacro{\PHISTAR}{1.5};
        \pgfmathsetmacro{\Dlt}{0.6};
        \pgfmathsetmacro{\Dltplus}{\Dlt+0.6}
        \pgfmathsetmacro{\WIDTH}{12};
        \tkzDefPoint(-\PHISTAR,0){LL};
        \tkzDefPoint(-\PHISTAR,\WIDTH){LU};
        \tkzDefPoint(\PHISTAR,0){RL};
        \tkzDefPoint(\PHISTAR,\WIDTH){RU};
        \tkzDefPoint(-\PHISTAR,\PHISTAR){X};
        \tkzDefPoint(\PHISTAR,\WIDTH-\PHISTAR){Y};
        \tkzDefPoint(-\PHISTAR,\WIDTH-\PHISTAR){IX};
        \tkzDefPoint(\PHISTAR,\PHISTAR){IY};
        \tkzDefPoint(-\PHISTAR,\PHISTAR+\Dlt){LLD};
        \tkzDefPoint(\PHISTAR,\PHISTAR+\Dlt){RLD};
        \tkzDefPoint(-\PHISTAR,\PHISTAR+\Dlt){LLD};
        \tkzDefPoint(\PHISTAR,\PHISTAR+\Dlt){RLD};
        \tkzDefPoint(-\PHISTAR,\WIDTH-\PHISTAR-\Dlt){LUD};
        \tkzDefPoint(\PHISTAR,\WIDTH-\PHISTAR-\Dlt){RUD};
        \tkzDefPoint(-\PHISTAR,\PHISTAR+\Dltplus){LLDp};
        \tkzDefPoint(\PHISTAR,\PHISTAR+\Dltplus){RLDp};
        \tkzDefPoint(-\PHISTAR,\WIDTH-\PHISTAR-\Dltplus){LUDp};
        \tkzDefPoint(\PHISTAR,\WIDTH-\PHISTAR-\Dltplus){RUDp};
        \tkzDefPoint(-\PHISTAR,\PHISTAR+\Dltplus+0.8){LLDdp};
        \tkzDefPoint(\PHISTAR,\PHISTAR+\Dltplus+0.8){RLDdp};
        \tkzDefPoint(-\PHISTAR,\WIDTH-\PHISTAR-\Dltplus-0.8){LUDdp};
        \tkzDefPoint(\PHISTAR,\WIDTH-\PHISTAR-\Dltplus-0.8){RUDdp};
        \tkzDefPoint (0.3*\PHISTAR,\PHISTAR+\Dltplus){montonestart};
        \tkzDefPoint (0.3*\PHISTAR,0.333*\WIDTH){montoneend};
        \draw [thin] (LL) --(LU);
        \draw [thin] (LU) --(RU);
        \draw [thin] (RL) --(RU);
        \draw [thin] (LL) --(RL);
        \draw [blue,ultra thin,dashed] (LL) --(IY);
        \draw [red,ultra thin,dashed]  (X)--(RL);
        \draw [blue,ultra thin,dashed] (RU) --(IX);
        \draw [red,ultra thin,dashed]  (LU)--(Y);
        \draw [purple,thick]  (LLD)--(RLD);
        \draw [purple,thick]  (LUD)--(RUD);
        \tkzLabelPoint[right](RU){$\pi$};
        \tkzLabelPoint[right](RL){$0$};
        \tkzLabelPoint[right](RLD){$\sin{\theta_1}=2r/R$, $d_1=2r$};
        \tkzLabelPoint[right](RUD){$\sin{\theta_1}=2r/R$, $d_1=2r$};
        \tkzLabelPoint[right](IY){$\theta_1=\Phistar$, $d_1=r\sin{\phistar}$};
        \tkzLabelPoint[right](Y){$\theta_1=\pi-\Phistar$, $d_1=r\sin{\phistar}$};
        \coordinate (RM) at ($(RL)!0.5!(RU)$);
        \tkzLabelPoint[right](RM){$\frac{\pi}{2}$};
        \coordinate (UpperM) at ($(LU)!0.5!(RU)$);
        \coordinate (LowerM) at ($(LL)!0.5!(RL)$);
        \coordinate (LM) at ($(LU)!0.5!(LL)$);
        \draw [ultra thin,dashed]  (LM)--(RM);
        \tkzLabelPoint[below](LL){$-\Phistar$};
        \tkzLabelPoint[below](RL){$+\Phistar$};
        \tkzLabelPoint[below](LowerM){$0$};
        
        \tkzDefPoint(0,0.5*\PHISTAR){LMend};
        \tkzDefPoint(0,\WIDTH-0.5*\PHISTAR){HMend};
        \tkzDefPoint(\PHISTAR,0.5*\PHISTAR){thetamin};
        \tkzDefPoint(\PHISTAR,\WIDTH-0.5*\PHISTAR){thetamax};
        \draw[ultra thin,dashed] (LMend)--(thetamin);
        \draw[ultra thin,dashed] (HMend)--(thetamax);
        \tkzLabelPoint[right](thetamin){$\theta_1=\frac{\Phistar}{2}$, $d_1=R\sin{(\Phistar/2)}=\frac{R\sin{\Phistar}}{2\cos{(\Phistar/2)}}=\frac{r\sin{\phistar}}{2\cos{(\Phistar/2)}}>\frac{r\sin{\phistar}}{2}$};
        \tkzLabelPoint[right](thetamax){$\theta_1=\pi-\frac{\Phistar}{2}$, $d_1=R\sin{(\Phistar/2)}=\frac{R\sin{\Phistar}}{2\cos{(\Phistar/2)}}=\frac{r\sin{\phistar}}{2\cos{(\Phistar/2)}}>\frac{r\sin{\phistar}}{2}$};
        \fill [green, opacity=10/30](LLD) -- (RLD) -- (IY) -- (LMend) -- (X) -- cycle;
        \fill [green, opacity=10/30](LUD) -- (RUD) -- (Y) -- (HMend) -- (IX) -- cycle;       \tkzDrawPoints(X,Y,IX,IY,LLD,LUD,RLD,RUD,LMend,HMend,thetamin,thetamax);
        \node[] at ([shift={(0,0.5*\PHISTAR)}]LMend) {\small $\mathcal{L}$};
        \node[] at ([shift={(0,-0.5*\PHISTAR)}]HMend) {\small $\mathcal{U}$};
    \end{tikzpicture}
    \caption{Two components $\mathcal{U}\sqcup\mathcal{L}=\big\{(\Phi_1,\theta_1)\in\MRin\cap\MRout\bigm|\theta_1<\sin^{-1}{(2r/R)}\big\}$ are regions in the context of \cref{TEcase_b}, $x_1\in\mathcal{U}\cup\mathcal{L}$.}\label{fig:A1contractionCase}
\end{center}
\end{figure}\COMMENT{The caption of \cref{fig:A1contractionCase} must be rewritten.}
\begin{proof}\strut\COMMENT{Are there any ideas in this proof that one could explain?}
\hidesubsection{Proof of \eqref{itemExpansionLowerBoundB1}} By \cref{contraction_region} in case (b) of \eqref{eqMainCases}, $x_0\in\Nout$ and \cref{Prop:Udisjoint} shows that $\Nout\cap\textcolor{red}{M^{\text{\upshape{out}}}_{r,0}}=\emptyset$. This means that it is impossible to have $x\in\Mrout\cap\Mrin=\textcolor{red}{M^{\text{\upshape{out}}}_{r,0}}$. Thus, $x\in\mathcal{F}^{-1}(\Mrout\setminus\Mrin)$, $x_0=\mathcal{F}(x)\in\Mrout$.
    
    Hence, $dx_0\dfn D\mathcal{F}_{x}(dx)$ is a tangent vector at $x_0$. Since $x$, $x_0$ are nonsingular on $M_r$ and $dx=(d\phi,d\theta)$ is a tangent vector in the quadrant cone, by \cref{proposition:ShearOneStepExpansion} we have $\frac{\|dx_0\|_\p}{\|dx\|_\p}\ge1$.
    
    We also have $x_1=\mathcal{F}(x_0)\in\MRin$, $dx_1\dfn D\mathcal{F}_{x_0}(dx_0)$. With $n_1=0$ in the orbit segment, we have $x_2=\mathcal{F}(x_1)\in\Mrin$, $dx_2\dfn D\mathcal{F}_{x_1}(dx_1)=D\mathcal{F}^2_{x_0}(dx_0)=D\mathcal{F}^3_{x}(dx)$, $dx_2$ is a tangent vector at $x_2$.

    \hidesubsection{Proof of \eqref{itemExpansionLowerBoundB2}} By \cref{contraction_region}, we have $x_2\in\Nin$ and by \cref{Prop:Udisjoint}, we have $\emptyset=\Nin\cap(\textcolor{blue}{M^{\text{\upshape{in}}}_{r,0}}\cup\textcolor{blue}{M^{\text{\upshape{in}}}_{r,1}}\cup\textcolor{blue}{M^{\text{\upshape{in}}}_{r,2}})$. Hence, $m(x_2)\ge3$. By definition of $m(x_2)$ and $\hat{M}$  (\cref{def:SectionSetsDefs}), $(\phi_3,\theta_3)=x_3\dfn\mathcal{F}(x_2),\mathcal{F}(x_3)$, $\cdots$, $\mathcal{F}^{m(x_2)-2}(x_3)=\mathcal{F}^{m(x_2)-1}(x_2)\in\hat{M}$ are nonsingular points on $M_r\setminus\Mrout$ since their $\Fc$-images are not on $M_R$. 
    
    Now $dx_3\dfn D\mathcal{F}^4_x(dx)=D\mathcal{F}_{x_2}(dx_2)$ is a tangent vector at $x_3$. By \cref{caseB}, $D\mathcal{F}^4_x$ is a negative matrix. Since $(d\phi_3,d\theta_3)=dx_3=D\mathcal{F}^4_x(dx)$ and $dx=(d\phi,d\theta)$ is in the quadrant with $\frac{d\theta}{d\phi}\ge0$, $(d\phi_3,d\theta_3)=dx_3$ is then in the interior of the quadrant, that is, $\frac{d\theta_3}{d\phi_3}>0$. Then \cref{proposition:ShearOneStepExpansion} implies $\frac{\|D\hat{F}_x(dx)\|_\p}{\|dx_3\|_\p}>1$.

    \hidesubsection{Proof of \eqref{itemExpansionLowerBoundB3}} Since $\frac{d\theta}{d\phi}\in[0,1]$, $(d\phi_0,d\theta_0)=dx_0=D\mathcal{F}_x(dx)$, $D\mathcal{F}_x=\begin{pmatrix}
        1 & 2\\
        0 & 1
    \end{pmatrix}$, that is $\begin{pmatrix} d\phi_0\\d\theta_0
    \end{pmatrix}=\begin{pmatrix}
        1 & 2\\
        0 & 1
    \end{pmatrix}\begin{pmatrix} d\phi\\d\theta
    \end{pmatrix}$, we obtain $\frac{d\theta_0}{d\phi_0}=\frac{d\theta}{d\phi+2d\theta}=\frac{1}{2+\frac{1}{\frac{d\theta}{d\phi}}}\in\big[0,\frac{1}{3}\big]$. Hence, by \eqref{eq:coordinateschange} and \eqref{eq:MirrorEquationExpansionFormula}, we have $\mathcal{B}^+_0\cos\varphi_0+\frac{1}{r}=\frac{d\varphi_0}{ds_0}=\frac{-1}{r}\frac{d\theta_0}{d\phi_0}\in[-\frac{1}{3r},0]$. Therefore, $\mathcal{B}^+_0=-\frac{1}{r\cos\varphi_0}(1+\frac{d\theta_0}{d\phi_0})$. For $k=1+\frac{1}{2+\frac{d\theta}{d\phi}}\in[1,\frac{4}{3}]$, $\mathcal{B}^+_0=\frac{-k}{r\cos\varphi_0}=\frac{-k}{r\sin\theta_0}\overbracket=^{\mathclap{\text{\cref{def:fpclf}}}}\frac{-k}{d_0}\in[-\frac{-4}{3d_0},\frac{-1}{d_0}]$. 
    
    Then 
    \begin{equation}\label{eq:Ik}
    \begin{aligned}
        \frac{\|dx_2\|_\p}{\|dx_0\|_\p}&\overset{\eqref{eq:dx2overdx0}}=|(\tau_0+\tau_1-\frac{2\tau_0\tau_1}{d_1})\Bcal^{+}_0+1-\frac{2\tau_1}{d_1}|=|(\tau_0+\tau_1-\frac{2\tau_0\tau_1}{d_1})(\frac{-k}{d_0})+1-\frac{2\tau_1}{d_1}|=|\mathcal{I}(k)|,
    \end{aligned}
    \end{equation}
     where $\mathcal{I}(k)=(\tau_0+\tau_1-\frac{2\tau_0\tau_1}{d_1})(\frac{-k}{d_0})+1-\frac{2\tau_1}{d_1}$.
\end{proof}




\begin{proposition}[Case (b): length restrictions]\label{proposition:lengthsrestrictionincaseB}
Under the standing assumptions for this section, the length functions  $\tau_0$, $\tau_1$, $d_0$, $d_1$, $d_2$ of  $x_1\in\big\{(\Phi_1,\theta_1)\in\MRin\cap\MRout\bigm|\sin\theta_1<2r/R\big\}$ (\cref{def:fpclf,fig:A1contractionCase}) satisfy
\COMMENT{All references to the following equations need to be fixed now, in each case a number 1,2,3, or 4 must be added to the label according to which specific equation is being used.\wentao{Work In progress}}
\begin{gather}\label{eq:lengthrestrictions1}
    \frac{1}{2}r\sin{\phistar}<\frac{r\sin{\phistar}}{2\cos{\frac{\Phistar}{2}}}=\frac{R\sin{\Phistar}}{2\cos{\frac{\Phistar}{2}}}=R\sin{\frac{\Phistar}{2}}\le d_1<2r,\\\label{eq:lengthrestrictions2}
    \tau_0<2d_1,\quad \tau_1<2d_1, \quad\text{and}\quad d_0+d_2<\tau_0+\tau_1<d_0+d_2+\frac{40.7r^2}{R\sin{\phistar}},\\\label{eq:lengthrestrictions3}
    0<\tau_0+\tau_1-2d_0<\frac{40.7r^2}{R\sin{\phistar}} \text{ and }0<\tau_0+\tau_1-2d_2<\frac{40.7r^2}{R\sin{\phistar}},\\\label{eq:lengthrestrictions4}
    |d_0-r\sin{\phistar}|<\frac{8.3r^2}{R\sin{\phistar}},\quad |d_2-r\sin{\phistar}|<\frac{8.3r^2}{R\sin{\phistar}}\quad\text{and}\quad|d_0-d_2|<\frac{16.6r^2}{R\sin{\phistar}}.
\end{gather}
\end{proposition}
\begin{proof}
    For $x_1=(\Phi_1,\theta_1)$ in case (b) of \eqref{eqMainCases} with $d_1<2r$, $n_1=0$, the smallest $\theta_1$ can be is $\frac{\Phistar}{2}$, which corresponds to the lowest point of the region $\mathcal{L}$ in \cref{fig:A1contractionCase}. The corresponding $d_1=R\sin(\Phistar/2)=\frac{R\sin\Phistar}{2\cos{(\Phistar/2)}}\overbracket{=}^{\mathllap{R\sin\Phistar=r\sin\phistar}}\frac{r\sin\phistar}{2\cos{(\Phistar/2)}}>\frac{r\sin\Phistar}{2}$, hence \eqref{eq:lengthrestrictions1}.

    In \eqref{eq:lengthrestrictions2},  $\tau_i<2d_1$, $i=0,1$ follows from \cref{remark:deffpclf}, and  
    $0<\tau_0+\tau_1-d_0-d_2<\frac{40.7r^2}{R\sin{\phistar}}$ is \cite[Equations (4.9)]{hoalb}.

    As to \eqref{eq:lengthrestrictions3}, $0<\tau_0+\tau_1-2d_i<\frac{40.7r^2}{R\sin{\phistar}}$ for $i=0,2$ is \cite[Equations (4.8)]{hoalb}.  
    
    Finally, $|d_i-r\sin{\phistar}|<\frac{8.3r^2}{R\sin{\phistar}}$, for $i=0,2$ is \cite[Equations (4.6)]{hoalb},
    and $|d_0-d_2|\le|d_0-r\sin{\phistar}|$+$|r\sin{\phistar}-d_2|<\frac{16.6r^2}{R\sin{\phistar}}$ by the triangle inequality.
\end{proof}
These simple estimates will be used throughout this section by referring to \cref{eq:lengthrestrictions1,eq:lengthrestrictions2,eq:lengthrestrictions3,eq:lengthrestrictions4}.\NEWPAGE
\subsection[\texorpdfstring{$x_2\rightarrow x_3$}{x2 to x3 expansion} expansion]{\texorpdfstring{$x_2\rightarrow x_3$}{x2 to x3} expansion}\label{subsec:x2Tox3Expansion} We now prove that $\frac{\|dx_3\|_\p}{\|dx_2\|_\p}\ge1$.
\begin{proposition}\label{Proposition:DF4ConeUpperBound}
\label{def:DF4matrix}
Under the standing assumptions for this section 
$G\dfn d_1d_2D\mathcal{F}^4_x=
    \begin{pmatrix}
        G_{11} & G_{12}\\
        G_{21} & G_{22}
    \end{pmatrix}$ is a negative matrix with 
$0<\frac{G_{22}}{G_{12}}<1$. 
\end{proposition}
\begin{proof}
By \cref{lemma:ExpansionLowerBoundB}\eqref{itemExpansionLowerBoundB1}, $x=\Fc^{-1}(x_0)$ and by \cref{caseB}, $D\Fc^{4}_{x}=D\Fc^{4}_{x_{-1}}$ is a $2$ by $2$ matrix with negative entries. Hence $G$ is a negative matrix. And equations \cite[Equations (4.11)(4.13)(4.15)(4.17)]{hoalb} give the following.
\begin{equation}\label{eq:Gmatrix}\begin{aligned}
            G_{11}&=6(\tau_1-\frac{2}{3}d_2)(\tau_0-d_0)-3d_1(\tau_0+\tau_1-d_0-\frac{2}{3}d_2)\text{ \cite[Equation 4.13]{hoalb}},\\
            G_{21}&=2(\tau_1-d_2)(\tau_0-d_0)-d_1(\tau_0+\tau_1-d_0-d_2)\text{ \cite[Equation 4.11]{hoalb}},\\
            G_{12}&=18(\tau_1-\frac{2}{3}d_2)(\tau_0-\frac{2}{3}d_0)-9d_1(\tau_0+\tau_1-\frac{2}{3}d_0-\frac{2}{3}d_2)\text{ \cite[Equation 4.17]{hoalb}},\\
            G_{22}&=6(\tau_1-d_2)(\tau_0-\frac{2}{3}d_0)-3d_1(\tau_0+\tau_1-\frac{2}{3}d_0-d_2)\text{ \cite[Equation 4.15]{hoalb}}.
\end{aligned}\end{equation} With $G_{12}<0$, $G_{22}<0$, it suffices to prove that $G_{12}-G_{22}<0$.
\begin{equation}\label{eq:G12G22diff}
        \begin{aligned}
            G_{12}-G_{22}&\overset{\mathclap{\eqref{eq:Gmatrix}}}=(12\tau_1-6d_2)(\tau_0-\frac{2}{3}d_0)-d_1(6\tau_0+6\tau_1-4d_0-3d_2)\\
            &=12(\tau_1-\frac{1}{2}d_2)(\tau_0-\frac{2}{3}d_0)-6d_1(\tau_0-\frac{2}{3}d_0+\tau_1-\frac{1}{2}d_2)\\
            & =6d_1(\tau_0-\frac{2}{3}d_0)(\tau_1-\frac{1}{2}d_2)(\frac{2}{d_1}-\frac{1}{\tau_0-\frac{2}{3}d_0}-\frac{1}{\tau_1-\frac{1}{2}d_2}).
        \end{aligned}   \end{equation}
We divide this into two cases.

\noindent\textbf{Case 1}: If $(\tau_0-\frac{2}{3}d_0)(\tau_1-\frac{1}{2}d_2)\le0$, then
\[G_{12}-G_{22}=12(\tau_1-\frac{1}{2}d_2)(\tau_0-\frac{2}{3}d_0)-6d_1(\tau_0-\frac{2}{3}d_0+\tau_1-\frac{1}{2}d_2)\le-6d_1(\underbracket{\tau_0+\tau_1}_{\mathclap{>d_0+d_2\text{ by \eqref{eq:lengthrestrictions2}}}}-\frac{2}{3}d_0-\frac{1}{2}d_2)\overbracket{<}^{\mathllap{\eqref{eq:lengthrestrictions2}}}-6d_1(\frac{1}{3}d_0+\frac{1}{2}d_2)<0.\] 

\noindent\textbf{Case 2}: If $(\tau_0-\frac{2}{3}d_0)(\tau_1-\frac{1}{2}d_2)>0$, then 
$\tau_0+\tau_1\overset{\eqref{eq:lengthrestrictions2}}>d_0+d_2$ implies \(\tau_0-\frac{2}{3}d_0>0\) and \(\tau_1-\frac{1}{2}d_2>0\), so once we prove $\frac{2}{d_1}<\frac{1}{\tau_0-\frac{2}{3}d_0}+\frac{1}{\tau_1-\frac{1}{2}d_2}$, \eqref{eq:G12G22diff} yields \[G_{12}-G_{22}=6d_1(\tau_0-\frac{2}{3}d_0)(\tau_1-\frac{1}{2}d_2)(\frac{2}{d_1}-\frac{1}{\tau_0-\frac{2}{3}d_0}-\frac{1}{\tau_1-\frac{1}{2}d_2})<0.\]

It remains to prove $\frac{2}{d_1}<\frac{1}{\tau_0-\frac{2}{3}d_0}+\frac{1}{\tau_1-\frac{1}{2}d_2}$. We do so in five subcases,  using \cref{proposition:lengthsrestrictionincaseB} frequently.

\noindent\textbf{Subcase 2.1}: If $\tau_0-\frac{2}{3}d_0\ge\frac{10}{9}d_1$, $\tau_1-\frac{1}{2}d_2>0$, then 
$2d_1\overset{\eqref{eq:lengthrestrictions2}}{>}\tau_0\ge\frac{2}{3}d_0+\frac{10}{9}d_1\Rightarrow \frac{8}{9}d_1>\frac{2}{3}d_2$, so $d_1\ge \frac{3}{4}d_0$. Therefore, we have
\(\displaystyle d_0+d_2+\frac{40.7r^2}{R\sin\phistar}\overbracket{>}^{\mathclap{\eqref{eq:lengthrestrictions2}}}\tau_0+\tau_1\overbracket{\ge}^{\mathclap{\tau_0\ge\frac{2}{3}d_0+\frac{10}{9}d_1}}\frac{2}{3}d_0+\frac{10}{9}d_1+\tau_1\), which implies \(\displaystyle\tau_1<\frac{1}{3}d_0+d_2-\frac{10}{9}d_1+\frac{40.7r^2}{R\sin\phistar}\) and hence
\begin{align*}
    0\overbracket<^{\mathclap{\text{Subcase 2.1}}}\tau_1-\frac{1}{2}d_2&<\frac{1}{3}d_0+\frac{1}{2}d_2-\frac{10}{9}d_1+\frac{40.7r^2}{R\sin\phistar}\overbracket{\le}^{d_1\ge\frac{3}{4}d_0}\frac{1}{3}d_0+\frac{1}{2}d_2-\frac{10}{9}\cdot\frac{3}{4}d_0+\frac{40.7r^2}{R\sin\phistar}=\frac{1}{2}(d_2-d_0)+\frac{40.7r^2}{R\sin\phistar}\\
    &\overbracket{<}^{\mathclap{\eqref{eq:lengthrestrictions4}}}\frac{16.6r^2}{2R\sin\phistar}+\frac{40.7r^2}{R\sin\phistar}=\frac{49r^2}{R\sin\phistar}\overbracket{<}^{\mathclap{R>\frac{196r}{\sin^2{\phistar}}}}\frac{r\sin\phistar}{4}\overset{\eqref{eq:lengthrestrictions1}}<\frac{d_1}{2}.
\end{align*}
thus, $\frac{2}{d_1}<\frac{1}{\tau_1-\frac{1}{2}d_2}<\frac{1}{\tau_1-\frac{1}{2}d_2}+\frac{1}{\tau_0-\frac{2}{3}d_0}$.

\noindent\textbf{Subcase 2.2}: If $0<\tau_0-\frac{2}{3}d_0\le\frac{1}{2}d_1$, $\tau_1-\frac{1}{2}d_2>0$, then we directly get $\frac{2}{d_1}\le\frac{1}{\tau_0-\frac{2}{3}d_0}<\frac{1}{\tau_0-\frac{2}{3}d_0}+\frac{1}{\tau_1-\frac{1}{2}d_2}$.

\noindent\textbf{Subcase 2.3}: $d_1\le\tau_0-\frac{2}{3}d_0<\frac{10}{9}d_1$, $\tau_1-\frac{1}{2}d_2>0$.

\NEWPAGE From \eqref{eq:lengthrestrictions2} we have $2d_1\overbracket{>}^{\mathclap{\eqref{eq:lengthrestrictions2}}}\tau_0\overbracket{\ge}^{\mathrlap{\text{Subcase 2.3}}}\frac{2d_0}{3}+d_1$, so $d_1> \frac{2}{3}d_0$. \cref{proposition:lengthsrestrictionincaseB} also gives 
\begin{align*}
    0\overbracket<^{\mathclap{\text{Subcase 2.3}}}\tau_1-\frac{1}{2}d_2\overbracket{<}^{\mathclap{\eqref{eq:lengthrestrictions2}}}&d_0+\frac{1}{2}d_2-\tau_0+\frac{40.7r^2}{R\sin\phistar}\overbracket{\le}^{\mathrlap{\text{Subcase 2.3}}\quad} d_0+\frac{1}{2}d_2-(\frac{2}{3}d_0+d_1)+\frac{40.7r^2}{R\sin\phistar}=\frac{1}{3}d_0+\frac{1}{2}d_2-d_1+\frac{40.7r^2}{R\sin\phistar}
    \\
    \overbracket{<}^{\quad\mathllap{d_1>\frac{2}{3}d_0}}&\frac{1}{3}d_0+\frac{1}{2}d_2-\frac{2}{3}d_0+\frac{40.7r^2}{R\sin\phistar}=\frac{1}{2}d_2-\frac{1}{3}d_0+\frac{40.7r^2}{R\sin\phistar}
    \\ \overbracket<^{\quad\mathllap{\eqref{eq:lengthrestrictions4}}}&\frac{1}{2}(r\sin{\phistar}+\frac{8.3r^2}{R\sin{\phistar}})-\frac{1}{3}(r\sin{\phistar}-\frac{8.3r^2}{R\sin{\phistar}})+\frac{40.7r^2}{R\sin{\phistar}}=\frac{1}{6}r\sin{\phistar}+\frac{(\frac{5}{6}\times8.3+40.7)r^2}{R\sin{\phistar}}
    \\
    \overset{\quad\mathllap{R>\frac{166r}{\sin^2\phistar}>\frac{(\frac{5}{6}\times8.3+40.7)r}{(\frac{5}{11}-\frac{1}{6})\sin^2\phistar}}}<&\frac{10}{11}\times\frac{1}{2}r\sin\phistar 
    \\ \overbracket<^{\quad\mathllap{\eqref{eq:lengthrestrictions4}: d_1>\frac{1}{2}r\sin\phistar}}&\frac{10}{11}d_1.
\end{align*}
Thus, $\frac{1}{\tau_1-\frac{1}{2}d_2}>\frac{11}{10d_1}$. The definition of \textbf{Subcase 2.3} includes $\frac{1}{\tau_0-\frac{2}{3}d_0}>\frac{9}{10d_1}$, so $\frac{2}{d_1}=\frac{11}{10d_1}+\frac{9}{10d_1}<\frac{1}{\tau_1-\frac{1}{2}d_2}+\frac{1}{\tau_0-\frac{2}{3}d_0}$.

\noindent\textbf{Subcase 2.4}: $0.7d_1\le\tau_0-\frac{2}{3}d_0<d_1, \tau_1-\frac{1}{2}d_2>0$.

We have $2d_1\overbracket{>}^{\mathclap{\eqref{eq:lengthrestrictions2}}}\tau_0\ge\frac{2}{3}d_0+0.7d_1$, so $d_1>\frac{2}{3.9}d_0$. \cref{proposition:lengthsrestrictionincaseB} also gives\COMMENT{This the first place where 1758 appears, and apparently the only one. We should refer to here from the discussion near \eqref{eqJZHypCond}. This is also the right time to comment on ``40.7,'' ``8.3,'' and ``39.''  How do they come about, what exact expressions do they replace, how close are they? Presumably $1758\approx39(40.7+\frac{8.3}{39}+\frac{8.3}{2})$.}
\begin{align*}
    0\overbracket{<}^{\mathclap{\text{Subcase 2.4}}}\tau_1-\frac{1}{2}d_2\overbracket{<}^{\mathclap{\eqref{eq:lengthrestrictions2}}}&d_0+\frac{1}{2}d_2-\tau_0+\frac{40.7r^2}{R\sin\phistar}\overbracket{\le}^{\mathclap{\substack{\text{Subcase 2.4}\\\tau_0\ge0.7d_1+\frac{2}{3}d_0\bigstrut}}} d_0+\frac{1}{2}d_2-(0.7d_1+\frac{2}{3}d_0)+\frac{40.7r^2}{R\sin\phistar}=\frac{1}{3}d_0+\frac{1}{2}d_2-0.7d_1+\frac{40.7r^2}{R\sin\phistar}\\
    \overbracket{<}^{\mathclap{d_1>\frac{2}{3.9}d_0}}&\frac{1}{2}d_2+(\frac{1}{3}-\frac{2\times0.7}{3.9})d_0+\frac{40.7r^2}{R\sin\phistar}=\frac{1}{2}d_2-\frac{1}{39}d_0+\frac{40.7r^2}{R\sin\phistar}\\ \overbracket{<}^{\mathclap{\eqref{eq:lengthrestrictions4}}}&\frac{1}{2}(r\sin{\phistar}+\frac{8.3r^2}{R\sin{\phistar}})-\frac{1}{39}(r\sin{\phistar}-\frac{8.3r^2}{R\sin{\phistar}})+\frac{40.7r^2}{R\sin{\phistar}}\\
    =&\big(\underbracket{-\frac{1}{39}+(40.7+\frac{8.3}{39}+\frac{8.3}{2})\frac{r}{R\sin^2\phistar}}_{<0\text{ when }R>\frac{1758r}{\sin^2{\phistar}}}\big)r\sin\phistar+\frac{1}{2}r\sin\phistar\overbracket{<}^{\mathclap{R>\frac{1758r}{\sin^2{\phistar}}}}\frac{1}{2}r\sin\phistar\overset{\eqref{eq:lengthrestrictions1}}<d_1.
\end{align*}
Together with $\frac{1}{\tau_0-\frac{2}{3}d_0}>\frac{1}{d_1}$ (from the definition of this subcase), this gives $\frac{2}{d_1}=\frac{1}{d_1}+\frac{1}{d_1}<\frac{1}{\tau_0-\frac{2}{3}d_0}+\frac{1}{\tau_1-\frac{1}{2}d_2}$.

\noindent\textbf{Subcase 2.5}: $\frac{1}{2}d_1<\tau_0-\frac{2}{3}d_0<0.7d_1$, $\tau_1-\frac{1}{2}d_2>0$.

$2d_1\overbracket{>}^{\mathllap{\eqref{eq:lengthrestrictions2}}}\tau_0\overbracket{>}^{\mathrlap{\text{Subcase 2.5}}}\frac{1}{2}d_1+\frac{2}{3}d_0$, so $d_1>\frac{4}{9}d_0$. \cref{proposition:lengthsrestrictionincaseB} also gives
\begin{align*}
    0\overbracket{<}^{\mathllap{\text{subcase 2.5}}}\tau_1-\frac{1}{2}d_2\overbracket{<}^{\eqref{eq:lengthrestrictions2}}&d_0+\frac{1}{2}d_2+\frac{40.7r^2}{R\sin\phistar}-\tau_0\overbracket{<}^{\mathclap{\tau_0>\frac{1}{2}d_1+\frac{2}{3}d_0}}d_0+\frac{1}{2}d_2+\frac{40.7r^2}{R\sin\phistar}-\frac{1}{2}d_1-\frac{2}{3}d_0\\
    \overbracket{<}^{\mathllap{d_1>\frac{4}{9}d_0}}&\frac{1}{3}d_0+\frac{1}{2}d_2-\frac{1}{2}\cdot\frac{4}{9}d_0+\frac{40.7r^2}{R\sin\phistar}=\frac{1}{9}d_0+\frac{1}{2}d_2+\frac{40.7r^2}{R\sin\phistar}\\ \overbracket{<}^{\mathllap{\eqref{eq:lengthrestrictions2}}}&\frac{1}{9}(r\sin\phistar+\frac{8.3r^2}{R\sin\phistar})+\frac{1}{2}(r\sin\phistar+\frac{8.3r^2}{R\sin\phistar})+\frac{40.7r^2}{R\sin\phistar}\\
    =&\frac{1}{9}r\sin{\phistar}+\frac{1}{2}r\sin{\phistar}+(\frac{1}{9}+\frac{1}{2})\frac{8.3r^2}{R\sin{\phistar}}+\frac{40.7r^2}{R\sin{\phistar}}\overbracket{<}^{\mathclap{R>\frac{174r}{\sin^2{\phistar}}}}\frac{7}{8}r\sin\phistar\overbracket{<}^{\mathllap{\eqref{eq:lengthrestrictions1}}}\frac{7}{4}d_1,
\end{align*}
where we used the fact that $R>\frac{174r}{\sin^2{\phistar}}\Rightarrow\big[(\frac{1}{9}+\frac{1}{2})\times8.3+40.7\big]\frac{r^2}{R\sin{\phistar}}<\big(\frac{7}{8}-\frac{1}{2}-\frac{1}{9}\big)r\sin{\phistar}$.

Together with $\frac{1}{\tau_0-\frac{2}{3}d_0}>\frac{10}{7d_1}$ (from the definition of the subcase), this gives $\frac{2}{d_1}=\frac{4}{7d_1}+\frac{10}{7d_1}<\frac{1}{\tau_0-\frac{2}{3}d_0}+\frac{1}{\tau_1-\frac{1}{2}d_2}$.
\end{proof}
\begin{proposition}\label{proposition:dx3dx2expansion}
Under the standing assumptions for this section we have $dx_2=(d\phi_2,d\theta_2)$ satisfies $\frac{d\theta_2}{d\phi_2}\in(-\infty,-1)\cup(0,+\infty)$ and  $\frac{\|dx_3\|_\p}{\|dx_2\|_\p}\ge1$.
\end{proposition}
\begin{proof}
By \cref{lemma:ExpansionLowerBoundB}\eqref{itemExpansionLowerBoundB1}, we have $x=\mathcal{F}^{-1}(x_0)=x_{-1}$. And by \cref{caseB}, $D\mathcal{F}^{4}_{x}=D\mathcal{F}^{4}_{x_{-1}}$ is a negative matrix. Then in \cref{lemma:ExpansionLowerBoundB}\eqref{itemExpansionLowerBoundB2}, for $dx=(d\phi,d\theta)$ in the cone $\big\{dx=(d\phi,d\theta)\bigm|\frac{d\theta}{d\phi}\in [0,1]\big\}$, we have $dx_3=D\mathcal{F}^4_{x}(dx)$ is a tangent vector at $x_3$ in the cone $\big\{dx_3=(d\phi_3,d\theta_3)\bigm|\frac{d\theta_3}{d\phi_3}\in(0,+\infty) \big\}$, $dx_2=D\mathcal{F}_{x}^{3}(dx)$ and $dx_3=D\mathcal{F}_{x_2}(dx_2)$ for $dx_2=(d\phi_2,d\theta_2)$ with $D\mathcal{F}_{x_2}=\begin{pmatrix}
1 & 2 \\
0 & 1 \end{pmatrix}$, that is, $\begin{pmatrix} d\phi_3\\d\theta_3
    \end{pmatrix}=\begin{pmatrix}
        1 & 2\\
        0 & 1
    \end{pmatrix}\begin{pmatrix} d\phi_2\\d\theta_2
    \end{pmatrix}$. $x_2$, $x_3$ are both nonsingular points in $M_r$ and we get a pair of equivalent ``continued-fraction'' equations:    \begin{equation}\label{eq:Twocontinuedfraction}
            \frac{d\theta_3}{d\phi_3}=\frac{d\theta_2}{d\phi_2+2d\theta_2} = \frac{1}{2}-\frac{1}{2}\frac{1}{1+2\frac{d\theta_2}{d\phi_2}},\qquad 
            \frac{d\theta_2}{d\phi_2}=\frac{d\theta_3}{d\phi_3-2d\theta_3} = -\frac{1}{2}+\frac{1}{2}\frac{1}{1-2\frac{d\theta_3}{d\phi_3}}.
    \end{equation}

For the negative matrix $G$ in \cref{def:DF4matrix} and with $\frac{d\theta}{d\phi}=k$ for some $k\in[0,1]$, 
\begin{align*}
    \begin{pmatrix} d\phi_3\\d\theta_3
    \end{pmatrix}=dx_3=D\Fc^4_x(dx)=\frac{1}{d_1d_2}\begin{pmatrix}
        G_{11} & G_{12}\\
        G_{21} & G_{22}
    \end{pmatrix}\begin{pmatrix} d\phi\\d\theta
    \end{pmatrix},
    \text{ so }\frac{d\theta_3}{d\phi_3}=\frac{G_{21}+kG_{22}}{G_{11}+kG_{12}}.
\end{align*}
Since the billiard map (and its iteration) is orientation-preserving \cite[Equation (2.29)]{cb}, the negative matrix $G$ in \cref{def:DF4matrix} has positive determinant: $G_{11}G_{22}-G_{12}G_{21}>0$. Hence $\frac{G_{22}}{G_{12}}-\frac{G_{21}}{G_{11}}=\frac{G_{11}G_{22}-G_{12}G_{21}}{G_{11}G_{12}}>0$ and\COMMENT{Format this to match the counterpart in Section 5.}
\begin{equation}\label{eq:casebdtheta3dphi3ratio}
\begin{aligned}
    \frac{d\theta_3}{d\phi_3}-\frac{G_{21}}{G_{11}}=\frac{G_{21}+kG_{22}}{G_{11}+kG_{12}}-\frac{G_{21}}{G_{11}}&=\frac{\overbracket{k}^{\mathllap{\ge0}}(\overbracket{G_{11}G_{22}-G_{21}G_{12}}^{\mathrlap{>0}})}{\underbracket{G_{11}(G_{11}+kG_{12})}_{>0}}\ge0,\\
    \frac{G_{22}}{G_{12}}-\frac{d\theta_3}{d\phi_3}=\frac{G_{22}}{G_{12}}-\frac{G_{21}+kG_{22}}{G_{11}+kG_{12}}&=\frac{(\overbracket{G_{11}G_{22}-G_{21}G_{12}}^{\mathrlap{>0}})}{\underbracket{G_{12}(G_{11}+kG_{12})}_{>0}}>0,\text{ hence}\\
    0<\frac{G_{21}}{G_{11}}\le\frac{d\theta_3}{d\phi_3}&=\frac{G_{21}+kG_{22}}{G_{11}+kG_{12}}<\frac{G_{22}}{G_{12}}\overbracket{<}^{\mathrlap{\text{\cref{def:DF4matrix}}}}1. 
\end{aligned}
\end{equation}
With $0<\frac{G_{21}}{G_{11}}\le\frac{d\theta_3}{d\phi_3}<\frac{G_{22}}{G_{12}}<1$, by \eqref{eq:Twocontinuedfraction} $\frac{d\theta_2}{d\phi_2}\in(-\infty,-1)\cup (0,+\infty]$ and we must analyze the following two cases:
\[
    \frac{d\theta_2}{d\phi_2}\in \big(0,+\infty\big]\overset{\eqref{eq:Twocontinuedfraction}}{\Longleftrightarrow} \frac{d\theta_3}{d\phi_3}\in\big(0,\frac{1}{2}\big]
    \quad\text{versus}\quad\frac{d\theta_2}{d\phi_2}\in \big[-\infty,-1\big) \overset{\eqref{eq:Twocontinuedfraction}}{\Longleftrightarrow} \frac{d\theta_3}{d\phi_3}\in\big[\frac{1}{2},1\big).
\]

In the first case, 
\cref{proposition:ShearOneStepExpansion} gives $\frac{\|dx_3\|_\p}{\|dx_2\|_\p}\ge1$.
%

In the second case, adapt the notations from \eqref{eq:MirrorEquationExpansionFormula}, denote by $\mathcal{B}_2^+$ the wave front curvature after collision at $x_2$, and note that the length of the free path between collisions $p(x_2)$ and $p(x_3)$ in the billiard table is $\tau_2=2r\sin\theta_2$. Then\COMMENT{Maybe overbracket and mathclap here...}
\[1\overbracket{<}^{\mathclap{\frac{d\theta_2}{d\phi_2}<-1}}-1-2\frac{d\theta_2}{d\phi_2}\overbracket{=}^{\mathclap{\eqref{eq:coordinateschange}}}-1-2r\frac{-d\varphi_2}{ds_2}\overbracket{=}^{\mathclap{\eqref{eq:MirrorEquationExpansionFormula}}}-1+\mathcal{B}_2^{+}\underbracket{2r\cos\varphi_2}_{\mathclap{=2r\sin\theta_2=\tau_2}}-2\underbracket{r\mathcal{K}}_{=-1}=1+\tau_2\mathcal{B}^+_2\overset{\eqref{eq:MirrorEquationExpansionFormula}}=\frac{\|dx_3\|_\p}{\|dx_2\|_\p}.\qedhere\]
\end{proof}

\subsection[\texorpdfstring{$x_0\longrightarrow x_1\longrightarrow x_2$}{x0 to x1 to x2} Transition Contraction]{\texorpdfstring{$x_0\longrightarrow x_1\longrightarrow x_2$}{x0 to x1 to x2} contraction}\label{subsec:A1contractionCase}\hfill

\begin{definition}\label{def:casebhandgfunctions}
Under the standing assumptions for this section, we define functions $g$, $h$ of $x_1\in\big\{(\Phi_1,\theta_1)\in\MRin\cap\MRout\bigm|\sin\theta_1<2r/R\big\}$  (\cref{fig:A1contractionCase}) by
\begin{equation}\label{eq:casebhandgfunctions}
        \begin{aligned}
            & h=\frac{\tau_0}{d_0}+\frac{\tau_1}{d_0}-2,\\
            & g=d_0-r\sin{\phistar}.
        \end{aligned}
    \end{equation}
\end{definition}
\begin{proposition}\label{proposition:casebhandgfunctions}
Under the standing assumptions for this section, \COMMENT{"From \cref{proposition:lengthsrestrictionincaseB}" is always weird because it is redundant. In this case it is not clear how the claims follow. Edit as before: omit the words that suggest that the reader memorize everything prior and make sure that every relation has a label with a reason.\wentao{Fixed}} \cref{proposition:lengthsrestrictionincaseB} and \eqref{eqJZHypCond}: $ R>\frac{30000r}{\sin^2\phistar}$ indicate\COMMENT{???????} $g$, $h$ from \cref{def:casebhandgfunctions} and $d_0$ satisfy the following.
       \COMMENT{The alignment is pointless, {\tt gathered} makes more sense. Refer to this often when 1.0003 shows up.\newline\boris Where do the numbers
       come from?\wentao{details added.}} 
       \begin{gather}\label{eq:casebghd0lengthrestrictions1}
    d_0-r\sin{\phistar}\overbracket{=}^{\mathllap{\eqref{eq:casebhandgfunctions}}}g\overbracket{\in}^{\mathrlap{\eqref{eq:lengthrestrictions4}}}\big(-\frac{8.3 r^2}{R\sin\phistar},\frac{8.3r^2}{R\sin\phistar}\big)\overbracket{\subset}^{\mathclap{R>\frac{30000r}{\sin^2\phistar}}} \big(-0.0003r\sin\phistar,0.0003r\sin\phistar\big),\\\label{eq:casebghd0lengthrestrictions2}
    0.9997r\sin\phistar\overbracket{<}^{\eqref{eq:casebghd0lengthrestrictions1}}d_0\overbracket{<}^{\eqref{eq:casebghd0lengthrestrictions1}}1.0003r\sin\phistar,\\\label{eq:casebghd0lengthrestrictions3}
    0\overbracket{<}^{\eqref{eq:lengthrestrictions3}}h\overbracket{=}^{\eqref{eq:casebhandgfunctions}}\frac{\tau_0}{d_0}+\frac{\tau_1}{d_0}-2=\frac{\tau_0+\tau_1-2d_0}{d_0}\overbracket{<}^{\eqref{eq:lengthrestrictions3}}\frac{\frac{40.7r^2}{R\sin\phistar}}{d_0}\overbracket{<}^{\eqref{eq:casebghd0lengthrestrictions2}}\frac{\frac{40.7r^2}{R\sin\phistar}}{0.9997r\sin\phistar}<\frac{40.8r}{R\sin^2{\phistar}},\\\label{eq:casebghd0lengthrestrictions4}
    d_2-r\sin{\phistar}\overbracket{\in}^{\mathrlap{\eqref{eq:lengthrestrictions4}}}\big(-\frac{8.3 r^2}{R\sin\phistar},\frac{8.3r^2}{R\sin\phistar}\big)\overbracket{\subset}^{\mathclap{R>\frac{30000r}{\sin^2\phistar}}} \big(-0.0003r\sin\phistar,0.0003r\sin\phistar\big),\\\label{eq:casebghd0lengthrestrictions5}
    0.9997r\sin\phistar\overbracket{<}^{\eqref{eq:casebghd0lengthrestrictions4}}d_2\overbracket{<}^{\eqref{eq:casebghd0lengthrestrictions4}}1.0003r\sin\phistar.
    \end{gather}
\end{proposition}
\begin{proposition}\label{proposition:IktransionInCaseb}
Under the standing assumptions for this section, $\mathcal{I}(k)=(\tau_0+\tau_1-\frac{2\tau_0\tau_1}{d_1})(\frac{-k}{d_0})+1-\frac{2\tau_1}{d_1}<-0.05$ for all $k\in[1,\frac{4}{3}]$ and all $x_1\in\big\{x_1=(\Phi_1,\theta_1)\in\MRin\cap\MRout\bigm|\sin\theta_1<2r/R\big\}$ (see \cref{fig:A1contractionCase}).
\end{proposition}
\begin{proof}
We first rewrite $\mathcal{I}(k)$ with $h$ from \cref{def:casebhandgfunctions}, that is,
\begin{equation}\label{eq:CaseBIkQuadraticPolynomial}
    \begin{aligned}
    \mathcal{I}(k)=&(\tau_0+\tau_1-\frac{2\tau_0\tau_1}{d_1})(\frac{-k}{d_0})+1-\frac{2\tau_1}{d_1}\overset{\eqref{eq:casebhandgfunctions}}=\big[(2+h)d_0-\frac{2\tau_0\tau_1}{d_1}\big]\big(\frac{-k}{d_0}\big)+1-\frac{2\tau_1}{d_1}\\
    =&-(2+h)k+\frac{2\tau_0\tau_1}{d_0d_1}k+1-\frac{2\tau_1}{d_1}=-2k+\frac{2\tau_0\tau_1}{d_0d_1}k+1-\frac{2\tau_1}{d_1}-hk\\
    =&-2k+1+\frac{2\tau_1}{d_1}(\frac{\tau_0}{d_0}k-1)-hk=-2k+1+\frac{2d_0}{d_1}\frac{\tau_1}{d_0}k(\frac{\tau_0}{d_0}-\frac{1}{k})-hk\\ \overbracket{=}^{\mathclap{\eqref{eq:casebhandgfunctions}: \frac{\tau_0}{d_0}=2+h-\frac{\tau_1}{d_0}}}& -2k+1+\frac{2d_0}{d_1}k\frac{\tau_1}{d_0}(2-\frac{1}{k}-\frac{\tau_1}{d_0})+\frac{2k\tau_1}{d_1}h-hk.
\end{aligned}
\end{equation}
Then, by symmetry, we assume that $x_1\in \mathcal{L}$ in \cref{fig:A1contractionCase} and consider two cases to analyze this: Case 1: $d_1\ge\frac{21}{32}r\sin{\phistar}$; Case 2: $d_1<\frac{21}{32}r\sin{\phistar}$.
\begin{figure}[h]
\begin{center}
    \begin{tikzpicture}
        \pgfmathsetmacro{\PHISTAR}{3};
        \pgfmathsetmacro{\Dlt}{0.9};
        \pgfmathsetmacro{\WIDTH}{5};
        \pgfmathsetmacro{\boudaryShift}{1.0};
        \pgfmathsetmacro{\LShift}{0.5};
        \tkzDefPoint(-\PHISTAR,0){LL};
        \tkzDefPoint(-\PHISTAR,\WIDTH){LU};
        \tkzDefPoint(\PHISTAR,0){RL};
        \tkzDefPoint(\PHISTAR,\WIDTH){RU};
        \tkzDefPoint(-\PHISTAR,\PHISTAR){X};
        \tkzDefPoint(\PHISTAR,\PHISTAR){IY};
        \tkzDefPoint(-\PHISTAR,\PHISTAR+\Dlt){LLD};
        \tkzDefPoint(\PHISTAR,\PHISTAR+\Dlt){RLD};
        \tkzDefPoint(0,0.5*\PHISTAR){LMD};
        \draw [ultra thick,red] (X) --(LLD);
        \draw [ultra thin,dashed] (LU) --(RU);
        \draw [ultra thin,dashed] (LU) --(LLD);
        \draw [ultra thin,dashed] (X) --(LL);
        \draw [ultra thin,dashed] (RU) --(RLD);
        \draw [ultra thin,dashed] (IY) --(RL);
        \draw [ultra thick,red] (IY) --(RLD);
        \draw [ultra thin,dashed] (LL) --(RL);
        \draw [red,ultra thick] (LMD) --(IY);
        \draw [red,ultra thick]  (X)--(LMD);
        \draw [red,ultra thick]  (LLD)--(RLD);
        \coordinate (LMDMark) at ([shift={(\PHISTAR,0)}]LMD);
        \coordinate (RBoundary) at ([shift={(-2*\boudaryShift,-\boudaryShift)}]IY);
        \coordinate (LBoundary) at ([shift={(2*\boudaryShift,-\boudaryShift)}]X);
         \coordinate (LlineBoundary) at ([shift={(2*\LShift,-\LShift)}]X);
         \coordinate (RlineBoundary) at ([shift={(-2*\LShift,-\LShift)}]IY);
        \coordinate (RlineBoundaryMark) at ([shift={(0,-\LShift)}]IY);
        \draw [ultra thick] (LBoundary) --(RBoundary);
        \draw [thick,blue] (LlineBoundary) --(RlineBoundary);
        \draw [ultra thin,blue,dashed] (RlineBoundary) --(RlineBoundaryMark);
        \coordinate (RBoundaryMark) at ([shift={(0,-\boudaryShift)}]IY);
        \draw [ultra thin,dashed]  (LMD)--(LMDMark);
        \draw [ultra thin,dashed] (RBoundary)--(RBoundaryMark);
        \fill [red, opacity=8/30](X) -- (LBoundary) -- (RBoundary) -- (IY) --(RLD)--(LLD)--cycle;
        \tkzDrawPoints(IY,RLD,LMD,RBoundaryMark,LMDMark,RlineBoundaryMark);
        \tkzLabelPoint[right](IY){$d_1=r\sin{\phistar}$};
        \tkzLabelPoint[right](RLD){$d_1=2r$};
        \tkzLabelPoint[right](RlineBoundaryMark){\textcolor{blue}{$d_1=\const$}};
        \tkzLabelPoint[right](RBoundaryMark){$d_1=\frac{21}{32}r\sin{\phistar}$};
        \tkzLabelPoint[right](LMDMark){$d_1=\frac{r\sin{\phistar}}{2\cos{\frac{\Phistar}{2}}}$};
        \tkzLabelPoint[below](LL){$-\Phistar$};
        \tkzLabelPoint[below](RL){$\Phistar$}
    \end{tikzpicture}
    \caption{
    \textbf{Case 1} Region of $\mathcal{L}$ in \cref{fig:A1contractionCase} with $d_1\ge\frac{21}{32}r\sin{\phistar}$
    }\label{fig:b_contractionCase1}
\end{center}
\end{figure}
\hidesubsection{Case 1} $x_1\in\mathcal{L}$ and $2r>d_1\ge\frac{21}{32}r\sin\phistar$ (see \cref{fig:b_contractionCase1}).

Suppose $x_1$ is on the line $d_1=\const\ge\frac{21}{32}r\sin{\phistar}$  inside $\mathcal{L}$ (the blue line in \cref{fig:b_contractionCase1}).
Then 
$\frac{(a+b)^2}{4}-ab=\frac{(a-b)^2}{4}\ge0\Longrightarrow ab\le\frac{(a+b)^2}{4}$ implies 
\begin{equation}\label{arithmeticgeometricmean}
\overbracket{\frac{\tau_1}{d_0}}^{=a}\overbracket{(2-\frac{1}{k}-\frac{\tau_1}{d_0})}^{=b}\le\frac{1}{4}(\overbracket{2-\frac{1}{k}}^{=a+b})^2.
\end{equation}
Thus
, we get 
\begin{align*}
    \mathcal{I}(k)&\overbracket{=}^{\mathclap{\eqref{eq:CaseBIkQuadraticPolynomial}}}-2k+1+\frac{2d_0}{d_1}k\times\frac{\tau_1}{d_0}(2-\frac{1}{k}-\frac{\tau_1}{d_0})
    +\frac{2k\tau_1}{d_1}h-hk\\
    &\overbracket{\le}^{\mathclap{\substack{d_0>0,d_1>0,k>0\\\text{and }\eqref{arithmeticgeometricmean}}}}-2k+1+\frac{2d_0}{d_1}k\times\frac{1}{4}(2-\frac{1}{k})^2+\frac{2k\tau_1}{d_1}h-hk=-2k+1+\frac{kd_0}{2d_1}(2-\frac{1}{k})^2+\frac{2k\tau_1}{d_1}h-hk\\
    &\overbracket{=}^{\mathclap{\eqref{eq:casebhandgfunctions}}}-2k+1+\frac{k}{2d_1}(2-\frac{1}{k})^2(g+r\sin\phistar)+\frac{2k\tau_1}{d_1}h-hk\\
    &= -2k+1+\frac{kr\sin{\phistar}}{2d_1}(2-\frac{1}{k})^2+\frac{kg}{2d_1}(2-\frac{1}{k})^2+\frac{2k\tau_1}{d_1}h-hk\\
    &=\underbracket{\frac{r\sin{\phistar}}{2d_1}}_{\mathclap{\text{(Case 1): }\le\frac12\frac{32}{21}}}\frac{(2k-1)^2}{k}+1-2k+\frac{kg}{2d_1}(2-\frac{1}{k})^2+(\frac{2\tau_1}{d_1}-1)hk\\
    &\le\frac{32}{21}\frac{(2k-1)^2}{2k}+1-2k+\frac{kg}{2d_1}(2-\frac{1}{k})^2+(\frac{2\tau_1}{d_1}-1)hk=\frac{32}{21}(1-\frac{1}{2k})(2k-1)-2k+1+\frac{kg}{2d_1}(2-\frac{1}{k})^2+(\frac{2\tau_1}{d_1}-1)hk\\
    &=-\frac{43}{21}+\underbracket{\frac{22k}{21}+\frac{16}{21k}}_{\mathclap{\text{increases with $k\in[1,\frac{4}{3}]$}}}+\frac{kg}{2d_1}(2-\frac{1}{k})^2+(\frac{2\tau_1}{d_1}-1)hk
    \le\underbracket{-\frac{43}{21}+\frac{22}{21}\times\frac{4}{3}+\frac{16}{21}\times\frac{3}{4}}_{=-5/63<-0.079}+\underbracket{\frac{kg}{2d_1}(2-\frac{1}{k})^2}_{<0.0005}+\underbracket{(\frac{2\tau_1}{d_1}-1)hk}_{\mathclap{<0.00544}}<-0.07,
\end{align*}
since 
\[(\frac{2\tau_1}{d_1}-1)hk\overbracket{<}^{\mathclap{\substack{1\le k\le4/3,\text{\eqref{eq:casebghd0lengthrestrictions3}: }h>0\\\text{\eqref{eq:lengthrestrictions2}: } \tau_1<2d_1}}}(2\cdot2-1)\cdot\frac{4}{3}h=4h\overbracket{<}^{\mathclap{\text{\eqref{eq:casebghd0lengthrestrictions3}}\ }}\frac{4\times40.8r}{R\sin^2{\phistar}}=\frac{163.2r}{R\sin^2{\phistar}}\overbracket{<}^{\mathclap{\text{\eqref{eqJZHypCond}: }R>\frac{30000r}{\sin^2{\phistar}}}}0.00544,\]
and\COMMENT{(6.11) below should be (6.10)\wentao{Thanks, fixed.}} \[\frac{kg}{2d_1}(2-\frac{1}{k})^2=\frac{\overbracket{g}^{\mathclap{\text{\eqref{eq:casebghd0lengthrestrictions1}: }g<0.0003r\sin{\phistar}}}}{\underbracket{d_1}_{\mathllap{\text{(Case 1): } d_1\ge\frac{21}{32}r\sin\phistar}}}\times\underbracket{(2k+\frac{1}{2k}-2)}_{\mathclap{\text{\quad increases with }k\in[1,\frac{4}{3}]}}<\frac{0.0003r\sin{\phistar}}{\frac{21}{32}r\sin{\phistar}}\times(2\times\frac{4}{3}+\frac{1}{2\times\frac{4}{3}}-2)<0.0005.\]
\begin{figure}[h]
\begin{minipage}{.45\textwidth}
\begin{center}
    \begin{tikzpicture}[xscale=0.9,yscale=0.9]
        \pgfmathsetmacro{\PHISTAR}{3};
        \pgfmathsetmacro{\Dlt}{0.9};
        \pgfmathsetmacro{\WIDTH}{5};
        \pgfmathsetmacro{\boudaryShift}{1.0};
        \pgfmathsetmacro{\LShift}{1.25};
        \tkzDefPoint(-\PHISTAR,0){LL};
        \tkzDefPoint(-\PHISTAR,\WIDTH){LU};
        \tkzDefPoint(\PHISTAR,0){RL};
        \tkzDefPoint(\PHISTAR,\WIDTH){RU};
        \tkzDefPoint(-\PHISTAR,\PHISTAR){X};
        \tkzDefPoint(\PHISTAR,\PHISTAR){IY};
        \tkzDefPoint(-\PHISTAR,\PHISTAR+\Dlt){LLD};
        \tkzDefPoint(\PHISTAR,\PHISTAR+\Dlt){RLD};
        \tkzDefPoint(0,0.5*\PHISTAR){LMD};
        \draw [ultra thick,red] (X) --(LLD);
        \draw [ultra thin,dashed] (LU) --(RU);
        \draw [ultra thin,dashed] (LU) --(LLD);
        \draw [ultra thin,dashed] (X) --(LL);
        \draw [ultra thin,dashed] (RU) --(RLD);
        \draw [ultra thin,dashed] (IY) --(RL);
        \draw [ultra thick,red] (IY) --(RLD);
        \draw [ultra thin,dashed] (LL) --(RL);
        \draw [red,ultra thick] (LMD) --(IY);
        \draw [red,ultra thick]  (X)--(LMD);
        \draw [red,ultra thick]  (LLD)--(RLD);
        \coordinate (LMDMark) at ([shift={(\PHISTAR,0)}]LMD);
        \coordinate (RBoundary) at ([shift={(-2*\boudaryShift,-\boudaryShift)}]IY);
        \coordinate (LBoundary) at ([shift={(2*\boudaryShift,-\boudaryShift)}]X);
         \coordinate (LlineBoundary) at ([shift={(2*\LShift,-\LShift)}]X);
         \coordinate (RlineBoundary) at ([shift={(-2*\LShift,-\LShift)}]IY);
        \coordinate (RlineBoundaryMark) at ([shift={(0,-\LShift)}]X);
        \draw [ultra thick] (LBoundary) --(RBoundary);
        \draw [thick,blue] (LlineBoundary) --(RlineBoundary);
        \draw [ultra thin,blue,dashed] (LlineBoundary) --(RlineBoundaryMark);
        \coordinate (RBoundaryMark) at ([shift={(0,-\boudaryShift)}]IY);
        \draw [ultra thin,dashed]  (LMD)--(LMDMark);
        \draw [ultra thin,dashed] (RBoundary)--(RBoundaryMark);
        \fill [red, opacity=8/30] (LBoundary) -- (RBoundary) -- (LMD) --cycle;
        \tkzDrawPoints(IY,RLD,LMD,RBoundaryMark,LMDMark,RlineBoundaryMark,RlineBoundary,LlineBoundary);
        \tkzLabelPoint[below right](RlineBoundaryMark){Line $\mathfrak{l}$};
        \tkzLabelPoint[above left](RBoundary){\small $d_1=\frac{21}{32}r\sin{\phistar}$};
        \tkzLabelPoint[below left](LMDMark){\small $d_1=\frac{r\sin{\phistar}}{2\cos{\frac{\Phistar}{2}}},\theta_1=\frac{1}{2}\Phistar$};
        \tkzLabelPoint[below](LL){$-\Phistar$};
        \tkzLabelPoint[below](RL){$\Phistar$}
         \tkzLabelPoint[left](LlineBoundary){$\mathcal{W}$};
        \tkzLabelPoint[right](RlineBoundary){$\mathcal{E}$};
    \end{tikzpicture}
    \captionsetup{width=0.85\textwidth}\caption{
    \textbf{Case 2}: $x_1\in\mathfrak{l}\subset\mathcal{L}$ in \cref{fig:A1contractionCase} 
    with $d_1<\frac{21}{32}r\sin{\phistar}$.
    The blue line segment $\overline{\mathcal{W}\mathcal{E}}=\mathfrak{l}\dfn\Big\{(\Phi_1,\theta_1=\sin^{-1}(\frac{\const}{R}))\bigm|\Phi_1\in\big[\Phistar-2\theta_1,2\theta_1-\Phistar\big],\frac{r\sin\phistar}{2}<\frac{r\sin{\phistar}}{2\cos{\frac{\Phistar}{2}}}<d_1=\const<\frac{21r\sin{\phistar}}{32}\Big\}$ has endpoints $\mathcal{W}$, $\mathcal{E}$ (``west'' and ``east'').
    }
\label{fig:b_contractionCase2}
\end{center}
\end{minipage}\hfill\vrule
\begin{minipage}{.5\textwidth}
\begin{center}
\begin{tikzpicture}[xscale=0.65,yscale=.65]
    \tkzDefPoint(0,0){Or};
    \tkzDefPoint(0,-4.8){Y};
    \pgfmathsetmacro{\rradius}{4.5};
    \pgfmathsetmacro{\bdist}{6.3};
    \tkzDefPoint(0,\bdist){OR};
    \pgfmathsetmacro{\phistardeg}{40};
    \pgfmathsetmacro{\XValueArc}{\rradius*sin(\phistardeg)};
    \pgfmathsetmacro{\YValueArc}{\rradius*cos(\phistardeg)};
    \tkzDefPoint(\XValueArc,-1.0*\YValueArc){B};    \tkzDefPoint(-1.0*\XValueArc,-1.0*\YValueArc){A};
    \pgfmathsetmacro{\Rradius}{veclen(\XValueArc,\YValueArc+\bdist)};
    \tkzDrawArc[name path=Cr,very thick](Or,B)(A);
    \tkzDrawArc[name path=CR,very thick](OR,A)(B);
    \pgfmathsetmacro{\PXValue}{\Rradius*sin(7)};
    \pgfmathsetmacro{\PYValue}{\bdist-\Rradius*cos(7)};
    \tkzDefPoint(\PXValue,\PYValue){P};
    \tkzDefPoint(0,\bdist-\Rradius){C};
    \draw[blue,ultra thin] (P)--(B);
    \draw[red,ultra thin,dashed](OR)--(P);
    \draw[red,ultra thin,dashed](OR)--(A);
    \draw[red,ultra thin,dashed](OR)--(B);
    \tkzDefPointBy[projection=onto OR--P](A)\tkzGetPoint{M};
    \coordinate (Q_far) at ($(A)!2!(M)$);
    \path[name path=linePQ] (P)--(Q_far);
    \path[name intersections={of=Cr and linePQ,by={Q}}];
    \tkzDrawPoints(A,B,P,Or,OR,Q);
    \begin{scope}[ultra thin,decoration={
    markings,
    mark=at position 0.5 with {\arrow{>}}}
    ] 
    \draw[blue,postaction={decorate}] (A)--(P);
    \draw[blue,postaction={decorate}] (P)--(Q);
    \end{scope}
    \tkzLabelPoint[below,rotate=0](P){$p(\mathcal{E})=P_{\mathcal{E}}=(x_\mathcal{E},y_\mathcal{E})$};
    \tkzLabelPoint[left,rotate=0](A){$A$};
    \tkzLabelPoint[right,rotate=0](Q){$Q_\mathcal{E}$};
    \tkzLabelPoint[below,rotate=0](B){$B$};
    \tkzLabelPoint[left,rotate=0](Or){$O_r$};
    \tkzLabelPoint[left,rotate=0](OR){$O_R$};
    \tkzMarkAngle[arc=lll,ultra thin,size=1.2](A,OR,P);
    \tkzMarkAngle[arc=lll,ultra thin,size=2.5](P,OR,B);
    \node[rotate=0] at ($(OR)+(-95:2.2)$)  {\small $2\theta_1$};
    \node[rotate=0] at ($(OR)+(-80:4)$)  {\small $2\Phistar-2\theta_1$};
\end{tikzpicture}
\captionsetup{width=0.85\textwidth}\caption{If $x_1=\mathcal{E}$, then $p(x_0)$ is at corner $A$, $|P_{\mathcal{E}}Q_{\mathcal{E}}|=\tilde{\tau}_1$ and $\tilde{\tau}_0=|AP_{\mathcal{E}}|=2d_1=2\cdot\const$}\label{fig:Fraction_tau1_d0_min}
\end{center}
\end{minipage}
\end{figure}
\hidesubsection{Case 2} $x_1\in\mathcal{L}$ in \cref{fig:A1contractionCase}\COMMENT{Refer to definition of $\mathcal{L}$ or a figure showing it.\wentao{referring the figure now}} and $d_1<\frac{21}{32}r\sin\phistar$ (see \cref{fig:b_contractionCase2}). Note that in this case $\Phistar/2<\theta_1<\Phistar$

Suppose $x_1=(\Phi_1,\theta_1)$ is on the line segment \[\mathfrak{l}\dfn\Big\{(\Phi_1,\theta_1=\sin^{-1}{(\frac{\const}{R})})\bigm|\Phi_1\in\big(\Phistar-2\theta_1,2\theta_1-\Phistar\big),\frac{r\sin{\phistar}}{2\cos{\frac{\Phistar}{2}}}<d_1=\const<\frac{21r\sin{\phistar}}{32}\Big\}\]in region $\mathcal{L}$. $\mathfrak{l}$ has a left end point $\mathcal{W}=(2\theta_1-\Phistar,\theta_1)$ and a right end point $\mathcal{E}=(\Phistar-2\theta_1,\theta_1)$. Note that $\mathcal{W}$, $\mathcal{E}$ are singular points, so $x_1$ cannot be $\mathcal{W}$ or $\mathcal{E}$.

Under the standing assumptions for this section, \COMMENT{This seems to be the wrong word. Are we really adding assumptions?\wentao{Much wording changed}} in the coordinate system given in \cref{def:StandardCoordinateTable} (\cref{fig:CaseAMonotone,fig:Fraction_tau1_d0_min}), $x_0=(\phi_0,\theta_0)$, $x_2=(\phi_2,\theta_2)$, $P=p(x_1)=(x_P,y_P)$, $Q=p(x_2)=(x_Q,y_Q)=(r\sin\phi_2,-r\cos\phi_2)$ and $T_0=p(x_0)=(x_{T_0},y_{T_0})=(r\sin\phi_0,-r\cos\phi_0)$.

Since we consider $\tau_0$, $\tau_1$, $d_1$, $d_2$ and $d_0$ as functions of $x_1=(\Phi_1,\theta_1)$ on line $\mathfrak{l}$, by \cref{monotone_p3}\eqref{itemmonotone_p3-1}\eqref{itemmonotone_p3-3}, we know 
\begin{align*}
    \frac{d\tau_1}{d\Phi_1}&=-\frac{bx_Q}{d_2},\\
    \frac{d\tau_0}{d\Phi_1}&=-\frac{bx_{T_0}}{d_0}.
\end{align*}
\COMMENT{\textcolor{red}{I have fixed a number of like errors before: in a context like this, the comma must be inside the display, not after it. Likewise with periods, etc. Please fix any occurrence that you see.\newline\wentao{Work in progress}}} By \cref{contraction_region}, \eqref{eq:casebcontractionregionx2_1} we have $|\phi_2-\phistar|\le\|x_2-Iy_*\|<\frac{14.6r}{R\sin{\phistar}}$ and $|\phi_0-(2\pi-\phistar)|\le\|x_0-x_*\|<\frac{14.6r}{R\sin{\phistar}}$ yielding the following inequalities.
\begin{align*}
    \phistar<&\phi_2<\phistar+\frac{14.6r}{R\sin{\phistar}}\overset{\eqref{eqJZHypCond}}<\tan^{-1}{(1/3)}+\frac{14.6r}{\sin{\phistar}}\frac{\sin^2\phistar}{30000r}<0.33<\pi/2,\\
    2\pi-\phistar>&\phi_0>2\pi-\phistar-\frac{14.6r}{R\sin{\phistar}}\overset{\eqref{eqJZHypCond}}>2\pi-\tan^{-1}{(1/3)}-\frac{14.6r}{\sin{\phistar}}\frac{\sin^2\phistar}{30000r}>5.96>\frac{3\pi}{2},
\end{align*} Therefore, with $x_Q=r\sin\phi_2>0$ and $x_{T_0}=r\sin\phi_0<0$, $\tau_1$ is a monotonically decreasing function of $\Phi_1\in(2\theta_1-\Phistar,\Phistar-2\theta_1)$ and $\tau_0$ is a monotonically increasing function of $\Phi_1\in(2\theta_1-\Phistar,\Phistar-2\theta_1)$. 

On the other hand, by \cref{monotone_p4} $d_0$ is a monotonically increasing function of $\Phi_1\in(2\theta_1-\Phistar,\Phistar-2\theta_1)$. Note that since $\mathfrak{l}\subset\MRin\cap\MRout$, the length functions in \cref{def:fpclf} are continuous on $\mathfrak{l}$. With $\theta_1$ fixed, we denote the length function limits at $\mathcal{E}$ by $\tilde{\tau}_0$, $\tilde{\tau}_1$, $\tilde{d}_0$ in the following.\COMMENT{Mention $\tau_0$ monotone\wentao{Done}}
\begin{equation}\label{eq:lengthfunc_Elimit}
\begin{aligned}
\tilde{\tau}_0&\dfn\lim_{\mathfrak{l}\ni x_1\to \mathcal{E}}\tau_0(x_1)=\lim_{\substack{x_1=(\Phi_1,\theta_1)\in\mathfrak{l}\\\Phi_1\to\Phistar-2\theta_1}}\tau_0,\\
    \tilde{\tau}_1&\dfn\lim_{\mathfrak{l}\ni x_1\to \mathcal{E}}\tau_1(x_1)=\lim_{\substack{x_1=(\Phi_1,\theta_1)\in\mathfrak{l}\\\Phi_1\to\Phistar-2\theta_1}}\tau_1,\\
    \tilde{d}_0&\dfn\lim_{\mathfrak{l}\ni x_1\to \mathcal{E}}d_0(x_1)=\lim_{\substack{x_1=(\Phi_1,\theta_1)\in\mathfrak{l}\\\Phi_1\to\Phistar-2\theta_1}}d_0.
\end{aligned}
\end{equation}
With $\theta_1$ fixed, $\frac{\tau_1}{d_0}$ is a decreasing function of $\Phi_1\in(2\theta_1-\Phistar,\Phistar-2\theta_1)$, therefore $\frac{\tau_1}{d_0}\ge\lim_{\mathfrak{l}\ni x_1\to \mathcal{E}}\frac{\tau_1}{d_0}(x_1)=\frac{\tilde{\tau}_1}{\tilde{d}_0}$. $\tilde{d}_0$ is the limit of $d_0(x_1)$ with $\mathfrak{l}\ni x_1\to\mathcal{E}$. Hence by \cref{proposition:casebhandgfunctions}\eqref{eq:casebghd0lengthrestrictions2}, $\tilde{d}_0$ and $\frac{\tilde{\tau}_1}{\tilde{d}_0}$ satisfy the following. \begin{equation}\label{eq:tau_d0_frac_min}
\begin{aligned}
    &0.9997r\sin\phistar\overbracket{\le}^{\text{\eqref{eq:casebghd0lengthrestrictions2}}}\tilde{d}_0\overbracket{\le}^{\text{\eqref{eq:casebghd0lengthrestrictions2}}}1.0003r\sin\phistar\\
    &\frac{\tau_1}{d_0}(x_1)\ge\frac{\tilde{\tau}_1}{\tilde{d}_0}\ge\frac{\tilde{\tau}_1}{1.0003r\sin\phistar}\\
\end{aligned}
\end{equation}

We observe that $\mathcal{E}$ in phase space represents the collision on $\Gamma_R$ with the billiard trajectory coming from corner $A$ (see \cref{fig:Fraction_tau1_d0_min}). In \cref{fig:Fraction_tau1_d0_min} and in the standard coordinate system of \cref{def:StandardCoordinateTable},\COMMENT{\wentao{One thing quite inconvenient to denote coordinates, have to always mention in the standard coordinate system of \cref{def:StandardCoordinateTable}.}\small\boris Not necessarily---in \cref{def:StandardCoordinateTable} we can say that these coordinates will henceforth be referred to simply as Cartesian coordinates. Then we can do that instead of "the standard coordinate system of \cref{def:StandardCoordinateTable}." Moreover, often the use of $x$ and $y$ is a strong hint in itself.} let 
\[P_\mathcal{E}\dfn p(\mathcal{E})=\big(x_\mathcal{E},y_\mathcal{E}\big)=\big(R\sin{(\Phistar-2\theta_1)},b-R\cos{(\Phistar-2\theta_1)}\big).\]  We can observe that for any $x_1=(\Phi_1,\theta_1)\in\mathfrak{l}$, that is, $P=p(x_1)=\big(x_P,y_P\big)=\big(R\sin\Phi_1,b-R\cos\Phi_1\big)$ with $2\theta_1-\Phistar<\Phi_1<\Phistar-2\theta_1$. Therefore, $x_P<x_{\mathcal{E}}$, $y_P<y_{\mathcal{E}}$. Moreover with corner $B=(x_B,y_B)$ in \cref{fig:Fraction_tau1_d0_min}, since $(x_Q,y_Q)=(r\sin\phi_2,-r\cos\phi_2)$ with $0<\phi_2<\pi/2$, we have 
\begin{align*}
x_P<&x_{\mathcal{E}}<x_{B}=r\sin\phistar<x_Q\\
y_P<&y_{\mathcal{E}}<y_{B}=-r\cos\phistar<y_Q
\end{align*}
Hence $\forall x_1\in\mathfrak{l}$, $\tau_1(x_1)=|PQ|=\sqrt{(x_Q-x_P)^2+(y_Q-y_P)^2}>\sqrt{(x_B-x_\mathcal{E})^2+(y_B-y_\mathcal{E})^2}=|P_\mathcal{E}B|$. Since $\tau_1(x_1)$ is a continuous function for $x_1\in\mathfrak{l}$, as a limit \[\tilde{\tau}_1=\lim_{\mathfrak{l}\ni x_1\to \mathcal{E}}\tau_1(x_1)\ge|P_\mathcal{E}B|.\]

We now compute a lower bound of $|P_\mathcal{E}B|$.

In \cref{fig:Fraction_tau1_d0_min}, $\measuredangle{AO_RP_{\mathcal{E}}}=2\theta_1$, $\measuredangle{P_{\mathcal{E}}O_RB}=2\Phistar-2\theta_1$. Using trigonometry, we get the following.

\begin{equation}\label{eq:PBLowerBound}
\begin{aligned}
 |P_{\mathcal{E}}B|=&2R\sin{(\Phistar-\theta_1)}\overbracket{=}^{\mathclap{\sin{(\Phistar-\theta_1)}=\sin{\Phistar}\cos{\theta_1}-\cos{\Phistar}\sin{\theta_1}}}2R\sin{\Phistar}\cos{\theta_1}-2R\cos{\Phistar}\sin{\theta_1}\overbracket=^{\mathclap{r\sin\phistar=R\sin\Phistar}}2r\sin\phistar\cdot\cos\theta_1-2R\sin\theta_1\cos\Phistar\\
 \overbracket=^{\quad\mathllap{d_1=R\sin\Phistar}}&2r\sin\phistar\cos\theta_1-2d_1\cos\Phistar\overbracket>^{\mathclap{\Phistar/2<\theta_1<\Phistar}}2r\sin\phistar\cos\Phistar-2d_1\cos\Phistar=2r\sin\phistar\cos\Phistar(1-\frac{d_1}{r\sin\phistar})\\
 \overbracket>^{\quad\mathllap{d_1<\frac{21}{32}r\sin\phistar}}&2r\sin\phistar\cos\Phistar(1-\frac{21}{32})=\frac{11}{16}r\sin\phistar\sqrt{1-\sin^2\Phistar}\overbracket=^{r\sin\phistar=R\sin\Phistar}\frac{11}{16}r\sin\phistar\sqrt{1-\frac{r^2\sin^2\phistar}{R^2}}.
 \end{aligned}
 \end{equation}
\eqref{eq:tau_d0_frac_min} and \eqref{eq:PBLowerBound} imply several useful inequalities:\COMMENT{Use overbracket and mathclap for the first one.}
\begin{equation}\label{eq:tau1d0lowerbound5over8}
     \frac{\tau_1}{d_0}\overset{\eqref{eq:tau_d0_frac_min}}{\ge}\frac{\tilde{\tau}_1}{\tilde{d}_0}\overset{\eqref{eq:tau_d0_frac_min}}{\ge}\frac{\tilde{\tau}_1}{1.0003r\sin\phistar}\overbracket{>}^{\mathclap{\substack{\tilde{\tau}_1\ge |P_{\mathcal{E}}B|\\ \text{and \eqref{eq:PBLowerBound}} }}}\frac{11r\sin{\phistar}}{16\times 1.0003r\sin{\phistar}}\sqrt{1-\frac{r^2\sin^2\phistar}{R^2}}\overbracket>^{\mathrlap{\eqref{eqJZHypCond}:R>1700r}\quad}0.6872>\frac{5}{8}.
 \end{equation}
On the other hand, if $\mathfrak{l}\ni x_1\to \mathcal{E}$, then $p(x_0)\to A$, $p(x_1)\to P_{\mathcal{E}}$. Therefore,
\[\tilde{\tau}_0\overbracket{=}^{\eqref{eq:lengthfunc_Elimit}}\lim_{\mathfrak{l}\ni x_1\to \mathcal{E}}\tau_1(x_1)\overbracket{=}^{\mathclap{\substack{\tau_1=\,\text{distance between }\\p(x_0),p(x_1)}}}\lim_{\mathfrak{l}\ni x_1\to \mathcal{E}}|p(x_0)p(x_1)|=|AP_{\mathcal{E}}|\overbracket{=}^{\mathclap{A,P_{\mathcal{E}}\text{ are on }\Gamma_R}}2d_1.\]

We extend the definitions of $g$, $h$ (\cref{def:casebhandgfunctions}) to $\tilde{h}$, $\tilde{g}$ on $\mathcal{E}$ by taking limits.\COMMENT{\LaTeX\ note: \textbackslash lim\textbackslash limits is unncessary in displays and bad in text, so it should never be used. \LaTeX\ does this right automatically.}\quad\COMMENT{\small\color{red}Sign error after "Therefore" below!\wentao{fixed now.}}\quad\COMMENT{The first two equations in \eqref{eq:casebcase2ghlimitatE} are backwards; they should start with the thing being defined.}
\begin{equation}\label{eq:casebcase2ghlimitatE}
    \begin{aligned}
    \frac{\tilde{\tau}_0}{\tilde{d}_0}+\frac{\tilde{\tau}_1}{\tilde{d}_0}-2&\overbracket{=}^{\mathclap{\eqref{eq:lengthfunc_Elimit}}}\lim_{\mathfrak{l}\ni x_1\to \mathcal{E}}{\big(\frac{\tau_0}{d_0}+\frac{\tau_1}{d_0}-2\big)}\overbracket{=}^{\mathclap{\text{\cref{def:casebhandgfunctions}}}}\lim_{\mathfrak{l}\ni x_1\to \mathcal{E}}{h}\nfd\tilde{h},\\
    \tilde{d}_0-r\sin{\phistar}&\overbracket{=}^{\mathclap{\eqref{eq:lengthfunc_Elimit}}}\lim_{\mathfrak{l}\ni x_1\to \mathcal{E}}{\big(d_0-r\sin\phistar\big)}\overbracket{=}^{\mathclap{\text{\cref{def:casebhandgfunctions}}}}\lim_{\mathfrak{l}\ni x_1\to \mathcal{E}}{g}\nfd\tilde{g},\\
    \text{Therefore, }-\tilde{h}&=2-\frac{\tilde{\tau}_0}{\tilde{d}_0}-\frac{\tilde{\tau}_1}{\tilde{d}_0}=2-\frac{2d_1}{\tilde{d}_0}-\frac{\tilde{\tau}_1}{\tilde{d}_0},\\
    \frac{\tilde{\tau}_1}{\tilde{d}_0}&=2-\frac{2d_1}{\tilde{d}_0}+\tilde{h}.
\end{aligned}
\end{equation}
Since $g$, $h$ are continuous functions of $x_1$ on $\mathfrak{l}$, by \cref{proposition:casebhandgfunctions}\eqref{eq:casebghd0lengthrestrictions1}\eqref{eq:casebghd0lengthrestrictions3}, $\tilde{g}$, $\tilde{h}$ as limits of $x_1\to\mathcal{E}$ satisfy the following.
\begin{equation}\label{eq:tildeghfunctionsrestriction}
    \begin{aligned}
\tilde{g}\overbracket{\in}^{\mathclap{\eqref{eq:casebghd0lengthrestrictions1}}}\big[-\frac{8.3 r^2}{R\sin\phistar},\frac{8.3r^2}{R\sin\phistar}\big]\overbracket{\subset}^{\mathclap{\text{\eqref{eqJZHypCond}: }R>\frac{30000r}{\sin^2{\phistar}}}}\big[-0.0003r\sin\phistar,0.0003r\sin\phistar\big],\qquad
0\overbracket{\le}^{\eqref{eq:casebghd0lengthrestrictions3}}\tilde{h}\overbracket{\le}^{\eqref{eq:casebghd0lengthrestrictions3}}\frac{40.8r}{R\sin^2\phistar}.
\end{aligned}
\end{equation}
For $k\in[1,\frac{3}{4}]$, the quadratic polynomial $\mathfrak{Q}(t)\dfn t(2-\frac{1}{k}-t)$\COMMENT{It is NOT necessary to use different funny Q's for the various quadratic polynomials; it is clear that each is "disposable" and will only be used briefly---so there can be no confusion.} is a decreasing function of $t\ge\frac{5}{8}\ge1-\frac{1}{2k}$, so 
\begin{equation}\label{eq:IkUpperBoundCasebCase2}
\begin{aligned}
    \mathcal{I}(k)\overbracket{=}^{\mathclap{\eqref{eq:CaseBIkQuadraticPolynomial}}}&-2k+1+\frac{2d_0}{d_1}k\times\underbracket{\frac{\tau_1}{d_0}(2-\frac{1}{k}-\frac{\tau_1}{d_0})}_{=\mathfrak{Q}(\frac{\tau_1}{d_0})}+\frac{2k\tau_1}{d_1}h-hk\\ \overbracket<^{\mathclap{\eqref{eq:tau1d0lowerbound5over8}: \frac{\tau_1}{d_0}\ge\frac{\tilde{\tau}_1}{\tilde{d}_0}>\frac{5}{8}}}&-2k+1+\frac{2d_0}{d_1}k\times\frac{\tilde{\tau}_1}{\tilde{d}_0}(2-\frac{1}{k}-\frac{\tilde{\tau}_1}{\tilde{d}_0})+\frac{2k\tau_1}{d_1}h-hk\\
    \overbracket{=}^{\mathclap{\text{\eqref{eq:casebcase2ghlimitatE}: }\frac{\tilde{\tau}_1}{\tilde{d}_0}=2-\frac{2d_1}{\tilde{d}_0}+\tilde{h}}}&-2k+1+\frac{2d_0}{d_1}k\times(2-\frac{2d_1}{\tilde{d_0}}+\tilde{h})(-\tilde{h}+\frac{2d_1}{\tilde{d_0}}-\frac{1}{k})+\frac{2k\tau_1}{d_1}h-hk\\
    =&-2k+1+\underbracket{\frac{2d_0}{d_1}k(2-\frac{2d_1}{\tilde{d_0}})(\frac{2d_1}{\tilde{d_0}}-\frac{1}{k})}_{\mathrlap{=\frac{2\tilde{d_0}}{d_1}k(2-\frac{2d_1}{\tilde{d_0}})(\frac{2d_1}{\tilde{d_0}}-\frac{1}{k})+2\frac{d_0-\tilde{d_0}}{d_1}k(2-\frac{2d_1}{\tilde{d_0}})(\frac{2d_1}{\tilde{d_0}}-\frac{1}{k})}}+\tilde{h}\frac{2d_0}{d_1}k(\frac{2d_1}{\tilde{d_0}}-\frac{1}{k})-\frac{2d_0}{d_1}\tilde{h}k(2-\frac{2d_1}{\tilde{d_0}})-\frac{2d_0}{d_1}k(\tilde{h})^2
    +h(\frac{2k\tau_1}{d_1}-k)\\
    =&\underbracket{-2k+1+\frac{2\tilde{d_0}}{d_1}k(2-\frac{2d_1}{\tilde{d_0}})(\frac{2d_1}{\tilde{d_0}}-\frac{1}{k})}_{<-0.0639\text{ by \eqref{eq:Ikfirst3items} below}}\\&\underbracket{+\tilde{h}\frac{2d_0}{d_1}k(\frac{2d_1}{\tilde{d_0}}-\frac{1}{k})-\frac{2d_0}{d_1}\tilde{h}k(2-\frac{2d_1}{\tilde{d_0}})-\frac{2d_0}{d_1}k(\tilde{h})^2+h(\frac{2k\tau_1}{d_1}-k)
    +2\frac{d_0-\tilde{d_0}}{d_1}k(2-\frac{2d_1}{\tilde{d_0}})(\frac{2d_1}{\tilde{d_0}}-\frac{1}{k})}_{\leq0\text{ when }R\text{ is large}}
\end{aligned}
\end{equation}
To bound this, we begin with the sum $-2k+1+\frac{2\tilde{d_0}}{d_1}k(2-\frac{2d_1}{\tilde{d_0}})(\frac{2d_1}{\tilde{d_0}}-\frac{1}{k})$ of the first three terms on the right-hand side of \eqref{eq:IkUpperBoundCasebCase2}.
Let $\mathfrak{e}=\frac{d_1}{\tilde{d}_0}>0$, which is close to $\frac{d_1}{r\sin{\phistar}}=\frac{\const}{r\sin\phistar}\in(\frac{1}{2},\frac{21}{32})$ in this case. Then 
\begin{equation}\label{eq:Ikfirst3items}
\begin{aligned}
    -2k+1+\frac{2\tilde{d_0}}{d_1}k(2-\frac{2d_1}{\tilde{d_0}})(\frac{2d_1}{\tilde{d_0}}-\frac{1}{k})=&-2k+1+\frac{2k}{\mathfrak{e}}(2-2\mathfrak{e})(2\mathfrak{e}-\frac{1}{k})=-2k+1+2k\frac{1}{\mathfrak{e}}(-4\mathfrak{e}^2+\frac{2\mathfrak{e}}{k}+4\mathfrak{e}-\frac{2}{k})\\
    =&-2k+1-8k\mathfrak{e}+4+8k-\frac{4}{\mathfrak{e}}=6k+5-\frac{4}{\mathfrak{e}}-8k\mathfrak{e}\\
    \underbracket{\le}_{\mathllap{8ke+\frac{4}{e}\ge4\times2\sqrt{2ke\times\frac{1}{e}}=8\sqrt{2k}}}&6k+5-8\sqrt{2k}\underbracket{=}_{\mathclap{\text{Let }\mu=\sqrt{2k}\in\big[\sqrt{2},\sqrt{8/3}\big]}}3\mu^2-8\mu+5\\
    \le& 3\times\frac{8}{3}-8\times\sqrt{\frac{8}{3}}+5=13-8\times\sqrt{\frac{8}{3}}<-0.0639.
\end{aligned}
\end{equation}
We now show that the last five terms of the inequality \eqref{eq:IkUpperBoundCasebCase2} can be arbitrarily small when $R$ is large. 
We reorganize this sum of five terms as follows \eqref{eq:IkUpperBoundCasebCase2MinorTerms}.
\begin{equation}\label{eq:IkUpperBoundCasebCase2MinorTerms}
\begin{aligned}
&+\tilde{h}\frac{2d_0}{d_1}k(\frac{2d_1}{\tilde{d_0}}-\frac{1}{k})-\frac{2d_0}{d_1}\tilde{h}k(2-\frac{2d_1}{\tilde{d_0}})-\frac{2d_0}{d_1}k(\tilde{h})^2+h(\frac{2k\tau_1}{d_1}-k)+2\frac{d_0-\tilde{d_0}}{d_1}k(2-\frac{2d_1}{\tilde{d_0}})(\frac{2d_1}{\tilde{d_0}}-\frac{1}{k})\\
=&\tilde{h}k\frac{4d_0}{\tilde{d}_0}-\tilde{h}\frac{2d_0}{d_1}-\frac{4d_0}{d_1}\tilde{h}k+\frac{4d_0}{\tilde{d}_0}\tilde{h}k-\frac{2d_0}{d_1}k(\tilde{h})^2+h(\frac{2k\tau_1}{d_1}-k)+2\frac{d_0-\tilde{d_0}}{d_1}k(2-\frac{2d_1}{\tilde{d_0}})(\frac{2d_1}{\tilde{d_0}}-\frac{1}{k})\\
=&\underbracket{-\tilde{h}d_0\big[k(-\frac{8}{\tilde{d}_0}+\frac{4}{d_1})+\frac{2}{d_1}\big]}_{\mathclap{\eqref{eq:casebcase2minorterm1}:<\frac{63r}{R\sin^2\phistar}}}\underbracket{-\frac{2d_0}{d_1}k(\tilde{h})^2}_{\le0}+\underbracket{hk(\frac{2\tau_1}{d_1}-1)}_{\mathclap{\eqref{eq:casebcase2minorterm2}:<\frac{164r}{R\sin^2\phistar}}}+\underbracket{2\frac{d_0-\tilde{d_0}}{d_1}(2-\frac{2d_1}{\tilde{d_0}})k(\frac{2d_1}{\tilde{d_0}}-\frac{1}{k})}_{\mathclap{\text{\eqref{eq:casebcase2minorterm3},\eqref{eq:casebcase2minorterm4}: }<\frac{50r}{R\sin^2\phistar}}}
\end{aligned}
\end{equation}
For $k\in[1,\frac{4}{3}]$, we have the following upper bounds for these four terms: $\displaystyle-\frac{2d_0}{d_1}k(\tilde{h})^2\le0$,
\begin{equation}\label{eq:casebcase2minorterm1}
\begin{aligned}\underbracket{-\tilde{h}d_0}_{\mathllap{\eqref{eq:tildeghfunctionsrestriction}:\le0}}\Big[\overbracket{k}^{\mathrlap{\ge1}}(\frac{-8}{\underbracket{\tilde{d}_0}_{\qquad\mathllap{\substack{\eqref{eq:tau_d0_frac_min}:\,0.9997r\sin\phistar\le\tilde{d}_0\\\le1.0003r\sin{\phistar}}}}}+\underbracket{\frac{4}{d_1})+\frac{2}{d_1}}_{\mathrlap{\text{case (2) and \eqref{eq:lengthrestrictions1}: }\frac{21r\sin{\phistar}}{32}>d_1>0.5r\sin\phistar}}\Big]\overbracket{<}^{\mathrlap{1\le k\le4/3}}&-\tilde{h}d_0\Big(\frac{-8k}{\tilde{d}_0}+\frac{4k+2}{d_1}\Big)<-\tilde{h}d_0\big[-8\frac{4}{3\tilde{d}_0}+\frac{4\cdot1+2}{(21/32)r\sin\phistar}\big]\\
<&-\tilde{h}d_0\Big(-\frac{32}{3\cdot0.9997r\sin\phistar}+\frac{6}{(21/32)r\sin\phistar}\Big)\\
<&\frac{-1.521\tilde{h}d_0}{r\sin\phistar}\\
\overbracket{<}^{\mathllap{\substack{\eqref{eq:casebghd0lengthrestrictions2}:\,0.9997r\sin{\phistar}<d_0\\<1.0003r\sin{\phistar}}}}&-1.521\cdot0.9997\tilde{h}\overbracket{<}^{\mathclap{\eqref{eq:tildeghfunctionsrestriction}:\tilde{h}\le\frac{40.8r}{R\sin^2\phistar}}}\qquad\frac{63r}{R\sin^2\phistar},\end{aligned}
\end{equation}
and
\begin{equation}\label{eq:casebcase2minorterm2}
\overbracket{h}^{\mathllap{\eqref{eq:casebghd0lengthrestrictions3}:h>0}}\underbracket{k}_{\mathclap{k\ge1}}(\frac{2\tau_1}{d_1}-1)\overbracket<^{\mathclap{\substack{\text{\eqref{eq:lengthrestrictions2}}:\\\tau_1<2d_1\quad}}}3hk\underbracket{\le}_{\mathclap{h>0,k\le4/3}}4h\overbracket{<}^{\mathclap{{\textstyle\strut}\eqref{eq:casebghd0lengthrestrictions3}:\, h<\frac{40.8r}{R\sin^2{\phistar}}}}\qquad\frac{164r}{R\sin^2\phistar}.
\end{equation}
%
From case 2 condition $d_1<\frac{21}{32}r\sin\phistar$, \eqref{eq:lengthrestrictions2}: $0.5r\sin\phistar<d_1$, \eqref{eq:tildeghfunctionsrestriction}, \eqref{eq:tau_d0_frac_min}, \cref{proposition:casebhandgfunctions}\eqref{eq:casebghd0lengthrestrictions2}, \eqref{eq:casebcase2ghlimitatE} and \eqref{eq:tildeghfunctionsrestriction}, we obtain the following.
\begin{equation}\label{eq:casebcase2minorterm3}
\begin{aligned}
k(\frac{2d_1}{\tilde{d_0}}-\frac{1}{k})=\frac{2d_1k}{\tilde{d_0}}-1&\overbracket{\in}^{\mathclap{\substack{0.5r\sin\phistar<d_1<(21/32)r\sin\phistar\\ \eqref{eq:tau_d0_frac_min}:\,0.9997r\sin\phistar<\tilde{d}_0<1.0003r\sin\phistar\\1\le k\le4/3}}}\qquad(2\cdot\frac{0.5r\sin{\phistar}}{1.0003r\sin{\phistar}}-1,2\cdot\frac{(21/32)r\sin{\phistar}}{0.9997r\sin{\phistar}}\frac{4}{3}-1)\subset(-0.003,0.751),\\
2-\frac{2d_1}{\tilde{d_0}}\qquad&\overbracket{\in}^{\mathclap{\substack{0.5r\sin\phistar<d_1<(21/32)r\sin\phistar\\ \eqref{eq:tau_d0_frac_min}:\,0.9997r\sin\phistar<\tilde{d}_0<1.0003r\sin{\phistar}\\\text{ }}}}
\qquad(2-2\cdot\frac{(21/32)r\sin{\phistar}}{0.9997r\sin{\phistar}},2-2\times\frac{0.5\times r\sin{\phistar}}{1.0003r\sin{\phistar}})\subset(0.687,1.0003),\\
d_0-\tilde{d_0}&\overbracket{\in}^{\mathllap{\substack{\eqref{eq:casebghd0lengthrestrictions1}\eqref{eq:casebcase2ghlimitatE}:g-\tilde{g}=d_0-\tilde{d}_0\\\eqref{eq:tildeghfunctionsrestriction}:\tilde{g}\in[-\frac{8.3 r^2}{R\sin\phistar},\frac{8.3r^2}{R\sin\phistar}]\\\eqref{eq:casebghd0lengthrestrictions1}: g\in(-\frac{8.3 r^2}{R\sin\phistar},\frac{8.3r^2}{R\sin\phistar})}}}\qquad(-\frac{16.6r^2}{R\sin{\phistar}}, \frac{16.6r^2}{R\sin{\phistar}}).
\end{aligned}
\end{equation} So, the conclusions of \eqref{eq:casebcase2minorterm3} consolidate into the following.
\begin{equation}\label{eq:casebcase2minorterm4}
\begin{aligned}
2\frac{d_0-\tilde{d_0}}{d_1}(2-\frac{2d_1}{\tilde{d_0}})k(\frac{2d_1}{\tilde{d_0}}-\frac{1}{k})&\le2\frac{|d_0-\tilde{d_0}|}{d_1}\bigm|(2-\frac{2d_1}{\tilde{d_0}})\bigm|\bigm|k(\frac{2d_1}{\tilde{d_0}}-\frac{1}{k})\bigm|<2\cdot\frac{\frac{16.6r^2}{R\sin{\phistar}}}{d_1}\cdot1.0003 \cdot0.751\\
&\overbracket{<}^{\mathllap{d_1>0.5r\sin\phistar}}2\cdot\frac{\frac{16.6r^2}{R\sin{\phistar}}}{0.5r\sin{\phistar}}\cdot1.0003\cdot0.751<\frac{50r}{R\sin^{2}{\phistar}}.
\end{aligned}
\end{equation}
This concludes \textbf{Case 2}, in which $d_1<\frac{21}{32}r\sin\phistar$, $k\in[1,\frac{4}{3}]$, $\mathcal{I}(k)$ in \eqref{eq:IkUpperBoundCasebCase2} satisfies the following. \[\mathcal{I}(k)\overbracket{<}^{\mathllap{\eqref{eq:Ikfirst3items},\eqref{eq:IkUpperBoundCasebCase2MinorTerms}}}-0.0639+\frac{63r}{R\sin^2\phistar}+\frac{164r}{R\sin^2\phistar}+\frac{50r}{R\sin^2\phistar}=-0.0639+\frac{277r}{R\sin^2\phistar}\overbracket<^{\mathclap{\eqref{eqJZHypCond}:R>\frac{30000r}{\sin^2\phistar}}}-0.05.\qedhere\]
\end{proof}
\NEWPAGE\section{Contraction control in case (c): more than one collision on \texorpdfstring{$\Gamma_R$}{flatter boundary}}\label{sec:ProofTECase_c}
This section is devoted to the proof of \cref{TEcase_c}, which involves the arduous work on the hardest case of those defined in \eqref{eqMainCases}. 

The first subsection defines functions which will be involved in lower bounds for the expansion in subsequent collisions, studies the return orbit segment involved in this case and basic formulas pertinent for the expansion computations, and begins estimates for the parameters involved. Later subsections then proceed according to the number $n_1$ of subsequent collisions on $\Gamma_R$  (\cref{def:MhatReturnOrbitSegment}) being at least 4, and equal to 3, 2, and 1, respectively. This becomes increasingly harder, and the last subsection alone takes up pages \pageref{SBSn1=1}--\pageref{sec:SACTWFELBWB}.

\begin{remark}[Standing assumptions]\label{REMCaseCAssumptions}
This section has as standing assumptions those of \cref{TEcase_c}: Fix $r,R$ satisfying the hyperbolicity condition \eqref{eqJZHypCond}, fix $\phistar\in(0,\tan^{-1}{(1/3)})$, $\hat F$ from \cref{def:SectionSetsDefs}, and the ``half-quadrant'' cone family $C_x\dfn\big\{(d\phi,d\theta)\bigm|\frac{d\theta}{d\phi}\in [0,1]\big\}$ from \cref{RMASUH}. We assume that $d_1<2r$ and that the orbit segment of an $x\in\hat M\aeq(\Mrin\cap\Mrout)\sqcup\mathcal{F}^{-1}(\Mrout\setminus\Mrin)$\COMMENT{\color{red}Do we need ``nonsingular''? I think not because singular ones can't return. They end in a corner.} from \eqref{def:Mhat} returning to $\hat M$ includes an orbit segment \eqref{eqPtsFromeqMhatOrbSeg} (page \pageref{eqPtsFromeqMhatOrbSeg}) with more than one collision with $\Gamma_R$.\COMMENT{Now that this is here, edit references to assumptions elsewhere} 
\end{remark}
We will prove:
\begin{enumerate}
    \item
    the cone family $C_x$ is strictly invariant under $D\hat{F}$, i.e., $D\hat{F}(C_x)\subset\text{\big\{interior of }C_{\hat{F}(x)}\text{\big\}}$.
    \item
    $\frac{\|D\hat{F}_{x}(dx)\|_\p}{\|dx\|_\p}>1-1743\big(r/R\big)-15450\big(r/R\big)^2>0.9$, for all $dx\in C_x$.
    \item
    when $n_1\ge4$, $\frac{\|D\hat{F}(dx)\|_\p}{\|dx\|_\p}>\frac{0.96}{1.01}n_1^2+5.655n_1-8.628>29.1$, for all $dx\in C_x$ (\cref{subsec:expansionLarge}).
\end{enumerate}
\subsection[A Formula for Expansion Lower Bound in \texorpdfstring{\cref{TEcase_c}}{case (c) transition}]{A candidate lower bound for expansion and parameter estimates}

\begin{notation}[Formulas for expansion in case (c) of \eqref{eqMainCases}]\label{def:AformularForLowerBoundOfExpansionForn1largerthan0}
For $x_0\in\Mrout$, $x_1\in\MRin$, $x_2\in\Mrin$, $n_1$ from \cref{def:MhatReturnOrbitSegment} and the free path and chord length functions $\tau_0$, $\tau_1$, $d_0$, $d_1$, $d_2$ from \cref{def:fpclf}. Denote by $\mathcal{B}^{\pm}_0$ the after/before wave front curvatures at $x_0$, as in \eqref{eq:MirrorEquationExpansionFormula} and let
\begin{equation}\label{eq:AformularForLowerBoundOfExpansionForn1largerthan0}
\begin{aligned}
\mathcal{II}_1(\mathcal{B}^+_0,n_1)&\dfn\big(1+\tau_0\mathcal{B}^+_0\big)+\big(\tau_1-\frac{2}{3}d_2\big)\Big[\big(2n_1+1-\frac{(2n_1+2)\tau_0}{d_1}\big)\mathcal{B}^+_0-\frac{2n_1+2}{d_1}\Big]\\
\mathcal{II}_2(\mathcal{B}^+_0,n_1)&\dfn 2n_1\big[1+(\tau_0-d_1)\mathcal{B}^+_0\big]
\end{aligned}
\end{equation}
\end{notation}

\begin{lemma}[Orbit configuration and expansion in \p-metric]\label{lemma:AformularForLowerBoundOfExpansionForn1largerthan0}
    Suppose a nonsingular $x\in \hat{M}\aeq(\Mrin\cap\Mrout)\sqcup\mathcal{F}^{-1}(\Mrout\setminus\Mrin)$ (\cref{def:SectionSetsDefs}) has a return orbit segment $x$, $\mathcal{F}(x)$, $\cdots$, $\mathcal{F}^{\sigma(x)}(x)=\hat{F}(x)\in \hat{M}$ as in \cref{def:MhatReturnOrbitSegment} with $x_{0}\in\Mrout$, $x_1=\mathcal{F}(x_0)=(\phi_1,\theta_1)\in \MRin$, $x_2\in\Mrin$, $\sin\theta_1< 2r/R$, $n_1\ge1$ (i.e., case (c) of \eqref{eqMainCases}),\COMMENT{Everything that precedes can be replaced by reference to \cref{REMCaseCAssumptions}, I believe, Check and implement.} and that $dx=(d\phi,d\theta)\in\Big\{(d\phi,d\theta)\bigm|\frac{d\theta}{d\phi}\in\big[0,1\big]\Big\}$ is a tangent vector at $x$, i.e., $dx$ is in the half-quadrant cone in the $\phi\theta$-coordinates. Then we have:
    \begin{enumerate}
\item\label{itemAformularForLowerBoundOfExpansionForn1largerthan01} $x\in\mathcal{F}^{-1}(\Mrout\setminus\Mrin)\subset M_r$, $(\phi_0,\theta_0)=x_0=\mathcal{F}(x)\in\Mrout$, $x_1=\mathcal{F}(x_0)\in\MRin$. If we set
        $x_{1,R}\dfn\mathcal{F}(x_1)\in M_R$, $\cdots$, $ x_{n_1,R}\dfn\mathcal{F}^{n_1}(x_1)\in M_R$, then $\mathcal{F}(x_{n_1,R})=x_2\in\Mrin$, and we have $\frac{\|dx_0\|_\p}{\|dx\|_\p}\ge1$ and $\Bcal^+_0\in\big[-\frac{4}{3d_0},-\frac{1}{d_0}\big]$, where $d_0=r\sin\theta_0$ is from \cref{def:fpclf}.
        \item\label{itemAformularForLowerBoundOfExpansionForn1largerthan02} $dx_0\dfn D\mathcal{F}_x(dx)$ is a tangent vector at $x_0$, $dx_1\dfn D\mathcal{F}_{x_0}(dx_0)$ is a tangent vector at $x_1$, $dx_{n_1,R}\dfn
        D\mathcal{F}^{n_1}_{x_1}(dx_1)$ is a tangent vector at $x_{n_1,R}$, $dx_2\dfn D\mathcal{F}_{x_{n_1,R}}(dx_{n_1,R})$ is a tangent vector at $x_2$.
        \item\label{itemAformularForLowerBoundOfExpansionForn1largerthan03} $m(x_2)\ge3$ hence $x_3\dfn\mathcal{F}(x_2)\in M_r$. $dx_3\dfn D\mathcal{F}_{x_2}(dx_2)=D\mathcal{F}^{4+n_1}_x(dx)$ and  $\mathcal{F}^{m(x_2)-2}(x_3)=\mathcal{F}^{m(x_2)-1}(x_2)\overbracket{=}^{\mathllap{\text{\cref{def:SectionSetsDefs}}}}\hat{F}(x)\in\hat{M}$, i.e. $x_3=\mathcal{F}(x_2)$, $\cdots$, $\mathcal{F}^{m(x_2)-1}(x_2)$ are all in $M_r\setminus\Mrout$. We have $\frac{\|D\hat{F}_x(dx)\|_\p}{\|dx_3\|_\p}>1$.
        \item\label{itemAformularForLowerBoundOfExpansionForn1largerthan04} Computation of $\frac{\|dx_3\|_\p}{\|dx_0\|_\p}$. We use the notation from \eqref{eq:MirrorEquationExpansionFormula} with indices: denote by $\mathcal{B}^{\pm}_1$, $\mathcal{B}^{\pm}_{n_1,R}$, $\mathcal{B}^{\pm}_2$ and $\mathcal{B}^{\pm}_3$ the infinitesimal wave front curvatures after/before collisions at $x_1$, $x_{n_1,R}$, $x_2$ and $x_3$, respectively, $\tau_{1,\text{\upshape im}}\dfn-2n_1d_1$ and $\Bcal^+_{1,\text{\upshape im}}\dfn\Bcal^-_1-\frac{1}{d_1}$ (``\emph{im}'' means imaginary billiard trajectory in \cite[Figure 8.8, Equation (8.12)]{cb}).\COMMENT{Is there a figure in this article that we can refer to?} Then 
        \begin{equation}\label{eq:Expansionx3x0}\begin{aligned}E(n_1,\Bcal^+_0)\dfn\frac{\|dx_3\|_{\p}}{\|dx_0\|_{\p}}&=\bigm|1+\tau_0\Bcal^+_0\bigm|\bigm|1-2n_1d_1(\Bcal^-_1-\frac{1}{d_1})\bigm|\bigm|1+\tau_1\Bcal^{+}_{n_1,R}\bigm|\bigm|1+2d_2\Bcal^{+}_{2}\bigm|\\
        &=\bigm|1+\tau_0\Bcal^+_0\bigm|\bigm|1+\tau_{1,\text{\upshape im}}\Bcal^+_{1,\text{\upshape im}}\bigm|\bigm|1+\tau_1\Bcal^{+}_{n_1,R}\bigm|\bigm|1+2d_2\Bcal^{+}_{2}\bigm|\\&
        %
        =3\big|\mathcal{II}_1(B^+_0,n_1)+\mathcal{II}_2(B^+_0,n_1)\big|,
        \end{aligned}\end{equation} where  $\mathcal{II}_1(\mathcal{B}^+_0,n_1)$, $\mathcal{II}_2(\mathcal{B}^+_0,n_1)$ and $\mathcal{B}^+_0$ are from \cref{def:AformularForLowerBoundOfExpansionForn1largerthan0}.
        \item\label{itemAformularForLowerBoundOfExpansionForn1largerthan05} 
         With $d_0$, $d_1$, $d_2$, $\tau_0$, $\tau_1$, $n_1$ as in \cref{def:fpclf,def:AformularForLowerBoundOfExpansionForn1largerthan0}, we have\begin{equation}\label{eq:B3minusEquation}
        \begin{multlined}[.85\textwidth]
            \mathcal{B}^{-}_3=\frac{-2}{d_2(1+2d_2\Bcal^+_2)}-\frac{1}{d_1(1+2d_2\Bcal_2^+)(1+\tau_1\Bcal^{+}_{n_1,R})}-\frac{1}{d_1(1+2d_2\Bcal^+_2)(1+\tau_1\Bcal^+_{n_1,R})(1+\tau_{1,\text{\upshape im}}\Bcal^+_{1,\text{\upshape im}})}\\
            +\frac{\Bcal^+_0}{(1+2d_2\Bcal^+_2)(1+\tau_1\Bcal^+_{n_1,R})(1+\tau_{1,\text{\upshape im}}\Bcal^+_{1,\text{\upshape im}})(1+\tau_0\Bcal_0^+)}.
        \end{multlined}\end{equation} 
    \end{enumerate}
\end{lemma}
\begin{proof}
\strut
   \hidesubsection{Proof of \eqref{itemAformularForLowerBoundOfExpansionForn1largerthan01}} For $x\in\hat{M}=\mathcal{F}^{-1}(\Mrout\setminus\Mrin)\sqcup(\Mrout\cap\Mrin)$ in \eqref{eqMhatOrbSeg} case (c), we first show that $x\notin \Mrout\cap\Mrin$. If we had  $x\in\Mrout\cap\Mrin$, then $x_0=x$ in \cref{def:MhatReturnOrbitSegment}. Since $x_0$ is nonsingular, we have $x_0=x\in\textcolor{red}{M^{\text{\upshape{out}}}_{r,0}}$. But by \cref{contraction_region} we have $x_0\in\Nout$ and \cref{Prop:Udisjoint} shows that $\Nout\cap\textcolor{red}{M^{\text{\upshape{out}}}_{r,0}}=\emptyset$. This means that it is impossible to have $x\in\Mrout\cap\Mrin$. Thus, $x\in\mathcal{F}^{-1}(\Mrout\setminus\Mrin)$, $x_0=\mathcal{F}(x)\in\Mrout$. 

    Hence, $dx_0\dfn D\mathcal{F}_{x}(dx)$ is a tangent vector at $x_0$. Since $x$, $x_0$ are nonsingular on $M_r$ and $dx=(d\phi,d\theta)$ is a tangent vector in the quadrant cone, by \cref{proposition:ShearOneStepExpansion} we have $\frac{\|dx_0\|_\p}{\|dx\|_\p}\ge1$. 
    
    Since $\frac{d\theta}{d\phi}\in[0,1]$, $(d\phi_0,d\theta_0)=dx_0=D\mathcal{F}_x(dx)$, $D\mathcal{F}_x=\begin{pmatrix}
        1 & 2\\
        0 & 1
    \end{pmatrix}$, that is $\begin{pmatrix} d\phi_0\\d\theta_0
    \end{pmatrix}=\begin{pmatrix}
        1 & 2\\
        0 & 1
    \end{pmatrix}\begin{pmatrix} d\phi\\d\theta
    \end{pmatrix}$, we have $\frac{d\theta_0}{d\phi_0}=\frac{d\theta}{d\phi+2d\theta}=\frac{1}{2+\frac{1}{\frac{d\theta}{d\phi}}}\in\big[0,\frac{1}{3}\big]$. Hence, by \eqref{eq:coordinateschange}, \eqref{eq:MirrorEquationExpansionFormula}, we have $\mathcal{B}^+_0\cos\varphi_0+\frac{1}{r}=\frac{d\varphi_0}{ds_0}=\frac{-1}{r}\frac{d\theta_0}{d\phi_0}\in[-\frac{1}{3r},0]$.
    
    Therefore, $\Bcal^+_0=-\frac{1}{r\cos\varphi_0}\big(1+\frac{d\theta_0}{d\phi_0}\big)\overbracket{\in}^{\mathclap{\substack{\frac{d\theta_0}{d\phi_0}\in[0,1/3]\text{, \eqref{eq:coordinateschange}, and \cref{def:fpclf}:}\\d_0=r\cos\varphi_0=r\sin\theta_0}}}\big[-\frac{4}{3d_0},-\frac{1}{d_0}\big]$. This proves item \eqref{itemAformularForLowerBoundOfExpansionForn1largerthan01}.\COMMENT{Fix reference above}
    
    \hidesubsection{Proof of \eqref{itemAformularForLowerBoundOfExpansionForn1largerthan02}}
    We also have $x_1=\mathcal{F}(x_0)\in\MRin$, $dx_1\dfn D\mathcal{F}_{x_0}(dx_0)$. By \cref{def:MhatReturnOrbitSegment}, $ x_{n_1,R}\dfn\mathcal{F}^{n_1}(x_1)\in M_R$, $\mathcal{F}(x_{n_1,R})=x_2\in\Mrin$ is a sub-segment of \eqref{eqMhatOrbSeg} so that we can define $dx_{1,n_1}\dfn
    D\mathcal{F}^{n_1}_{x_1}(dx_1)$, $dx_2\dfn D\mathcal{F}_{x_{1,n_1}}(dx_{1,n_1})$, where $dx_{n_1,R}$ is a tangent vector at $x_{n_1,R}$, $dx_2$ is a tangent vector at $x_2$.
    This proves item \eqref{itemAformularForLowerBoundOfExpansionForn1largerthan02}.
    
    \hidesubsection{Proof of \eqref{itemAformularForLowerBoundOfExpansionForn1largerthan03}} By \cref{contraction_region}, we have $x_2\in\Nin$ and by \cref{Prop:Udisjoint}, we have $\emptyset=\Nin\cap(\textcolor{blue}{M^{\text{\upshape{in}}}_{r,0}}\cup\textcolor{blue}{M^{\text{\upshape{in}}}_{r,1}}\cup\textcolor{blue}{M^{\text{\upshape{in}}}_{r,2}})$. Hence, $x_2$ is nonsingular and $m(x_2)\ge3$. By the definition of $m(x_2)$ and $\hat{M}$ in \eqref{def:Mhat}, $(\phi_3,\theta_3)=x_3\dfn\mathcal{F}(x_2),\mathcal{F}(x_3)$, $\cdots$, $\mathcal{F}^{m(x_2)-2}(x_3)=\mathcal{F}^{m(x_2)-1}(x_2)\in\hat{M}$ are nonsingular points on $M_r\setminus\Mrout$ since their $\mathcal{F}$-images are not on $M_R$. Now $dx_3\dfn D\mathcal{F}_{x_2}(dx_2)=D\mathcal{F}^{4+n_1}_x(dx)$ and  $\mathcal{F}^{m(x_2)-2}(x_3)=\mathcal{F}^{m(x_2)-1}(x_2)\overset{\text{\cref{def:SectionSetsDefs}}}{\scalebox{7}[1]{$=$}}\hat{F}(x)\in\hat{M}$. 
    
    By \cref{caseB}, $D\mathcal{F}^{4+n_1}_x$ is a matrix with positive entries, since $(d\phi_3,d\theta_3)=dx_3=D\mathcal{F}^4_x(dx)$ and $dx=(d\phi,d\theta)$ is in the quadrant with $\frac{d\theta}{d\phi}\ge0$, $(d\phi_3,d\theta_3)=dx_3$ is then in the interior of the quadrant, that is, $\frac{d\theta_3}{d\phi_3}>0$. Since $D\hat{F}_x(dx)=D\mathcal{F}^{m(x_2)-2}(dx_3)$ and $x_3,\mathcal{F}(x_3),\cdots,\mathcal{F}^{m(x_2)-2}(x_3)$ are all in $M_r$, by \cref{proposition:ShearOneStepExpansion} we have $\frac{\|D\hat{F}_x(dx)\|_\p}{\|dx_3\|_\p}=\frac{\|D\mathcal{F}^{m(x_2)-2}(dx_3)\|_\p}{\|dx_3\|_\p}>1$. This proves item \eqref{itemAformularForLowerBoundOfExpansionForn1largerthan03}.
    
    \hidesubsection{Proof of \eqref{itemAformularForLowerBoundOfExpansionForn1largerthan04}}
    By applying \cite[exercises 8.28, 8.29]{cb} to $\Gamma_R$, we get\COMMENT{\textcolor{red}{\boris All the \textbackslash text and \textbackslash upshape commands must be removed and their content taken out of math mode. Only formulas should be in math mode.} Some fixing is needed, and that is too hard until this is cleaned up.}
    \begin{enumerate}[label={(\roman*)}]\label{Property:WavefrontsProofOf4}
        \item\label{item01ContinuedfractionResults} $\frac{\|dx_2\|_\p}{\|dx_{n_1,R}\|_\p} \overset{\eqref{eq:MirrorEquationExpansionFormula}}=\big|1+\tau_1\Bcal^{+}_{n_1,R}\big|$, since $\tau_1$ is the distance between  $p(x_{n_1,R})=p(\Fcal^{-1}(x_2))$ and $p(x_2)$,\COMMENT{I fixed this line in case you were unclear on what is needed; all others on various pages should be edited in the same way.}
        \item\label{item02ContinuedfractionResults} $\frac{\|dx_3\|_\p}{\|dx_2\|_\p}\overset{\eqref{eq:MirrorEquationExpansionFormula}}=\big|1+2d_2\Bcal^+_2\big|\text{ , since }2d_2 \text{ is the distance between } p(x_2)=p(\Fcal^{-1}(x_3))\text{ and }p(x_3)$,
        \item\label{item03ContinuedfractionResults} $\frac{1}{\Bcal^{+}_{n_1,R}}+\tau_1 \overset{\text{\eqref{eqBBr}, \cite[equation (3.31)]{cb}}}{\scalebox{16}[1]{$=$}}\frac{1}{\Bcal^{-}_2} \text{ , thus } \Bcal^{-}_{2}=\frac{\Bcal^{+}_{n_1,R}}{1+\tau_1
        \Bcal^{+}_{n_1,R}}$,
        \item\label{item04ContinuedfractionResults} $\Bcal^{+}_2\overset{\eqref{eq:MirrorEquationExpansionFormula} \text{ with } \mathcal{K}=\frac{-1}{r}}{\scalebox{11}[1]{$=$}} \Bcal_2^{-}-\frac{2}{r\cos{\varphi_2}}=\Bcal^{-}_2-\frac{2}{d_2}\text{ , since }\eqref{eq:coordinateschange}\text{: }\theta_2+\varphi_2=\pi/2 \text{ and \cref{def:fpclf}: }d_2=r\sin\theta_2$,
        \item\label{item05ContinuedfractionResults} $\Bcal^{+}_{n_1,R} \overset{\text{\cite[exercise 8.28]{cb}}}{\scalebox{12}[1]{$=$}}-\frac{1}{d_1}+\frac{1}{-2n_1d_1+\frac{1}{-\frac{1}{d_1}+\Bcal^{-}_1}}  = -\frac{1}{d_1}+\frac{\Bcal^{-}_1-\frac{1}{d_1}}{2n_1+1-2n_1d_1\Bcal^{-}_1}$,
        \item\label{item06ContinuedfractionResults} $\frac{\|dx_{n_1,R}\|_\p}{\|dx_1\|_\p}\overset{\text{\cite[exercise 8.29]{cb}}}{\scalebox{12}[1]{$=$}}\big|1-2n_1d_1(\Bcal^{-}_1+\frac{-1}{d_1})\big| =\big|1+2n_1-2n_1d_1\Bcal^{-}_1\big|$,
        \item\label{item07ContinuedfractionResults} $\frac{1}{\Bcal_{1}^{-}}\overset{\text{\eqref{eqBBr}, \cite[equation (3.31)]{cb}}}{\scalebox{16}[1]{$=$}}\frac{1}{\Bcal_{0}^{+}}+\tau_0\text{ , thus }\Bcal^-_1=\frac{\Bcal^+_0}{1+\tau_0\Bcal^+_0}$,
        \item\label{item08ContinuedfractionResults} $\frac{\|dx_1\|_\p}{\|dx_0\|_\p} \overset{\eqref{eq:MirrorEquationExpansionFormula}}=\big|1+\tau_0\Bcal^{+}_0\big|\text{ , since }\tau_0 \text{ is the distance between } p(x_0)\text{ and }p(x_1)$.
        \end{enumerate}
        Therefore, \begin{align*}E(n_1,\Bcal^+_0)\dfn&\frac{\|dx_3\|_{\p}}{\|dx_0\|_{\p}}=\frac{\|dx_3\|_{\p}}{\|dx_2\|_{\p}}\frac{\|dx_2\|_{\p}}{\|dx_{n_1,R}\|_{\p}}\frac{\|dx_{n_1,R}\|_{\p}}{\|dx_1\|_{\p}}\frac{\|dx_1\|_{\p}}{\|dx_0\|_{\p}}\\
        \overbracket{=}^{\quad\mathllap{\text{\ref{item01ContinuedfractionResults}\ref{item02ContinuedfractionResults}\ref{item04ContinuedfractionResults}\ref{item07ContinuedfractionResults}}}}&\bigm|1+2d_2\Bcal^+_2\bigm|\bigm|1+\tau_1\Bcal^{+}_{n_1,R}\bigm|\bigm|1-2n_1d_1(\Bcal^-_1-\frac{1}{d_1})\bigm|\bigm|1+\tau_0\Bcal^+_0\bigm|.
        \end{align*}        \begin{align*}\mathllap{\text{Next,\qquad}}
        \frac{\|dx_3\|_\p}{\|dx_{n_1,R}\|_\p}&=\frac{\|dx_2\|_\p}{\|dx_{n_1,R}\|_\p}\frac{\|dx_3\|_\p}{\|dx_2\|_\p}\overset{\text{(i) and (ii)}}{\scalebox{7}[1]{$=$}}\big|1+\tau_1\Bcal^{+}_{n_1,R}\big|\big|1+2d_2\Bcal^+_2\big|\overbracket{=}^{\mathclap{\text{(iv)}}}\big|1+\tau_1\Bcal^{+}_{n_1,R}\big|\big|1+2d_2(\Bcal^-_2-\frac{2}{d_2})\big|\\
        &\overset{\text{(iii)}}{=}\big|1+\tau_1\Bcal^{+}_{n_1,R}\big|\big|1+2d_2\big(\frac{\Bcal^{+}_{n_1,R}}{1+\tau_1\Bcal^{+}_{n_1,R}}-\frac{2}{d_2}\big)\big|=\big|1+\tau_1 \Bcal^{+}_{n_1,R}\big|\big|3-2d_2\frac{\Bcal_{n_1,R}^{+}}{1+\tau_1\Bcal^{+}_{n_1,R}}\big|\\
        &=\big|3+(3\tau_1-2d_2)\Bcal^{+}_{n_1,R}\big|\overset{\text{(v)}}{=}3\big|1+(\tau_1-\frac{2}{3}d_2)(\frac{-1}{d_1}+\frac{\Bcal^-_1-\frac{1}{d_1}}{2n_1+1-2n_1d_1\Bcal^{-}_1})\big|.
    \end{align*} 
    \begin{align*}\mathllap{\text{Therefore, }}
    \frac{\|dx_3\|_\p}{\|dx_1\|_\p}&=\frac{\|dx_3\|_\p}{\|dx_{n_1,R}\|_\p}\frac{\|dx_{n_1,R}\|_\p}{\|dx_1\|_\p}\overset{\text{(vi)}}=3\big|1+(\tau_1-\frac{2d_2}{3})(\frac{-1}{d_1}+\frac{\Bcal^-_1-\frac{1}{d_1}}{2n_1+1-2n_1d_1\Bcal^{-}_1})\big|\big|1+2n_1-2n_1d_1\Bcal^{-}_1\big|\\
    &=3\big|1+2n_1-2n_1d_1\Bcal^{-}_1+(\tau_1-\frac{2}{3}d_2)(\Bcal^{-}_1-\frac{1}{d_1}-\frac{2n_1+1}{d_1}+2n_1\Bcal^{-}_1)\big|\\
    &=3\Big|1+2n_1-2n_1d_1\Bcal^{-}_1+\big(\tau_1-\frac{2}{3}d_2\big)\big[(2n_1+1)\Bcal^-_1-\frac{2n_1+2}{d_1}\big]\Big|\\
    &\overbracket{=}^{\mathclap{{\text{(vii)}}}}3\Big|1+2n_1-2n_1d_1\frac{\Bcal^+_0}{1+\tau_0\Bcal^+_0}+\big(\tau_1-\frac{2}{3}d_2\big)\big[(2n_1+1)\frac{\Bcal^+_0}{1+\tau_0\Bcal^+_0}-\frac{2n_1+2}{d_1}\big]\Big|\text{. }
    \end{align*}
\begin{equation}\label{eq:caseC0to3Expansion}
    \begin{aligned}\mathllap{\text{Hence \eqref{itemAformularForLowerBoundOfExpansionForn1largerthan04}:\ }}
    \frac{\|dx_3\|_\p}{\|dx_0\|_\p}&=\frac{\|dx_1\|_\p}{\|dx_0\|_\p}\frac{\|dx_3\|_\p}{\|dx_1\|_\p}\\
    &\overbracket{=}^{\mathclap{\text{(viii)}}}3\Big|1+\tau_0\Bcal^+_0\Big|\Big|1+2n_1-2n_1d_1\frac{\Bcal^+_0}{1+\tau_0\Bcal^+_0}+\big(\tau_1-\frac{2}{3}d_2\big)\big[(2n_1+1)\frac{\Bcal^+_0}{1+\tau_0\Bcal^+_0}-\frac{2n_1+2}{d_1}\big]\Big|\\
    & =3\Big|(1+\tau_0\Bcal^+_0)+2n_1(1+\tau_0\Bcal^+_0)-2n_1d_1\Bcal^+_0+\big(\tau_1-\frac{2}{3}d_2\big)\big[(2n_1+1)\Bcal^+_0-\frac{2n_1+2}{d_1}(1+\tau_0\Bcal^+_0)\big]\Big|\\
    & =3\Big|\underbracket{(1+\tau_0\Bcal^+_0)+(\tau_1-\frac{2}{3}d_2)\big[[(2n_1+1)-\frac{(2n_1+2)\tau_0}{d_1}]\Bcal^+_0-\frac{2n_1+2}{d_1}\big]}_{=\mathcal{II}_1(\Bcal^+_0,n_1)}+\underbracket{2n_1\big[1+(\tau_0-d_1)\Bcal^+_0\big]}_{=\mathcal{II}_2(\Bcal^+_0,n_1)}\Big|.
    \end{aligned}
\end{equation}

\hidesubsection{Proof of \eqref{itemAformularForLowerBoundOfExpansionForn1largerthan05}}\strut\COMMENT{Note the overbracket fixes below; implement further, \color{red}such as, adding mathclap}\COMMENT{Does ``\eqref{eqBBr}, \cite[equation (3.31)]{cb}'' below mean that both are needed? I suspect only \eqref{eqBBr} is. If correct, then delete ``\cite[equation (3.31)]{cb}''!}\COMMENT{Does ``since $2d_2$ is the distance between\dots'' explain the ``='' before or after? Both of those have added explanations on top anyway, so this seems confusing.}
\begin{align*}
\Bcal^{-}_3\overbracket{=}^{\Bigstrut\mathclap{\text{\eqref{eqBBr}, \cite[equation (3.31)]{cb}}}}&\frac{\Bcal^+_2}{1+2d_2\Bcal^+_2}\text{ , since }2d_2 \text{ is the distance between } p(x_2)=p(\Fcal^{-1}(x_3))\text{ and }p(x_3),\\
\overbracket{=}^{\bigstrut\mathclap{\text{\eqref{eq:MirrorEquationExpansionFormula} with } \mathcal{K}={-1}/{r},\ d_2=r\cos{\varphi_2}}}&\frac{-2/d_2}{1+2d_2\Bcal^+_2}+\frac{\Bcal^-_{2}}{1+2d_2\Bcal_2^+},\\
\overbracket{=}^{\bigstrut\mathclap{\text{\eqref{eqBBr}, \cite[equation (3.31)]{cb}}}}&\frac{-2/d_2}{1+2d_2\Bcal^+_2}+\frac{\overbracket{\Bcal^+_{n_1,R}}^{\mathclap{\text{\ref{item05ContinuedfractionResults} in proof of (4)}}}}{(1+2d_2\Bcal_2^+)(1+\tau_1\Bcal^+_{n_1,R})},\\
\overbracket{=}^{\mathclap{\text{\ref{item05ContinuedfractionResults} in proof of (4)}}}&\qquad\qquad\frac{-2/d_2}{1+2d_2\Bcal^+_2}+\frac{-\frac{1}{d_1}+\frac{\Bcal^{-}_1-\frac{1}{d_1}}{2n_1+1-2n_1d_1\Bcal^{-}_1}}{(1+2d_2\Bcal_2^+)(1+\tau_1\Bcal^+_{n_1,R})},\\
=&\frac{-2}{d_2(1+2d_2\Bcal^+_2)}-\frac{1}{d_1(1+2d_2\Bcal_2^+)(1+\tau_1\Bcal^+_{n_1,R})}-\frac{1}{d_1(1+2d_2\Bcal_2^+)(1+\tau_1\Bcal^+_{n_1,R})[1-2n_1d_1(\Bcal^-_1-\frac{1}{d_1})]}\\
&+\frac{\Bcal^-_1}{(1+2d_2\Bcal_2^+)(1+\tau_1\Bcal^+_{n_1,R})[1-2n_1d_1(\Bcal^-_1-\frac{1}{d_1})]}\\
\overbracket{=}^{\mathclap{\text{\eqref{eqBBr}, \cite[equation (3.31)]{cb}}}}&\qquad\qquad\frac{-2}{d_2(1+2d_2\Bcal^+_2)}-\frac{1}{d_1(1+2d_2\Bcal_2^+)(1+\tau_1\Bcal^+_{n_1,R})}-\frac{1}{d_1(1+2d_2\Bcal_2^+)(1+\tau_1\Bcal^+_{n_1,R})[1-2n_1d_1(\Bcal^-_1-\frac{1}{d_1})]}\\
&+\frac{\Bcal^+_0}{(1+2d_2\Bcal_2^+)(1+\tau_1\Bcal^+_{n_1,R})[1-2n_1d_1(\Bcal^-_1-\frac{1}{d_1})](1+\tau_0\Bcal^+_0)}.
\end{align*} This gives \eqref{eq:B3minusEquation}.
\end{proof}



\begin{figure}[h]
\begin{center}
    \begin{tikzpicture}[xscale=0.5,yscale=.5]
        \pgfmathsetmacro{\PHISTAR}{1.5};
        \pgfmathsetmacro{\Dlt}{0.6};
        \pgfmathsetmacro{\Dltplus}{\Dlt+0.6}
        \pgfmathsetmacro{\WIDTH}{8};
        \tkzDefPoint(-\PHISTAR,0){LL};
        \tkzDefPoint(-\PHISTAR,\WIDTH){LU};
        \tkzDefPoint(\PHISTAR,0){RL};
        \tkzDefPoint(\PHISTAR,\WIDTH){RU};
        \tkzDefPoint(-\PHISTAR,\PHISTAR){X};
        \tkzDefPoint(\PHISTAR,\WIDTH-\PHISTAR){Y};
        \tkzDefPoint(-\PHISTAR,\WIDTH-\PHISTAR){IX};
        \tkzDefPoint(\PHISTAR,\PHISTAR){IY};
        \tkzDefPoint(-\PHISTAR,\PHISTAR+\Dlt){LLD};
        \tkzDefPoint(\PHISTAR,\PHISTAR+\Dlt){RLD};
        \tkzDefPoint(-\PHISTAR,\PHISTAR+\Dlt){LLD};
        \tkzDefPoint(\PHISTAR,\PHISTAR+\Dlt){RLD};
        \tkzDefPoint(-\PHISTAR,\WIDTH-\PHISTAR-\Dlt){LUD};
        \tkzDefPoint(\PHISTAR,\WIDTH-\PHISTAR-\Dlt){RUD};
        \tkzDefPoint(-\PHISTAR,\PHISTAR+\Dltplus){LLDp};
        \tkzDefPoint(\PHISTAR,\PHISTAR+\Dltplus){RLDp};
        \tkzDefPoint(-\PHISTAR,\WIDTH-\PHISTAR-\Dltplus){LUDp};
        \tkzDefPoint(\PHISTAR,\WIDTH-\PHISTAR-\Dltplus){RUDp};
        \tkzDefPoint(-\PHISTAR,\PHISTAR+\Dltplus+0.8){LLDdp};
        \tkzDefPoint(\PHISTAR,\PHISTAR+\Dltplus+0.8){RLDdp};
        \tkzDefPoint(-\PHISTAR,\WIDTH-\PHISTAR-\Dltplus-0.8){LUDdp};
        \tkzDefPoint(\PHISTAR,\WIDTH-\PHISTAR-\Dltplus-0.8){RUDdp};
        \tkzDefPoint (0.3*\PHISTAR,\PHISTAR+\Dltplus){montonestart};
        \tkzDefPoint (0.3*\PHISTAR,0.333*\WIDTH){montoneend};
        \draw [thin] (LL) --(LU);
        \draw [thin] (LU) --(RU);
        \draw [thin] (RL) --(RU);
        \draw [thin] (LL) --(RL);
        \draw [blue,ultra thin,dashed] (LL) --(IY);
        \draw [red,ultra thin,dashed]  (X)--(RL);
        \draw [blue,ultra thin,dashed] (RU) --(IX);
        \draw [red,ultra thin,dashed]  (LU)--(Y);
        \tkzLabelPoint[right](RU){$\pi$};
        \tkzLabelPoint[right](RL){$0$};
        \tkzLabelPoint[right](IY){$\theta_1=\Phistar,d_1=r\sin{\phistar}$};
        \tkzLabelPoint[right](Y){$\theta_1=\pi-\Phistar,d_1=r\sin{\phistar}$};
        \coordinate (RM) at ($(RL)!0.5!(RU)$);
        \tkzLabelPoint[right](RM){${\pi}/{2}$};
        \coordinate (UpperM) at ($(LU)!0.5!(RU)$);
        \coordinate (LowerM) at ($(LL)!0.5!(RL)$);
        \coordinate (LM) at ($(LU)!0.5!(LL)$);
        \tkzLabelPoint[below](LL){$-\Phistar$};
        \tkzLabelPoint[below](RL){$+\Phistar$};
        \tkzLabelPoint[below](LowerM){$0$};
        \tkzDefPoint(0,0.5*\PHISTAR){LMend};
        \tkzDefPoint(0,\WIDTH-0.5*\PHISTAR){HMend};
        \fill [green, opacity=10/30] (LL) -- (LMend) -- (X) -- cycle;
        \fill [green, opacity=10/30] (RU) -- (HMend) -- (Y) -- cycle;
        \tkzDrawPoints(X,Y,IX,IY,LMend,HMend);
        \node[] at (0,0.5*\WIDTH){\Large $n_1=0$};
        \node[] at (-0.7*\PHISTAR,0.5*\PHISTAR){\small $n_1\ge1$};
        \node[] at (0.7*\PHISTAR,\WIDTH-0.5*\PHISTAR){\small $n_1\ge1$};
    \end{tikzpicture}
    \caption{Two components of \MRin\ where $x_1\in\MRin$ with $n_1+1\ge2$}\label{fig:x1_multi_collisions_caseC}
\end{center}
\end{figure}
\begin{proposition}[Case (c) 
parameter estimates
]\label{proposition:casecparametersestimate} In the context of \cref{TEcase_c},\COMMENT{Edit if we state a standing assumption} for the length functions $\tau_0$, $\tau_1$, $d_0$, $d_1$, $d_2$ from \cref{def:AformularForLowerBoundOfExpansionForn1largerthan0}, under condition \eqref{eqJZHypCond}, we have the following parameter estimates.

If $n_1\ge1$, then\COMMENT{Eliminate all the \textbackslash text and \textbackslash upshape below. Text should not be part of inline mathematics.\Hrule Ideally replace all $\Y<$ by a reason, but this is less urgent.}
\begin{enumerate}[label={(\roman*)}]
        \item\label{itemcasecparametersestimate1} $\big|d_0-r\sin\phistar\big|\overbracket{<}^{\mathllap{\text{\cite[Inequality (3.20)]{hoalb}}}}\frac{17r^2\sin\phistar}{4R}$ and $\big|d_2-r\sin\phistar\big|\overbracket{<}^{\mathllap{\text{\cite[Inequality (3.20)]{hoalb}}}}\frac{17r^2\sin\phistar}{4R},$ (also see \cite[Inequality (5.11)]{hoalb}),
        \item\label{itemcasecparametersestimate2} $2r\sin{\phistar}\overbracket{<}^{\mathclap{\text{\cite[Inequalities (3.16) and (4.7)]{hoalb}}}}\tau_0+\tau_1+2n_1d_1\overbracket{<}^{\mathrlap{\text{\cite[Inequalities (3.16) and (4.7)]{hoalb}}}}2r\sin{\phistar}+\frac{8.1r^2\sin{\phistar}}{R},$ (also see \cite[Inequality (5.12)]{hoalb}),
        \item\label{itemcasecparametersestimate3} $d_1\overbracket{>}^{\mathclap{\substack{\text{\upshape{By (ii) and \cref{remark:deffpclf}: }}\\\tau_1<2d_1,\tau_0<2d_1}}}\frac{1}{n_1+2}r\sin{\phistar}$,
        \item\label{itemcasecparametersestimate4} $d_1<\frac{r\sin{\phistar}}{n_1}+\frac{8.1r}{2n_1R}r\sin{\phistar}$ by \ref{itemcasecparametersestimate2}: $2n_1d_1<\tau_0+\tau_1+2n_1d_1<2r\sin{\phistar}+\frac{8.1r^2\sin{\phistar}}{R},$
        \item\label{itemcasecparametersestimate5} $\frac{r\sin\phistar}{d_0}\in(0.9975,1.0026)\text{ \upshape{and} } \frac{r\sin\phistar}{d_2}\in(0.9975,1.0026)$ by \ref{itemcasecparametersestimate1} and \eqref{eqJZHypCond}: $R>1700r$.\COMMENT{Edit! ("divide 9i0 by $r\sin\phistar$")\wentao{explained in remark now.}}
        \end{enumerate}
        
        If $n_1\ge2$, then
        \begin{enumerate}[label={(\roman*)}]
        \setcounter{enumi}{5}
        \item\label{itemcasecparametersestimate6} $\tau_0<\frac{2}{3}d_0+\frac{2}{3}d_1$ by \cite[Proposition 5.2]{hoalb},
        \item\label{itemcasecparametersestimate7} $\frac{2}{3}d_2+\frac{2}{3}d_1-\tau_1>\frac{1}{6}r\sin{\phistar}-\frac{16.6r}{3R}r\sin{\phistar}$ and $\frac{2}{3}d_0+\frac{2}{3}d_1-\tau_0>\frac{1}{6}r\sin{\phistar}-\frac{16.6r}{3R}r\sin{\phistar}$ \COMMENT{\textcolor{red}{Needs editing, but first remove all the \textbackslash text and the like.\wentao{clarified in remark now}}} by \cite[(5.15)]{hoalb} with its following statement.
    \end{enumerate} 
\end{proposition}
\begin{remark}

\textbf{\cref{proposition:casecparametersestimate}\ref{itemcasecparametersestimate1}} with $n_1\ge2$ is the claim of \cite[Inequalities (5.11)]{hoalb}. Its proof is based on the Mean-Value Theorem and \cite[Inequalities (3.20)]{hoalb}. This proof argument also works for $n_1=1$. 

\textbf{\cref{proposition:casecparametersestimate}\ref{itemcasecparametersestimate2}} with $n_1\ge2$ is the claim of \cite[(5.12)]{hoalb}. Its proof relies on \cite[Inequalities (3.16)(4.7)]{hoalb}. This proof argument works the same for $n_1=1$. 

In \cite{hoalb}, we also note that the statement following \cite[equation (5.19)]{hoalb} regards \cite[Inequalities (5.11)(5.12)]{hoalb} as established for $n_1=1$ to argue that $\tau_0\ge\frac{2}{3}d_0+\frac{1}{2}d_1$ and $\tau_1\ge\frac{1}{2}d_1+\frac{2}{3}d_2$ cannot simultaneously hold for $n_1=1$.

\textbf{\cref{proposition:casecparametersestimate}\ref{itemcasecparametersestimate5}} follows \textbf{\cref{proposition:casecparametersestimate}\ref{itemcasecparametersestimate1}}. Dividing by $r\sin\phistar$ on both sides of the two inequalities in \ref{itemcasecparametersestimate1} gives $\bigm|\frac{d_i}{r\sin\phistar}-1\bigm|<\frac{17r}{4R}$ for $i=0,2$. Then \eqref{eqJZHypCond}: $R>1700$ gives the conclusion of \ref{itemcasecparametersestimate5}.

\textbf{\cref{proposition:casecparametersestimate}\ref{itemcasecparametersestimate7}} is the claim of \cite[Inequalities (5.15)]{hoalb} and its following statement in \cite{hoalb}.
    
\end{remark}


\subsection[The expansion in Case (c) of \eqref{eqMainCases} with \texorpdfstring{$n_1\ge4$}{n1>3}]{The expansion in Case (c) of \texorpdfstring{\eqref{eqMainCases}}{orbit transition} with \texorpdfstring{$n_1\ge4$}{n1>3}}\label{subsec:expansionLarge}\hfill

\begin{proposition}\label{proposition:casecexpansionn1ge4}
    In the context of \cref{TEcase_c}, that is, Case (c) of \eqref{eqMainCases}, if $n_1\ge4$, then $\mathcal{II}_1(\mathcal{B}^+_0,n_1)$, $\mathcal{II}_2(\mathcal{B}^+_0,n_1)$ in \cref{def:AformularForLowerBoundOfExpansionForn1largerthan0,lemma:AformularForLowerBoundOfExpansionForn1largerthan0} satisfy $\mathcal{II}_1(\Bcal^+_0,n_1)+\mathcal{II}_2(\Bcal^+_0,n_1)>\frac{0.32}{1.01}n_1^2+1.885n_1-2.876>9.73$. 
    
    Hence,$E(n_1,\Bcal^+_0)=\frac{\|dx_3\|_\p}{\|dx_0\|_\p}=3\big|\mathcal{II}_1(\mathcal{B}^+_0,n_1)+\mathcal{II}_2(\mathcal{B}^+_0,n_1)\big|>\frac{0.96}{1.01}n_1^2+5.655n_1-8.628>29.1$ in \cref{lemma:AformularForLowerBoundOfExpansionForn1largerthan0}(4).
\end{proposition}
\begin{proof}
        We first estimate $\mathcal{II}_2(\mathcal{B}^+_0,n_1)=2n_1\big[1+(\tau_0-d_1)\Bcal^+_0\big]$ with $n_1\ge1$.\COMMENT{Fix alignments. \textbackslash mathclap is better than \textbackslash mathllap here.}\COMMENT{\textcolor{red}{Punctuation is missing virtually throughout. Learn how to use it. It will save you time and help make your work readable.}}
        \begin{align*}
            \mathcal{II}_2(\mathcal{B}^+_0,n_1)&=2n_1\big[1+(\tau_0-d_1)\Bcal^+_0\big]\overbracket{>}^{\mathclap{\substack{\text{\cref{lemma:AformularForLowerBoundOfExpansionForn1largerthan0} (1): } \Bcal^+_0<0 \\\text{ \cref{remark:deffpclf}: } \tau_0-d_1<d_1}}}2n_1+2n_1d_1\Bcal_0^+\overbracket{>}^{\mathclap{\substack{\Bcal^+_0<0\text{ and}\\\text{\cref{proposition:casecparametersestimate} \ref{itemcasecparametersestimate4}}}}}2n_1+2n_1\big(\frac{r\sin\phistar}{n_1}+\frac{8,1r}{2n_1R}r\sin\phistar\big)\Bcal_0^+\\
            &=2n_1+\big(2+\frac{8.1r}{R}\big)r\sin\phistar\Bcal_0^+\overbracket{\ge}^{\mathclap{\Bcal_0^+\ge\frac{-4}{3d_0}}}2n_1-\big(2+\frac{8.1r}{R}\big)\frac{4r\sin\phistar}{3d_0}\\
            &\overbracket{>}^{\quad\mathllap{\text{\cref{proposition:casecparametersestimate}  }\ref{itemcasecparametersestimate5}}}2n_1-\frac{8\times1.0026}{3}-\frac{32.4r\times1.0026}{3R}\overbracket{>}^{\mathrlap{\eqref{eqJZHypCond}:R>1700r}}2n_1-2.693,
        \end{align*}
        that is,
        \begin{equation}\label{eq:II2B4estimate1}
        \begin{aligned}
            \text{if }n_1\ge1\text{, then }&\mathcal{II}_2(\mathcal{B}^+_0,n_1)>2n_1-2.693.
        \end{aligned}
        \end{equation}
        Next we estimate $\mathcal{II}_1(\mathcal{B}^+_0,n_1)=\big(1+\tau_0\mathcal{B}^+_0\big)+\big(\tau_1-\frac{2}{3}d_2\big)\Big[\big(2n_1+1-\frac{(2n_1+2)\tau_0}{d_1}\big)\mathcal{B}^+_0-\frac{2n_1+2}{d_1}\Big]$, and we first get the following.\COMMENT{Reword "We can see".\newline Misalignment. \textbackslash mathclap or \textbackslash mathllap should help. Omit all the \textbackslash qquads you put in these things. Almost always they only make things worse.}\COMMENT{\textcolor{red}{Some reasoning labels seem wrong. Check!!}}
        \begin{align*}
        \tau_1-\frac{2}{3}d_2&\overbracket{<}^{\mathclap{\substack{\text{\cref{remark:deffpclf}: }\\\tau_1<2d_1}}}2d_1-\frac{2}{3}d_2\overbracket{<}^{\mathclap{\strut\text{\cref{proposition:casecparametersestimate} \ref{itemcasecparametersestimate1}}}}2d_1-\frac{2}{3}\big(r\sin\phistar-\frac{17r}{4R}r\sin\phistar\big)\\
        &\overbracket{<}^{\quad\mathllap{\text{\cref{proposition:casecparametersestimate}\ref{itemcasecparametersestimate4}}}}\big(\frac{2r\sin{\phistar}}{n_1}+\frac{8.1r\sin{\phistar}}{n_1}\frac{r}{R}\big)-\frac{2}{3}\big(r\sin\phistar-\frac{17r}{4R}r\sin\phistar\big)\\
        &\overbracket{\le}^{\quad\mathllap{n_1\ge4}}\frac{1}{2}r\sin{\phistar}+\frac{8.1r}{4R}r\sin{\phistar}-\frac{2}{3}r\sin{\phistar}+\frac{17r}{6R}r\sin{\phistar}\\
        &=\big(-\frac{1}{6}+\frac{58.3}{12}\cdot\frac{r}{R}\big)r\sin{\phistar}\overbracket{<}^{\mathclap{\eqref{eqJZHypCond}:R>1700r}}-0.16r\sin{\phistar}.
        \end{align*}
        Thus, we have shown the following.\COMMENT{Remember that displays need punctuation. Remember this as you write so there is not so much time wasted fixing that.} \begin{equation}\label{eq:II1B4estimate1}
        \begin{aligned}
            \text{If }n_1\ge1\text{, then }\tau_1-\frac{2}{3}d_2<&\big(\frac{2r\sin\phistar}{n_1}+\frac{8.1r\sin\phistar}{n_1}\frac{r}{R}\big)-\frac{2}{3}(r\sin\phistar-\frac{17R}{4r}r\sin\phistar).\\
            \text{If }n_1\ge4\text{, then }\tau_1-\frac{2}{3}d_2<&-0.16r\sin\phistar.
        \end{aligned}
        \end{equation}
    We next use that $0<\frac{\tau_0}{d_1}<2$ (\cref{remark:deffpclf}), hence $(2n_1+1)-\frac{(2n_1+2)\tau_0}{d_1}> 2n_1+1-2(2n_1+2)=-2n_1-3$, and  $\Bcal^+_0<0$ (\cref{lemma:AformularForLowerBoundOfExpansionForn1largerthan0} (1)) to get
\begin{align*}
    \overbracket{\big[(2n_1+1)-\frac{(2n_1+2)\tau_0}{d_1}\big]}^{>-2n_1-3}\overbracket{\Bcal^+_0}^{<0}-\frac{2n_1+2}{d_1}<&-(2n_1+3)\Bcal^+_0-\frac{2n_1+2}{d_1}\overbracket\le^{\mathclap{\Bcal_0^+\ge{-4}/{3d_0}}}\frac{4}{3d_0}(2n_1+3)-\frac{2n_1+2}{d_1}\\
    \overbracket{<}^{\quad\mathllap{\text{\eqref{eqJZHypCond} and \cref{proposition:casecparametersestimate} \ref{itemcasecparametersestimate1} }: 0<r\sin\phistar-\frac{17r}{4R}r\sin\phistar<d_0}}&\frac{4(2n_1+3)}{3(r\sin{\phistar}-\frac{17r}{4R}r\sin{\phistar})}-\frac{2n_1+2}{d_1}\\
    \overbracket{<}^{\quad\mathllap{\text{\cref{proposition:casecparametersestimate}\ref{itemcasecparametersestimate1}: }0<d_1<\frac{r}{n_1}\sin{\phistar}+\frac{8.1r}{2n_1R}r\sin{\phistar}}}&\frac{4(2n_1+3)}{3(r\sin{\phistar}-\frac{17r}{4R}r\sin{\phistar})}-\frac{(2n_1+2)}{\frac{r}{n_1}\sin{\phistar}+\frac{8.1r}{2n_1R}r\sin{\phistar}}\\
    \overbracket{<}^{\quad\mathllap{\eqref{eqJZHypCond}\text{: R>1700r}}}&\big[\frac{4}{3\times0.99}(2n_1+3)-\frac{n_1(2n_1+2)}{1.01}\big]\frac{1}{r\sin{\phistar}}\\
    =&\big[-\frac{2}{1.01}n_1^2+(\frac{8}{3\times0.99}-\frac{2}{1.01})n_1+\frac{4}{0.99}\big]\frac{1}{r\sin\phistar}\\
    <&\underbracket{\big[-\frac{2}{1.01}n_1^2+0.7135n_1+\frac{4}{0.99}\big]}_{\text{decreases with } n_1\ge1}\frac{1}{r\sin\phistar}\overbracket<^{\mathclap{n_1\ge4}}\frac{-24}{r\sin\phistar}<0.
\end{align*}
Therefore:
\begin{equation}\label{eq:II1B4estimate2}
\begin{aligned}
\text{If }n_1\ge1\text{, then }&\big[(2n_1+1)-\frac{(2n_1+2)\tau_0}{d_1}\big]\Bcal^+_0-\frac{2n_1+2}{d_1}<\big[-\frac{2}{1.01}n_1^2+0.7135n_1+\frac{4}{0.99}\big]\frac{1}{r\sin\phistar}.\\
\text{If }n_1\ge4\text{, then }&\big[(2n_1+1)-\frac{(2n_1+2)\tau_0}{d_1}\big]\Bcal^+_0-\frac{2n_1+2}{d_1}<\frac{-24}{r\sin\phistar}<0.
\end{aligned}
\end{equation}
Furthermore, since $(n_1+1)\tau_0=\tau_0+n_1\tau_0\overbracket{<}^{\mathclap{\substack{\text{\cref{remark:deffpclf}:}\\\tau_0<2d_1}}}\tau_0+2n_1d_1+\tau_1\overbracket{<}^{\mathclap{\substack{\text{\cref{proposition:casecparametersestimate}\ref{itemcasecparametersestimate2}}}}}2r\sin\phistar+\frac{8.1r}{R}r\sin\phistar$, we have \[0<\tau_0<\frac{2r\sin\phistar}{n_1+1}+\frac{8.1r}{(n_1+1)R}r\sin\phistar=\big[\frac{2}{n_1+1}+\frac{8.1r}{(n_1+1)R}\big]r\sin\phistar,\] as well as\COMMENT{Use overbracket and mathclap on the last inequality}
    \begin{align*}
    1+\underbracket{\tau_0}_{\mathclap{<[\frac{2}{n_1+1}+\frac{8.1r}{(n_1+1)R}]r\sin\phistar}}\overbracket{\Bcal_0^+}^{\mathclap{\text{By \cref{lemma:AformularForLowerBoundOfExpansionForn1largerthan0} \ref{itemAformularForLowerBoundOfExpansionForn1largerthan01}, }\Bcal^+_0<0}}>&1+\big[\frac{2}{n_1+1}+\frac{8.1r}{(n_1+1)R}\big]r\sin\phistar\Bcal^+_0\\
    \overbracket{>}^{\quad\mathllap{\text{\cref{lemma:AformularForLowerBoundOfExpansionForn1largerthan0} (1), }\Bcal^+_0\in[\frac{-4}{3d_0},\frac{-1}{d_0}]}}&1+\big[\frac{2}{n_1+1}+\frac{8.1r}{(n_1+1)R}\big]\frac{-4r\sin\phistar}{3d_0}=1-\frac{4}{3}\big[\frac{2}{n_1+1}+\frac{8.1r}{(n_1+1)R}\big]\frac{r\sin\phistar}{d_0}\\
    \overbracket{>}^{\quad\mathllap{\text{\cref{proposition:casecparametersestimate} \ref{itemcasecparametersestimate5}}}}&1-\frac{4}{3}\big[\frac{2}{1+n_1}+\frac{8.1r}{(1+n_1)R}\big]\times1.0026\overbracket>^{\mathclap{\substack{\text{\eqref{eqJZHypCond}: R>1700r}\\\text{ and }n_1\ge4}}}1-\frac{4}{3}(\frac{2}{5}+\frac{8.1}{1700\cdot 5})\times1.0026>0.4640.
    \end{align*}
Thus:\COMMENT{Check last inequality in previous display}
    \begin{equation}\label{eq:II1B4estimate3}
    \begin{aligned}
        \text{If }n_1\ge1\text{, then }&1+\tau_0\Bcal^+_0>1+\big[\frac{2}{n_1+1}+\frac{8.1r}{(n_1+1)R}\big]\times\big(\frac{-4r\sin\phistar}{3d_0}\big).\\
        \text{If }n_1\ge4\text{, then }&1+\tau_0\Bcal^+_0>0.464.
    \end{aligned}
    \end{equation}
    Therefore, 
    if $n_1\ge4$, then \begin{align*}\mathcal{II}_1(\mathcal{B}^+_0,n_1)+\mathcal{II}_2(\mathcal{B}^+_0,n_1)&=\big(\overbracket{1+\tau_0\mathcal{B}^+_0}^{\mathclap{\text{\eqref{eq:II1B4estimate3}: }>0.464}}\big)+\big(\underbracket{\tau_1-\frac{2}{3}d_2}_{\mathclap{\text{\eqref{eq:II1B4estimate1}:}<-0.16r\sin\phistar<0}}\big)\Big[\overbracket{\big(2n_1+1-\frac{(2n_1+2)\tau_0}{d_1}\big)\mathcal{B}^+_0-\frac{2n_1+2}{d_1}}^{\mathclap{\text{\eqref{eq:II1B4estimate2}: }<\big(-\frac{2}{1.01}n_1^2+0.7135n_1+\frac{4}{0.99}\big)\frac{1}{r\sin\phistar}<\frac{-24}{r\sin\phistar}<0}}\Big]
    +\overbracket{\mathcal{II}_2(\mathcal{B}^+_0,n_1)}^{\mathclap{\text{\eqref{eq:II2B4estimate1}: }>2n_1-2.693}}\\
    &>0.464+(-0.16)\times(-\frac{2}{1.01}n_1^2+0.7135n_1+\frac{4}{0.99})+2n_1-2.693\\
    &>\frac{0.32}{1.01}n_1^2+1.885n_1-2.876.\qedhere
    \end{align*}
\end{proof}

\subsection[The expansion in Case (c) of \eqref{eqMainCases} with \texorpdfstring{$n_1=3$}{n1=3}]{The expansion in Case (c) of \texorpdfstring{\eqref{eqMainCases}}{orbit transition} with \texorpdfstring{$n_1=3$}{n1=3}}\hfill

By symmetry, assume that $x_1$ with $n_1=3$ is on the lower component in \cref{fig:x1_multi_collisions_caseC} with $\theta_1\in (0,\frac{\pi}{2})$. More specifically, $x_1=(\Phi_1,\theta_1)$ is in the \cref{fig:x1_n1cases_caseC} labeled\COMMENT{Wrong word; lots of things are labeled in \cref{fig:x1_n1cases_caseC}.} region\COMMENT{Wrong word order. See the fixed word order above, learn it, and fix it here. Ask for help if needed.} with $\theta_1$ varying in $\big(\frac{1}{5}\Phistar,\frac{1}{3}\Phistar\big)$.
We use a horizontal line $\theta_1=0.235\Phistar$ to separate this region into two parts: upper part $\theta_1\ge0.235\Phistar$, lower part $\theta_1<0.235\Phistar$, as shown in \cref{fig:x1_n1equal3case_subcases}.\COMMENT{Nice figure placement}

\begin{figure}[h]
\begin{minipage}{.54\textwidth}
\begin{center}
    \begin{tikzpicture}[xscale=0.6, yscale=0.6]
        \pgfmathsetmacro{\PHISTAR}{6};
        \pgfmathsetmacro{\WIDTH}{8};
        \tkzDefPoint(-\PHISTAR,0){LL};
        \tkzDefPoint(-\PHISTAR,\WIDTH){LU};
        \tkzDefPoint(\PHISTAR,0){RL};
        \tkzDefPoint(\PHISTAR,\WIDTH){RU};
        \tkzDefPoint(-\PHISTAR,\PHISTAR){X};
        \tkzDefPoint(\PHISTAR,\WIDTH-\PHISTAR){Y};
        \tkzDefPoint(-\PHISTAR,\WIDTH-\PHISTAR){IX};
        \tkzDefPoint(\PHISTAR,\PHISTAR){IY};
        \tkzDefPoint(-\PHISTAR,0.5*\PHISTAR){S1};
        \tkzDefPoint(-\PHISTAR,\PHISTAR/3){S2};
        \tkzDefPoint(-\PHISTAR,\PHISTAR/4){S3};
        \tkzDefPoint(-\PHISTAR,\PHISTAR/5){L4};
        \draw [thin] (LL) --(LU);
        \draw [thin] (LU) --(RU);
        \draw [thin] (RL) --(RU);
        \draw [thin] (LL) --(RL);
        \draw [name path=lineLLIY, blue,ultra thin,dashed] (LL) --(IY);
        \draw [red,ultra thin,dashed]  (X)--(RL);
        \draw [red,ultra thin,dashed]  (S1)--(RL);
        \draw [name path=lineS2RL, red,ultra thin,dashed]  (S2)--(RL);
        \draw [name path=lineS3RL, red,ultra thin,dashed]  (S3)--(RL);
        \path[name intersections= 
        {of=lineLLIY and lineS3RL,by={LS}}];
        \path[name intersections=
        {of=lineLLIY and lineS2RL,by={HS}}];
        \tkzDrawPoints(S2,S3,HS,LS,L4);
        \draw [ultra thin,dashed] (L4) --(LS);
        \tkzLabelPoint[above left](S2){\small$\theta_1=\frac{1}{3}\Phistar$};
        \tkzLabelPoint[left](S3){\small$\theta_1=0.25\Phistar$};
        \tkzLabelPoint[below left](L4){\small $\theta_1=0.2\Phistar$};
        \node[] at (0,1.1*\PHISTAR){\Large $n_1=0$};
        \node[] at (-0.7*\PHISTAR,0.6*\PHISTAR){\Large $n_1=1$};
        \node[] at (-0.83*\PHISTAR,0.38*\PHISTAR){\Large $n_1=2$};
        \node[] at (-0.85*\PHISTAR,0.27*\PHISTAR){\small $n_1=3$};
        \node[] at (-0.90*\PHISTAR,0.13*\PHISTAR){ $n_1\ge4$};
        \fill [green, opacity=5/30](LS) -- (HS) -- (S2) -- (S3) -- cycle;
        \tkzLabelPoint[below](LL){$-\Phistar$};
        \tkzLabelPoint[below](RL){$+\Phistar$};
    \end{tikzpicture}
    \captionsetup{{width=.9\linewidth}}
    \caption{Quadrilateral of $n_1=3$, $\theta_1\in[0.2\Phistar,\frac{1}{3}\Phistar]$}\label{fig:x1_n1cases_caseC}
\end{center}
\end{minipage}\hfill
\begin{minipage}{0.45\textwidth}
\begin{center}
    \begin{tikzpicture}[xscale=2.2, yscale=2.2]
        \pgfmathsetmacro{\PHISTAR}{6};
        \pgfmathsetmacro{\WIDTH}{8};
        \tkzDefPoint(-\PHISTAR,0){LL};
        \tkzDefPoint(-\PHISTAR,\WIDTH){LU};
        \tkzDefPoint(\PHISTAR,0){RL};
        \tkzDefPoint(\PHISTAR,\WIDTH){RU};
        \tkzDefPoint(-\PHISTAR,\PHISTAR){X};
        \tkzDefPoint(\PHISTAR,\WIDTH-\PHISTAR){Y};
        \tkzDefPoint(-\PHISTAR,\WIDTH-\PHISTAR){IX};
        \tkzDefPoint(\PHISTAR,\PHISTAR){IY};
        \tkzDefPoint(-\PHISTAR,0.5*\PHISTAR){S1};
        \tkzDefPoint(-\PHISTAR,\PHISTAR/3){S2};
        \tkzDefPoint(-\PHISTAR,\PHISTAR/4){S3};
        \tkzDefPoint(-\PHISTAR,\PHISTAR/5){L4};
        \clip(-\PHISTAR-0.05,0.15*\PHISTAR) rectangle (-\PHISTAR+\PHISTAR*0.52, \PHISTAR*0.35);
        \draw [thin] (LL) --(LU);
        \draw [thin] (LU) --(RU);
        \draw [thin] (RL) --(RU);
        \draw [thin] (LL) --(RL);
        \draw [name path=lineLLIY, blue,ultra thin,dashed] (LL) --(IY);
        \draw [red,ultra thin,dashed]  (X)--(RL);
        \draw [red,ultra thin,dashed]  (S1)--(RL);
        \draw [name path=lineS2RL, red,ultra thin,dashed]  (S2)--(RL);
        \draw [name path=lineS3RL, red,ultra thin,dashed]  (S3)--(RL);
        \path[name intersections= 
        {of=lineLLIY and lineS3RL,by={LS}}];
        \path[name intersections=
        {of=lineLLIY and lineS2RL,by={HS}}];
        \tkzDefPoint(-\PHISTAR,0.235*\PHISTAR){L5};
        \tkzDefPoint(-\PHISTAR,0.222*\PHISTAR){L6};
        \path[name path=HorSeg] (L5)-- +(5,0);
        \path[name path=HorSeg1] (L6)-- +(5,0);
        \path[name intersections={of=lineS3RL and HorSeg,by={Seg1}}];
        \path[name intersections={of=lineLLIY and HorSeg,by={Seg2}}];
        \draw[very thin] (L5)--(Seg2);
        \tkzLabelPoint[above](Seg1){\small$\theta_1\ge0.235\Phistar$};
        \tkzLabelPoint[below](Seg1){\small$\theta_1<0.235\Phistar$};
        \path[name intersections={of=lineS3RL and HorSeg1,by={Seg3}}];
        \path[name intersections={of=lineLLIY and HorSeg1,by={Seg4}}];
        \draw[very thin,blue] (Seg3)--(Seg4);
        \tkzLabelPoint[below](Seg3){\small$\mathcal{W}$};
        \tkzLabelPoint[below](Seg4){\small$\mathcal{E}$};
         \tkzLabelSegment[below](Seg3,Seg4){\small{Line segment $\mathfrak{l}$}};  \tkzDrawPoints(S2,S3,HS,LS,L5,Seg1,Seg2,Seg3,Seg4);
        \fill [green, opacity=5/30](LS) -- (HS) -- (S2) -- (S3) -- cycle;
        \tkzLabelPoint[below](LL){$-\Phistar$};
        \tkzLabelPoint[below](RL){$+\Phistar$};
    \end{tikzpicture}
    \caption{Quadrilateral $n_1=3$ divided into two parts: upper part with $\theta_1\ge0.235\Phistar$ and lower part with $\theta_1<0.235\Phistar$. Line segment $\mathfrak{l}=\big\{(\Phi_1,\theta_1=\const)\big\}$ with endpoints $\mathcal{W}$, $\mathcal{E}$ in lower part.}\label{fig:x1_n1equal3case_subcases}
\end{center}
\end{minipage}\hfill
\end{figure}

\begin{proposition}\label{proposition:casecexpansionn1eq3}
    In the context of \cref{TEcase_c}, that is, Case (c) of \eqref{eqMainCases}, if $n_1=3$, then $\mathcal{II}_1(\mathcal{B}^+_0,3)$, $\mathcal{II}_2(\mathcal{B}^+_0,3)$ in \cref{def:AformularForLowerBoundOfExpansionForn1largerthan0,lemma:AformularForLowerBoundOfExpansionForn1largerthan0} have\COMMENT{Wrong word. Learn this and fix. "satisfy" is likely what you mean.} $\mathcal{II}_1(\Bcal^+_0,3)+\mathcal{II}_2(\Bcal^+_0,3)>4.28$. Hence in (4) of \cref{lemma:AformularForLowerBoundOfExpansionForn1largerthan0}, $E(3,\Bcal^+_0)=\frac{\|dx_3\|_\p}{\|dx_0\|_\p}=3\big|\mathcal{II}_1(\mathcal{B}^+_0,3)+\mathcal{II}_2(\mathcal{B}^+_0,3)\big|>12.8$.
\end{proposition}
\begin{proof}
    We first estimate $\tau_1-\frac{2}{3}d_2$ for $x_1=(\Phi_1,\theta_1)$ with $n_1=3$. 
    \begin{equation}\label{eq:tau1minus2thirdsd2n1is3}
        \begin{aligned}
        &\left.\begin{aligned}
        &-\frac{2}{3}d_2\overset{\text{\cref{proposition:casecparametersestimate} (i)}}<-\frac{2}{3}r\sin\phistar+\frac{2}{3}\cdot\frac{17r}{4R}r\sin\phistar\\
        &\tau_1\overset{\text{\cref{proposition:casecparametersestimate} \ref{itemcasecparametersestimate2}}}<2r\sin\phistar-6d_1-\tau_0+\frac{8.1r}{R}r\sin\phistar
        \end{aligned}
        \right\}\Rightarrow\tau_1-\frac{2}{3}d_2<\frac{4r\sin\phistar}{3}-6d_1-\tau_0+\frac{65.6r}{6R}r\sin\phistar
    \end{aligned}
    \end{equation}
    
    Then we analyze the two cases based on $\theta_1$ i.e. $x_1 $in two parts of \cref{fig:x1_n1equal3case_subcases}.\NEWPAGE
    
    \textbf{Case 1}: $0.235\Phistar\le\theta_1<\Phistar/3$ and $x_1$ is in the upper part of the $n_1=3$ quadrilateral  in \cref{fig:x1_n1equal3case_subcases}.\COMMENT{Below, $sinc$ does not work. It is the product of variables $x$, $i$, $n$ and $c$!\Hrule\small Use \textbackslash sinc throughout!\Hrule Also, not everyone knows $\sinc$.\Hrule Inside cases use underbracket and mathclap}
    \begin{align*}
        &\left.\begin{aligned}
        &0.235\Phistar\overset{\text{\cref{def:BasicNotations0}}}=0.235\sin^{-1}\big(\frac{r\sin\phistar}{R}\big)\overset{\text{\eqref{eqJZHypCond}: }R>1700r}<0.2\\
        &\text{If }0<z<0.2\text{, then } \mathrm{sinc}(z)>1/1.01\\
        \end{aligned}
        \right\} \Rightarrow \frac{\sin{(0.235\Phistar)}}{0.235\Phistar}> \frac{1}{1.01} \\
        &\Rightarrow d_1=R\sin{\theta_1}\overset{0.235\Phistar\le\theta_1<\pi/2}\ge R\sin{(0.235\Phistar)}>\frac{0.235}{1.01}R\Phistar>\frac{0.235}{1.01}R\sin{\Phistar}=\frac{0.235}{1.01}r\sin{\phistar}\\
    \end{align*} 
    
    Hence in \textbf{case 1} we have\COMMENT{Use overbracket and mathclap\wentao{Done}} \begin{equation}\label{eq:tau1minus2thirdsd2n1is3case1}
        \begin{aligned}
        \tau_1-\frac{2}{3}d_2&\overbracket{<}^{\mathclap{\eqref{eq:tau1minus2thirdsd2n1is3}}}\frac{4}{3}r\sin{\phistar}-6d_1-\tau_0+\frac{65.6r}{6R}r\sin{\phistar}\\
        &\overbracket{<}^{\mathllap{d_1>\frac{0.235}{1.01}r\sin\phistar}}\frac{4}{3}r\sin{\phistar}-\frac{6\times0.235}{1.01}r\sin{\phistar}+\frac{65.6r}{6R}r\sin{\phistar}\\
        &<(-0.0625+\frac{65.6r}{6R})r\sin{\phistar}\overset{\text{\eqref{eqJZHypCond}:} R>1700r}<-0.056r\sin{\phistar} 
        \end{aligned}
    \end{equation}
\textbf{Case 2: } $\theta_1\in \big(0.2\Phistar,0.235\Phistar\big)$, $x_1$ is in the lower part of the $n_1=3$ quadrilateral in \cref{fig:x1_n1equal3case_subcases}.

In \cref{fig:x1_n1equal3case_subcases}, suppose that $x_1$ is on a line segment $\mathfrak{l}=\Big\{(
\Phi_1,\theta_1=\const
=\nu\Phistar)\bigm|\Phi_1\in\big[\Phi_{\mathcal{W}},\Phi_{\mathcal{E}}\big)\Big\}$ with end points $\mathcal{W}=(\Phi_{\mathcal{W}},\nu\Phistar)$ and $\mathcal{E}=(\Phi_{\mathcal{E}},\nu\Phistar)$ on the boundary of the lower part of the $n_1=3$ quadrilateral, $\Phi_{\mathcal{W}}\le\Phi_1<\Phi_{\mathcal{E}}<0$ and $\nu\in\big(0.2,0.235\big)$. Note that $\tau_0$ as a function of $x_1=(\Phi_1,\theta_1=\const)\in\mathfrak{l}$ and thus of $\Phi_1\in\big[\Phi_{\mathcal{W}},\Phi_{\mathcal{E}}\big)$.

In the coordinates of \cref{def:StandardCoordinateTable} for \cref{fig:n1_equal3_leftendpoint_case}, $(x_P,y_P)=P=p(x_1)=(R\sin{\Phi_1},b-R\cos{\Phi_1})$, $A=(x_A,y_A)=(-r\sin\phistar,-r\cos\phistar)$, $P_\mathcal{W}\dfn p(\mathcal{W})=(x_\mathcal{W},y_\mathcal{W})=(R\sin{(\Phi_{\mathcal{W}})},b-R\cos{(\Phi_{\mathcal{W}})})$. Since $-\pi/2<-\Phistar<\Phi_\mathcal{W}\le\Phi_1<0$ and by checking \cref{fig:n1_equal3_leftendpoint_case}, we have $x_A<x_\mathcal{W}\le x_P$ and $y_A>y_\mathcal{W}\ge y_P$.

Suppose that $(x_{T_0},y_{T_0})=T_0=p(x_0)=p(\Fcal^{-1}(x_1))$ with $x_0=(\phi_0,\theta_0)$. By \cref{contraction_region}, $x_0\in\Nin$. More specifically, we see $T_0$ is in some neighborhood of corner $A$ (see \cref{fig:n1_equal3_leftendpoint_case} for a belief). Furthermore, by the definition of \Nin in \cref{def:SectionSetsDefs} we have 
\[
0<2\pi-\phistar-\phi_0\overbracket{<}^{\mathclap{\strut\text{definition of \Nin}}}\frac{17}{\sin{(\phistar/2)}}\sqrt{\frac rR}\Rightarrow 2\pi-\phistar>\phi_0>2\pi-\phistar-\frac{17}{\sin(\phistar/2)}\sqrt{\frac rR}\overbracket{>}^{\mathclap{\substack{\phistar\in(0,\tan^{-1}(1/3))\ \&\text{ \eqref{eqJZHypCond}: }R>\frac{30000r}{\sin^2(\phistar/2)}}}}\frac{3\pi}{2}
\] 
Then $(x_{T_0},y_{T_0})=(r\sin\phi_0,-r\cos\phi_0)$ has $x_{T_0}<x_A$ and $y_{T_0}>y_A$.\COMMENT{Edit\wentao{done, better sentence now.}} will yield the following.
\begin{align*}
    &\left.\begin{aligned}
        &x_{T_0}<x_A<x_{\mathcal{W}}\le x_P\\
        &y_{T_0}>y_A>y_{\mathcal{W}}\ge y_P
        \end{aligned}
    \right\} \Rightarrow |AP_{\mathcal{W}}|=\sqrt{(x_A-x_\mathcal{W})^2+(y_A-y_\mathcal{W})^2}<\sqrt{(x_P-x_{T_0})^2+(y_{T_0}-y_P)^2}=|T_0P|=\tau_0
\end{align*}
\NEWPAGE So $|AP_{\mathcal{W}}|$ is a lower bound of $\tau_0(x_1)$ for all $x_1\in\mathfrak{l}$. Since $p(\Fcal^4{\mathcal{W}})=B$, from \cref{fig:n1_equal3_leftendpoint_case} we get $|AP_\mathcal{W}|=2R\sin{\big(\frac{2\Phistar-8\theta_1}{2}\big)}$ and\COMMENT{{\color{red}\small This display needs to be disassembled.\Hrule}Do not forget punctuation.\Hrule The purpose of the large overbracket is unclear. It's a jumble without indication of what is used how and where.\Hrule Fix $sinc$ below and throughout.\wentao{Done} Use \textbackslash sinc instead of sinc everywhere and define sinc when it first appears.\Hrule\small Use mathclap (plus strut or bigstrut as needed to fix text collisions).}
\begin{align*}
        &\left.\begin{aligned}
        &\overbracket{0<(1-4\nu)\Phistar<0.2\Phistar=0.2\sin^{-1}(\frac{r\sin\phistar}{R})<0.2\text{ and }0<\nu\Phistar<0.235\Phistar<0.2}^{\text{By }\nu\in(0.2,0.235)\text{, \cref{def:BasicNotations0}: } r\sin\phistar=R\sin\Phistar \text{ and \eqref{eqJZHypCond}: }R>1700r}\\
        &\text{If }0\le z<0.2\text{, then }\sinc(z)\dfn\frac{\sin z}z>\frac{1}{1.01}
        \end{aligned}
        \right\}\Longrightarrow \\
        \Longrightarrow 
        &\left\{\begin{aligned}
        \frac{\sin{\big[(1-4\nu)\Phistar\big]}}{(1-4\nu)\Phistar}=\sinc\big[(1-4\nu)\Phistar\big]&>\frac{1}{1.01}\\
        \frac{\sin(\nu\Phistar)}{\nu\Phistar}=\sinc(\nu\Phistar)&>\frac{1}{1.01}
        \end{aligned}
        \right.\\
        \Longrightarrow&\left\{\begin{aligned}\tau_0>&|AP_\mathcal{W}|=2R\sin{\big(\frac{2\Phistar-8\theta_1}{2}\big)}\overbracket{=}^{\theta_1=\nu\Phistar}2R\sin{\big[(1-4\nu)\Phistar\big]}>2R\frac{(1-4v)\Phistar}{1.01}\\
        >&2R\frac{1-4v}{1.01}\sin{\Phistar}\overset{R\sin{\Phistar}=r\sin{\phistar}}=\frac{1-4v}{1.01}(2r\sin{\phistar})\\
        d_1=&R\sin\theta_1\overbracket{=}^{\theta_1=\nu\Phistar}R\sin{(\nu\Phistar)}>\frac{R\nu\Phistar}{1.01}>\frac{\nu R\sin\Phistar}{1.01}\overset{R\sin{\Phistar}=r\sin{\phistar}}=\frac{\nu}{1.01}r\sin\phistar
        \end{aligned}
        \right.
\end{align*}
Hence in \textbf{case 2} we have \begin{equation}\label{eq:tau1minus2thirdsd2n1is3case2}
        \begin{aligned}
        \tau_1-\frac{2}{3}d_2&\overbracket{<}^{\mathclap{\eqref{eq:tau1minus2thirdsd2n1is3}}}\frac{4}{3}r\sin{\phistar}-6d_1-\tau_0+\frac{65.6r}{6R}r\sin{\phistar}\\
        &\overbracket{<}^{\mathllap{\tau_0>\frac{(1-4\nu)2r\sin\phistar}{1.01}\text{ and }d_1>\frac{\nu r\sin\phistar}{1.01}}}\frac{4}{3}r\sin{\phistar}-\frac{6\nu}{1.01}r\sin{\phistar}-\frac{1-4\nu}{1.01}2r\sin{\phistar}+\frac{65.6r}{R}r\sin{\phistar}\\
        &=(\frac{4}{3}-\frac{2}{1.01}+\frac{2\nu}{1.01})r\sin{\phistar}+\frac{65.6r}{6R}r\sin{\phistar}\\
        &\overbracket{<}^{\mathllap{\nu<0.235\text{ and \eqref{eqJZHypCond}: }R>1700r}}-0.181r\sin\phistar+0.007r\sin\phistar<-0.17r\sin\phistar
        \end{aligned}
    \end{equation}

To summarize, both in \textbf{case 1} \eqref{eq:tau1minus2thirdsd2n1is3case1} and in \textbf{case 2} \eqref{eq:tau1minus2thirdsd2n1is3case2}, we get $\tau_1-\frac{2}{3}d_2<-0.056r\sin{\phistar}$. 
Furthermore,\COMMENT{Fix this alignment jumble. Useless alignment, needlessly repeated arguments create confusion. Pay more attention to readability. Note that the last underbracket pretends that the first equation is not there. You made more work for yourself by not planning ahead. Instead there should be a reference to the top equation.\Hrule mathclap will help} \begin{align*}\big[(2n_1+1)-\frac{(2n_1+2)\tau_0}{d_1}\big]\Bcal^+_0-\frac{2n_1+2}{d_1}&\overbracket{<}^{\mathclap{\eqref{eq:II1B4estimate2}}}\big[-\frac{2}{1.01}n_1^2+0.7135n_1+\frac{4}{0.99}\big]\frac{1}{r\sin\phistar}\overset{n_1=3}{<}\frac{-11.64}{r\sin\phistar}<0\\
1+\tau_0\Bcal^+_0&\underbracket{>}_{\mathclap{\text{\eqref{eq:II1B4estimate3}}}}1+\big[\frac{2}{n_1+1}+\frac{8.1r}{(n_1+1)R}\big]\times\big(\frac{-4r\sin\phistar}{3d_0}\big)\underbracket{=}_{n_1=3}1+\big[0.5+\frac{8.1r}{4R}\big]\big(\frac{-4r\sin\phistar}{3d_0}\big)\\
&\underbracket{>}_{\mathllap{\text{\cref{proposition:casecparametersestimate}\ref{itemcasecparametersestimate5}}}}1+\big(0.5+\frac{8.1r}{4R}\big)\frac{-4\times1.0026}{3}\underbracket{>}_{\text{\eqref{eqJZHypCond}: }R>1700r}0.33\\
\mathcal{II}_1(\mathcal{B}^+_0,n_1=3)=\underbracket{\big(1+\tau_0\mathcal{B}^+_0\big)}_{>0.33}+&\underbracket{\big(\tau_1-\frac{2}{3}d_2\big)}_{\text{\eqref{eq:II1B4estimate1}:}<-0.056r\sin\phistar<0}\underbracket{\Big[\big(2n_1+1-\frac{(2n_1+2)\tau_0}{d_1}\big)\mathcal{B}^+_0-\frac{2n_1+2}{d_1}\Big]\big|_{n_1=3}}_{\text{\eqref{eq:II1B4estimate2} with }n_1=3\text{: }<\frac{-11.64}{r\sin\phistar}<0}\\
>0.33+0.056\times&11.64>0.98\\
\mathcal{II}_2(\mathcal{B}^+_0,n_1=3)\overbracket{>}^{\text{\eqref{eq:II2B4estimate1}}}(2n_1-2.693)\bigm|_{n_1=3}&=3.307
\end{align*}
Therefore, $\mathcal{II}_1(\mathcal{B}^+_0,3)+\mathcal{II}_2(\mathcal{B}^+_0,3)>0.98+3.307=4.287$. This proves the claim.
\end{proof}
\begin{figure}[h]
\begin{center}
\begin{tikzpicture}[]
    \tkzDefPoint(0,0){Or};
    \tkzDefPoint(0,-4.8){Y};
    \pgfmathsetmacro{\rradius}{4.5};
    \pgfmathsetmacro{\phistardeg}{40};
    \pgfmathsetmacro{\XValueArc}{\rradius*sin(\phistardeg)};
    \pgfmathsetmacro{\YValueArc}{\rradius*cos(\phistardeg)};
    \pgfmathsetmacro{\Rradius}{
        \rradius*sin(\phistardeg)/sin(18)
    };
    \pgfmathsetmacro{\bdist}{(-\rradius*cos(\phistardeg))+cos(18)*\Rradius}
    \tkzDefPoint(0,\bdist){OR};
    \tkzDefPoint(\XValueArc,-1.0*\YValueArc){B};    \tkzDefPoint(-1.0*\XValueArc,-1.0*\YValueArc){A};
    \draw[name path=Cr,green,ultra thin] (Or) circle (\rradius);
    \tkzDrawArc[name path=CR](OR,A)(B);
    \pgfmathsetmacro{\PXValue}{-\Rradius*sin(14)};
    \pgfmathsetmacro{\PYValue}{\bdist-\Rradius*cos(14)};
    \tkzDefPoint(\PXValue,\PYValue){P};
    \pgfmathsetmacro{\PXValue}{-\Rradius*sin(6)};
    \pgfmathsetmacro{\PYValue}{\bdist-\Rradius*cos(6)};
    \tkzDefPoint(\PXValue,\PYValue){P1};
    \pgfmathsetmacro{\PXValue}{\Rradius*sin(2)};
    \pgfmathsetmacro{\PYValue}{\bdist-\Rradius*cos(2)};
    \tkzDefPoint(\PXValue,\PYValue){P2};
    \pgfmathsetmacro{\PXValue}{\Rradius*sin(10)};
    \pgfmathsetmacro{\PYValue}{\bdist-\Rradius*cos(10)};
    \tkzDefPoint(\PXValue,\PYValue){P3};
    \path[name path=Line_farT_P](P)-- +(162:3.5);
    \path[name intersections={of=Cr and Line_farT_P,by={T}}];
    \tkzDefPoint(0,\bdist-\Rradius){C};
    \draw[red,ultra thin,dashed](OR)--(P);
    \draw[red,ultra thin,dashed](OR)--(A);
    \draw[red,ultra thin,dashed](OR)--(B);
    \begin{scope}[ultra thin,decoration={markings,mark=at position 0.95 with {\arrow{>}}}] 
    \draw[blue,postaction={decorate}] (T)--(P);
    \draw[blue,postaction={decorate}] (P)--(P1);
    \draw[blue,postaction={decorate}] (P1)--(P2);
    \draw[blue,postaction={decorate}] (P2)--(P3);
    \draw[blue,postaction={decorate}] (P3)--(B);
    \end{scope}
    \tkzLabelPoint[above](P){$P_\mathcal{W}$};
    \tkzLabelPoint[below](A){$A$};
    \tkzLabelPoint[above left](T){$T_\mathcal{W}$};
    \tkzLabelPoint[above](P1){\small $p(\mathcal{F}(\mathcal{W}))$};
    \tkzLabelPoint[below](P2){\small $p(\mathcal{F}^2(\mathcal{W}))$};
    \tkzLabelPoint[above](P3){\small $p(\mathcal{F}^3(\mathcal{W}))$};
    \tkzLabelPoint[below right](B){\small $B=p(\mathcal{F}^4(\mathcal{W}))$};
    \tkzLabelPoint[below](Or){$O_r$};
    \tkzLabelPoint[above](OR){$O_R$};
    \tkzMarkAngle[arc=lll,ultra thin,size=4.5](A,OR,P);
    \tkzMarkAngle[arc=lll,ultra thin,size=2](P,OR,B);
    \tkzDrawPoints(A,B,P,Or,OR,T,P1,P2,P3);
    \node[] at ($(OR)+(-88:3.2)$)  {\Large $8\theta_1$};
    \node[] at ($(OR)+(-110:5.5)$)  { $2\Phistar-8\theta_1$};
\end{tikzpicture}
\caption{ For $\mathcal{W}=(\Phi_{\mathcal{W}},\theta_1)$ in \cref{fig:x1_n1equal3case_subcases},
$P_\mathcal{W}=p(\mathcal{W})$, $T_\mathcal{W}=p(\mathcal{F}^{-1}(\mathcal{W}))$, $p(\mathcal{F}^4(\mathcal{W}))=B$. $\measuredangle{AO_RP_\mathcal{W}}=2\Phistar-8\theta_1$, $\measuredangle{P_\mathcal{W}O_RB}=8\theta_1$, $|AP_\mathcal{W}|=2R\sin{(\frac{2\Phistar-8\theta_1}{2})}$}\label{fig:n1_equal3_leftendpoint_case}
\end{center}
\end{figure}
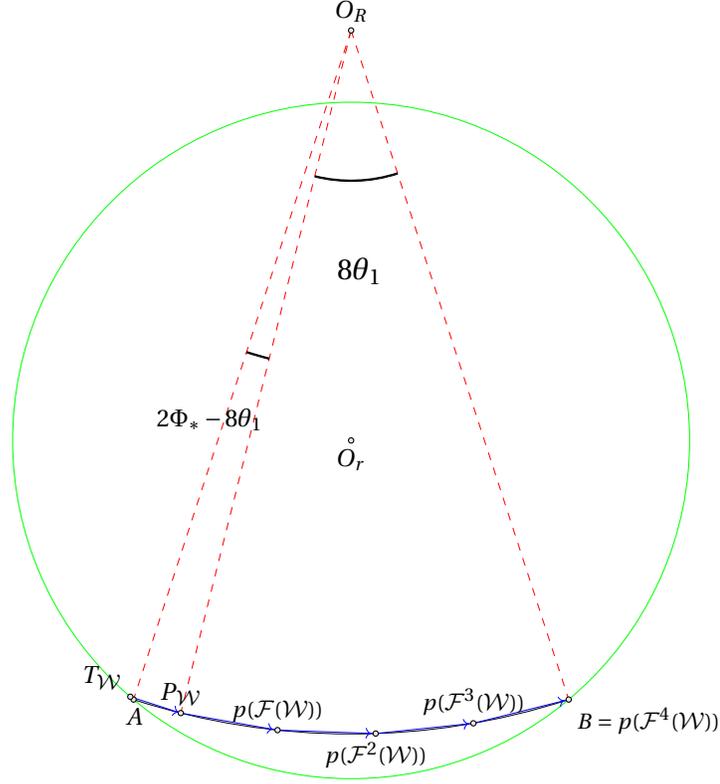

\subsection[Expansion in Case (c) of \eqref{eqMainCases} with \texorpdfstring{$n_1=2$}{n1=2}]{The expansion in Case (c) of \texorpdfstring{\eqref{eqMainCases}}{orbit transition} with \texorpdfstring{$n_1=2$}{n1=2}}\hfill
\begin{proposition}\label{proposition:casecexpansionn1eq2}
    In the context of \cref{TEcase_c}, that is, case (c) of \eqref{eqMainCases},\COMMENT{Not needed once we have a standing assumption, right?} if $n_1=2$, then 
    $\mathcal{II}_1(\Bcal^+_0,2)+\mathcal{II}_2(\Bcal^+_0,2)>1.86$ (see \cref{def:AformularForLowerBoundOfExpansionForn1largerthan0,lemma:AformularForLowerBoundOfExpansionForn1largerthan0}). 
    
    Hence in \cref{lemma:AformularForLowerBoundOfExpansionForn1largerthan0}(4), $E(2,\Bcal^+_0)=\frac{\|dx_3\|_\p}{\|dx_0\|_\p}=3\big|\mathcal{II}_1(\mathcal{B}^+_0,2)+\mathcal{II}_2(\mathcal{B}^+_0,2)\big|>5.58$.
\end{proposition}
\begin{proof}
We consider two cases: $\tau_1-\frac{2}{3}d_2>0$ and $\tau_1-\frac{2}{3}d_2\le0$.

\hidesubsection[$\tau_1-\frac{2}{3}d_2>0$]{Case 1} We\COMMENT{\color{red}This use of \textbackslash hidesection is good for cases.} first show that, in fact, $0<\tau_1-\frac{2}{3}d_2<0.06r\sin\phistar$:
\begin{align*}
        &\frac{2}{3}(r\sin{\phistar}-\frac{17r}{4R}r\sin{\phistar})+4d_1\overbracket{<}^{\mathclap{\strut\text{\cref{proposition:casecparametersestimate} \ref{itemcasecparametersestimate1}}}}\frac{2}{3}d_2+4d_1\overbracket{<}^{\mathclap{\tau_1>2d_2/3}}\tau_1+4d_1<\tau_1+4d_1+\tau_0\overbracket{<}^{\mathclap{\text{\cref{proposition:casecparametersestimate} \ref{itemcasecparametersestimate2} with } n_1=2}}2r\sin\phistar+\frac{8.1r}{R}r\sin\phistar\\
        \Rightarrow&\frac{2}{3}r\sin{\phistar}-\frac{17r}{6R}r\sin{\phistar}+4d_1<2r\sin\phistar+\frac{8.1r}{R}r\sin\phistar\\
        \Rightarrow&4d_1<\frac{4}{3}r\sin{\phistar}+(\frac{17r}{6R}+\frac{8.1r}{R})r\sin{\phistar}\Rightarrow\frac{2}{3}d_1<\frac{2}{9}r\sin{\phistar}+(\frac{17r}{36R}+\frac{8.1r}{6R})r\sin{\phistar}\\
        \Rightarrow & 0<\tau_1-\frac{2}{3}d_2\underbracket{<}_{\mathclap{\strut\text{\cref{proposition:casecparametersestimate} \ref{itemcasecparametersestimate7}}}}\frac{2}{3}d_1-\frac{1}{6}r\sin{\phistar}+\frac{16.6r}{3R}r\sin{\phistar}\underbracket{<}_{\strut\mathclap{\frac{2}{3}d_1<\frac{2}{9}r\sin{\phistar}+(\frac{17r}{36R}+\frac{8.1r}{6R})r\sin{\phistar}}}\big(\frac{2}{9}-\frac{1}{6}\big)r\sin{\phistar}+\big(\frac{16.6r}{3R}+\frac{17r}{36R}+\frac{8.1r}{6R}\big)r\sin{\phistar}\underbracket{<}_{\mathclap{\text{\eqref{eqJZHypCond}: }R>1700r}}0.06r\sin\phistar,
    \end{align*}
as claimed. 
We now turn to $\mathcal{II}_1(\mathcal{B}^+_0,2)$ (see \eqref{eq:AformularForLowerBoundOfExpansionForn1largerthan0}). One term in it is
\begin{align*}
        \big[(2n_1+1-\frac{(2n_1+2)\tau_0}{d_1})\Bcal^+_0-\frac{2n_1+2}{d_1}\big]\bigm|_{n_1=2}=&5\Bcal^+_0-\Bcal^+_0\frac{6\tau_0}{d_1}-\frac{6}{d_1}\overbracket{>}^{\mathclap{\text{\cref{lemma:AformularForLowerBoundOfExpansionForn1largerthan0}(1): }\Bcal_0^+<0}}5\Bcal_0^+-\frac{6}{d_1}\\
        \overbracket{>}^{\mathclap{\substack{\text{\cref{lemma:AformularForLowerBoundOfExpansionForn1largerthan0}(1): }\Bcal_0^+\ge{-4}/{3d_0}\\\text{\cref{proposition:casecparametersestimate}\ref{itemcasecparametersestimate3} with }n_1=2\text{, }d_1>\frac{1}{4}r\sin\phistar\strut}}}&-\frac{5\times4}{3d_0}-\frac{6\times4}{r\sin\phistar}\underbracket{>}_{\strut\mathclap{\text{\cref{proposition:casecparametersestimate} \ref{itemcasecparametersestimate5}: }d_0>{r\sin\phistar}/{1.0026}}}\frac{-20\times1.0026}{3r\sin\phistar}-\frac{24}{r\sin\phistar}=\frac{-30.684}{r\sin\phistar}.
    \end{align*}
So, \(\displaystyle\big(\tau_1-\frac{2}{3}d_2\big)\big[(2n_1+1-\frac{(2n_1+2)\tau_0}{d_1})\Bcal^+_0-\frac{2n_1+2}{d_1}\big]\bigm|_{n_1=2}\underbracket{>}_{\quad\mathclap{\tau_1-\frac{2}{3}d_2>0}}\big(\tau_1-\frac{2}{3}d_2\big)\big(\frac{-30.684}{r\sin\phistar}\big)\underbracket{>}_{\bigstrut\mathclap{\tau_1-\frac{2}{3}d_2<0.06r\sin\phistar}}0.06\times(-30.684)>-1.842.\)
The\COMMENT{Edit the next display the way I edited case 2 below} remaining term in $\mathcal{II}_1(\mathcal{B}^+_0,2)$  involves $\tau_0$, and we have\COMMENT{Use overbracket and mathclap below.} 
\[
\tau_0+2d_2\overbracket{<}^{\mathclap{\tau_1>\frac{2}{3}d_2}}\tau_0+3\tau_1=\tau_0+\tau_1+2\tau_1\overbracket{<}^{\mathclap{\text{\cref{remark:deffpclf}: }\tau_1<2d_1}}\tau_1+\tau_0+4d_1\underbracket{<}_{\mathclap{\text{\cref{proposition:casecparametersestimate}\ref{itemcasecparametersestimate2} with }n_1=2}}2r\sin{\phistar}+\frac{8.1r}{R}r\sin{\phistar},
\]
hence
\[
\tau_0<2(r\sin\phistar-d_2)+\frac{8.1r}{R}r\sin{\phistar}\underbracket{<}_{\strut\mathclap{\text{\cref{proposition:casecparametersestimate}\ref{itemcasecparametersestimate1}}}}2\times\frac{17r}{4R}r\sin\phistar+\frac{8.1r}{R}r\sin\phistar=\frac{16.6r}{R}r\sin\phistar,
\]
so
\[
1+\tau_0\Bcal^+_0\overset{\Bcal^+_0<0}>1+\frac{16.6r}{R}r\sin\phistar\Bcal_0^+\overset{\Bcal^+_0\ge\frac{-4}{3d_0}}\ge1-\frac{4\times16.6r}{3R}\frac{r\sin{\phistar}}{d_0}\overbracket{>}^{\strut\mathclap{\text{\cref{proposition:casecparametersestimate}\ref{itemcasecparametersestimate5}}}}1-\frac{4\times16.6r}{3R}1.0026\overbracket{>}^{\mathclap{\eqref{eqJZHypCond}:R>1700r}}0.986
\]
and
\begin{multline*}
\mathcal{II}_2(B^+_0,2)\overbracket{=}^{\mathclap{\text{\cref{def:AformularForLowerBoundOfExpansionForn1largerthan0}}}}4\big[1+(\tau_0-d_1)\Bcal^+_0\big]=4(1+\tau_0\Bcal^+_0)-4d_1\Bcal^+_0>4\times0.986-4d_1\Bcal^+_0
        \\
        \overbracket{>}^{\mathclap{\substack{\Bcal^+_0<0\text{ and \cref{proposition:casecparametersestimate}\ref{itemcasecparametersestimate3}:}\\n_1=2,\ d_1>(r\sin\phistar)/4}}}3.944-\Bcal^+_0r\sin\phistar
        \overbracket{\ge}^{\mathclap{-\Bcal^+_0\ge\frac{1}{d_0}}}3.944+\frac{r\sin\phistar}{d_0}\overbracket{>}^{\strut\mathclap{\text{\cref{proposition:casecparametersestimate}\ref{itemcasecparametersestimate5}}}}3.944+0.9975=4.9415.
\end{multline*}
Thus, in case 1, we conclude
\begin{multline*}
\mathcal{II}_1(\mathcal{B}^+_0,2)+\mathcal{II}_2(\mathcal{B}^+_0,2)
=\overbracket{\mathcal{II}_2(\mathcal{B}^+_0,2)}^{>4.9415}+\overbracket{\big(1+\tau_0\mathcal{B}^+_0\big)}^{>0.986}+\overbracket{\big(\tau_1-\frac{2}{3}d_2\big)\Big[\big(2n_1+1-\frac{(2n_1+2)\tau_0}{d_1}\big)\mathcal{B}^+_0-\frac{2n_1+2}{d_1}\Big]\bigm|_{n_1=2}}^{>-1.842}
\\
>4.9415+0.986-1.842=4.0855>1.8698.
\end{multline*}
\hidesubsection[$\tau_1-\frac{2}{3}d_2\le0$]{Case 2}\strut\COMMENT{Use overbracket below.}
\begin{align*}
    \big([(2n_1+1)-\frac{(2n_1+2)\tau_0}{d_1}]\Bcal^+_0-\frac{2n_1+2}{d_1}\big)&\overbracket=^{\strut\mathclap{n_1=2}}5\Bcal^+_0-\Bcal^+_0\frac{6\tau_0}{d_1}-\frac{6}{d_1}\overbracket{<}^{\mathclap{\substack{\text{\cref{remark:deffpclf}: }\tau_0<2d_1\\\text{and }\Bcal^+_0<0}}}-7\Bcal_0^+-\frac{6}{d_1}\overbracket{\le}^{\mathclap{-\Bcal^+_0\le\frac{4}{3d_0}}}-\frac{28}{3d_0}-\frac{6}{d_1}\\
    &\overbracket{<}^{\mathllap{\text{\cref{proposition:casecparametersestimate}\ref{itemcasecparametersestimate5}: }\frac{1}{d_0}<\frac{1.0026}{r\sin\phistar}\text{ and \ref{itemcasecparametersestimate3} with }n_1=2}}\frac{28\times1.0026}{3r\sin{\phistar}}-\frac{6}{(\frac{1}{2}+\frac{8.1r}{4R}r\sin{\phistar})}\overset{\eqref{eqJZHypCond}:R>1700r}<\frac{-2.61}{r\sin\phistar}<0
\end{align*} Multiplying by $\tau_1-\frac{2}{3}d_2\le0$, we have 
\[\overbracket{\big(\tau_1-\frac{2}{3}d_2\big)}^{\le0}\overbracket{\big([(2n_1+1)-\frac{(2n_1+2)\tau_0}{d_1}]\Bcal^+_0-\frac{2n_1+2}{d_1}\big)}^{<0}\ge0.\] 
Next,\COMMENT{This will look much better with mathclap throughout.\Hrule Add punctuation.\Hrule You are getting carried away with alignment acrobatics and it does not always help. I fixed most of the jumble, maybe you can tidy up the corresponding part in case 1.} 
\(\displaystyle3\tau_0\overbracket{<}^{\mathclap{\text{\cref{def:fpclf}}:\ \tau_0<2d_1}}\tau_0+4d_1<\tau_0+4d_1+\tau_1\overbracket{<}^{\mathclap{\text{\cref{proposition:casecparametersestimate}\ref{itemcasecparametersestimate2}}}}2r\sin\phistar+\frac{8.1r}{R}r\sin\phistar\) implies \(\tau_0<\frac{2}{3}r\sin\phistar+\frac{2.7r}{R}r\sin\phistar\), so
\[
1+\tau_0\Bcal^+_0\overbracket{\ge}^{\mathclap{\Bcal^+_0\ge\frac{-4}{3d_0}}}1-\frac{4\tau_0}{3d_0}>1-\frac{4}{3}\frac{\frac{2}{3}r\sin{\phistar}+\frac{2.7r}{R}r\sin{\phistar}}{d_0}\overbracket{>}^{\bigstrut\mathclap{\text{\cref{proposition:casecparametersestimate}\ref{itemcasecparametersestimate5}}}}1-\frac{4}{3}(\frac{2}{3}+\frac{2.7r}{R})1.0026\overbracket{>}^{\mathclap{\text{\eqref{eqJZHypCond}}: R>1700r}}0.1066,
\]
and, again using \(\tau_0<\frac{2}{3}r\sin\phistar+\frac{2.7r}{R}r\sin\phistar\),
\begin{align*}\mathcal{II}_2(\Bcal^+_0,2)\overbracket{=}^{\mathclap{\text{\cref{def:AformularForLowerBoundOfExpansionForn1largerthan0}}}}&4\big(1+(\tau_0-d_1)\Bcal^+_0\big)\overbracket{>}^{\mathclap{\Bcal^+_0<0}}4\big(1+(\frac{2}{3}r\sin\phistar+\frac{2.7r}{R}r\sin\phistar-d_1)\Bcal^+_0\big)\\
\overbracket{>}^{\mathclap{\substack{\text{\cref{proposition:casecparametersestimate}\ref{itemcasecparametersestimate3}: }\\n_1=2, d_1>(r\sin\phistar)/4\\-\Bcal^+_0>0}}}&4\big(1+(\frac{2}{3}r\sin\phistar+\frac{2.7r}{R}r\sin\phistar-\frac{r\sin\phistar}{4})\Bcal^+_0\big)=4\Big[1+\big(\frac{5}{12}+\frac{2.7r}{R}\big)r\sin\phistar\Bcal^+_0\Big]\\
\overbracket{\ge}^{\mathclap{\Bcal^+_0\ge\frac{-4}{3d_0}}}&4\Big[1+\frac{-4}{3}\big(\frac{5}{12}+\frac{2.7r}{R}\big)\frac{r\sin\phistar}{d_0}\Big]\overbracket{>}^{\mathclap{\text{\cref{proposition:casecparametersestimate}\ref{itemcasecparametersestimate5}}}}4\cdot\big[1-\frac{4\times1.0026}{3}\big(\frac{5}{12}+\frac{2.7r}{R}\big)\big]\overbracket{>}^{\hskip-2em\mathrlap{\text{\eqref{eqJZHypCond}: } R>1700r}}4\times0.4408=1.7632.\\
\end{align*}
Combining these, we get
\begin{multline*}
\mathcal{II}_1(\mathcal{B}^+_0,2)+\mathcal{II}_2(\mathcal{B}^+_0,2)=\overbracket{\mathcal{II}_2(\mathcal{B}^+_0,2)}^{>1.7632}+\overbracket{\big(1+\tau_0\mathcal{B}^+_0\big)}^{>0.1066}+\overbracket{\big(\tau_1-\frac{2}{3}d_2\big)\Big[\big(2n_1+1-\frac{(2n_1+2)\tau_0}{d_1}\big)\mathcal{B}^+_0-\frac{2n_1+2}{d_1}\Big]\bigm|_{n_1=2}}^{\ge0}
\\
>1.7632+0.1066+0=1.8698
\end{multline*}

Thus, $\mathcal{II}_1(\mathcal{B}^+_0,2)+\mathcal{II}_2(\mathcal{B}^+_0,2)>1.8698$ in both cases. This proves the claim.
\end{proof}

\subsection[The expansion in Case (c) of \eqref{eqMainCases} with \texorpdfstring{$n_1=1$}{n1=1}]{The expansion in Case (c) of \texorpdfstring{\eqref{eqMainCases}}{orbit transition} with \texorpdfstring{$n_1=1$}{n1=1}}\label{SBSn1=1}

\begin{definition}\label{def:II1kAndII2k}
    In \cref{def:AformularForLowerBoundOfExpansionForn1largerthan0},  with $\Bcal^+_0=\frac{-k}{d_0}$ for $k\in[1,\frac{4}{3}]$,\COMMENT{Do we need any of the preceding to introduce this notation?} we let $\mathcal{II}_1(\frac{-k}{d_0},1)\nfd \mathcal{II}_1(k)$ and $\mathcal{II}_2(\frac{-k}{d_0},1)\nfd \mathcal{II}_1(k)$.\COMMENT{What follows is not definition; take it out of the definition.} From \eqref{eq:AformularForLowerBoundOfExpansionForn1largerthan0} we have the formula \begin{equation}\label{eq:II1kAndII2k}
        \begin{aligned}
            \mathcal{II}_1(k)&=\mathcal{II}_1(\frac{-k}{d_0},1)=\Big[\big(1+\tau_0\Bcal^+_0\big)+\big(\tau_1-\frac{2}{3}d_2\big)\big[(2n_1+1-\frac{(2n_1+2)\tau_0}{d_1})B^+_0-\frac{2n_1+2}{d_1}\big]\Big]\bigm|_{n_1=1,\Bcal_0^+=\frac{-k}{d_0}}\\
            &=1+\frac{-k\tau_0}{d_0}+\big(\tau_1-\frac{2}{3}d_2\big)\big[(3-\frac{4\tau_0}{d_1})\Bcal^+_0-\frac{4}{d_1}\big]=1-\frac{k\tau_0}{d_0}+\big(\tau_1-\frac{2}{3}d_2\big)\big[(3-\frac{4\tau_0}{d_1})\frac{-k}{d_0}-\frac{4}{d_1}\big]\\
            &=1-k\frac{\tau_0}{d_0}+\big(\frac{\tau_1}{d_0}-\frac{2}{3}\frac{d_2}{d_0}\big)\big(\frac{3d_1}{d_0}-\frac{4\tau_0}{d_0}\big)\big(\frac{-kd_0}{d_1}\big)-\frac{4d_0}{d_1}\big(\frac{\tau_1}{d_0}-\frac{2}{3}\frac{d_2}{d_0}\big).\\
            \mathcal{II}_2(k)&=\mathcal{II}_2(\frac{-k}{d_0},1)=2n_1\big[1+(\tau_0-d_1)\Bcal^+_0\big]\bigm|_{n_1=1,\Bcal^+_0=\frac{-k}{d_0}}=2-\frac{2k\tau_0}{d_0}+\frac{2kd_1}{d_0}.
        \end{aligned}
    \end{equation}
\end{definition}

\begin{definition}[Length fraction functions]\label{def:lengthfractionfunctionsn1eq1}
    In the context of \cref{TEcase_c}, for $x_1\in\MRin$ with $n_1=1$ and the length functions $\tau_0$, $\tau_1$, $d_0$, $d_1$, $d_2$ of \cref{def:fpclf},\COMMENT{Do we need any of the preceding to define these notations??} we define the length fraction functions to be \begin{enumerate}
        \item $h_1\dfn\frac{\tau_0}{d_0}+\frac{\tau_1}{d_0}+\frac{2d_1}{d_0}-2$\label{item01lengthfractionfunctionsn1eq1}
        \item $g_0\dfn d_0-r\sin{\phistar}$\label{item02lengthfractionfunctionsn1eq1}
        \item $g_1\dfn\frac{1}{d_0}-\frac{1}{r\sin{\phistar}}$\label{item03lengthfractionfunctionsn1eq1}
        \item $g_2\dfn d_2-r\sin{\phistar}$\label{item04lengthfractionfunctionsn1eq1}
        \item $g_3\dfn\frac{1}{d_2}-\frac{1}{r\sin{\phistar}}$.\label{item05lengthfractionfunctionsn1eq1}
    \end{enumerate}
 These are defined for $x_1$ in the quadrilateral with $n_1=1$ in \cref{fig:x1_n1cases_caseC}, also shown as \textcolor{blue}{$M^{\text{\upshape{in}}}_{R,1}$} in \cref{fig:x1_n1equal1case}.\COMMENT{It sems that instead, they are defined whenever the denominators are not zero. We only USE them under some assumptions; those should be stated when needed.} 
\end{definition}
\begin{lemma}[Estimates for length fraction functions]\label{lemma:h1g1g0g1g2g3lengthfractionestimate}\COMMENT{\color{red}Here and in many other places you put \textbackslash label before the optional argument of the environment. It must always come after. There are more instances than I can fix; please find and fix them all, starting with this one.\wentao{All fixed now.}}
    In the context of \cref{TEcase_c}, under condition \eqref{eqJZHypCond}, for the functions in \cref{def:lengthfractionfunctionsn1eq1} defined for $x_1$ in the $n_1=1$ region in \cref{fig:x1_n1cases_caseC} or equivalently the \textcolor{blue}{$M^{\text{\upshape{in}}}_{R,1}$} component in \cref{fig:x1_n1equal1case},\COMMENT{Edit to remove any standing assumptions.} we have \begin{enumerate}
        \item\label{item01h1g1g0g1g2g3lengthfractionestimate} $-8.53r/R<h_1<16.65r/R$,
        \item\label{item02h1g1g0g1g2g3lengthfractionestimate} $-\frac{17r}{4R}r\sin\phistar<g_0<\frac{17r}{4R}r\sin\phistar$,
        \item\label{item03h1g1g0g1g2g3lengthfractionestimate} $-\frac{17r}{4R}\frac{1}{0.9975r\sin\phistar}<g_1<\frac{17r}{4R}\frac{1}{0.9975r\sin\phistar}$,
        \item\label{item04h1g1g0g1g2g3lengthfractionestimate} $-\frac{17r}{4R}r\sin\phistar<g_2<\frac{17r}{4R}r\sin\phistar$,
        \item\label{item05h1g1g0g1g2g3lengthfractionestimate} $-\frac{17r}{4R}\frac{1}{0.9975r\sin\phistar}<g_3<\frac{17r}{4R}\frac{1}{0.9975r\sin\phistar}$
        \item\label{item06h1g1g0g1g2g3lengthfractionestimate} $-\frac{17r}{2R}r\sin\phistar<d_0-d_2=g_0-g_2<\frac{17r}{2R}r\sin\phistar$
    \end{enumerate}
\end{lemma}
\begin{proof} We provide item-by-item proofs.\COMMENT{Lots of mathclap! Always with overbracket or underbracket. Sometimes with strut or bigstrut\Hrule\small\color{red}Replace enumerate by hidesubsection as was done in the proof of \cref{proposition:QUpperLowerRegionsLowerBoundsChatCtilde}} 
\begin{enumerate}
    \item $h_1\overbracket{=}^{\quad\mathllap{\text{\cref{def:lengthfractionfunctionsn1eq1}\eqref{item01lengthfractionfunctionsn1eq1}}}}\frac{\overbracket{\tau_0+\tau_1+2d_1}^{\mathrlap{\substack{\text{\cref{proposition:casecparametersestimate}\ref{itemcasecparametersestimate2} with }n_1=1: \\\in(2r\sin\phistar,\quad 2r\sin\phistar+(\frac{8.1r}{R})r\sin\phistar)}}\qquad\qquad\qquad\qquad}}{d_0}-2\in\big(\frac{2r\sin\phistar}{d_0}-2,\frac{2r\sin\phistar}{d_0}+\frac{8.1r}{R}\frac{r\sin\phistar}{d_0}-2\big),$
    
    Since 
    \begin{align*}
        \frac{2r\sin\phistar}{d_0}-2=&2\frac{\overbracket{r\sin\phistar-d_0}^{\mathllap{\substack{\text{\cref{proposition:casecparametersestimate}\ref{itemcasecparametersestimate1}: }\\ r\sin\phistar-d_0>-\frac{17r}{4R}r\sin\phistar}}}}{d_0}>-\frac{17r}{2R}\frac{r\sin\phistar}{d_0}\overbracket{>}^{\qquad\qquad\mathllap{\text{\cref{proposition:casecparametersestimate}\ref{itemcasecparametersestimate5}}}}-\frac{17r}{2R}\frac{1}{0.9975}>-8.53r/R,\\
        \frac{2r\sin\phistar}{d_0}+\frac{8.1r}{R}\frac{r\sin\phistar}{d_0}-2=&2\frac{\overbracket{r\sin\phistar-d_0}^{\mathllap{\substack{\text{\cref{proposition:casecparametersestimate}\ref{itemcasecparametersestimate1}: }\\ r\sin\phistar-d_0<\frac{17r}{4R}r\sin\phistar}}}}{d_0}+\frac{8.1r}{R}\frac{r\sin\phistar}{d_0}<\frac{17r}{2R}\frac{r\sin\phistar}{d_0}+\frac{8.1r}{R}\frac{r\sin\phistar}{d_0}\\
        \overbracket{<}^{\mathllap{\text{\cref{proposition:casecparametersestimate}\ref{itemcasecparametersestimate5}}}}&\frac{17r}{2\times0.9975R}+\frac{8.1r}{0.9975R}<16.65r/R.
    \end{align*}
    Therefore, $-8.53r/R<h<16.65r/R$.
    \item It follows\COMMENT{Fix!} \cref{proposition:casecparametersestimate}\ref{itemcasecparametersestimate1} that $|g_0|=|d_0-r\sin\phistar|<\frac{17r}{4R}r\sin\phistar$. Likewise for (4).
    \item $|g_1|\overbracket{=}^{\mathclap{\text{\cref{def:lengthfractionfunctionsn1eq1}\eqref{item01lengthfractionfunctionsn1eq1}}\quad}}\frac{|d_0-r\sin\phistar|}{d_0r\sin\phistar}\overbracket{<}^{\mathclap{\quad\text{\cref{proposition:casecparametersestimate}\ref{itemcasecparametersestimate1}}}}\frac{17r}{4R}\frac{r\sin\phistar}{d_0r\sin\phistar}=\frac{17r}{4R}\frac{1}{d_0}<\frac{17r}{4R}\frac{1}{0.9975r\sin\phistar}$. Likewise for (5).
    \item The proof is the same as (2).\COMMENT{Omit}
    \item The proof is the same as (3).\COMMENT{Omit}
    \item $g_0-g_2\overbracket{=}^{\mathclap{\text{\cref{def:lengthfractionfunctionsn1eq1}\eqref{item02lengthfractionfunctionsn1eq1}\eqref{item04lengthfractionfunctionsn1eq1}}}}d_0-d_2$ and by the triangle inequality 
    \[|d_0-d_2|\le|d_0-r\sin\phistar|+|r\sin\phistar-d_2|\overbracket{<}^{\mathclap{\text{\cref{proposition:casecparametersestimate}\ref{itemcasecparametersestimate1}}}}\frac{17r}{4R}r\sin\phistar+\frac{17r}{4R}r\sin\phistar=\frac{17r}{2R}r\sin\phistar.\qedhere\]
\end{enumerate}
\end{proof}
\begin{proposition}\label{proposition:II1II2C1C2}
    In the context of \cref{TEcase_c,lemma:AformularForLowerBoundOfExpansionForn1largerthan0}, under condition \eqref{eqJZHypCond}, for the $\mathcal{II}_1(k)$, $\mathcal{II}_2(k)$ in \cref{def:II1kAndII2k} with $k\in[1,4/3]$ and $h_1$, $g_0$, $g_1$, $g_2$, $g_3$ from \cref{def:lengthfractionfunctionsn1eq1} defined for $x_1$ with $n_1=1$ in \cref{fig:x1_n1cases_caseC} or equivalently in the \textcolor{blue}{$M^{\text{\upshape{in}}}_{R,1}$} component in \cref{fig:x1_n1equal1case}, we have\COMMENT{Needs editing, but only after we decided about any standing assumptions.}
    \begin{equation}\label{eq:II1II2C1C2}
        \begin{aligned}
            \mathcal{II}_1(k)=&\frac{3+16k}{3}+\big(\frac{8-16k}{3}\big)\frac{r\sin{\phistar}}{d_1}+\frac{2kd_1}{r\sin{\phistar}}+\Big[\frac{(32k-12)r\sin{\phistar}}{3d_1}-10k\Big]\big(\frac{\tau_1}{d_0}\big)\\
            &-\frac{4kr\sin{\phistar}}{d_1}\big(\frac{\tau_1}{d_0}\big)^2+\mathcal{C}_1(k)\\
            \mathcal{II}_2(k)=&2-4k+6k\frac{d_1}{r\sin{\phistar}}+2k\frac{\tau_1}{d_0}+\mathcal{C}_2(k)\\
            \mathcal{II}_1(k)+\mathcal{II}_2(k)=&\frac{9+4k}{3}+\big(\frac{8-16k}{3}\big)\frac{r\sin\phistar}{d_1}+\frac{8kd_1}{r\sin\phistar}+\Big[\frac{(32k-12)r\sin{\phistar}}{3d_1}-8k\Big]\big(\frac{\tau_1}{d_0}\big)\\
            &-\frac{4kr\sin{\phistar}}{d_1}\big(\frac{\tau_1}{d_0}\big)^2+\mathcal{C}_1(k)+\mathcal{C}_2(k)
        \end{aligned}
    \end{equation}
    where \begin{equation}\label{eq:C1C2}
    \begin{aligned}
        \mathcal{C}_1(k)=&2kd_1g_1-\frac{2(g_2-g_0)}{3d_0}\frac{8kg_0}{d_1}+(\frac{\tau_1}{d_0}-\frac{2}{3})8k\frac{g_0}{d_1}-\frac{2(g_2-g_0)}{3d_0}(-11k+8k\frac{r\sin{\phistar}}{d_1})\\
        &+(-4k\frac{\tau_1}{d_0})[-\frac{r\sin{\phistar}}{d_1}\frac{2(g_2-g_0)}{3d_0}+\frac{g_0}{d_1}(\frac{\tau_1}{d_0}-\frac{2}{3})-\frac{g_0}{d_1}\frac{2(g_2-g_0)}{3d_0}]-\frac{4g_0}{d_1}(\frac{\tau_1}{d_0}-\frac{2}{3})\\
        &+\frac{8(g_2-g_0)}{3d_1}-kh_1+4h_1(\frac{\tau_1}{d_0}-\frac{2d_2}{3d_0})\frac{kd_0}{d_1}\\
        \mathcal{C}_2(k)=&6kd_1g_1-2kh_1
    \end{aligned}
    \end{equation}
\end{proposition}
\begin{remark}\label{remark:II1II2Quadratic}
    For fixed $k\in[1,4/3]$, $\frac{d_1}{r\sin\phistar}=\const$, we define a quadratic function by\COMMENT{Use the same Q as previously.} \begin{equation}\label{eq:II1II2Quadratic}
        \mathtt{Q}(t)\dfn-\frac{4kr\sin\phistar}{d_1}t^2+\big[\frac{(32k-12)r\sin\phistar}{3d_1}-8k\big]t+\big(\frac{8-16k}{3}\big)\frac{r\sin\phistar}{d_1}+\frac{8kd_1}{r\sin{\phistar}}+\frac{4k+9}{3}.
    \end{equation}
\emph{\bfseries Alternate version; choose one:}    For fixed $k\in[1,4/3]$, define a quadratic function by \begin{equation}\label{eq:II1II2Quadratic2}
        \mathfrak{Q}(t)\dfn-12kDt^2+\big[(32k-12)D-8k\big]t+(8-16k)D+\frac{8k}{3D}+\frac{4k+9}{3},\quad\text{where}\quad D\dfn\frac{r\sin\phistar}{3d_1}.
    \end{equation}
    Then in \cref{proposition:II1II2C1C2}, \begin{equation}\label{eq:II1II2QuadraticAndC1C2}\mathcal{II}_1(k)+\mathcal{II}_2(k)=\mathtt{Q}\big(\frac{\tau_1}{d_0}\big)+\mathcal{C}_1(k)+\mathcal{C}_2(k).\end{equation} 
    And we will see in \cref{proposition:C1C2lowerbound} that the lower bound of $\mathcal{C}_1(k)$ and $\mathcal{C}_2(k)$ can be arbitrarily close to $0$ as $R/r$ is large enough, and the same holds for the upper bound. 
\end{remark}
\begin{proof}
    We first compute $\mathcal{II}_1(k)$.\COMMENT{Lots of mathclap!\Hrule\small Always with overbracket or underbracket. Almost always without \textbackslash quad or 
    \textbackslash qquad} 
    \begin{align*}
        \mathcal{II}_1(k)\overbracket{=}^{\mathclap{\text{\eqref{eq:II1kAndII2k}}}}&1-k\frac{\tau_0}{d_0}+\big(\frac{\tau_1}{d_0}-\frac{2}{3}\frac{d_2}{d_0}\big)\big(\frac{3d_1}{d_0}-\frac{4\tau_0}{d_0}\big)\big(\frac{-kd_0}{d_1}\big)-\frac{4d_0}{d_1}\big(\frac{\tau_1}{d_0}-\frac{2}{3}\frac{d_2}{d_0}\big)\\
        \overbracket{=}^{\mathclap{\substack{\text{\cref{def:lengthfractionfunctionsn1eq1} \eqref{item01lengthfractionfunctionsn1eq1}:}\\\frac{\tau_0}{d_0}=2+h_1-\frac{\tau_1}{d_0}-\frac{2d_1}{d_0}}}}&1-k\big(2+h_1-\frac{\tau_1}{d_0}-\frac{2d_1}{d_0}\big)+\big(\frac{\tau_1}{d_0}-\frac{2}{3}\frac{d_2}{d_0}\big)\big[\frac{3d_1}{d_0}-4(2+h_1-\frac{\tau_1}{d_0}-\frac{2d_1}{d_0})\big]\big(\frac{-kd_0}{d_1}\big)-\frac{4d_0}{d_1}\big(\frac{\tau_1}{d_0}-\frac{2}{3}\frac{d_2}{d_0}\big)\\
        =&1-k\big(2+h_1-\frac{\tau_1}{d_0}-\frac{2d_1}{d_0}\big)+\big(\frac{\tau_1}{d_0}-\frac{2}{3}\frac{d_2}{d_0}\big)\big(\frac{3d_1}{d_0}-8-4h_1+\frac{4\tau_1}{d_0}+\frac{8d_1}{d_0}\big)\big(\frac{-kd_0}{d_1}\big)-\frac{4d_0}{d_1}\big(\frac{\tau_1}{d_0}-\frac{2}{3}\frac{d_2}{d_0}\big)\\
        =&1-k\big(2-\frac{\tau_1}{d_0}-\frac{2d_1}{d_0}\big)-kh_1+\big(\frac{\tau_1}{d_0}-\frac{2d_2}{3d_0}\big)\big(\frac{11d_1}{d_0}-8+\frac{4\tau_1}{d_0}\big)\big(\frac{-kd_0}{d_1}\big)\\
        &+4h_1\big(\frac{\tau_1}{d_0}-\frac{2d_2}{3d_0}\big)\frac{kd_0}{d_1}-\frac{4d_0}{d_1}\big(\frac{\tau_1}{d_0}-\frac{2}{3}\frac{d_2}{d_0}\big),\\
    \end{align*}
    Using the functions from \cref{def:lengthfractionfunctionsn1eq1}, we split and recombine this:
    \begin{align*}
        \mathcal{II}_1(k)=&1-2k+k\frac{\tau_1}{d_0}+2k\frac{d_1}{d_0}+\big(\frac{\tau_1}{d_0}-\frac{2d_2}{3d_0}\big)\big(-11k+8k\frac{d_0}{d_1}\big)+\big(\frac{\tau_1}{d_0}-\frac{2d_2}{3d_0}\big)\big(\frac{-4k\tau_1}{d_0}\big)\big(\frac{d_0}{d_1}\big)-\frac{4d_0}{d_1}\big(\frac{\tau_1}{d_0}-\frac{2}{3}\frac{d_2}{d_0}\big)\\
        &-kh_1+4h_1\big(\frac{\tau_1}{d_0}-\frac{2d_2}{3d_0}\big)
        \\
        =&1-2k+k\frac{\tau_1}{d_0}+2kd_1\big(\overbracket{\frac{1}{r\sin{\phistar}}+g_1}^{\mathclap{\text{\cref{def:lengthfractionfunctionsn1eq1}\eqref{item03lengthfractionfunctionsn1eq1}: }=\frac{1}{d_0}}}\big)+\big[\frac{\tau_1}{d_0}-\frac{2}{3}-\frac{2\overbracket{(g_2-g_0)}^{\mathclap{\substack{\text{\cref{def:lengthfractionfunctionsn1eq1}\eqref{item02lengthfractionfunctionsn1eq1}, \eqref{item04lengthfractionfunctionsn1eq1}: }\\ g_2-g_0=d_2-d_0}}}}{3d_0}\big]\big(-11k+\overbracket{8k\frac{r\sin{\phistar}}{d_1}+8k\frac{g_0}{d_1}}^{\mathclap{\text{\cref{def:lengthfractionfunctionsn1eq1}\eqref{item02lengthfractionfunctionsn1eq1}: }=8k\frac{d_0}{d_1}}}\big)
        \\
        &+\big[\frac{\tau_1}{d_0}-\frac{2}{3}-\frac{2\overbracket{(g_2-g_0)}^{\mathclap{\text{\cref{def:lengthfractionfunctionsn1eq1}\eqref{item02lengthfractionfunctionsn1eq1}\eqref{item04lengthfractionfunctionsn1eq1}: } g_2-g_0=d_2-d_0\quad}}}{3d_0}\big]\big(-4k\frac{\tau_1}{d_0}\big)\big(\overbracket{\frac{r\sin{\phistar}}{d_1}+\frac{g_0}{d_1}}^{\mathclap{\text{\cref{def:lengthfractionfunctionsn1eq1}\eqref{item02lengthfractionfunctionsn1eq1}: }=\frac{d_0}{d_1}}}\big)-\big(\frac{4g_0}{d_1}+\frac{4r\sin{\phistar}}{d_1}\big)\big(\frac{\tau_1}{d_0}-\frac{2}{3}-\frac{2}{3}\frac{\overbracket{g_2-g_0}^{\mathclap{\text{\cref{def:lengthfractionfunctionsn1eq1}\eqref{item02lengthfractionfunctionsn1eq1}, \eqref{item04lengthfractionfunctionsn1eq1}: } g_2-g_0=d_2-d_0}}}{d_0}\big)
        \\
        &-kh_1+4h_1\big(\frac{\tau_1}{d_0}-\frac{2d_2}{3d_0}\big)
        \\
        \overbracket{=}^{\mathclap{\substack{\text{Split terms into terms}\\\text{with and without}\\ g_0,g_1,g_2,h_1}}\qquad}&1-2k+k\frac{\tau_1}{d_0}+\frac{2kd_1}{r\sin{\phistar}}+\big(\frac{\tau_1}{d_0}-\frac{2}{3}\big)\big(-11k+8k\frac{r\sin{\phistar}}{d_1}\big)+\big(\frac{\tau_1}{d_0}-\frac{2}{3}\big)\big(-4k\frac{\tau_1}{d_0}\big)\frac{r\sin{\phistar}}{d_1}-\frac{4r\sin{\phistar}}{d_1}\big(\frac{\tau_1}{d_0}-\frac{2}{3}\big)\\
        &+2kd_1g_1-\frac{2(g_2-g_0)}{3d_0}\big(-11k+8k\frac{r\sin\phistar}{d_1}+8k\frac{g_0}{d_1}\big)+\big(\frac{\tau_1}{d_0}-\frac{2}{3}\big)8k\frac{g_0}{d_1}+\frac{2(g_2-g_0)}{3d_0}\frac{4k\tau_1}{d_0}\big(\frac{r\sin\phistar}{d_1}+\frac{g_0}{d_1}\big)\\
        &+\big(\frac{\tau_1}{d_0}-\frac{2}{3}\big)\big(-4k\frac{\tau_1}{d_0}\big)\frac{g_0}{d_1}-\frac{4g_0}{d_1}\big(\frac{\tau_1}{d_0}-\frac{2}{3}-\frac{2}{3}\frac{g_2-g_0}{d_0}\big)+\frac{4r\sin\phistar}{d_1}\frac{2(g_2-g_0)}{3d_0}\\
        &-kh_1+4h_1\big(\frac{\tau_1}{d_0}-\frac{2d_2}{3d_0}\big)\\
        =&1-2k+\frac{2kd_1}{r\sin\phistar}+\frac{22k}{3}-\frac{16kr\sin\phistar}{3d_1}+\frac{8r\sin\phistar}{3d_1}+\big(k-11k+8k\frac{r\sin\phistar}{d_1}+\frac{8kr\sin\phistar}{3d_1}-\frac{4r\sin\phistar}{d_1}\big)\big(\frac{\tau_1}{d_0}\big)\\
        &-\frac{4kr\sin\phistar}{d_1}\big(\frac{\tau_1}{d_0}\big)^2+\mathcal{C}_1(k)\\
        =&\frac{3+16k}{3}+\big(\frac{8-16k}{3}\big)\frac{r\sin{\phistar}}{d_1}+\frac{2kd_1}{r\sin{\phistar}}+\Big[\frac{(32k-12)r\sin{\phistar}}{3d_1}-10k\Big]\frac{\tau_1}{d_0}-\frac{4kr\sin{\phistar}}{d_1}\big(\frac{\tau_1}{d_0}\big)^2+\mathcal{C}_1(k),
    \end{align*} where\COMMENT{\small Isn't this definition in \cref{eq:C1C2}? Then it should be implemented the FIRST (not last) time this expression appears above. This means that ALL occurrences above of the formula for $\mathcal{C}_1(k)$ should be replaced by ``$\mathcal{C}_1(k)$''---with reference to \cref{eq:C1C2} the first time, of course. This will reduce clutter a lot!}
    \begin{multline*}
    \mathcal{C}_1(k)=2kd_1g_1-\frac{2(g_2-g_0)}{3d_0}\big(-11k+8k\frac{r\sin\phistar}{d_1}+8k\frac{g_0}{d_1}\big)+\big(\frac{\tau_1}{d_0}-\frac{2}{3}\big)8k\frac{g_0}{d_1}+\frac{2(g_2-g_0)}{3d_0}\frac{4k\tau_1}{d_0}\big(\frac{r\sin\phistar}{d_1}+\frac{g_0}{d_1}\big)
    \\
    +\big(\frac{\tau_1}{d_0}-\frac{2}{3}\big)\big(-4k\frac{\tau_1}{d_0}\big)\frac{g_0}{d_1}-\frac{4g_0}{d_1}\big(\frac{\tau_1}{d_0}-\frac{2}{3}-\frac{2}{3}\frac{g_2-g_0}{d_0}\big)+\frac{4r\sin\phistar}{d_1}\frac{2(g_2-g_0)}{3d_0}-kh_1+4h_1\big(\frac{\tau_1}{d_0}-\frac{2d_2}{3d_0}\big).
        \end{multline*}

Similarly by introducing functions in \cref{def:lengthfractionfunctionsn1eq1}, we also have 
\begin{align*}
    \mathcal{II}_2(k)&\overbracket{=}^{\mathclap{\text{\eqref{eq:II1kAndII2k}}}}2-\frac{2k\tau_0}{d_0}+\frac{2kd_1}{d_0}\\
    &\overbracket{=}^{\mathclap{\substack{\text{\cref{def:lengthfractionfunctionsn1eq1} \eqref{item01lengthfractionfunctionsn1eq1}: }\\\frac{\tau_0}{d_0}=2+h_1-\frac{\tau_1}{d_0}-\frac{2d_1}{d_0}}}}2-2k\big(2-\frac{\tau_1}{d_0}-\frac{2d_1}{d_0}+h_1\big)+\frac{2kd_1}{d_0}=2-4k+\frac{6kd_1}{d_0}+\frac{2k\tau_1}{d_0}-2kh_1
    \\
    &\overbracket{=}^{\mathclap{\substack{\text{\cref{def:lengthfractionfunctionsn1eq1}\eqref{item03lengthfractionfunctionsn1eq1}}:\\\frac{1}{d_0}=\frac{1}{r\sin\phistar}+g_1}}}2-4k+6kd_1\big(\frac{1}{r\sin{\phistar}}+g_1\big)+\frac{2k\tau_1}{d_0}-2kh_1\\
    &=2-4k+6k\frac{d_1}{r\sin\phistar}+2k\frac{\tau_1}{d_0}+\underbracket{6kd_1g_1-2kh_1}_{\mathcal{C}_2(k)}=2-4k+6k\frac{d_1}{r\sin{\phistar}}+2k\frac{\tau_1}{d_0}+\mathcal{C}_2(k).\qedhere
\end{align*}
\end{proof}
\begin{proposition}\label{proposition:C1C2lowerbound}
    For $k\in[1,4/3]$, $\mathcal{C}_1(k)$ and $\mathcal{C}_2(k)$ from \eqref{eq:C1C2} 
    satisfy 
    \[\mathcal{C}_1(k)+\mathcal{C}_2(k)>-389\frac{r}{R}-483\big(\frac{r}{R}\big)^2.\] 
\end{proposition}
\begin{proof}
    Since  $\frac{1}{3}<\frac{d_1}{r\sin\phistar}<1$ by \cref{proposition:QUpperLowerRegionsLowerBoundsChatCtilde} \eqref{item01Qlowerboundn1eq1}, when $x_1$ is in the $n_1=1$ region of \cref{fig:x1_n1cases_caseC}, we have\COMMENT{\color{red}Bad \TeX\ below and probably throughout the paper: Note that all the ``='' are too close to the 2 or 8 to the right. This is because of ``=\&'' This should ALWAYS be ``\&=''---and then it should always be ``\&\textbackslash quad+'' to move the ``+'' to the right!}
    \begin{equation}\label{eq:C1C2reformed}
    \begin{aligned}
    \mathcal{C}_1(k)+\mathcal{C}_2(k)\overbracket{=}^{\mathclap{\text{\eqref{eq:C1C2}}}}&2kd_1g_1-\frac{2(g_2-g_0)}{3d_0}\big(-11k+8k\frac{r\sin\phistar}{d_1}+8k\frac{g_0}{d_1}\big)+\big(\frac{\tau_1}{d_0}-\frac{2}{3}\big)8k\frac{g_0}{d_1}+\frac{2(g_2-g_0)}{3d_0}\frac{4k\tau_1}{d_0}\big(\frac{r\sin\phistar}{d_1}+\frac{g_0}{d_1}\big)
    \\
        &+\big(\frac{\tau_1}{d_0}-\frac{2}{3}\big)\big(-4k\frac{\tau_1}{d_0}\big)\frac{g_0}{d_1}-\frac{4g_0}{d_1}\big(\frac{\tau_1}{d_0}-\frac{2}{3}-\frac{2}{3}\frac{g_2-g_0}{d_0}\big)+\frac{4r\sin\phistar}{d_1}\frac{2(g_2-g_0)}{3d_0}-kh_1+4h_1\big(\frac{\tau_1}{d_0}-\frac{2d_2}{3d_0}\big)
        \\
        &+6kd_1g_1-2kh_1
        \\
        =&8kd_1g_1+h_1\big(-3k+\frac{\tau_1}{d_0}-\frac{2d_2}{3d_0}\big)+\frac{2(g_2-g_0)}{3d_0}\big(\frac{4r\sin\phistar}{d_1}+\frac{4k\tau_1}{d_0}\frac{r\sin\phistar}{d_1}+\frac{4k\tau_1}{d_0}\frac{g_0}{d_1}\big)+\frac{2(g_2-g_0)}{3d_0}\frac{4g_0}{d_1}
        \\
        &-\frac{2(g_2-g_0)}{3d_0}\big(-11k+8k\frac{r\sin\phistar}{d_1}+8k\frac{g_0}{d_1}\big)
        -\frac{4g_0}{d_1}\big(\frac{\tau_1}{d_0}-\frac{2}{3}\big)+\big(\frac{\tau_1}{d_0}-\frac{2}{3}\big)\big(-4k\frac{\tau_1}{d_0}\big)\frac{g_0}{d_1}+\big(\frac{\tau_1}{d_0}-\frac{2}{3}\big)8k\frac{g_0}{d_1}
        \\
        =&8kd_1g_1+h_1\big(-3k+\frac{\tau_1}{d_0}-\frac{2d_2}{3d_0}\big)+\big(\frac{\tau_1}{d_0}-\frac{2}{3}\big)\frac{g_0}{d_1}\big(8k-4-4k\frac{\tau_1}{d_0}\big)
        \\
        &+\frac{2(g_2-g_0)}{3d_0}\big(\frac{4r\sin\phistar}{d_1}+\frac{4k\tau_1}{d_0}\frac{r\sin\phistar}{d_1}+\frac{4k\tau_1+4d_0}{d_0}\frac{g_0}{d_1}\big)-\frac{2(g_2-g_0)}{3d_0}\big(-11k+8k\frac{r\sin\phistar}{d_1}+8k\frac{g_0}{d_1}\big)\\
        =&8kd_1g_1+h_1\big(-3k+\frac{\tau_1}{d_0}-\frac{2d_2}{3d_0}\big)+\big(\frac{\tau_1}{d_0}-\frac{2}{3}\big)\frac{g_0}{d_1}\big(8k-4-4k\frac{\tau_1}{d_0}\big)\\
        &+\frac{2(g_2-g_0)}{3d_0}\Big[11k-\frac{8kr\sin\phistar}{d_1}+\frac{4r\sin\phistar}{d_1}+\frac{4k\tau_1}{d_0}\frac{r\sin\phistar}{d_1}+\frac{4k\tau_1+(4-8k)d_0}{d_0}\frac{g_0}{d_1}\Big].
    \end{aligned}
    \end{equation}
Note that\COMMENT{Lots of mathclap with overbracket or underbracket. With strut or bigstrut if needed.} 
    \begin{equation}\label{eq:g0d1ratioestimate}
        -\frac{51r}{4R}\overbracket{<}^{d_1>(r\sin\phistar)/3}-\frac{17r}{4R}\frac{r\sin\phistar}{d_1}\overbracket{<}^{\text{\cref{lemma:h1g1g0g1g2g3lengthfractionestimate}\eqref{item02h1g1g0g1g2g3lengthfractionestimate}}}\frac{g_0}{d_1}\overbracket{<}^{\text{\cref{lemma:h1g1g0g1g2g3lengthfractionestimate}\eqref{item02h1g1g0g1g2g3lengthfractionestimate}}}\frac{17r}{4R}\frac{r\sin\phistar}{d_1}\overbracket{<}^{d_1>(r\sin\phistar)/3}\frac{51r}{4R},
    \end{equation} 
    and \begin{equation}\label{eq:g2minusg0overd0estimate}
        -\frac{17r}{3R}\cdot1.0026\overbracket{<}^{\qquad\mathllap{\text{\cref{proposition:casecparametersestimate}\ref{itemcasecparametersestimate5}}}}-\frac{2}{3d_0}\frac{17r}{2R}r\sin\phistar\overbracket{<}^{\qquad\mathllap{\text{\cref{lemma:h1g1g0g1g2g3lengthfractionestimate}\eqref{item06h1g1g0g1g2g3lengthfractionestimate}}}}\frac{2(g_2-g_0)}{3d_0}\overbracket{<}^{\mathrlap{\text{\cref{lemma:h1g1g0g1g2g3lengthfractionestimate}\eqref{item06h1g1g0g1g2g3lengthfractionestimate}}}\qquad\qquad}\frac{2}{3d_0}\frac{17r}{2R}r\sin\phistar\overbracket{<}^{\qquad\qquad\mathllap{\text{\cref{proposition:casecparametersestimate}\ref{itemcasecparametersestimate5}}}}\frac{17r}{3R}\cdot1.0026.
    \end{equation}
Therefore,\COMMENT{Below are examples of how to get rid of unneeded ``alignments''} $\displaystyle -\frac{20}{3}\overset{k\le4/3}\le4-8k\overset{\tau_1>0}<\frac{4k\tau_1+(4-8k)d_0}{d_0}\overbracket{<}^{\qquad\qquad\mathllap{\substack{\text{\cref{remark:deffpclf}}:\\\tau_1<2d_0}}}2\times4k+(4-8k)=4$. This and \eqref{eq:g0d1ratioestimate} imply
\begin{equation}\label{eq:4thtermg0d1ration}
    \begin{aligned}
    -\frac{85r}{R}=-\frac{20\times51r}{3\times4R}<\frac{4k\tau_1+(4-8k)d_0}{d_0}\frac{g_0}{d_1}<\frac{20\times51r}{3\times4R}=\frac{85r}{R}.
\end{aligned}
\end{equation}
This and $\displaystyle11k-\frac{8kr\sin\phistar}{d_1}+\frac{4r\sin\phistar}{d_1}+\frac{4k\tau_1}{d_0}\frac{r\sin\phistar}{d_1}=\underbracket{\overbracket{11k+\big(4+\frac{4k\tau_1}{d_0}-8k\big)\frac{r\sin\phistar}{d_1}}^{\overbracket{\scriptstyle<}^{\mathclap{\text{\cref{remark:deffpclf}}:\ \tau_1<2d_0}}11k+\frac{4r\sin\phistar}{d_1}\overbracket{\scriptstyle<}^{\mathclap{d_1>\frac{r\sin\phistar}{3}}}11k+12}}_{\underbracket{\scriptstyle>}_{\mathclap{\tau_1>0}}11k+(4-8k)\frac{4r\sin\phistar}{d_1}
    \underbracket{\scriptstyle>}_{\mathclap{4-8k<0\text{ and }d_1<r\sin\phistar}}11k+4-8k=3k+4}
    $ in turn give
\begin{equation}\label{eq:4thtermanalysis}
\begin{aligned}
    0\underbracket{<}_{\mathclap{\qquad\eqref{eqJZHypCond}:\ R>1700r,\ k\ge1}}3k+4-\frac{85r}{R}<\Big[11k-\frac{8kr\sin\phistar}{d_1}+\frac{4r\sin\phistar}{d_1}+&\frac{4k\tau_1}{d_0}\frac{r\sin\phistar}{d_1}+\frac{4k\tau_1+(4-8k)d_0}{d_0}\frac{g_0}{d_1}\Big]<11k+12+\frac{85r}{R}.\\
\end{aligned}
\end{equation}
Then, assuming \eqref{eqJZHypCond} (i.e., $R>1700r$), by \eqref{eq:g2minusg0overd0estimate}, \eqref{eq:4thtermg0d1ration}, we find\COMMENT{Use multline below}
\begin{equation}\label{eq:4thterminC1C2}
    \begin{aligned}
        &-\frac{17\times1.0026r}{3R}\big(11k+12+\frac{85r}{R}\big)\\
        \overbracket{<}^{\mathllap{\eqref{eq:g2minusg0overd0estimate}, \eqref{eq:4thtermanalysis}}}&\frac{2(g_2-g_0)}{3d_0}\Big[11k-\frac{8kr\sin\phistar}{d_1}+\frac{4r\sin\phistar}{d_1}+\frac{4k\tau_1}{d_0}\frac{r\sin\phistar}{d_1}+\frac{4k\tau_1+(4-8k)d_0}{d_0}\frac{g_0}{d_1}\Big]\\\overbracket{<}^{\mathllap{\eqref{eq:g2minusg0overd0estimate}, \eqref{eq:4thtermanalysis}}}&\frac{17\times1.0026r}{3R}\big(11k+12+\frac{85r}{R}\big).
    \end{aligned}
\end{equation}Also,\COMMENT{Hard-to-read jumbles follow. Maybe OK...\Hrule More mathclap}
\begin{equation}\label{eq:3rdterminC1C2}
\begin{aligned}
&\left.
\begin{aligned}
&-4\overbracket{<}^{\mathllap{\text{\cref{remark:deffpclf}: }\tau_1<2d_0}}\big(8k-4-4k\frac{\tau_1}{d_0}\big)\overbracket{<}^{\qquad\qquad\mathllap{\tau_1>0}}8k-4\overset{k\le4/3}{\le}\frac{20}{3}\\
&-\frac{2}{3}\overbracket{<}^{\qquad\mathllap{\tau_1>0}}\frac{\tau_1}{d_0}-\frac{2}{3}\overbracket{<}^{\qquad\quad\qquad\qquad\mathllap{\text{\cref{remark:deffpclf}: }\tau_1<2d_0}}\frac{4}{3}
\end{aligned}
\right\}\Rightarrow-\frac{16}{3}<\big(\frac{\tau_1}{d_0}-\frac{2}{3}\big)\big(8k-4-4k\frac{\tau_1}{d_0}\big)<\frac{80}{9}\\
\overbracket{\Rightarrow}^{\text{\eqref{eq:g0d1ratioestimate}}}&-\frac{1020r}{9R}=-\frac{80\times51r}{9\times4R}<\big(\frac{\tau_1}{d_0}-\frac{2}{3}\big)\frac{g_0}{d_1}\big(8k-4-4k\frac{\tau_1}{d_0}\big)<\frac{80\times51r}{9\times4R}=\frac{1020r}{9R}.\\
\end{aligned}
\end{equation}
And
\begin{equation}\label{eq:2ndterminC1C2}
    \begin{aligned}
    &\left.
    \begin{aligned}
    -3k-\frac{2d_2}{3d_0}\overbracket{<}^{\qquad\mathllap{\tau_1>0}}&\big(-3k+\frac{\tau_1}{d_0}-\frac{2d_2}{3d_0}\big)\overbracket{<}^{\qquad\quad\qquad\qquad\mathllap{\text{\cref{remark:deffpclf}: }\tau_1<2d_0}}\big(-3k+2-\frac{2d_2}{3d_0}\big)\\
    \text{\cref{proposition:casecparametersestimate}\ref{itemcasecparametersestimate5}: }&0.9975<\frac{d_2}{r\sin\phistar}<1.0026,\quad\frac{1}{1.0026}<\frac{r\sin\phistar}{d_0}<\frac{1}{0.9975}
    \end{aligned}
    \right\}\\
    \Longrightarrow& -3k-\frac{2\times1.0026}{3\times0.9975}<-3k+\frac{\tau_1}{d_0}-\frac{2d_2}{3d_0}<-3k+2-\frac{2\times0.9975}{3\times1.0026}\overset{k\ge1}<0\\
    \overbracket{\Longrightarrow}^{\mathllap{\text{\cref{lemma:h1g1g0g1g2g3lengthfractionestimate}\eqref{item01h1g1g0g1g2g3lengthfractionestimate}}}}& \big(-3k-\frac{2\times1.0026}{3\times0.9975}\big)\big(16.65r/R\big)<h_1\big(-3k+\frac{\tau_1}{d_0}-\frac{2d_2}{3d_0}\big)<\big(-3k+2-\frac{2\times0.9975}{3\times1.0026}\big)\big(-8.53r/R\big)
\end{aligned}
\end{equation}
\begin{equation}\label{eq:1stterminC1C2}
    \begin{aligned}-8k\frac{17r}{4R}\frac{1.0026}{0.9975}\overbracket{<}^{\mathrlap{\text{\cref{proposition:casecparametersestimate}\ref{itemcasecparametersestimate5}}}\qquad}&\qquad-8k\frac{17r}{4R}\frac{d_1}{0.9975r\sin\phistar}\\
    \overbracket{<}^{\mathllap{\text{\cref{lemma:h1g1g0g1g2g3lengthfractionestimate}\eqref{item03h1g1g0g1g2g3lengthfractionestimate}}}}&8kd_1g_1\overset{\text{\cref{lemma:h1g1g0g1g2g3lengthfractionestimate}\eqref{item03h1g1g0g1g2g3lengthfractionestimate}}}<8k\frac{17r}{4R}\frac{d_1}{0.9975r\sin\phistar}\overbracket{<}^{\qquad\qquad\mathllap{\text{\cref{proposition:casecparametersestimate}\ref{itemcasecparametersestimate5}}}\qquad}8k\frac{17r}{4R}\frac{1.0026}{0.9975}
    \end{aligned}
\end{equation}
In the end, from \eqref{eq:C1C2reformed}, by \eqref{eq:1stterminC1C2}, \eqref{eq:2ndterminC1C2}, \eqref{eq:3rdterminC1C2}, \eqref{eq:4thterminC1C2} we get
\begin{align*}
\mathcal{C}_1(k)+\mathcal{C}_2(k)\overbracket{=}^{\mathllap{\eqref{eq:C1C2reformed}}}&\overbracket{8kd_1g_1}^{\eqref{eq:1stterminC1C2}}+\overbracket{h_1\big(-3k+\frac{\tau_1}{d_0}-\frac{2d_2}{3d_0}\big)}^{\eqref{eq:2ndterminC1C2}}+\overbracket{\big(\frac{\tau_1}{d_0}-\frac{2}{3}\big)\frac{g_0}{d_1}\big(8k-4-4k\frac{\tau_1}{d_0}\big)}^{\eqref{eq:3rdterminC1C2}}\\
&+\underbracket{\frac{2(g_2-g_0)}{3d_0}\Big[11k-\frac{8kr\sin\phistar}{d_1}+\frac{4r\sin\phistar}{d_1}+\frac{4k\tau_1}{d_0}\frac{r\sin\phistar}{d_1}+\frac{4k\tau_1+(4-8k)d_0}{d_0}\frac{g_0}{d_1}\Big]}_{\eqref{eq:4thterminC1C2}}\\
>&-8k\frac{17r}{4R}\frac{1.0026}{0.9975}+\big(-3k-\frac{2\times1.0026}{3\times0.9975}\big)\big(\frac{16.65r}{R}\big)-\frac{1020r}{9R}-\frac{17\times1.0026r}{3R}\big(11k+12+\frac{85r}{R}\big)\\
>&\frac{r}{R}\big(-146.7k-192.7-483\frac{r}{R}\big)\overbracket{>}^{k\le4/3}-389\frac{r}{R}-483\big(\frac{r}{R}\big)^2\qedhere
\end{align*}
\end{proof}

\begin{figure}[h]
\begin{center}    
    \begin{tikzpicture}[xscale=1.0,yscale=1.0]
        \pgfmathsetmacro{\PHISTAR}{6};
        \pgfmathsetmacro{\WIDTH}{8};
        \tkzDefPoint(-\PHISTAR,0){LL};
        \tkzDefPoint(\PHISTAR,0){RL};
        \tkzDefPoint(-\PHISTAR,\PHISTAR){X};
        \tkzDefPoint(\PHISTAR,\PHISTAR){IY};
        \tkzDefPoint(-\PHISTAR,0.5*\PHISTAR){S1};
        \tkzDefPoint(\PHISTAR,0.5*\PHISTAR){RS1};
        \tkzDefPoint(-\PHISTAR,\PHISTAR/3){L4};
        \path[name path=lineLLIY] (LL)--(IY);
        \path[name path=lineS1RL] (S1)--(RL);
        \path[name path=lineRS1LL] (RS1)--(LL);
        \path[name path=lineXRL] (X)--(RL);
        \path[name intersections= 
        {of=lineLLIY and lineS1RL,by={LS}}];
        \path[name intersections=
        {of=lineLLIY and lineXRL,by={HS}}];
        \path[name intersections=
        {of=lineRS1LL and lineXRL,by={RLS}}];
        \draw[name path=IYRS1, dashed, thin] (IY)--(RS1);
        \draw[name path=RS1RLS, dashed, thin] (RS1)--(RLS);
        \draw[name path=RLSHS, dashed, thin] (RLS)--(HS);
        \draw[name path=HSIY, dashed, thin] (HS)--(IY);

        \tkzLabelPoint[left](L4){\small $\theta_1=\frac{1}{3}\Phistar$};
        \draw[name path=lineLSHS, very thick](LS)--(HS);
        \draw[name path=lineHSX, very thick](HS)--(X);
        \draw[name path=lineXS1, very thick](X)--(S1);
        \draw[name path=lineS1LS, very thick](S1)--(LS);
        \draw[ultra thin,dashed](L4)--(LS);
        \fill [orange, opacity=2/30](LS) -- (HS) -- (X) -- (S1) -- cycle;
        \tkzDefPoint(-\PHISTAR,0.47*\PHISTAR){L5};
        \path[name path=HorSeg] (L5)-- +(25,0);
        \path[name intersections={of=lineS1RL and HorSeg,by={Seg1}}];
        \path[name intersections={of=lineLSHS and HorSeg,by={Seg2}}];
        \draw[ultra thin,red,dashed](Seg1)--(Seg2);
        \draw[ultra thin,red,dashed](L5)--(Seg1);
        \tkzDefPoint(-\PHISTAR,0.4*\PHISTAR){L6};
        \path[name path=HorSeg2] (L6)-- +(15,0);
        \path[name intersections={of=lineS1LS and HorSeg2,by={Seg3}}];
        \path[name intersections={of=lineLSHS and HorSeg2,by={Seg4}}];
        \draw[ultra thin,blue](Seg3)--(Seg4);
        \path[name intersections={of=RLSHS and HorSeg2,by={Seg7}}];
        \path[name intersections={of=RS1RLS and HorSeg2,by={Seg8}}];
        \draw[blue, ultra thin] (Seg7)--(Seg8);
        \tkzLabelSegment[above](Seg7,Seg8){\small $\Fcal(\mathfrak{l}_L)$}; 
        \tkzLabelPoint[left, xshift=-0.7](Seg7){\small $\Fcal(\mathcal{W}_L)$};
        \tkzLabelPoint[right, xshift=2.7](Seg8){\small $\Fcal(\mathcal{E}_L)$};
        \tkzLabelPoint[above](L5){\small $\textcolor{red}{d_1=\frac{1}{2}r\sin{\phistar}}$};
        \tkzLabelPoint[above right](X){\small $d_1=r\sin{\phistar}$, $\theta_1=\Phistar$};
        \tkzDefPoint(-\PHISTAR,0.6*\PHISTAR){Seg5};
        \path[name path=HorSeg3] (Seg5)-- +(15,0);
        \path[name intersections={of=lineHSX and HorSeg3,by={Seg6}}];
        \draw[ultra thin,blue](Seg5)--(Seg6);
        \path[name intersections={of=HSIY and HorSeg3,by={Seg9}}];
        \path[name intersections={of=IYRS1 and HorSeg3,by={Seg10}}];
        \draw[blue, ultra thin] (Seg9)--(Seg10);
        \tkzLabelSegment[above](Seg9,Seg10){\small $\Fcal(\mathfrak{l}_U)$}; 
        \tkzLabelPoint[left, xshift=-0.7](Seg9){\small $\Fcal(\mathcal{W}_U)$};
        \tkzLabelPoint[right, xshift=2.7](Seg10){\small $\Fcal(\mathcal{E}_U)$};
        \tkzDrawPoints(X,LS,HS,S1,L4,Seg1,Seg2,L5,Seg3,Seg4,Seg5,Seg6,Seg7,Seg8,Seg9,Seg10,RL,LL,IY);
        \tkzLabelPoint[right, xshift=2.7](Seg4){\small $\mathcal{E}_L$};
        \tkzLabelPoint[right, xshift=2.7](Seg6){\small $\mathcal{E}_U$};
        \tkzLabelPoint[left, xshift=-0.7](Seg5){\small $\mathcal{W}_U$};
        \tkzLabelSegment[above](Seg5,Seg6){\small $\mathfrak{l}_U$}; 
        \tkzLabelPoint[left, xshift=-0.7](Seg3){\small $\mathcal{W}_L$};
        \tkzLabelSegment[above](Seg3,Seg4){\small{$\mathfrak{l}_L$}};
        \node[] at ($(X)+(-58:1)$)  { \textcolor{blue}{$M^{\text{\upshape{in}}}_{R,1}$}};
        \node[] at ($(IY)+(-122:1)$)  { \textcolor{red}{$M^{\text{\upshape{out}}}_{R,1}$}};
        \draw[dashed,ultra thin](LL)--(RL);
        \draw[dashed,ultra thin](RL)--(RS1);
        \draw[dashed,ultra thin](LL)--(S1);
        \draw[dashed,ultra thin](RLS)--(LL);
        \draw[dashed,ultra thin](RLS)--(RL);
        \draw[dashed,ultra thin](LS)--(LL);
        \draw[dashed,ultra thin](LS)--(RL);
        \tkzLabelPoint[below right](LL){\small $\Phi=-\Phistar, SW=(-\Phistar,0)$};
        \tkzLabelPoint[below](RL){\small $\Phi=+\Phistar$};
        \tkzLabelPoint[above](IY){\small $NE=(\Phistar,\Phistar)$};
        
    \end{tikzpicture}
    \caption{The \textcolor{blue}{$M^{\text{\upshape in}}_{R,1}$} component  from the $n_1=1$-region of \cref{fig:x1_n1cases_caseC}.}\label{fig:x1_n1equal1case}
    \end{center}
\end{figure}

\begin{proposition}\label{proposition:QUpperLowerRegionsLowerBoundsChatCtilde}
    Let $x_1=(\Phi_1,\theta_1)\in \textcolor{blue}{M^{\text{\upshape{in}}}_{R,1}}$  (\cref{fig:x1_n1equal1case}), i.e., $n_1=1$ (see \cref{fig:x1_n1cases_caseC}). Then this indeed falls into case (c) of \eqref{eqMainCases}, which is the hypothesis of \cref{TEcase_c}. Suppose $x_1$ is on the line $\mathfrak{l}=\Big\{(\Phi,\theta_1)\bigm|d_1=R\sin\theta_1=\const\Big\}\subset\textcolor{blue}{M^{\text{\upshape{in}}}_{R,1}}$. With $\mathtt{Q}(t)$ from \eqref{eq:II1II2Quadratic},\COMMENT{Remember that we are considering a change in notation} and fixed $k\in[1,4/3]$, $\nu=\frac{d_1}{r\sin\phistar}$ we have the following.
    \begin{enumerate}
        \item\label{item01Qlowerboundn1eq1} $\nu\in(1/3,1)$.
        \item\label{item02Qlowerboundn1eq1} The length functions $\tau_0$, $\tau_1$, $d_0$, $d_1$ $d_2$ defined in \cref{def:fpclf} are functions of $x\in\mathfrak{l}$. $d_0$ as a function of $x=(\Phi,\theta_1)\in\mathfrak{l}$ with constant $\theta_1$, is thus a monotone increasing function of $\Phi$. And $\tau_1$ as a function of $x=(\Phi,\theta_1)\in\mathfrak{l}$ with constant $\theta_1$, is thus a monotonically decreasing function of $\Phi$.\COMMENT{Why ``thus''?}
        \item\label{item03Qlowerboundn1eq1} If  $\nu>0.5$ that is $x_1\in\mathfrak{l}=\mathfrak{l}_U=
        \Big\{x=(\Phi,\theta_1)\bigm|d_1=\nu r\sin\phistar>0.5r\sin\phistar\Big\}$ with endpoints $\mathcal{W}_U=(\Phi_{\mathcal{W}_U},\theta_1)$ (think "west"), $\mathcal{E}_U=(\Phi_{\mathcal{E}_U},\theta_1)$ (think "east") in \cref{fig:x1_n1equal1case} and $\Phi_{\mathcal{W}_U}<\Phi_{\mathcal{E}_U}$, then \[\mathtt{Q}\big(\frac{\tau_1}{d_0}\bigm|_{x=x_1}\big)>\min{\Big\{\frac{11}{9},\frac{1}{3}+\mathcal{C}^{\mathcal{W}}_{U}(\nu,k)\Big\}},\] where, with $h_1,g_1$ from \cref{def:lengthfractionfunctionsn1eq1}\eqref{item01lengthfractionfunctionsn1eq1}, \eqref{item03lengthfractionfunctionsn1eq1}, $\hat{h}_1(U)\dfn\lim_{\mathfrak{l}_U\ni x\to \mathcal{W}_U}h_1(x)$, $\hat{g}_1(U)\dfn\lim_{\mathfrak{l}_U\ni x\to \mathcal{W}_U}g_1(x)$,
        \begin{equation}\label{eq:CUWnuk}
        \begin{aligned}
        \mathcal{C}^{\mathcal{W}}_U(\nu,k)=&\Big[\frac{-16k-12}{3\nu}+8k\Big]\big(\hat{h}_1(U)-2d_1\hat{g}_1(U)\big)-\frac{4k}{\nu}\big(\hat{h}_1(U)-2d_1\hat{g}_1(U)\big)^2.
        \end{aligned}    
        \end{equation}
        
        \item\label{item04Qlowerboundn1eq1} If  $\nu\le 0.5$ that is $x_1\in\mathfrak{l}=\mathfrak{l}_L=
        \Big\{x=(\Phi,\theta_1)\bigm|d_1=\nu r\sin\phistar\le0.5r\sin\phistar\Big\}$ with endpoints $\mathcal{W}_L=(\Phi_{\mathcal{W}_L},\theta_1)$, $\mathcal{E}_L=(\Phi_{\mathcal{E}_L},\theta_1)$ in \cref{fig:x1_n1equal1case} and $\Phi_{\mathcal{W}_L}<\Phi_{\mathcal{E}_L}$, then \[\mathtt{Q}\big(\frac{\tau_1}{d_0}\Bigm|_{x=x_1}\big)>\min{\Big\{\frac{11}{9}+\mathcal{C}^{\mathcal{E}}_{L}(\nu,k),\frac{1}{3}+\mathcal{C}^{\mathcal{W}}_{L}(\nu,k)\Big\}},
        \] 
        where
        \begin{equation}\label{eqL:CLWnukandCLEnukand}
        \begin{aligned}
            \mathcal{C}^{\mathcal{W}}_{L}(\nu,k)=&2d_1\hat{g}_1(L)\big(-24k+\frac{32k-12}{3\nu}\big)-\frac{16kd_1^2}{\nu}\hat{g}^2_1(L),
            \\
            \mathcal{C}^{\mathcal{E}}_{L}(\nu,k)=&\big(-8k+\frac{32k-12}{3\nu}\big)\big(\tilde{h}_1(L)-4d_1\tilde{g}_1(L)\big)-\frac{8k}{\nu}\big(\tilde{h}_1(L)-4d_1\tilde{g}_1(L)\big)\big(2-4\nu\big)-\frac{4k}{\nu}\big(\tilde{h}_1(L)-4d_1\tilde{g}_1(L)\big)^2,
        \end{aligned}   
        \end{equation}
        with (see \cref{def:lengthfractionfunctionsn1eq1}\eqref{item01lengthfractionfunctionsn1eq1},\eqref{item03lengthfractionfunctionsn1eq1})
        \begin{align*}
        \hat{h}_1(L)\dfn&\lim_{\mathfrak{l}_L\ni x\to \mathcal{W}_L}h_1(x),&
            \hat{g}_1(L)\dfn&\lim_{\mathfrak{l}_L\ni x\to \mathcal{W}_L}g_1(x),\\
            \tilde{h}_1(L)\dfn&\lim_{\mathfrak{l}_L\ni x\to \mathcal{E}_L}h_1(x),&
            \tilde{g}_1(L)\dfn&\lim_{\mathfrak{l}_L\ni x\to \mathcal{E}_L}g_1(x).
        \end{align*}
    \item\label{item05Qlowerboundn1eq1}  $\mathcal{C}_U^\mathcal{W}(\nu,k)$, $\mathcal{C}_L^\mathcal{E}(\nu,k)$, and $\mathcal{C}_L^\mathcal{W}(\nu,k)$  from \eqref{eq:CUWnuk},  \eqref{eqL:CLWnukandCLEnukand} satisfy \begin{align*}
        \mathcal{C}_U^\mathcal{W}(\nu,k)>&-192\big(\frac{r}{R}\big)-4667\big(\frac{r}{R}\big)^2,\\
        \mathcal{C}_L^\mathcal{W}(\nu,k)>&-50\big(\frac{r}{R}\big)-194\big(\frac{r}{R}\big)^2,\\
        \mathcal{C}_L^\mathcal{E}(\nu,k)>&-92\big(\frac{r}{R}\big)-7984\big(\frac{r}{R}\big)^2.
    \end{align*}
    \end{enumerate} 
\end{proposition}
\begin{proof}
These computations will run to page \pageref{proposition:casecexpansionn1eq1}.
        \hidesubsection{Proof of (1)}\cref{proposition:casecparametersestimate}\eqref{item03lengthfractionfunctionsn1eq1} with $n_1=1$ implies $d_1>\frac{1}{3}r\sin\phistar$. Also, since $p(x_1)=P$, $p(\mathcal{F}(x_1))=P_1$ are two consecutive collisions on $\Gamma_R$ with distance $2d_1$, $2d_1=|p(x_1)p(\mathcal{F}(x_1))|<|AB|=2r\sin\phistar$ (as shown in \cref{fig:n1_equal1_upper_leftendpoint_case,fig:n1_equal1_lower_leftendpoint_case,fig:n1_equal1_upper_leftendpoint_case2,fig:n1_equal1_lower_rightendpoint_case,fig:n1_equal1_upper_rightendpoint_case,fig:n1_equal1_upper_rightendpoint_case2}). Hence $\frac{1}{3}r\sin\phistar<d_1<r\sin\phistar$, $\nu=\frac{d_1}{r\sin\phistar}\in\big(\frac{1}{3},1\big)$.
        \hidesubsection{Proof of (2)}By \cref{def:fpclf}, $\tau_0=|p(x)p(\Fcal^{-1}(x))|$, $d_0$ are determined by the collision angle $\theta$ coordinate of $\Fcal^{-1}(x)\in\Mrout$, $\tau_1=|p(\Fcal(x))p(\Fcal^{2}(x))|$, $d_2$ are determined by the collision angle $\theta$ coordinate of $\Fcal^{2}(x)\in \Mrin$. Hence length functions $\tau_0$, $\tau_1$, $d_0$, $d_2$ defined in \cref{def:fpclf} are functions of $x\in\mathfrak{l}$.
        
        Since in \cref{fig:x1_n1equal1case} either $\mathfrak{l}=\mathfrak{l}_U\subset\MRin$ or $\mathfrak{l}=\mathfrak{l}_L\subset\MRin$ with some constant $\theta_1<\Phistar=\sin^{-1}(r\sin\phistar/R)\overset{R>r}<\sin^{-1}{(\sqrt{4r/R})}$, \cref{monotone_p4} implies that $d_0$ is a monotone increasing function of $\Phi$.

        Let $x\in\mathfrak{l}$, $\Fcal(x)\in\Fcal(\mathfrak{l})\subset\MRout$ hence $\Fcal^2(x)\in\Mrin$ then we can suppose $\Gamma_r\ni p(\Fcal^2(x))=(x_Q,y_Q)$ in the standard coordinate system defined by \cref{def:StandardCoordinateTable} for the billiard table.

        We note that $\Fcal\bigm|_\mathfrak{l}(x)=\Fcal(\Phi,\theta_1=\const)=(\Phi+2\theta_1,\theta_1=\const)$ is the bijection/diffeomorphism\COMMENT{Why ``bijection/diffeomorphism'' instead of one of the two? This kind of a/b should be avoided whenever possible.} between $\mathfrak{l}_U$ and $\Fcal(\mathfrak{l}_U)\subset\MRout$, between $\mathfrak{l}_L$ and $\Fcal(\mathfrak{l}_L)\subset\MRout$ in \cref{fig:x1_n1equal1case}. By \cref{monotone_p3}\eqref{itemmonotone_p3-1}, \eqref{eq1monotone_p3} and chain rule we get\COMMENT{The annotations in the next equation are unclear---and have too large vertical bars.} 
        \[
        \frac{d\tau_1}{d\Phi}\Bigm|_{x=(\Phi,\theta_1=\const)}\overbracket{=}^{\strut\mathclap{\text{chain rule}}}\frac{d\tau_1}{d\Phi\bigm|_{\Fcal(\mathfrak{l})}}\cdot\overbracket{\frac{d\Phi\bigm|_{\Fcal(\mathfrak{l})}}{d\Phi\bigm|_{\mathfrak{l}}}\Bigm|_{x}}^{\mathclap{\text{since }\Phi|_{\Fcal(\mathfrak{l})}=(\Phi+2\theta_1)|_{\mathfrak{l}},\ \frac{d(\Phi+2\const)}{d\Phi}=1}}=\frac{d\tau_1}{d\Phi\bigm|_{\Fcal(\mathfrak{l})}}\overbracket{=}^{\bigstrut\mathclap{\text{\cref{monotone_p3}\eqref{itemmonotone_p3-1}, \eqref{eq1monotone_p3}}}}\frac{-bx_Q}{d_2}
        \]
        We also recall that $I y_*=(\phistar,\phistar)$ in \cref{fig:MrN}, $\Fcal^2(x)\nfd(\phi_Q,\theta_Q)\in\Nin$ by \cref{contraction_region} and more precisely in \eqref{eq:caseCNinNout1}\COMMENT{?} we have 
        \[
        0<\phi_Q-\phistar<\sqrt{|\theta_Q-\phistar|^2+|\phi_Q-\phistar|^2}=\|x_Q-Iy_*\|\overbracket{<}^{\strut\mathclap{\text{\cite[equation (3.21)]{hoalb}}}}\frac{5.84r\sin\phistar}{R}.
        \] 
        Therefore, with conditions $\eqref{eqJZHypCond}$ and $\phistar<\tan^{-1}(1/3)$, we get\COMMENT{overbracket and mathclap in next line} 
        \[
        \phistar<\phi_Q<\phistar+\frac{5.84r\sin\phistar}{R}\overset{\text{\eqref{eqJZHypCond}:} R>1700r\text{ and } \phistar<\tan^{-1}\big(1/3\big)}<\tan^{-1}\big(1/3\big)+\frac{5.84}{1700}<\frac{\pi}{2}.
        \] 
        Hence in coordinate system \cref{def:StandardCoordinateTable}, $x_Q=r\sin{\phi_Q}>0$, $\frac{d\tau_1}{d\Phi}\Bigm|_{x=(\Phi,\theta_1=\const)}=-\frac{bx_Q}{d_2}<0$, that is, $\tau_1$ as a function of $x=(\Phi,\theta_1)\in\mathfrak{l}$ with constant $\theta_1$, is thus a monotonically decreasing function of $\Phi$.
\begin{figure}[h]
\begin{center}
\begin{tikzpicture}[xscale=1.4,yscale=1.4]
    \tkzDefPoint(0,0){Or};
    \tkzDefPoint(0,-4.8){Y};
    \clip(-5.7,-4.2) rectangle (5.7,-2.8);
    \pgfmathsetmacro{\rradius}{4.5};
    \pgfmathsetmacro{\phistardeg}{40};
    \pgfmathsetmacro{\XValueArc}{\rradius*sin(\phistardeg)};
    \pgfmathsetmacro{\YValueArc}{\rradius*cos(\phistardeg)};
    \pgfmathsetmacro{\Rradius}{
        \rradius*sin(\phistardeg)/sin(18)
    };
    \pgfmathsetmacro{\bdist}{(-\rradius*cos(\phistardeg))+cos(18)*\Rradius}
    \tkzDefPoint(0,\bdist){OR};
    \tkzDefPoint(\XValueArc,-1.0*\YValueArc){B};    \tkzDefPoint(-1.0*\XValueArc,-1.0*\YValueArc){A};
    \tkzDrawArc[name path=Cr,very thick](Or,B)(A);
    \tkzDrawArc[name path=CR,very thick](OR,A)(B);
    \pgfmathsetmacro{\PXValue}{\Rradius*sin(14)};
    \pgfmathsetmacro{\PYValue}{\bdist-\Rradius*cos(14)};
    \tkzDefPoint(\PXValue,\PYValue){P2};
    \tkzDefPointBy[projection=onto OR--P2](A)\tkzGetPoint{M};
    \coordinate (Q_far) at ($(A)!2!(M)$);
    \path[name path=Line_P2_Q_far](P2)-- (Q_far);
    \path[name intersections={of=Cr and Line_P2_Q_far,by={Q}}];
    \tkzDefPoint(0,\bdist-\Rradius){C};
    \begin{scope}[very thin,decoration={markings,mark=at position 0.5 with {\arrow{>}}}] 
    \draw[blue,postaction={decorate}] (P2)--(Q);
    \draw[blue,postaction={decorate}] (A)--(P2);
    \end{scope}
    \tkzLabelPoint[above right](A){\small $A=P=p(x=\mathcal{W}_{U})$};
    \tkzLabelPoint[below left](A){\small $\hat{\tau}_0(U)=|AP|=0$};
    \tkzLabelSegment[above right=0.2pt and 0.6pt](A,P2){$2\hat{d}_1(U)=2d_1$};
    \tkzLabelPoint[below, yshift=-3.5](P2){\small $P_1=p(\mathcal{F}(x=\mathcal{W}_{U}))$};
    \tkzLabelPoint[right](Q){\small $Q=p(\mathcal{F}^2(x=\mathcal{W}_{U}))$};
    \tkzLabelPoint[right](B){$B$};
    \tkzLabelPoint[below](Or){$O_r$};
    \tkzDrawPoints(A,B,Or,Q,P2);
    \tkzLabelSegment[above left=0.5pt and 3pt](P2,Q){$\hat{\tau}_1(U)=|P_1Q|$};
    
\end{tikzpicture}
\caption{With constant $\theta_1\ge0.5\Phistar$, $(\Phi_{\mathcal{W}_U},\theta_1)=x=\mathcal{W}_U$ is on the trajectory from corner $A$  }\label{fig:n1_equal1_upper_leftendpoint_case}

\end{center}
\end{figure}
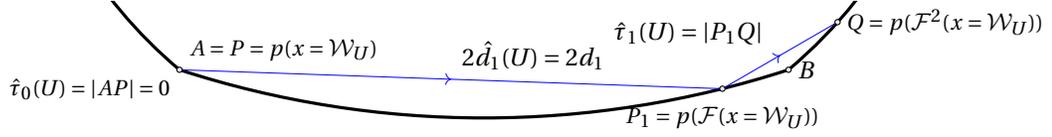

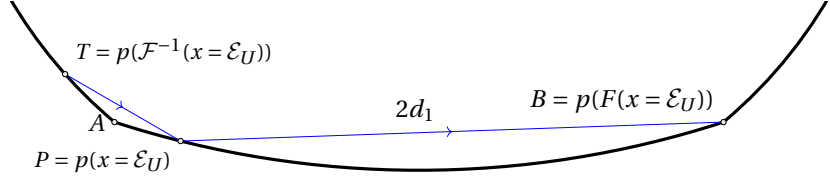
\begin{figure}[h] 
\begin{center}
\begin{tikzpicture}[xscale=1.4,yscale=1.4]
    \tkzDefPoint(0,0){Or};
    \tkzDefPoint(0,-4.8){Y};
    \clip(-5.5,-4.3) rectangle (5.5,-2.3);
    \pgfmathsetmacro{\rradius}{4.5};
    \pgfmathsetmacro{\phistardeg}{40};
    \pgfmathsetmacro{\XValueArc}{\rradius*sin(\phistardeg)};
    \pgfmathsetmacro{\YValueArc}{\rradius*cos(\phistardeg)};
    \pgfmathsetmacro{\Rradius}{
        \rradius*sin(\phistardeg)/sin(18)
    };
    \pgfmathsetmacro{\bdist}{(-\rradius*cos(\phistardeg))+cos(18)*\Rradius}
    \tkzDefPoint(0,\bdist){OR};
    \tkzDefPoint(\XValueArc,-1.0*\YValueArc){B};    \tkzDefPoint(-1.0*\XValueArc,-1.0*\YValueArc){A};
    \tkzDrawArc[name path=Cr,very thick](Or,B)(A);
    \tkzDrawArc[name path=CR,very thick](OR,A)(B);
    \pgfmathsetmacro{\PXValue}{-\Rradius*sin(14)};
    \pgfmathsetmacro{\PYValue}{\bdist-\Rradius*cos(14)};
    \tkzDefPoint(\PXValue,\PYValue){P1};
    
    \tkzDefPointBy[projection=onto OR--P1](B)\tkzGetPoint{M};
    \coordinate (T_far) at ($(B)!2!(M)$);
    \path[name path=Line_P1_T_far](P1)-- (T_far);
    \path[name intersections={of=Cr and Line_P1_T_far,by={T}}];
    \tkzDefPoint(0,\bdist-\Rradius){C};
    \begin{scope}[ultra thin,decoration={markings,mark=at position 0.5 with {\arrow{>}}}] 
    \draw[blue,postaction={decorate}] (T)--(P1);
    \draw[blue,postaction={decorate}] (P1)--(B);
    \end{scope}
    \tkzLabelPoint[left](A){$A$};
    \tkzLabelPoint[below left](P1){\small $P=p(x=\mathcal{E}_U)$};
    \tkzLabelPoint[above right](T){\small $T=p(\mathcal{F}^{-1}(x=\mathcal{E}_U))$};
    \tkzLabelPoint[above left](B){$B=p(F(x=\mathcal{E}_U))$};
    \tkzLabelPoint[below](Or){$O_r$};
    \tkzLabelSegment[above left=0.5pt and 3pt](P1,B){$2d_1$};
    \tkzDrawPoints(A,B,Or,P1,T);
\end{tikzpicture}
\caption{With constant $\theta_1\ge0.5\Phistar$, $(\Phi_{\mathcal{E}_U},\theta_1)=x=\mathcal{E}_U$ is on a trajectory to corner $B$}\label{fig:n1_equal1_upper_rightendpoint_case}
\end{center}
\end{figure}

\begin{figure}[h]
\begin{center}
\begin{tikzpicture}[xscale=1.8,yscale=1.8]
    \tkzDefPoint(0,0){Or};
    \tkzDefPoint(0,-4.8){Y};
    \clip(-4.3,-4.2) rectangle (4.3,-2.9);
    \pgfmathsetmacro{\rradius}{4.5};
    \pgfmathsetmacro{\phistardeg}{40};
    \pgfmathsetmacro{\XValueArc}{\rradius*sin(\phistardeg)};
    \pgfmathsetmacro{\YValueArc}{\rradius*cos(\phistardeg)};
    \pgfmathsetmacro{\Rradius}{
        \rradius*sin(\phistardeg)/sin(18)
    };
    \pgfmathsetmacro{\bdist}{(-\rradius*cos(\phistardeg))+cos(18)*\Rradius}
    \tkzDefPoint(0,\bdist){OR};
    \tkzDefPoint(\XValueArc,-1.0*\YValueArc){B};    \tkzDefPoint(-1.0*\XValueArc,-1.0*\YValueArc){A};
    \tkzDrawArc[name path=Cr,very thick](Or,B)(A);
    \tkzDrawArc[name path=CR,very thick](OR,A)(B);
    \pgfmathsetmacro{\PXValue}{-\Rradius*sin(1)};
    \pgfmathsetmacro{\PYValue}{\bdist-\Rradius*cos(1)};
    \tkzDefPoint(\PXValue,\PYValue){P1};
    
    \tkzDefPointBy[projection=onto OR--P1](A)\tkzGetPoint{M};
    \coordinate (P2) at ($(A)!2!(M)$);
    
    \tkzDefPointBy[projection=onto OR--P2](P1)\tkzGetPoint{M2};
    \coordinate (Q_far) at ($(P1)!2!(M2)$);
    \path[name path=Line_P2_Q_far](P2)-- (Q_far);
    \path[name intersections={of=Cr and Line_P2_Q_far,by={Q}}];
    \tkzDefPoint(0,\bdist-\Rradius){C};
    \begin{scope}[ultra thin,decoration={markings,mark=at position 0.5 with {\arrow{>}}}] 
    \draw[blue,postaction={decorate}] (A)--(P1);
    \draw[blue,postaction={decorate}] (P1)--(P2);
    \draw[blue,postaction={decorate}] (P2)--(Q);
    \end{scope}
    \tkzLabelPoint[above right](A){\small $A=p\big(\mathcal{F}^{-1}(x=\mathcal{E}_U)\big)$};
    \tkzLabelPoint[below](P1){\small $P=p\big(x=\mathcal{E}_U\big)$};
    \tkzLabelPoint[below](P2){\small $p\big(\mathcal{F}(x=\mathcal{E}_U)\big)$};
    \tkzLabelPoint[above left](Q){\small $Q=p\big(\mathcal{F}^{2}(x=\mathcal{E}_U)\big)$};
    \tkzLabelPoint[right](B){$B$};
    \tkzLabelPoint[below](Or){$O_r$};
    \tkzLabelSegment[above left=0.5pt and 0.1pt](P1,P2){$2d_1$};
    \tkzDrawPoints(A,B,Or,P1,P2,Q);
    
\end{tikzpicture}
\caption{With constant $\theta_1<0.5\Phistar$ and $d_1>0.5r\sin\phistar$, $(\Phi_{\mathcal{E}_U},\theta_1)=x=\mathcal{E}_U$ is on a trajectory from corner $A$
} \label{fig:n1_equal1_upper_rightendpoint_case2}
\end{center}
\end{figure}

\begin{figure}[h]
\begin{center}
\begin{tikzpicture}[xscale=1.5,yscale=1.5]
    \tkzDefPoint(0,0){Or};
    \tkzDefPoint(0,-4.8){Y};
    \clip(-4.7,-4.2) rectangle (4.7,-2.7);
    \pgfmathsetmacro{\rradius}{4.5};
    \pgfmathsetmacro{\phistardeg}{40};
    \pgfmathsetmacro{\XValueArc}{\rradius*sin(\phistardeg)};
    \pgfmathsetmacro{\YValueArc}{\rradius*cos(\phistardeg)};
    \pgfmathsetmacro{\Rradius}{
        \rradius*sin(\phistardeg)/sin(18)
    };
    \pgfmathsetmacro{\bdist}{(-\rradius*cos(\phistardeg))+cos(18)*\Rradius}
    \tkzDefPoint(0,\bdist){OR};
    \tkzDefPoint(\XValueArc,-1.0*\YValueArc){B};    \tkzDefPoint(-1.0*\XValueArc,-1.0*\YValueArc){A};
    \tkzDrawArc[name path=Cr, very thick](Or,B)(A);
    \tkzDrawArc[name path=CR, very thick](OR,A)(B);
    \pgfmathsetmacro{\PXValue}{\Rradius*sin(1)};
    \pgfmathsetmacro{\PYValue}{\bdist-\Rradius*cos(1)};
    \tkzDefPoint(\PXValue,\PYValue){P2};
    
    \tkzDefPointBy[projection=onto OR--P2](B)\tkzGetPoint{M};
    \coordinate (P1) at ($(B)!2!(M)$);
    
    \tkzDefPointBy[projection=onto OR--P1](P2)\tkzGetPoint{M2};
    \coordinate (T_far) at ($(P2)!2!(M2)$);
    \path[name path=Line_P1_T_far](P1)-- (T_far);
    \path[name intersections={of=Cr and Line_P1_T_far,by={T}}];
    \tkzDefPoint(0,\bdist-\Rradius){C};
    \begin{scope}[ultra thin,decoration={markings,mark=at position 0.5 with {\arrow{>}}}] 
    \draw[blue,postaction={decorate}] (T)--(P1);
    \draw[blue,postaction={decorate}] (P1)--(P2);
    \draw[blue,postaction={decorate}] (P2)--(B);
    \end{scope}
    \tkzLabelPoint[left](A){$A$};
    \tkzLabelPoint[below =2.5pt](P1){\small $P=p\big(x_1=\mathcal{W}_U\big)$};
    \tkzLabelPoint[below](P2){\small $p\big(\mathcal{F}(x_1=\mathcal{W}_U)\big)$};
    \tkzLabelPoint[ left =0.8pt](T){\small $T=p(\mathcal{F}^{-1}(x_1=\mathcal{W}_U))$};
    \tkzLabelPoint[above right](B){\small $B=p\big(\mathcal{F}^2(x_1=\mathcal{W}_U)\big)$};
    \tkzLabelPoint[below](Or){$O_r$};
    \tkzLabelSegment[above =0.2pt](P1,P2){$2\hat{d}_1(U)=2d_1$};
    \tkzLabelSegment[above =0.5pt](P2,B){$2d_1=\hat{\tau}_1(U)$};
    \tkzLabelSegment[above right =0pt and 0.3pt](P1,T){\small $\hat{\tau}_0(U)=|PT|>0$};
    \tkzDrawPoints(A,B,P1,P2,T);
\end{tikzpicture}
\caption{
With constant $\theta_1<0.5\Phistar$ and $d_1>0.5r\sin\phistar$, $(\Phi_{\mathcal{W}_U},\theta_1)=x=\mathcal{W}_U$ is on a trajectory to corner $B$
}\label{fig:n1_equal1_upper_leftendpoint_case2}

\end{center}
\end{figure}

\hidesubsection{Proof of (3)} For $x_1\in \mathfrak{l}_U$ in \cref{fig:x1_n1equal1case} with $d_1>0.5r\sin\phistar$, we define\COMMENT{Why not ``$\hat{\tau}_i(U)\dfn\lim_{\mathfrak{l}_U\ni x\to \mathcal{W}_U}\tau_i\ge0$, $\hat{d}_1(U)\dfn\lim_{\mathfrak{l}_U\ni x\to \mathcal{W}_U}d_1$''? Alternatively, fit the pair of $\tau$ and the triple of $d$ side by side.}
\begin{equation}\label{eq:lengthlimitWU}
    \begin{aligned}
        \hat{\tau}_0(U)\dfn&\lim_{\mathfrak{l}_U\ni x\to \mathcal{W}_U}\tau_0\ge0
        \\
        \hat{\tau}_1(U)\dfn&\lim_{\mathfrak{l}_U\ni x\to \mathcal{W}_U}\tau_1\ge0
        \\
        \hat{d}_0(U)\dfn&\lim_{\mathfrak{l}_U\ni x\to \mathcal{W}_U}d_0>0
        \\
        \hat{d}_1(U)\dfn&\lim_{\mathfrak{l}_U\ni x\to \mathcal{W}_U}d_1=\const=\nu\cdot r\sin{\phistar}
        \\
        \hat{d}_2(U)\dfn&\lim_{\mathfrak{l}_U\ni x\to \mathcal{W}_U}d_2>0
    \end{aligned}
    \end{equation}
    
Note that these limits exist.\COMMENT{Because...} In \cref{fig:x1_n1equal1case}, $\mathfrak{l}_U$ and its image $\Fcal(\mathfrak{l}_U)$ are line segments with interiors containing no 1st order singularity of $\Fcal$ (See \cref{def:BasicNotations0}). The line segment $\mathfrak{l}_U$ has endpoints $\mathcal{W}_U$, $\mathcal{E}_U$ possibly being on the boundary of $M_R$.  The line segment $\Fcal(\mathfrak{l}_U)$ has endpoints that are possible on the boundary $M_R$.

By \cref{contraction_region}$, \Fcal^{-1}(\mathfrak{l}_U)$ is a curve in $\Nout$ and $\Fcal^{2}(\mathfrak{l}_U)$ is a curve in $\Nin$ (see \cref{fig:MrN}) with endpoints possibly being on the boundary of $M_r$.
Then these collision angle and position limits for the endpoints can be defined by continuous extension from $\Fcal^{-1}(\mathfrak{l}_U)$, $\Fcal^{2}(\mathfrak{l}_U)$ the interior of $M_r$ to endpoints possibly on the boundary of $M_r$. These defined lengths continuously extended for $x$ at the end points $\mathcal{W}_U$, $\mathcal{E}_U$ can be visualized in \cref{fig:n1_equal1_upper_leftendpoint_case,fig:n1_equal1_upper_rightendpoint_case,fig:n1_equal1_upper_leftendpoint_case2,fig:n1_equal1_upper_rightendpoint_case2}.

Also, since by \eqref{eqJZHypCond} and the \Nin,\Nout defined in \cref{def:SectionSetsDefs} and the $\delta_*$ in \cref{def:SectionSetsDefs} satisfies\COMMENT{This sentence is broken, and so is the next line.}
\[\delta_*=\frac{17}{\sin{(\phistar/2)}}\sqrt{r/R}\overbracket{<}^{\eqref{eqJZHypCond}}\frac{17}{\sin{(\phistar/2)}}\sqrt{\frac{\sin^2{(\phistar/2)}\sin^2{\phistar}}{324}}=\frac{17\sin\phistar}{18}<\frac{17}{18}\phistar\text{ bounded away from }\phistar,\] the $\Nin$, $\Nout$ region $\theta$ coordinates are bounded away from $0$ and $\pi$ (see \cref{fig:MrN}). Therefore, the $\theta$ coordinates of $\Fcal^{-1}(x)\bigm|_{x\in\mathfrak{l}_U}$, $\Fcal^{2}(x)\bigm|_{x\in\mathfrak{l}_U}$ are bounded away from the boundaries $\theta=0$ and $\theta=\pi$ of $M_r$. Thus, $d_0(x)\overset{\text{\cref{def:fpclf}}}=r\sin{\big(\theta \text{ coordinate of }\Fcal^{-1}(x)\big)}$ and $d_2(x)\overset{\text{\cref{def:fpclf}}}=r\sin{\big(\theta \text{ coordinate of }\Fcal^{2}(x)\big)}$ are bounded away from $0$ for all $x\in\mathfrak{l}_U$. These limits $\hat{d}_0>0$, $\hat{d}_2>0$.\COMMENT{How does this relate to the ``$>0$'' in the previous equation?}
Then 
define
\begin{equation}\label{eq:hhatUg1hatU}
    \begin{aligned}
    \hat{h}_1(U)\dfn&\lim_{\mathfrak{l}_U\ni x\to \mathcal{W}_U}h_1(x)\overbracket{=}^{\mathclap{\text{\cref{def:lengthfractionfunctionsn1eq1}\eqref{item01lengthfractionfunctionsn1eq1}}}}\lim_{\mathfrak{l}_U\ni x\to \mathcal{W}_U}\big(\frac{\tau_0}{d_0}+\frac{\tau_1}{d_0}+\frac{2d_1}{d_0}-2\big)\overbracket{=}^{\mathclap{\strut\eqref{eq:lengthlimitWU}\ }}\frac{\hat{\tau}_0(U)}{\hat{d}_0(U)}+\frac{\hat{\tau}_1(U)}{\hat{d}_0(U)}+\frac{2d_1}{\hat{d}_0(U)}-2
    \\
    \hat{g}_1(U)\dfn&\lim_{\mathfrak{l}_U\ni x\to \mathcal{W}_U}g_1(x)\overbracket{=}^{\mathclap{\text{\cref{def:lengthfractionfunctionsn1eq1}\eqref{item03lengthfractionfunctionsn1eq1}}}}\lim_{\mathfrak{l}_U\ni x\to \mathcal{W}_U}\big(\frac{1}{d_0}-\frac{1}{r\sin\phistar}\big)\overbracket=^{\mathclap{\eqref{eq:lengthlimitWU}}}\frac{1}{\hat{d}_0(U)}-\frac{1}{r\sin\phistar}
\end{aligned}
\end{equation}
On the other hand, by the monotone property \cref{proposition:QUpperLowerRegionsLowerBoundsChatCtilde}\eqref{item02Qlowerboundn1eq1}, we have 
\begin{align*}
    0<\frac{\tau_1}{d_0}\Bigm|_{x=x_1}\overbracket{<}^{\mathclap{\text{\cref{proposition:QUpperLowerRegionsLowerBoundsChatCtilde}\eqref{item02Qlowerboundn1eq1}}\strut\qquad\quad}}&\frac{\hat{\tau}_1(U)}{d_0\bigm|_{x=x_1}}\overbracket<^{\mathclap{\qquad\quad\strut\text{\cref{proposition:QUpperLowerRegionsLowerBoundsChatCtilde}\eqref{item02Qlowerboundn1eq1}}}}\frac{\hat{\tau}_1(U)}{\hat{d}_0(U)}\overset{\text{\eqref{eq:hhatUg1hatU}}}=2+\hat{h}_1(U)-\frac{\hat{\tau}_0(U)}{\hat{d}_0(U)}-\frac{2\hat{d}_1(U)}{\hat{d}_0(U)}\\
    \overbracket{\le}^{\mathclap{\text{\eqref{eq:lengthlimitWU}: }\hat{\tau}_0(U)\ge0, d_1=\nu r\sin\phistar}}&2+\hat{h}_1(U)-\frac{2\nu r\sin\phistar}{\hat{d}_0}\overset{\text{\eqref{eq:hhatUg1hatU}}}=2+\hat{h}_1(U)-2\nu r\sin\phistar\big(\hat{g}_1(U)+\frac{1}{r\sin\phistar}\big)\\
    =&2-2\nu+\hat{h}_1(U)-2\nu r\sin\phistar\hat{g}_1(U)\overbracket=^{\mathclap{d_1=\nu r\sin\phistar}}2-2\nu+\hat{h}_1(U)-2d_1\hat{g}_1(U)
\end{align*}

With fixed $k\in[1,4/3]$ and constant $\nu=\frac{d_1}{r\sin\phistar}$, the function $\mathtt{Q}(t)$ from \eqref{eq:II1II2Quadratic} for $t\in (-\infty,+\infty)$ has only one extreme point, which is the global maximum point. Therefore, for $t\in[0,2-2\nu+\hat{h}_1(U)-2d_1\hat{g}_1(U)]$, $\mathtt{Q}\big(t\big)$ cannot reach the minimum in the interior of $[0,2-2\nu+\hat{h}_1(U)-2d_1\hat{g}_1(U)]$. Hence, \begin{equation}\label{eq:UpperQlowerbound}
    \mathtt{Q}\Big(\frac{\tau_1}{d_0}\bigm|_{x=x_1}\Big)>\min_{t\in\big[0,\,\,2-2\nu+\hat{h}_1(U)-2d_1\hat{g}_1(U)\big]}\mathtt{Q}\big(t\big)=\min\Big\{\mathtt{Q}\big(0\big),\,\,\mathtt{Q}\big(2-2\nu+\hat{h}_1(U)-2d_1\hat{g}_1(U)\big)\Big\}.
\end{equation}
For $\frac{d_1}{r\sin\phistar}=\nu\in(0.5,1)$ we estimate \begin{equation}\label{eq:UpperEQlowerbound}
\begin{aligned}
\mathtt{Q}\big(0\big)\overbracket{=}^{\mathclap{\strut\text{\eqref{eq:II1II2Quadratic}}}}&\frac{4k+9}{3}+\big(\frac{8-16k}{3}\big)\frac{r\sin{\phistar}}{d_1}+\frac{8kd_1}{r\sin{\phistar}}\\
    \overbracket{=}^{\mathclap{\frac{d_1}{r\sin\phistar}=\nu}}&3+\frac{4k}{3}-\big(\frac{16k}{3}-\frac{8}{3}\big)\frac{1}{\nu}+8k\nu\\
    \overbracket{>}^{\mathclap{\substack{ \text{For fixed }k\in[1,4/3],\ \frac{16k}{3}-\frac{8}{3}\ge0\text{, so }\\\nu\mapsto-\big(\frac{16k}{3}-\frac{8}{3}\big)\frac{1}{\nu}+8k\nu\text{ is}\\ 
    \text{increasing with }\nu\in[0.5,1)\strut}}}& 3+\frac{4k}{3}-\big(\frac{16k}{3}-\frac{8}{3}\big)\frac{1}{0.5}+8k\cdot0.5=\frac{25}{3}-\frac{16k}{3}\overbracket{\ge}^{\strut\mathclap{k\le4/3}\strut}\frac{25}{3}-\frac{64}{9}=\frac{11}{9}\\
    \end{aligned}
    \end{equation}Then we further compute
    \begin{equation}\label{eq:UpperWQlowerbound}
    \begin{aligned}
    &\mathtt{Q}\big(2-2\nu+\hat{h}_1(U)-2d_1\hat{g}_1(U)\big)\\
    \overbracket{=}^{\mathllap{\text{\eqref{eq:II1II2Quadratic}}}}&\frac{9+4k}{3}+\big(-\frac{16}{3}k+\frac{8}{3}\big)\frac{r\sin\phistar}{d_1}+\frac{8kd_1}{r\sin\phistar}-\frac{4kr\sin\phistar}{d_1}\big(2-2\nu+\hat{h}_1(U)-2d_1\hat{g}_1(U)\big)^2\\
    &+\big[\frac{(32k-12)r\sin\phistar}{3d_1}-8k\big]\big(2-2\nu+\hat{h}_1(U)-2d_1\hat{g}_1(U)\big)\\
    \overbracket{=}^{\mathllap{\nu=\frac{d_1}{r\sin\phistar}}}&\frac{9+4k}{3}+\big(\frac{-16k}{3}+\frac{8}{3}\big)\frac{1}{\nu}+8k\nu-\frac{4k}{\nu}(2-2\nu)^2+\big(-8k+\frac{32k-12}{3}\frac{1}{\nu}\big)\big(2-2\nu\big)\\
    &+\underbracket{\big[\frac{32k-12}{3\nu}-8k\big]\big(\hat{h}_1(U)-2d_1\hat{g}_1(U)\big)-\frac{8k}{\nu}\big(2-2\nu\big)\big(\hat{h}_1(U)-2d_1\hat{g}_1(U)\big)-\frac{4k}{\nu}\big(\hat{h}_1(U)-2d_1\hat{g}_1(U)\big)^2}_{=\mathcal{C}^{\mathcal{W}}_U(\nu,k)}\\    
    =&\frac{9+4k}{3}+32k-16k-\frac{64k-24}{3}+\big(\frac{-16k}{3}+\frac{8}{3}-16k+\frac{64k-24}{3}\big)\frac{1}{\nu}+\big(8k-16k+16k\big)\nu\\
    +&\underbracket{\big[\frac{-16k-12}{3\nu}+8k\big]\big(\hat{h}_1(U)-2d_1\hat{g}_1(U)\big)-\frac{4k}{\nu}\big(\hat{h}_1(U)-2d_1\hat{g}_1(U)\big)^2}_{=\mathcal{C}^{\mathcal{W}}_U(\nu,k)}\\
    =&-4k+11-\frac{16}{3}\frac{1}{\nu}+8k\nu+\mathcal{C}^{\mathcal{W}}_U(\nu,k)\underbracket{>}_{\mathrlap{\substack{\text{Given a } k\in[1,4/3],\\-\frac{16}{3}\frac{1}{\nu}+8k\nu \text{ increases}\\ \text{with }\nu\in[0.5,1]}}\qquad}-4k+11-\frac{16}{3\cdot0.5}+\cdot8k\cdot0.5+\mathcal{C}^{\mathcal{W}}_U(\nu,k)=\frac{1}{3}+\mathcal{C}^{\mathcal{W}}_U(\nu,k),
\end{aligned}\end{equation} where \begin{align*}
    \mathcal{C}^{\mathcal{W}}_U(\nu,k)=\Big[\frac{-16k-12}{3\nu}+8k\Big]\big(\hat{h}_1(U)-2d_1\hat{g}_1(U)\big)-\frac{4k}{\nu}\big(\hat{h}_1(U)-2d_1\hat{g}_1(U)\big)^2
\end{align*}
Hence, by \eqref{eq:UpperQlowerbound}, \eqref{eq:UpperEQlowerbound} and \eqref{eq:UpperWQlowerbound} we have proved \eqref{item03Qlowerboundn1eq1}.
\hidesubsection{Proof of (4)} With $\nu\in(1/3,0.5]$, for $\mathfrak{l}_L$ in \cref{fig:x1_n1equal1case}, we define length functions limits on its endpoints $\mathcal{W}_L$ and $\mathcal{E}_L$ \newline as follows.\COMMENT{Why not ``$\hat{\tau}_i(L)\dfn\lim_{\mathfrak{l}_L\ni x\to \mathcal{W}_L}\tau_i\ge0$'' etc. below?}

\smallskip\begin{minipage}{.47\textwidth}
    Limit at $\mathcal{W}_L$:
    \begin{equation}\label{eq:lengthlimitWL}
        \begin{aligned}
        \hat{\tau}_0(L)&\dfn\lim_{\mathfrak{l}_L\ni x\to \mathcal{W}_L}\tau_0\ge0\\
        \hat{\tau}_1(L)&\dfn\lim_{\mathfrak{l}_L\ni x\to \mathcal{W}_L}\tau_1\ge0\\
        \hat{d}_0(L)&\dfn\lim_{\mathfrak{l}_L\ni x\to \mathcal{W}_L}d_0>0\\
        \hat{d}_1(L)&\dfn\lim_{\mathfrak{l}_L\ni x\to \mathcal{W}_L}d_1=\const=\nu\cdot r\sin{\phistar}\\
        \hat{d}_2(L)&\dfn\lim_{\mathfrak{l}_L\ni x\to \mathcal{W}_L}d_2>0
    \end{aligned}
    \end{equation}
    \end{minipage}\vrule\quad
    \begin{minipage}{0.47\textwidth}
    Limit at $\mathcal{E}_L$:
    \begin{equation}\label{eq:lengthlimitEL}
        \begin{aligned}
        \tilde{\tau}_0(L)\dfn&\lim_{\mathfrak{l}_L\ni x\to \mathcal{E}_L}\tau_0\ge0\\
        \tilde{\tau}_1(L)\dfn&\lim_{\mathfrak{l}_L\ni x\to \mathcal{E}_L}\tau_1\ge0\\
        \tilde{d}_0(L)\dfn&\lim_{\mathfrak{l}_L\ni x\to \mathcal{E}_L}d_0>0\\
        \tilde{d}_1(L)\dfn&\lim_{\mathfrak{l}_L\ni x\to \mathcal{E}_L}d_1=\nu\cdot r\sin{\phistar}\\
        \tilde{d}_2(L)\dfn&\lim_{\mathfrak{l}_L\ni x\to \mathcal{E}_L}d_2>0
    \end{aligned}
    \end{equation}
    \end{minipage}
    
    These limits in \eqref{eq:lengthlimitWL}, \eqref{eq:lengthlimitEL} exist for the same reason as for the existence of limits in \eqref{eq:lengthlimitWU}. And these defined lengths continuously extended for $x$\COMMENT{Edit} at the end points $\mathcal{W}_U$, $\mathcal{E}_U$ can be visualized in \cref{fig:n1_equal1_lower_leftendpoint_case,fig:n1_equal1_lower_rightendpoint_case}. 
    
    For the continuous extension to $\mathcal{W}_L$,  we observe that for the following. 
    
    $x\in\mathfrak{l}_L$, if $x\rightarrow\mathcal{W}_L$, then $p(\mathcal{F}^2(x))\rightarrow B$ which is the corner and $\Fcal(x)\rightarrow\Fcal(\mathcal{W}_L)$ which is an interior point of $M_R$ (see \cref{fig:x1_n1equal1case}), that is, $p(\Fcal(\mathcal{W}_L))$ is an interior point of $\Gamma_R$ (see \cref{fig:n1_equal1_lower_leftendpoint_case}).
    
    If $x\rightarrow\mathcal{E}_L$, then $p(\mathcal{F}^{-1}(x))\rightarrow A$ which is the corner and $x\rightarrow\mathcal{E}_L$ which is an interior point of $M_R$ (see \cref{fig:x1_n1equal1case}), that is, $p(\mathcal{E}_L)$ is an interior point of $\Gamma_R$ (see \cref{fig:n1_equal1_lower_rightendpoint_case}).
    Hence \begin{equation}\label{eq:tau0Ltau1Llimits}
        \begin{aligned}
        \hat{\tau}_1(L)=&\lim_{\mathfrak{l}_L\ni x\to \mathcal{W}_L}\overbracket{|p(\Fcal(x))p(\Fcal^2(x))|}^{=\tau_1\bigm|_x}=|p(\Fcal(\mathcal{W}_L))B|\overbracket{=}^{\qquad\qquad\mathllap{B\text{, } p(\Fcal(\mathcal{W}_L))\text{ are on }\Gamma_R}} 2\hat{d}_1(L)=2d_1\\
        \tilde{\tau}_0(L)=&\lim_{\mathfrak{l}_L\ni x\to \mathcal{E}_L}\overbracket{|p(x)p(\Fcal^{-1}(x))|}^{{=\tau_0\bigm|_x}}=|p(\mathcal{E}_L)A|\overbracket{=}^{\qquad\qquad\mathllap{A\text{, } p(\mathcal{E}_L)\text{ are on }\Gamma_R}} 2\tilde{d}_1(L)=2d_1
        \end{aligned}
    \end{equation}
    
    Similarly to \eqref{eq:hhatUg1hatU}, we define
    \begin{equation}\label{eq:htildehatLg1tildehatL}
    \begin{aligned}
    \hat{h}_1(L)\dfn&\lim_{\mathfrak{l}_U\ni x\to \mathcal{W}_L}h_1(x)\overbracket{=}^{\mathllap{\text{\cref{def:lengthfractionfunctionsn1eq1}\eqref{item01lengthfractionfunctionsn1eq1}}}}\lim_{\mathfrak{l}_U\ni x\to \mathcal{W}_L}\big(\frac{\tau_0}{d_0}+\frac{\tau_1}{d_0}+\frac{2d_1}{d_0}-2\big)\overset{\eqref{eq:lengthlimitWL}}=\frac{\hat{\tau}_0(L)}{\hat{d}_0(L)}+\frac{\hat{\tau}_1(L)}{\hat{d}_0(L)}+\frac{2d_1}{\hat{d}_0(L)}-2\\
    \hat{g}_1(L)\dfn&\lim_{\mathfrak{l}_L\ni x\to \mathcal{W}_L}g_1(x)\overbracket{=}^{\mathllap{\text{\cref{def:lengthfractionfunctionsn1eq1}\eqref{item03lengthfractionfunctionsn1eq1}}}}\lim_{\mathfrak{l}_U\ni x\to \mathcal{W}_L}\big(\frac{1}{d_0}-\frac{1}{r\sin\phistar}\big)\overset{\eqref{eq:lengthlimitWL}}=\frac{1}{\hat{d}_0(L)}-\frac{1}{r\sin\phistar}\\
    \tilde{h}_1(L)\dfn&\lim_{\mathfrak{l}_L\ni x\to \mathcal{E}_L}h_1(x)\overbracket{=}^{\mathllap{\text{\cref{def:lengthfractionfunctionsn1eq1}\eqref{item01lengthfractionfunctionsn1eq1}}}}\lim_{\mathfrak{l}_U\ni x\to \mathcal{E}_L}\big(\frac{\tau_0}{d_0}+\frac{\tau_1}{d_0}+\frac{2d_1}{d_0}-2\big)\overset{\eqref{eq:lengthlimitEL}}=\frac{\tilde{\tau}_0(L)}{\tilde{d}_0(L)}+\frac{\tilde{\tau}_1(L)}{\tilde{d}_0(L)}+\frac{2d_1}{\tilde{d}_0(L)}-2\\
    \tilde{g}_1(L)\dfn&\lim_{\mathfrak{l}_L\ni x\to \mathcal{E}_L}g_1(x)\overbracket{=}^{\mathllap{\text{\cref{def:lengthfractionfunctionsn1eq1}\eqref{item03lengthfractionfunctionsn1eq1}}}}\lim_{\mathfrak{l}_L\ni x\to \mathcal{E}_L}\big(\frac{1}{d_0}-\frac{1}{r\sin\phistar}\big)\overset{\eqref{eq:lengthlimitEL}}=\frac{1}{\tilde{d}_0(L)}-\frac{1}{r\sin\phistar}
    \end{aligned}
\end{equation}

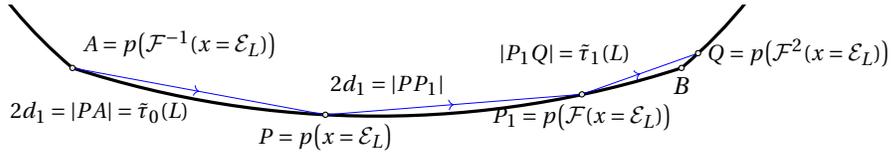
\begin{figure}[h]
\begin{center}
\begin{tikzpicture}[xscale=1.4,yscale=1.4]
    \tkzDefPoint(0,0){Or};
    \tkzDefPoint(0,-4.8){Y};
    \clip(-5.8,-4.23) rectangle (5.8,-2.85);
    \pgfmathsetmacro{\rradius}{4.5};
    \pgfmathsetmacro{\phistardeg}{40};
    \pgfmathsetmacro{\XValueArc}{\rradius*sin(\phistardeg)};
    \pgfmathsetmacro{\YValueArc}{\rradius*cos(\phistardeg)};
    \pgfmathsetmacro{\Rradius}{
        \rradius*sin(\phistardeg)/sin(18)
    };
    \pgfmathsetmacro{\bdist}{(-\rradius*cos(\phistardeg))+cos(18)*\Rradius}
    \tkzDefPoint(0,\bdist){OR};
    \tkzDefPoint(\XValueArc,-1.0*\YValueArc){B};    \tkzDefPoint(-1.0*\XValueArc,-1.0*\YValueArc){A};
    \tkzDrawArc[name path=Cr,very thick](Or,B)(A);
    \tkzDrawArc[name path=CR,very thick](OR,A)(B);
    \pgfmathsetmacro{\PXValue}{-\Rradius*sin(3)};
    \pgfmathsetmacro{\PYValue}{\bdist-\Rradius*cos(3)};
    \tkzDefPoint(\PXValue,\PYValue){P1};
    
    \tkzDefPointBy[projection=onto OR--P1](A)\tkzGetPoint{M};
    \coordinate (P2) at ($(A)!2!(M)$);
    
    \tkzDefPointBy[projection=onto OR--P2](P1)\tkzGetPoint{M2};
    \coordinate (Q_far) at ($(P1)!2!(M2)$);
    \path[name path=Line_P2_Q_far](P2)-- (Q_far);
    \path[name intersections={of=Cr and Line_P2_Q_far,by={Q}}];
    \tkzDefPoint(0,\bdist-\Rradius){C};
    \begin{scope}[ultra thin,decoration={markings,mark=at position 0.5 with {\arrow{>}}}] 
    \draw[blue,postaction={decorate}] (A)--(P1);
    \draw[blue,postaction={decorate}] (P1)--(P2);
    \draw[blue,postaction={decorate}] (P2)--(Q);
    \end{scope}
    \tkzLabelPoint[above right](A){\small $A=p\big(\mathcal{F}^{-1}(x=\mathcal{E}_L)\big)$};
    \tkzLabelPoint[below](P1){\small $P=p\big(x=\mathcal{E}_L\big)$};
    \tkzLabelPoint[below](P2){\small $P_1=p\big(\mathcal{F}(x=\mathcal{E}_L)\big)$};
    \tkzLabelPoint[right](Q){\small $Q=p\big(\mathcal{F}^{2}(x=\mathcal{E}_L)\big)$};
    \tkzLabelPoint[below](B){$B$};
    \tkzLabelPoint[below](Or){$O_r$};
    \tkzLabelSegment[above left=0.05pt and 0.03pt](P1,P2){\small $2d_1=|PP_1|$};
    \tkzLabelSegment[below left=0.5pt and 0.3pt](A,P1){\small $2d_1=|PA|=\tilde{\tau}_0(L)$};
     \tkzLabelSegment[above left=0.1pt and 0.2pt](P2,Q){\small $|P_1Q|=\tilde{\tau}_1(L)$};
    \tkzDrawPoints(A,B,Or,P1,P2,Q);
    
\end{tikzpicture}
\caption{With $d_1\le0.5r\sin\phistar$, $(\Phi_{\mathcal{E}_L},\theta_1)=x=\mathcal{E}_L$ is on a trajectory from corner $A$}\label{fig:n1_equal1_lower_rightendpoint_case}
\end{center}
\end{figure}

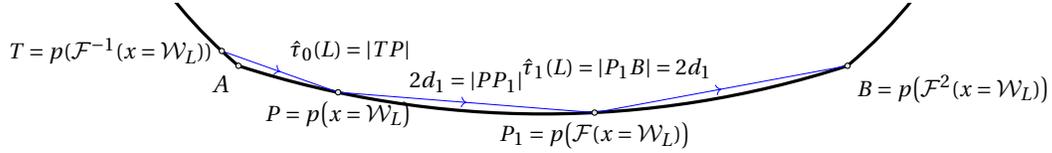
\begin{figure}[h]
\begin{center}
\begin{tikzpicture}[xscale=1.4,yscale=1.4]
    \tkzDefPoint(0,0){Or};
    \tkzDefPoint(0,-4.8){Y};
    \clip(-5.7,-4.23) rectangle (5.7,-2.85);
    \pgfmathsetmacro{\rradius}{4.5};
    \pgfmathsetmacro{\phistardeg}{40};
    \pgfmathsetmacro{\XValueArc}{\rradius*sin(\phistardeg)};
    \pgfmathsetmacro{\YValueArc}{\rradius*cos(\phistardeg)};
    \pgfmathsetmacro{\Rradius}{
        \rradius*sin(\phistardeg)/sin(18)
    };
    \pgfmathsetmacro{\bdist}{(-\rradius*cos(\phistardeg))+cos(18)*\Rradius}
    \tkzDefPoint(0,\bdist){OR};
    \tkzDefPoint(\XValueArc,-1.0*\YValueArc){B};    \tkzDefPoint(-1.0*\XValueArc,-1.0*\YValueArc){A};
    \tkzDrawArc[name path=Cr, very thick](Or,B)(A);
    \tkzDrawArc[name path=CR, very thick](OR,A)(B);
    \pgfmathsetmacro{\PXValue}{\Rradius*sin(3)};
    \pgfmathsetmacro{\PYValue}{\bdist-\Rradius*cos(3)};
    \tkzDefPoint(\PXValue,\PYValue){P2};
    
    \tkzDefPointBy[projection=onto OR--P2](B)\tkzGetPoint{M};
    \coordinate (P1) at ($(B)!2!(M)$);
    
    \tkzDefPointBy[projection=onto OR--P1](P2)\tkzGetPoint{M2};
    \coordinate (T_far) at ($(P2)!2!(M2)$);
    \path[name path=Line_P1_T_far](P1)-- (T_far);
    \path[name intersections={of=Cr and Line_P1_T_far,by={T}}];
    \tkzDefPoint(0,\bdist-\Rradius){C};
    \begin{scope}[ultra thin,decoration={markings,mark=at position 0.5 with {\arrow{>}}}] 
    \draw[blue,postaction={decorate}] (T)--(P1);
    \draw[blue,postaction={decorate}] (P1)--(P2);
    \draw[blue,postaction={decorate}] (P2)--(B);
    \end{scope}
    \tkzLabelPoint[below left](A){$A$};
    \tkzLabelPoint[below](P1){\small $P=p\big(x=\mathcal{W}_L\big)$};
    \tkzLabelPoint[below](P2){\small $P_1=p\big(\mathcal{F}(x=\mathcal{W}_L)\big)$};
    \tkzLabelPoint[left](T){\small $T=p(\mathcal{F}^{-1}(x=\mathcal{W}_L))$};
    \tkzLabelPoint[below right](B){\small $B=p\big(\mathcal{F}^2(x=\mathcal{W}_L)\big)$};
    \tkzLabelPoint[below](Or){$O_r$};
    \tkzLabelSegment[above =0.5pt](P1,P2){\small $2d_1=|PP_1|$};
    \tkzLabelSegment[above right= 0.5pt and 0.2pt](P1,T){\small $\hat{\tau}_0(L)=|TP|$};
    \tkzLabelSegment[above left= 0.5pt and 0.02pt](P2,B){\small $\hat{\tau}_1(L)=|P_1B|=2d_1$};
    \tkzDrawPoints(A,B,P1,P2,T);
\end{tikzpicture}
\caption{With $d_1\le0.5r\sin\phistar$, $(\Phi_{\mathcal{W}_L},\theta_1)=x=\mathcal{W}_L$ is on a trajectory to corner $B$}\label{fig:n1_equal1_lower_leftendpoint_case}


\end{center}
\end{figure}

By the previous monotone property item \eqref{item02Qlowerboundn1eq1}, we have
\[\frac{\tilde{\tau}_1(L)}{\tilde{d}_0(L)}\overbracket{<}^{\qquad\qquad\qquad\mathllap{\substack{\text{item \eqref{item02Qlowerboundn1eq1}: }\\\tilde{d}_0(L)>d_0|_{x=x_1}}}}\frac{\tilde{\tau}_1(L)}{d_0|_{x=x_1}}\overbracket{<}^{\qquad\qquad\qquad\mathllap{\substack{\text{item \eqref{item02Qlowerboundn1eq1}: }\\\tilde{\tau}_1(L)<\tau_1|_{x=x_1}}}}\frac{\tau_1}{d_0}\Bigm|_{x=x_1}\overbracket{<}^{\qquad\qquad\qquad\mathllap{\substack{\text{item \eqref{item02Qlowerboundn1eq1}: }\\\hat{\tau}_1(L)>\tau_1|_{x=x_1}}}}\frac{\hat{\tau}_1(L)}{d_0\bigm|_{x=x_1}}\overbracket{<}^{\qquad\qquad\qquad\mathllap{\substack{\text{item \eqref{item02Qlowerboundn1eq1}: }\\\hat{d}_0(L)<d_0|_{x=x_1}}}}\frac{\hat{\tau}_1(L)}{\hat{d}_0(L)}\]
And for these two bounds we have 
\begin{align*}
    \frac{\hat{\tau}_1(L)}{\hat{d}_0(L)}\overbracket{=}^{\eqref{eq:tau0Ltau1Llimits}}&\frac{2d_1}{\hat{d}_0(L)}\overbracket{=}^{\eqref{eq:htildehatLg1tildehatL}}2\big(\frac{d_1}{r\sin{\phistar}}+d_1\hat{g}_1(L)\big)\overbracket{=}^{d_1=\nu r\sin\phistar}2\nu+2d_1\hat{g}_1(L)\\
    \frac{\tilde{\tau}_1(L)}{\tilde{d}_0(L)}\overbracket{=}^{\eqref{eq:htildehatLg1tildehatL}}&2+\tilde{h}_1(L)-\frac{2d_1}{\tilde{d}_0(L)}-\frac{\tilde{\tau}_0(L)}{\tilde{d}_0(L)}\\
    \overbracket{=}^{\qquad\mathllap{\text{\eqref{eq:tau0Ltau1Llimits}: }\tilde{\tau}_0(L)=2d_1}}&2+\tilde{h}_1(L)-\frac{4d_1}{\tilde{d}_0(L)}\\
    \overbracket{=}^{\eqref{eq:htildehatLg1tildehatL}}&2+\tilde{h}_1(L)-4\big(\frac{d_1}{r\sin{\phistar}}+d_1\tilde{g}_1(L)\big)\overbracket{=}^{\qquad\mathllap{d_1=\nu r\sin\phistar}}2-4\nu+\tilde{h}_1(L)-4d_1\tilde{g}_1(L)\\
\end{align*}

With fixed $k\in[1,4/3]$ and constant $\nu=\frac{d_1}{r\sin\phistar}$, the function $\mathtt{Q}(t)$ from \eqref{eq:II1II2Quadratic} for $t\in (-\infty,+\infty)$ has only one extreme point, which is the global maximum point. Therefore, for $t\in\big[2-4\nu+\tilde{h}_1(L)-4d_1\tilde{g}_1(L),\,\,2\nu+2d_1\hat{g}_1(L)\big]$, $\mathtt{Q}\big(t\big)$ cannot reach the minimum in the interior of $\big[2-4\nu+\tilde{h}_1(L)-4d_1\tilde{g}_1(L),\,\,2\nu+2d_1\hat{g}_1(L)\big]$. Hence,\begin{equation}\label{eq:LowerQlowerbound}
    \mathtt{Q}\Big(\frac{\tau_1}{d_0}\bigm|_{x=x_1}\Big)>\min_{t\in\big[2-4\nu+\tilde{h}_1(L)-4d_1\tilde{g}_1(L),\,\,2\nu+2d_1\hat{g}_1(L)\big]}\mathtt{Q}\big(t\big)=\min{\Big\{\mathtt{Q}\big(2-4\nu+\tilde{h}_1(L)-4d_1\tilde{g}_1(L)\big),\,\,\mathtt{Q}\big(2\nu+2d_1\hat{g}_1(L)\big)\Big\}}
\end{equation}
With $\frac{d_1}{r\sin\phistar}=\nu\in(1/3,0.5]$ we firstly compute $\mathtt{Q}\big(2-4\nu+\tilde{h}_1(L)-4d_1\tilde{g}_1(L)\big)$.\COMMENT{Use multlined. Ask how if needed}
\begin{equation}\label{eq:LowerEQlowerboundSimplified}
\begin{aligned}
&\mathtt{Q}\big(2-4\nu+\tilde{h}_1(L)-4d_1\tilde{g}_1(L)\big)\\
\overbracket{=}^{\mathllap{\text{\eqref{eq:II1II2Quadratic}}}}&\frac{4k+9}{3}+\big(\frac{8-16k}{3}\big)\frac{r\sin{\phistar}}{d_1}+\frac{8kd_1}{r\sin{\phistar}}-\frac{4kr\sin\phistar}{d_1}\big(2-4\nu+\tilde{h}_1(L)-4d_1\tilde{g}_1(L)\big)^2\\
&+\big[\frac{(32k-12)r\sin\phistar}{3d_1}-8k\big]\big(2-4\nu+\tilde{h}_1(L)-4d_1\tilde{g}_1(L)\big)\\
=&\frac{9+4k}{3}+\frac{8-16k}{3\nu}+8k\nu-\frac{4k(2-4\nu)^2}{\nu}-\frac{8k(2-4\nu)}{\nu}\big(\tilde{h}_1(L)-4d_1\tilde{g}_1(L)\big)\\
&-\frac{4k\big(\tilde{h}_1(L)-4d_1\tilde{g}_1(L)\big)^2}{\nu}+\big(\frac{32k-12}{3\nu}-8k\big)\big(2-4\nu\big)+\big(\frac{32k-12}{3\nu}-8k\big)\big(\tilde{h}_1(L)-4d_1\tilde{g}_1(L)\big)\\
=&\frac{9+4k}{3}+64k-16k-\frac{128k-48}{3}+\big(\frac{8-16k}{3}-16k+\frac{64k-24}{3}\big)\frac{1}{\nu}+\big(8k-64k+32k\big)\nu\\
&+\big(-8k+\frac{32k-12}{3\nu}\big)\big(\tilde{h}_1(L)-4d_1\tilde{g}_1(L)\big)-\frac{8k}{\nu}\big(\tilde{h}_1(L)-4d_1\tilde{g}_1(L)\big)\big(2-4\nu\big)-\frac{4k}{\nu}\big(\tilde{h}_1(L)-4d_1\tilde{g}_1(L)\big)^2\\
=&19+\frac{20k}{3}-\frac{16}{3\nu}-24k\nu+\big(24k+\frac{-16k-12}{3\nu}\big)\big(\tilde{h}_1(L)-4d_1\tilde{g}_1(L)\big)-\frac{4k}{\nu}\big(\tilde{h}_1(L)-4d_1\tilde{g}_1(L)\big)^2\\
=&19+\frac{20k}{3}-\frac{16}{3\nu}-24k\nu+ \mathcal{C}^{\mathcal{E}}_{L}(\nu,k),\\
\end{aligned}
\end{equation}
where $\mathcal{C}^{\mathcal{E}}_{L}(\nu,k)=\big(24k+\frac{-16k-12}{3\nu}\big)\big(\tilde{h}_1(L)-4d_1\tilde{g}_1(L)\big)-\frac{4k}{\nu}\big(\tilde{h}_1(L)-4d_1\tilde{g}_1(L)\big)^2$.

For fixed $k\in[1,4/3]$, since $19+\frac{20k}{3}-\frac{16}{3\nu}-24k\nu$ as a function of $\nu\in(0,+\infty)$ has only one extreme point, which is a maximum point. Therefore, for $\nu\in[1/3,1/2]$, $19+\frac{20k}{3}-\frac{16}{3\nu}-24k\nu$ cannot reach a minimum at an interior point. So, it reaches the minimum when $\nu=1/2$ or $\nu=1/3$. We get the following.\COMMENT{Is the $\min_{k\in[1,4/3]}$ needed in the next computation? It does not seem so. Cut the clutter. Also, remove all \textbackslash limits.}
\begin{equation}\label{eq:LowerEQlowerbound}
\begin{aligned}
&\min_{\substack{k\in[1,4/3]\\\nu\in[1/3,1/2]}}{\Big\{\overbracket{19+\frac{20k}{3}-\frac{16}{3\nu}-24k\nu}^{\qquad\qquad\qquad\qquad\qquad\mathllap{\text{minimum attained }\text{at } \nu=1/2\text{ or }1/3}}\Big\}}\\
=&\min_{k\in[1,4/3]}{\Big\{\min{\big\{19+\frac{20k}{3}-\frac{16}{3\cdot0.5}-24k\cdot0.5,19+\frac{20k}{3}-\frac{16}{3\cdot(1/3)}-24k\cdot(1/3)\big\}}\Big\}}\\
=&\min_{k\in[1,4/3]}{\Big\{\min{\big\{\frac{25}{3}-\frac{16k}{3},3-\frac{4k}{3}\big\}}\Big\}}\overset{k\le4/3}=\frac{11}{9}.\\\text{Hence }&
\mathtt{Q}\big(2-4\nu+\tilde{h}_1(L)-4d_1\tilde{g}_1(L)\big)\overbracket{=}^{\mathllap{\eqref{eq:LowerEQlowerboundSimplified}}}19+\frac{20k}{3}-\frac{16}{3\nu}-24k\nu+\mathcal{C}^{\mathcal{E}}_{L}(\nu,k)\overbracket{\ge}^{\mathllap{k\in[1,4/3],\nu\in(1/3,1/2]}}\frac{11}{9}+\mathcal{C}^{\mathcal{E}}_{L}(\nu,k).
\end{aligned}
\end{equation}

We then compute $\mathtt{Q}\big(2\nu+2d_1\hat{g}_1(L)\big)$.
\begin{equation}\label{eq:LowerWQlowerboundSimplified}
    \begin{aligned}
    \mathtt{Q}\big(2\nu+2d_1\hat{g}_1(L)\big)\overbracket{=}^{\mathllap{\text{\eqref{eq:II1II2Quadratic}}}}&\frac{4k+9}{3}+\big(\frac{8-16k}{3}\big)\frac{r\sin{\phistar}}{d_1}+\frac{8kd_1}{r\sin{\phistar}}-\frac{4kr\sin\phistar}{d_1}\big(2\nu+2d_1\hat{g}_1(L)\big)^2\\
        &+\big[\frac{(32k-12)r\sin\phistar}{3d_1}-8k\big]\big(2\nu+2d_1\hat{g}_1(L)\big)\\
        \overbracket{=}^{\mathllap{d_1=\nu r\sin\phistar}}&\frac{9+4k}{3}+\big(\frac{8-16k}{3\nu}\big)+8k\nu+\big(\frac{32k-12}{3\nu}-8k\big)\big(2\nu+2d_1\hat{g}_1(L)\big)-\frac{4k}{\nu}\big(2\nu+2d_1\hat{g}_1(L)\big)^2\\
        =&\frac{9+4k}{3}+\frac{64k-24}{3}+\big(\frac{8-16k}{3\nu}\big)+(8k-16k-16k)\nu\\
        &+\underbracket{\big(\frac{32k-12}{3\nu}-8k\big)\cdot2d_1\hat{g}_1(L)-\frac{4k\cdot8\nu d_1\hat{g}_1(L)}{\nu}-\frac{16kd_1^2(\hat{g}_1(L))^2}{\nu}}_{=\mathcal{C}^{\mathcal{W}}_{L}(\nu,k)}\\
        =&\frac{68k}{3}-5-\big(\frac{16k-8}{3\nu}\big)-24k\nu+\underbracket{\big(\frac{32k-12}{3\nu}-24k\big)\cdot2d_1\hat{g}_1(L)-\frac{16kd_1^2(\hat{g}_1(L))^2}{\nu}}_{\mathcal{C}^{\mathcal{W}}_{L}(\nu,k)}\\
        =&\frac{68k}{3}-5-\big(\frac{16k-8}{3\nu}\big)-24k\nu+\mathcal{C}^{\mathcal{W}}_{L}(\nu,k),\\
    \end{aligned}
\end{equation}
where $\mathcal{C}^{\mathcal{W}}_{L}(\nu,k)=\big(\frac{32k-12}{3\nu}-24k\big)\cdot2d_1\hat{g}_1(L)-\frac{16kd_1^2(\hat{g}_1(L))^2}{\nu}$

For fixed $k\in[1,4/3]$, since $\frac{68k}{3}-5-\big(\frac{16k-8}{3\nu}\big)-24k\nu$ as a function of $\nu\in(0,+\infty)$ has only one extreme point, which is a maximum point. Therefore, for $\nu\in[1/3,1/2]$, $\frac{68k}{3}-5-\big(\frac{16k-8}{3\nu}\big)-24k\nu$ cannot reach a minimum at an interior point. So, it reaches the minimum when $\nu=1/2$ or $\nu=1/3$. We get the following.
\begin{equation}\label{eq:LowerWQlowerbound}
\begin{aligned}
&\min_{\substack{k\in[1,4/3]\\\nu\in[1/3,1/2]}}{\Big\{\overbracket{\frac{68k}{3}-5-\big(\frac{16k-8}{3\nu}\big)-24k\nu}^{{\mathllap{\text{minimum attained }\text{at } \nu=1/2\text{ or }1/3}}}\Big\}}\\
=&\min_{k\in[1,4/3]}{\Big\{\min{\big\{\frac{68k}{3}-5-\frac{16k-8}{3\cdot0.5}-24k\cdot0.5,\frac{68k}{3}-5-\big(\frac{16k-8}{3\cdot(1/3)}\big)-24k\cdot(1/3)\big\}}\Big\}}\\
=&\min_{k\in[1,4/3]}{\Big\{\min{\big\{\frac{1}{3},3-\frac{4k}{3}\big\}}\Big\}}\overbracket{=}^{k\le4/3}\frac{1}{3}.\\
\text{Hence }&\mathtt{Q}\big(2\nu+2d_1\hat{g}_1(L)\big)\overbracket{=}^{\mathllap{\eqref{eq:LowerWQlowerboundSimplified}}}\frac{68k}{3}-5-\big(\frac{16k-8}{3\nu}\big)-24k\nu+\mathcal{C}^{\mathcal{W}}_{L}(\nu,k)\overbracket{\ge}^{\mathllap{k\in[1,4/3],\nu\in(1/3,1/2]}}\frac{1}{3}+\mathcal{C}^{\mathcal{W}}_{L}(\nu,k).
\end{aligned}    
\end{equation}

By \eqref{eq:LowerQlowerbound}, \eqref{eq:LowerEQlowerbound} and \eqref{eq:LowerWQlowerbound} we have proved \eqref{item04Qlowerboundn1eq1}.
\hidesubsection{Proof of (5)} We first look at $\mathcal{C}^{\mathcal{W}}_{U}(\nu,k)$. Note that its parameter ranges are $k\in[1,4/3]$, $\nu\in(1/2,1)$.
Therefore, we have \begin{equation}\label{eq:2ndproducttermCUWnuk}
    \begin{aligned}
    &\frac{-8k-12}{3}=\frac{-16k-12}{3\times0.5}+8k\overbracket{<}^{\qquad\mathllap{\substack{-16k-12<0,\\\nu>0.5}}}\frac{-16k-12}{3\nu}+8k\overbracket{<}^{\qquad\mathllap{\substack{-16k-12<0,\\\nu<1}}}\frac{-16k-12}{3\times1}+8k=\frac{8k-12}{3}\overset{k\le4/3}<0 
\end{aligned}
\end{equation}
On the other hand, by \eqref{eq:lengthlimitWU}, \eqref{eq:hhatUg1hatU} $\hat{h}_1(U)$, $\hat{g}_1(U)$ are the limits of $h_1(x)$, $g_1(x)$ which are continuous functions of $x\in \mathfrak{l}_U$. Hence we find
\begin{equation}\label{eq:hat1Ug1Ubounds}
\begin{aligned}
&\text{\cref{lemma:h1g1g0g1g2g3lengthfractionestimate}\eqref{item01h1g1g0g1g2g3lengthfractionestimate}: } \forall x\in\mathfrak{l}_U\text{, }h_1(x)\in\big(-8.53r/R,\quad16.65r/R\big)\\
\Longrightarrow&\hat{h}_1(U)\in\big[-8.53r/R,\quad 16.65r/R\big]\\
&\text{and \cref{lemma:h1g1g0g1g2g3lengthfractionestimate}\eqref{item03h1g1g0g1g2g3lengthfractionestimate}: } \forall x\in\mathfrak{l}_U\text{, }g_1(x)\in\big(-\frac{17r}{4R}\frac{1}{0.9975r\sin\phistar},\quad \frac{17r}{4R}\frac{1}{0.9975r\sin\phistar}\big)\\
\Longrightarrow&\hat{g}_1(U)\in\big[-\frac{17r}{4R}\frac{1}{0.9975r\sin\phistar},\quad \frac{17r}{4R}\frac{1}{0.9975r\sin\phistar}\big],\\
\Longrightarrow&-2d_1\hat{g}_1(U)\in\big[-\frac{17r}{2R}\frac{d_1}{0.9975r\sin\phistar},\quad \frac{17r}{2R}\frac{d_1}{0.9975r\sin\phistar}\big]\overset{\nu=\frac{d_1}{r\sin\phistar}}=\big[-\frac{8.5\nu}{0.9925}\frac{r}{R},\quad +\frac{8.5\nu}{0.9975}\frac{r}{R}\big].\\
\end{aligned}\end{equation}Hence 
\begin{equation}\label{eq:1stproducttermCUWnuk}
\hat{h}_1(U)-2d_1\hat{g}_1(U)\in\Big[-\big(8.53+\frac{8.5\nu}{0.9975}\big)\frac{r}{R},\quad \big(16.65+\frac{8.5\nu}{0.9975}\big)\frac{r}{R}\Big]
\end{equation}
Therefore, 
\begin{equation}\label{eq:CWUnukbound}
\begin{aligned}
    \mathcal{C}^{\mathcal{W}}_{U}(\nu,k)=&\Big(\frac{-16k-12}{3\nu}+8k\Big)\big(\hat{h}_1(U)-2d_1\hat{g}_1(U)\big)-\frac{4k}{\nu}\big(\hat{h}_1(U)-2d_1\hat{g}_1(U)\big)^2\\
    \overbracket{>}^{\mathllap{\text{\eqref{eq:1stproducttermCUWnuk},\eqref{eq:2ndproducttermCUWnuk}}}}&\big(\underbracket{\frac{-8k-12}{3}}_{<0}\big)\big(16.65+\frac{8.5\nu}{0.9975}\big)\frac{r}{R}-\frac{4k}{\nu}\Big[\big(16.65+\frac{8.5\nu}{0.9975}\big)\frac{r}{R}\Big]^2\\
    \overbracket{>}^{\mathllap{\nu<1}}&\big(\frac{-8k-12}{3}\big)\big(16.65+\frac{8.5}{0.9975}\big)\frac{r}{R}-4k\underbracket{\big(\frac{16.65}{\nu}+\frac{8.5}{0.9975}\big)\big(16.65+\frac{8.5\nu}{0.9975}\big)}_{<875 \text{ with }0.5\le\nu\le1}\big(\frac{r}{R}\big)^2\\
    >&-(68k+101)\frac{r}{R}-3500k\big(\frac{r}{R}\big)^2\overset{k\le4/3}>-192\frac{r}{R}-4667\big(\frac{r}{R}\big)^2
\end{aligned}
\end{equation}
We secondly look at $\mathcal{C}^{\mathcal{W}}_L(\nu,k)$. Note that its parameter ranges are $k\in[1,4/3]$, $\nu\in(1/3,1/2]$. Therefore, we have \begin{equation}\label{eq:1stproducttermCWLnuk}
\begin{aligned}
\frac{-8k-24}{3}=-24k+\frac{32k-12}{3\cdot0.5}\overbracket{\le}^{\mathrlap{32k-12>0,\nu\le1/2}\qquad}-24k+\frac{32k-12}{3\nu}\overbracket{<}^{\mathrlap{32k-12>0,\nu>1/3}}-24k+32k-12=8k-12\overset{k\le4/3}<0    
\end{aligned}    
\end{equation}
by \eqref{eq:lengthlimitWL}, \eqref{eq:htildehatLg1tildehatL} $\hat{g}_1(L)$ is the limit of $h_1(x)$,  which is a continuous function of $x\in \mathfrak{l}_L$. Hence we find
\begin{equation}\label{eq:2ndproducttermCWLnuk}
\begin{aligned}&\text{\cref{lemma:h1g1g0g1g2g3lengthfractionestimate}\ref{item03h1g1g0g1g2g3lengthfractionestimate}: }\forall x\in\mathfrak{l}_L\text{, } g_1(x)\in \big(-\frac{17r}{4R}\frac{1}{0.9975r\sin\phistar},\quad\frac{17r}{4R}\frac{1}{0.9975r\sin\phistar}\big)\\
\Longrightarrow&\hat{g}_1(L)\in\big[-\frac{17r}{4R}\frac{1}{0.9975r\sin\phistar},\quad\frac{17r}{4R}\frac{1}{0.9975r\sin\phistar}\big]\\
\Longrightarrow&2d_1\hat{g}_1(L)\in\big[-\frac{17r}{2R}\frac{d_1}{0.9975r\sin\phistar},\quad\frac{17r}{2R}\frac{d_1}{0.9975r\sin\phistar}\big]\overset{\nu=\frac{d_1}{r\sin\phistar}}=\big[-\frac{8.5\nu}{0.9975}\frac{r}{R},\quad\frac{8.5\nu}{0.9975}\frac{r}{R}\big],
\end{aligned}
\end{equation}
Therefore \begin{equation}\label{eq:CWLnukbound}
    \begin{aligned}
        \mathcal{C}^{\mathcal{W}}_{L}(\nu,k)=&\underbracket{\big(\frac{32k-12}{3\nu}-24k\big)\cdot2d_1\hat{g}_1(L)}_{\eqref{eq:1stproducttermCWLnuk}: <0}-\frac{16kd_1^2(\hat{g}_1(L))^2}{\nu}\\
        \overbracket{>}^{\mathllap{\text{\eqref{eq:1stproducttermCWLnuk},\eqref{eq:2ndproducttermCWLnuk}}}}&\frac{-8k-24}{3}\frac{8.5\nu}{0.9975}\frac{r}{R}-\frac{4k}{\nu}\big(\frac{8.5}{0.9975}\frac{r}{R}\big)^2\nu^2\\
        \overbracket{>}^{\mathllap{1/3<\nu\le0.5}}&\frac{\overbracket{-8k-24}^{<0}}{3}\frac{8.5\times0.5}{0.9975}\frac{r}{R}-4k\cdot0.5\cdot\big(\frac{8.5}{0.9975}\frac{r}{R}\big)^2\overset{k\le4/3}>-50\big(\frac{r}{R}\big)-194\big(\frac{r}{R}\big)^2
    \end{aligned}
\end{equation}
We lastly look at $\mathcal{C}^{\mathcal{E}}_L(\nu,k)$. Note that its parameter ranges are $k\in[1,4/3]$, $\nu\in(1/3,1/2]$. Therefore,
\begin{equation}\label{eq:2ndproducttermCLEnuk}
    -4\overbracket{\le}^{k\ge1}8k-12\overbracket{<}^{\nu>1/3}24k+\frac{\overbracket{-16k-12}^{<0}}{3\nu}\overbracket{\le}^{\nu\le0.5}24k+\frac{-16k-12}{3\cdot0.5}=\frac{40k-24}{3}\overbracket{\le}^{k\le4/3}\frac{88}{9}
\end{equation}

On the other hand, by the same reason for \eqref{eq:hat1Ug1Ubounds}, we have
\begin{align*}
    \tilde{h}_1(L)\in&\big[-8.53r/R,\quad 16.65r/R\big]\\
    -2d_1\tilde{g}_1(L)\in&\big[-\frac{8.5\nu}{0.9975}\frac{r}{R},\quad\frac{8.5\nu}{0.9975}\frac{r}{R}\big]
\end{align*}Hence 
\begin{equation}\label{eq:1stproducttermCLEnuk}
\tilde{h}_1(L)-4d_1\tilde{g}_1(L)\in\Big[-\big(8.53+\frac{17\nu}{0.9975}\big)\frac{r}{R},\quad \big(16.65+\frac{17\nu}{0.9975}\big)\frac{r}{R}\Big]
\end{equation} Therefore\COMMENT{Does the next equation need to be numbered? I am pretty sure that it does not.}
\begin{equation}\label{eq:CELnukbound}
\begin{aligned}
    \mathcal{C}^{\mathcal{E}}_{L}(\nu,k)=&\big(24k+\frac{-16k-12}{3\nu}\big)\big(\tilde{h}_1(L)-4d_1\tilde{g}_1(L)\big)-\frac{4k}{\nu}\big(\tilde{h}_1(L)-4d_1\tilde{g}_1(L)\big)^2\\
    \overbracket{\ge}^{\mathllap{\text{\eqref{eq:1stproducttermCLEnuk},\eqref{eq:2ndproducttermCLEnuk}}}}&-\big(8.53+\frac{17\nu}{0.9975}\big)\frac{r}{R}\cdot\frac{88}{9}-\frac{4k}{\nu}\big(16.65+\frac{17\nu}{0.9975}\big)^2\big(\frac{r}{R}\big)^2\\
    \overbracket{\ge}^{\mathllap{\nu\le0.5}}&-\big(8.53+\frac{17\cdot0.5}{0.9975}\big)\frac{r}{R}\cdot\frac{88}{9}-4\overbracket{k}^{\mathllap{k\le4/3}}\underbracket{\big(16.65+\frac{17\nu}{0.9975}\big)\big(\frac{16.65}{\nu}+\frac{17}{0.9975}\big)}_{<1497\text{ with }\nu\in(1/3,1/2]}\big(\frac{r}{R}\big)^2\\
    >&-92\big(\frac{r}{R}\big)-7984\big(\frac{r}{R}\big)^2.
\end{aligned}
\end{equation}\eqref{eq:CWUnukbound},\eqref{eq:CELnukbound}, \eqref{eq:CWLnukbound} are the claim of \eqref{item05Qlowerboundn1eq1}\qedhere
\end{proof}
Finally we get the following conclusion.
\begin{proposition}\label{proposition:casecexpansionn1eq1}
        In the context of \cref{TEcase_c}, that is, case (c) of \eqref{eqMainCases}, if $n_1=1$, then $\mathcal{II}_1(\mathcal{B}^+_0,1)$, $\mathcal{II}_2(\mathcal{B}^+_0,1)$ in \cref{def:AformularForLowerBoundOfExpansionForn1largerthan0} and in \eqref{itemAformularForLowerBoundOfExpansionForn1largerthan04} of \cref{lemma:AformularForLowerBoundOfExpansionForn1largerthan0} satisfy \[E(1,\Bcal^+_0)=\frac{\|dx_3\|_\p}{\|dx_0\|_\p}=3\big|\mathcal{II}_1(\mathcal{B}^+_0,1)+\mathcal{II}_2(\mathcal{B}^+_0,1)\big|>1-1743\big(\frac{r}{R}\big)-15450\big(\frac{r}{R}\big)^2>0.9.\]
\end{proposition}
\begin{proof}
With $n_1=1$, $\frac{\|dx_3\|_\p}{\|dx_0\|_\p}\overbracket{=}^{\mathrlap{\text{\cref{lemma:AformularForLowerBoundOfExpansionForn1largerthan0}\eqref{itemAformularForLowerBoundOfExpansionForn1largerthan04}}}}3\big|\mathcal{II}_1(\mathcal{B}^+_0,1)+\mathcal{II}_2(\mathcal{B}^+_0,1)\big|\overbracket{=}^{\mathllap{\text{\cref{def:II1kAndII2k}}}}3\big|\mathcal{II}_1(k)+\mathcal{II}_2(k)\big|\overbracket{=}^{\mathllap{\text{\cref{remark:II1II2Quadratic},\eqref{eq:II1II2QuadraticAndC1C2}}}}3\big|\mathtt{Q}\big(\frac{\tau_1}{d_0}\big)+\mathcal{C}_1(k)+\mathcal{C}_2(k)\big|$ with $\Bcal^+_0=\frac{-k}{d_0}$ for some $k\in[1,\frac{4}{3}]$.

With $\nu=\frac{d_1}{r\sin\phistar}$, by \cref{proposition:QUpperLowerRegionsLowerBoundsChatCtilde}\eqref{item03Qlowerboundn1eq1},\eqref{item04Qlowerboundn1eq1} and \eqref{item05Qlowerboundn1eq1}\begin{align*}
    \mathtt{Q}\big(\frac{\tau_1}{d_0}\big)&\overbracket{>}^{\mathclap{\text{\eqref{item03Qlowerboundn1eq1},\eqref{item04Qlowerboundn1eq1}}}}\min{\Big\{\frac{11}{9},\frac{1}{3}+\mathcal{C}^{\mathcal{W}}_{U}(\nu,k),\frac{11}{9}+\mathcal{C}^{\mathcal{E}}_{L}(\nu,k),\frac{1}{3}+\mathcal{C}^{\mathcal{W}}_{L}(\nu,k)\Big\}}\\
    &\overbracket{>}^{\mathclap{\text{\eqref{item05Qlowerboundn1eq1}}}}\min{\Big\{\frac{11}{9},\frac{1}{3}-192\big(\frac{r}{R}\big)-4667\big(\frac{r}{R}\big)^2,\frac{11}{9}-92\big(\frac{r}{R}\big)-7984\big(\frac{r}{R}\big)^2,\frac{1}{3}-50\big(\frac{r}{R}\big)-194\big(\frac{r}{R}\big)^2\Big\}}\\
    &\overbracket{=}^{\mathllap{\text{\eqref{eqJZHypCond}:}R>1700r}}\frac{1}{3}-192\big(\frac{r}{R}\big)-4667\big(\frac{r}{R}\big)^2
\end{align*}
Hence \begin{align*}
    3\Big(\mathtt{Q}\big(\frac{\tau_1}{d_0}\big)+\overbracket{\mathcal{C}_1(k)+\mathcal{C}_2(k)}^{\mathclap{\text{\cref{proposition:C1C2lowerbound}:}>-389(\frac{r}{R})-483(\frac{r}{R})^2}}\Big)>&3\Big(\frac{1}{3}-192\big(\frac{r}{R}\big)-4667\big(\frac{r}{R}\big)^2-389\big(\frac{r}{R}\big)-483\big(\frac{r}{R}\big)^2\Big)=3\Big(\frac{1}{3}-581\big(\frac{r}{R}\big)-5150\big(\frac{r}{R}\big)^2\Big)\\
    =&1-1743\big(\frac{r}{R}\big)-15450\big(\frac{r}{R}\big)^2\overbracket{>}^{\mathrlap{\eqref{eqJZHypCond}:R>30000r}}0.9,\\
    \frac{\|dx_3\|_\p}{\|dx_0\|_\p}>&1-1743\big(\frac{r}{R}\big)-15450\big(\frac{r}{R}\big)^2>0.9.\qedhere
\end{align*}
\end{proof}
\begin{theorem}[In the context of \cref{lemma:AformularForLowerBoundOfExpansionForn1largerthan0}, thus of \cref{TEcase_c}]\label{theorem:B3minusdx3dphi3ratiobound}
     The wave front curvature $\Bcal^-_3$ at $x_3$ in \cref{lemma:AformularForLowerBoundOfExpansionForn1largerthan0}\eqref{itemAformularForLowerBoundOfExpansionForn1largerthan05} and \eqref{eq:B3minusEquation} is uniformly bounded. i.e., $|\Bcal^-_3|<const$ in case (c) of \eqref{eqMainCases} for some uniform constant.
     
     Furthermore, since \eqref{eq:MirrorEquationExpansionFormula}: $\frac{d\theta_3}{d\phi_3}=-r(\Bcal^-_3\cos\varphi_3-\frac{1}{r})$, thus in \cref{lemma:AformularForLowerBoundOfExpansionForn1largerthan0}\eqref{itemAformularForLowerBoundOfExpansionForn1largerthan03} $(d\phi_3,d\theta_3)=dx_3=D\Fc^{n_1+4}(dx_0)$ has $\frac{d\theta_3}{d\phi_3}$ also uniformly bounded in case (c). Thus, by \cref{caseC}, there is a uniform upper bound $\lambda_2(r,R,\phistar)>0$ such that $0<\frac{d\theta_3}{d\phi_2}<\lambda_2(r,R,\phistar)$, where constant $\lambda_2(r,R,\phistar)$ is determined by lemon billiard configuration.\COMMENT{This reads like a proof rather than a Theorem. That makes is harder to follow.}
\end{theorem}
\begin{proof}
    It suffices to prove each term in \eqref{eq:B3minusEquation} is uniformly bounded. The first term: 
    \begin{equation}\label{eq:Bmins3firstterm}
    \begin{aligned}\bigm|\frac{-2}{d_2(1+2d_2\Bcal^+_2)}\bigm|\overbracket{=}^{\mathclap{\eqref{eq:Expansionx3x0}}}&\frac{2}{d_2}\frac{|1+\tau_0\Bcal^+_0||1+\tau_{1,\text{\upshape im}}\Bcal^+_{1,\text{\upshape im}}||1+\tau_1\Bcal^+_{n_1,R}|}{E(n_1,\Bcal^+_0)}\\
    \overbracket{=}^{\mathllap{\text{\cref{lemma:AformularForLowerBoundOfExpansionForn1largerthan0}\eqref{itemAformularForLowerBoundOfExpansionForn1largerthan04}}}}&\frac{|1+\tau_0\Bcal^+_0|}{E(n_1,\Bcal^+_0)}\bigm|1+(-2n_1d_1)(\Bcal^{-}_1-\frac{1}{d_1})\bigm||1+\tau_1\Bcal^+_{n_1,R}|\\
    \overbracket{=}^{\mathllap{\text{Proof of \cref{lemma:AformularForLowerBoundOfExpansionForn1largerthan0}\eqref{itemAformularForLowerBoundOfExpansionForn1largerthan04}\ref{item04ContinuedfractionResults}}}}&\frac{|1+\tau_0\Bcal^+_0|}{E(n_1,\Bcal^+_0)}\bigm|1+(-2n_1d_1)(\Bcal^{-}_1-\frac{1}{d_1})\bigm|\bigm|1+\tau_1[-\frac{1}{d_1}+\frac{\Bcal^{-}_1-\frac{1}{d_1}}{1-2n_1d_1(\Bcal^{-}_1-\frac{1}{d_1})}]\bigm|\\
    =&\frac{|1+\tau_0\Bcal^+_0|}{E(n_1,\Bcal^+_0)}\bigm|1-\frac{\tau_1}{d_1}+2n_1(\tau_1-d_1)(\Bcal^-_1-\frac{1}{d_1})\bigm|\\
    =&\frac{|1+\tau_0\Bcal^+_0|}{E(n_1,\Bcal^+_0)}\bigm|1+2n_1-(2n_1+1)\frac{\tau_1}{d_1}+2n_1(\tau_1-d_1)\Bcal^-_1\bigm|\\
    \le&\frac{|1+\tau_0\Bcal^+_0|}{E(n_1,\Bcal^+_0)}\bigm|1+2n_1-(2n_1+1)\frac{\tau_1}{d_1}\bigm|+\frac{|\Bcal^-_1||1+\tau_0\Bcal^+_0|}{E(n_1,\Bcal^+_0)}\bigm|2n_1(\tau_1-d_1)\bigm|\\
    \overbracket{<}^{\mathllap{\substack{\text{\cref{remark:deffpclf}: }0<\tau_1<2d_1\text{ and }\\\text{\cite[equation (3.31)]{cb}: }\\\Bcal^-_1(1+\tau_0\Bcal^+_0)=\Bcal_0^+}}}&\frac{|1+\tau_0\Bcal^+_0|(2n_1+1)}{E(n_1,\Bcal^+_0)}+\frac{2n_1d_1|\Bcal^+_0|}{E(n_1,\Bcal^+_0)}
    \end{aligned}
    \end{equation}
Now we notice the following.
\begin{itemize}
    \item By \cref{lemma:AformularForLowerBoundOfExpansionForn1largerthan0}\eqref{itemAformularForLowerBoundOfExpansionForn1largerthan01}: $\Bcal^+_0\in[\frac{-4}{3d_0},\frac{-1}{d_0}]$ and \cref{proposition:casecparametersestimate}\ref{itemcasecparametersestimate1}: $|d_0-r\sin\phistar|<\frac{17R}{4R}r\sin\phistar$, thus $\Bcal^+_0$ is uniformly bounded.
    \item $1+\tau_0\Bcal^+_0\in(1-\frac{4\tau_0}{3d_0},1-\frac{\tau_0}{d_0})\subset(-\frac{5}{3},1)$ since $0<\tau_0<2d_0$. And $2n_1d_1<2R\Phistar$ the arc length of $\Gamma_R$.
    \item By \cref{proposition:casecexpansionn1ge4,proposition:casecexpansionn1eq3,proposition:casecexpansionn1eq2,proposition:casecexpansionn1eq1}, for every $n_1>0$, $E(n_1,\Bcal^+_0)>0.9$. And especially by \cref{proposition:casecexpansionn1ge4}, $\frac{2n_1+1}{E(n_1,\Bcal^+_0)}\rightarrow0$ as $n_1\rightarrow+\infty$.
\end{itemize}
Therefore, we see that the right terms in the last line of \eqref{eq:Bmins3firstterm} are uniformly bounded. $\bigm|\frac{-2}{d_2(1+2d_2\Bcal^+_2)}\bigm|$ is uniformly bounded.

The second term:
\begin{equation}\label{eq:Bmins3secondterm}
    \begin{aligned}
        \bigm|-\frac{1}{d_1(1+2d_2\Bcal_2^+)(1+\tau_1\Bcal^{+}_{n_1,R})}\bigm|\overbracket{=}^{\mathclap{\eqref{eq:Expansionx3x0}}}&\frac{1}{d_1}\frac{|1+\tau_0\Bcal^+_0||1+\tau_{1,\text{\upshape im}}\Bcal^+_{1,\text{\upshape im}}|}{E(n_1,\Bcal^+_0)}\\
        \overbracket{=}^{\mathllap{\text{\cref{lemma:AformularForLowerBoundOfExpansionForn1largerthan0}\eqref{itemAformularForLowerBoundOfExpansionForn1largerthan04}}}}&\frac{1}{d_1}\frac{|1+\tau_0\Bcal_0^+||1-2n_1d_1(\Bcal^-_1-\frac{1}{d_1})|}{E(n_1,\Bcal^+_0)}\\
        =&\frac{1}{d_1}\frac{|(1+2n_1)(1+\tau_0\Bcal^+_0)-2n_1d_1\Bcal^-_1(1+\tau_0\Bcal^+_0)|}{E(n_1,\Bcal_0^+)}\\
        \le&\frac{|1+\tau_0\Bcal^+_0|(2n_1+1)}{d_1E(n_1,\Bcal_0^+)}+\frac{2n_1|\Bcal^+_0|}{E(n_1,\Bcal^+_0)}\\
        \overbracket{<}^{\mathllap{\text{\cref{proposition:casecparametersestimate}\ref{itemcasecparametersestimate3}: } d_1>\frac{1}{n_1+2}r\sin\phistar}}&\frac{|1+\tau_0\Bcal^+_0|(2n_1+1)(n_1+2)}{r\sin\phistar E(n_1,\Bcal^+_0)}+\frac{2n_1|\Bcal^+_0|}{E(n_1,\Bcal^+_0)}
    \end{aligned}
\end{equation}
The same as in the first term estimate, both $|\Bcal^+_0|$, $|1+\tau_0\Bcal^+_0|$ are uniformly bounded. \cref{proposition:casecexpansionn1ge4,proposition:casecexpansionn1eq3,proposition:casecexpansionn1eq2,proposition:casecexpansionn1eq1} gives for each $n_1>0$, $E(n_1,\Bcal^+_0)>0.9$. By \cref{proposition:casecexpansionn1ge4} $\limsup_{n_1\rightarrow\infty}{\frac{(2n_1+1)(n_1+2)}{E(n_1,\Bcal^+_0)}}< C_{c}$ for some constant $C_{c}>0$ and $\lim_{n_1\rightarrow\infty}\frac{n_1}{E(n_1,\Bcal^+_0)}=0$. Therefore, the two terms in the last line of \eqref{eq:Bmins3secondterm} are uniformly bounded. $-\frac{1}{d_1(1+2d_2\Bcal_2^+)(1+\tau_1\Bcal^{+}_{n_1,R})}$ is uniformly bounded.

The third term: 
\begin{equation}\label{eq:B3minusthirdterm}
    \begin{aligned}
            \bigm|-\frac{1}{d_1(1+2d_2\Bcal^+_2)(1+\tau_1\Bcal^+_{n_1,R})(1+\tau_{1,\text{\upshape im}}\Bcal^+_{1,\text{\upshape im}})}\bigm|\overbracket{=}^{\mathclap{\eqref{eq:Expansionx3x0}}}&\frac{|1+\tau_0\Bcal_0^+|}{d_1E(n_1,\Bcal_0^+)}\\
            \overbracket{<}^{\mathllap{\text{\cref{proposition:casecparametersestimate}\ref{itemcasecparametersestimate3}: } d_1>\frac{1}{n_1+2}r\sin\phistar}}&\frac{|1+\tau_0\Bcal^+_0|(n_1+2)}{r\sin\phistar E(n_1,\Bcal_0^+)}.
    \end{aligned}
\end{equation}
The same as in the first term estimate, $|1+\tau_0\Bcal^+_0|$ is uniformly bounded. \cref{proposition:casecexpansionn1ge4,proposition:casecexpansionn1eq3,proposition:casecexpansionn1eq2,proposition:casecexpansionn1eq1} give for each $n_1>0$, $E(n_1,\Bcal^+_0)>0.9$. By \cref{proposition:casecexpansionn1ge4}  $\lim_{n_1\rightarrow\infty}\frac{n_1+2}{E(n_1,\Bcal^+_0)}=0$. Therefore, the last line term of \eqref{eq:B3minusthirdterm} is uniformly bounded. $|-\frac{1}{d_1(1+2d_2\Bcal^+_2)(1+\tau_1\Bcal^+_{n_1,R})(1+\tau_{1,\text{\upshape im}}\Bcal^+_{1,\text{\upshape im}})}|$ is uniformly bounded.

The fourth term $\bigm|\frac{\Bcal^+_0}{(1+2d_2\Bcal^+_2)(1+\tau_1\Bcal^+_{n_1,R})(1+\tau_{1,\text{\upshape im}}\Bcal^+_{1,\text{\upshape im}})(1+\tau_0\Bcal_0^+)}\bigm|=\frac{|\Bcal^+_0|}{E(n_1,\Bcal^+_0)}.$
Same proof as in the first item gives $|\Bcal^+_0|$ is uniformly bounded. \cref{proposition:casecexpansionn1ge4,proposition:casecexpansionn1eq3,proposition:casecexpansionn1eq2,proposition:casecexpansionn1eq1} imply that for each $n_1>0$, $E(n_1,\Bcal^+_0)>0.9$. Therefore $|\frac{\Bcal^+_0}{(1+2d_2\Bcal^+_2)(1+\tau_1\Bcal^+_{n_1,R})(1+\tau_{1,\text{\upshape im}}\Bcal^+_{1,\text{\upshape im}})(1+\tau_0\Bcal_0^+)}|$ is uniformly bounded.\COMMENT{This section started on page \pageref{sec:ProofTECase_c}---maybe it deserves a brief summary at the end.}
\end{proof}

\section{Approximation of 1-petal billiard by lemon billiard}\label{sec:SACTWFELBWB}

In this section, we describe which values of the parameter $\phistar$ must be excluded, and how to determine $R_{\text{\upshape HF}}(r,\phistar)$ in \eqref{eqJZHypCond}. We prove \cref{RMASUH,thm:MroutreturnUniformExpansion,thm:UniformExpansionOnMhat1}\COMMENT{Edit} by comparing the Lemon billiard expansion and trajectory with the 1-petal billiard.

\subsection[Exceptional \texorpdfstring{$\phistar$}{configuration angle phi*} in \texorpdfstring{main theorems: \cref{RMASUH,thm:UniformExpansionOnMhat1,thm:MroutreturnUniformExpansion}}{main theorems}]{Exceptional \texorpdfstring{$\phistar$}{configuration angle phi*} in \texorpdfstring{main theorems \cref{RMASUH,thm:UniformExpansionOnMhat1,thm:MroutreturnUniformExpansion}}{main theorems}}\hfill

\begin{figure}[ht]
\strut\hfill
    \begin{tikzpicture}[xscale=0.65,yscale=.65]
    \tkzDefPoint(0,0){Or};
    \tkzDefPoint(0,-4.8){Y};
    \pgfmathsetmacro{\rradius}{4.5};
    \pgfmathsetmacro{\bdist}{6.3};
    \clip(-1.0*\rradius-0.2,-1.4*\rradius) rectangle (\rradius*1.0+0.2, \rradius*1.5);
    \tkzDefPoint(0,\bdist){OR};
    \pgfmathsetmacro{\phistardeg}{18.1};
    \pgfmathsetmacro{\XValueArc}{\rradius*sin(\phistardeg)};
    \pgfmathsetmacro{\YValueArc}{\rradius*cos(\phistardeg)};    
    \tkzDefPoint(\XValueArc,-1.0*\YValueArc){B};    \tkzDefPoint(-1.0*\XValueArc,-1.0*\YValueArc){A};
    \pgfmathsetmacro{\Rradius}{veclen(\XValueArc,\YValueArc+\bdist)};
    \pgfmathsetmacro{\XValueArcBOne}{\rradius*sin(3*\phistardeg)};
    \pgfmathsetmacro{\YValueArcBOne}{\rradius*cos(3*\phistardeg)};
    \pgfmathsetmacro{\XValueArcAOne}{-\rradius*sin(3*\phistardeg)};
    \pgfmathsetmacro{\YValueArcAOne}{\rradius*cos(3*\phistardeg)};
    
    \tkzDefPoint(\XValueArcBOne,-1.0*\YValueArcBOne){B1};
    \tkzDefPoint(-\XValueArcBOne,-1.0*\YValueArcBOne){A1};
    \tkzDrawArc[name path=Cr,very thick](Or,B)(A);
    \draw[very thick](A)--(B);
    \draw[very thin,blue,arrows={-Triangle[length=0.2cm]}](B) -- (B1);
    \draw[very thin,blue,arrows={-Triangle[length=0.2cm]}](B1) -- (B);
    \draw[very thin,blue,arrows={-Triangle[length=0.2cm]}](A) -- (A1);
    \draw[very thin,blue,arrows={-Triangle[length=0.2cm]}](A1) -- (A);
    \tkzLabelPoint[below,rotate=0](A){$A$};
    \tkzLabelPoint[left,rotate=0](A1){$A_1$};
    \tkzLabelPoint[below,rotate=0](B){$B$};
    \tkzLabelPoint[right,rotate=0](B1){$B_1$};
    \tkzLabelPoint[above,rotate=0](Or){$O_r$};
    \tkzDrawPoints(A,B,Or);
    \draw[dashed,ultra thin](Or)--+(0,-\rradius)coordinate(negYaxis);
    \draw[dashed,ultra thin](Or)--(A);
    \draw[dashed,ultra thin](Or)--(B);
    \tkzMarkAngle[arc=lll,ultra thin,size=1.1,mark=|](A,Or,Y);
    \tkzMarkAngle[arc=lll,ultra thin,size=1.2,mark=|](Y,Or,B);
    \node[rotate=0] at ($(Or)+(-80:1.8)$)  {\small $\phistar$};
    \node[rotate=0] at ($(Or)+(-0.6*\rradius,0.93*\rradius)$)  {\Large $\Gamma_r$};
    \node[rotate=0] at ($(Or)+(0,-1.18*\rradius)$)  {\Large $\PET(r,\phistar)$};
    \end{tikzpicture}\hfill
    \begin{tikzpicture}
        [xscale=0.65,yscale=.65]
    \tkzDefPoint(0,0){Or};
    \tkzDefPoint(0,-4.8){Y};
    \pgfmathsetmacro{\rradius}{4.5};
    \pgfmathsetmacro{\bdist}{6.3};
    \clip(-1.0*\rradius-0.2,-1.4*\rradius) rectangle (\rradius*1.0+0.2, \rradius*1.5);
    \tkzDefPoint(0,\bdist){OR};
    \pgfmathsetmacro{\phistardeg}{18.1};
    \pgfmathsetmacro{\XValueArc}{\rradius*sin(\phistardeg)};
    \pgfmathsetmacro{\YValueArc}{\rradius*cos(\phistardeg)};    
    \tkzDefPoint(\XValueArc,-1.0*\YValueArc){B};    \tkzDefPoint(-1.0*\XValueArc,-1.0*\YValueArc){A};
    \pgfmathsetmacro{\Rradius}{veclen(\XValueArc,\YValueArc+\bdist)};
    \pgfmathsetmacro{\XValueArcBOne}{\rradius*sin(3*\phistardeg)};
    \pgfmathsetmacro{\YValueArcBOne}{\rradius*cos(3*\phistardeg)};
    \pgfmathsetmacro{\XValueArcAOne}{-\rradius*sin(3*\phistardeg)};
    \pgfmathsetmacro{\YValueArcAOne}{\rradius*cos(3*\phistardeg)};
    
    \tkzDefPoint(\XValueArcBOne,-1.0*\YValueArcBOne){B1};
    \tkzDefPoint(-\XValueArcBOne,-1.0*\YValueArcBOne){A1};
    \tkzDrawArc[name path=Cr,very thick](Or,B)(A);
    \tkzDrawArc[name path=CR,very thick](OR,A)(B);
    \draw[very thin,blue,arrows={-Triangle[length=0.2cm]}](B) -- (B1);
    \draw[very thin,blue,arrows={-Triangle[length=0.2cm]}](B1) -- (B);
    \draw[very thin,blue,arrows={-Triangle[length=0.2cm]}](A) -- (A1);
    \draw[very thin,blue,arrows={-Triangle[length=0.2cm]}](A1) -- (A);
    \tkzLabelPoint[below,rotate=0](A){$A$};
    \tkzLabelPoint[left,rotate=0](A1){$A_1$};
    \tkzLabelPoint[below,rotate=0](B){$B$};
    \tkzLabelPoint[right,rotate=0](B1){$B_1$};
    \tkzLabelPoint[above,rotate=0](Or){$O_r$};
    \tkzLabelPoint[left,rotate=0](OR){$O_R$};
    \tkzDrawPoints(A,B,Or,OR);
    \draw[dashed,ultra thin](Or)--+(0,-\rradius)coordinate(negYaxis);
    \draw[dashed,ultra thin](Or)--(A);
    \draw[dashed,ultra thin](Or)--(B);
    \draw[dashed,ultra thin](OR)--(A);
    \draw[dashed,ultra thin](OR)--(B);
    \tkzMarkAngle[arc=lll,ultra thin,size=1.1,mark=|](A,Or,Y);
    \tkzMarkAngle[arc=lll,ultra thin,size=1.2,mark=|](Y,Or,B);
    \node[rotate=0] at ($(Or)+(-80:1.8)$)  {\small $\phistar$};
    \node[rotate=0] at ($(Or)+(0.6*\rradius,0.93*\rradius)$)  {\Large $\Gamma_r$};
    \node[rotate=0] at ($(Or)+(0,-1.18*\rradius)$)  {\Large $\LEM(r,R,\phistar)$};
    \end{tikzpicture}
    \hfill\strut
    \caption{A 1-petal billiard table $\PET(r,\phistar)$ and a lemon billiard table $\LEM(r,R,\phistar)$ with same parameters $r$, $\phistar$.}\label{fig:1-petalExtensionFromCornerAndLemon}
\end{figure}

In this subsection, we note that the exceptional $\phistar\in\big(0,\tan^{-1}{(1/3)}\big)$ in \cref{thm:MroutreturnUniformExpansion,RMASUH} are $\pi/n$ for $n\in\mathbb{N}$. This is related to two comparison billiards for the lemon billiard as follows. In the lemon billiard, the boundary arc that connects $A$ and $B$ in \cref{fig:1-petalExtensionFromCornerAndLemon} has curvature $1/R$. If the curvature is zero instead, then we obtain the 1-petal billiard on the left of \cref{fig:1-petalExtensionFromCornerAndLemon}, and if the curvature is $1/r$, then the billiard is perfectly circular. The exceptional $\phistar$s are exactly those for which $A,B$ are adjacent points of a periodic orbit in this circular billiard. 

Subsequent subsections make quantitative comparisons between the two billiards in \cref{fig:1-petalExtensionFromCornerAndLemon}, and we now explicitly adopt our standing notations for work on the 1-petal billiard.
\begin{notation}[Lemon billiard notations used for 1-petal billiard]\label{def:LemonAndPetTogether}\hfill
    \begin{enumerate}
        \item\label{item01LemonAndPetTogether} We denote by $\PET(r,\phistar)$ a Wojtowkski 1-petal billiard table with corners $A$ and $B$, a flat boundary $\overline{AB}$ and a boundary $\Gamma_r$ that is an arc of circle $C_r$ centered on $O_r$, and $\measuredangle{AO_rB}=2\phistar<2\tan^{-1}\big(1/3\big)$ shown in \cref{fig:WojtkowskiBilliardTable,fig:1-petalExtensionFromCornerAndLemon}. We denote by $\overline{F}_{\phistar}$  the 1-petal billiard map of $\PET(r,\phistar)$ as defined by \cite[Equation (2.18)]{cb}. (This is the return map of the 1-petal billiard flow to the boundary \cite[Sections 2.9--2.12]{cb}.)\COMMENT{Are these references necessary? If the beginning of this paper has a clear definition of "billiard map" at the start, then one can refer to it instead of elsewhere. Even that may not be needed.}
        \item\label{item02LemonAndPetTogether} The Wojtowkski 1-petal billiard collisions on $\Gamma_r$ are represented by $(\phi,\theta)\in M_r=[\phistar,2\pi-\phistar]\times[0,\pi]$, which is in \cref{fig:Mr} the same as the lemon billiard phase space $M_r$ in \cref{def:BasicNotations0}, with $\phi$ being the arc angle parameter of $C_r$ for the collision position and $\theta$ being the collision angle with the tangent direction of $\Gamma_r$. The collisions on $\overline{AB}$ are represented by $(s,\theta_1)\in M_{\text{\upshape{f}}}\dfn[-r\sin\phistar,r\sin\phistar]\times[0,\pi]$ with $s$ being the length parameter on $\overline{AB}$ and $\theta_1$ being the collision angle with tangential direction. The billiard map $\overline{F}_{\phistar}$ is a map from phase space $M_r\sqcup M_{\text{\upshape{f}}}\nfd\overline{M}$ to itself.
        \item\label{item03LemonAndPetTogether} For $x\in M_{\text{\upshape{f}}}$, we define $p(x)$ as the collision position in $\PET(r,\phistar)$. 
        \item\label{item04LemonAndPetTogether} The $\Mrout$, $\Mrin$ in \cref{def:BasicNotations} defined for a $\LEM(r,R,\phistar)$ billiard map $\mathcal{F}$ have the similar meaning for the $\PET(r,\phistar)$ billiard map $\overline{F}_{\phistar}$. That is, they
 are the collisions that occur on the boundary $\Gamma_r$ with the next(previous) collision on another boundary. i.e. \begin{align*}
            \Int M_r\cap\mathcal{F}^{-1}(\Int M_R)\overbracket{=}^{\mathllap{\text{\cref{def:BasicNotations}}}}\Mrout&=\Int M_r\cap\overline{F}_{\phistar}^{-1}(\Int M_{\text{\upshape{f}}}),\\
            \Int M_r\cap\mathcal{F}(\Int M_R)\overbracket{=}^{\mathllap{\text{\cref{def:BasicNotations}}}}\Mrin&=\Int M_r\cap\overline{F}_{\phistar}(\Int M_{\text{\upshape{f}}}), 
        \end{align*}
        
        For $x\in \Int M_r$, $p(\Fcal(x))\in \Int\Gamma_r$ is equivalent to $p(\overline{F}_{\phistar}(x))\in\Int\Gamma_r$. Therefore for $x\in\Mrin$, \[
        \max{\big\{m\ge0\bigm|p(\overline{F}_{\phistar}^i(x))\in\operatorname{int}\Gamma_r\text{ for }0\le i\le m\big\}}=\max{\big\{m\ge0\bigm|p(\mathcal{F}^i(x))\in\operatorname{int}\Gamma_r\text{ for }0\le i\le m\big\}}\overbracket{=}^{\mathllap{\text{\cref{def:SectionSetsDefs}}}}m(x)
        \]
        \item\label{item05LemonAndPetTogether} For $\overline{x}_0=(\overline{\phi}_0,\overline{\theta}_0)\in\Mrout$ with $\overline{x}_1=\overline{F}_{\phi_*}(\overline{x}_0)\in M_{\text{\upshape f}}$ and $(\overline{\phi}_2,\overline{\theta}_2)=\overline{x}_2=\overline{F}^2_{\phi_*}(\overline{x}_0)\in\Mrin$ define $\overline{d}_0\dfn r\sin{\overline{\theta}_0}$, $\overline{d}_2\dfn r\sin{\overline{\theta}_2}$ and let $\overline{\tau}_0\dfn|p(\overline{x}_0)p(\overline{x}_1)|$, $\overline{\tau}_1\dfn|p(\overline{x}_1)p(\overline{x}_2)|$ be the Euclidean distance between the two consecutive collisions on $\mathbb{L}(r,\phi_*)$ billiard table in \cref{fig:1-petalExtensionFromCornerAndLemon}.
        \item\label{item06LemonAndPetTogether} With same $\Min_{r,n}$, $\Mout_{r,n}$ and global section $\hat{M}= \big(\bigsqcup_{n\ge1}\mathcal{F}^{n-1}\Min_{r,n}\big)\sqcup\Min_{r,0}=\mathcal{F}^{-1}\big(\Mrout\smallsetminus{\Mrin}\big)\sqcup\big(\Mrout\cap\Mrin\big)$ in \cref{def:SectionSetsDefs} and \eqref{def:Mhat}, we define $\hat{\overline{F}}_{\phistar}$ the return map on $\hat{M}$ in $\PET(r,\phistar)$ table under $\overline{F}_{\phistar}$ dynamics. And we define $\tilde{\overline{F}}_{\phistar}$ the return map on $\Mrout$ in $\PET(r,\phistar)$ table under $\overline{F}_{\phistar}$ dynamics.
        \item\label{item07LemonAndPetTogether} In a table $\PET(r,\phistar)$, we consider $\overrightarrow{BA}$ as a collision on the corner $A$ with direction $\overrightarrow{BA}$. It is a singular/boundary point of $M_r$ that is $(2\pi-\phistar,\pi-\phistar)=Ix_*$ in \cref{fig:Mr}.
        
        Resuming $\overrightarrow{BA}$ to have a reflection on $\Gamma_r$ at the corner $A$, in \cref{fig:1-petalExtensionFromCornerAndLemon} we define $\overrightarrow{AA_1}$ as its next collision direction at $A_1\in\Gamma_r$. Then in phase space \cref{fig:Mr}, $(2\pi-3\phistar,\pi-\phistar)\in M_r$ is the extended definition of $\overrightarrow{BA}$'s next collision at $A_1$. 
        
        Symmetrically, we consider $\overrightarrow{AB}$ as a collision on the corner $B$ with direction $\overrightarrow{AB}$. It is a singular / boundary point of $M_r$, that is, $(\phistar,\phistar)=Iy_*$ in \cref{fig:Mr}.
        
        Resuming\COMMENT{??? Look this word up.} $\overrightarrow{AB}$ to have a reflection on $\Gamma_r$ at the corner $B$, in \cref{fig:1-petalExtensionFromCornerAndLemon} We define $\overrightarrow{BB_1}$ as its next collision direction at $B_1\in\Gamma_r$. Then $(3\phistar,\phistar)\in M_r$ in \cref{fig:Mr} is the extended definition of $\overrightarrow{AB}$'s next collision at $B_1$.
        
        \COMMENT{This is excessively complicated and virtually incomprehensible. Therefore there is nothing obvious in the next line. Is there a good reason for the complications? I suspect not.\newline\wentao{I moved this item to be become a new definition \cref{def:ExceptionPhistar}.}}
    \end{enumerate}
\end{notation}

\begin{definition}[exceptional and generalized exceptional $\phistar$]\label{def:ExceptionPhistar}
    For lemon billiard $\mathbb{L}(r,\phistar)$ with $0<\phistar<\tan^{-1}(1/3)$, we define \emph{$\phi_*$ as exceptional} if in the 1-petal billiard $\mathbb{P}(r,\phistar)$ the image of $(3\phistar,\phistar)$ under some positive iteration of $\overline{F}_{\phistar}$ is $(\phistar,\pi-\phistar)$ or $(2\pi-\phistar,\phistar)$. i.e. $\overline{F}^n_{\phistar}(3\phistar,\phistar)=(\phistar,\pi-\phistar)$ or $(2\pi-\phistar,\phistar)$ for some $n\ge1$.

    The more \emph{generalized exceptional $\phi_*$} are the $0<\phistar<\tan^{-1}(1/3)$ configurations for $\mathbb{L}(r,\phi_*)$ such that the image of $(3\phistar,\phistar)$ under some positive iteration of $\overline{F}_{\phistar}$ hits the corner. i.e., the $\phi$-coordinate  of $\overline{F}^n_{\phistar}(3\phistar,\phistar)$ is either $\phistar$ or $2\pi-\phistar$ for some $n\ge1$.
    
    Especially we call this $(3\phi_*,\phi_*)$'s images under $\overline{F}^n_{\phi_*}|_{n\ge0}$ the \emph{positive trajectory of $(3\phi_*,\phi_*)$}.

\end{definition}
\begin{remark}\label{remark:exceptionalphistar} In \cref{fig:1-petalExtensionFromCornerAndLemon}, we can envision the following fact and we will make some claims in future work.
\begin{itemize}
    \item The exceptional $\phi_*$ are those tables $\mathbb{P}(r,\phi_*)$ with a trajectory $\overrightarrow{BB_1}$ (or equivalently $\overrightarrow{AA_1}$) in the future hitting the corner $A$ in the direction of $\overrightarrow{A_1A}$ or hitting the corner $B$ in the direction of $\overrightarrow{B_1B}$.
    \item The obvious exceptional $\phistar$s are $\pi/n$, $n\in\mathbb{N}$.\COMMENT{Since the last item above is incomprehensible, this is not obvious.}
    \item By some differential map checking and applying Sard theorem, the generalized exceptional $\phi_*$ in $(0,\tan^{-1}(1/3))$ has measure $0$. Note that exceptional $\phistar$ is a subset of generalized exceptional $\phistar$, therefore it is also measure $0$. We further conjecture the exceptional $\phi_*$ is countable. 
    \item If $\tan^{-1}\big(1/3\big)>\phistar\ne\frac{\pi}{n} \forall n\in\mathbb{N}$, then the lower bound for the expansion of the cone vectors along the orbit of $\overrightarrow{BB_1}$ in $\PET(r,\phistar)$ exponentially increases with the return times (on the section). So, within finite returns the expansion of the cone vector in the $\Nin$ orbit of $\LEM(r,R,\phistar)$ will be greater than $22.5$ so that it overcomes the contraction of $0.05$ within finite many steps (of the return map).
    \end{itemize}
\end{remark}

\begin{lemma}
[Approximation for consecutive collisions on $\Gamma_r$]\label{lemma:ApproximationForConsecutiveCollisions}
    For a given nonsingular point $\overline{x}=(\overline{\phi},\overline{\theta})\in 
    \Mrin$, 
    \begin{itemize}
        \item $\exists\,\overline{\delta}_0(\overline{x})>0$ such that if $x=(\phi,\theta)\in\Mrin$ satisfies the following:
    \begin{equation}\label{eq1:ApproximationForConsecutiveCollisions}
        \begin{aligned}
            |\phi-\overline{\phi}|<\overline{\delta}_0\\
            |\theta-\overline{\theta}|<\overline{\delta}_0,
        \end{aligned}
    \end{equation} then the \cref{def:LemonAndPetTogether}\eqref{item04LemonAndPetTogether} defined $m(x)=m(\overline{x})$
    \item For $\forall\varepsilon>0$, $\exists\,\delta_0(\overline{x},\varepsilon)\dfn\min\big\{\overline{\delta}_0(\overline{x}),\,\frac{\varepsilon}{2},\,\frac{\varepsilon}{2m(\overline{x})}\big\}$ such that if $x=(\phi,\theta)\in\Mrin$ satisfies the following:
        \begin{equation}\label{eq2:ApproximationForConsecutiveCollisions}
        \begin{aligned}
            |\phi-\overline{\phi}|<\delta_0\\
            |\theta-\overline{\theta}|<\delta_0,
        \end{aligned}
    \end{equation} then the $\phi\theta$ coordinates of $\Fc^{m(x)}(\phi,\theta)=(\phi+2m(x)\theta,\quad\theta)\in\Mrout$ and $\overline{F}_{\phistar}^{m(\overline{x})}(\overline{\phi},\overline{\theta})=(\overline{\phi}+2m(\overline{x})\overline{\theta},\quad\overline{\theta})\in\Mrout$ are $\varepsilon$ close. i.e. $\bigm|\overline{\theta}-\theta\bigm|<\varepsilon$, $\bigm|\phi+2m(x)\theta-(\overline{\phi}+2m(\overline{x})\overline{\theta})\bigm|<\varepsilon$.
    \end{itemize}
    Under condition $0<\phistar<\tan^{-1}{(1/3)}$ and $\phistar\ne\frac{\pi}{n},\forall n\ge1$, for special point $x_*=(\phistar,\phistar)$, $\exists$ finite $m_*\dfn\inf\big\{ k>0\bigm|\phistar+2k\phistar\in(-\phistar,+\phistar)(\text{ mod }2\pi)\big\}$
    \begin{itemize}
        \item $\exists\,\overline{\delta}_0(x_*)>0$ such that if nonsingular $x=(\phi,\theta)\in\Mrin$ satisfies the following.
        \begin{equation}\label{eq3:ApproximationForConsecutiveCollisions}
            \begin{aligned}
                |\phi-\overline{\phi}|<\overline{\delta}_0(x_*)\\
                |\theta-\overline{\theta}|      <\overline{\delta}_0(x_*),
            \end{aligned}
        \end{equation} 
        then from the \cref{def:LemonAndPetTogether}\eqref{item04LemonAndPetTogether} $m(x)=m_*$.\COMMENT{Edit\wentao{Better organized sentence now. pls check}}
        \item For $\forall\varepsilon>0$, $\exists\,\delta_0(x_*,\varepsilon)\dfn\min\big\{\overline{\delta}_0(x_*),\,\frac{\varepsilon}{2},\,\frac{\varepsilon}{2m_*}\big\}$ such that if $x=(\phi,\theta)\in\Mrin$ satisfies the following:
        \begin{equation}\label{eq4:ApproximationForConsecutiveCollisions}
        \begin{aligned}
            |\phi-\overline{\phi}|<\delta_0\\
            |\theta-\overline{\theta}|<\delta_0,
        \end{aligned}
    \end{equation} then the $\phi\theta$ coordinates of $\Fc^{m(x)}(\phi,\theta)=(\phi+2m(x)\theta,\,\theta)\in\Mrout$ and $\overline{F}_{\phistar}^{m_*}(\phistar,\phistar)=(\phistar+2m_*\phistar,\,\phistar)\in\Mrout$ are $\varepsilon$ close. i.e. $\bigm|\overline{\theta}-\theta\bigm|<\varepsilon$, $\bigm|\phi+2m(x)\theta-(\phistar+2m_*\phistar)\bigm|<\varepsilon$.
    \end{itemize}
\end{lemma}
\begin{proof}
The first conclusion follows the fact that the nonsingular $\overline{x}$ is in the interior of $\Min_{r,m(\overline{x})}$ which is a finite number of disconnected quadrilateral components (see \cref{fig:MrN}). The second conclusion follows $m(x)=m(\overline{x})$.

The third conclusion follows the fact that $x_*$ is on the boundary of $\Min_{r,m_*}$, which is a finite number of disconnected quadrilateral components. The fourth conclusion follows $m(x)=m_*$.
\end{proof}

\subsection{$\mathbb{P}(r,\phi_*)$ 1-petal billiard trajectory}
\COMMENT{think about the subsection title}
\begin{definition}[1-petal  billiard trajectory details]

If $x\in\partial M_r$, which represents the collisions in a corner, then $x=(\phistar,\theta)$ with $\theta\in(0,\pi-\phistar)$, or $x=(2\pi-\phistar,\theta)\in\partial M_r$ with $\theta\in(\phistar,\pi)$, i.e. $x$ is on the open line segment $D_1E_2$ or $D_3E_4$ in \cref{fig:MrBoundaryExtension1}. Then we set $$\overline{F}_{\phistar}\big(x=(\phi(\text{mod}2\pi),\theta)\big)=(\phi+2\theta(\text{mod}2\pi),\theta).$$ So, they are considered representing reflections on $\Gamma_r$ and they are the continuous extension of $\overline{F}_{\phistar}$ from $\Int M_r$.
\end{definition}

For $x\in\partial\Mrout$, we consider $p(x)$ as the collision position on $\Gamma_r$ (maybe $A$ or $B$). And at $p(x)$, the billiard flow direction represented by $x$ points to the corner $A$ or $B$. To resume the billiard flow when it hits the corner, there are two options: reflecting on the boundary $\Gamma_r$ or on the boundary $\overline{AB}$. This will extend the $\PET(r,\phistar)$ billiard map to be multivalued since for $x\in\partial \Mrout$ there are two options(branches) of value $\overline{F}_{1,\phistar}$ and $\overline{F}_{2,\phistar}$.\COMMENT{What follows looks incomprehensible. We need to talk about this---what is the purpose and how can it be achieved?} 


\begin{definition}[Multi-valued definitions for $\overline{F}_{\phistar}$'s iteration]\strut\COMMENT{This runs several pages!!! This needs to be made readable}

    \textbf{Branch (1)} 
    
    For $x\in\partial\Mrout$, the following defined $\overline{F}_{1,\phistar}(x)$ on $\partial \Mrout\sqcup\partial M_{\text{\upshape{f}}}$ resumes the collision hitting corner as hitting / reflecting on the boundary $\overline{AB}$. That is $\overline{F}_{1,\phistar}(\partial \Mrout)\subset \partial M_{\text{\upshape{f}}}$ and $\overline{F}_{1,\phistar}(\partial M_{\text{\upshape{f}}})=\partial\Mrin$
    
    The branch (1) values of $\overline{F}_{\phistar}(x)$, $\overline{F}^2_{\phistar}(x)$ are described as

\begin{align*}
    \overline{F}_{\phistar}(x)=\overline{F}_{1,\phistar}(x)\dfn&\lim_{\Mrout\ni x_0\to x\in\partial \Mrout}\overline{F}_{\phistar}(x_0)\in\partial M_{\text{\upshape{f}}}\\
    \text{ and }\overline{F}^2_{\phistar}(x)=\overline{F}_{1,\phistar}(\overline{F}_{1,\phistar}(x))\dfn&\lim_{\Mrout\ni x_0\to x\in\partial \Mrout}\overline{F}^2_{\phistar}(x_0)\in\partial\Mrin\\
\end{align*} 

By definition, it is the continuous extension from the points of $\Int(\Mrout)$ to $\partial\Mrout$. It resumes the hitting corner trajectory by first hitting / reflecting at the corner with boundary $\overline{AB}$ and then hitting on boundary $\Gamma_r$ which are the same boundaries in the order hit / reflecting on in the next two steps by its nearby points in $\Int(\Mrout)$.

$F^2_{1,\phistar}$ is the diffeomorphism between the same number labeled piece-wise segments of $\partial\Mrout$ and $\partial\Mrin$ in \cref{fig:MrBoundaryExtension1}. Hence, in branch (1), $\overline{F}^2_{\phistar}$ is by definition the homeomorphism between $\text{closure}(\Mrout)$ and $\text{closure}(\Mrin)$. 

\textbf{Branch (2)}

In the following case (i) and (ii) $x\in\partial\Mrout$, the defined $\overline{F}_{2,\phistar}(x)$ resumes the collision hitting corner as hitting /reflecting on the boundary $\Gamma_r$. 

Case (i) If $x=(\phi,\theta)\in\partial\Mrout$ with $\pi<\theta<\pi-\phistar$, more specifically $x$ is on line segment $F_2C_4$ or $F_4C_2$ in \cref{fig:MrBoundaryExtension1}, then the branch (2) values of $\overline{F}_{\phistar}(x)$, $\overline{F}^2_{\phistar}(x)$ are described as
\begin{equation}\label{eq:extdefF1}
\begin{aligned}
    \overline{F}_{\phistar}(x)=\overline{F}_{2,\phistar}(x)\dfn&\lim_{\Int(M_r\setminus\Mrout)\ni x_0\to x\in\partial \Mrout}\overline{F}_{\phistar}(x_0)\in\partial M_{r}\setminus\partial \Mrout\\
    \text{ and }\overline{F}^2_{\phistar}(x)=\overline{F}^2_{2,\phistar}(x)\dfn&\lim_{\Int(M_r\setminus\Mrout)\ni x_0\to x\in\partial \Mrout}\overline{F}^2_{\phistar}(x_0)\in\Int(M_r),
\end{aligned}
\end{equation} which makes it as a continuous extension from $\Int(M_r\setminus\Mrout)$ to $\partial\Mrout$. It resumes the hitting corner trajectory by first hitting / reflecting at the corner with boundary $\Gamma_r$ and then again hitting on $\Int(\Gamma_r)$ which are the same boundaries in the order hit / reflecting on in the next two steps by its nearby points in $\Int(M_r\setminus\Mrout)$.

Case (ii) If $x=(\phi,\theta)\in\partial\Mrout$ with  $\theta\in(0,\phistar)\sqcup(\pi-\phistar,\pi)$ and $\phi\ne\phistar,2\pi-\phistar$, more specifically $x$ is on line segment $D_4F_4$ or $D_2F_2$ in \cref{fig:MrBoundaryExtension2}, then we define $$\overline{F}_{\phistar}(x)=\overline{F}_{2,\phistar}(x)\dfn\quad\lim_{\Int(M_r\setminus\Mrout)\ni x_0\to x\in\partial \Mrout}\overline{F}_{\phistar}(x_0)\in\partial \Mrout$$
And let $y=\overline{F}_{2,\phistar}(x)$, then $y\in\partial\Mrout$ and more specifically $y$ is on line segment $E_4D_4$ or $E_2D_2$ in \cref{fig:MrBoundaryExtension2}. We define\COMMENT{There is a lot going on here, and this can't all be definitions. Those should use :=\Hrule The alignment is unhelpful.}\begin{align*}
    \overline{F}^2_{\phistar}(x)&=\overline{F}_{1,\phistar}(\overline{F}_{2,\phistar}(x))=\overline{F}_{1,\phistar}(y)\overbracket{=}^{\mathllap{\text{by equation }\eqref{eq:extdefF1}\in\partial M_{\text{\upshape{f}}}}}\lim_{\Mrout\ni y_0\to y\in\partial \Mrout}\overline{F}_{\phistar}(y_0)\\
    &\overbracket{=}^{\mathllap{\substack{\text{In small enough nbd of }x\text{,}\\\overline{F}_{\phistar}\text{ is a homeomorphism.}}}}\lim_{\Mrout\ni y_0\to y\in\partial \Mrout}\overline{F}^2_{\phistar}(\overbracket{\overline{F}^{-1}_{\phistar}(y_0)}^{\mathrlap{\substack{\text{For }y_0\text{ in small enough nbd of }y,\\x_0=\overline{F}^{-1}_{\phistar}(y_0)\in\Int(M_r\setminus\Mrout),\\\text{and since }y_0\to y,\text{ it implies } x_0\to x}}})=\lim_{\Int(M_r\setminus\Mrout)\ni x_0\to x}\overline{F}^2_{\phistar}(x_0)\in\partial M_{\text{\upshape{f}}},
\end{align*}


And then we define \begin{align*}
    \overline{F}^3_{\phistar}(x)&=\overline{F}_{1,\phistar}(\overline{F}_{1,\phistar}(\overline{F}_{2,\phistar}(x)))=\overline{F}^2_{1,\phistar}(y)\overbracket{=}^{\mathllap{\text{by equation }\eqref{eq:extdefF1}\in\partial\Mrin}}\lim_{\Mrout\ni y_0\to y\in\partial \Mrout}\overline{F}^2_{\phistar}(y_0)\\
    &\overbracket{=}^{\mathllap{\substack{\text{In small enough nbd of }x\text{,}\\\overline{F}_{\phistar}\text{ is a homeomorphism.}}}}\lim_{\Mrout\ni y_0\to y\in\partial \Mrout}\overline{F}^3_{\phistar}(\overbracket{\overline{F}^{-1}_{\phistar}(y_0)}^{\mathrlap{\substack{\text{For }y_0\text{ in small enough nbd of }y,\\x_0=\overline{F}^{-1}_{\phistar}(y_0)\in\Int(M_r\setminus\Mrout),\\\text{and since }y_0\to y,\text{ it implies } x_0\to x}}})=\lim_{\Int(M_r\setminus\Mrout)\ni x_0\to x}\overline{F}^3_{\phistar}(x_0)\in\partial\Mrin,
\end{align*}
which makes it as a continuous extension from $\Int(M_r\setminus\Mrout)$ to $\partial\Mrout$. It resumes the hitting corner trajectory by first hitting / reflecting at the corner with boundary $\Gamma_r$, next hitting on boundary $\overline{AB}$, and third hitting on boundary $\Gamma_r$ which are the same boundaries in the order hit / reflecting in the next three steps by its nearby points in $\Int(M_r\setminus\Mrout)$.

\end{definition}\strut\COMMENT{Can this truly be the end of this definition?}

\begin{figure}[h]
    \begin{center}
    \begin{tikzpicture}[xscale=1.08,yscale=1.08]
        \pgfmathsetmacro{\PHISTAR}{0.8}
        \pgfmathsetmacro{\WIDTH}{7}
        \pgfmathsetmacro{\LENGTH}{7}
        \tkzDefPoint(-\LENGTH,0.5*\WIDTH){LLM};
        \tkzDefPoint(-\LENGTH,0){LLL};
        \tkzDefPoint(-\LENGTH+\PHISTAR,0){LL};
        \tkzDefPoint(-\LENGTH+\PHISTAR,\WIDTH){LU};
        \tkzDefPoint(-\LENGTH,\WIDTH){LLU};
        \tkzDefPoint(\LENGTH-\PHISTAR,0){RL};
        \tkzDefPoint(\LENGTH-\PHISTAR,\WIDTH){RU};
        \tkzDefPoint(\LENGTH,0){RRL};
        \tkzDefPoint(\LENGTH,\WIDTH){RRU};
        \tkzDefPoint(-\LENGTH+\PHISTAR,\PHISTAR){X};
        \tkzDefPoint(-\LENGTH+\PHISTAR,\WIDTH-\PHISTAR){IX};
        \tkzDefPoint(\LENGTH-\PHISTAR,\PHISTAR){IY};
        \tkzDefPoint(\LENGTH-\PHISTAR,\WIDTH-\PHISTAR){Y};

        \tkzDefPoint(-\LENGTH+3*\PHISTAR,2*\PHISTAR){C1};
        \tkzDefPoint(-\LENGTH+3*\PHISTAR,\WIDTH-2*\PHISTAR){C2};
        \tkzDefPoint(\LENGTH-3*\PHISTAR,\WIDTH-2*\PHISTAR){C3};
        \tkzDefPoint(\LENGTH-3*\PHISTAR,2*\PHISTAR){C4};
        
        \tkzLabelPoint[left](LLL){0};
        \tkzLabelPoint[left](LLM){$\frac{\pi}{2}$};
        \tkzLabelPoint[left](LLU){$\pi$};
        \tkzLabelPoint[below](LLL){0};
        \tkzLabelPoint[below](RRL){\ $2\pi$};
        \tkzLabelPoint[right](IY){$x_*=E_4$};
        \tkzLabelPoint[left](IX){$E_2=y_*$};
        \tkzLabelPoint[right](Y){\small $I x_*=E_3$};
        \tkzLabelPoint[left](X){\small $E_1=Iy_*$};
        
        \tkzLabelPoint[above left](C1){\small $C_1=(3\phistar,2\phistar)$};
        \tkzLabelPoint[below left](C2){\small $C_2=(3\phistar,\pi-2\phistar)$};
        \tkzLabelPoint[below right](C3){\small $C_3=(2\pi-3\phistar,\pi-2\phistar)$};
        \tkzLabelPoint[above right](C4){\small $C_4=(2\pi-3\phistar,2\phistar)$};
        \tkzLabelPoint[below right](LL){\small $D_1=(\phistar,0)$};
        \tkzLabelPoint[above right](LU){\small $D_2=(\phistar,\pi)$};
        \tkzLabelPoint[above left](RU){\small $D_3=(2\pi-\phistar,\pi)$};
        \tkzLabelPoint[below left](RL){\small $D_4=(2\pi-\phistar,0)$};

        \tkzLabelSegment[above right](C4,IY){\small $1$};
        \tkzLabelSegment[right](IY,RL){\small $2$};
        \tkzLabelSegment[below left](RL,C2){\small $3$};
        \tkzLabelSegment[below left](C2,IX){\small $4$};
        \tkzLabelSegment[left](IX,LU){\small $5$};
        \tkzLabelSegment[above right](LU,C4){\small $6$};

        \tkzLabelSegment[left](LL,X){\small $1$};
        \tkzLabelSegment[above left](X,C1){\small $2$};
        \tkzLabelSegment[above left](C1,RU){\small $3$};
        \tkzLabelSegment[right](RU,Y){\small $4$};
        \tkzLabelSegment[below right](Y,C3){\small $5$};
        \tkzLabelSegment[below right](C3,LL){\small $6$};


        \draw [thick,-triangle 45](C4)--(IY);
        \draw [thick,-triangle 45](IY)--(RL);
        \draw [thick,-triangle 45](RL)--(C2);
        \draw [thick,-triangle 45](C2)--(IX);
        \draw [thick,-triangle 45](IX)--(LU);
        \draw [thick,-triangle 45](LU)--(C4);

        \draw [thick,-triangle 45](LL)--(X);
        \draw [thick,-triangle 45](X)--(C1);
        \draw [thick,-triangle 45](C1)--(RU);
        \draw [thick,-triangle 45](RU)--(Y);
        \draw [thick,-triangle 45](Y)--(C3);
        \draw [thick,-triangle 45](C3)--(LL);
        
        \draw [ultra thin,dashed] (LL) --(LU);
        \draw [ultra thin,dashed] (LU) --(RU);
        \draw [ultra thin,dashed] (RL) --(RU);
        \draw [ultra thin,dashed] (LL) --(RL);
        \draw [ultra thin,dashed](LLL)--(LL);
        \draw [ultra thin,dashed](LLU)--(LU);
        \draw [ultra thin,dashed](LLU)--(LLL);
        \draw [ultra thin,dashed](RRL)--(RL);
        \draw [ultra thin,dashed](RRU)--(RU);
        \draw [ultra thin,dashed](RRU)--(RRL);
        \fill [blue, opacity=3/30](LL) -- (Y) -- (RU) -- (X) -- cycle;
        \fill [red, opacity=3/30](LU) -- (IY) -- (RL) -- (IX) -- cycle;
        
        \coordinate (little_x) at ($(RL)!0.65!(IY)$);
        \begin{scope}
            \clip (RL) -- (IY) -- (LU) -- cycle;
            \draw[fill=orange] circle[at=(little_x),radius=0.12];
        \end{scope}
        \coordinate (F1little_x) at ($(C1)!0.35!(X)$);
        \begin{scope}
            \clip (LL) -- (X) -- (RU) -- cycle;
            \draw[fill=orange] ellipse[rotate=27, at=(F1little_x),x radius=0.201, y radius =0.102];
        \end{scope}
        \tkzLabelPoint[left](little_x){$U_1$};
        \tkzLabelPoint[below right](F1little_x){$F^2_{\phistar}(U_1)$};
        
        \coordinate (little_x_2) at ($(C4)!0.2!(LU)$);
        \begin{scope}
            \clip (IX) -- (LU) -- (IY) -- (RL) -- cycle;
            \draw[fill=yellow] ellipse[rotate=0, at=(little_x_2),radius=0.12];
        \end{scope}

        \coordinate (F1little_x_2) at ($(C3)!0.68!(LL)$);
        \begin{scope}
            \clip (LL) -- (X) -- (RU) -- (Y) -- cycle;
            \draw[fill=yellow] ellipse[rotate=30, at=(F1little_x_2),x radius=0.36, y radius =0.09];
        \end{scope}
        \tkzLabelPoint[right](little_x_2){$U_2$};
        \tkzLabelPoint[below right](F1little_x_2){$F^2_{\phistar}(U_2)$};

        \coordinate (little_x_3) at ($(C2)!0.9!(RL)$);
        \begin{scope}
            \clip (IX) -- (LU) -- (IY) -- (RL) -- cycle;
            \draw[fill=green] ellipse[rotate=0, at=(little_x_3),radius=0.12];
        \end{scope}

        \coordinate (F1little_x_3) at ($(C1)!0.13!(RU)$);
        \begin{scope}
            \clip (LL) -- (X) -- (RU) -- (Y) -- cycle;
            \draw[fill=green] ellipse[rotate=15, at=(F1little_x_3),x radius=0.23, y radius =0.06];
        \end{scope}
        \tkzLabelPoint[left](little_x_3){$U_3$};
        \tkzLabelPoint[above left](F1little_x_3){$F^2_{\phistar}(U_3)$};
        \tkzDrawPoints(X,IX,Y,IY,RL,RU,LL,LU,C1,C2,C3,C4,little_x,F1little_x,little_x_2,F1little_x_2,little_x_3,F1little_x_3);

    \end{tikzpicture}
    \end{center}
    \caption{$\PET(r,\phistar)$ billiard branch (1) extended definition for $x\in\partial\Mrout$.}\label{fig:MrBoundaryExtension1}    
\end{figure}

\COMMENT{What is this?}Line segment 1: $C_4E_4\subset\partial\Mrout$ is $\overline{F}^2_{1,\phistar}-$ diffeomorphic to line segment 1: $D_1E_1\subset\partial\Mrin$.
    
    Line segment 2: $E_4D_4\subset\partial\Mrout$ is $\overline{F}^2_{1,\phistar}-$ diffeomorphic to line segment 2: $E_1C_1\subset\partial\Mrin$.
    
    Line segment 3: $D_4C_2\subset\partial\Mrout$ is $\overline{F}^2_{1,\phistar}-$ diffeomorphic to line segment 3: $C_1D_3\subset\partial\Mrin$.
    
    Line segment 4: $C_2E_2\subset\partial\Mrout$ is $\overline{F}^2_{1,\phistar}-$ diffeomorphic to line segment 4: $D_3E_3\subset\partial\Mrin$.
    
    Line segment 5: $E_2D_2\subset\partial\Mrout$ is $\overline{F}^2_{1,\phistar}-$ diffeomorphic to $\text{line segment } 5: E_3C_3\subset\partial\Mrin$.
    
    Line segment 6: $D_2C_4\subset\partial\Mrout$ is $\overline{F}^2_{1,\phistar}-$ diffeomorphic to line segment 6: $C_3D_1\subset\partial\Mrin$.

In \cref{fig:MrBoundaryExtension1}, for $x\in$ one of the number-labeled line segments of $\partial\Mrout$, \newline $\overline{F}^2_{\phistar}(\text{nbd}(x)\cap\Mrout)= \text{some nbd}(\overline{F}_{1,\phistar}(x))\cap\Mrin$ for which $\overline{F}_{1,\phistar}(x)$ is on the same number-labeled line segment of $\partial\Mrin$.

For example, the $U_3$ (a nbd of the $x\in$ line segment $3$) in $\Mrout$ is diffeomorphically mapped by $\overline{F}^2_{\phistar}$ to the green set (a nbd of $\overline{F}_{1,\phistar}(x)\in$ line segment $3$) in $\Mrin$.

The $U_2$ (a nbd of the $x\in$ line segment $6$) in $\Mrout$ is diffeomorphically mapped by $\overline{F}^2_{\phistar}$ to the yellow set (a nbd of $\overline{F}_{1,\phistar}(x)\in$ line segment $6$) in $\Mrin$.

The $U_1$ (a nbd of $x\in$ line segment $2$) in $\Mrout$ is diffeomorphically mapped by $\overline{F}^2_{\phistar}$ to the orange set (a nbd of $\overline{F}_{1,\phistar}(x)\in$ line segment $2$) in $\Mrin$.

\begin{figure}[ht]
    \begin{center}
    \begin{tikzpicture}[xscale=1.08,yscale=1.08]
        \pgfmathsetmacro{\PHISTAR}{0.8}
        \pgfmathsetmacro{\WIDTH}{7}
        \pgfmathsetmacro{\LENGTH}{7}
        \tkzDefPoint(-\LENGTH,0.5*\WIDTH){LLM};
        \tkzDefPoint(-\LENGTH,0){LLL};
        \tkzDefPoint(-\LENGTH+\PHISTAR,0){LL};
        \tkzDefPoint(-\LENGTH+\PHISTAR,\WIDTH){LU};
        \tkzDefPoint(-\LENGTH,\WIDTH){LLU};
        \tkzDefPoint(\LENGTH-\PHISTAR,0){RL};
        \tkzDefPoint(\LENGTH-\PHISTAR,\WIDTH){RU};
        \tkzDefPoint(\LENGTH,0){RRL};
        \tkzDefPoint(\LENGTH,\WIDTH){RRU};
        \tkzDefPoint(-\LENGTH+\PHISTAR,\PHISTAR){X};
        \tkzDefPoint(-\LENGTH+\PHISTAR,\WIDTH-\PHISTAR){IX};
        \tkzDefPoint(\LENGTH-\PHISTAR,\PHISTAR){IY};
        \tkzDefPoint(\LENGTH-\PHISTAR,\WIDTH-\PHISTAR){Y};

        \tkzDefPoint(-\LENGTH+3*\PHISTAR,2*\PHISTAR){C1};
        \tkzDefPoint(-\LENGTH+3*\PHISTAR,\WIDTH-2*\PHISTAR){C2};
        \tkzDefPoint(\LENGTH-3*\PHISTAR,\WIDTH-2*\PHISTAR){C3};
        \tkzDefPoint(\LENGTH-3*\PHISTAR,2*\PHISTAR){C4};

        \tkzDefPoint(-\LENGTH+3*\PHISTAR,\PHISTAR){F1};
        \tkzDefPoint(-\LENGTH+3*\PHISTAR,\WIDTH-\PHISTAR){F2};
        \tkzDefPoint(\LENGTH-3*\PHISTAR,\WIDTH-\PHISTAR){F3};
        \tkzDefPoint(\LENGTH-3*\PHISTAR,\PHISTAR){F4};

        \tkzLabelPoint[left](LLL){0};
        \tkzLabelPoint[left](LLM){$\frac{\pi}{2}$};
        \tkzLabelPoint[left](LLU){$\pi$};
        \tkzLabelPoint[below](LLL){0};
        \tkzLabelPoint[below](RRL){\ $2\pi$};
        \tkzLabelPoint[right](IY){$x_*=E_4$};
        \tkzLabelPoint[left](IX){$E_2=y_*$};
        \tkzLabelPoint[right](Y){\small $I x_*=E_3$};
        \tkzLabelPoint[left](X){\small $E_1=Iy_*$};
        
        \tkzLabelPoint[below right](LL){\small $D_1=(\phistar,0)$};
        \tkzLabelPoint[above right](LU){\small $D_2=(\phistar,\pi)$};
        \tkzLabelPoint[above left](RU){\small $D_3=(2\pi-\phistar,\pi)$};
        \tkzLabelPoint[below left](RL){\small $D_4=(2\pi-\phistar,0)$};

        \tkzLabelPoint[below right](F1){\small $F_1$};
        \tkzLabelPoint[above right](F2){\small $F_2$};
        \tkzLabelPoint[above left](F3){\small $F_3$};
        \tkzLabelPoint[below left](F4){\small $F_4$};

        \tkzLabelPoint[below right](C1){\small $C_1$};
        \tkzLabelPoint[below left](C2){\small $C_2$};
        \tkzLabelPoint[above left](C3){\small $C_3$};
        \tkzLabelPoint[above right](C4){\small $C_4$};
        





        \draw [thick] (X) --(IX);
        \draw [ultra thin,dashed] (LU) --(RU);
        \draw [thick] (IY) --(Y);
        \draw [ultra thin,dashed] (LL) --(RL);
        \draw [ultra thin,dashed](LLL)--(LL);
        \draw [ultra thin,dashed](LLU)--(LU);
        \draw [ultra thin,dashed](LLU)--(LLL);
        \draw [ultra thin,dashed](RRL)--(RL);
        \draw [ultra thin,dashed](RRU)--(RU);
        \draw [ultra thin,dashed](RRU)--(RRL);

        \draw [thick](LU)--(IY);
        \draw [thick](IX)--(RL);
        \draw [thick](X)--(RU);
        \draw [thick](LL)--(Y);

        \draw [thick](LL)--(X);
        \draw [thick](IX)--(LU);
        \draw [thick](Y)--(RU);
        \draw [thick](RL)--(IY);
        
        \fill [blue, opacity=2/30](LL) -- (Y) -- (RU) -- (X) -- cycle;
        \fill [red, opacity=2/30](LU) -- (IY) -- (RL) -- (IX) -- cycle;
        
        \coordinate (little_x) at ($(RL)!0.45!(F4)$);
        \begin{scope}
            \clip (RL) -- (LL) -- (IX) -- cycle;
            \draw[fill=green] circle[at=(little_x),radius=0.12];
        \end{scope}
        \tkzLabelPoint[above](little_x){\small $U$}
        \coordinate (F2little_x) at ($(RL)!0.45!(IY)$);
        \begin{scope}
            \clip (RL) -- (IX) -- (LU) -- (IY)--cycle;
            \draw[fill=green] ellipse[rotate=25, at=(F2little_x),x radius=0.14, y radius =0.097];
        \end{scope}
        \tkzLabelPoint[right](F2little_x){\small $F_{\phistar}(U)$};
        \coordinate (F1F2little_x) at ($(C1)!0.35!(X)$);
        \begin{scope}
            \clip (X) -- (RU) -- (Y) --(LL) -- cycle;
            \draw[fill=green] ellipse[rotate=20, at=(F1F2little_x),x radius=0.26, y radius =0.07];
        \end{scope}
        \tkzLabelPoint[left](F1F2little_x){\small $F^3_{\phistar}(U)$};
        \coordinate (little_x2) at ($(IY)!0.25!(F2)$);
        \begin{scope}
            \clip (IY) -- (LU) -- (RU) -- cycle;
            \draw[fill=yellow] circle[at=(little_x2),radius=0.12];
        \end{scope}
        \tkzLabelPoint[above right](little_x2){\small $V$};
        \coordinate (F2little_x2) at ($(X)!0.25!(IX)$);
        \begin{scope}
            \clip (X) -- (RU) -- (LU)--cycle;
            \draw[fill=yellow] ellipse[rotate=70, at=(F2little_x2),x radius=0.158, y radius =0.096];
        \end{scope}
        \tkzLabelPoint[above right](F2little_x2){\small $F_{\phistar}(V)$};
        \coordinate (F2F2little_x2) at ($(F1)!0.25!(Y)$);
        \begin{scope}
            \clip (Y) -- (LL) -- (RL) -- cycle;
            \draw[fill=yellow] ellipse[rotate=18, at=(F2F2little_x2),x radius=0.29, y radius =0.075];
        \end{scope}
        \tkzLabelPoint[below right](F2F2little_x2){\small $F^2_{\phistar}(V)$};

        

        

        \tkzDrawPoints(X,IX,Y,IY,RL,RU,LL,LU,F1,F2,F3,F4,C1,C2,C3,C4,little_x,F2little_x,F1F2little_x,little_x2,F2little_x2,F2F2little_x2);

    \end{tikzpicture}
    \end{center}
    \caption{$\PET$ billiard branch (2) for $x\in\partial\Mrout$ in case (i) with neighborhood $V$ and in case (ii) with neighborhood $U$.
    }\label{fig:MrBoundaryExtension2}
\end{figure}
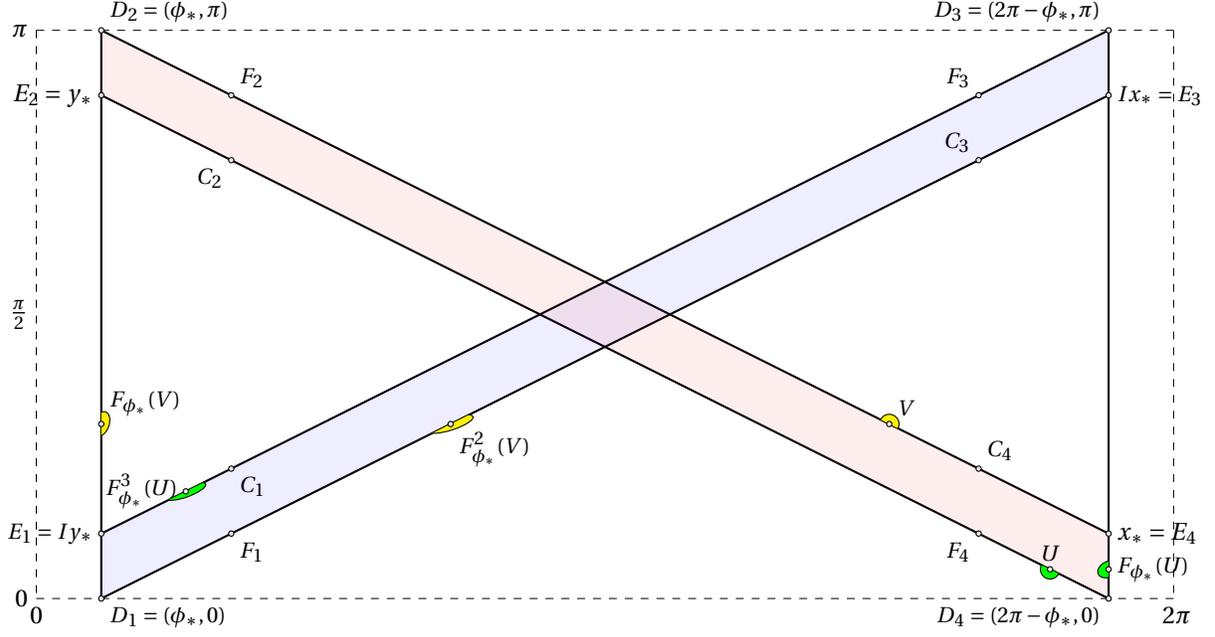\hfill
    
    If $x$ is in case (i) and in \cref{fig:MrBoundaryExtension2} yellow $V$ is the $nbd(x)\cap\Int(M_r\setminus\Mrout)$, then $\overline{F}_{\phistar}(V)$ and $\overline{F}^2_{\phistar}(V)$ are the neighborhood of branch (2) points $\overline{F}_{\phistar}(x)\in\partial M_r$, $\overline{F}^2_{\phistar}(x)\in\partial \Mrin$.
    
    If $x$ is in case (ii) and in \cref{fig:MrBoundaryExtension2} green $U$ is the $nbd(x)\cap\Int(M_r\setminus\Mrout)$, then $\overline{F}_{\phistar}(U)$, $\overline{F}^2_{\phistar}(U)$ and $\overline{F}^3_{\phistar}(U)$ are the neighborhood of branch (2) points $\overline{F}_{\phistar}(x)\in\partial M_r\cap \partial\Mrout$, $\overline{F}^2_{\phistar}(x)\in \partial M_{\text{\upshape{f}}}$ and $\overline{F}^3_{\phistar}(x)\in\partial\Mrin$.

\begin{lemma}\label{lemma:PETExpExpansion}
For a $\PET(r,\phistar)$ billiard table with $\phistar\ne\frac{\pi}{n}$, $\forall n\in\mathbb{N}$, if the trajectory of $x\in\Mrout$: $\hat{\overline{F}}_{\phistar}^k(x)|_{k\ge0}\nfd \hat{\overline{x}}_k|_{k\ge0}$ never hits the corner under the iteration of the return map $\hat{\overline{F}}_{\phistar}$ (\cref{def:LemonAndPetTogether}\eqref{item06LemonAndPetTogether}) on $\Mrout$, then with $\overline{\tau}_{0,k}$, $\overline{\tau}_{1,k}$, $\overline{d}_{0,k}$ defined in \cref{def:LemonAndPetTogether}\eqref{item06LemonAndPetTogether} for each $\hat{\overline{x}}_k$ and for all $n\ge0$, $\prod_{k=1}^{n}\bigm|-1-\frac{\overline{\tau}_{0,k}+\overline{\tau}_{1,k}-2\overline{d}_{0,k}}{\overline{d}_{0,k}}\bigm|>\overline{C}\cdot\overline{E}^n$ with some constant $\overline{C}>0$ and $\overline{E}>1$. i.e. it has an exponential expansion growth.
\end{lemma}
\begin{proof}
    By the definitions of length functions $\overline{\tau}_0$, $\overline{\tau}_1$, $\overline{d}_0$ in \cref{def:LemonAndPetTogether}\eqref{item06LemonAndPetTogether}, in phase space $\overline{\tau}_{0,k}+\overline{\tau}_{1,k}-2\overline{d}_{0,k}>0$ in \cref{fig:WojtkowskiBilliardTable}, which is also the root cause of the defocusing mechanism in \cite[section 8.1,8.2]{cb}.
    
    In addition, we note that $\overline{\tau}_{0,k}+\overline{\tau}_{1,k}-2\overline{d}_{0,k}$ is bounded away from $0$ if $\hat{\overline{x}}_{k}\in\Mrout$ and bounded away from the boundary segments $C_4E_4$ and $C_2E_2$ (segments 1 and 4 in \cref{fig:MrBoundaryExtension1}). And the small enough neighborhood of $C_4E_4$ ($C_2E_2$) is diffeomorphically mapped by $\overline{F}_{\phi_*}$ to the small enough neighborhood of $D_1E_1$ ($D_3E_3$). 
    
    Now we are to provide proof arguments based on the notions for phase space \cref{fig:MrBoundaryExtension1}. Note that in the phase space $E_1=(\phistar,\phistar)$ is the endpoint for the line segment $D_1E_1=\big\{(\phistar,\theta)\in M_r\bigm|\theta\in[0,\phi_*]\big\}$. For $\forall x=(\phi,\theta)\in D_1E_1$, let $k(x)=\max{\big\{k\bigm|p(\overline{F}_{\phi_*}^{i}(x))\in\Gamma_r,i=0,1,\cdots,k\big\}}$. Since $x,\cdots,\overline{F}_{\phi_*}^{k(x)}(x)$ are consecutive collisions on $\Gamma_r$, all the $\theta$ coordinates for $x,\cdots,\overline{F}_{\phi_*}^{k(x)}(x)$ are the same as for the $x$ and in the range $[0,\phi_*]$. Therefore, if $x\in D_1E_1$ is bounded away from $E_1=Iy_*=(\phi_*,\phi_*)$, then $\overline{F}_{\phi_*}^{k(x)}(x)$ is bounded away from $C_4E_4$ and $E_2C_2$, which are line segments with $\theta\in[\phistar,2\phistar]\sqcup[\pi-2\phistar,\pi-\phistar]$.

    Since $\phistar\ne\pi/n$ $\forall n\ge1$, $\overline{F}_{\phistar}^{k(Iy_*)}(Iy_*)\in\Mrout$ and is bounded away from line segments $C_4E_4$ and $C_2E_2$. By continuity of $\overline{F}_{\phi_*}$, there is $N(Iy_*)$: a neighborhood of $Iy_*$ such that $\overline{F}_{\phistar}^{k(Iy_*)}(N(Iy_*))\subset\Mrout$ and $\overline{F}_{\phistar}^{k(Iy_*)}(N(Iy_*))$ is bounded away from $2\pi-\phi_*$. Thus, $\overline{F}_{\phistar}^{k(Iy_*)}(N(Iy_*)\cap D_1E_1)$ is bounded away from $C_4C_4$ (also $C_2E_2$).

    Hence, for $D_1E_1\ni x$ with $x$ belonging to a neighborhood of $Iy_*$ or being outside of this neighborhood of $Iy*$, $\overline{F}_{\phi_*}^{k(x)}(x)$ is bounded away from line segments $C_4E_4$ and $E_2C_2$. Therefore, for every two consecutive $\hat{\overline{x}}_k, \hat{\overline{x}}_{k+1}$, at least one of $\overline{\tau}_{0,k}+\overline{\tau}_{1,k}-2\overline{d}_{0,k}$, $\overline{\tau}_{0,k+1}+\overline{\tau}_{1,k+1}-2\overline{d}_{0,k+1}$ is larger than some constant $\overline{\delta}>0$. On the other hand for all $k$, $\overline{d}_{0,k}\le r$. We get the following.
    \[\bigm|-1-\frac{\overline{\tau}_{0,k}+\overline{\tau}_{1,k}-2\overline{d}_{0,k}}{\overline{d}_{0,k}}\bigm|\bigm|-1-\frac{\overline{\tau}_{0,k+1}+\overline{\tau}_{1,k+1}-2\overline{d}_{0,k+1}}{\overline{d}_{0,k+1}}\bigm|>\bigm|-1-\frac{\overline{\delta}}{r}\bigm|.\]
Hence, if we set $\overline{E}=(1+\frac{\overline{\delta}}{r})^{0.5}$ and $\overline{C}=(1+\frac{\overline{\delta}}{r})^{-0.5}$, then $\prod_{k=1}^{n}\bigm|-1-\frac{\overline{\tau}_{0,k}+\overline{\tau}_{1,k}-2\overline{d}_{0,k}}{\overline{d}_{0,k}}\bigm|>(1+\frac{\overline{\delta}}{r})^{\left\lfloor\frac{n}{2}\right\rfloor}>(1+\frac{\overline{\delta}}{r})^{\cdot{\frac{n}{2}-1}}=\overline{C}\cdot \overline{E}^n$ for $n\ge0$.
\end{proof}

\subsection[Trajectory position and direction closeness  between Lemon billiard and 1-petal billiard]{Trajectory position and direction closeness  between Lemon billiard and 1-petal billiard}\COMMENT{Subsection 8.3 needs explanations and guided tour what's going on, nit just computations.}

\begin{definition}[\cref{fig:Trajectory_LemonCompareWithWojtkowski} flow collision position and angle]\label{def:flowcollsionpositionandangle}\label{Def9_2}\label{Def9_3} With $O_r=(0,0)$, $O_R=(0,b)$ in coordinate system of \cref{def:StandardCoordinateTable}, the $\PET(r,\phistar)$ billiard and the $\LEM(r,R,\phistar)$ billiard in \cref{def:LemonAndPetTogether,fig:1-petalExtensionFromCornerAndLemon} share the same corners $A$, $B$ and the same boundary $\Gamma_r$. 

For $x\in M_R$, we denote by $(\mathcal{X},\mathcal{Y},\mathcal{\Theta})$ of $x$ the position and direction of the $\LEM(r,R,\phistar)$ billiard flow after collision on boundary $\Gamma_R$. 

For $\overline{x}\in M_{\text{\upshape{f}}}$, we denote by $(\overline{\mathcal{X}},\overline{\mathcal{Y}},\overline{\mathcal{\Theta}})$ of $\overline{x}$ the position and direction of the $\PET(r,\phistar)$ billiard flow after collision on boundary $\overline{AB}$.
\end{definition}
\COMMENT{Explain what the whole section is about and how it is needed. Right now it looks like random factoids. It would likely be far better if the theorems were each stated at the earliest time possible in their respective subsections, with some CLEAR explanations what they are saying.\newline\wentao{rewritten/with conclusions put in the beginning}}

\begin{figure}[ht]
\begin{center}
\begin{tikzpicture}[xscale=0.8,yscale=0.8]
    \tkzDefPoint(0,0){Or};
    \tkzDefPoint(0,-4.8){Y};
    \pgfmathsetmacro{\rradius}{4.5};
    \clip(-2.4*\rradius,-1.1*\rradius) rectangle (1.3*\rradius, 0.03*\rradius);
    \pgfmathsetmacro{\phistardeg}{40};
    \pgfmathsetmacro{\XValueArc}{\rradius*sin(\phistardeg)};
    \pgfmathsetmacro{\YValueArc}{\rradius*cos(\phistardeg)};
    \pgfmathsetmacro{\Rradius}{
        \rradius*sin(\phistardeg)/sin(12)
    };
    \pgfmathsetmacro{\bdist}{(-\rradius*cos(\phistardeg))+cos(12)*\Rradius}
    \tkzDefPoint(0,\bdist){OR};
    \tkzDefPoint(\XValueArc,-1.0*\YValueArc){B};    \tkzDefPoint(-1.0*\XValueArc,-1.0*\YValueArc){A};
    \tkzDrawArc[name path=Cr, very thick](Or,B)(A);
    \tkzDrawArc[name path=CR](OR,A)(B);
    \pgfmathsetmacro{\PXValue}{-\Rradius*sin(3)};
    \pgfmathsetmacro{\PYValue}{\bdist-\Rradius*cos(3)};
    \tkzDefPoint(\PXValue,\PYValue){P};
    
    
    \path[name path=Line_P_T0_far](P)-- (176:6);
    \path[name path=Line_P_T_far](P)-- (172:6);
    \path[name intersections={of=Cr and Line_P_T0_far,by={T0}}];
    \path[name intersections={of=Cr and Line_P_T_far,by={T}}];
    \tkzDefPoint(0,\bdist-\Rradius){C};
    \path[name intersections={of=Cr and Line_P_T0_far,by={T0}}];
    \draw[name path=LineAB,thin] (A)--(B);
     \coordinate (PHorizon) at ([shift={(3,0)}]P);
    \draw[red,ultra thin,dashed] (P)--(PHorizon);
    \path[name intersections={of=Line_P_T0_far and LineAB,by={P_bar}}];
    \tkzDefPointBy[projection=onto Or--P](M)\tkzGetPoint{M2};
    
    \begin{scope}[ultra thin,decoration={markings,mark=at position 0.5 with {\arrow{>}}}] 
    \draw[red,postaction={decorate}] (T0)--(P_bar);
    \draw[blue,postaction={decorate}] (T)--(P);
    \end{scope}
    \tkzLabelPoint[left](A){$A$};
    \tkzLabelPoint[below](P){\small $P=p(\mathcal{F}(x))=(\mathcal{X}_1,\mathcal{Y}_1)$};
    \tkzLabelPoint[below left](T0){\small $T_0=p(\overline{x}_0)=(r\sin{\overline{\phi}_0},-r\cos{\overline{\phi}_0})$};
    \tkzLabelPoint[left](T){\small $T=p(x)=(r\sin{\phi},-r\cos{\phi})$};
    \tkzLabelPoint[below](Or){$O_r$};
    \tkzLabelPoint[below](B){$B$};
    \tkzLabelPoint[above right]([shift={(0.15,0.23)}]P_bar){\small $\overline{P}=p(\overline{F}_{\phistar}(\overline{x}_0))=(\overline{\mathcal{X}}_1,\overline{\mathcal{Y}}_1)$};
    \tkzLabelPoint[right]([shift={(0.08,0.15)}]P_bar){\small $\overline{\Theta}_1$};
    \tkzLabelPoint[right]([shift={(0.15,0.15)}]P){\small $\Theta_1$};
    \tkzDefPointBy[projection=onto B--A](T0)\tkzGetPoint{M0};
    \coordinate (T0Mirror) at ($(T0)!2!(M0)$);
    \coordinate (PbarArrow) at ($(T0Mirror)!1.16!(P_bar)$);
    \tkzDefPointBy[projection=onto P--OR](T)\tkzGetPoint{M};
    \coordinate (TMirror) at ($(T)!2!(M)$);
    \coordinate (PArrow) at ($(P)!0.18!(TMirror)$);
    \begin{scope}[ultra thin,decoration={markings,mark=at position 1.0 with {\arrow{>}}}] 
    \draw[red,postaction={decorate},ultra thin] (P_bar)--(PbarArrow);
    \draw[blue,postaction={decorate},ultra thin] (P)--(PArrow);
    \end{scope}
    \tkzMarkAngle[arc=lll,ultra thin,size=0.15](B,P_bar,PbarArrow);
    \tkzMarkAngle[arc=lll,ultra thin,size=0.18](PHorizon,P,PArrow);
    \tkzDrawPoints(A,B,Or,P,T0,T,P_bar);
\end{tikzpicture}
\caption{
$\overline{P}=(\overline{\mathcal{X}}_1,\overline{\mathcal{Y}}_1)=p(\overline{F}_{\phistar}(x_0))$ on $\PET(r,\phistar)$.
\newline $P=(\mathcal{X}_1,\mathcal{Y}_1)=p(\mathcal{F}(x))$ on $\LEM(r,R,\phistar)$.
} \label{fig:Trajectory_LemonCompareWithWojtkowski}
\end{center}
\end{figure}


\COMMENT{Much rewritten.}
\COMMENT{Think about $\bar F$ versus $\overline{F}$. I think it looks much better.}
We now embark on the proof of our first trajectory approximation lemma, \cref{lemma:Trajectory_approximation_1}.\COMMENT{Explain here and elsewhere what this theorem says and what role it plays.\Hrule Is there a good reason why the theorem is not stated right here?\newline\wentao{Agree}}

\begin{lemma}[Shadowing lemma 1: approximation for $\overline{F}_{\phistar}$ trajectory in \cref{fig:Trajectory_LemonCompareWithWojtkowski}]\label{lemma:Trajectory_approximation_1}\hfill\newline
Let $\overline{x}_0=(\overline{\phi}_0,\overline{\theta}_0)\in\Mrout$ and $x=(\phi,\theta)\in\Mrout$. Using the notation from \cref{def:LemonAndPetTogether} for the collision maps $\overline{F}_{\phistar}$ for $\PET(r,\phistar)$ and $\mathcal{F}$ for $\LEM(r,R,\phistar)$, we suppose that $x_1=(\Phi_1,\theta_1)\nfd\mathcal{F}(x)\in M_R$. Suppose that $(\mathcal{X}_1,\mathcal{Y}_1,\Theta_1)$ is the $\LEM(r,R,\phistar)$ billiard flow after collision on boundary $\Gamma_R$ position and direction angle (in \cref{def:flowcollsionpositionandangle}) of $\mathcal{F}(x)\in M_R$, $(\overline{\mathcal{X}}_1,\overline{\mathcal{Y}}_1,\overline{\Theta}_1)$ is the $\PET(r,\phistar)$ billiard flow after collision on boundary $\overline{AB}$ position and direction angle (in \cref{def:flowcollsionpositionandangle}) of $\overline{F}_{\phistar}(\overline{x}_0)\nfd \overline{x}_1\in M_{\text{\upshape{f}}}$. We have the following.

For $\forall\, \varepsilon,\varepsilon_1>0$, if $x$, $\overline{x}_0$ and $R$ satisfy the following \begin{equation}\label{eq:x0barx0ClosingCondition}
    \begin{aligned}
        |\theta-\overline{\theta}_0|<\mathtt{\delta}_1(\overline{\phi}_0,\overline{\theta}_0,\varepsilon,\varepsilon_1)\dfn\min&{\Big\{0.01\sin^2{(\overline{\phi}_0+\overline{\theta}_0)},-0.05\frac{\varepsilon}{r}\sin^3{(\overline{\phi}_0+\overline{\theta}_0)},0.4\varepsilon_1,\frac{\pi}{2}-\frac{\overline{\phi}_0+\overline{\theta}_0}{4},\frac{\overline{\phi}_0+\overline{\theta}_0}{4}\Big\}}\\
        |\phi-\overline{\phi}_0|<\mathtt{\delta}_2(\overline{\phi}_0,\overline{\theta}_0,\varepsilon,\varepsilon_1)=\min&{\Big\{0.01\sin^2{(\overline{\phi}_0+\overline{\theta}_0)},-0.05\frac{\varepsilon}{r}\sin^3{(\overline{\phi}_0+\overline{\theta}_0)},0.4\varepsilon_1,\frac{\pi}{2}-\frac{\overline{\phi}_0+\overline{\theta}_0}{4},\frac{\overline{\phi}_0+\overline{\theta}_0}{4}\Big\}}\\
        R>\mathtt{R}_1(\overline{\phi}_0,\overline{\theta}_0,\varepsilon,\varepsilon_1)\dfn\max&\Big\{r+\frac{23r^2}{-\varepsilon\sin^3{(\overline{\phi}_0+\overline{\theta}_0)}},r+\frac{500r}{\sin^2{(\overline{\phi}_0+\overline{\theta}_0)}},\frac{2r^2\sin^2{\phistar}}{\varepsilon},\frac{2.02r\sin{\phistar}}{0.4\varepsilon_1},\\
        &1700r
        ,\frac{4r}{\sin^2{(\frac{2\pi-\overline{\phi}_0-\overline{\theta}_0}{4})}},\frac{4r}{\sin^2{(\frac{\overline{\theta}_0+\overline{\phi}_0-\pi}{4})}}\Big\},
    \end{aligned}
\end{equation}\newline then 
\begin{equation}\label{eq:eq:x0barx0Closing}
        \begin{aligned}
            |\mathcal{X}_1-\overline{\mathcal{X}}_1|<&\varepsilon\\
            |\mathcal{Y}_1-\overline{\mathcal{Y}}_1|<&\varepsilon\\
            |\Theta_1-\overline{\Theta}_1|<&\varepsilon_1.\\
        \end{aligned}
    \end{equation} And especially if \begin{equation}\label{eq:xorbitA0Condition}
        \begin{aligned}
            \bigm|\theta-\overline{\theta}_0\bigm|<&\min\big\{\frac{\pi}{2}-\frac{\overline{\phi}_0+\overline{\theta}_0}{4},\frac{\overline{\phi}_0+\overline{\theta}_0}{4}\big\}\\
            \bigm|\phi-\overline{\phi}_0\bigm|<&\min\big\{\frac{\pi}{2}-\frac{\overline{\phi}_0+\overline{\theta}_0}{4},\frac{\overline{\phi}_0+\overline{\theta}_0}{4}\big\}\\
            R>&\max\Big\{\frac{4r}{\sin^2{(\frac{2\pi-\overline{\phi}_0-\overline{\theta}_0}{4})}},\frac{4r}{\sin^2{(\frac{\overline{\theta}_0+\overline{\phi}_0-\pi}{4})}}\Big\}
        \end{aligned}
    \end{equation} then $\sin{\theta_1}>\sqrt{4r/R}$, which means that if the $\LEM(r,R,\phistar)$ billiard $\hat{M}$ return orbit in \cref{def:MhatReturnOrbitSegment} with $x_0=x=(\phi,\theta)$ satisfying condition \eqref{eq:x0barx0ClosingCondition}, then $x_1=\Fcal(x)=\Fcal(x_0)$ is in the case (a0) in \eqref{eqMainCases}.
\end{lemma}
\begin{proof}
    First note that if $x=(\phi,\theta)\in\Mrout$ and $\overline{x}_0=(\overline{\phi}_0,\overline{\theta}_0)\in\Mrout$, then by checking \cref{fig:Mr} to see that $\Mrout$ region is below line $\big\{(\phi,\theta)\bigm|\phi+\theta=2\pi\big\}$ and above line $\big\{(\phi,\theta)\bigm|\phi+\theta=\pi\big\}$. So, we have $\pi<\phi+\theta<2\pi$ and $\pi<\overline{\phi}_0+\overline{\theta}_0<2\pi$, thus $\sin(\overline{\phi}_0+\overline{\theta}_0)<0$, $\sin(\phi+\theta)<0$.

    In \cref{fig:Trajectory_LemonCompareWithWojtkowski} the coordinate system with $O_r=(0,0)$, $(0,b)$ and $\overline{AB}$ is on line $y=-r\cos\phistar$, we have the following coordinates of points.
    \begin{equation}\label{eq:xT0xTclose}
        \begin{aligned}
        (x_T,y_T)=T\dfn&p(x)=(r\sin\phi,-r\cos\phi)\in\Gamma_r\\ (x_{T_0},y_{T_0})=T_0\dfn&p(\overline{x}_0)=(r\sin{\overline{\phi}_0},-r\cos{\overline{\phi}_0})\in\Gamma_r\\
        P\dfn&p(\Fcal(x))=(\mathcal{X}_1,\mathcal{Y}_1)\in\Gamma_R\\
        \overline{P}\dfn&p(\overline{F}_{\phistar}(\overline{x}_0))=(\overline{\mathcal{X}}_1,\overline{\mathcal{Y}}_1)\in\overline{AB}\text{ with }\overline{\mathcal{Y}}_1=-r\cos\phistar.\\
        \end{aligned}
    \end{equation}

    On the other hand, by the geometric meaning of $\phi$ and $\theta$, given counterclockwise as the positive orientation for the billiard plane, the vector $\overrightarrow{TP}$ has an angle $\phi+\theta$ with respect to the positive axis $x$ and the vector $\overrightarrow{T_0\overline{P}}$ has an angle $\overline{\phi}_0+\overline{\theta}_0$ with respect to the positive axis $x$. Note that $2\pi-\overline{\Theta}_1$ is also the angle between $\overrightarrow{T_0\overline{P}}$ and the positive axis $x$ given the counterclockwise orientation. Therefore, we get the following equations.
    \begin{equation}\label{eq:X1barY1barTheta1bar}
        \begin{aligned}
            2\pi-\overline{\Theta}_1&=\overline{\phi}_0+\overline{\theta}_0\\
            \frac{-r\cos\phi+r\cos\phistar}{r\sin\phi-\overline{\mathcal{X}}_1}\overset{\eqref{eq:xT0xTclose}}=\frac{y_T-\overline{\mathcal{Y}}_1}{x_T-\overline{\mathcal{X}}_1}&\overbracket{=}^{\mathrlap{\text{slope of }T\overline{P}}}\tan{(\overline{\phi}_0+\overline{\theta}_0)}=\frac{\sin{(\overline{\phi}_0+\overline{\theta}_0)}}{\cos{(\overline{\phi}_0+\overline{\theta}_0)}}\\
            \Longrightarrow \overline{P}=(\overline{\mathcal{X}}_1,\overbracket{\overline{\mathcal{Y}}_1}^{\mathllap{\eqref{eq:xT0xTclose}:=-r\cos\phistar}})&=(r\sin{\overline{\phi}_0}+\frac{r\cos{\overline{\phi}_0}-r\cos{\phistar}}{\sin{(\overline{\phi}_0+\overline{\theta}_0)}}\cos{(\overline{\phi}_0+\overline{\theta}_0)},-r\cos{\phistar}).
        \end{aligned}
    \end{equation} Let $t=|PT|$, since $P$ is on the circle $C_R$, we get the following equations for $(\mathcal{X}_1,\mathcal{Y}_1)$.
\begin{equation}\label{eq:X1Y1Polynomial}
\begin{aligned}
    &\left.\begin{aligned}\mathcal{X}_1=&r\sin{\phi}+t\cos{(\phi+\theta)}\\
        \mathcal{Y}_1=&-r\cos{\phi}+t\sin{(\phi+\theta)}\\
       R^2=&\mathcal{X}^2_1+(\mathcal{Y}_1-b)^2\\
    \end{aligned}\right\}\Longrightarrow R^2=\big(\overbracket{r\sin{\phi}+t\cos{(\phi+\theta)}}^{=\mathcal{X}_1}\big)^2+\big(\overbracket{-r\cos{\phi}+t\sin{(\phi+\theta)}-b}^{=\mathcal{Y}_1-b}\big)^2\\
    \Longrightarrow&R^2=r^2+t^2+b^2+\overbracket{2rt\sin\phi\cos{(\phi+\theta)}-2rt\cos\phi\sin{(\phi+\theta)}}^{-2rt\sin\theta}-2bt\sin{(\phi+\theta)}+2br\cos\phi\\
    \Longrightarrow&t^2-2rt\sin\theta-2bt\sin{(\phi+\theta)}+2br\cos\phi+r^2+b^2-R^2=0\\
    \overbracket{\Longrightarrow}^{\mathrlap{\substack{\text{Cosine Law: }\\R^2=r^2+b^2-2br\cos{\angle{O_RO_rB}}\\
    \measuredangle{O_RO_rB}=\pi-\phistar}}}&t^2-2\big(r\sin\theta+b\sin{(\phi+\theta)}\big)t+2br\cos\phi-2br\cos\phistar=0.
    \end{aligned} 
\end{equation}
    Since $\pi\in(\phistar,2\pi-\phistar)$, quadratic polynomial $t^2-2\big(r\sin\theta+b\sin{(\phi+\theta)}\big)t+2br\cos\phi-2br\cos\phistar$ has $|PT|$ as its positive root $t_+$
    \begin{equation}\label{eq:PTposRoot}
    \begin{aligned}|PT|=t_+=&b\sin{(\phi+\theta)}+r\sin\theta+\sqrt{\big(b\sin{(\phi+\theta)}+r\sin\theta\big)^2+\big(2br\cos\phistar-2br\cos\phi\big)}\\
    =&\frac{2br\cos\phistar-2br\cos\phi}{-b\sin{(\phi+\theta)}-r\sin\theta+\sqrt{\big(b\sin{(\phi+\theta)}+r\sin\theta\big)^2+\big(2br\cos\phistar-2br\cos\phi\big)}}\\
    =&\frac{2r\cos\phistar-2r\cos\phi}{-\sin{(\phi+\theta)}-(r/b)\sin\theta+\sqrt{\big(\sin{(\phi+\theta)}+(r/b)\sin\theta\big)^2+\big(2(r/b)\cos\phistar-2(r/b)\cos\phi\big)}}.
    \end{aligned}
    \end{equation}
\newline Note that with $r$ fixed, as $R\rightarrow\infty$, $b\overbracket{>}^{\mathllap{\substack{b,r,R\text{ are edge lengths of }\\ \triangle O_rO_RB\text{ in \cref{fig:1-petalExtensionFromCornerAndLemon}}}}}R-r\rightarrow\infty$. So, with $(r/b)\rightarrow 0$, $\phi\rightarrow\overline{\phi}_0$, $\theta\rightarrow\overline{\theta}_0$, the right-hand side of the last equality in \eqref{eq:PTposRoot} has the limit $t_+\longrightarrow\frac{2r\cos\phistar-2r\cos{\overline{\phi}_0}}{-\sin{(\overline{\phi}_0+\overline{\theta}_0)}+|\sin{(\overline{\phi}_0+\overline{\theta}_0)}|}\overbracket{=}^{\qquad\qquad\mathllap{\substack{\pi<\overline{\phi}_0+\overline{\theta}_0<2\pi\\\sin{(\overline{\phi}_0+\overline{\theta}_0)}<0}}}\frac{r\cos\phistar-r\cos{\overline{\phi}_0}}{-\sin{(\overline{\phi}_0+\overline{\theta}_0)}}$. Therefore, if we write 
\begin{equation}\label{eq:xietaphitheta}
    \begin{aligned}
        \xi(\phi,\theta)=&(2r\cos{\phistar}-2r\cos{\phi})\cos{(\phi+\theta)},\\
        \eta(\phi,\theta)=&\sqrt{\sin^2{(\theta+\phi)}+(r/b)^2\sin^2{\theta}+2(r/b)\sin{\theta}\sin{(\theta+\phi)}+2(r/b)\cos{\phistar}-2(r/b)\cos{\phi}},\\
        &-\sin{(\theta+\phi)}-(r/b)\sin{\theta}\\
        \frac{\xi(\phi,\theta)}{\eta(\phi,\theta)}=t_+\cos{(\phi+\theta)}\overbracket{\longrightarrow}^{\mathllap{\substack{(r/b)\to0,\\\phi\to\overline{\phi}_0,\\\theta\to\overline{\theta}_0}}}&\frac{r\cos\phistar-r\cos{\overline{\phi}_0}}{-\sin{(\overline{\phi}_0+\overline{\theta}_0)}}\cos{(\overline{\phi}_0+\overline{\theta}_0)},
    \end{aligned}
\end{equation}
\newline then by \eqref{eq:PTposRoot}, \eqref{eq:xietaphitheta}, \eqref{eq:X1Y1Polynomial}, \eqref{eq:X1barY1barTheta1bar} we conclude
\begin{equation}\label{eq:X1X1bardiff}
\begin{aligned}\mathcal{X}_1\overbracket{=}^{\mathllap{\text{\eqref{eq:X1Y1Polynomial}}}}&r\sin{\phi}+\overbracket{t_+\cos{(\phi+\theta)}}^{\text{\eqref{eq:xietaphitheta}: }=\frac{\xi(\phi,\theta)}{\eta(\phi,\theta)}}\overbracket{\longrightarrow}^{\qquad\qquad\mathllap{\substack{(r/b)\to0,\\\phi\to\overline{\phi}_0,\\\theta\to\overline{\theta}_0}}}r\sin\overline{\phi}_0+\frac{r\cos\phistar-r\cos{\overline{\phi}_0}}{-\sin{(\overline{\phi}_0+\overline{\theta}_0)}}\cos{(\overline{\phi}_0+\overline{\theta}_0)}\overbracket{=}^{\qquad\mathllap{\text{\eqref{eq:X1barY1barTheta1bar}}}}\overline{\mathcal{X}}_1,\\
    \mathcal{X}_1-\overline{\mathcal{X}}_1\overbracket{=}^{\mathllap{\text{\eqref{eq:X1barY1barTheta1bar}\eqref{eq:X1Y1Polynomial}}}}&r\sin\phi-r\sin{\overline{\phi}_0}+\frac{\xi(\phi,\theta)}{\eta(\phi,\theta)}-\frac{r\cos\phistar-r\cos{\overline{\phi}_0}}{-\sin{(\overline{\phi}_0+\overline{\theta}_0)}}\cos{(\overline{\phi}_0+\overline{\theta}_0)}.\\
    \end{aligned}
\end{equation} For the $\xi(\phi,\theta)$ limit approximation, we observe that \begin{equation}\label{eq:xiapproximation}
    \begin{aligned}
        &\bigm|\xi(\phi,\theta)-(2r\cos{\phistar}-2r\cos{\overline{\phi}_0})\cos{(\overline{\phi}_0+\overline{\theta}_0)}\bigm|\\
        \overbracket{=}^{\mathllap{\text{\eqref{eq:xietaphitheta}}}}&\bigm|\big(2r\cos{\phistar}-2r\cos{\phi}\big)\cos{(\phi+\theta)}-\big(2r\cos{\phistar}-2r\cos{\overline{\phi}_0}\big)\cos{(\overline{\phi}_0+\overline{\theta}_0)}\bigm|\\
        =&2r\bigm|\cos{\phistar}\big(\cos{(\phi+\theta)}-\cos{(\overline{\phi}_0+\overline{\theta}_0)}\big)+\cos{\overline{\phi}_0}\cos{(\overline{\phi}_0+\overline{\theta}_0)}-\cos{\phi}\cos{(\phi+\theta)}\bigm|\\
        =&2r\bigm|\cos{\phistar}(\cos{(\phi+\theta)}-\cos{(\overline{\phi}_0+\overline{\theta}_0)})+\cos{\overline{\phi}_0}\cos{(\overline{\phi}_0+\overline{\theta}_0)}-\cos{\overline{\phi}_0}\cos{(\phi+\theta)}\\
        &+\cos{\overline{\phi}_0}\cos{(\phi+\theta)}-\cos{\phi}\cos{(\phi+\theta)}\bigm|\\
        \le&2r\cos{\phistar}\bigm|\overbracket{\cos{(\phi+\theta)}-\cos{(\overline{\phi}_0+\overline{\theta}_0)}}^{\mathllap{\substack{\text{By mean value thm }\exists\text{ some }w_1\\  \cos{(\phi+\theta)}-\cos{(\overline{\phi}_0+\overline{\theta}_0)}=-(\phi+\theta-\overline{\phi}_0-\overline{\theta}_0)\sin{w_1}}}}\bigm|+2r\bigm|\cos{\overline{\phi}_0}\bigm|\bigm|\overbracket{\cos{(\phi+\theta)}-\cos{(\overline{\phi}_0+\overline{\theta}_0)}}^{\mathllap{\substack{\text{By mean value thm }\exists\text{ some }w_1\\  \cos{(\phi+\theta)}-\cos{(\overline{\phi}_0+\overline{\theta}_0)}=-(\phi+\theta-\overline{\phi}_0-\overline{\theta}_0)\sin{w_1}}}}\bigm|\\
&+2r\bigm|\cos{(\phi+\theta)}\bigm|\bigm|\overbracket{\cos{\overline{\phi}_0}-\cos\phi}^{{\mathllap{\substack{\text{By mean value thm }\exists\text{ some }w_2\\  \cos{\overline{\phi}_0}-\cos{\phi}=(\overline{\phi}_0-\phi)\sin{w_2}}}}}\bigm|\\
        \overbracket{\le}^{\mathllap{\text{By mean value thm}}}&2r\bigm|\phi+\theta-\overline{\phi}_0-\overline{\theta}_0\bigm|+2r\bigm|\phi+\theta-\overline{\phi}_0-\overline{\theta}_0\bigm|+2r\bigm|\overline{\phi}_0-\phi\bigm|\\
        \le&6r\bigm|\overline{\phi}_0-\phi\bigm|+4r\bigm|\overline{\theta}_0-\theta\bigm|.
    \end{aligned}
\end{equation} For the $\eta(\phi,\theta)$ limit approximation, we see that
\begin{equation}\label{eq:etaaproximation}
\begin{aligned}
    &\bigm|\eta(\phi,\theta)-(-2\sin{(\overline{\phi}_0+\overline{\theta}_0)})\bigm|\\\overbracket{=}^{\mathllap{\text{\eqref{eq:xietaphitheta}}}}&\bigm|\sqrt{\sin^2{(\theta+\phi)}+(r/b)^2\sin^2{\theta}+2(r/b)\sin{\theta}\sin{(\theta+\phi)}+2(r/b)\cos{\phistar}-2(r/b)\cos{\phi}}\\
    &-\sin{(\theta+\phi)}-(r/b)\sin{\theta}+2\sin{(\overline{\phi}_0+\overline{\theta}_0)}\bigm|\\
    \overbracket{\le}^{\mathllap{\substack{\text{Triangle inequality}}}}& \bigm|\sin{(\overline{\phi}_0+\overline{\theta}_0)}-\sin{(\theta+\phi)}\bigm|+\bigm|\frac{r\sin{\theta}}{b}\bigm|\\
    &+\Bigm|\sqrt{\sin^2{(\theta+\phi)}+\frac{r^2\sin^2{\theta}}{b^2}+\frac{2r\sin{\theta}\sin{(\theta+\phi)}}{b}+\frac{2r\cos{\phistar}-2r\cos{\phi}}{b}}+\sin{(\overline{\phi}_0+\overline{\theta}_0)}\Bigm|\\
    \overbracket{\le}^{\mathllap{\text{Mean value thm}}}&|\theta+\phi-\overline{\theta}_0-\overline{\phi}_0|+\frac{r}{b}\\
    &+\Bigm|\frac{(r^2/b^2)\sin^2{\theta}+(2r/b)\sin{\theta}\sin{(\theta+\phi)}+(2r/b)(\cos{\phistar}-\cos{\phi})+\sin^2{(\theta+\phi)}-\sin^2{(\overline{\phi}_0+\overline{\theta}_0)}}{\underbracket{\sqrt{\sin^2{(\theta+\phi)}+\frac{r^2\sin^2{\theta}}{b^2}+\frac{2r\sin{\theta}\sin{(\theta+\phi)}}{b}+\frac{2r\cos{\phistar}-2r\cos{\phi}}{b}}-\sin{(\overline{\phi}_0+\overline{\theta}_0)}}_{\ge-\sin{(\overline{\phi}_0+\overline{\theta}_0)}>0\text{, by }\pi<\overline{\phi}_0+\overline{\theta}_0<2\pi}}\Bigm|\\
    \le&|\theta+\phi-\overline{\theta}_0-\overline{\phi}_0|+\frac{r}{b}+\frac{(r^2/b^2+6r/b)}{-\sin{(\overline{\phi}_0+\overline{\theta}_0)}}+\frac{|\sin^2{(\phi+\theta)}-\sin^2{(\overline{\phi}_0+\overline{\theta}_0)}|}{-\sin{(\overline{\phi}_0+\overline{\theta}_0)}}\\
    \overbracket{\le}^{\mathllap{\text{Mean value thm}}}&|\theta+\phi-\overline{\theta}_0-\overline{\phi}_0|+\overbracket{\frac{r}{b}+\frac{(r^2/b^2+6r/b)}{-\sin{(\overline{\phi}_0+\overline{\theta}_0)}}}^{<(7.1r)/b}+\frac{2|\theta+\phi-\overline{\theta}_0-\overline{\phi}_0|}{-\sin{(\overline{\phi}_0+\overline{\theta}_0)}}\\
    \overbracket{<}^{\mathllap{\substack{b>R-r>1699r\\\text{and }\frac{1}{-\sin{(\overline{\phi}_0+\overline{\theta}_0)}}\ge1}}}&\frac{7.1r}{b}\frac{1}{-\sin{(\overline{\phi}_0+\overline{\theta}_0)}}+\frac{3|\theta-\overline{\theta}_0|}{-\sin{(\overline{\phi}_0+\overline{\theta}_0)}}+\frac{3|\phi-\overline{\phi}_0|}{-\sin{(\overline{\phi}_0+\overline{\theta}_0)}}.\\
\end{aligned}
\end{equation}\newline Hence, if  
\begin{equation}\label{eq:thetaphiclose}
    \begin{aligned}
    \bigm|\theta-\overline{\theta}_0\bigm|<&0.01\sin^2{(\overline{\phi}_0+\overline{\theta}_0)},\bigm|\phi-\overline{\phi}_0\bigm|<0.01\sin^2{(\overline{\phi}_0+\overline{\theta}_0)}, \\
    b>R-r>&\frac{500r}{\sin^2(\overline{\phi}_0+\overline{\theta}_0)},\\
    \textit{therefore }\frac{r}{b}<&\frac{1}{500}\sin^2{(\overline{\phi}_0+\overline{\theta}_0)},\\
    \end{aligned}
\end{equation}\newline then since we have the following estimate for $\eta(\phi,\theta)$,
\begin{align*}
    \bigm|\eta(\phi,\theta)-(-2\sin{(\overline{\phi}_0+\overline{\theta}_0)})\bigm|\overbracket{<}^{\mathllap{\eqref{eq:etaaproximation}}}&\frac{7.1r}{b}\frac{1}{-\sin{(\overline{\phi}_0+\overline{\theta}_0)}}+\frac{3|\theta-\overline{\theta}_0|}{-\sin{(\overline{\phi}_0+\overline{\theta}_0)}}+\frac{3|\phi-\overline{\phi}_0|}{-\sin{(\overline{\phi}_0+\overline{\theta}_0)}}\\
    \overbracket{<}^{\mathllap{\eqref{eq:thetaphiclose}}}&
    \frac{7.1}{500}\frac{\sin^2{(\overline{\phi}_0+\overline{\theta}_0)}}{-\sin{(\overline{\phi}_0+\overline{\theta}_0)}}+\frac{6\times0.01\sin^2{(\overline{\phi}_0+\overline{\theta}_0)}}{-\sin{(\overline{\phi}_0+\overline{\theta}_0)}}=-0.0742\sin{(\overline{\phi}_0+\overline{\theta}_0)}.\\
\end{align*}, we find a lower bound for $\eta(\phi,\theta)$ that is 
\begin{equation}\label{eq:etalowerbound}
    \begin{aligned}
    \eta(\phi,\theta)-(-2\sin{(\overline{\phi}_0+\overline{\theta}_0)})>&0.0742\sin{(\overline{\phi}_0+\overline{\theta}_0)}\\
    \eta(\phi,\theta)>&-1.9258\sin{(\overline{\phi}_0+\overline{\theta}_0)}>0.
    \end{aligned}
\end{equation}

For $\frac{\xi(\phi,\theta)}{\eta(\phi,\theta)}$, we see that
\begin{equation}\label{eq:xietaratioApproximate}
    \begin{aligned}
        &\Bigm|\frac{\xi(\phi,\theta)}{\eta(\phi,\theta)}-\frac{2r(\cos\phistar-r\cos{\overline{\phi}_0})\cos{(\overline{\phi}_0+\overline{\theta}_0)}}{-2\sin{(\overline{\phi}_0+\overline{\theta}_0)}}\Bigm|\\
        =&\Bigm|\frac{\xi(\phi,\theta)(-2\sin{(\overline{\phi}_0+\overline{\theta}_0)})-\eta(\phi,\theta)[2r(\cos{\phistar}-\cos{\overline{\phi}_0})\cos{(\overline{\phi}_0+\overline{\theta}_0)}]}{(-2\sin{(\overline{\phi}_0+\overline{\theta}_0)})\eta(\phi,\theta)}\Bigm|\\
        =&\Bigm|\frac{\xi(\phi,\theta)(-2\sin{(\overline{\phi}_0+\overline{\theta}_0)})-\xi(\phi,\theta)\eta(\phi,\theta)+\xi(\phi,\theta)\eta(\phi,\theta)-\eta(\phi,\theta)[2r(\cos{\phistar}-\cos{\overline{\phi}_0})\cos{(\overline{\phi}_0+\overline{\theta}_0)}]}
        {-2\sin{(\overline{\phi}_0+\overline{\theta}_0)}\eta(\phi,\theta)}\Bigm|\\
        \le&\frac{\overbracket{|\xi(\phi,\theta)|}^{\le4r}}{-2\sin{(\overline{\phi}_0+\overline{\theta}_0)}\underbracket{|\eta(\phi,\theta)|}_{\mathllap{\eqref{eq:etalowerbound}:>-1.9258\sin{(\overline{\phi}_0+\overline{\theta}_0)}>0}}}\underbracket{\Bigm|-2\sin{(\overline{\phi}_0+\overline{\theta}_0)}-\eta(\phi,\theta)\Bigm|}_{\text{LHS of \eqref{eq:etaaproximation}}}+\frac{\overbracket{\bigm|\xi(\phi,\theta)-(2r\cos{\phistar}-2r\cos{\overline{\phi}_0})\cos{(\overline{\phi}_0+\overline{\theta}_0)}\bigm|}^{\text{LHS of \eqref{eq:xiapproximation}}}}{-2\sin{(\overline{\phi}_0+\overline{\theta}_0)}}\\
        \overbracket{<}^{\mathllap{\text{\eqref{eq:xiapproximation}, \eqref{eq:etaaproximation} and \eqref{eq:etalowerbound}}}}&\frac{2r}{1.9258\sin^2{(\overline{\phi}_0+\overline{\theta}_0)}}\Big[\frac{7.1r}{b}\frac{1}{-\sin{(\overline{\phi}_0+\overline{\theta}_0)}}+\frac{3|\theta-\overline{\theta}_0|}{-\sin{(\overline{\phi}_0+\overline{\theta}_0)}}+\frac{3|\phi-\overline{\phi}_0|}{-\sin{(\overline{\phi}_0+\overline{\theta}_0)}}\Big]+\frac{6r\bigm|\overline{\phi}_0-\phi\bigm|+4r\bigm|\overline{\theta}_0-\theta\bigm|}{-2\sin{(\overline{\phi}_0+\overline{\theta}_0)}}\\
    \end{aligned}
\end{equation}

We summarize if
\begin{equation}\label{eq:X1closeCondition}
\begin{aligned}
    |\phi-\overline{\phi}_0|<&\min{\big\{\underbracket{0.01\sin^2{(\overline{\phi}_0+\overline{\theta}_0)}}_{\mathllap{\text{required in \eqref{eq:thetaphiclose}}}},\underbracket{-0.05\frac{\varepsilon}{r}\sin^3{(\overline{\phi}_0+\overline{\theta}_0)}}_{\text{used in \eqref{eq:X1X1BarClose}}}\big\}},\\ |\theta-\overline{\theta}_0|<&\min{\big\{\underbracket{0.01\sin^2{(\overline{\phi}_0+\overline{\theta}_0)}}_{\mathllap{\text{required in \eqref{eq:thetaphiclose}}}},\underbracket{-0.05\frac{\varepsilon}{r}\sin^3{(\overline{\phi}_0+\overline{\theta}_0)}}_{\text{used in \eqref{eq:X1X1BarClose}}}\big\}}\\
    \text{And }R>&\max{\big\{\underbracket{1700r}_{\mathllap{\text{required by \eqref{eq:etaaproximation}}}},\underbracket{r+\frac{500r}{\sin^2(\overline{\phi}_0+\overline{\theta}_0)}}_{\mathllap{\text{required in \eqref{eq:thetaphiclose}}}},\underbracket{r+\frac{23r^2}{-\varepsilon\sin^3{(\overline{\phi}_0+\overline{\theta}_0)}}}_{\text{used in \eqref{eq:X1X1BarClose}}}\big\}},
\end{aligned} 
\end{equation} then
\begin{equation}\label{eq:X1X1BarClose}
    \begin{aligned}
        \bigm|\mathcal{X}_1-\overline{\mathcal{X}}_1\bigm|\overbracket{=}^{\mathllap{\eqref{eq:X1X1bardiff}}}&\bigm|r\sin\phi-r\sin{\overline{\phi}_0}+\frac{\xi(\phi,\theta)}{\eta(\phi,\theta)}-\frac{r\cos\phistar-r\cos{\overline{\phi}_0}}{-\sin{(\overline{\phi}_0+\overline{\theta}_0)}}\cos{(\overline{\phi}_0+\overline{\theta}_0)}\bigm|\\
        \le&r\underbracket{\bigm|\sin\phi-\sin{\overline{\phi}_0}\bigm|}_{\mathllap{\text{ apply mean value thm: }\le|\phi-\overline{\phi}_0|}}+\underbracket{\bigm|\frac{\xi(\phi,\theta)}{\eta(\phi,\theta)}-\frac{r\cos\phistar-r\cos{\overline{\phi}_0}}{-\sin{(\overline{\phi}_0+\overline{\theta}_0)}}\cos{(\overline{\phi}_0+\overline{\theta}_0)}\bigm|}_{\mathllap{\text{\eqref{eq:xietaratioApproximate}}}}\\
        \overbracket{<}^{{\mathllap{\text{\eqref{eq:xietaratioApproximate}}}}}&r|\phi-\overline{\phi}_0|+\frac{2r}{1.9258\sin^2{(\overline{\phi}_0+\overline{\theta}_0)}}\Big[\frac{7.1r}{b}\frac{1}{-\sin{(\overline{\phi}_0+\overline{\theta}_0)}}+\frac{3|\theta-\overline{\theta}_0|}{-\sin{(\overline{\phi}_0+\overline{\theta}_0)}}+\frac{3|\phi-\overline{\phi}_0|}{-\sin{(\overline{\phi}_0+\overline{\theta}_0)}}\Big]\\
        &+\frac{6r\bigm|\overline{\phi}_0-\phi\bigm|+4r\bigm|\overline{\theta}_0-\theta\bigm|}{-2\sin{(\overline{\phi}_0+\overline{\theta}_0)}}\\
        =&\big(\frac{6r}{-1.9258\sin^3{(\overline{\phi}_0+\overline{\theta}_0)}}+\frac{3r}{-\sin{(\overline{\phi}_0+\overline{\theta}_0)}}+r\big)\big|\phi-\overline{\phi}_0\big|\\
        &+\big(\frac{6r}{-1.9258\sin^3{(\overline{\phi}_0+\overline{\theta}_0)}}+\frac{2r}{-\sin{(\overline{\phi}_0+\overline{\theta}_0)}}+r\big)\big|\theta-\overline{\theta}_0\big|+\frac{14.2r^2}{-1.9258b\sin^3{(\overline{\phi}_0+\overline{\theta}_0)}}\\
        \overbracket{<}^{\mathllap{\frac{1}{-\sin{(\overline{\phi}_0+\overline{\theta}_0)}}\ge1}}&\frac{7.2r}{-\sin^3{(\overline{\phi}_0+\overline{\theta}_0)}}\big|\phi-\overline{\phi}_0\big|+\frac{6.2r}{-\sin^3{(\overline{\phi}_0+\overline{\theta}_0)}}\big|\theta-\overline{\theta}_0\big|+\frac{7.4r^2}{\underbracket{-b\sin^3{(\overline{\phi}_0+\overline{\theta}_0)}}_{\mathllap{\text{\eqref{eq:X1closeCondition}:} b>R-r>\frac{23r^2}{-\varepsilon\sin^3{(\overline{\phi}_0+\overline{\theta}_0)}}}}}\\
        \overbracket{<}^{\mathllap{\eqref{eq:X1closeCondition}}}&\frac{7.2r}{-\sin^3{(\overline{\phi}_0+\overline{\theta}_0)}}\frac{-0.05\varepsilon}{r}\sin^3{(\overline{\phi}_0+\overline{\theta}_0)}+\frac{6.2r}{-\sin^3{(\overline{\phi}_0+\overline{\theta}_0)}}\frac{-0.05\varepsilon}{r}\sin^3{(\overline{\phi}_0+\overline{\theta}_0)}+\frac{7.6\varepsilon}{23}<\varepsilon\\
    \end{aligned}
\end{equation}

In table \cref{fig:Trajectory_LemonCompareWithWojtkowski} $\LEM(r,R,\phistar)$, we denote by $M$ the midpoint $\arc{AB}$. We have the following. If \begin{equation}\label{eq:Y1closecondition}
    R>\frac{2r^2}{\varepsilon}\sin^2\phistar
\end{equation} then
\begin{equation}
    \begin{aligned}
        \bigm|\mathcal{Y}_1-\overline{\mathcal{Y}}_1\bigm|=&\bigm|\mathcal{Y}_1-(-r\cos\phistar)\bigm|\le\text{distance}(M,|AB|)=2R\sin^2{(\frac{\Phistar}{2})}\\<&2R\sin^2{\Phistar}\overbracket{=}^{\qquad\mathllap{r\sin\phistar=R\sin\Phistar}}2r\sin{\phistar}\sin{\Phistar}\overbracket{=}^{\qquad\mathllap{r\sin\phistar=R\sin\Phistar}}\frac{2r^2}{R}\sin^2{\phistar}\overbracket{<}^{\mathllap{\text{\eqref{eq:Y1closecondition}}}}\varepsilon.
    \end{aligned}
\end{equation}\newline Since $\mathcal{F}(x)=(\Phi_1,\theta_1)$, $P=p(\mathcal{F}(x))\in\Gamma_R$, in \cref{fig:Trajectory_LemonCompareWithWojtkowski} the geometric meanings of $\phi$ and $\theta$ imply that $\overrightarrow{TP}$ has angle $2\pi-(\phi+\theta)\in(0,\pi)$ with positive axis $x$ given the clockwise orientation of the billiard plane. By the geometric meanings of $\Phi_1$ the tangential direction at $P$ of $\Gamma_R$ has angle $-\Phi_1$ with the positive axis $x$ given the clockwise orientation of the billiard plane. So, $2\pi+\Phi_1-(\phi+\theta)=\theta_1$ is the incident angle between $\overrightarrow{TP}$ and the tangential direction at $P$. By the geometric meaning of $\Theta_1$ in \cref{def:LemonAndPetTogether}, $\Theta_1-\Phi_1$ is the reflection angle $=\theta_1$. Hence we see \begin{equation}\label{eq:relectionAndincident}\begin{aligned}
    \Phi_1&\in\big(-\Phistar,\quad +\Phistar\big)\\
    \Theta_1-\Phi_1&=\theta_1=2\pi+\Phi_1-(\phi+\theta)\\
    \Theta_1&=2\pi+2\Phi_1-(\phi+\theta)
\end{aligned}\end{equation}\newline So, if \begin{equation}\label{eq:Theta1CloseCondition}
    \begin{aligned}
        R>&\frac{2.02r}{0.4\varepsilon_1}\\
        |\phi-\overline{\phi}_0|<&0.3\varepsilon_1\\
        |\theta-\overline{\theta}_0|<&0.3\varepsilon_1
    \end{aligned}
\end{equation}, then
\begin{equation}\label{eq:theta1theta1barApproximation}
    \begin{aligned}
        \overbracket{2\pi-\overline{\Theta}_1}^{\mathllap{\eqref{eq:X1barY1barTheta1bar}: =\overline{\phi}_0+\overline{\theta}_0}}+\overbracket{\Theta_1-\Phi_1}^{\mathrlap{\eqref{eq:relectionAndincident}: =2\pi+\Phi_1-(\phi+\theta)}}=&2\pi+\Phi_1-(\phi+\theta)+\overline{\phi}_0+\overline{\theta}_0\\
        \big|\Theta_1-\overline{\Theta}_1\big|=&\big|2\Phi_1-(\phi+\theta)+\overline{\phi}_0+\overline{\theta}_0\big|\\
        \le&2\big|\Phi_1\big|+\big|\phi-\overline{\phi}_0\big|+\big|\theta-\overline{\theta}_0\big|\\
        <&2|\Phistar|+\big|\phi-\overline{\phi}_0\big|+\big|\theta-\overline{\theta}_0\big|\\
        \overbracket{<}^{\mathllap{r\sin\phistar=R\sin\Phistar}}&2\sin^{-1}\big((r/R)\sin\phistar\big)+\big|\phi-\overline{\phi}_0\big|+\big|\theta-\overline{\theta}_0\big|\\
        \overbracket{<}^{\mathllap{\substack{R>1700r: 0<\sin^{-1}\big((r/R)\sin\phistar\big)<0.02\\ \text{If }0<z<0.2\text{, then } \mathrm{sinc}(z)>1/1.01}}}&2\cdot1.01\sin{\big(\sin^{-1}{\big((r/R)\sin\phistar\big)}\big)}+\big|\phi-\overline{\phi}_0\big|+\big|\theta-\overline{\theta}_0\big|\\
        <&2.02(r/R)+\big|\phi-\overline{\phi}_0\big|+\big|\theta-\overline{\theta}_0\big|\\
        \overbracket{<}^{\mathllap{\eqref{eq:Theta1CloseCondition}}}&2.02\frac{0.4\varepsilon_1}{2.02}+0.3\varepsilon_1+0.3\varepsilon_1=\varepsilon_1
    \end{aligned}
\end{equation}

Note that if 
\begin{equation}\label{eq:theta1a0condition1}
    \begin{aligned}
    \overline{\phi}_0-\phi>&-(2\pi-\overline{\phi}_0-\overline{\theta}_0)/4\\
    \overline{\theta}_0-\theta>&-(2\pi-\overline{\phi}_0-\overline{\theta}_0)/4\\
    R>&4r/\sin^2{(\frac{\pi}{2}-\frac{\overline{\phi}_0+\overline{\theta}_0}{4})},
    \end{aligned}
\end{equation}then 
\begin{equation}\label{eq:theta1a0bound1}
    \begin{aligned}
        \theta_1-\sin^{-1}{(\sqrt{4r/R})}\overbracket{=}^{\mathrlap{\text{\eqref{eq:relectionAndincident}}}}&2\pi+\Phi_1-(\phi+\theta)-\sin^{-1}{(\sqrt{4r/R})}\\
        =&2\pi-(\overline{\phi}_0+\overline{\theta}_0)+(\overline{\phi}_0-\phi)+(\overline{\theta}_0-\theta)+\overbracket{\Phi_1}^{\mathrlap{>-\Phistar=-\sin^{-1}{\big((r/R)\sin{\phistar}\big)}}}-\sin^{-1}{(\sqrt{4r/R})}\\
        >&2\pi-(\overline{\phi}_0+\overline{\theta}_0)+\underbracket{(\overline{\phi}_0-\phi)}_{\mathllap{\overline{\phi}_0-\phi>-\frac{2\pi-\overline{\phi}_0-\overline{\theta}_0}{4}}}+\underbracket{(\overline{\theta}_0-\theta)}_{\overline{\theta}_0-\theta>-\frac{2\pi-\overline{\phi}_0-\overline{\theta}_0}{4}}\underbracket{-\sin^{-1}{\big((r/R)\sin{\phistar}\big)}}_{>-\sin^{-1}{\sqrt{4r/R}}}-\sin^{-1}{(\sqrt{4r/R})}\\
        >&\pi-\frac{\overline{\phi}_0+\overline{\theta}_0}{2}-2\sin^{-1}(\sqrt{4r/R})\underbracket{>}_{\mathllap{R>4r/\sin^2{(\frac{\pi}{2}-\frac{\overline{\phi}_0+\overline{\theta}_0}{4})}}}0.
    \end{aligned}
\end{equation}
\newline If \begin{equation}\label{eq:theta1a0condition2}
    \begin{aligned}
        -\overline{\phi}_0+\phi>&-\frac{-\pi+\overline{\phi}_0+\overline{\theta}_0}{4},\\
        -\overline{\theta}_0+\theta>&-\frac{-\pi+\overline{\phi}_0+\overline{\theta}_0}{4}\\
        R>&\frac{4r}{\sin^2{(\frac{\overline{\phi}_0+\overline{\theta}_0}{4}-\frac{\pi}{4})}},
    \end{aligned}
\end{equation}\newline then
\begin{equation}\label{eq:theta1a0bound2}
    \begin{aligned}
        \pi-\sin^{-1}{(\sqrt{4r/R})}-\theta_1\overbracket{=}^{\mathrlap{\text{\eqref{eq:relectionAndincident}}}}&\pi-\big(2\pi+\Phi_1-(\phi+\theta)\big)-\sin^{-1}{(\sqrt{4r/R})}\\
        =&-\pi+(\overline{\phi}_0+\overline{\theta}_0)+(-\overline{\phi}_0+\phi)+(-\overline{\theta}_0+\theta)\overbracket{-\Phi_1}^{\mathrlap{>-\Phistar=-\sin^{-1}{\big((r/R)\sin{\phistar}\big)}}}-\sin^{-1}{(\sqrt{4r/R})}\\
        >&-\pi+(\overline{\phi}_0+\overline{\theta}_0)+\underbracket{(-\overline{\phi}_0+\phi)}_{\mathllap{-\overline{\phi}_0+\phi>-\frac{-\pi+\overline{\phi}_0+\overline{\theta}_0}{4}}}+\overbracket{(-\overline{\theta}_0+\theta)}^{\mathllap{-\overline{\theta}_0+\theta>-\frac{-\pi+\overline{\phi}_0+\overline{\theta}_0}{4}}}\underbracket{-\sin^{-1}{\big((r/R)\sin{\phistar}\big)}}_{\mathllap{>-\sin^{-1}{\sqrt{4r/R}}}}-\sin^{-1}{(\sqrt{4r/R})}\\
        >&-\frac{\pi}{2}+\frac{\overline{\phi}_0+\overline{\theta}_0}{2}-2\sin^{-1}(\sqrt{4r/R})\underbracket{>}_{\mathllap{R>4r/\sin^2{(\frac{\overline{\phi}_0+\overline{\theta}_0}{4}-\frac{\pi}{4})}}}0.
    \end{aligned}
\end{equation}\newline Therefore, conditions \eqref{eq:theta1a0condition1} and \eqref{eq:theta1a0condition2} imply \eqref{eq:theta1a0bound1} and  \eqref{eq:theta1a0bound2} thus $\sin{\theta_1}>\sqrt{4r/R}$.\newline Conditions \eqref{eq:etalowerbound}, \eqref{eq:X1closeCondition}, \eqref{eq:Y1closecondition}, 
\eqref{eq:Theta1CloseCondition}, \eqref{eq:theta1a0condition1} and \eqref{eq:theta1a0condition2} collected become condition \eqref{eq:x0barx0ClosingCondition}.
\end{proof}

\begin{figure}[h]
\begin{center}
\begin{tikzpicture}[xscale=0.6,yscale=0.6]
    \tkzDefPoint(0,0){Or};
    \tkzDefPoint(0,-4.8){Y};
    \pgfmathsetmacro{\rradius}{4.5};
    \clip(-1.1*\rradius,-1.1*\rradius) rectangle (2.8*\rradius, 0.02*\rradius);
    \pgfmathsetmacro{\phistardeg}{40};
    \pgfmathsetmacro{\XValueArc}{\rradius*sin(\phistardeg)};
    \pgfmathsetmacro{\YValueArc}{\rradius*cos(\phistardeg)};
    \pgfmathsetmacro{\Rradius}{
        \rradius*sin(\phistardeg)/sin(18)
    };
    \pgfmathsetmacro{\bdist}{(-\rradius*cos(\phistardeg))+cos(18)*\Rradius}
    \tkzDefPoint(0,\bdist){OR};
    \tkzDefPoint(\XValueArc,-1.0*\YValueArc){B};    \tkzDefPoint(-1.0*\XValueArc,-1.0*\YValueArc){A};
    \tkzDrawArc[name path=Cr, thick](Or,B)(A);
    \tkzDrawArc[name path=CR, ultra thin](OR,A)(B);
    \pgfmathsetmacro{\PXValue}{-\Rradius*sin(3)};
    \pgfmathsetmacro{\PYValue}{\bdist-\Rradius*cos(3)};
    \tkzDefPoint(\PXValue,\PYValue){P};
    



    \path[name path=Line_P_Q0_far](P)-- (3:8);
    \path[name path=Line_P_Q_far](P)-- (1:8);
    \path[name intersections={of=Cr and Line_P_Q0_far,by={Q0}}];
    \path[name intersections={of=Cr and Line_P_Q_far,by={Q}}];
    \tkzDefPoint(0,\bdist-\Rradius){C};
    \path[name intersections={of=Cr and Line_P_Q0_far,by={Q0}}];
    \draw[name path=LineAB,thin] (A)--(B);
    \coordinate (PHorizon) at ([shift={(3,0)}]P);
    \draw[red,ultra thin,dashed] (P)--(PHorizon);
    
    \path[name intersections={of=Line_P_Q0_far and LineAB,by={P_bar}}];
    

    \tkzLabelPoint[below](A){$A$};
    \tkzLabelPoint[below](P){\small $p(x_1)=P$};
    \tkzLabelPoint[below](Or){$O_r$};
    \tkzLabelPoint[below](B){$B$};
    \tkzLabelPoint[above left]([shift={(-0.01,0.02)}]P_bar){\small $p(\overline{x}_1)=\overline{P}$};
    \tkzLabelPoint[right]([shift={(0.42,0.3)}]P_bar){\small $\overline{\Theta}_1$};
    \tkzLabelPoint[right]([shift={(0.2,0.13)}]P){\small $\Theta_1$};
    
    \tkzDefPointBy[projection=onto B--A](T0)\tkzGetPoint{M0};
    \coordinate (T0Mirror) at ($(T0)!2!(M0)$);
    \coordinate (PbarArrow) at ($(T0Mirror)!1.06!(P_bar)$);
    \tkzDefPointBy[projection=onto P--OR](T)\tkzGetPoint{M};
    \coordinate (TMirror) at ($(T)!2!(M)$);
    \coordinate (PArrow) at ($(P)!0.07!(TMirror)$);

    

    \tkzMarkAngle[arc=lll,ultra thin,size=0.55](B,P_bar,Q0);
    \tkzMarkAngle[arc=lll,ultra thin,size=0.22](PHorizon,P,Q);
    \tkzDefPoint(0,\bdist-\Rradius){C};
    \begin{scope}[ultra thin,decoration={markings,mark=at position 0.7 with {\arrow{>}}}] 
    \draw[red,postaction={decorate},ultra thin] (P_bar)--(Q0);
    \draw[blue,postaction={decorate},ultra thin] (P)--(Q);
    \end{scope}
    \tkzLabelPoint[above right](Q0){$\overline{Q}=p(\overline{x}_2)=(r\sin{\overline{\phi}_2},-r\cos{\overline{\phi}_2})$};
    \tkzLabelPoint[below right](Q){$Q=p(x_2)=(r\sin{\phi_2},-r\cos{\phi_2})$};
    \tkzDrawPoints(A,B,P,P_bar,Q0,Q,Or);
\end{tikzpicture}
\caption{
$\overline{P}=(\overline{\mathcal{X}}_1,\overline{\mathcal{Y}}_1)=p(\overline{x}_1)$ and $\overline{x}_2=\overline{F}_{\phistar}(\overline{x}_1),\overline{Q}=p(\overline{x}_2)$.\newline $P=(\mathcal{X}_1,\mathcal{Y}_1)=p(x_1)$, and $x_2=\mathcal{F}(x_1),Q=p(x_2)$
} \label{fig:Trajectory_LemonCompareWithWojtkowski2}
\end{center}
\end{figure}\hfill
\begin{lemma}[Shadowing lemma 2, trajectory approximation in \cref{fig:Trajectory_LemonCompareWithWojtkowski2}]\label{lemma:Trajectory_approximation_2} Let $(\mathcal{X}_1,\mathcal{Y}_1,\Theta_1)$ be the position and direction of the $\LEM(r,R,\phistar)$ billiard flow after collision on the boundary $\Gamma_R$ (in \cref{def:flowcollsionpositionandangle}) of $x_1\in \MRout\subset M_R$ with $(\phi_2,\theta_2)=x_2=\Fcal(x_1)\in\Mrin$. In the coodinate system \cref{def:StandardCoordinateTable} for both petal and lemon billiard table in \cref{fig:1-petalExtensionFromCornerAndLemon}, let $P=p(x_1)=(\mathcal{X}_1,\mathcal{Y}_1)\in\Gamma_R$ and let $(\overline{\mathcal{X}}_1,\overline{\mathcal{Y}}_1,\overline{\Theta}_1)$ be the position and direction of the $\PET(r,\phistar)$ billiard flow after collision on the boundary $\overline{AB}$  (in \cref{def:flowcollsionpositionandangle}) of $\overline{x}_1\in M_{\text{\upshape{f}}}$. Suppose $\overline{P}=p(\overline{x}_1)=(\overline{\mathcal{X}}_1,\overline{\mathcal{Y}}_1)\in\Int(\overline{AB})$. i.e. $\big|\overline{\mathcal{X}}_1\big|<r\sin\phistar$ and $(\overline{\phi}_2,\overline{\theta}_2)=\overline{x}_2\nfd\overline{F}_{\phistar}(\overline{x}_1)\in\Mrin$. We have the following conclusion. For $\varepsilon>0$ if \begin{equation}\label{eq:X1Y1Theta1CloseCondition}
    \begin{aligned}
        \bigm|\mathcal{X}_1-\overline{\mathcal{X}}_1\bigm|<\delta_3(\varepsilon,\overline{\mathcal{X}}_1)\dfn\min&\Big\{\frac{0.14r}{\pi}\big(1+\frac{6r}{\sqrt{r^2\sin^2{\phistar}-\overline{\mathcal{X}}^2_1}}\big)^{-1}\varepsilon,\big(1+\frac{6r}{\sqrt{r^2\sin^2{\phistar}-\overline{\mathcal{X}}^2_1}}\big)^{-1}0.015r\sin{\phistar},\\
        &\frac{0.2r\varepsilon}{3\pi}\Big\}\\
        \bigm|\mathcal{Y}_1-\overbracket{\overline{\mathcal{Y}}_1}^{\mathrlap{=-r\cos\phistar}}\bigm|<\delta_4(\varepsilon,\overline{\mathcal{X}}_1)\dfn\min&\Big\{\frac{0.14r}{\pi}\big(1+\frac{6r}{\sqrt{r^2\sin^2{\phistar}-\overline{\mathcal{X}}^2_1}}\big)^{-1}\varepsilon, \big(1+\frac{6r}{\sqrt{r^2\sin^2{\phistar}-\overline{\mathcal{X}}^2_1}}\big)^{-1}0.015r\sin{\phistar},\\
        &\frac{0.2r\varepsilon}{3\pi}\Big\},\\
         \bigm|\Theta_1-\overline{\Theta}_1\bigm|<\delta_5(\varepsilon,\overline{\mathcal{X}}_1)\dfn\min&\Big\{\frac{0.28}{3\pi}\big(1+\frac{4r}{\sqrt{r^2\sin^2{\phistar}-\overline{\mathcal{X}}^2_1}}\big)^{-1}\varepsilon,\big(1+\frac{4r}{\sqrt{r^2\sin^2{\phistar}-\overline{\mathcal{X}}^2_1}}\big)^{-1}0.01\sin{\phistar},\\
         &\frac{0.1\varepsilon}{4\pi+1}\Big\},
    \end{aligned}
\end{equation} then \begin{equation}\label{eq:phi2theta2close}
    \begin{aligned}
        |\overline{\phi}_2-\phi_2|<&\varepsilon,\\
        |\overline{\theta}_2-\theta_2|<&\varepsilon.
    \end{aligned}
\end{equation}
\end{lemma}
\begin{proof}
In \cref{fig:Trajectory_LemonCompareWithWojtkowski2} with coordinate system of \cref{def:StandardCoordinateTable},\newline $O_r=(0,0)$, $O_R=(0,b)$, $Q=p(x_2)=(r\sin{\phi_2},-r\cos{\phi_2})=(x_Q,y_Q)$, $\overline{Q}=(r\sin{\overline{\phi}_2},-r\cos{\overline{\phi}_2})=(x_{\overline{Q}},y_{\overline{Q}})$.
    
    Denote $\overline{t}=|\overline{P}\,\overline{Q}|$ and $t=|PQ|$. Note that the vector $\overrightarrow{\overline{P}\,\overline{Q}}$ has an angle $\overline{\Theta}_1$ with respect to the positive axis $x$, the vector $\overrightarrow{PQ}$ has an angle $\Theta_1$ with respect to the positive axis $x$ and $Q$, $\overline{Q}$ are on the circle $C_r$. We get the following.
\begin{equation}\label{eq:tau1lengthQuadraticEq}
    \begin{aligned}
        &\left.
        \begin{aligned}
            x_Q=\mathcal{X}_1+t\cos{\Theta_1}&=r\sin\phi_2\\
            y_Q=\mathcal{Y}_1+t\sin{\Theta_1}&=-r\cos\phi_2\\
        \end{aligned}\right\}\Longrightarrow(\overbracket{\mathcal{X}_1+t\cos{\Theta_1}}^{\mathllap{=r\sin{\phi_2}}})^2+(\overbracket{\mathcal{Y}_1+t\sin{\Theta_1}}^{=-r\cos{\phi_2}})^2=r^2\\
        &\Longrightarrow t^2+(2\mathcal{X}_1\cos{\Theta_1}+2\mathcal{Y}_1\sin{\Theta_1})t+\mathcal{X}_1^2+\mathcal{Y}_1^2-r^2=0.\\
    \end{aligned}
\end{equation} 
Since $(\mathcal{X}_1,\mathcal{Y}_1)=P$ is on the circle $C_R$ centered at $(0,b)$ and $\measuredangle{O_RO_rB}=\pi-\phistar$, we can see the following. \begin{equation}\label{eq:tau1lengthQuadraticEq2}
    \begin{aligned}
        &\left.
        \begin{aligned}
            \mathcal{X}^2_1+\mathcal{Y}^2_1-2b\mathcal{Y}_1+b^2=\mathcal{X}_1^2+(\mathcal{Y}_1-b)^2=&R^2\\
            b^2+r^2-2br\cos{(\angle{O_RO_rB})}=b^2+r^2+2br\cos\phistar=&R^2\\
        \end{aligned}
        \right\}\Rightarrow \mathcal{X}^2_1+\mathcal{Y}^2_1-2b\mathcal{Y}_1+b^2=b^2+r^2+2br\cos\phistar\\
        \Longrightarrow&\mathcal{X}^2_1+\mathcal{Y}^2_1-r^2=2b\mathcal{Y}_1+2br\cos\phistar\\
        &\text{Therefore, the last quadratic equation in \eqref{eq:tau1lengthQuadraticEq} becomes}\\
        &t^2+(2\mathcal{X}_1\cos{\Theta_1}+2\mathcal{Y}_1\sin{\Theta_1})t+\underbracket{2b\mathcal{Y}_1+2br\cos\phistar}_{=\mathcal{X}_1^2+\mathcal{Y}_1^2-r^2}=0
    \end{aligned}
\end{equation}
 
Since $2b(\mathcal{Y}_1+r\cos\phistar)=\mathcal{X}_1^2+\mathcal{Y}_1^2-r^2<0$, the last quadratic equation in \eqref{eq:tau1lengthQuadraticEq2} has one positive root and one negative root. $|PQ|=t$ is its positive root. 
\begin{equation}\label{eq:tplus}
    \begin{aligned}
        |PQ|=t=&-\mathcal{X}_1\cos{\Theta_1}-\mathcal{Y}_1\sin{\Theta_1}+\sqrt{\Delta}\\
        \text{Where: }\Delta=&(\mathcal{X}_1\cos{\Theta_1}+\mathcal{Y}_1\sin{\Theta_1})^2-2b(\mathcal{Y}_1+r\cos{\phistar})>0\text{, since }\mathcal{Y}_1+r\cos{\phistar}<0\\
    \end{aligned}
\end{equation}
\begin{equation}\label{eq:tau1barlengthQuadraticEq}
    \begin{aligned}
        &\left.
        \begin{aligned}
            x_{\overline{Q}}=\overline{\mathcal{X}}_1+\overline{t}\cos{\overline{\Theta}_1}&=r\sin\overline{\phi}_2\\
            y_{\overline{Q}}=\overline{\mathcal{Y}}_1+\overline{t}\sin{\overline{\Theta}_1}&=-r\cos\overline{\phi}_2\\
            \overline{\mathcal{Y}}_1&=-r\cos{\phistar}\end{aligned}\right\}\Longrightarrow(\overbracket{\overline{\mathcal{X}}_1+\overline{t}\cos{\overline{\Theta}_1}}^{\mathllap{=r\sin{\overline{\phi}_2}}})^2+(\overbracket{-r\cos\phistar+\overline{t}\sin{\overline{\Theta}_1}}^{=-r\cos{\overline{\phi}_2}})^2=r^2\\
        \Longrightarrow& \overline{t}^2+(2\overline{\mathcal{X}}_1\cos{\overline{\Theta}_1}-2r\cos{\phistar}\sin{\overline{\Theta}_1})\overline{t}+\overline{\mathcal{X}}_1^2-r^2\sin^2{\phistar}=0\\
    \end{aligned}
\end{equation}
Since $\overline{\mathcal{X}}_1^2-r^2\sin^2{\phistar}=\overline{\mathcal{X}}^2_1+\overline{\mathcal{Y}}^2_1-r^2<0$, the last quadratic equation \eqref{eq:tau1barlengthQuadraticEq} has one positive and one negative root. $|\overline{P}\,\overline{Q}|=\overline{t}$ is its positive root. So \begin{equation}\label{eq:tau1barlengthQuadraticEq2}
    \begin{aligned}
    &|\overline{P}\,\overline{Q}|=\overline{t}=-\overline{\mathcal{X}}_1\cos{\overline{\Theta}_1}+r\cos\phistar\sin{\overline{\Theta}_1}+\sqrt{\overline{\Delta}},\\
        &\text{where }\overline{\Delta}=(\overline{\mathcal{X}}_1\cos{\overline{\Theta}_1}-r\cos\phistar\sin{\overline{\Theta}_1})^2-(\overline{\mathcal{X}}^2_1-r^2\sin^2{\phistar})>0,\text{ since }\overline{\mathcal{X}}_1\in(-r\sin{\phistar},\quad r\sin{\phistar})\\
    \end{aligned}
\end{equation}\newline In order to compare $t$ with $\overline{t}$ and to compare $\sqrt{\Delta}$ with $\sqrt{\overline{\Delta}}$, we will derive/use the following seven inequalities. \begin{equation}\label{eq:DeltaCompareIneq1}
    \begin{aligned}
        |2b(\mathcal{Y}_1+r\cos{\phistar})-(\mathcal{X}^2_1-r^2\sin^2{\phistar})|=&|2b\mathcal{Y}_1-\mathcal{X}^2_1+2br\cos{\phistar}+r^2\sin^2{\phistar}|\\
        \overbracket{=}^{\mathllap{\eqref{eq:tau1lengthQuadraticEq2}:\mathcal{X}^2_1+\mathcal{Y}^2_1-2b\mathcal{Y}_1+b^2=R^2}}&|\mathcal{Y}^2_1+b^2-R^2+2br\cos{\phistar}+r^2\sin^2{\phistar}|\\
        \overbracket{=}^{\mathllap{b^2+r^2+2br\cos\phistar=R^2}}&|\mathcal{Y}^2_1-r^2+r^2\sin^2{\phistar}|=|\mathcal{Y}^2_1-r^2\cos^2{\phistar}|\\
        =&|\mathcal{Y}_1-r\cos{\phistar}||\mathcal{Y}_1+r\cos{\phistar}|\overbracket{<}^{\mathllap{|\mathcal{Y}_1-r\cos{\phistar}|<|\mathcal{Y}_1|+r<2r}}2r|\mathcal{Y}_1+r\cos{\phistar}|\\
    \end{aligned}
\end{equation}
\begin{equation}\label{eq:DeltaCompareIneq2}
    \begin{aligned}
        &|\mathcal{X}^2_1-r^2\sin^2{\phistar}-(\overline{\mathcal{X}}^2_1-r^2\sin^2{\phistar})|=|\mathcal{X}^2_1-\overline{\mathcal{X}}^2_1|=|\overline{\mathcal{X}}_1-\mathcal{X}_1|\overbracket{|\overline{\mathcal{X}}_1+\mathcal{X}_1|}^{\mathllap{|\overline{\mathcal{X}}_1+\mathcal{X}_1|<|\overline{\mathcal{X}}_1|+|\mathcal{X}_1|<2r}}<2r|\overline{\mathcal{X}}_1-\mathcal{X}_1|\\
    \end{aligned}
\end{equation}
\begin{equation}\label{eq:DeltaCompareIneq3}
    \begin{aligned}
        &\bigm|2b(\mathcal{Y}_1+r\cos{\phistar})-(\overline{\mathcal{X}}^2_1-r^2\sin^2{\phistar})\bigm|\\
        =&\bigm|2b(\mathcal{Y}_1+r\cos{\phistar})-(\mathcal{X}^2_1-r^2\sin^2{\phistar})+\mathcal{X}^2_1-r^2\sin^2{\phistar}-(\overline{\mathcal{X}}^2_1-r^2\sin^2{\phistar})\bigm|\\
        \le&\bigm|2b(\mathcal{Y}_1+r\cos{\phistar})-(\mathcal{X}^2_1-r^2\sin^2{\phistar})\bigm|+\bigm|\mathcal{X}^2_1-r^2\sin^2{\phistar}-(\overline{\mathcal{X}}^2_1-r^2\sin^2{\phistar})\bigm|\\
        \overbracket{<}^{\mathllap{\eqref{eq:DeltaCompareIneq1},\eqref{eq:DeltaCompareIneq2}}}&2r\bigm|\mathcal{Y}_1+r\cos{\phistar}\bigm|+2r\bigm|\overline{\mathcal{X}}_1-\mathcal{X}_1\bigm|\\
    \end{aligned}
\end{equation}
\begin{equation}\label{eq:DeltaCompareIneq4}
    \begin{aligned}
    &\bigm|\overline{\mathcal{X}}_1\cos{\overline{\Theta}_1}-r\cos{\phistar}\sin{\overline{\Theta}_1}-(\mathcal{X}_1\cos{\Theta_1}+\mathcal{Y}_1\sin{\Theta_1})\bigm|\\
    =&\bigm|\overline{\mathcal{X}}_1\cos{\overline{\Theta}_1}-\mathcal{X}_1\cos{\overline{\Theta}_1}+\mathcal{X}_1\cos{\overline{\Theta}_1}-\mathcal{X}_1\cos{\Theta_1}-r\cos{\phistar}\sin{\overline{\Theta}_1}-\mathcal{Y}_1\sin{\overline{\Theta}_1}+\mathcal{Y}_1\sin{\overline{\Theta}_1}-\mathcal{Y}_1\sin{\Theta_1}\bigm|\\
    \le& |\overline{\mathcal{X}}_1\cos{\overline{\Theta}_1}-\mathcal{X}_1\cos{\overline{\Theta}_1}|+\underbracket{|\mathcal{X}_1\cos{\overline{\Theta}_1}-\mathcal{X}_1\cos{\Theta_1}|}_{\qquad\qquad\mathllap{\text{By Mean-Value Theorem: }\le|\mathcal{X}_1||\overline{\Theta}_1-\Theta_1|< r|\overline{\Theta}_1-\Theta_1|}}+|-r\cos{\phistar}\sin{\overline{\Theta}_1}-\mathcal{Y}_1\sin{\overline{\Theta}_1}|+\underbracket{|\mathcal{Y}_1\sin{\overline{\Theta}_1}-\mathcal{Y}_1\sin{\Theta_1}|}_{\qquad\qquad\mathllap{\text{By Mean-Value Theorem: }\le|\mathcal{Y}_1||\overline{\Theta}_1-\Theta_1|<r|\overline{\Theta}_1-\Theta_1|}}\\
    \underbracket{<}_{\mathrlap{\substack{\text{By Mean-Value Theorem and }\\ |\mathcal{X}_1|<r, |\mathcal{Y}_1|<r}}}&\bigm|\overline{\mathcal{X}}_1-\mathcal{X}_1\bigm|+\bigm|r\cos{\phistar}+\mathcal{Y}_1\bigm|+2r\bigm|\Theta_1-\overline{\Theta}_1\bigm|
    \end{aligned}
\end{equation}

\begin{equation}\label{eq:DeltaCompareIneq5}
    \begin{aligned}
        &\bigm|(\overline{\mathcal{X}}_1\cos{\overline{\Theta}_1}-r\cos{\phistar}\sin{\overline{\Theta}_1})^2-(\mathcal{X}_1\cos{\Theta_1}+\mathcal{Y}_1\sin{\Theta_1})^2\bigm|\\
        =&\underbracket{\bigm|\overline{\mathcal{X}}_1\cos{\overline{\Theta}_1}-r\cos{\phistar}\sin{\overline{\Theta}_1}+\mathcal{X}_1\cos{\Theta_1}+\mathcal{Y}_1\sin{\Theta_1}\bigm|}_{<4r\text{ since } |\overline{\mathcal{X}}_1|<r,|\mathcal{X}_1|<r,|\mathcal{Y}_1|<r}\bigm|\overline{\mathcal{X}}_1\cos{\overline{\Theta}_1}-r\cos{\phistar}\sin{\overline{\Theta}_1}-(\mathcal{X}_1\cos{\Theta_1}+\mathcal{Y}_1\sin{\Theta_1})\bigm|\\<&4r\bigm|\overline{\mathcal{X}}_1\cos{\overline{\Theta}_1}-r\cos{\phistar}\sin{\overline{\Theta}_1}-(\mathcal{X}_1\cos{\Theta_1}+\mathcal{Y}_1\sin{\Theta_1})\bigm|\\
        \underbracket{<}_{\mathrlap{\eqref{eq:DeltaCompareIneq4}}}&4r\Big(\bigm|\overline{\mathcal{X}}_1-\mathcal{X}_1\bigm|+\bigm|r\cos{\phistar}+\mathcal{Y}_1\bigm|+2r\bigm|\Theta_1-\overline{\Theta}_1\bigm|\Big)
    \end{aligned}
\end{equation}

Since \eqref{eq:tau1barlengthQuadraticEq2}: $\overline{\Delta}=(\overline{\mathcal{X}}_1\cos{\overline{\Theta}_1}-r\cos\phistar\sin{\overline{\Theta}_1})^2-(\overline{\mathcal{X}}^2_1-r^2\sin^2{\phistar})\ge r^2\sin^2{\phistar}-\overline{\mathcal{X}}^2_1>0$, we can compare/estimate the difference of $\sqrt{\Delta}$ and $\sqrt{\overline{\Delta}}$ in the following.
\begin{equation}\label{eq:DeltaCompareIneq6}
\begin{aligned}
    &\bigm|\sqrt{\Delta}-\sqrt{\overline{\Delta}}\bigm|=\frac{\bigm|\Delta-\overline{\Delta}\bigm|}{\sqrt{\Delta}+\sqrt{\overline{\Delta}}}<\frac{1}{\underbracket{\sqrt{\overline{\Delta}}}_{\mathrlap{\ge\sqrt{r^2\sin^2{\phistar}-\overline{\mathcal{X}}^2_1}}}}\bigm|\Delta-\overline{\Delta}\bigm|\\
    \le&\frac{1}{\sqrt{r^2\sin^2{\phistar}-\overline{\mathcal{X}}^2_1}}\bigm|\Delta-\overline{\Delta}\bigm|\\
    \underbracket{=}_{\mathrlap{\text{\eqref{eq:tplus}, \eqref{eq:tau1barlengthQuadraticEq2}}}}&\frac{\bigm|(\mathcal{X}_1\cos{\Theta_1}+\mathcal{Y}_1\sin{\Theta_1})^2-2b(\mathcal{Y}_1+r\cos{\phistar})-\big[(\overline{\mathcal{X}}_1\cos{\overline{\Theta}_1}-r\cos\phistar\sin{\overline{\Theta}_1})^2-(\overline{\mathcal{X}}^2_1-r^2\sin^2{\phistar})\big]\bigm|}{\sqrt{r^2\sin^2{\phistar}-\overline{\mathcal{X}}^2_1}}\\
    \le&\frac{\bigm|(\mathcal{X}_1\cos{\Theta_1}+\mathcal{Y}_1\sin{\Theta_1})^2-(\overline{\mathcal{X}}_1\cos{\overline{\Theta}_1}-r\cos\phistar\sin{\overline{\Theta}_1})^2\bigm|}{\sqrt{r^2\sin^2{\phistar}-\overline{\mathcal{X}}^2_1}}+\frac{\bigm|-2b(\mathcal{Y}_1+r\cos{\phistar})+(\overline{\mathcal{X}}^2_1-r^2\sin^2{\phistar})\bigm|}{\sqrt{r^2\sin^2{\phistar}-\overline{\mathcal{X}}^2_1}}\\
    \overbracket{<}^{\mathllap{\eqref{eq:DeltaCompareIneq5},\eqref{eq:DeltaCompareIneq3}}}&\frac{4r\Big(\bigm|\overline{\mathcal{X}}_1-\mathcal{X}_1\bigm|+\bigm|r\cos{\phistar}+\mathcal{Y}_1\bigm|+2r\bigm|\Theta_1-\overline{\Theta}_1\bigm|\Big)+2r\bigm|\mathcal{Y}_1+r\cos{\phistar}\bigm|+2r\bigm|\overline{\mathcal{X}}_1-\mathcal{X}_1\bigm|}{\sqrt{r^2\sin^2{\phistar}-\overline{\mathcal{X}}^2_1}}\\
    =&\frac{2r\big(3\bigm|\overline{\mathcal{X}}_1-\mathcal{X}_1\bigm|+3\bigm|\mathcal{Y}_1+r\cos{\phistar}\bigm|+4r\bigm|\Theta_1-\overline{\Theta}_1\bigm|\big)}{\sqrt{r^2\sin^2{\phistar}-\overline{\mathcal{X}}^2_1}}.
\end{aligned}
\end{equation} And we can estimate the difference between $t$ and $\overline{t}$ by the following.\begin{equation}
    \begin{aligned}
        \bigm|t-\overline{t}\bigm|\overbracket{=}^{\mathllap{\text{\eqref{eq:tplus}, \eqref{eq:tau1barlengthQuadraticEq2}}}}&\bigm|-\mathcal{X}_1\cos{\Theta_1}-\mathcal{Y}_1\sin{\Theta_1}+\sqrt{\Delta}-(-\overline{\mathcal{X}}_1\cos{\overline{\Theta}_1}+r\cos\phistar\sin{\overline{\Theta}_1}+\sqrt{\overline{\Delta}})\bigm|\\
        \le&\underbracket{\bigm|-\mathcal{X}_1\cos{\Theta_1}-\mathcal{Y}_1\sin{\Theta_1}-(-\overline{\mathcal{X}}_1\cos{\overline{\Theta}_1}+r\cos\phistar\sin{\overline{\Theta}_1})\bigm|}_{\qquad\qquad\qquad\mathllap{\text{\eqref{eq:DeltaCompareIneq4}: }<|\overline{\mathcal{X}}_1-\mathcal{X}_1|+|r\cos{\phistar}+\mathcal{Y}_1|+2r|\Theta_1-\overline{\Theta}_1|}}+\underbracket{\bigm|\sqrt{\Delta}-\sqrt{\overline{\Delta}}\bigm|}_{\qquad\qquad\qquad\qquad\qquad\qquad\qquad\mathllap{\text{\eqref{eq:DeltaCompareIneq6}: }<\frac{2r(3|\overline{\mathcal{X}}_1-\mathcal{X}_1|+3|\mathcal{Y}_1+r\cos{\phistar}|+4r|\Theta_1-\overline{\Theta}_1|)}{\sqrt{r^2\sin^2{\phistar}-\overline{\mathcal{X}}^2_1}}}}\\
        \overbracket{<}^{\mathllap{\text{\eqref{eq:DeltaCompareIneq4},\eqref{eq:DeltaCompareIneq6}}}}&\big(1+\frac{6r}{\sqrt{r^2\sin^2{\phistar}-\overline{\mathcal{X}}^2_1}}\big)\big(\bigm|\overline{\mathcal{X}}_1-\mathcal{X}_1\bigm|+|\mathcal{Y}_1+\underbracket{r\cos{\phistar}}_{=-\mathcal{Y}_1}|\big)+\big(2r+\frac{8r^2}{\sqrt{r^2\sin^2{\phistar}-\overline{\mathcal{X}}^2_1}}\big)\bigm|\Theta_1-\overline{\Theta}_1\bigm|.
    \end{aligned}
\end{equation}\newline For $\forall \varepsilon>0$, if \begin{equation}\label{eq:ttbarvarepsilonCondition}
    \begin{aligned}
        \bigm|\mathcal{X}_1-\overline{\mathcal{X}}_1\bigm|<&\frac{1.4r}{3\pi}\cdot\big(1+\frac{6r}{\sqrt{r^2\sin^2{\phistar}-\overline{\mathcal{X}}^2_1}}\big)^{-1}\cdot0.3\varepsilon\\
        \bigm|\mathcal{Y}_1-\overline{\mathcal{Y}}_1\bigm|<&\frac{1.4r}{3\pi}\cdot\big(1+\frac{6r}{\sqrt{r^2\sin^2{\phistar}-\overline{\mathcal{X}}^2_1}}\big)^{-1}\cdot0.3\varepsilon\\
        \bigm|\Theta_1-\overline{\Theta}_1\bigm|<&\frac{1.4r}{3\pi}\cdot\big(2r+\frac{8r^2}{\sqrt{r^2\sin^2{\phistar}-\overline{\mathcal{X}}^2_1}}\big)^{-1}\cdot0.4\varepsilon,\\
        \text{then    }\bigm|t-\overline{t}\bigm|<&\frac{1.4r}{3\pi}\varepsilon.
    \end{aligned}
\end{equation} Now we compare $\phi_2$ with $\overline{\phi}_2$ and compare $\theta_2$ with $\overline{\theta}_2$.
\begin{equation}\label{eq:sindiffphi2phibar2}
    \begin{aligned}
        r^2\bigm|\sin{(\phi_2-\overline{\phi}_2)}\bigm|=&r^2\bigm|\sin{(\overline{\phi}_2-\phi_2)}\bigm|=\bigm|-\underbracket{(r\sin{\phi_2})}_{\qquad\mathllap{\text{\eqref{eq:tau1lengthQuadraticEq}: }=\mathcal{X}_1+t\cos{\Theta_1}}}\overbracket{(r\cos{\overline{\phi}_2})}^{\mathllap{\text{\eqref{eq:tau1barlengthQuadraticEq}: }=r\cos{\phistar}-\overline{t}\sin{\overline{\Theta}_1}}}+\underbracket{(r\cos{\phi_2})}_{\qquad\mathllap{\text{\eqref{eq:tau1lengthQuadraticEq}: }=-\mathcal{Y}_1-t\sin{\Theta_1}}}\overbracket{(r\sin{\overline{\phi}_2})}^{\mathrlap{\text{\eqref{eq:tau1barlengthQuadraticEq}: }\overline{\mathcal{X}}_1+\overline{t}\cos{\overline{\Theta}_1}}\qquad}\bigm|\\
        =&\bigm|-
        (\mathcal{X}_1+t\cos{\Theta_1})(r\cos{\phistar}-\overline{t}\sin{\overline{\Theta}_1})-(\mathcal{Y}_1+t\sin{\Theta_1})(\overline{\mathcal{X}}_1+\overline{t}\cos{\overline{\Theta}_1})\bigm|\\
        =&\big|-\mathcal{X}_1r\cos{\phistar}-tr\cos{\phistar}\cos{\Theta_1}+\mathcal{X}_1\overline{t}\sin{\overline{\Theta}_1}+\cos{\Theta_1}\sin{\overline{\Theta}_1}t\overline{t}\\
        &-\overline{\mathcal{X}}_1\mathcal{Y}_1-\overline{\mathcal{X}}_1\sin{\Theta_1}t-\mathcal{Y}_1\cos{\overline{\Theta}_1}\overline{t}-\sin{\Theta_1}\cos{\overline{\Theta}_1}t\overline{t}\bigm|\\
        \le&\bigm|\underbracket{-\mathcal{X}_1r\cos{\phistar}-\overline{\mathcal{X}}_1\mathcal{Y}_1}_{\qquad\qquad\mathllap{=-\mathcal{X}_1r\cos{\phistar}-\mathcal{X}_1\mathcal{Y}_1+\mathcal{X}_1\mathcal{Y}_1-\overline{\mathcal{X}}_1\mathcal{Y}_1}}\bigm|+\bigm|\underbracket{\mathcal{X}_1\overline{t}\sin{\overline{\Theta}_1}-\overline{\mathcal{X}}_1\sin{\Theta_1}t}_{\mathrlap{=\mathcal{X}_1\overline{t}\sin{\overline{\Theta}_1}-\overline{\mathcal{X}}_1\overline{t}\sin{\overline{\Theta}_1}+\overline{\mathcal{X}}_1\overline{t}\sin{\overline{\Theta}_1}-\overline{\mathcal{X}}_1\sin{\Theta_1}t}}\bigm|\\
        &+\bigm|-tr\cos{\phistar}\cos{\Theta_1}-\mathcal{Y}_1\cos{\overline{\Theta}_1}\overline{t}\bigm|+\bigm|\cos{\Theta_1}\sin{\overline{\Theta}_1}-\sin{\Theta_1}\cos{\overline{\Theta}_1}\bigm|t\overline{t}\\
        \le&\underbracket{\bigm|-\mathcal{X}_1r\cos{\phistar}-\mathcal{X}_1\mathcal{Y}_1\bigm|}_{< r|\mathcal{Y}_1+r\cos\phistar|}+\underbracket{\bigm|\mathcal{X}_1\mathcal{Y}_1-\overline{\mathcal{X}}_1\mathcal{Y}_1\bigm|}_{<r|\mathcal{X}_1-\overline{\mathcal{X}}_1|}+\underbracket{\bigm|\mathcal{X}_1\overline{t}\sin{\overline{\Theta}_1}-\overline{\mathcal{X}}_1\overline{t}\sin{\overline{\Theta}_1}\bigm|}_{\text{By }\overline{t}\le2r\text{, }<2r|\mathcal{X}_1-\overline{\mathcal{X}}_1|}\\
        &+\underbracket{\bigm|\overline{\mathcal{X}}_1\overline{t}\sin{\overline{\Theta}_1}-\overline{\mathcal{X}}_1\sin{\Theta_1}t\bigm|}_{\mathllap{\substack{=|\overline{\mathcal{X}}_1\overline{t}\sin{\overline{\Theta}_1}-\overline{\mathcal{X}}_1t\sin{\overline{\Theta}_1}+\overline{\mathcal{X}}_1t\sin{\overline{\Theta}_1}-\overline{\mathcal{X}}_1\sin{\Theta_1}t|\\ <r|t-\overline{t}|+2r^2|\Theta_1-\overline{\Theta}_1|}}}+\underbracket{\bigm|-tr\cos{\phistar}\cos{\Theta_1}+\overline{t}r\cos{\phistar}\cos{\Theta_1}\bigm|}_{\mathrlap{r|t-\overline{t}|}}\\
        &+\underbracket{\bigm|-tr\cos{\phistar}\cos{\Theta_1}-\mathcal{Y}_1\cos{\overline{\Theta}_1}\overline{t}\bigm|}_{\qquad\qquad\mathllap{\substack{=|-tr\cos{\phistar}\cos{\Theta_1}+tr\cos{\phistar}\cos{\overline{\Theta}_1}-tr\cos{\phistar}\cos{\overline{\Theta}_1}\\+\overline{t}r\cos{\phistar}\cos{\overline{\Theta}_1}-\overline{t}r\cos{\phistar}\cos{\overline{\Theta}_1}-\mathcal{Y}_1\cos{\overline{\Theta}_1}\overline{t}|}}}+\underbracket{t\overline{t}\bigm|\sin{(\Theta_1-\overline{\Theta}_1)}\bigm|}_{\mathrlap{\substack{\text{By Mean-Value Theorem and }t\le2r,\overline{t}\le2r: \\\le4r^2|\Theta_1-\overline{\Theta}_1|}}}\\
        <&3r\bigm|\mathcal{X}_1-\overline{\mathcal{X}}_1\bigm|+r\bigm|\mathcal{Y}_1+r\cos\phistar\bigm|+6r^2\bigm|\Theta_1-\overline{\Theta}_1\bigm|+2r\bigm|t-\overline{t}\bigm|\\
        &+\underbracket{\bigm|-tr\cos{\phistar}\cos{\Theta_1}+tr\cos{\phistar}\cos{\overline{\Theta}_1}\bigm|}_{\text{By} |t|\le2r \text{ and Mean-Value Theorem }<2r^2|\Theta_1-\overline{\Theta}_1|}+\underbracket{\bigm|-tr\cos{\phistar}\cos{\overline{\Theta}_1}+\overline{t}r\cos{\phistar}\cos{\overline{\Theta}_1}\bigm|}_{<r|t-\overline{t}|}\\
        &+\underbracket{\bigm|-\overline{t}r\cos{\phistar}\cos{\overline{\Theta}_1}-\mathcal{Y}_1\cos{\overline{\Theta}_1}\overline{t}\bigm|}_{<2r|\mathcal{Y}_1+r\cos\phistar|}\\
        <&3r\bigm|\mathcal{X}_1-\overline{\mathcal{X}}_1\bigm|+3r\bigm|\mathcal{Y}_1+r\cos\phistar\bigm|+8r^2\bigm|\Theta_1-\overline{\Theta}_1\bigm|+3r\bigm|t-\overline{t}\bigm|\\
    \end{aligned}
\end{equation}

\begin{equation}\label{eq:cosdiffphi2phibar2_1}
    \begin{aligned}
        r^2\cos{(\overline{\phi}_2-\phi_2)}=&\overbracket{(\mathcal{Y}_1+t\sin{\Theta_1})}^{\mathllap{\text{\eqref{eq:tau1lengthQuadraticEq}: }=-r\cos{\phi_2}}}\overbracket{(\overline{\mathcal{Y}}_1+\overline{t}\sin{\overline{\Theta}_1})}^{\text{\eqref{eq:tau1barlengthQuadraticEq}: }=-r\cos{\overline{\phi}_2}}+\overbracket{(\mathcal{X}_1+t\cos{\Theta_1})}^{\text{\eqref{eq:tau1lengthQuadraticEq}: }=r\sin{\phi_2}}\overbracket{(\overline{\mathcal{X}}_1+\overline{t}\cos{\overline{\Theta}_1})}^{\text{\eqref{eq:tau1barlengthQuadraticEq}: }=r\sin{\phi_2}}\\
        =&\overbracket{(\overline{\mathcal{Y}}_1+\overline{t}\sin{\overline{\Theta}_1})^2}^{=r^2\cos^2{\overline{\phi}_2}}+(\mathcal{Y}_1+t\sin{\Theta_1}-\overline{\mathcal{Y}}_1-\overline{t}\sin{\overline{\Theta}_1})(\overline{\mathcal{Y}}_1+\overline{t}\sin{\overline{\Theta}_1})\\
        &+\overbracket{(\overline{\mathcal{X}}_1+\overline{t}\cos{\overline{\Theta}_1})^2}^{=r^2\sin^2{\overline{\phi}_2}}+(\mathcal{X}_1+t\cos{\Theta_1}-\overline{\mathcal{X}}_1-\overline{t}\cos{\overline{\Theta}_1})(\overline{\mathcal{X}}_1+\overline{t}\cos{\overline{\Theta}_1})\\
        =&r^2-(-\mathcal{Y}_1-t\sin{\Theta_1}+\overline{\mathcal{Y}}_1+\overline{t}\sin{\overline{\Theta}_1})(\overline{\mathcal{Y}}_1+\overline{t}\sin{\overline{\Theta}_1})\\
        &-(-\mathcal{X}_1-t\cos{\Theta_1}+\overline{\mathcal{X}}_1+\overline{t}\cos{\overline{\Theta}_1})(\overline{\mathcal{X}}_1+\overline{t}\cos{\overline{\Theta}_1}),\\
    \end{aligned}
\end{equation}
\begin{equation}\label{eq:cosdiffphi2phibar2_2}
    \begin{aligned}
        (-\mathcal{Y}_1-t\sin{\Theta_1}+\overline{\mathcal{Y}}_1+\overline{t}\sin{\overline{\Theta}_1})(\overline{\mathcal{Y}}_1+\overline{t}\sin{\overline{\Theta}_1})\le&\bigm|-\mathcal{Y}_1-t\sin{\Theta_1}+\overline{\mathcal{Y}}_1+\overline{t}\sin{\overline{\Theta}_1}\bigm|\overbracket{\bigm|\overline{\mathcal{Y}}_1+\overline{t}\sin{\overline{\Theta}_1}\bigm|}^{=r|\cos{\overline{\phi}_1}|\le r}\\
        \le&r\bigm|\overline{\mathcal{Y}}_1-\mathcal{Y}_1-t\sin{\Theta_1}+\overline{t}\sin{\Theta_1}-\overline{t}\sin{\Theta_1}+\overline{t}\sin{\overline{\Theta}_1}\bigm|\\
        \overbracket{\le}^{\mathllap{\text{By Mean-Value Theorem and }t,\overline{t}\le2r}}&r\bigm|\overline{\mathcal{Y}}_1-\mathcal{Y}_1\bigm|+r\bigm|t-\overline{t}\bigm|+2r^2\bigm|\Theta_1-\overline{\Theta}_1\bigm|,
    \end{aligned}
\end{equation}
\begin{equation}\label{eq:cosdiffphi2phibar2_3}
    \begin{aligned}
        (-\mathcal{X}_1-t\cos{\Theta_1}+\overline{\mathcal{X}}_1+\overline{t}\cos{\overline{\Theta}_1})(\overline{\mathcal{X}}_1+\overline{t}\cos{\overline{\Theta}_1})\le&\bigm|-\mathcal{X}_1-t\cos{\Theta_1}+\overline{\mathcal{X}}_1+\overline{t}\cos{\overline{\Theta}_1}\bigm|\overbracket{\bigm|\overline{\mathcal{X}}_1+\overline{t}\cos{\overline{\Theta}_1}\bigm|}^{=r|\sin{\overline{\phi}_1}|\le r}\\
        \le&r\bigm|\overline{\mathcal{X}}_1-\mathcal{X}_1-t\cos{\Theta_1}+\overline{t}\cos{\Theta_1}-\overline{t}\cos{\Theta_1}+\overline{t}\cos{\overline{\Theta}_1}\bigm|\\
        \overbracket{\le}^{\mathllap{\text{By Mean-Value Theorem and }t,\overline{t}\le2r}}&r\bigm|\overline{\mathcal{X}}_1-\mathcal{X}_1\bigm|+r\bigm|\overline{t}-t\bigm|+2r^2\bigm|\Theta_1-\overline{\Theta}_1\bigm|
    \end{aligned}
\end{equation}
By \eqref{eq:cosdiffphi2phibar2_1}, \eqref{eq:cosdiffphi2phibar2_2}, \eqref{eq:cosdiffphi2phibar2_3}, we see that if \begin{equation}\label{eq:X1Y1Theta1closeCond1}
\begin{aligned}
    \bigm|\overline{\mathcal{X}}_1-\mathcal{X}_1\bigm|&<0.1r\\
    \bigm|\overline{\mathcal{Y}}_1-\mathcal{Y}_1\bigm|&<0.1r\\
    \bigm|t-\overline{t}\bigm|&<0.1r\\
    \bigm|\Theta_1-\overline{\Theta}_1\bigm|&<0.1
\end{aligned}
\end{equation} then \begin{align*}
    r^2-r^2\cos{(\phi_2-\overline{\phi}_2)}\overbracket{=}^{\mathllap{\eqref{eq:cosdiffphi2phibar2_1}}}&(-\mathcal{Y}_1-t\sin{\Theta_1}+\overline{\mathcal{Y}}_1+\overline{t}\sin{\overline{\Theta}_1})(\overline{\mathcal{Y}}_1+\overline{t}\sin{\overline{\Theta}_1})\\
    &+(-\mathcal{X}_1-t\cos{\Theta_1}+\overline{\mathcal{X}}_1+\overline{t}\cos{\overline{\Theta}_1})(\overline{\mathcal{X}}_1+\overline{t}\cos{\overline{\Theta}_1})\\
    \overbracket{<}^{\mathllap{\eqref{eq:cosdiffphi2phibar2_2}, \eqref{eq:cosdiffphi2phibar2_3}}}&r\bigm|\overline{\mathcal{X}}_1-\mathcal{X}_1\bigm|+r\bigm|\overline{\mathcal{Y}}_1-\mathcal{Y}_1\bigm|+2r\bigm|\overline{t}-t\bigm|+4r^2|\Theta_1-\overline{\Theta}_1|\\
    <&r^2,
\end{align*} which implies $\cos{(\phi_2-\overline{\phi}_2)}>0$. Since $\phi_2,\overline{\phi}_2\in\big[\phistar,2\pi-\phistar\big]$, so 
\begin{equation}\label{eq:phi2phibar2diffrange}
\phi_2-\overline{\phi}_2\in[-2\pi+2\phistar,\,-\frac{3\pi}{2})\cup(-\frac{\pi}{2},\,\frac{\pi}{2})\cup(\frac{3\pi}{2},\,2\pi-2\phistar].\end{equation}\newline  With $\phistar<\pi/4$ and by \eqref{eq:sindiffphi2phibar2} if \begin{equation}\label{eq:X1Y1Theta1closeCond2}
    \begin{aligned}
        \bigm|\mathcal{X}_1-\overline{\mathcal{X}}_1\bigm|&<0.05r\sin{(2\phistar)}\\
        \bigm|\mathcal{Y}_1-\overbracket{\overline{\mathcal{Y}}_1}^{\mathrlap{=-r\cos\phistar}}\bigm|&<0.05r\sin{(2\phistar)}\\
        \bigm|\Theta_1-\overline{\Theta}_1\bigm|&<0.05\sin{(2\phistar)}\\
        \bigm|t-\overline{t}\bigm|&<0.05r\sin{(2\phistar)}
    \end{aligned}
\end{equation} then \begin{align*}
    &r^2\bigm|\sin{(\phi_2-\overline{\phi}_2)}\bigm|\overbracket{<}^{\mathllap{\eqref{eq:sindiffphi2phibar2}}}3r\bigm|\mathcal{X}_1-\overline{\mathcal{X}}_1\bigm|+3r\bigm|\mathcal{Y}_1+r\cos\phistar\bigm|+8r^2\bigm|\Theta_1-\overline{\Theta}_1\bigm|+3r\bigm|t-\overline{t}\bigm|<r^2\sin(2\phistar)\\
    \Rightarrow& \bigm|\sin{(\phi_2-\overline{\phi}_2)}\bigm|<\sin{(2\phistar)}\overbracket{\Rightarrow}^{2\phistar<\pi/2} \phi_2-\overline{\phi}_2\notin[-2\pi+2\phistar,\,-\frac{3\pi}{2})\cup(\frac{3\pi}{2},\,2\pi-2\phistar]\overbracket{\Rightarrow}^{\eqref{eq:phi2phibar2diffrange}} (\phi_2-\overline{\phi}_2)\in(-\pi/2,\pi/2).
\end{align*} Therefore, since if $|z|<\pi/2$, then $|\mathrm{sinc}(z)|=\frac{|\sin{z}|}{|z|}>\frac{2}{\pi}$. we get the following estimate for $|\phi_2-\overline{\phi}_2|$. \begin{equation}\label{eq:phi2phibar2diff}
    \begin{aligned}
        \bigm|\phi_2-\overline{\phi}_2\bigm|&<\frac{\pi}{2}\bigm|\sin{(\phi_2-\overline{\phi}_2)}\bigm|\overbracket{<}^{\eqref{eq:sindiffphi2phibar2}}\frac{\pi}{2}\big(\frac{3}{r}\bigm|\mathcal{X}_1-\overline{\mathcal{X}}_1\bigm|+\frac{3}{r}\bigm|\mathcal{Y}_1-\overline{\mathcal{Y}}_1\bigm|+8\bigm|\Theta_1-\overline{\Theta}_1\bigm|+\frac{3}{r}\bigm|t-\overline{t}\bigm|\big).\\
    \end{aligned} 
\end{equation} 
To satisfy $\bigm|t-\overline{t}\bigm|<0.05r\sin{(2\phistar)}$ in \eqref{eq:X1Y1Theta1closeCond2},  
by \eqref{eq:ttbarvarepsilonCondition}, it suffices to make $\bigm|\mathcal{X}_1-\overline{\mathcal{X}}_1\bigm|$, $\bigm|\mathcal{Y}_1-\overline{\mathcal{Y}}_1\bigm|$ and $\bigm|\Theta_1-\overline{\Theta}_1\bigm|$ satisfy
\begin{equation}\label{eq:X1Y1Theta1closeCond3}
    \begin{aligned}
        \bigm|\mathcal{X}_1-\overline{\mathcal{X}}_1\bigm|<&\big(1+\frac{6r}{\sqrt{r^2\sin^2{\phistar}-\overline{\mathcal{X}}^2_1}}\big)^{-1}0.015r\sin{\phistar}\\
        \bigm|\mathcal{Y}_1-\overline{\mathcal{Y}}_1\bigm|<&\big(1+\frac{6r}{\sqrt{r^2\sin^2{\phistar}-\overline{\mathcal{X}}^2_1}}\big)^{-1}0.015r\sin{\phistar}\\
        \bigm|\Theta_1-\overline{\Theta}_1\bigm|<&\big(2r+\frac{8r^2}{\sqrt{r^2\sin^2{\phistar}-\overline{\mathcal{X}}^2_1}}\big)^{-1}\cdot0.02r\sin{\phistar}=\big(1+\frac{4r}{\sqrt{r^2\sin^2{\phistar}-\overline{\mathcal{X}}^2_1}}\big)^{-1}0.01\sin{\phistar}
    \end{aligned}
\end{equation}

In \cref{fig:Trajectory_LemonCompareWithWojtkowski2} the tangential direction at $\overline{P}$ of $\Gamma_r$ has angle $\overline{\phi}_2$ and $\overrightarrow{\overline{P}\,\overline{Q}}$ has angle $\overline{\Theta}_1$, the collision angle $\overline{\theta}_2=\overline{\phi}_2-\overline{\Theta}_1$. The tangential direction at $P$ of $\Gamma_r$ has angle $\phi_2$ and $\overrightarrow{PQ}$ has angle $\Theta_1$, the collision angle $\theta_2=\phi_2-\Theta_1$. we get the following estimate for $|\theta_2-\overline{\theta}_2|$. \begin{equation}\label{eq:theta2thetabar2diff}
    \begin{aligned}
    \bigm|\theta_2-\overline{\theta}_2\bigm|&=\bigm|\phi_2-\Theta_1-(\overline{\phi}_2-\overline{\Theta}_1)\bigm|=\bigm|\phi_2-\overline{\phi}_2+\overline{\Theta}_1-\Theta_1\bigm|\\
    &\le\bigm|\phi_2-\overline{\phi}_2\bigm|+\bigm|\overline{\Theta}_1-\Theta_1\bigm|\\
    &\overbracket{\le}^{\mathllap{\eqref{eq:phi2phibar2diff}}}\frac{3\pi}{2r}\bigm|\mathcal{X}_1-\overline{\mathcal{X}}_1\bigm|+\frac{3\pi}{2r}\bigm|\mathcal{Y}_1-\overline{\mathcal{Y}}_1\bigm|+(4\pi+1)\bigm|\Theta_1-\overline{\Theta}_1\bigm|+\frac{3\pi}{2r}\bigm|t-\overline{t}\bigm|
    \end{aligned}
    \end{equation}
For $\forall \varepsilon>0$ if $\bigm|\mathcal{X}_1-\overline{\mathcal{X}}_1\bigm|$, $\bigm|\mathcal{Y}_1-\overline{\mathcal{Y}}_1\bigm|$ and $\bigm|\Theta_1-\overline{\Theta}_1\bigm|$ satisfy \eqref{eq:ttbarvarepsilonCondition} and the following
\begin{equation}\label{eq:X1X1barY1Y1barTheta1Thetabar1Varepsilon}
\begin{aligned}
    \bigm|\mathcal{X}_1-\overline{\mathcal{X}}_1\bigm|&<\frac{2r}{3\pi}\cdot0.1\varepsilon\\
    \bigm|\mathcal{Y}_1-\overline{\mathcal{Y}}_1\bigm|&<\frac{2r}{3\pi}\cdot0.1\varepsilon\\
    \bigm|\Theta_1-\overline{\Theta}_1\bigm|&<\frac{1}{4\pi+1}\cdot0.1\varepsilon,
\end{aligned}
\end{equation}then \begin{align*}
    \max{\big\{\bigm|\phi_2-\overline{\phi}_2\bigm|,\bigm|\theta_2-\overline{\theta}_2\bigm|\big\}}\overbracket{\le}^{\mathllap{\eqref{eq:theta2thetabar2diff},\eqref{eq:phi2phibar2diff}}}&\frac{3\pi}{2r}\bigm|\mathcal{X}_1-\overline{\mathcal{X}}_1\bigm|+\frac{3\pi}{2r}\bigm|\mathcal{Y}_1-\overline{\mathcal{Y}}_1\bigm|+(4\pi+1)\bigm|\Theta_1-\overline{\Theta}_1\bigm|+\underbracket{\frac{3\pi}{2r}\bigm|t-\overline{t}\bigm|}_{\mathllap{\eqref{eq:ttbarvarepsilonCondition}: <\frac{3\pi}{2r}\frac{1.4r}{3\pi}\varepsilon=0.7\varepsilon}}\\
    <&0.1\varepsilon+0.1\varepsilon+0.1\varepsilon+0.7\varepsilon=\varepsilon.
\end{align*} So we conclude the collected conditions from \eqref{eq:ttbarvarepsilonCondition}, \eqref{eq:X1X1barY1Y1barTheta1Thetabar1Varepsilon}, \eqref{eq:X1Y1Theta1closeCond3}, \eqref{eq:X1Y1Theta1closeCond2} and \eqref{eq:X1Y1Theta1closeCond1} to be condition \eqref{eq:X1Y1Theta1CloseCondition}.\qedhere
\end{proof}
\subsection{Closeness of expansion in Lemon billiard and 1-petal billiard}
\begin{lemma}[Expansion of $\PET(r,\phistar)$ and $\LEM(r,R,\phistar)$ trajectories]\label{lemma:expansionbetweenPETandLEM} With notations from \cref{def:LemonAndPetTogether}\eqref{item05LemonAndPetTogether}, $\overline{x}_0=(\overline{\phi}_0,\overline{\theta}_0)\in\Mrout$ with $\overline{x}_1=F_{\phi_*}(\overline{x}_0)\in M_{\text{\upshape f}}$, $\overline{d}_0= r\sin{\overline{\theta}_0}$, $\overline{\tau}_0=|p(\overline{x}_0)p(\overline{x}_1)|$, $\overline{d}_2= r\sin{\overline{\theta}_2}$, $\overline{\tau}_1=|p(\overline{x}_1)p(\overline{x}_2)|$, then in the $\mathbb{P}(r,\phi_*)$ billiard table coordinate \cref{def:StandardCoordinateTable} in \cref{fig:1-petalExtensionFromCornerAndLemon} let $p(\overline{x}_1)=(\overline{\mathcal{X}}_1,\overline{\mathcal{Y}}_1)$ and suppose that $p(\overline{x}_1)=(\overline{\mathcal{X}}_1,\overline{\mathcal{Y}}_1)$ and $(\overline{\mathcal{X}}_1,\overline{\mathcal{Y}}_1,\overline{\Theta}_1)$ (in \cref{def:flowcollsionpositionandangle}) be the position and direction of the $\PET(r,\phistar)$ billiard flow after collision at $p(\overline{x}_1)$.

Also, with the notation of \cref{def:fpclf}: $x=(\phi,\theta)\in\Mrout$, $(\Phi_1,\theta_1)=x_1=\Fcal(x)\in M_R$ and $d_0= r\sin{\theta}$, $d_1=R\sin{\theta_1}$, $\tau_0=|p(x_1)p(x)|$,then in the $\mathbb{L}(r,R,\phi_*)$ billiard table of the same coordinate \cref{def:StandardCoordinateTable} in \cref{fig:1-petalExtensionFromCornerAndLemon} let $p(x_1)=(\mathcal{X}_1,\mathcal{Y}_1)$ and suppose $(\mathcal{X}_1,\mathcal{Y}_1,\Theta_1)$ to be the position and direction of the $\LEM(r,R,\phistar)$ (in \cref{def:flowcollsionpositionandangle}) billiard flow after collision at $p(x_1)$. We have the following.

$\overline{x}_0\in\Mrout$ implies $p(\overline{x}_1)\in\Int\overline{AB}$ in \cref{fig:1-petalExtensionFromCornerAndLemon}, that is, $|\overline{\mathcal{X}}_1|<r\sin\phistar$. For $\forall\varepsilon_2>0$ if \begin{equation}\label{eq:ExpansionClosingCondition_i}
        \begin{aligned}
            \bigm|\theta-\overline{\theta}_0\bigm|<&\delta_6(\overline{\phi}_0,\overline{\theta}_0,\varepsilon_2,\overline{\mathcal{X}}_1)\dfn\min\Bigm\{0.5\sin{\overline{\theta}_0},\,\frac{1}{32}\varepsilon_2\sin^2{\overline{\theta}_0},\,\frac{-\varepsilon_2\sin^3{(\overline{\phi}_0+\overline{\theta}_0)}\sin{\overline{\theta}_0}}{64},\\
            &\mathtt{\delta}_1\big(\overline{\phi}_0,\,\overline{\theta}_0,\,[1+\frac{6r}{\sqrt{r^2\sin^2{\phistar}-\overline{\mathcal{X}}^2_1}}]^{-1}0.075\varepsilon_2r\sin{\overline{\theta}_0},\,[2r+\frac{8r^2}{\sqrt{r^2\sin^2{\phistar}-\overline{\mathcal{X}}^2_1}}]^{-1}0.1\varepsilon_2r\sin{\overline{\theta}_0}\big)\Bigm\},\\
            \bigm|\phi-\overline{\phi}_0\bigm|<&\delta_7(\overline{\phi}_0,\overline{\theta}_0,\varepsilon_2,\overline{\mathcal{X}}_1)\dfn\min\Bigm\{\frac{-\varepsilon_2\sin^3{(\overline{\phi}_0+\overline{\theta}_0)}\sin{\overline{\theta}_0}}{64},\,\\
            &\,\mathtt{\delta}_2\big(\overline{\phi}_0,\,\overline{\theta}_0,\,[1+\frac{6r}{\sqrt{r^2\sin^2{\phistar}-\overline{\mathcal{X}}^2_1}}]^{-1}0.075\varepsilon_2 r\sin{\overline{\theta}_0},[2r+\frac{8r^2}{\sqrt{r^2\sin^2{\phistar}-\overline{\mathcal{X}}^2_1}}]^{-1}0.1\varepsilon_2r\sin{\overline{\theta}_0}\Big)\Bigm\},\\
            R>&\mathtt{R}_2(\overline{\phi}_0,\overline{\theta}_0,\varepsilon_2,\overline{\mathcal{X}}_1)\dfn\max\Bigm\{\,\frac{1}{64\varepsilon^2_2}r,\,r+\frac{64r}{-\varepsilon_2\sin^3{(\overline{\phi}_0+\overline{\theta}_0)}\sin{\overline{\theta}_0}},\\
            &\mathtt{R}_1\Big(\overline{\phi}_0,\overline{\theta}_0,\,[1+\frac{6r}{\sqrt{r^2\sin^2{\phistar}-\overline{\mathcal{X}}^2_1}}]^{-1}0.075\varepsilon_2r\sin{\overline{\theta}_0},[2r+\frac{8r^2}{\sqrt{r^2\sin^2{\phistar}-\overline{\mathcal{X}}^2_1}}]^{-1}0.1\varepsilon_2r\sin{\overline{\theta}_0}\Big)\Bigm\},
        \end{aligned}
    \end{equation} where $\delta_1$, $\mathtt{R}_1$ are defined in \eqref{eq:x0barx0ClosingCondition}, then $R$, $(\overline{\phi}_0,\overline{\theta}_0)$ and $(\phi,\theta)$ satisfy condition \eqref{eq:xorbitA0Condition}, i.e. $\sin{\theta_1}>\sqrt{4r/R}$. \newline Hence the $\LEM(r,R,\phistar)$ billiard return orbit on $\hat{M}$ in \cref{def:MhatReturnOrbitSegment} with $(\phi_0,\theta_0)=x_0=x=(\phi,\theta)$, $x_1=\Fcal(x)=\Fcal(x_0)$ is in the case (a0) in \eqref{eqMainCases} and $\mathcal{I}$ from \cref{def:AformularForLowerBoundOfExpansion} as function of $x_1$ thus also of $x$ satisfies 
    \begin{equation}
    \begin{aligned}
        \Bigm|\mathcal{I}-\big(-1-\frac{\overline{\tau}_0+\overline{\tau}_1-2\overline{d}_0}{\overline{d}_0}\big)\Big|<\varepsilon_2
    \end{aligned}
    \end{equation}

\end{lemma}
\begin{proof}
    Conditions in \eqref{eq:ExpansionClosingCondition_i} imply $\overline{\phi}_0$, $\overline{\theta}_0$, $R$ satisfy condition \eqref{eq:xorbitA0Condition} thus $d_1=R\sin{\theta_1}>\sqrt{4rR}$.
        \begin{equation}\label{eq:LemonexpansionIcompareWithPET}
        \begin{aligned}
            \Big|\mathcal{I}-\big(-1-\frac{\overline{\tau}_0+\overline{\tau}_1-2\overline{d}_0}{\overline{d}_0}\big)\Big|&=\Big|\overbracket{-1+\frac{\tau_1}{d_0}\Big[\frac{2(\tau_0-d_0)}{d_1}-\frac{\tau_0+\tau_1-2d_0}{\tau_1}\Big]}^{{\text{\cref{def:AformularForLowerBoundOfExpansion}}}\nfd\mathcal{I}}-\big(-1-\frac{\overline{\tau}_0+\overline{\tau}_1-2\overline{d}_0}{\overline{d}_0}\big)\Big|\\
            &=\Big|\frac{\overline{\tau}_0+\overline{\tau}_1}{\overline{d}_0}-\frac{\tau_0+\tau_1}{d_0}+\frac{\tau_1}{d_0}\frac{2(\tau_0-d_0)}{d_1}\Big|\\
            &\le\bigm|\frac{\overline{\tau}_0+\overline{\tau}_1}{\overline{d}_0}-\frac{\tau_0+\tau_1}{d_0}\bigm|+\underbracket{\bigm|\frac{\tau_1}{d_1}\frac{2(\tau_0-d_0)}{d_0}\bigm|}_{\substack{<2\sqrt{r/R}\\\text{since }\tau_1\le2r,d_1=R\sin\theta_1>\sqrt{4rR},\\\text{\cref{remark:deffpclf}:}-d_0<\tau_0-d_0<d_0}}\\
            &<\frac{1}{d_0\overline{d}_0}\bigm|\overline{d}_0(\tau_0+\tau_1)-d_0(\overline{\tau}_0+\overline{\tau}_1)\bigm|+2\sqrt{r/R}\\
            &=\frac{1}{d_0\overline{d}_0}\bigm|\overline{d}_0(\tau_0+\tau_1)-d_0(\tau_0+\tau_1)+d_0(\tau_0+\tau_1)-d_0(\overline{\tau}_0+\overline{\tau}_1)\bigm|+2\sqrt{r/R}\\
            &\le\frac{1}{d_0\overline{d}_0}\bigm|\tau_0+\tau_1\bigm|\bigm|d_0-\overline{d}_0\bigm|+\frac{1}{\overline{d}_0}\bigm|\tau_0-\overline{\tau}_0\bigm|+\frac{1}{\overline{d}_0}\bigm|\tau_1-\overline{\tau}_1\bigm|+2\sqrt{r/R}
        \end{aligned}                    
        \end{equation} Note that the same as in \cref{lemma:Trajectory_approximation_1}, $(\phi_0,\theta_0)$ satisfies $\phi_0+\theta_0\in(\pi,2\pi)$.
    
    For the first term in the rightmost of \eqref{eq:LemonexpansionIcompareWithPET}, if $x=(\phi,\theta)$ satisfies the following \begin{equation}\label{eq:FirstTermInLemonexpansionIcompareWithPET}
        \begin{aligned}
            \bigm|\theta-\overline{\theta}_0\bigm|<&\min\Bigm\{0.5\sin{\overline{\theta}_0},\,\frac{1}{32}\varepsilon_2\sin^2{\overline{\theta}_0}\Bigm\},
        \end{aligned}
    \end{equation} then
    \begin{equation}\label{eq:1sttermexpansioncomparison}
    \begin{aligned}
        \bigm|d_0-\overline{d}_0\bigm|=&\bigm|r\sin{\theta}-r\sin{\overline{\theta}_0}\bigm|\overbracket{\le}^{\mathllap{\text{Mean-Value Theorem}}}r\bigm|\theta-\overline{\theta}_0\bigm|\overbracket{<}^{\mathllap{\eqref{eq:FirstTermInLemonexpansionIcompareWithPET}}}0.5r\sin{\overline{\theta}_0}\Longrightarrow r\sin{\theta}>0.5r\sin{\overline{\theta}_0}\\
        \Longrightarrow\frac{\bigm|\tau_0+\tau_1\bigm|\bigm|d_0-\overline{d}_0\bigm|}{d_0\overline{d}_0}=&\frac{\tau_0+\tau_1}{d_0\overline{d}_0}\bigm|d_0-\overline{d}_0\bigm|\le\frac{4r}{\underbracket{d_0}_{\mathllap{>0.5r\sin{\overline{\theta}_0}}}\underbracket{\overline{d}_0}_{\mathrlap{=r\sin{\overline{\theta}_0}}}}\bigm|d_0-\overline{d}_0\bigm|<\frac{8r}{r^2\sin^2{\overline{\theta}_0}}\bigm|d_0-\overline{d}_0\bigm|\\
        \overbracket{<}^{\mathllap{\text{Mean-Value Theorem}}}&\frac{8r^2}{r^2\sin^2{\overline{\theta}_0}}\bigm|\theta-\overline{\theta}_0\bigm|\overbracket{<}^{\mathllap{\eqref{eq:FirstTermInLemonexpansionIcompareWithPET}}}\frac{\varepsilon_2}{4}\\
    \end{aligned}
    \end{equation}

    For the second term in the rightmost of \eqref{eq:LemonexpansionIcompareWithPET}, if the following holds  
    \begin{equation}\label{eq:secondtermtau0taubar0closecondition}
        \begin{aligned}
            \bigm|\theta-\overline{\theta}_0\bigm|<&\min\big\{0.01\sin^2{(\overline{\phi}_0+\overline{\theta}_0)},\,\frac{-\varepsilon_2\sin^3{(\overline{\phi}_0+\overline{\theta}_0)}\sin{\overline{\theta}_0}}{64}\big\}\\
            \bigm|\phi-\overline{\phi}_0\bigm|<&\min\big\{0.01\sin^2{(\overline{\phi}_0+\overline{\theta}_0)},\,\frac{-\varepsilon_2\sin^3{(\overline{\phi}_0+\overline{\theta}_0)}\sin{\overline{\theta}_0}}{64}\big\}\\
            R>&\max\big\{r+\frac{500r}{\sin^2(\overline{\phi}_0+\overline{\theta}_0)},1700r,r+\frac{64r}{-\varepsilon_2\sin^3{(\overline{\phi}_0+\overline{\theta}_0)}\sin{\overline{\theta}_0}}\big\},\\
    \end{aligned}
    \end{equation} then \eqref{eq:thetaphiclose} and thus \eqref{eq:etalowerbound},\eqref{eq:etaaproximation} hold. Therefore from the results in \eqref{eq:PTposRoot} and from \cref{fig:Trajectory_LemonCompareWithWojtkowski}, we get the following. \begin{align*}
        \tau_0=&t_+=\frac{2r\cos\phistar-2r\cos\phi}{-\sin{(\phi+\theta)}-(r/b)\sin\theta+\sqrt{\big(\sin{(\phi+\theta)}+(r/b)\sin\theta\big)^2+\big(2(r/b)\cos\phistar-2(r/b)\cos\phi\big)}}\nfd\frac{\varsigma(\phi,\theta)}{\eta(\phi,\theta)},\\
        \overline{\tau}_0=&\frac{r\cos\phistar-r\cos{\overline{\phi}_0}}{-\sin{(\overline{\phi}_0+\overline{\theta}_0)}}=\frac{2r\cos\phistar-2r\cos{\overline{\phi}_0}}{-2\sin{(\overline{\phi}_0+\overline{\theta}_0)}},
    \end{align*}\newline where $\varsigma(\phi,\theta)=2r\cos\phistar-2r\cos\phi$, and from \eqref{eq:xietaphitheta}\[\eta(\phi,\theta)=-\sin{(\phi+\theta)}-(r/b)\sin\theta+\sqrt{\big(\sin{(\phi+\theta)}+(r/b)\sin\theta\big)^2+\big(2(r/b)\cos\phistar-2(r/b)\cos\phi\big)}.\] And by Mean-Value Theorem we have
    \begin{equation}\label{eq:mvtforvarsigma}
        \begin{aligned}
            &\bigm|\varsigma(\phi,\theta)-(2r\cos{\phistar}-2r\cos{\overline{\phi}_0})\bigm|=2r\bigm|\cos\phi-\cos{\overline{\phi}_0}\bigm|\overset{\text{Mean-Value Theorem}}\le2r\bigm|\phi-\overline{\phi}_0\bigm|.
        \end{aligned}
    \end{equation} And \begin{equation}\label{eq:2ndtermexpansioncomparison}
    \begin{aligned}
        \bigm|\tau_0-\overline{\tau}_0\bigm|=&\Bigm|\frac{\varsigma(\phi,\theta)}{\eta(\phi,\theta)}-\frac{2r(\cos\phistar-\cos{\overline{\phi}_0})}{-2\sin{(\overline{\phi}_0+\overline{\theta}_0)}}\Bigm|\\
        =&\Bigm|\frac{\varsigma(\phi,\theta)\big[-2\sin{(\overline{\phi}_0+\overline{\theta}_0)}\big]-\eta(\phi,\theta\big[2r(\cos{\phistar}-\cos{\overline{\phi}_0})\big]}{(-2\sin{(\overline{\phi}_0+\overline{\theta}_0)})\eta(\phi,\theta)}\Bigm|\\
        =&\Bigm|\frac{\varsigma(\phi,\theta)\big[-2\sin{(\overline{\phi}_0+\overline{\theta}_0)}\big]-\varsigma(\phi,\theta)\eta(\phi,\theta)+\varsigma(\phi,\theta)\eta(\phi,\theta)-\eta(\phi,\theta)[2r(\cos{\phistar}-\cos{\overline{\phi}_0})]}
        {-2\sin{(\overline{\phi}_0+\overline{\theta}_0)}\eta(\phi,\theta)}\Bigm|\\
        \le&\frac{\overbracket{|\varsigma(\phi,\theta)|}^{\le4r}}{-2\sin{(\overline{\phi}_0+\overline{\theta}_0)}\underbracket{|\eta(\phi,\theta)|}_{\mathllap{\eqref{eq:etalowerbound}:>-1.9258\sin{(\overline{\phi}_0+\overline{\theta}_0)}>0}}}\underbracket{\Bigm|-2\sin{(\overline{\phi}_0+\overline{\theta}_0)}-\eta(\phi,\theta)\Bigm|}_{\text{LHS of \eqref{eq:etaaproximation}}}+\frac{\overbracket{\bigm|\varsigma(\phi,\theta)-(2r\cos{\phistar}-2r\cos{\overline{\phi}_0})\bigm|}^{\text{LHS of \eqref{eq:mvtforvarsigma}}}}{-2\sin{(\overline{\phi}_0+\overline{\theta}_0)}}\\
        \underbracket{<}_{\mathllap{\text{\eqref{eq:mvtforvarsigma}, \eqref{eq:etaaproximation} and \eqref{eq:etalowerbound}}}}&\frac{2r}{1.9258\sin^2{(\overline{\phi}_0+\overline{\theta}_0)}}\Big[\frac{7.1r}{-b\sin{(\overline{\phi}_0+\overline{\theta}_0)}}+\frac{3|\theta-\overline{\theta}_0|}{-\sin{(\overline{\phi}_0+\overline{\theta}_0)}}+\frac{3|\phi-\overline{\phi}_0|}{-\sin{(\overline{\phi}_0+\overline{\theta}_0)}}\Big]+\frac{2r\bigm|\overline{\phi}_0-\phi\bigm|}{-2\sin{(\overline{\phi}_0+\overline{\theta}_0)}}\\
        \underbracket{<}_{\mathllap{\frac{1}{-\sin{(\overline{\phi}_0+\overline{\theta}_0)}}\ge1}}&\underbracket{\frac{5.2r}{-\sin^3{(\overline{\phi}_0+\overline{\theta}_0)}}\big|\phi-\overline{\phi}_0\big|}_{\text{\eqref{eq:secondtermtau0taubar0closecondition}:}<\frac{5.2\varepsilon_2}{64}r\sin{\overline{\theta}_0}}+\underbracket{\frac{3.2r}{-\sin^3{(\overline{\phi}_0+\overline{\theta}_0)}}\big|\theta-\overline{\theta}_0\big|}_{\text{\eqref{eq:secondtermtau0taubar0closecondition}:}<\frac{3.2\varepsilon_2}{64}r\sin{\overline{\theta}_0}}+\frac{7.4r^2}{\underbracket{-b\sin^3{(\overline{\phi}_0+\overline{\theta}_0)}}_{\mathrlap{\text{\eqref{eq:secondtermtau0taubar0closecondition}:} b>R-r>\frac{64r}{-\varepsilon_2\sin^3{(\overline{\phi}_0+\overline{\theta}_0)}\sin{\overline{\theta}_0}}}}}\\\overbracket{<}^{\mathllap{\text{\eqref{eq:secondtermtau0taubar0closecondition}:}}}&\big(5.2/64+3.2/64+7.4/64\big)r\sin{\overline{\theta}_0}<\frac{\varepsilon_2}{4}\overline{d}_0.\\
    \end{aligned}
    \end{equation}\newline Therefore, $\frac{|\tau_0-\overline{\tau}_0|}{\overline{d}_0}<\frac{\varepsilon_2}{4}$.
    
    For the third term in the rightmost of \eqref{eq:LemonexpansionIcompareWithPET}, if the following condition holds that is,
    \begin{equation}\label{eq:thirdtermclosingcondition}
        \begin{aligned}
            |\theta-\overline{\theta}_0|<&\mathtt{\delta}_1\Big(\overline{\phi}_0,\,\overline{\theta}_0,\,\big(1+\frac{6r}{\sqrt{r^2\sin^2{\phistar}-\overline{\mathcal{X}}^2_1}}\big)^{-1}0.075\varepsilon_2r\sin{\overline{\theta}_0},\,\big(2r+\frac{8r^2}{\sqrt{r^2\sin^2{\phistar}-\overline{\mathcal{X}}^2_1}}\big)^{-1}0.1\varepsilon_2r\sin{\overline{\theta}_0}\Big)\\
            |\phi-\overline{\phi}_0|<&\mathtt{\delta}_2\Big(\overline{\phi}_0,\,\overline{\theta}_0,\,\big(1+\frac{6r}{\sqrt{r^2\sin^2{\phistar}-\overline{\mathcal{X}}^2_1}}\big)^{-1}0.075\varepsilon_2r\sin{\overline{\theta}_0},\,\big(2r+\frac{8r^2}{\sqrt{r^2\sin^2{\phistar}-\overline{\mathcal{X}}^2_1}}\big)^{-1}0.1\varepsilon_2r\sin{\overline{\theta}_0}\Big)\\
            R>&\mathtt{R}_1\Big(\overline{\phi}_0,\overline{\theta}_0,\,\big(1+\frac{6r}{\sqrt{r^2\sin^2{\phistar}-\overline{\mathcal{X}}^2_1}}\big)^{-1}0.075\varepsilon_2r\sin{\overline{\theta}_0},\,\big(2r+\frac{8r^2}{\sqrt{r^2\sin^2{\phistar}-\overline{\mathcal{X}}^2_1}}\big)^{-1}0.1\varepsilon_2r\sin{\overline{\theta}_0}\Big),
        \end{aligned}
    \end{equation}then \cref{lemma:Trajectory_approximation_1}, \eqref{eq:x0barx0ClosingCondition} and \eqref{eq:eq:x0barx0Closing} conclude that the conditions in \eqref{eq:thirdtermclosingcondition} imply the following condition. \begin{equation}\label{eq:X1Y1Theta1closingforTau1}\begin{aligned}
        \bigm|\mathcal{X}_1-\overline{\mathcal{X}}_1\bigm|<&\big(1+\frac{6r}{\sqrt{r^2\sin^2{\phistar}-\overline{\mathcal{X}}^2_1}}\big)^{-1}\cdot0.075\varepsilon_2r\sin{\overline{\theta}_0},\\
        \bigm|\mathcal{Y}_1-\overline{\mathcal{Y}}_1\bigm|<&\big(1+\frac{6r}{\sqrt{r^2\sin^2{\phistar}-\overline{\mathcal{X}}^2_1}}\big)^{-1}\cdot0.075\varepsilon_2r\sin{\overline{\theta}_0},\\
        \bigm|\Theta_1-\overline{\Theta}_1\bigm|<&\big(2r+\frac{8r^2}{\sqrt{r^2\sin^2{\phistar}-\overline{\mathcal{X}}^2_1}}\big)^{-1}\cdot0.1\varepsilon_2r\sin{\overline{\theta}_0}.\\
    \end{aligned}\end{equation}\newline And
\eqref{eq:ttbarvarepsilonCondition} concludes that the conditions of \eqref{eq:X1Y1Theta1closingforTau1} imply \begin{equation}\label{eq:3rdtermexpansioncomparison}
    \begin{aligned}
        |\tau_1-\overline{\tau}_1|=\overbracket{|t-\overline{t}|}^{\text{in \eqref{eq:ttbarvarepsilonCondition}}}<&\frac{\varepsilon_2}{4}r\sin{\overline{\theta}_0}\\
        \frac{|\tau_1-\overline{\tau}_1|}{\overline{d}_0}<\frac{\varepsilon_2}{4}.
    \end{aligned}
\end{equation}
Therefore, the conditions in \eqref{eq:thirdtermclosingcondition} imply that $\frac{|\tau_1-\overline{\tau}_1|}{\overline{d}_0}<\frac{\varepsilon_2}{4}$.

It is clear that if \begin{equation}\label{eq:fourthtermclosingcondition}
    R>\frac{1}{64\varepsilon^2_2}r, 
\end{equation} then the fourth term in the rightmost of \eqref{eq:LemonexpansionIcompareWithPET}: $2\sqrt{r/R}<\frac{\varepsilon_2}{4}$. 

The conditions in \eqref{eq:FirstTermInLemonexpansionIcompareWithPET}, \eqref{eq:secondtermtau0taubar0closecondition}, \eqref{eq:thirdtermclosingcondition}, \eqref{eq:fourthtermclosingcondition} all together become conditions in \eqref{eq:ExpansionClosingCondition_i} to make \begin{align*}
            \Bigm|\mathcal{I}-\big(-1-\frac{\overline{\tau}_0+\overline{\tau}_1-2\overline{d}_0}{\overline{d}_0}\big)\Big|\overset{\eqref{eq:LemonexpansionIcompareWithPET}}{<}&\underbracket{\frac{1}{d_0\overline{d}_0}\bigm|\tau_0+\tau_1\bigm|\bigm|d_0-\overline{d}_0\bigm|}_{\eqref{eq:1sttermexpansioncomparison}:<0.25\varepsilon_2}+\underbracket{\frac{1}{\overline{d}_0}\bigm|\tau_0-\overline{\tau}_0\bigm|}_{\eqref{eq:2ndtermexpansioncomparison}:<0.25\varepsilon_2}+\underbracket{\frac{1}{\overline{d}_0}\bigm|\tau_1-\overline{\tau}_1\bigm|}_{\eqref{eq:3rdtermexpansioncomparison}:<0.25\varepsilon_2}+\underbracket{2\sqrt{r/R}}_{\eqref{eq:fourthtermclosingcondition}:<0.25\varepsilon_2}<\varepsilon_2\qedhere\\
\end{align*}

\end{proof}

\subsection[Uniform expansion on section sets (proof of the main theorems)]{Uniform expansion on section sets (proof of the main theorems)}\label{subsec:UniformExpansionOnSection}\hfill

Now we start to derive $R_{\text{\upshape HF}}$ in \cref{def:AcomputationForRHF} to suffice for uniform hyperbolicity. It depends on how close we need to approximate a $\mathbb{L}(r,R,\phistar)$ trajectory by a $\mathbb{P}(r,\phistar)$ trajectory, and definitely it is not the optimal/smallest $R$ to suffice for uniform hyperbolicity.\COMMENT{Need to define petal trajectory and $\bar l$ and then derive $R$\wentao{Now given in the two definitions following.}}

\begin{definition}[$\overline{l}$: sufficient times of return steps to overcome contraction]\label{def:ApettrajectoryExpansion}
    Suppose $\phistar\in(0,\tan^{-1}{(1/3)})$ to be not exceptional in \cref{def:ExceptionPhistar}. For the trajectory of $x_{
    B_1}=(3\phistar,\phistar)$, $\overline{F}^n_{\phistar}(3\phistar,\phistar)$ in \cref{def:ExceptionPhistar} for $\mathbb{P}(r,\phistar)$, suppose for $k\ge1$, the subsequence  $\overline{y}_{n_k>0}=\overline{F}_{\phistar}^{n_k>0}(x_{B_1})\in \Mrout$ and we initialize $\overline{y}_{n_0=0}=(3\phi_*,\phi_*)=x_{B_1}$ and have $n_{k+1}>n_{k}$, $\forall k\ge0$. Then per \cref{lemma:expansionbetweenPETandLEM} with definitions of length functions from \cref{def:LemonAndPetTogether}\eqref{item05LemonAndPetTogether}, for $k>0$, $\overline{\tau}_{0,k}\dfn\overline{\tau}_0(\overline{y}_{n_k})$, $\overline{\tau}_{1,k}\dfn\overline{\tau}_0(\overline{y}_{n_k})$
    $\overline{d}_{0,k}\dfn\overline{d}_0(\overline{y}_{n_k})$,
    $\overline{d}_{1,k}\dfn\overline{d}_1(\overline{y}_{n_k})$,
    $\overline{d}_{2,k}\dfn\overline{d}_2(\overline{y}_{n_k})$. With $\phistar\neq\frac{\pi}{n},\forall n\in\mathbb{N}$, \cref{lemma:PETExpExpansion} ensures that $\prod_{k=1}^{l}\bigm|-1-\frac{\overline{\tau}_{0,k}+\overline{\tau}_{1,k}-2\overline{d}_{0,k}}{\overline{d}_{0,k}}\bigm|$ grows exponentially\COMMENT{Add explanation using figure\wentao{Now providing a lemma for explanation}}.There exists a constant $\mathbb{N}\ni\overline{l}>0$ such that $\prod_{k=1}^{\overline{l}}\bigm|-1-\frac{\overline{\tau}_{0,k}+\overline{\tau}_{1,k}-2\overline{d}_{0,k}}{\overline{d}_{0,k}}\bigm|>25$.\COMMENT{Why is this called a definition?} We call $\overline{l}$ is the sufficient times of return steps to overcome contraction.
\end{definition}

\begin{definition}[Computation of $R_{\text{\upshape HF}}$]\label{def:AcomputationForRHF} From \cref{def:ApettrajectoryExpansion}, with its order-reverted orbit section segment: $\overline{y}_{n_{\overline{l}}}$, $\cdots$, $\overline{y}_{n_0=0}=x_{B_1}$. When $1\le k\le\overline{l}$, let $(\overline{\mathcal{X}}_k,\overline{\mathcal{Y}}_k)=p(\overline{F}_{\phistar}(\overline{y}_{n_k}))$ in \cref{fig:Trajectory_LemonCompareWithWojtkowski2} with coordinate system \cref{def:StandardCoordinateTable}, suppose $(\overline{\phi}_k,\overline{\theta}_k)=\overline{y}_{n_k}$ and let $\overline{E}_k=-1-\frac{\overline{\tau}_0+\overline{\tau}_1-2\overline{d}_0}{\overline{d}_0}$. We initially set $R_{\text{\upshape HF}}=r\cdot\max{\big\{\frac{30000}{\sin^2\phistar},\frac{324}{\sin^2{(\phistar/2)\sin^2(\phistar)}}\big\}}$. We have the following finite steps to iteratively update $R_{\text{\upshape HF}}$ and compute $\varepsilon_k$ by iterating $k$ decreasingly from $\overline{l}$ to $1$.
\begin{equation}\label{eq:RHFcomputation}
\begin{aligned}
        \text{For }&k=\overline{l},\cdots,1\\
        &\text{set }\varepsilon_{k-1,1}=\left\{\begin{aligned}&\min\Big\{\delta_6(\overline{\phi}_k,\overline{\theta}_k,\bigm|1-0.9^{\frac{1}{\overline{l}}}\bigm|\bigm|\overline{E}_k\bigm|,\overline{\mathcal{X}}_k),\delta_7(\overline{\phi}_k,\overline{\theta}_k,\bigm|1-0.9^{\frac{1}{\overline{l}}}\bigm|\bigm|\overline{E}_k\bigm|,r,\overline{\mathcal{X}}_k),\varepsilon_k\Big\}\text{ if }k\ne\overline{l},\\&\min\Big\{\delta_6(\overline{\phi}_k,\overline{\theta}_k,\bigm|1-0.9^{\frac{1}{\overline{l}}}\bigm|\bigm|\overline{E}_k\bigm|,\overline{\mathcal{X}}_k),\delta_7(\overline{\phi}_k,\overline{\theta}_k,\bigm|1-0.9^{\frac{1}{\overline{l}}}\bigm|\bigm|\overline{E}_k\bigm|,r,\overline{\mathcal{X}}_k)\Big\}\text{ if }k=\overline{l},
        \end{aligned}\right.\\
        &\text{set }R_{\text{\upshape HF}}=\max\big\{R_{\text{\upshape HF}},\,r+\mathtt{R}_2(\overline{\phi}_k,\overline{\theta}_k,\bigm|1-0.9^{\frac{1}{\overline{l}}}\bigm|\bigm|\overline{E}_k\bigm|,\overline{\mathcal{X}}_k)\big\},\\
        &\text{set } \varepsilon_{k-1,2}=\left\{\begin{aligned}&\min\big\{\delta_3(\varepsilon_{k},\overline{\mathcal{X}}_{k}),\delta_4(\varepsilon_{k},\overline{\mathcal{X}}_{k})\big\},\text{ if }k\ne\overline{l},\\&\varepsilon_{k-1,1},\text{ if }k=\overline{l},\end{aligned}\right.\\
        &\text{set }\varepsilon_{k-1,3}=\left\{\begin{aligned}&\delta_5(\varepsilon_{k},\overline{\mathcal{X}}_{k}),\text{ if }k\ne\overline{l},\\&\varepsilon_{k-1,1},\text{ if }k=\overline{l},\end{aligned}\right.\\
        &\text{set }\varepsilon_{k-1,4}=\left\{\begin{aligned}&\min\Big\{\delta_1(\overline{\phi}_k,\overline{\theta}_k,\varepsilon_{k-1,2},\varepsilon_{k-1,3}),\delta_2(\overline{\phi}_k,\overline{\theta}_k,\varepsilon_{k-1,2},\varepsilon_{k-1,3})\Big\},\text{ if }k\ne\overline{l},\\&\varepsilon_{k-1,1},\text{ if }k=\overline{l},\end{aligned}\right.\\
        &\text{set }R_{\text{\upshape HF}}=\max\big\{R_{\text{\upshape HF}},r+\mathtt{R}_1(\overline{\phi}_k,\overline{\theta}_k,\varepsilon_{k-1,2},\varepsilon_{k-1,3})\big\},\\
        &\text{set }\varepsilon_{k-1}=\left\{\begin{aligned}&\delta_0\Big(\overline{F}^2_{\phistar}(\overline{y}_{n_{k-1}}),\min\big\{\varepsilon_{k-1,1},\varepsilon_{k-1,4}\big\}\Big)\text{, if }k>1,\\
        &\delta_0\Big(x_*,\min\big\{\varepsilon_{k-1,1},\varepsilon_{k-1,4}\big\}\Big) \text{, if }k=1,
        \end{aligned}\right.
    \end{aligned}
\end{equation}
where $\delta_6,\delta_7,\mathtt{R}_2$ are defined in \cref{lemma:expansionbetweenPETandLEM}\eqref{eq:ExpansionClosingCondition_i}, $\delta_3,\delta_4,\delta_5$ are defined in \cref{lemma:Trajectory_approximation_2}\eqref{eq:X1Y1Theta1CloseCondition}, $\delta_1,\delta_2,\mathtt{R}_1$ are defined in \cref{lemma:Trajectory_approximation_1}\eqref{eq:x0barx0ClosingCondition}. $\delta_0$ is defined in \cref{lemma:ApproximationForConsecutiveCollisions}\eqref{eq2:ApproximationForConsecutiveCollisions}\eqref{eq4:ApproximationForConsecutiveCollisions}. The procedure \eqref{eq:RHFcomputation} terminates with computed $\varepsilon_0$ and the updated $R_{\text{\upshape HF}}$. We then finally update $R_{\text{\upshape HF}}=\max\big\{R_{\text{\upshape HF}}, \frac{289r}{\varepsilon_0\sin^2(\phistar/2)}\big\}$.
\end{definition}
\begin{remark}\label{remark:GaranteeLargeTheta}
    By making $R_{\text{\upshape HF}}(r,\phistar)$ as large as in \cref{def:AcomputationForRHF}, by \cref{lemma:Trajectory_approximation_1}\eqref{eq:xorbitA0Condition} we see that we have ensured that if $R>R_{\text{\upshape HF}}(r,\phistar)$ for each point nongsingular $x\in\Nin$ the neighborhoods of $(\phistar,\phistar)$ and $(2\pi-2\phistar,\pi-\phistar)$ $\exists j(x)\ge0$ s.t. $p(\Fc^i(x))\in\Gamma_r$, $i=0,\cdots,j(x)$ and $p(\Fc^{l(x)}(x))\in\Gamma_R$ with $\Fc^{l(x)}(x)=(\Phi,\theta)\in M_R$ having $\theta>\sin^{-1}{(\sqrt{r/R})}$. That is, orbits starting from $\Nin$ ( $\Nout$ ) cannot have 2 consecutive returns to $\Mrin$ ( $\Mrout$ ) being returns to $\Nin$ ( $\Nout$ ).
\end{remark}
\begin{proposition}[$R_{\text{\upshape{HF}}}$ ensures expansion for orbits starting from $\Fc^{-1}(\Nout)$]\label{proposition:ExpansionWins}
    Suppose that the lemon billiard $\mathbb{L}(r,R,\phistar)$ has $\phistar\in(0,\tan^{-1}{(1/3)})$ and $\phistar$ not being generalized exceptional $\phistar$ in \cref{def:ExceptionPhistar}. For $R_{\text{\upshape HF}}(r,\phistar)>0$ defined in \cref{def:AcomputationForRHF}, the lemon billiard $\mathbb{L}(r,R,\phistar)$ with $R$ satisfying \eqref{eqJZHypCond} ensures the following:\newline With $dx_3=D\Fc_{x_2}(dx_2)$ for $dx_2$ defined in \cref{lemma:ExpansionLowerBoundA0A1,lemma:ExpansionLowerBoundB,lemma:AformularForLowerBoundOfExpansionForn1largerthan0} corresponding to (a1)(b)(c)cases in \eqref{eqMainCases}, \[\frac{\|D\hat{F}^{\overline{l}}_x(dx)\|_{\p}}{\|dx_3\|_{\p}}>25\times0.9=22.5\]
\end{proposition}
\begin{proof}
     \cref{def:ApettrajectoryExpansion,def:AcomputationForRHF,lemma:Trajectory_approximation_1,lemma:Trajectory_approximation_2,lemma:expansionbetweenPETandLEM}
     ensure that for $\hat{F}$ iterations, after the cases (a1)(b)(c) orbit segment (in \eqref{eqMainCases}) there will have at least $\overline{l}$ times case (a0) return orbit segment following. Suppose that in the $j$-th following return orbit segment, the $x_1$ in \cref{def:MhatReturnOrbitSegment} is $x_{1,j}=(\Phi_{1,j},\theta_{1,j})$, for $j=1,\cdots,\overline{l}$.

     In \cref{def:ApettrajectoryExpansion,def:AcomputationForRHF}, \cref{lemma:expansionbetweenPETandLEM} ensures that each $j$-th following orbit segment of $\hat{M}$ return has $\mathcal{I}_j$ (in \cref{def:AformularForLowerBoundOfExpansion}) satisfying $\bigm|\mathcal{I}_j-\overline{E}_j\bigm|<\bigm|1-0.9^{\frac{1}{\overline{l}}}\bigm|\bigm|\overline{E}_j\bigm|$ so that $\bigm|\mathcal{I}_j\bigm|\ge\bigm|\overline{E}_j\bigm|-\bigm|\overline{E}_j-\mathcal{I}_j\bigm|>(1-|1-0.9^{\frac{1}{\overline{l}}}|)\bigm|\overline{E}_j\bigm|=0.9^{\frac{1}{\overline{l}}}\bigm|\overline{E}_j\bigm|$. Therefore, \[\frac{\|D\hat{F}^{\overline{l}}_x(dx)\|_{\p}}{\|dx_3\|_{\p}}\overbracket{>}^{\text{\cref{lemma:ExpansionLowerBoundA0A1}}}\prod_{j=1}^{\overline{l}}\bigm|\mathcal{I}_j\bigm|>\prod_{j=1}^{\overline{l}}\bigm|(0.9)^{\frac{1}{\overline{l}}}E_j\bigm|\overbracket{>}^{\text{\cref{def:ApettrajectoryExpansion}}}25\times0.9=22.5\qedhere\]
\end{proof}
\begin{definition}[uniform hyperbolic lemon billiard]\label{def:UniformHyperbolicLemonBilliard}
    Now we define our \emph{uniform hyperbolic lemon billiard} to be the $\mathbb{L}(r,R,\phistar)$ with $\phistar\in(0,\tan^{-1}{(1/3)})$, $\phistar$ not being the \emph{generalized exceptional $\phistar$} defined in \cref{def:ExceptionPhistar}, and with $R_{\text{\upshape HF}}(r,\phistar)>0$ given in \cref{def:AcomputationForRHF}.
\end{definition}
\begin{theorem}[\cref{RMASUH} uniform expansion of $\hat F$]\label{RMASUHwithProof}\label{thm:UniformExpansionOnMhat}
For the lemon billiards $\LEM(r,R,\phistar)$ satisfying \cref{def:UniformHyperbolicLemonBilliard} and for the return map $\hat{F}$ defined on $\hat{M}$ in \cref{def:Mhatdefinition}, we have the following. 
\begin{enumerate}
    \item\label{item01RMASUHwithProof} The cone family $C_x\dfn\big\{(d\phi,d\theta)\bigm|\frac{d\theta}{d\phi}\in [0,1]\big\}$ is strictly invariant under the return map $\hat{F}$ from \cref{def:SectionSetsDefs}.
    \item\label{item02RMASUHwithProof} There exist $c>0$ and $\Lambda>1$ such that $\dfrac{\|D\hat{F}_x^n(dx)\|_\p}{\|dx\|_\p}>c\Lambda^n$ for nonsingular $x\in\hat{M}$, $dx\in C_x$ and $n\in\mathbb{N}$.\COMMENT{Edit. Do not refer to assumptions in proposition. Say that this the main theorem restated. Or state that theorem exactly and prove here.}
\end{enumerate}
\end{theorem}
\begin{proof}
    The first conclusion that $C_x$ is strictly invariant is proved by \cref{corollary:InvQuadrant}, \eqref{eq:strictlyinvcone_casea1bc}, and \cref{corollary:halfquadrantcone}. 
    
    For the second conclusion, it remains to show the uniform exponential expansion for nonsingular $x\in\hat{M}$ with its return orbit segment cases \eqref{eqMainCases}. 
    
    In case (a0) of  \eqref{eqMainCases}, $\frac{\|D\hat{F}_x(dx)\|_{\p}}{\|dx\|_{\p}}>1+\lambda_c$ with constant $\lambda_c>0$ from \cref{TEcase_a0}.

    In cases (a1)(b)(c) of \eqref{eqMainCases}, for the $x_2$ and $dx_2$ defined in \cref{lemma:ExpansionLowerBoundA0A1,lemma:ExpansionLowerBoundB,lemma:AformularForLowerBoundOfExpansionForn1largerthan0}, with $dx_3=D\Fc_{x_2}(dx_2)$, \cref{lemma:ExpansionLowerBoundA0A1,lemma:ExpansionLowerBoundB,lemma:AformularForLowerBoundOfExpansionForn1largerthan0} give $\frac{\|dx_3\|_{\p}}{\|dx\|_{\p}}>0.05$.  Then with our $R_{\text{\upshape HF}}$ in \cref{def:AcomputationForRHF}\eqref{eq:RHFcomputation} and the $\mathbb{L}(r,\phistar)$ trajectory orbit section segment $\big\{y_{n_k}\big\}_{k=0}^{\overline{l}}$ in \cref{def:ApettrajectoryExpansion}, \cref{proposition:ExpansionWins} gives $\frac{\|D\hat{F}^{\overline{l}}_x(dx)\|_{\p}}{\|dx_3\|_{\p}}>25\times0.9=22.5$. Therefore, in cases (a1)(b)(c), $\frac{\|D\hat{F}^{\overline{l}+1}_x(dx)\|_{\p}}{\|dx\|_{\p}}>0.05\times22.5=1.125$.

    Let $\overline{\Lambda}=\min\{1+\lambda_c,1.125\}$, where $\lambda_c$ is from \cref{TEcase_a0}. \cref{def:ApettrajectoryExpansion,def:AcomputationForRHF,lemma:Trajectory_approximation_1}\eqref{eq:xorbitA0Condition} ensure that among every consecutive $\overline{l}+1$ times $\hat{M}$ return orbit segment defined in \eqref{eqMhatOrbSeg} there can exist at most $1$ time return orbit segment not in case (a0) of \eqref{eqMainCases}. Therefore, We can conclude that for all nonsingular $x$ and $\forall n\ge1$, \[\frac{\|D\hat{F}^{n}_x(dx)\|_{\p}}{\|dx\|_{\p}}>(\overline{\Lambda})^{\left\lfloor\frac{n}{\overline{l}+1}\right\rfloor}(0.05).\] 
    Therefore, $\frac{\|D\hat{F}^{n}_x(dx)\|_{\p}}{\|dx\|_{\p}}>(\overline{\Lambda})^{\frac{n}{\overline{l}+1}-1}(0.05)=c\Lambda^n$, where $c=(0.05)/\overline{\Lambda}$, $\Lambda=(\overline{\Lambda})^{\frac{1}{1+\overline{l}}}$.\COMMENT{Move to beginning}
\end{proof}
\begin{remark}
    \cref{RMASUHwithProof} is \cref{RMASUH} since the generalized exceptional $\phi_*\in(0,\tan^{-1}(1/3))$ has $0$ measure per \cref{remark:exceptionalphistar}.
\end{remark}

\begin{corollary}[In the context of \cref{RMASUHwithProof}]\label{corollary:sinpmetricexpansion}\label{corollary:sinpmetricexpansionMhat}
There exist constant $\hat{c}>0$ and $\Lambda>1$ such that $\dfrac{\sin\theta\|D\hat{F}_x^n(dx)\|_\p}{\|dx\|_\p}>\hat{c}\Lambda^n$ for all nonsingular $x=(\phi,\theta)\in\hat{M}$, $dx\in C_x$ and $n\in\mathbb{N}$.
\end{corollary}
\begin{proof}
    Let $x=\hat{x}_0=(\phi,\theta)\in\hat{M}$, $(d\phi,d\theta)=dx\in C_{x}$, $\hat{x}_n=\hat{F}^n(x)$, $d\hat{x}_n=D\hat{F}_x^n(dx)$, $n\ge0$.
    \begin{equation}\label{eq:pexpansion_sin}
\displaystyle\sin\theta\frac{\|d\hat{x}_n\|_\p}{\|dx\|_\p}=\frac{\|d\hat{x}_n\|_\p}{\|d\hat{x}_1\|_\p}\frac{\|d\hat{x}_1\|_\p\sin{\theta}}{\|dx\|_\p}
\end{equation}
For the orbit segment of $x$, defined in \cref{def:MhatReturnOrbitSegment} and \eqref{eqMhatOrbSeg} with points $x_0=(\phi_0,\theta_0)\in\Mrout$, $x_1=(\Phi_1,\theta_1)\in\MRin$ and $x_2\in\Mrin$ in \eqref{eqPtsFromeqMhatOrbSeg}. Note that since for either $x\in\Mrin\cap\Mrout$ or $\Fc^{-1}(\Mrout\smallsetminus\Mrin)$, we have $\theta=\theta_0$ since either $x=x_0$ or $x$, $x_0$ are consecutive collisions on $\Gamma_r$. 

We analyze the two cases classified by $\theta_1$.

case (i) $\sin{\theta_1}<\sqrt{4r/R}$ that are cases (a1)(b)(c) of \eqref{eqMainCases}, \cref{TEcase_a1,TEcase_b,TEcase_c} imply that $\frac{\|d\hat{x}_1\|_\p}{\|dx\|_\p}>0.05$.

And by \cref{contraction_region} and \eqref{eqJZHypCond}, $|\sin{\theta_0}-\sin{\phistar}|\overset{x_0\in\Nout}<\frac{17}{\sin{(\frac{\phistar}{2})}}\sqrt{\frac{r}{R}}\overset{\eqref{eqJZHypCond}}<\frac{17\sin\phistar}{18}$, thus $\sin\theta_0$ is bounded away from 0. Therefore in this case there exists $C_0>0$ such that $\frac{\|d\hat{x}_1\|_\p\sin{\theta}}{\|dx\|_\p}=\frac{\|d\hat{x}_1\|_\p\sin{\theta_0}}{\|dx\|_\p}>C_0$. 

case (ii) $\sin{\theta_1}\ge\sqrt{4r/R}$ that are cases (a0) of \eqref{eqMainCases}, therefore, \cref{TEcase_a0} gives $\frac{\|d\hat{x}_1\|_\p\sin{\theta_0}}{\|dx\|_\p}>\sin\theta_0$. 

On the other hand, by \cref{def:AformularForLowerBoundOfExpansion} and computation in \cref{subsec:AFELB},
\begin{equation}\label{eq:casea0onestepsinexpansion}
\begin{aligned}
\frac{\|d\hat{x}_1\|_\p\sin{\theta_0}}{\|dx\|_\p}&\underbracket{>}_{\mathclap{\substack{\text{\cref{lemma:ExpansionLowerBoundA0A1}}\\\text{\cref{def:AformularForLowerBoundOfExpansion}}}}}\Big|-1+\frac{\tau_1}{d_1}\frac{2(\tau_0-d_0)}{d_0}-\frac{\tau_0+\tau_1-2d_0}{d_0}\Big|\sin\theta_0\\
    & =\Big|\underbracket{-\sin\theta_0}_{\mathllap{\rightarrow0}}+\underbracket{\frac{\tau_1}{d_1}\frac{2(\tau_0-d_0)}{d_0}\sin\theta_0}_{\mathllap{\text{\cref{remark:deffpclf}: }\rightarrow0}}-\underbracket{\frac{\tau_0+\tau_1-2d_0}{d_0}\sin\theta_0}_{\mathclap{\substack{\quad=\frac{\tau_0+\tau_1-2d_0}{r}\text{ since }d_0=r\sin{\theta_0}\\
    \text{Since }d_0\rightarrow0,\rho\rightarrow r\text{, \eqref{eq:NewCoordinateLegnth}: }\theta_1\rightarrow\phistar-\Phistar\text{ or }\pi-\phistar+\Phistar\\\rho\rightarrow r\text{ in \eqref{eq:continuousExtension}: }\tau_0+\tau_1\rightarrow 2r\sin{(2\phistar-2\Phistar)}}}}\Big|\xrightarrow{\sin\theta_0\rightarrow0}2\sin{(2\phistar-2\Phistar)}
\end{aligned}
\end{equation}
If $\sin{\theta_0}\rightarrow0$, then $d_0=r\sin\theta_0\rightarrow0$ and $0\le\tau_0<2d_0\rightarrow 0$. This implies the collision position $p(x_1)$ will be approximating the corners, that is, $p(x_1)\rightarrow A \cup B$.

Thus in equations \eqref{eq:lengthextensionA}, \eqref{eq:lengthextensionB}, $\theta_0\longrightarrow \{0\} \cup \{\pi\}$ will yield $\rho\rightarrow r$, $d_0\rightarrow 0$, $\theta_1\rightarrow\{\phistar-\Phistar\}\cup\{ \pi-\phistar+\Phistar\}$,  $\tau_0\rightarrow0$ and $\tau_1\rightarrow 2r\sin{(2\phistar-2\Phistar)}$.

Since $\frac{\tau_1}{d_1}\in[0,2],\frac{2(\tau_0-d_0)}{d_0}\in[-2,2]$, the right hand side of the last equality in \eqref{eq:casea0onestepsinexpansion} $\rightarrow 2\sin{(2\phistar-2\Phistar)}$ and is bounded away from 0 as $\sin\theta_0\rightarrow0$. Hence, there exists $C_1>0$ such that $\frac{\|d\hat{x}_1\|_\p\sin{\theta}}{\|dx\|_\p}>C_1$ in case (ii). 

in both case (i) and case (ii), there exists $C=\min\big\{C_0,C_1\big\}>0$ such that $\frac{\|d\hat{x}_1\|_\p\sin{\theta}}{\|dx\|_\p}=\frac{\|d\hat{x}_1\|_\p\sin{\theta_0}}{\|dx\|_\p}>C$. Thus \eqref{eq:pexpansion_sin} gives \[\sin\theta\frac{\|d\hat{x}_n\|_\p}{\|dx\|_\p}=\frac{\|d\hat{x}_n\|_\p}{\|d\hat{x}_1\|_\p}\frac{\|d\hat{x}_1\|_\p\sin{\theta}}{\|dx\|_\p}\overbracket{>}^{\mathrlap{\text{\cref{RMASUHwithProof}}}}C\cdot c\Lambda^{n-1}=\underbracket{cC\Lambda^{-1}}_{\mathclap{=\hat{c}}}\Lambda^{n}\qedhere\]
\end{proof}

\begin{corollary}[of \cref{RMASUHwithProof}, uniform expansion in Euclidean Metric] \label{Corollary:EuclideanExpansionMhat}\label{Col:EuclideanExpansionMhat} Under the return map $\hat{F}$ on $\hat{M}$ the vectors in the unstable cone $\mathcal{HQ}_x(\mathrm{I,III})=C_x$ defined in \cref{RMASUHwithProof} expand uniformly with respect to the Euclidean metric $\|\,\|$, i.e., for nonsingular $x\in\hat{M}$, $dx=(d\phi,d\theta)\in C_x$ and all $n\ge0$, there exist $c_0>0,\Lambda>1$ such that $\dfrac{\|D\hat{F}_x^n(dx)\|}{\|dx\|}>c_0\Lambda^n$.
\end{corollary}
\begin{proof}
    Let $x=(\phi,\theta)\in\hat{M}$, $(d\phi,d\theta)=(d\hat{\phi}_0,d\hat{\theta}_0)=d\hat{x}\in C_{x}$, $\hat{x}_n=\hat{F}^n(\hat{x}_0)$, $(d\hat{\phi}_n,d\hat{\theta}_n)=d\hat{x}_n=D\hat{F}_x^n(dx),\mathcal{V}_n=\frac{d\hat{\theta}_n}{d\hat{\phi}_n}\in[0,1]$, $n\ge0$.
Reasoning as in \cite[equation (8.21)]{cb}, it gives
\begin{equation}\label{eq:pmetricConvertEuclideanmetric}
\displaystyle\frac{\|d\hat{x}_n\|}{\|dx\|}=\frac{\|d\hat{x}_n\|_\p}{\|dx\|_\p}\frac{\sin{\theta}}{\sin{\hat{\theta}_n}}\frac{\sqrt{1+\mathcal{V}^2_n}}{\sqrt{1+\mathcal{V}^2_0}}\ge\sin{\theta}\frac{\|d\hat{x}_n\|_\p}{\|dx\|_\p}\frac{1}{\sqrt{2}}\underbracket{>}_{\text{\cref{corollary:sinpmetricexpansion}}}\frac{\hat{c}}{\sqrt{2}}\Lambda^n\nfd c_0\Lambda^n.\qedhere
\end{equation}
\end{proof}

\begin{theorem}[Uniform expansion on section $\hat{M}_1$]\label{thm:UniformExpansionOnMhat1} For the lemon billiards $\LEM(r,R,\phistar)$ satisfying \cref{def:UniformHyperbolicLemonBilliard} and for the return map $\hat{F}_1$ on $\hat{M}_1\dfn(\Mrin\cap\Mrout)\sqcup\Fc(\Mrin\smallsetminus\Mrout)\aeq \underbracket{\Min_{r,0}\sqcup \Fc(\Min_{r,1})}_{\subset\Mrout}\sqcup\underbracket{\Fc(\bigsqcup_{i\ge2}\Min_{r,i})}_{\mathrlap{\subset M_r\smallsetminus\Mrout}}$ we have the following.

\begin{enumerate}[label={(\roman*)}]
    \item\label{item01:UniformExpansionOnotherSections}
    For nonsingular $x\in\hat{M}_1$, if $dx\in \hat{C}^{u}_{1,x}\dfn\begin{cases}
        \big\{(d\phi,d\theta)\in T_x\hat{M}_1\bigm|\frac{d\theta}{d\phi}\in[0,1]\big\}\text{, if } x\in \Mrout,\\
            \big\{(d\phi,d\theta)\in T_x\hat{M}_1\bigm|\frac{d\theta}{d\phi}\in[0,\hat{\lambda}]\big\} \text{ elsewhere,}
    \end{cases}$\newline where $\hat{\lambda}=\max\{1,\lambda_2\}$ with $\lambda_2$ from \cref{theorem:B3minusdx3dphi3ratiobound}, then $D\hat{F}_{1,x}(dx)\in\big\{\text{interior of }\hat{C}^{u}_{1,\hat{F}_1(x)}\big\}$.
    \item\label{item02:UniformExpansionOnotherSections} For nonsingular $x\in\hat{M}_1$, $dx\in \hat{C}^u_{1,x}$, $\exists c_1>0$, $\Lambda_1>1$ such that $\frac{\|D\hat{F}^n_{1,x}(dx)\|_\p}{\|dx\|_\p}>c_1\Lambda_1^n$, $\forall n\ge1$.
    \end{enumerate}
\end{theorem}
\begin{proof}
    We use positive quadrant/half quadrant notation from \cref{corollary:halfquadrantcone}: $\mathcal{Q}_x(\mathrm{I,III})=\big\{(d\phi,d\theta)\in T_x M\bigm|\frac{d\theta}{d\phi}\in[0,+\infty]\big\}$, $\mathcal{HQ}_x(\mathrm{I,III})=\big\{(d\phi,d\theta)\in T_x M\bigm|\frac{d\theta}{d\phi}\in[0,1]\big\}$.

    \textbf{proof of (i)} is similar to \cref{corollary:InvQuadrant}. For nonsingular $x\in\hat{M}_1$, $\hat{F}_1$ return orbit segment (similar to \eqref{eqMhatOrbSeg}) is 
    \begin{equation}\label{eqhatF1returnorbit}\underbracket{x}_{\in\hat{M}_1},\underbracket{\Fc(x),\cdots}_{\notin\hat{M}_1},\underbracket{\Fc^{\sigma_1(x)}(x)=\hat{F}_1(x)}_{\in\hat{M}_1}\text{ where }\sigma_1(x)\dfn\inf\big\{k>0\mid \mathcal{F}^{k}(x)\in\hat{M}_1\big\}.\end{equation} 
    Note that it is not hard see $\sigma_1(x)\ge2$. For $\forall dx\in\mathcal{Q}_x(\mathrm{I,III})$, 
    if $x\in\Min_{r,0}\sqcup\overbracket{\Fc(\Min_{r,1})}^{\mathclap{=\Mout_{r,1}}}$, then $x_0\dfn x$, $dx_0\dfn dx$. And if $x\in\Fc(\Min_{r,i})$ for some $i\ge2$, then $\Mrout\ni x_0\dfn\Fc^{m(\Fc^{-1}(x))-1}(x)$ and $dx_0\dfn D\Fc_{x}^{m(\Fc^{-1}(x))-1}(dx)\in \mathcal{Q}_{x_0}(\mathrm{I,III})$ since $D\Fc_{x}^{m(\Fc^{-1}(x))-1}=\begin{pmatrix}
1 & 2(m(\Fc^{-1}(x))-1) \\
0 & 1 
\end{pmatrix}$.

Therefore, the same as \eqref{eqPtsFromeqMhatOrbSeg}, the orbit segment \eqref{eqhatF1returnorbit} contains $x_0$, $x_1=\Fc(x_0)\in\MRin$, $x_2=\Fc^{n_1+1}(x_1)\in\Mrin$ with $n_1\ge0$. We define $dx_1=D\Fc_{x_0}(dx_0)$, $dx_2=D\Fc^{n_1+1}_{x_1}(dx_1)$. There are length functions $\tau_0,\tau_1,d_0,d_1,d_2$ defined in \cref{def:fpclf}, as well as three cases of $d_1$, $n_1$ by \eqref{eqMainCases}.

    For case (a): $d_1\ge2r$, \cref{caseA}: $D\Fc^2_{x_0}$ is a negative matrix. Since $dx_0\in\mathcal{Q}_{x_0}(\mathrm{I,III})$, 
    \begin{equation}\label{eq:Fhat1dx2caseA}
    dx_2=D\hat{F}_1(dx)=D\Fc_{x_0}^2(D\Fc_{x}^{m(\Fc^{-1}(x))-1}(dx))=D\Fc_{x_0}^2(dx_0)\in\big\{\text{ interior of }\mathcal{Q}_{x_2}\big\}
    \end{equation} 
    
    And especially if $\hat{F}_1(x)\in\Min_{r,0}$, then $x_2=\hat{F}_1(x)=\Fc^2(x_0)\in\Min_{r,0}$. And if $\hat{F}_1(x)\in\Fc(\Min_{r,1})$, then $x_2\in\Min_{r,1}$, $\hat{F}_1(x)=\Fc(x_2)$. 
    
    Therefore, if $x_2\in\Mrin\cap\Mrout$, then from \cref{def:DF2matrix} and \eqref{eq:DF2matrix} $d_1d_2D\Fc^2_{x_0}=\begin{bmatrix}
    D_{11} & D_{12}\\
    D_{21} & D_{22}
    \end{bmatrix}$ satisfies $0<\frac{D_{21}}{D_{11}}<\frac{D_{22}}{D_{12}}<1$, thus by \eqref{eq:narrowconereason1}, \eqref{eq:narrowconereason2} and \eqref{eq:narrowconereason3} 
    \begin{equation}\label{eq:Mhat1DF2}
        D\Fc^2_{x_0}(dx_0)\in\big\{\text{interior of }\mathcal{HQ}_{x_2}(\mathrm{I,III})\big\}
    \end{equation}
    In the case $x_2=\hat{F}_{1,x}(x)$, $D\hat{F}_{1,x}(dx)=D\Fc_{x_0}^2(dx_0)\in \mathcal{HQ}_{\hat{F}_1(x)}(\mathrm{I,III})$.
    
    In the case $\Fc(x_2)=\hat{F}_1(x)$, $D\hat{F}_{1,x}(dx)=D\Fc_{x_2}(D\Fc_{x_0}^2(dx_0))\in \mathcal{HQ}_{\hat{F}_1(x)}(\mathrm{I,III})$ since $D\Fc_{x_2}=\begin{bmatrix}
    1 & 2\\
    0 & 1
    \end{bmatrix}$.
    
    Otherwise $x_2\in\Fc(\Mrin\smallsetminus\Mrout)$, then by \eqref{eq:Fhat1dx2caseA}, 
    \begin{equation}
    \begin{aligned}
        D\Fc_{x_2}=\begin{bmatrix}
    1 & 2\\
    0 & 1
    \end{bmatrix}\\
        D\hat{F}_{1,x}(dx)=D\Fc_{x_2}(\underbracket{dx_2}_{\mathclap{\eqref{eq:Fhat1dx2caseA}: \in \mathcal{Q}_{x_2}(\mathrm{I,III})}})&\in\big\{\text{interior of }\mathcal{HQ}_{\hat{F}_1(x)}(\mathrm{I,III})\big\}\subset\big\{\text{interior of }\hat{C}^u_{1,\hat{F}_1(x)}\big\}\\
    \end{aligned}\end{equation}

    For cases (b) and (c), per \cref{contraction_region}, $x_0\in\Nout$, $x_2\in\Nin$. Since \cref{Prop:Udisjoint}, $\Nout\cap(\Mout_{r,0}\sqcup\Mout_{r,1}\sqcup\Mout_{r,2})=\emptyset$ and $\Nin\cap(\Min_{r,0}\sqcup\Min_{r,1}\sqcup\Min_{r,2})=\emptyset$. Therefore, $x\neq x_0$ and $x_2\notin\Min_{r,0}$, thus $\hat{F}_1(x)=\Fc(x_2)$. Since $x_0=\Fc^{m(\Fc^{-1}(x))-1}(x)$, $m(\Fc^{-1}(x))\ge3$, $\exists x_{-1}\dfn \Fc^{-1}(x_0)$ s.t. $x_{-1}=\Fc^{m(\Fc^{-1}(x))-2}(x)$, $dx_{-1}=D\Fc_{x}^{m(\Fc^{-1}(x))-2}(x)$ with $D\Fc_{x}^{m(\Fc^{-1}(x))-2}=\begin{pmatrix}
1 & 2(m(\Fc^{-1}(x))-2) \\
0 & 1 
\end{pmatrix}.$ Therefore, $dx_{-1}\in\mathcal{Q}_{x_{-1}}(\mathrm{I,III})$.

In case (b) \cref{caseB}: $G=D\Fc^4_{x_{-1}}$ is a negative matrix. \cref{Proposition:DF4ConeUpperBound} and \eqref{eq:casebdtheta3dphi3ratio} shows for $x_3=\Fc^4(x_{-1})=\hat{F}_{1}(x)$. And $(d\phi_3,d\theta_3)\nfd dx_3=D\hat{F}_1(dx)=D\Fc_{x_{-1}}^4(D\Fc_{x}^{m(\Fc^{-1}(x))-2}(dx))=D\Fc_{x_{-1}}^4(dx_{-1})$ satisfies $0<\frac{d\theta_3}{d\phi_3}<1$ therefore is the interior of $\hat{C}^u_{1,\hat{F}_1(x)}$.

In case (c) \cref{caseC}: $D\Fc^{n_1+4}_{x_{-1}}$ is a positive matrix, $x_3=\Fc^4(x_{-1})=\hat{F}_{1}(x)$. Hence \[D\hat{F}_1(dx)=D\Fc_{x_{-1}}^{4+n_1}(D\Fc_{x}^{m(\Fc^{-1}(x))-2}(dx))=D\Fc_{x_{-1}}^{4+n_1}(dx_{-1})=dx_3\nfd(d\phi_3,d\theta_3)\] satisfies $0<\frac{d\theta_3}{d\phi_3}<\lambda_2$ by \cref{theorem:B3minusdx3dphi3ratiobound}\eqref{eq:casebdtheta3dphi3ratio}. Thus, $D\hat{F}_1(dx)\in\Int(\hat{C}^u_{1,\hat{F}_1(x)})$.



\textbf{proof of (ii)}: The proof argument is the same as in \cref{thm:MroutreturnUniformExpansion}. The return orbit segment \eqref{eqhatF1returnorbit} contains $x_0\in\Mrout$, $x_1=\Fc(x_0)\in\MRin$, $x_2=\Fc^{n_1+1}(x_1)\in\Mrin$. And if $x_0\in\Nin$, $\hat{F}_1(x)=\Fc(x_2)\nfd x_3$ and $x_{-1}=\Fc^{-1}(x_0)$ are also contained in the return orbit segment \eqref{eqhatF1returnorbit}.

The same conclusions \cref{TEcase_a0,TEcase_a1,TEcase_b,TEcase_c} as in the previous sections hold: 

In case (a0) and (a1) $d_1\ge2r$, $\frac{\|D\hat{F}_{1,x}(dx)\|_{\p}}{\|dx\|_{\p}}\ge\frac{\|dx_2\|_{\p}}{\|dx_0\|_{\p}}>0.26$.

In case (b) $d_1<2r$ and $n_1=0$, $\frac{\|D\hat{F}_{1,x}(dx)\|_{\p}}{\|dx\|_{\p}}\ge\frac{\|dx_3\|_{\p}}{\|dx_{-1}\|_{\p}}>0.05$.

In case (c) $d_1<2r$ and $n_1\ge1$, $\frac{\|D\hat{F}_{1,x}(dx)\|_{\p}}{\|dx\|_{\p}}\ge\frac{\|dx_3\|_{\p}}{\|dx_{-1}\|_{\p}}>0.9$.

\cref{caseA,caseB,caseC,TEcase_a0,TEcase_a1,TEcase_b,TEcase_c,corollary:halfquadrantcone,Corollary:CollisionOnMrin0Mrin1} imply the following:
\begin{itemize}
    \item For nongsingular $x\in\hat{M}_1$, we denote by $k_1(x)\dfn\inf\big\{k>0\bigm|\Fc^{k_1(x)}(x)\in\hat{M}\big\}$ and $d\hat{x}\dfn D\Fc_x^{k_1(x)}(dx)$. For $\forall dx\in C^u_{1,x}$, $d\hat{x}_1=D\Fc_x^{k(x)}(dx)\in C_x$ and $\frac{\|d\hat{x}_1\|_{\p}}{\|dx\|_{\p}}>0.05$,
    \item For nongsingular $x\in\hat{M}$, we denote by $k(x)\dfn\inf\big\{k>0\bigm|\Fc^{k(x)}(x)\in\hat{M}_1\big\}$ and $d\hat{x}_1\dfn D\Fc_x^{k(x)}(dx)$. For $\forall dx\in C_{x}$, $d\hat{x}_1=D\Fc_x^{k(x)}(dx)\in \hat{C}^u_{1,x}$ and $\frac{\|d\hat{x}_1\|_{\p}}{\|dx\|_{\p}}>0.05$,
\end{itemize}
where $C_x$ is the cone defined in \cref{RMASUHwithProof}.
    
Since we can find that every 2 consecutive return orbit segments \eqref{eqhatF1returnorbit} contain at least one element in $\hat{M}$. The $\p$-metric uniform exponential expansion on $\hat{M}$ by \cref{RMASUHwithProof} implies the $\p$-metric uniform exponential expansion on $\hat{M}_1$.
\end{proof}

\begin{corollary}[In the context of \cref{thm:UniformExpansionOnMhat1}]\label{corollary:sinpmetricexpansionMhat1}
There exist constant $\hat{c}_1>0$ and $\Lambda_1>1$ such that \[\dfrac{\sin\theta\|D\hat{F}_{1,x}^n(dx)\|_\p}{\|dx\|_\p}>\hat{c}_1\Lambda_1^n\] for all nonsingular $x=(\phi,\theta)\in\hat{M}_1$, $dx\in \hat{C}^u_{1,x}$ and $n\in\mathbb{N}$.
\end{corollary}
\begin{proof}
The proof is the same as \cref{corollary:sinpmetricexpansion}.
\end{proof}

\begin{corollary}[In the context of \cref{thm:UniformExpansionOnMhat1}, uniform expansion in Euclidean Metric] \label{Col:EuclideanExpansionMhat1}
Under the return map $\hat{F}_1$ on $\hat{M}_1$ the vectors in the unstable cone $\hat{C}^u_{1,x}$ defined in \cref{thm:UniformExpansionOnMhat1} expand uniformly with respect to the Euclidean metric $\|\,\|$, i.e., for nonsingular $x\in\hat{M}_1$, $dx=(d\phi,d\theta)\in \hat{C}^u_{1,x}$ and all $n\ge0$, there exist $c_2>0,\Lambda_1>1$ such that $\dfrac{\|D\hat{F}_{1,x}^n(dx)\|}{\|dx\|}>c_2\Lambda_1^n$.
\end{corollary}
\begin{proof}
    The proof is the same as \cref{Corollary:EuclideanExpansionMhat}.
\end{proof}

\begin{corollary}[Trajectory Expansion between $\hat{M}$ and $\hat{M}_1$ in $\p$-metric and Euclidean metric]\label{corollary:ExpansionBetweenMhatminusMhatplus} The following hold.

\begin{enumerate}
    \item\label{item01ExpansionBetweenMhatminusMhatplus} $\exists$ constant $c_3>0$ such that for all nonsingular $(\phi,\theta)=x\in\hat{M}$, $dx\in C_x$ as defined in \cref{RMASUHwithProof}, let $j(x)=\inf\big\{i>0\bigm|\Fc^i(x)\in\hat{M}_1\big\}$, $\frac{\sin\theta\|D\Fc^{j(x)(dx)}_x\|_{\p}}{\|dx\|_{\p}}>c_3$ and $\frac{\|D\Fc^{j(x)}_x(dx)\|}{\|dx\|}>c_3$.
    \item\label{item02ExpansionBetweenMhatminusMhatplus} $\exists$ constant $c_4>0$ such that for all nonsingular $(\phi,\theta)=x\in\hat{M}_1$, $dx\in \hat{C}^u_{1,x}$ as defined in \cref{thm:UniformExpansionOnMhat1}, let $j_1(x)=\inf\big\{i>0\bigm|\Fc^i(x)\in\hat{M}\big\}$, $\frac{\|D\Fc^{j_1(x)}_x(dx)\|_{\p}}{\|dx\|_{\p}}>c_4$ and $\frac{\|D\Fc^{j_1(x)}_x(dx)\|}{\|dx\|}>c_4$.
\end{enumerate}
\end{corollary}
\begin{proof}
    For the first conclusion, note that the trajectory segment $x\in\hat{M}_1,\Fc(x),\cdots,\Fc^{j_1(x)}(x)\in\hat{M}_1$ contains $x_0\in\Mrout$, $x_1=\Fc(x_0)\in \MRin$ and $x_2=\Fc^{n_1+1}(x_1)\in\Mrin$ with $n_1=\inf\big\{k\ge0:p(\Fc^i(x_1)),i=0,\cdots,k\in\Gamma_R\big\}$ similar to that in \cref{def:MhatReturnOrbitSegment}. Thus, the same analysis for \cref{TEcase_a0,TEcase_a1,TEcase_b,TEcase_c,corollary:sinpmetricexpansion,Corollary:EuclideanExpansionMhat} yields the first conclusion, 

    For the second conclusion, if the trajectory segment $x\in\hat{M}_1,\Fc(x),\cdots,\Fc^{j(x)}(x)\in\hat{M}$ are all collisions on $\Gamma_r$, then since $\hat{C}^u_{1,x}$ is a subset of positive quadrant $\mathcal{Q}_x(\mathrm{I,III})$ and $D\Fc^{j(x)}_x=\begin{pmatrix}
1 & 2j(x) \\
0 & 1 
\end{pmatrix}$ in $\phi\theta$ coordinate, $\frac{\|D\Fc^{j_1(x)}_x\|_{\p}}{\|dx\|_{\p}}>1$ and $\frac{\|D\Fc^{j_1(x)}_x\|}{\|dx\|}>1$. Otherwise, $x\in\hat{M}$, $\Fc(x)$, $\cdots$, $\Fc^{j(x)}(x)\in\hat{M}$ contains $x_0\in\Mrout$, $x_1=\Fc(x_0)\in \MRin$ and $x_2=\Fc^{n_1+1}(x_1)\in\Mrin$ with $n_1=\max\big\{k\ge0:p(\Fc^i(x_1))\in\Gamma_R\text{ for } i=0,\cdots,k\big\}$, then the same analysis for \cref{TEcase_a0,TEcase_a1,TEcase_b,TEcase_c,RMASUHwithProof,corollary:sinpmetricexpansion,Corollary:EuclideanExpansionMhat} yields the second conclusion. 
\end{proof}

\begin{theorem}[Uniform expansion on \Mrout]\label{thm:MroutreturnUniformExpansion}
    For the lemon billiards $\LEM(r,R,\phistar)$ satisfying \cref{def:UniformHyperbolicLemonBilliard} and for the return map $\hat{F}_1$ on, the return map $\tilde{F}$ on $\tilde{M}=\Mrout$ satisfies the following.
    \begin{enumerate}[label={(\roman*)}]
    \item\label{item01:MroutreturnUniformExpansion} For nonsingular $x\in\Mrout$, if $dx\in \tilde{C}^{u}_x\dfn\begin{cases}
            \big\{(d\phi,d\theta)\in T_x\Mrout\bigm|\frac{d\theta}{d\phi}\in[0,1]\big\}\text{, if } x\in \Mrout\smallsetminus\Nout,\\
            \big\{(d\phi,d\theta)\in T_x\Mrout\bigm|\frac{d\theta}{d\phi}\in[0,\frac{1}{3}]\big\} \text{, if } x\in\Nout, 
        \end{cases}$ \newline then $D\tilde{F}_x(dx)\in\big\{\text{interior of }\tilde{C}^{u}_{\tilde{F}(x)}\big\}$.
    \item\label{item02:MroutreturnUniformExpansion} For nonsingular $x\in\Mrout$, $dx\in \tilde{C}^{u}_x$, $\exists \tilde{c}>0$, $\tilde{\Lambda}>1$ such that $\frac{\|D\tilde{F}^n_x(dx)\|_\p}{\|dx\|_\p}>\tilde{c}\tilde{\Lambda}^n$, $\forall n\ge1$.
    \end{enumerate}
\end{theorem}
\begin{proof}
\textbf{Proof of \ref{item01:MroutreturnUniformExpansion}}: 
For any nonsingular $y\in\Mrout$, the $y$ return to $\Mrout$ orbit segment is defined as follows. \begin{equation}\label{eq:yMRoutreturnorbitsegment}
\underbracket{y=x_0}_{\mathclap{\in\Mrout}},\underbracket{\mathcal{F}(y)=x_1,\cdots}_{\notin\Mrout},\underbracket{\mathcal{F}^{\tilde{\sigma}(y)}(y)=\tilde{F}(y)\nfd\tilde{y}_1}_{\in \Mrout},\text{ where }\tilde{\sigma}(y)\dfn\inf\{k>0\mid \mathcal{F}^{k}(y)\in\Mrout\}.
\end{equation} and similar as \eqref{eqPtsFromeqMhatOrbSeg}, \eqref{eq:yMRoutreturnorbitsegment} also contains $x_1=\Fc(y)=\Fc(x_0)\in\MRin$ and $x_2=\Fc^{n_1+1}(x_1)\in\Mrin$. In \eqref{eq:yMRoutreturnorbitsegment} there are also the same (a0)(a1)(b)(c) cases in \eqref{eqMainCases} based on $x_1=(\Phi_1,\theta_1)$ and $n_1$. \cref{corollary:InvQuadrant,corollary:halfquadrantcone}, \eqref{eq:strictlyinvcone_casea1bc} and the same argument as in \cref{thm:UniformExpansionOnMhat} imply the strict invariance of $\tilde{C}^{u}_x$.

\textbf{Proof of \ref{item02:MroutreturnUniformExpansion}}: 
In \eqref{eq:yMRoutreturnorbitsegment}, there is a $\hat{y}_1\in\hat{M}$ with 
\begin{equation}\label{eq:hatydefinitionInOrbit}
\begin{aligned}
    \hat{y}_1&=\left\{
    \begin{aligned}
    &\tilde{y}_1,\text{ if }\tilde{y}_1\in\Mrin\cap\Mrout\\ 
    &\Fc^{-1}(\tilde{y}_1).\text{ if }\tilde{y}_1\in\Mrout\smallsetminus\Mrin.
    \end{aligned}\right.
\end{aligned}
\end{equation}
Suppose $T_{\hat{y}_1}M\ni\hat{y}_1=\Fc^{t_1}(y)$ with $t_1\in\mathbb{N}$. Then for $dy\in \tilde{C}^u_{y}$, let $\hat{y}_1\dfn D\Fc_y^{t_1}(dy)$ and $d\tilde{y}_1\dfn D\tilde{F}_y(dy)$.

For all cases (a0)(a1)(b)(c) of \eqref{eq:yMRoutreturnorbitsegment}, the analysis as in \cref{lemma:ExpansionLowerBoundA0A1,theorem:UniformExpansionA0,theorem:proofofTEcase_b,theorem:contractioncontrolA1,lemma:ExpansionLowerBoundB,lemma:AformularForLowerBoundOfExpansionForn1largerthan0,proposition:casecexpansionn1eq1,proposition:casecexpansionn1eq2,proposition:casecexpansionn1eq3,proposition:casecexpansionn1ge4} implies that $\frac{\|d\hat{y}_1\|_{\p}}{\|dy\|_{\p}}>0.05$. 

The $n$ times $\tilde{F}$ return orbit segment for $y$ is \[\underbracket{y\dfn\tilde{y}_0}_{\in\Mrout},\cdots,\underbracket{\tilde{y}_1\dfn\tilde{F}(y)}_{\in\Mrout},\cdots,\underbracket{\tilde{y}_n\dfn\tilde{F}^n(y)}_{\in\Mrout}\] and it contains $\underbracket{\hat{y}_1}_{\in\hat{M}},\cdots,\underbracket{\hat{y}_n}_{\in\hat{M}}$ with each $\hat{y}_k$ defined similarly as in \eqref{eq:hatydefinitionInOrbit}. \begin{align*}
    \hat{y}_k&=\left\{
    \begin{aligned}
    &\tilde{y}_k,\text{ if }\tilde{y}_k\in\Mrin\cap\Mrout\\ 
    &\Fc^{-1}(\tilde{y}_k).\text{ if }\tilde{y}_k\in\Mrout\smallsetminus\Mrin.
    \end{aligned}\right.
\end{align*}

With $d\tilde{y}_k=D\tilde{F}^{k}_y(dy)$, $T_{\hat{y}_k}M\ni d\hat{y}_k=D\Fc^{t_k}_{\tilde{y}_{k-1}}(d\tilde{y}_{k-1})$, $t_k>0$, $k=1,\cdots,n$. 

By \cref{corollary:halfquadrantcone}, $d\hat{y}_n\in C_{\hat{y}_n}=\mathcal{HQ}_{\hat{y}_n}(\mathrm{I,III})$. And since either $d\hat{y}_n=d\tilde{y}_n$ or $D\Fc_{\hat{y}_n}(d\hat{y}_n)=d\tilde{y}_n$ with $D\Fc_{\hat{y}_n}=\begin{pmatrix}
1 & 2 \\
0 & 1 
\end{pmatrix}$, $\|d\tilde{y}_n\|_{\p}\ge \|d\hat{y}_n\|_{\p}$ by \cref{proposition:ShearOneStepExpansion}. We get \begin{align*}
    \frac{\|D\tilde{F}^n_y(dy)\|}{\|dy\|_{\p}}=&\overbracket{\frac{\|d\tilde{y}_n\|_{\p}}{\|d\hat{y}_n\|_{\p}}}^{\ge1}\cdot\qquad\overbracket{\frac{\|d\hat{y}_n\|_{\p}}{\|d\hat{y}_1\|_{\p}}}^{\mathclap{\text{\cref{thm:UniformExpansionOnMhat}: }>c\Lambda^{n-1}}}\qquad\cdot\overbracket{\frac{\|d\hat{y}_1\|_{\p}}{\|dy\|_{\p}}}^{>0.05}\\
    >&1\cdot c\Lambda^{n-1}\cdot0.05\dfn\tilde{c}\tilde{\Lambda}^n\qedhere
\end{align*}
\end{proof}
\begin{corollary}\label{corollary:MroutreturnUniformExpansion}In the context of \cref{thm:MroutreturnUniformExpansion} the following hold. Under the return map $\tilde{F}$ on $\Mrout$ the vectors in the unstable cone $\tilde{C}_x^u$ defined in \cref{thm:MroutreturnUniformExpansion} expand uniformly with respect to the Euclidean metric $\|\,\|$. i.e., For nonsingular $(\phi,\theta)=x\in\MRout$, $dx=(d\phi,d\theta)\in \tilde{C}^u_x$ and all $n\ge0$, there exist $\tilde{c}_0,\tilde{c}_1>0,\tilde{\Lambda}_0>1$ such that \begin{align*}
        \dfrac{\sin\theta\|D\tilde{F}_x^n(dx)\|_{\p}}{\|dx\|_{\p}}&>\tilde{c}_0\tilde{\Lambda}_0^n,\\
        \dfrac{\|D\tilde{F}_x^n(dx)\|}{\|dx\|}&>\tilde{c}_1\tilde{\Lambda}_0^n.
    \end{align*}
\end{corollary}
\begin{proof}
    The proof is the same as \cref{corollary:sinpmetricexpansionMhat,Col:EuclideanExpansionMhat} by using the estimate in \eqref{eq:casea0onestepsinexpansion}.
\end{proof}

\begin{lemma}[symmetries between $\hat{M}_1$ and $\hat{M}$ from \cref{RMASUHwithProof} (see \cref{fig:hatUplusUminus}) ]\label{lemma:symmetriesbetweenMhatMhat1} Both symmetries $I$ and $J$ from \eqref{eq:InversionSymmetryDef} and  \eqref{eq:FoldSymmetryDef} are conjugates between $\hat{F}$ on $\hat{M}$ and $\hat{F}^{-1}$ on $\hat{M}_1$. i.e. $I\circ \hat{F}_1=\hat{F}^{-1}\circ I$, $J\circ \hat{F}_1=\hat{F}^{-1}\circ J$.
\end{lemma}
\begin{proof}
    Note that both $I$, $J$ are isometries between $\hat{M}$ and $\hat{M}_1$ since $J\Mrin=\Mrout=I\Mrin$ and $J\Mout_{r,i}=\Min_{r,i}=I\Mout_{r,i}$, $\forall i\ge0$, therefore $I(\hat{M})=\hat{M}_1=J(\hat{M})$ (also see \cref{fig:MrN}). 
    
    $\forall$ nonsingular point $x\in\hat{M}_1$, let $y=I(x)\in\hat{M}$ and then $x=I(y)$. Since $I$, $J$ are conjugates between $\Fc$ and $\Fc^{-1}$ and $I\circ I=Id=J\circ J$. 
    
    Therefore, $\forall k>0$, $\Fc^k(x)=\Fc^k(I(y))\overbracket{=}^{\mathclap{\eqref{eq:InversionSymmetryDef}}}I\Fc^{-k}(y)$, and thus \begin{equation}\label{eq:F1hatFhatConjugate}
        \Fc^{k}(x)\in\hat{M}_1\text{ if only if }\Fc^{-k}(y)\in I^{-1}(\hat{M}_1)=\hat{M}. 
    \end{equation}
    Note that $\sigma_1(x)\ge2$ and by \eqref{eqhatF1returnorbit}: $\sigma_1(x)=\inf\{k>0\mid\Fc^k(x)\in\hat{M}_1\}$. This means that for $k=1,\cdots,\sigma_1(x)-1$, $\Fc^{k}(x)\notin\hat{M}_1$ and $\Fc^{\sigma_1(x)}(x)\in\hat{M}_1$.
    
    Therefore, by \eqref{eq:F1hatFhatConjugate} for $k=1,\cdots,\sigma_1(x)-1$, $\Fc^{-k}(y)\notin\hat{M}$ and $\Fc^{-\sigma_1(x)}(y)\in\hat{M}$. With $y\in\hat{M}$ and by definition \eqref{eq:F1hatFhatConjugate}, we see $\sigma_1(x)=\sigma(\hat{F}^{-1}(y))$. We further get the following. \begin{align*}I\hat{F}_1(x)&\overbracket{=}^{\mathclap{\eqref{eq:F1hatFhatConjugate}}}I\circ \Fc^{\sigma_1(x)}(x)\overbracket{=}^{\mathclap{\eqref{eq:InversionSymmetryDef}}}\Fc^{-\sigma_1(x)}\circ I(x)\overbracket{=}^{\mathclap{I(x)=y}}\Fc^{-\sigma_1(x)}(y)\overbracket{=}^{\mathclap{\sigma_1(x)=\sigma(\hat{F}^{-1}(y))}}\hat{F}^{-1}(y)\underbracket{=}_{\mathclap{y=I(x)}}\hat{F}^{-1}I(x)\\
    \end{align*}
    
    Let $z=J(x)$, by the exactly same preceding reasoning and using \eqref{eq:FoldSymmetryDef}, we get $J\hat{F}_1(x)=\hat{F}^{-1}J(x)$.
\end{proof}

Using symmetries $I$ and $J$, we can also prove that the stable cone on $\hat{M}$ is strictly invariant under $\hat{F}^{-1}$.

\section{Local ergodicity conditions}\label{sec:LET}
\subsection[Notations for singularity curves and local ergodic conditions L1, L2]{Notations for singularity curves and local ergodic conditions L1, L2}\label{subsec:NotationL1L2}

\begin{notation}\label{notation:symplecticformmeasurenbhd}[{\cite[Definition 2.1, 2.3 and 2.8]{MR3082535}}] From here on, we introduce and use notations from \cite{{MR3082535}} and make standing assumption for lemon billiard configuration. 
\begin{itemize}
    \item $(M,\omega,\mu,\Fc)$ denotes our billiard map $\Fc$ on phase space $M$ with the invariant symplectic form $\omega=\sin\theta d\phi \wedge d\theta$ and $\mu$ is its induced Liouville measure on $M$.\COMMENT{This is redundant with prior notation, which makes us look stupid.}
    \item $d(x,y)$ is the Euclidean distance between points $x,y\in M$.
    \item Given a subset $\mathcal{A}\subset M$ and $\varepsilon>0$, let $\mathcal{A}(\varepsilon)=\big\{y\in M:d(y,\mathcal{A})<\varepsilon\big\}$ be the $\varepsilon$-neighborhood of $\mathcal{A}$.
    \item We also assume that our billiard is a lemon billiard $\mathbb{L}(r,R,\phistar)$ that satisfies \cref{def:UniformHyperbolicLemonBilliard}.
    \item With singularity curves notations from \cref{def:BasicNotations0}, we further define $\mathcal{S}^{-}_{1}\dfn\mathcal{S}_{-1}\smallsetminus{\partial M}$ and $\mathcal{S}^{+}_{1}\dfn\mathcal{S}_{1}\smallsetminus{\partial M}$. Similarly, let $\mathcal{S}^{-}_{k}=\mathcal{S}_{-k}\smallsetminus\mathcal{S}_{-(k-1)}$ and let $\mathcal{S}^{+}_{k}=\mathcal{S}_{k}\smallsetminus\mathcal{S}_{k-1}$ for $k\ge2$.
    \item Let $\mathcal{L}_{-}(\cdot)$ be the Lebesgue measure restricted on the line segments $\mathcal{S}^{-}_{1}$. Let $\mathcal{L}_{+}(\cdot)$ be the Lebesgue measure restricted on the line segments $\mathcal{S}^{+}_{1}$.
\end{itemize}
\begin{definition}[Local ergodic conditions L1, L2]\hfill

    \begin{itemize}
        \item L1 (\cite[Theorem 4.1]{MR3082535}). The sets $\mathcal{S}_k$ and $\mathcal{S}_{-k}$ are regular (\cref{def:RegularSingularityCurve}) for every $k>0$,

        \item L2 (\cite[Theorem 4.1]{MR3082535}). For every $k\ge1$: 
    \end{itemize}
\end{definition}
\end{notation}

\subsection{Cones, quadratic forms and sufficient points}
\begin{definition}[{Cones and quadratic forms \cite[Definition 2.11, 2.12, 2.14, 2.16]{MR3082535}}]\label{def:QuadraticFormForCone}\label{def:ConeWithQuadraticForm}\label{def:ExpansionInCone}
    Consider two transverse Lagrangian subspaces $A_x$, $B_x$ for all $x=(\phi,\theta)\in U$ in an open set $U$. For $u\in T_xM$ define the \emph{quadratic forms}  
    \begin{align*}
         Q_x(u)&\dfn\omega_x(u_1,u_2),\\
         Q'_x(u)&\dfn-\omega_x(u_1,u_2),
    \end{align*}where $u_1\in A_x$, $u_2\in B_x$ are uniquely defined by $u=u_1+u_2$. For every $x\in U$, the \emph{cone field} $\mathcal{C}=\big\{\mathcal{C}(x)\big\}_{x\in U}$ and its complement $\mathcal{C}'$ associated to $A$ and $B$ is the family of closed cones given by \begin{align*}
    \mathcal{C}(x)&=Q_x^{-1}([0,+\infty))\subset T_x(M),\\
    \mathcal{C}'(x)&=Q_x^{-1}((-\infty,0)).
\end{align*} For every $u\in\mathcal{C}(x)$, the quadratic form $Q$ induces a norm on cone $\|u\|_{Q}\dfn\sqrt{Q_x(u)}$.

If $x \mapsto A_x$, $x\mapsto B_x$ are continuous, then we say that $Q$ and $C$ are continuous.
    \begin{itemize}
        \item $Q$ is \emph{monotone} (with respect to $\Fc$) if $Q_{\Fc^k(x)}(D\Fc^k_xu)\ge Q_x(u)$ whenever $u\in T_xM$, $k>0$, $x,\Fc^k(x)\in U$.
    \item $\mathcal{C}$ is \emph{invariant} (with respect to $\Fc$) if $D\Fc^k_x(\mathcal{C}(x))\subset\mathcal{C}(\Fc^k(x))$ whenever $k>0$, $x,\Fc^k(x)\in U$.
\end{itemize}
Whenever $k>0$, $x,\Fc^k(x)\in U$, write
\begin{equation}\label{eq:ExpansionInCone}
        \sigma_{\mathcal{C}}(D\Fc_x^k)\dfn\inf_{u\in\mathcal{C}(x)}\sqrt{\frac{Q_{\Fc^k(x)}(D\Fc_x^k(u))}{Q_x(u)}}\quad\text{and}\quad
        \sigma^{*}_{\mathcal{C}}(D\Fc_x^k)\dfn\inf_{u\in\mathcal{C}(x)}\frac{\sqrt{Q_{\Fc^k(x)}(D\Fc_x^k(u))}}{\|u\|}.
    \end{equation}

Similarly, for every $u'\in\mathcal{C}'(x)$, the quadratic form $Q'$ induces a norm on cone $\mathcal{C}'$, $\|u'\|_{Q'}=\sqrt{Q'_x(u')}$.
\end{definition}

\begin{remark}[{\cite[Proposition 2.13 and 2.15, Remark 2.17]{MR3082535}}]\label{remark:MonotoneCone} $\mathcal{C}$ is invariant if and only if the quadratic form $Q$ is monotone \cite[Theorem 4.4]{MR1346498} if and only if $D\Fc_{x}^{-k}C'(x)\subset \mathcal{C}'(\Fc^{-k}(x))$ for every $x\in U$ and $k>0$ such that $\Fc^{-k}(x)\in U$ \cite[Proposition 6.2]{MR1346498}. $D\Fc^k_x(\mathcal{C}(x))\subset\mathcal{C}(\Fc^k(x))$ implies  $\sigma_{\mathcal{C}}(D\Fc_x^k)\ge 1$ \cite[Proposition 6.1]{MR1346498}.\COMMENT{Check!\wentao{looks OK. Just too squeezed}}

Furthermore, if $D\Fc_{x}^k(\mathcal{C}(x))\subset\Int\mathcal{C}(\Fc^k(x))$ (strict invariance), then $\sigma_{\mathcal{C}}(D\Fc_x^k)> 1$.
\end{remark}

\begin{definition}[{Joint invariance \cite[Definition 2.18]{MR3082535}}]\label{def:JointInvariant}
Cone fields $\mathcal{C}_0,\mathcal{C}_1$ on open sets $U_0,U_1$, respectively, are said to be \emph{jointly invariant} if $D\Fc^k_x(\mathcal{C}_i)\subset\mathcal{C}_{1-i}(\Fc^kx)$ for every $x\in U_i$ and $k>0$ such that $\Fc^kx\in U_{1-i}$.
\end{definition}
\begin{remark}
    The cone $C_x$ on $\hat{M}$ is jointly invariant with $\hat{C}^u_{1,x}$ on $\hat{M}_1$.
\end{remark}
\begin{definition}[Sufficient point {\cite[Definition 3.1]{MR3082535}}]\label{def:sufficientPoint}
A point $x\in M\smallsetminus\partial M$ is said to be sufficient if there exist:
\begin{enumerate}[label={(\roman*)}]
\item an integer $l$ such that $\Fc^l$ is a local diffeomorphism at $x$;
\item a neighborhood $U$ of $\Fc^l x$ and an integer $N>0$ such that $U\cap\mathcal{S}_{-N}=\emptyset$ ($\mathcal{S}_{-N}$ from singularity point defintion \cref{def:BasicNotations0});
\item an invariant continuous cone field $\mathcal{C}$ on $U\cup\Fc^{-N}(U)$ such that $\sigma_{\mathcal{C}}(D\Fc^N_y)>3$ for every $y\in\Fc^{-N}U$.
\end{enumerate} We say that $x$ is a sufficient point with quadruple $(l,N,U,\mathcal{C})$ (from \cref{def:sufficientPoint}).
\end{definition}

\subsection[Lemon billiard cones, sufficient points on reference sets]{Lemon billiard cones, sufficient points on reference sets}\strut\COMMENT{Replace $\hat M_1$ by $\check M$?}
\begin{notation}\label{def:UIntRegionofD} Define the open set $\mathfrak{U}\dfn\big\{(\phi,\theta)\in \Int(M_R) \bigm |\sin{\theta} >\sqrt{4r/R} \big\}\subset\MRout\cap\MRin$ (interior of the orange rectangle in \cref{fig:A_compactRegion}) and $\mathfrak{U}_{\pm}\dfn\Fc^{\pm1}(\mathfrak{U})$, $\hat{\mathfrak{U}}_{-}\dfn(\Mrin\cap\Mrout)\sqcup\Fc^{-1}(\mathfrak{U}_{-}\smallsetminus\Mrin)$, and $\hat{\mathfrak{U}}_{+}\dfn(\Mrin\cap\Mrout)\sqcup\Fc(\mathfrak{U}_{+}\smallsetminus\Mrout)$ (see \cref{fig:hatUplusUminus}). 
\end{notation}

\begin{remark}\label{remark:UIntRegionofD}\hfill
    \begin{itemize}
        \item  $I\mathfrak{U}=\mathfrak{U}$ since $\sin\theta>\sqrt{4r/R}$ is equivalent to $\sin{(\pi-\theta)>\sqrt{4r/R}}$.
        \item $\mathfrak{U}_{-}$ is the set of $x_0$, $\mathfrak{U}_{+}$ is the set of $x_2$ in \cref{def:fpclf} and \eqref{eqPtsFromeqMhatOrbSeg} for case (a0) of \eqref{eqMainCases}. By \cref{contraction_region}, $(\Mrout\smallsetminus \Nout)\subset\mathfrak{U}_{-}\subset\Mrout$ and $(\Mrin\smallsetminus \Nin)\subset\mathfrak{U}_{+}\subset\Mrin$ (see \cref{fig:MrN}).
        \item By \cref{remark:GaranteeLargeTheta,lemma:Trajectory_approximation_1} with $R_{HF}(r,\phistar)$ as in \eqref{eqJZHypCond}, we have $\Nin\subset\bigcup_{k\in\mathbb{Z}}(\Fc^k(\mathfrak{U}))$. So $\Nout=I\Nin\subset \bigcup_{k\in\mathbb{Z}}(I\Fc^k(\mathfrak{U}))\overbracket{=}^{\mathclap{\eqref{eq:InversionSymmetryDef}}}\bigcup_{k\in\mathbb{Z}}(\Fc^{-k}(I\mathfrak{U}))\overbracket{=}^{\mathclap{I\mathfrak{U}=\mathfrak{U}}}\bigcup_{k\in\mathbb{Z}}(\Fc^{k}(\mathfrak{U}))$. This means $\Mrout\subset\bigcup_{k\in\mathbb{Z}}(\Fc^{k}(\mathfrak{U}))$ and $\Mrin\subset\bigcup_{k\in\mathbb{Z}}(\Fc^{k}(\mathfrak{U}))$, hence $M\aeq \bigcup_{k\in\mathbb{Z}}(\Fc^{k}(\mathfrak{U}))\aeq \bigcup_{k\in\mathbb{Z}}(\Fc^{k}(\hat{\mathfrak{U}}_{-}))\aeq \bigcup_{k\in\mathbb{Z}}(\Fc^{k}(\hat{\mathfrak{U}}_{+}))$.
    \end{itemize}
\end{remark}

\begin{definition}[Lemon billiard cones and quadratic norm]\label{def:billiardQuadrantConeAndForm}
With $A_x=\frac{\partial}{\partial\phi}$, $B_x=\frac{\partial}{\partial\theta}$, and $\omega_(x)=\sin\theta d\phi\wedge d\theta$ we obtain $\|a\frac{\partial}{\partial\phi}+b\frac{\partial}{\partial\theta}\|^2_{Q}=Q_x(a\frac{\partial}{\partial\phi}+b\frac{\partial}{\partial\theta})=\omega_x(a\frac{\partial}{\partial\phi},b\frac{\partial}{\partial \phi})=ab\sin\theta$ (\cref{def:QuadraticFormForCone}). We write $\|dx\|^2_{Q}=d\phi d\theta\sin\theta$ and set\begin{equation}\label{eq:PositiveNegativeConeOnLemonSets}\begin{aligned}
        \mathcal{C}(x)=&Q_x^{-1}([0,+\infty))=\big\{dx=(d\phi,d\theta)\in T_xM\bigm|d\phi d\theta\ge0\big\}=\big\{dx=(d\phi,d\theta)\in T_xM\bigm|\frac{d\theta}{d\phi}\in[0,+\infty]\big\},\\
        \mathcal{C}'(x)=&Q_x^{-1}((-\infty,0 ))=\big\{dx=(d\phi,d\theta)\in T_xM\bigm|\frac{d\theta}{d\phi}\in(-\infty,0)\big\}.
    \end{aligned}\end{equation}
\end{definition}

\begin{figure}[ht]
\begin{center}
    \begin{tikzpicture}[xscale=1.0,yscale=1.0]
        \pgfmathsetmacro{\PHISTAR}{0.68}
        \pgfmathsetmacro{\WIDTH}{7}
        \pgfmathsetmacro{\LENGTH}{7}
        \tkzDefPoint(-\LENGTH,0.5*\WIDTH){LLM};
        \tkzDefPoint(-\LENGTH,0){LLL};
        \tkzDefPoint(-\LENGTH+\PHISTAR,0){LL};
        \tkzDefPoint(-\LENGTH+\PHISTAR,\WIDTH){LU};
        \tkzDefPoint(-\LENGTH,\WIDTH){LLU};
        \tkzDefPoint(\LENGTH-\PHISTAR,0){RL};
        \tkzDefPoint(\LENGTH-\PHISTAR,\WIDTH){RU};
        \tkzDefPoint(\LENGTH,0){RRL};
        \tkzDefPoint(\LENGTH,\WIDTH){RRU};
        \tkzDefPoint(-\LENGTH+\PHISTAR,\PHISTAR){X};
        \tkzDefPoint(-\LENGTH+\PHISTAR,\WIDTH-\PHISTAR){IX};
        \tkzDefPoint(\LENGTH-\PHISTAR,\PHISTAR){IY};
        \tkzDefPoint(\LENGTH-\PHISTAR,\WIDTH-\PHISTAR){Y};
        \tkzDefPoint(0,0.5*\WIDTH+0.5*\PHISTAR){UMrin01};
        \tkzDefPoint(\PHISTAR,0.5*\WIDTH){RightMrin01};
        \tkzDefPoint(0,0.5*\WIDTH-0.5*\PHISTAR){LMrin01};
        \tkzDefPoint(-1.0*\PHISTAR,0.5*\WIDTH){LeftMrin01};
        \fill [red, opacity=10/30](UMrin01) -- (RightMrin01) -- (LMrin01) -- (LeftMrin01) -- cycle;
        \fill [blue, opacity=10/30](UMrin01) -- (RightMrin01) -- (LMrin01) -- (LeftMrin01) -- cycle;
        \tkzDrawPoints(UMrin01,RightMrin01,LMrin01,LeftMrin01);
        \tkzLabelPoint[left](LeftMrin01){K1};
        \tkzLabelPoint[below](LMrin01){K2};
        \tkzLabelPoint[right](RightMrin01){K3};
        \tkzLabelPoint[above](UMrin01){K4};
 
        \tkzDefPoint(\LENGTH-\PHISTAR,0.5*\WIDTH){RUhatUplus};
        \tkzDefPoint(\LENGTH-\PHISTAR,0.5*\WIDTH-0.5*\PHISTAR){RLhatUplus};
        \tkzDefPoint(-1.0*\LENGTH+3.0*\PHISTAR,\PHISTAR){LUhatUplus};
        \fill [blue, opacity=10/30](RUhatUplus) -- (RLhatUplus) -- (LL) -- (LUhatUplus) -- cycle;
        \begin{scope}
            \clip (LL) -- (LUhatUplus) -- (RUhatUplus) -- cycle;
            \draw[fill=blue] circle[at=(LUhatUplus),radius=0.15];
        \end{scope}
        \node[] at ([shift={(0.24,-0.33)}]LUhatUplus) {\small $\Fc(N^{\text{\upshape in}})$};
        \node[] at ([shift={(-0.03,0.2)}]LUhatUplus) {\tiny $LF=(3\phistar,\phistar)$};

        \tkzDefPoint(-\LENGTH+\PHISTAR,0.5*\WIDTH){LLhatUminus};
        \tkzDefPoint(-\LENGTH+\PHISTAR,0.5*\WIDTH+0.5*\PHISTAR){LUhatUminus};
        \tkzDefPoint(1.0*\LENGTH-3.0*\PHISTAR,\WIDTH-\PHISTAR){RLhatUminus};
        \fill [blue, opacity=10/30](LLhatUminus) -- (LUhatUminus) -- (RU) -- (RLhatUminus) -- cycle;
        \begin{scope}
            \clip (RU) -- (RLhatUminus) -- (LLhatUminus) -- cycle;
            \draw[fill=blue] circle[at=(RLhatUminus),radius=0.15];
        \end{scope}
         \begin{scope}
            \clip (RU) -- (RLhatUminus) -- (LLhatUminus) -- cycle;
            \draw[fill=blue] circle[at=(RLhatUminus),radius=0.15];
        \end{scope}
         \node[] at ([shift={(0.1,-0.13)}]RLhatUminus) {\small $\Fc(N^{\text{\upshape in}})$};

        \tkzLabelPoint[below](LL){$LL=(\phistar,0)$};
        \tkzLabelPoint[below]([shift={(-0.5,0)}]RL){$2\pi-\phistar$};
        \tkzLabelPoint[left](LLL){0};
        \tkzLabelPoint[left](LLM){$\frac{\pi}{2}$};
        \tkzLabelPoint[left](LLU){$\pi$};
        \tkzLabelPoint[below](RRL){$2\pi$};
        \tkzLabelPoint[right](IY){$x_*$};
        \tkzLabelPoint[left](IX){$y_*$};
        \tkzLabelPoint[right](Y){$I x_*$};
        \tkzLabelPoint[left](X){$IY=I y_*$};
        
        \draw [thick] (LL) --(LU);
        \draw [thick] (LU) --(RU);
        \draw [thick] (RL) --(RU);
        \draw [thick] (LL) --(RL);
        \draw [thin,dashed](LLL)--(LL);
        \draw [thin,dashed](LLU)--(LU);
        \draw [thin,dashed](LLU)--(LLL);
        \draw [thin,dashed](RRL)--(RL);
        \draw [thin,dashed](RRU)--(RU);
        \draw [thin,dashed](RRU)--(RRL);
        \draw [thin,dashed](LL)--(Y);
        \draw [thin,dashed](RU)--(X);
        \draw [thin,dashed](LU) -- (IY);
        \draw [thin,dashed](RL) -- (IX);
        \tkzDefPoint(-\LENGTH+\PHISTAR,\WIDTH/2-\PHISTAR/4-\PHISTAR/4){M1LL};
        \tkzDefPoint(-\LENGTH+\PHISTAR,\WIDTH/2-\PHISTAR/4+\PHISTAR/4){M1LU};
        \tkzDefLine[parallel=through M1LU](M1LL,RL)\tkzGetPoint{M11RU};
        \path[name path=lineM1U] (M11RU) -- (M1LU);
        \path[name path=lineM0L] (IX) -- (RL);
        \path[name intersections={of=lineM1U and lineM0L,by={M11RU}}];
        \fill [red, opacity=10/30](RL) -- (M1LL) -- (M1LU) -- (M11RU) -- cycle;
        \begin{scope}
            \clip (RL) -- (M11RU) -- (M1LU) -- cycle;
            \draw[fill=red] circle[at=(M11RU),radius=0.15];
        \end{scope}
        \node[] at ([shift={(-0.14,0.18)}]M11RU) {\small $\Fc^{-1}(N^{\text{\upshape out}})$};
        \path[name path=RBoundry](RL) -- (RU);
        \tkzDefLine[parallel=through LU](M1LL,RL)\tkzGetPoint{UM11RU};
        \path[name path=lineUM1U](LU)--(UM11RU);
        \path[name intersections={of=lineUM1U and RBoundry,by={UM1RU}}];
        \tkzDefShiftPoint[UM1RU](0,-\PHISTAR/2){UM1RL};
        \tkzDefLine[parallel=through UM1RL](UM1RU,LU)\tkzGetPoint{UM11LL};
        \path[name path=lineLUIY](LU) -- (IY);
        \path[name path=UM11LLUM1RL](UM11LL)--(UM1RL);
        \path[name intersections={of=lineLUIY and UM11LLUM1RL,by={UM1LL}}];
        \fill [red, opacity=10/30](LU) -- (UM1RU) -- (UM1RL) -- (UM1LL) -- cycle;
        \begin{scope}
            \clip (LU) -- (UM1LL) -- (UM1RL) -- cycle;
            \draw[fill=red] circle[at=(UM1LL),radius=0.15];
        \end{scope}
        \node[] at ([shift={(-0.12,-0.18)}]UM1LL) {\small $\Fc^{-1}(N^{\text{\upshape out}})$};
        \tkzDefPoint(-\LENGTH+\PHISTAR,\WIDTH/3){L_LUM2};
        \tkzDefShiftPoint[L_LUM2](0,-\PHISTAR/3){L_LLM2};
        \tkzDefLine[parallel=through L_LUM2](L_LLM2,RL)\tkzGetPoint{Far_L_RUM2};
        \path[name path=line_L_LUM2_Far_L_RUM2] (L_LUM2)--(Far_L_RUM2);
        \path[name path=line_RL_M1LL] (RL)--(M1LL);
        \path[name intersections={of=line_L_LUM2_Far_L_RUM2 and line_RL_M1LL,by={L_RUM2}}];
        \tkzDefPoint(\LENGTH-\PHISTAR,2*\WIDTH/3){R_RLM2};
        \tkzDefShiftPoint[R_RLM2](0,\PHISTAR/3){R_RUM2};
        \tkzDefLine[parallel=through R_RLM2](R_RUM2,LU)\tkzGetPoint{Far_R_LLM2};
        \path[name path=line_R_RLM2_Far_R_LLM2] (R_RLM2)--(Far_R_LLM2);
        \path[name path=line_LU_UM1RU] (LU)--(UM1RU);
        \path[name intersections={of=line_R_RLM2_Far_R_LLM2 and line_LU_UM1RU,by={R_LLM2}}];

        \tkzDefPoint(\LENGTH-\PHISTAR,1*\WIDTH/3){MR_RLM2};
        \tkzDefShiftPoint[MR_RLM2](0,\PHISTAR/3){MR_RUM2};
        \path[name path=line_IY_LU] (IY)--(LU);
        \tkzDefLine[parallel=through MR_RLM2](R_RLM2,R_LLM2)\tkzGetPoint{far_MR_LLM2};
        \tkzDefLine[parallel=through MR_RUM2](R_RLM2,R_LLM2)\tkzGetPoint{far_MR_LUM2};
        \path[name path=line_MR_RLM2_far_MR_LLM2] (far_MR_LLM2)--(MR_RLM2);
        \path[name path=line_MR_RUM2_far_MR_LUM2] (far_MR_LUM2)--(MR_RUM2);
        \path[name path=line_LU_IY] (LU)--(IY);
        \path[name intersections={of=line_LU_IY and line_MR_RLM2_far_MR_LLM2,by={MR_LLM2}}];
        \path[name intersections={of=line_LU_IY and line_MR_RUM2_far_MR_LUM2,by={MR_LUM2}}];

        \tkzDefPoint(-\LENGTH+\PHISTAR,2*\WIDTH/3){ML_LUM2};
        \tkzDefShiftPoint[ML_LUM2](0,-\PHISTAR/3){ML_LLM2};
        \tkzDefLine[parallel=through ML_LUM2](R_LLM2,R_RLM2)\tkzGetPoint{far_ML_RUM2};
        \tkzDefLine[parallel=through ML_LLM2](R_LLM2,R_RLM2)\tkzGetPoint{far_ML_RLM2};
        \path[name path=line_ML_LUM2_far_ML_RUM2] (ML_LUM2)--(far_ML_RUM2);
        \path[name path=line_ML_LLM2_far_ML_RLM2] (ML_LLM2)--(far_ML_RLM2);
        \path[name path=line_RL_IX] (RL)--(IX);
        \path[name intersections={of=line_RL_IX and line_ML_LUM2_far_ML_RUM2,by={ML_RUM2}}];
        \path[name intersections={of=line_RL_IX and line_ML_LLM2_far_ML_RLM2,by={ML_RLM2}}];
    \begin{scope}
      \clip (LL) -- (X) -- (RU) -- cycle;
      \draw[fill=green] circle[at=(X),radius=0.2];
    \end{scope}
    \begin{scope}
      \clip (RU) -- (Y) -- (LL) -- cycle;
      \draw[fill=green] circle[at=(Y),radius=0.2];
    \end{scope}
    \begin{scope}
      \clip (RL) -- (IX) -- (LU) -- cycle;
      \draw[fill=yellow] circle[at=(IX),radius=0.2];
    \end{scope}
    \begin{scope}
      \clip (LU) -- (IY) -- (RL) -- cycle;
      \draw[fill=yellow] circle[at=(IY),radius=0.2];
    \end{scope}
    
    \node[] at (0,0.5*\WIDTH)  {$M^{\text{\upshape in}}_{r,0}$};
    \node[] at (0,0.25*\WIDTH)  {\tiny $\Fc(M^{\text{\upshape in}}_{r,2})$};
    \node[] at (0,0.75*\WIDTH)  {\tiny $\Fc(M^{\text{\upshape in}}_{r,2})$};
    
    \path[name path=line_RU_X] (RU)--(X);
    \path[name intersections={of=line_RU_X and lineM1U,by={Mr1Label}}];

    \path[name intersections={of=line_RU_X and line_L_LUM2_Far_L_RUM2,by={Mr2Label}}];
    \node[] at ([shift={(0.2,-0.05)}]X) {\tiny $N^{\text{\upshape in}}$};
    \node[] at ([shift={(-0.14,0.13)}]Y) {\tiny $N^{\text{\upshape in}}$};
    \node[] at ([shift={(0.25,0.1)}]IX) {\tiny $N^{\text{\upshape out}}$};
    \node[] at ([shift={(-0.08,-0.1)}]IY) {\tiny $N^{\text{\upshape out}}$};
    \tkzDrawPoints(LL,X,LUhatUplus);
    
    \end{tikzpicture}
    \caption{ The two strips with blue one representing $\hat{M}_{+}=\hat{M}_1=(\Mrout\cap\Mrin)\sqcup\Fc(\Mrin\smallsetminus\Mrout)$ and red one representing $\hat{M}_{-}=\hat{M}=(\Mrout\cap\Mrin)\sqcup\Fc^{-1}(\Mrout\smallsetminus\Mrin)$, with marked $\Fc(\Nin)$ and $\Fc^{-1}(\Nout)$. $(\hat{M}_1\smallsetminus\Fc{(\Nin)})\subset\hat{\mathfrak{U}}_{+}\subset \hat{M}_1=\hat{M}_{+}\,,\qquad (\hat{M}\smallsetminus\Fc^{-1}{(\Nout)})\subset\hat{\mathfrak{U}}_{-}\subset \hat{M}=\hat{M}_{-}$. \newline We consider the line components $IY-K1$, $LL-K3$ of $\mathcal{S}_{-1}$ and their $\xi$ neighborhoods in \cref{theorem:L4condition}.
    }\label{fig:hatUplusUminus}.
\end{center}
\end{figure}
\begin{theorem}[Sufficient points on $\hat{\mathfrak{U}}_{+}$]\label{theorem:sufficentptsonMhatplus} For every $x\in\hat{\mathfrak{U}}_{+}\smallsetminus\mathcal{S}_{-\infty}$, the followings hold.
    \begin{enumerate}[label={(\roman*)}]
        \item\label{item01sufficentptsonMhatplus}  
        $\exists$ some integer $l<0$ such that $\Fc^{l}$ is a local diffeomorphism at $x$ such that $\Fc^l(x)\in\hat{\mathfrak{U}}_{+}$. 
        \item\label{item02sufficentptsonMhatplus} And there is an open neighborhood $U\subset \hat{\mathfrak{U}}_{+}$ of $\Fc^l(x)$ and an integer $N>0$ such that $U\cap\mathcal{S}_{-N}=\emptyset$ and $\Fc^{-N}(U)\subset \hat{\mathfrak{U}}_{+}$. The cone field $\mathcal{C}$ given in \cref{def:billiardQuadrantConeAndForm} restricted on $U\cup\Fc^{-N}(U)$ is invariant and $\sigma_{\mathcal{C}}(D\Fc^N_y)>3$ for every $y\in\Fc^{-N}U$.
    \end{enumerate}
\end{theorem}
\begin{proof}
    \textbf{Proof of (i)}. Since $x\in\hat{\mathfrak{U}}_{+}\subset\hat{M}_1=(\Mrout\cap\Mrin)\sqcup\Fc(\Mrin\smallsetminus\Mrout)$, $x\in\Mrin$ or $\Fc^{-1}(x)\in\Mrin$. Hence, recall \cref{def:UIntRegionofD} the definition of $\mathfrak{U}=\big\{(\phi,\theta)\in \Int(M_R) \bigm |\sin{\theta} >\sqrt{4r/R} \big\}$ and \cref{remark:UIntRegionofD}, we can let 
    \begin{equation}\label{eq:x0forUplushatinversetreturn}
        \begin{aligned}
        \Mrout\ni x_0\dfn\left\{\begin{aligned}
        \Fc^{-2}(x)&\text{ if }x\in\Mrout\cap\Mrin\\
        \Fc^{-3}(x)&\text{ otherwise.}
    \end{aligned}\right.
    \end{aligned}
    \end{equation}
    Since $x_0\notin\mathcal{S}_{-\infty}$, similar to the $m(x)$ definition in \cref{def:SectionSetsDefs} we can define the following $m'(x_0)$ \[\exists\quad 0\le m'(x_0)\dfn\max\big\{m'\ge 0\mid p(\mathcal{F}^{-i}(x_0))\in\operatorname{int}\Gamma_r\text{ for }0\le i\le m'\big\}.\]  
    If $m'(x_0)=0$, i.e. $x_0\in\Mrin\cap\Mrout$, then we set $\Fc^{-l}(x)=x_0\in\Mrin\cap\Mrout\subset\hat{\mathcal{U}}_{+}$. If $m'(x_0)\ge1$, then $\Fc^{-m'(x_0)+1}(x_0)\in\Fc(\Mrin\smallsetminus\Mrout)$. If $\Fc^{-m'(x_0)+1}(x_0)\in\hat{\mathcal{U}}_{+}$, then we set $\Fc^{-l}(x)=\Fc^{-m'(x_0)+1}(x_0)$. Either $l=-m'(x_0)-1$ or $l=-m'(x_0)-2$ is determined by \eqref{eq:x0forUplushatinversetreturn}.
    
    Now, if $\Fc^{-m'(x_0)+1}(x_0)\notin\hat{\mathcal{U}}_{+}$, then in \cref{fig:hatUplusUminus}, $\Fc^{-m'(x_0)}(x_0)\in\Nin$. Hence $\exists$ $j(x_0)>0$ to be the smallest integer $i>0$ such that $\Fc^{-i}(\Fc^{-m'(x_0)}(x_0))\in\Nout$. Denote $x'_0\dfn \Fc^{-j}(\Fc^{-m'(x_0)}(x_0))$. Then $\Fc^{-m'(x'_0)+1}(x'_0)\in\Fc(\Mrin\smallsetminus\Mrout)$ and $\Fc^{-m'(x'_0)}(x'_0)\in\Mrin$. 
    
    Note that by \cref{lemma:Trajectory_approximation_1,remark:GaranteeLargeTheta} for our chosen $R>R_{HF}(r,\phistar)$, the orbit of $\Nin$ point cannot have two consecutive returns to $\Mrin$ being the returns to $\Nin$. We have $\Fc^{-m'(x'_0)}(x'_0)\in(\Mrin\smallsetminus\Nin)\subset\mathfrak{U}_{+}$. Therefore, $\Fc^{-m'(x'_0)+1}(x'_0)\in\Fc(\mathfrak{U}_{+}\smallsetminus\Mrout)\overbracket{\subset}^{\mathclap{\text{\cref{def:UIntRegionofD}}}}\hat{\mathfrak{U}}_{+}$. We set \[\Fc^{l}(x)=\overbracket{\Fc^{-m'(x'_0)+1}(x'_0)}^{\mathclap{\in\hat{\mathfrak{U}}_{+}}}=\Fc^{-m'(x'_0)+1}(\Fc^{-j}(\Fc^{-m'(x_0)}(x_0)))=\Fc^{-m'(x'_0)-m'(x_0)-j+1}(x_0)\in\hat{\mathfrak{U}}_{+},\] where $l=-m'(x'_0)-m'(x_0)-j(x_0)-1$ or $l=-m'(x'_0)-m'(x_0)-j(x_0)-2$ depending on \cref{eq:x0forUplushatinversetreturn}, so $x$. Then by \cite[equation (2.26) and theorem 2.33]{cb}, $\Fc^l$ is a local diffeomorphism between $x,\Fc^l(x)\in\mathfrak{U}_{+}$.

    \textbf{Proof of (ii)}. It is clear by \cref{caseA,caseB,caseC} the positive cone $\mathcal{C}$ given in \cref{def:billiardQuadrantConeAndForm} restricted on $U\cup\Fc^{-N}(U)$ is invariant for any $N>0$ and open $U\subset\hat{\mathfrak{U}}_{+}$. For any $N>0$,  as long as $U$ is small enough, $\Fc^{-N}$ is a diffeomorphism by \cite[equation (2.26) and theorem 2.33]{cb} on $U$ such that $U\cap\mathcal{S}_{-N}=\emptyset$.
    
    Since $\hat{M}_1\ni\Fc^l(x)\notin\mathcal{S}_{-\infty}$, there are infinite return times of $\Fc^l(x)$ on $\hat{M}_1$ under iterations of $\hat{F}^{-1}_1$. 
    
    Let $\hat{x}_L\dfn\hat{F}^{L}_1(\Fc^l(x))$ for $L<0$. Note that there are no consecutive $\hat{x}_{L-1},\hat{x}_{L}\in\Fc(\Nin)$. Otherwise,  $\Fc^{-1}(\hat{x}_{L-1})\in\Nin$ would be the point in $\Nin$ that returns $\Mrin$ landing at $\Fc^{-1}(\hat{x}_L)\in\Nin$. This is impossible under our chosen $R_{\text{\upshape HF}}(r,\phistar)$ in \cref{lemma:Trajectory_approximation_1,remark:GaranteeLargeTheta}\eqref{eq:xorbitA0Condition}.
    
    Therefore, there are infinitely many $\hat{x}_L$, $L<0$ such that $\hat{x}_L\in\hat{M}_1\smallsetminus\Fc(\Nin)\subset\hat{\mathfrak{U}}_+$. 
    
    We choose $L_2<L_1<0$ with $k_2<k_1<0$ such that $\hat{x}_{L_1}=\Fc^{k_1}(\Fc^l(x))\in\hat{\mathfrak{U}}_{+}$, $\hat{x}_{L_2}=\Fc^{k_2}(\Fc^l(x))\in\hat{\mathfrak{U}}_{+}$. Note that since there are infinitely many such $L_2<0$, $L_1<0$, so we can choose $L_1<<0$ arbitrarily far away from $0$. Here we first assume $L_1<-10$.
    
    \begin{equation}\label{eq:sigmaExpansionCone}
    \begin{aligned}\sigma_{\mathcal{C}}(D\Fc_{\hat{x}_{L_2}}^{-k_2})&\overbracket{=}^{\mathclap{\eqref{eq:ExpansionInCone}}}\inf_{u\in\mathcal{C}(\hat{x}_{L_2})}\sqrt{\frac{Q_{\Fc^{-k_2}({\hat{x}_{L_2}})}(D\Fc_{\hat{x}_{L_2}}^{-k_2}(u))}{Q_{\hat{x}_{L_2}}(u)}}\\
    &=\inf_{u\in\mathcal{C}(\hat{x}_{L_2})}\sqrt{\frac{Q_{\Fc^{-k_2}({\hat{x}_{L_2}})}(D\Fc_{\hat{x}_{L_2}}^{-k_2}(u))}{Q_{\Fc^{k_1-k_2}({\hat{x}_{L_2}})}(D\Fc_{\hat{x}_{L_2}}^{k_1-k_2}(u))}\cdot\frac{Q_{\Fc^{k_1-k_2}({\hat{x}_{L_2}})}(D\Fc_{\hat{x}_{L_2}}^{k_1-k_2}(u))}{Q_{\hat{x}_{L_2}}(u)}}\\
    &\ge\inf_{u\in\mathcal{C}(\hat{x}_{L_2})}\sqrt{\frac{Q_{\Fc^{-k_2}({\hat{x}_{L_2}})}(D\Fc_{\hat{x}_{L_2}}^{-k_2}(u))}{Q_{\Fc^{k_1-k_2}({\hat{x}_{L_2}})}(D\Fc_{\hat{x}_{L_2}}^{k_1-k_2}(u))}}\cdot\underbracket{\inf_{u\in\mathcal{C}(\hat{x}_{L_2})}\sqrt{\frac{Q_{\Fc^{k_1-k_2}({\hat{x}_{L_2}})}(D\Fc_{\hat{x}_{L_2}}^{k_1-k_2}(u))}{Q_{\hat{x}_{L_2}}(u)}}}_{\mathclap{\substack{\text{\cref{remark:MonotoneCone}: }>1\text{ since }D\Fc_{\hat{x}_{L_2}}^{k_1-k_2}(u)\in\Int\mathcal{C}(\Fc^{k_1-k_2}({\hat{x}_{L_2}}))=\Int\mathcal{C}({\hat{x}_{L_1}})\\\text{ the cone is strictly invariant on } \hat{\mathfrak{U}}_{+}\cup\hat{\mathfrak{U}}_{-}}}}\\
    &\overbracket{>}^{\mathllap{\Fc^{k_1-k_2}(\hat{x}_{L_2})=\hat{x}_{L_1}}}\inf_{u\in\mathcal{C}(\hat{x}_{L_2})}\sqrt{\frac{Q_{\Fc^{-k_2}({\hat{x}_{L_2}})}(D\Fc_{\hat{x}_{L_2}}^{-k_2}(u))}{Q_{\hat{x}_{L_1}}(D\Fc_{\hat{x}_{L_2}}^{k_1-k_2}(u))}}=\inf_{u\in\mathcal{C}(\hat{x}_{L_2})}\sqrt{\frac{Q_{\Fc^{l}(x)}(D\Fc_{\hat{x}_{L_2}}^{-k_2}(u))}{Q_{\hat{x}_{L_1}}(D\Fc_{\hat{x}_{L_2}}^{k_1-k_2}(u))}}
    \end{aligned}   
    \end{equation}
    
    Denote by $\Fc^l(x)=(\phi,\theta)$ and $(d\phi,d\theta)=dx\dfn D\Fc_{\hat{x}_{L_2}}^{-k_2}(u)=D\Fc_{\hat{x}_{L_1}}^{-k_1}(d\hat{x}_1)\in T_{\Fc^l(x)}M$.
    
    Suppose $(\hat{\phi}_2,\hat{\theta}_2)=\hat{x}_{L_2}$ and $(\hat{\phi}_1,\hat{\theta}_1)=\hat{x}_{L_1}$. Let $u=(d\hat{\phi}_2,d\hat{\theta}_2)\dfn d\hat{x}_2\in \mathcal{C}(\hat{x}_{L_2})=\big\{dx=(d\phi,d\theta)\in T_{\hat{x}_{L_2}}M\bigm|\frac{d\theta}{d\phi}\in[0,\infty]\big\}$, then $(d\hat{\phi}_1,d\hat{\theta}_1)=d\hat{x}_1\dfn D\Fc^{k_1-k_2}_{\hat{x}_{L_2}}(u)=D\Fc^{k_1-k_2}_{\hat{x}_{L_2}}(d\hat{x}_2)\in\big\{dx=(d\phi,d\theta)\in T_{\hat{x}_{L_1}}M\bigm|\frac{d\theta}{d\phi}\in[0,\infty]\big\}$ by \cref{caseA,caseB,caseC}. Especially if $\hat{x}_1\in\Mrin\cap\Mrout$, by \cref{eq:narrowconereason1,eq:narrowconereason2,eq:narrowconereason3}, $d\hat{x}_1$ is in the half-quadrant $\big\{dx=(d\phi,d\theta)\in T_{\hat{x}_{L_1}}M\bigm|\frac{d\theta}{d\phi}\in[0,1]\big\}$. 
    
    In orbit segment: $\underbracket{\hat{\mathfrak{U}}_{+}\ni\hat{x}_{L_1},\Fc(\hat{x}_{L_1}),\cdots,\Fc^{-k_2}(\hat{x}_{L_1})=\Fc^l(x)\in\hat{\mathfrak{U}}_{+}}_{\mathclap{\text{contains }\hat{x}_{L_1},\cdots,\hat{x}_{-1},\Fc^{l}(x)\in\hat{\mathfrak{U}_1}}}$, there exist $-L_1+1$ elements in $\hat{\mathfrak{U}}_{+}$, therefore in the orbit segment there exist at least $-L_1+1$ elements in $\hat{M}_1$. Hence, \cref{thm:UniformExpansionOnMhat1} gives $\frac{\|dx\|_{\p}}{\|d\hat{x}_1\|_{\p}}>c_2\Lambda_2^{-L_1}$ with constants $c_2>0,\Lambda_2>1$. \cref{corollary:sinpmetricexpansionMhat1} gives $\dfrac{\sin{\hat{\theta}_1}\|dx\|_\p}{\|d\hat{x}_1\|_\p}>\hat{c}_1\Lambda_1^{-L_1}$.

    \begin{equation}\label{eq:Qnormratio}
    \begin{aligned}
        \frac{Q_{\Fc^{l}(x)}(D\Fc_{\hat{x}_{L_2}}^{-k_2}(u))}{Q_{\hat{x}_{L_1}}(D\Fc_{\hat{x}_{L_2}}^{k_1-k_2}(u))}\qquad\overbracket{=}^{\mathclap{\text{\cref{def:QuadraticFormForCone}}}}\qquad&\frac{\|dx\|^2_Q}{\|d\hat{x}_1\|^2_Q}\qquad\overbracket{=}^{\mathclap{\text{\cref{def:billiardQuadrantConeAndForm}}}}\qquad\frac{|\sin\theta d\phi d\theta|}{|\sin{\hat{\theta}_1}d\hat{\phi}_1d\hat{\theta}_1|}=\frac{\sin{\hat{\theta}_1}}{\sin\theta}\frac{\sin^2\theta d\phi^2 |d\theta/d\phi|}{\sin^2{\hat{\theta}_1}d\hat{\phi}_1^2 |d\hat{\theta}_1/d \hat{\phi}_1|}\\
        \overbracket{=}^{\mathllap{\p\text{-metric definition in \cref{def:SectionSetsDefs}}}}\qquad&\frac{\sin{\hat{\theta}_1}}{\sin\theta}\frac{\|dx\|^2_{\p}}{\|d\hat{x}_1\|^2_{\p}}\frac{|d\theta/d\phi|}{|d\hat{\theta}_1/d \hat{\phi}_1|}\\
        \overbracket{>}^{\mathllap{\text{\cref{thm:UniformExpansionOnMhat1,corollary:sinpmetricexpansionMhat1}}}}\qquad&\hat{c}_1c_1(\Lambda_1)^{-2L_1}\frac{|d\theta/d\phi|}{|d\hat{\theta}_1/d \hat{\phi}_1|}\text{ for all }u\in\mathcal{C}(\hat{x}_{L_2})\\
        >\qquad&\hat{c}_1c_1(\Lambda_1)^{-2L_1}\Upsilon>9\qquad\text{, if }L_1\text{ satisfies \eqref{eq:L2largeConditionTobesufficient}.}
    \end{aligned}
    \end{equation}
Explanation for $\Upsilon$: Since $\mathfrak{U}$ is a subset of $\mathfrak{D}$ in \cref{fig:A_compactRegion}, by the conclusions of \cref{corollary:cone-lowerbound},  \eqref{eq:narrowconereason2} and \eqref{eq:narrowconereason3}, we can conclude the following. 

If $x\in\Mrin\cap\Mrout$, then $d\theta/d\phi\in\big[\lambda_0(r,R,\phistar),\lambda_1(r,R,\phistar)\big]$.

If $x\in\Fc(\Mrin\smallsetminus\Mrout)$, then $d\theta/d\phi\in\big[\frac{1}{2+1/\lambda_0(r,R,\phistar)},\frac{1}{2+1/\lambda_1(r,R,\phistar)}\big]$.

Similarly, if $\hat{x}_{L_1}\in\Mrin\cap\Mrout$, then $d\hat{\theta}_1/d\hat{\phi}_1\in\big[\lambda_0(r,R,\phistar),\lambda_1(r,R,\phistar)\big]$.

If $\hat{x}_{L_1}\in\Fc(\Mrin\smallsetminus\Mrout)$, then $d\hat{\theta}_1/d\hat{\phi}_1\in\big[\frac{1}{2+1/\lambda_0(r,R,\phistar)},\frac{1}{2+1/\lambda_1(r,R,\phistar)}\big]$.\newline The $\lambda_0(r,R,\phistar)$, $\lambda_1(r,R,\phistar)$ are constants from \cref{corollary:cone-lowerbound} determined by the billiard configuration. Hence, $d\theta/d\phi$ is bounded above $0$. 

Thus, $d\hat{\theta}_1/d\hat{\phi}_1$ is bounded away from $\infty$. Therefore, we can conclude that $\frac{|d\theta/d\phi|}{|d\hat{\theta}_1/d \hat{\phi}_1|}>\Upsilon$, that is, it is uniformly bounded above some constant $\Upsilon(r,R,\phistar)>0$, depending on the billiard configuration. 

Therefore, if \begin{equation}\label{eq:L2largeConditionTobesufficient}
-L_2>-L_1>\frac{\log(9/(\hat{c}_1c_1\Upsilon))}{2\log(\Lambda_1)},
\end{equation} then by \eqref{eq:sigmaExpansionCone}\eqref{eq:Qnormratio} $\sigma_{\mathcal{C}}(D\Fc_{\hat{x}_{L_2}}^{-k_2})>3$. Note that there must exist $L_2<L_1<<0$ to ensure \eqref{eq:L2largeConditionTobesufficient} since there are infinitely many $\hat{x}_L$, $L<0$ satisfying \eqref{eq:Qnormratio}.

Hence, with $N=-k_2$, $\Fc^{-N}(\Fc^l(x))=\hat{x}_{L_2}$, since $\Fc^{-N}$ is also a local diffeomorphism by \cite[equation (2.26) and theorem 2.33]{cb}, we can choose a sufficiently small $U\subset\hat{\mathfrak{U}}_{+}$ an open neighborhood of $\Fc^l(x)$ that satisfies $U\cap\mathcal{S}_{-N}=\emptyset$ and $\Fc^{-N}$ is a diffeomorphism on $U$.

For $y\in U$, let $\hat{y}_{L}=\hat{F}^{L}_1(y)$, $L<0$, then for $\hat{y}_{L_2}=\Fc^{k_2}(\Fc^l(x))=\Fc^{-N}(\Fc^l(x))$. For the same reasons as in \eqref{eq:Qnormratio}, $\sigma_{\mathcal{C}}(D\Fc_{\hat{y}_{L_2}}^{-k_2})>3$, as long as $L_2<L_1<<0$ satisfies \eqref{eq:Qnormratio}. Since $\Fc^{-N}$ is a diffeomorphism on $U$, $\hat{y}_{L_2}$ can be any point of $\Fc^{-N}U$. Therefore, $\sigma_{\mathcal{C}}(D\Fc^N_y)>3$ for every $y\in\Fc^{-N}U$.
\end{proof}
\begin{theorem}[Sufficient points on $\hat{\mathfrak{U}}_{-}$]\label{theorem:sufficentptsonMhatminus} For every $x\in\hat{\mathfrak{U}}_{-}\smallsetminus\mathcal{S}_{-\infty}$, the followings hold.
    \begin{enumerate}[label={(\roman*)}]
        \item\label{item01sufficentptsonMhatminus}  
        $\exists$ some integer $l<0$ such that $\Fc^{l}$ is a local diffeomorphism at $x$ such that $\Fc^l(x)\in\hat{\mathfrak{U}}_{-}$. 
        \item\label{item02sufficentptsonMhatminus} And there is an open neighborhood $U\subset \hat{\mathfrak{U}}_{-}$ of $\Fc^l(x)$ and an integer $N>0$ such that $U\cap\mathcal{S}_{-N}=\emptyset$ and $\Fc^{-N}(U)\subset \hat{\mathfrak{U}}_{-}$. The cone field $\mathcal{C}$ given in \cref{def:billiardQuadrantConeAndForm} restricted on $U\cup\Fc^{-N}(U)$ is invariant and $\sigma_{\mathcal{C}}(D\Fc^N_y)>3$ for every $y\in\Fc^{-N}U$.
    \end{enumerate}
\end{theorem}
\begin{proof}
    The proof of \ref{item01sufficentptsonMhatplus} is the same as \cref{theorem:sufficentptsonMhatplus}\ref{item01sufficentptsonMhatplus}.

    Proof of (ii): It is clear by \cref{caseA,caseB,caseC} the positive cone $\mathcal{C}$ given in \cref{def:billiardQuadrantConeAndForm} restricted on $U\cup\Fc^{-N}(U)$ is invariant for any $N>0$ and open $U\subset\hat{\mathfrak{U}}_{-}$. For any $N>0$,  as long as $U\ni\Fc^l(x)$ is small enough, $\Fc^{-N}$ is a diffeomorphism by \cite[equation (2.26) and theorem 2.33]{cb} on $U$ such that $U\cap\mathcal{S}_{-N}=\emptyset$.
    
    For $x\in\hat{\mathfrak{U}}_{-}\smallsetminus\mathcal{S}_{-\infty}$, $l<0$ given in (i), $\Fc^{l}(x)\in\hat{\mathfrak{U}}_{-}\smallsetminus\mathcal{S}_{-\infty}$. 
    
    For each $x\in\hat{M}\smallsetminus\mathcal{S}_{-\infty}$, suppose that $j(x)\dfn\inf\big\{i>0,\Fc^{-i}(x)\in\hat{M}_1\big\}$.
    
    For each $x\in\hat{M}_1\smallsetminus\mathcal{S}_{-\infty}$, suppose that $j_1(x)\dfn\inf\big\{i>0,\Fc^{-i}(x)\in\hat{M}\big\}$.
    
    Then we denote $x_{-}\dfn\Fc^{-j(\Fc^{l}(x))}(\Fc^{l}(x))\in\hat{M}_1$. And we further define $\hat{M}_1\ni\hat{x}_{L}\dfn\hat{F}^{L}_{1}(x_{-})$ for $L<0$ as the backward trajectory of $x_{-}$ under $\hat{F}^{-1}_{1}$ iterations. Then we pick $L_2<L_1<<0$ satisfying \eqref{eq:L2largeConditionTobesufficient} with $k_2<k_1<0$ such that $\hat{x}_{L_2}=\Fc^{k_2}(x_{-})$, $\hat{x}_{L_1}=\Fc^{k_1}(x_{-})$.

    For all $u\in\mathcal{C}(\hat{x}_{L_2})$, since $\hat{x}_{L_2},x_{-}\in\hat{M}_1$,
    by the same reason as in \eqref{eq:sigmaExpansionCone} \eqref{eq:Qnormratio} $\frac{Q_{x_{-}}(D\Fc_{\hat{x}_{L_2}}^{-k_2}(u))}{Q_{\hat{x}_{L_2}}(u)}>9$.

    Then we let $\hat{x}=\Fc^{-j_1(\hat{x}_{L_2})}(\hat{x}_{L_2})\in\hat{M}$, for all $u\in\mathcal{C}(\hat{x})$ by \cref{remark:MonotoneCone}, since $\mathcal{C}$ is strictly invariant on $\hat{\mathfrak{U}}_+\cup\hat{\mathfrak{U}_-}$, $D\Fc^{j_1(\hat{x}_{L_2})}u\in\Int(\mathcal{C}(\hat{x}_{L_2}))$, $\frac{Q_{\hat{x}_{L_2}}(D\Fc_{\hat{x}}^{j_1(\hat{x}_{L_2})}u)}{Q_{\hat{x}}(u)}>1$.
    
    The same reason as in \cref{remark:MonotoneCone} gives for all $u\in\mathcal{C}(x_{-})$ since $\mathcal{C}$ is strictly invariant in $\hat{\mathfrak{U}}_+\cup\hat{\mathfrak{U}_-}$, \[\frac{Q_{\Fc^l(x)}(D\Fc_{x_{-}}^{j(\Fc^l(x))}(u))}{Q_{x_{-}}(u)}>1.\] 
    
    Therefore, we choose $N=j_1(\hat{x}_{L_2})-k_2+j(\Fc^{l}(x))$ so that $\Fc^N(\hat{x})=\Fc^l(x)$. By chain rule and the invariance of the cone $\mathcal{C}$ it gives:
    \begin{equation}\label{eq:sigmaconeexpansionNmhat1}
        \begin{aligned}
            \sigma_{\mathcal{C}}(D\Fc^N_{\hat{x}})&=\inf_{u\in\mathcal{C}(\hat{x})}\sqrt{\frac{Q_{\Fc^{l}(x)}(D\Fc_{\hat{x}}^{N}(u))}{Q_{\hat{x}}(u)}}\\
            &=\inf_{u\in\mathcal{C}(\hat{x})}\sqrt{\underbracket{\frac{Q_{\Fc^{l}(x)}(D\Fc_{\hat{x}}^{N}(u))}{Q_{x_{-}}(D\Fc_{\hat{x}}^{j_1(\hat{x}_{L_2})-k_2}(u))}}_{>1}\cdot\underbracket{\frac{Q_{x_{-}}(D\Fc_{\hat{x}}^{j_1(\hat{x}_{L_2})-k_2}(u))}{Q_{\hat{x}_{L_2}}(D\Fc_{\hat{x}}^{j_1(\hat{x}_{L_2})}(u))}}_{>9}\cdot\underbracket{\frac{Q_{\hat{x}_{L_2}}(D\Fc_{\hat{x}}^{j_1(\hat{x}_{L_2})}(u))}{Q_{\hat{x}}(u)}}_{>1}}>3.
        \end{aligned}
    \end{equation}
Since $\Fc^{-N}$ is also a local diffeomorphism by \cite[equation (2.26) and theorem 2.33]{cb}, we can choose an open neighborhood small enough $U\subset\hat{\mathfrak{U}}_{-}$ of $\Fc^l(x)$ that satisfies $U\cap\mathcal{S}_{-N}=\emptyset$ and $\Fc^{-N}$ is a diffeomorphism on $U$. For each $y\in U$, let $\hat{y}=\Fc^{-N}(y)$. Then for the same reason as in \eqref{eq:sigmaconeexpansionNmhat1}, \[\sigma_{\mathcal{C}}(D\Fc^N_{\hat{y}})>3.\] Since $\Fc^{-N}$ is a diffeomorphism on $U$, $\hat{y}$ can be any point of $\Fc^{-N}U$. Therefore, $\sigma_{\mathcal{C}}(D\Fc^N_y)>3$ for every $y\in\Fc^{-N}U$.
\end{proof}
\begin{corollary}\label{corollary:sufficientpointInMhatPlusMinus}
    Every point $x\in(\hat{\mathfrak{U}}_{-}\cup\hat{\mathfrak{U}}_{+})\smallsetminus\mathcal{S}_{-\infty}$ is a sufficient point with some quadruple $(l,N,U,\mathcal{C})$ (\cref{def:sufficientPoint}).
\end{corollary}
\begin{proof}
    Every $x\in\hat{\mathfrak{U}}_{-}$ is a sufficient point with quadruple $(l,N,U,\mathcal{C})$ given in \cref{theorem:sufficentptsonMhatminus}. Every $x\in\hat{\mathfrak{U}}_{+}$ is a sufficient point with quadruple $(l,N,U,\mathcal{C})$ given in \cref{theorem:sufficentptsonMhatplus}.
\end{proof}

\subsection[Local (un)stable manifolds and L4 contraction condition]{Local (un)stable manifold and L4 contraction condition}

\begin{proposition}\label{proposition:stablemanifoldthm}\cite[Proposition 3.4 and Definition 3.5]{MR3082535} Let $x\in M\smallsetminus\partial M$ be a sufficient point (\cref{def:sufficientPoint}) with quadruple $(l,N,U,\mathcal{C})$. Then there exists an invariant measurable set $\Omega\subset\bigcup_{k\in\mathbb{Z}}\Fc^kU$ with $\mu(\bigcup_{k\in\mathbb{Z}}U\smallsetminus\Omega)=0$ and two families of $C^2$ submanifolds $V^s=\big\{V^s_y\big\}_{y\in\Omega}$ and $V^u=\big\{V^u_y\big\}_{y\in\Omega}$ such that for every $y\in\Omega$, the following hold:
    \begin{enumerate}
    \item\label{item01stablemanifoldthm} $V^s_y\cap V^u_y=\{y\}$;
    \item\label{item02stablemanifoldthm} $V^s_y$ and $V^u_y$ are embedded $1$-dimensional intervals;
    \item\label{item03stablemanifoldthm} $T_y V^s_y\subset \mathcal{C}'(y)$ and $T_y V^u_y\subset \mathcal{C}(y)$ provided that $y\in U\cup\Fc^{-N}U$;
    \item\label{item04stablemanifoldthm} $\Fc V^s_y\subset V^s_{\Fc y}$ and $\Fc^{-1}V^u_y\subset V^u_{\Fc^{-1}y}$;
    \item\label{item05stablemanifoldthm} $d(\Fc^n y,\Fc^n z)\rightarrow 0$ exponentially as $n\rightarrow+\infty$ for every $z\in V^s_y$, and the same is true as $n\rightarrow-\infty$ for every $z\in V^u_y$.
    \end{enumerate} Furthermore,  $V^s_y$ and $V^u_y$ vary measurably with $y\in\Omega$, and the families $V^s$ and $V^u$ have the absolute continuity property. 
    
    \cite[Definition 3.5]{MR3082535}: The submanifolds forming the families $V^s$ and $V^u$ are called \emph{local stable manifolds} and \emph{local unstable manifolds}, respectively. For every $y\in\Omega$, denote by $W^u_y$ the connected component of $\bigcup_{k\ge0}\Fc^kV^u_{\Fc^{-k}y}$ containing $y$. Analogously, denote by $W^s_y$ the set obtained by replacing $\Fc$ with $\Fc^{-1}$ and $V^u$ with $V^s$ in definition of $W^u_y$. The sets $W_y^s$ and $W_y^u$ are immersed submanifolds of $M$.
\end{proposition}
\begin{proof}
    The proof is in \cite[Proposition 3.4]{MR3082535} and is an application of \cite[Katok-Strelcyn theory]{MR872698}.
\end{proof}
\begin{lemma}[Unstable manifold and the cone]\label{lemma:UnstableManifoldInCone} We have the following conclusions for unstable manifolds $W^u_y$.
\begin{enumerate}
    \item\label{item01UnstableManifoldInCone} For $z\in W^u_y\cap\hat{M}_1$, $(d\phi,d\theta)\nfd dz=T_z(W^u_y)$ satisfies $dz\in \hat{C}^u_{1,z}$ (cone defined \cref{thm:UniformExpansionOnMhat1}).
    \item\label{item02UnstableManifoldInCone} For $z\in W^u_y\cap\hat{M}_1$, $(d\phi,d\theta)\nfd dz=T_z(W^u_y)$ and let $z_k\dfn\hat{F}^{-k}(z)$, then $dz=\bigcap_{k\ge0}(D\hat{F}^{k}_{1,z_k}(\hat{C}^u_{1,{z_k}}))$, where $\hat{C}^u_{1,{z_k}}$ is the cone defined in \cref{thm:UniformExpansionOnMhat1}.

    \item\label{item03UnstableManifoldInCone} For $z\in W^u_y\cap\hat{M}$, $(d\phi,d\theta)\nfd dz=T_z(W^u_y)$ satisfies $dz\in C_z$ (cone defined in \cref{RMASUHwithProof}).
    \item\label{item04UnstableManifoldInCone} For $z\in W^u_y\cap\hat{M}$, $(d\phi,d\theta)\nfd dz=T_z(W^u_y)$ and let $z_k\dfn\hat{F}^{-k}(z)$, then $dz=\bigcap_{k\ge0}(D\hat{F}^{k}_{z_k}(C_{z_k}))$, where $C_z$ is the cone defined in \cref{RMASUHwithProof}.
\end{enumerate}
\end{lemma}
\begin{proof}
    \eqref{item02UnstableManifoldInCone} implies \eqref{item01UnstableManifoldInCone}. \eqref{item04UnstableManifoldInCone} implies \eqref{item03UnstableManifoldInCone}. And both \eqref{item02UnstableManifoldInCone} and \eqref{item04UnstableManifoldInCone} are the conclusion of the cone technique in \cite[Equation (3.55)]{cb}.
\end{proof}

\begin{theorem}[Lemon billiard satisfies L4 condition in \cref{thm:lemonLET}]\label{theorem:L4condition}\hfill

For the lemon billiard $\mathbb{L}(r,R,\phistar)$ satisfying the condition \cref{def:UniformHyperbolicLemonBilliard} with $\Phistar=\sin^{-1}(r/R\sin\phistar)$ (see \cref{fig:1-petalExtensionFromCornerAndLemon}), suppose that $x\in\hat{\mathfrak{U}}_{+}$ is a sufficient point with quadruple $(l,N,U,\mathcal{C})$ given in \cref{theorem:sufficentptsonMhatplus} or $x\in\hat{\mathfrak{U}}_{-}$ is a sufficient point with quadruple $(l,N,U,\mathcal{C})$ given in \cref{theorem:sufficentptsonMhatminus}. 

And let $\Omega$ be the subset of $\bigcup_{k\in\mathbb{Z}}\Fc^k U$ given in \cref{proposition:stablemanifoldthm}. 

Then there exist constants $\beta>0$ and $\xi>0$ such that: \begin{enumerate}
    \item\label{item01L4condition} if $y\in\Omega\cap U$, $z\in W^u_y$ and $\Fc^{-k}z\in\mathcal{S}^{-}_{1}(\xi)$ with $k>0$, then $\|D\Fc_z^{-k}\bigm|_{T_zW^u_y}\|\le\beta$,
    \item\label{item02L4condition} if $y\in\Omega\cap U$, $z\in W^s_y$ and $\Fc^kz\in\mathcal{S}^{+}_1(\xi)$ with $k>0$, then $\|D\Fc^{k}_z\bigm|_{T_zW^s_y}\|\le\beta$,
\end{enumerate} where $\mathcal{S}^{+}_1(\xi)$, $\mathcal{S}^{-}_{1}(\xi)$ are the $\xi-$ neighborhoods of $\mathcal{S}^{+}_1$ and $\mathcal{S}_{1}^{-}$ per \cref{notation:symplecticformmeasurenbhd}.
\end{theorem}
\begin{proof}
    Since by symmetry $I\,\mathcal{S}^{-}_{1}(\xi)=\mathcal{S}^{+}_1(\xi)$, $I\,W^u_y=W^s_{I(y)}$, $I\,\hat{\mathfrak{U}}_{+}=\hat{\mathfrak{U}}_{-}$, $I\Fc=\Fc^{-1}I$, it suffices to only prove the first conclusion \eqref{item01L4condition} i.e. for $x$ and $U$ to be either
    \begin{itemize}
        \item a sufficient point $x\in\hat{\mathfrak{U}}_{+}$ with quadruple $(l,N,U,\mathcal{C})$ given in \cref{theorem:sufficentptsonMhatplus} that is with $l<0$, $U$ as a neighborhood of $\Fc^{l}(x)\subset U\subset\hat{\mathfrak{U}}_{+}\subset \hat{M}_{+}=\hat{M}_{1}$.
        \item a sufficient point $x\in\hat{\mathfrak{U}}_{-}$ with quadruple $(l,N,U,\mathcal{C})$ given in \cref{theorem:sufficentptsonMhatminus} that is with $l<0$, $U$ as a neighborhood of $\Fc^{l}(x)\subset U\subset\hat{\mathfrak{U}}_{-}\subset\hat{M}_{-}\subset\hat{M}$.
    \end{itemize}
    
    On the other hand, if $z\in \mathcal{S}_{-\infty}$, then $z\in\partial W^u_y$ and $D\Fc_z^{-k}\bigm|_{T_zW^u_y}$ is the one-sided limit of the interior points of $W^u_y$: $\lim_{\Int{(W^u_y)\ni x\rightarrow\partial(W^u_y)}} D\Fc_{x}^{-k}\bigm|_{\partial(W^u_y)}.$ Because by \cref{proposition:stablemanifoldthm} $W^u_y$ is $C^2$ and by \cite[Equation (2.26) and Theorem 2.33]{cb} $\Fc^{-1}$ and its iteration is smooth, it suffices to prove conclusion \eqref{item01L4condition} for $z\notin\mathcal{S}_{-\infty}$. Now we assume $z\notin\mathcal{S}_{-\infty}$.

     For $y\in W^u_{y}\subset\Omega\cap U$ and $y'$ with $k>0$ satisfying $\Fc^{k}(y')=y\in\Omega\cap U$, the chain rule and the local unstable mainfold's invariance under $\Fc^{-1}$ (\cref{proposition:stablemanifoldthm}\eqref{item04stablemanifoldthm}) imply that for every $z\in W^{u}_y$, $z\notin\mathcal{S}_{-\infty}$, there is a unique $z'=\Fc^{-k}(z)\in W^u_{\Fc^{-k}(y)}=W^u_{y'}$. Therefore, the conclusion \eqref{item01L4condition} is equivalent to $\exists$  uniform $\beta>0$, $\xi>0$ such that 
     \begin{equation}\label{eq:L4AlternativeClaim}
        \begin{aligned}
        \text{for all }z'=\Fc^{-k}(z)\in W^{u}_{\Fc^{-k}(y)}\cap\mathcal{S}^{-}_{1}(\xi)&=W^{u}_{y'}\cap\mathcal{S}^{-}_{1}(\xi),\\
        \|D\Fc^{k}_{z'}\bigm|_{T_{z'}W_{y'}^u}\|&\ge 1/\beta. 
        \end{aligned}
     \end{equation}
    
    Since the connected components of $\mathcal{S}^{-}_{1}$ are the line segments: $\partial\Mrin\smallsetminus\partial M_r$ and $\partial\MRin\smallsetminus\partial M_R$ (see \cref{fig:Mr,fig:MR}), we make proofs based on the two cases for singularity curve component $\mathbb{S}$ to be $\partial\Mrin$ or $\partial\MRin$.
    
    \textbf{Singularity curve case (1)}: Let $\mathbb{S}=\partial\Mrin\smallsetminus\partial M_r$ i.e. be the dashed line segment as the boundary of $\Mrin$ in \cref{fig:hatUplusUminus}. Suppose $z'\in W^u_{y'}\cap\mathbb{S}(\xi)$ in the $\xi-$ neighborhood of $\mathbb{S}$. And by symmetry we can assume that $z'$ is in the $\xi$-neighborhood of the two closed line segments $LL-K3$ or $IY-K1$ in \cref{fig:hatUplusUminus}. 
    
    We make sure that $\xi$ is small enough to satisfy \begin{equation}\label{eq:xicondition1}
        \begin{aligned}
            &\mathbb{S}(\xi)\cap\Nout=\emptyset,\,\Fc\big((\mathbb{S}\setminus\Mrout)(\xi)\big)\cap\Nout=\emptyset  
            \text{ and }\Fc^2\big((\mathbb{S}\setminus(\Mrout\cup\Fc^{-1}(\Mrout\smallsetminus\Mrin)))(\xi)\big)\cap\Nout=\emptyset.
        \end{aligned}
    \end{equation}

    This $\xi$ restriction is a not hard to satisfy by the same proof for \cref{Prop:Udisjoint}. Also note that
    $\Fc(\mathbb{S}(\xi)\setminus\Mrout)$ and 
    are the boundary of blue strips in \cref{fig:hatUplusUminus}. 
    
     Local unstable manifold $W^u_{y'}$ is a connected line segment and cannot intersect with $(\mathcal{S}_{-1}\smallsetminus\partial M)\subset\partial\Mrin$ at interior point of $W^u_{y'}$ by the definition of $W^u_{y'}$ in \cref{proposition:stablemanifoldthm}. Hence either $z'\in W^u_{y'}\subset\Mrin$ or $z'\in W^u_{y'}\subset(M_r\smallsetminus\Mrin)$. We have the following two subcases of the \textbf{ singularity curve case (1)}.

    \textbf{Singularity curve subcase (1.1):} $z'\in W^u_{y'}\subset(M_r\smallsetminus\Mrin)\cap \mathbb{S}(\xi)$, $dz'\in T_{z'}W^u_{y'}$. 

    First we note that by the local unstable manifolds definition in \cref{proposition:stablemanifoldthm}, $y\notin\mathcal{S}_{-\infty}$ and $z\notin\mathcal{S}_{-\infty}$. And by $F^k(z')=z\in W^u_y\subset U\subset\hat{\mathfrak{U}}_{+}\cup\hat{\mathfrak{U}}_{-}$, $z'$ cannot be in a periodic orbit with points that only collide with $\Gamma_r$. Hence, there exists a finite $j_1(z')\dfn\inf\big\{j'\ge0\bigm|p(\Fc^{-j'}(z'))\in\Gamma_R\big\}$. Since $z'\in M_r\smallsetminus\Mrin$, $j_1(z')\ge2$ and $\Fc^{-j_1(z')+1}(z')\in\Mrin\smallsetminus\Mrout$. 
    
    Also note that it is also true that $\Fc^{-j_1(z')+1}(y')\in\Mrin\smallsetminus\Mrout$. Otherwise, if $\Fc^{-j_1(z')+1}(y')\notin \Mrin$, then  $\Fc^{-j_1(z')}(y')$ and $\Fc^{-j_1(z')}(z')$ are not both in $M_R$. Or if $\Fc^{-j_1(z')+1}(y')\in \Mrout$, then $\Fc^{-j_1(z')+2}(y')$ and $\Fc^{-j_1(z')+2}(z')$ are not both in $M_r$. This means that the $y$'s and $z$'s preimages under $\Fc$ are separate/do not belong to the same local unstable manifold. This contradicts \cref{proposition:stablemanifoldthm}\eqref{item04stablemanifoldthm}.
    
    Again By \cref{proposition:stablemanifoldthm}\eqref{item04stablemanifoldthm}: the $\Fc^{-1}$ invariance of unstable manifold, $\Fc^{-j_1(z')+1}(z')$ is on the local unstable manifold of $\Fc^{-j_1(z')+1}(y')$, we can define $\hat{z}_{+}\dfn\Fc^{-j_1(z')+2}(z')$ and $\hat{y}_{+}\dfn\Fc^{-j_1(z')+2}(y')$ so that $\hat{z}_{+}\in W^u_{\hat{y}_{+}}\subset\Fc(\Mrin\smallsetminus\Mrout)\subset \hat{M}_1=\hat{M}_{+}$ with $d\hat{z}_{+}\in T_{\hat{z}_{+}} W^u_{\hat{y}_{+}}$ and $D\Fc_{\hat{z}_{+}}^{j_1(z')-2}(d\hat{z}_{+})=dz'$. By \cref{lemma:UnstableManifoldInCone}\eqref{item02UnstableManifoldInCone} $d\hat{z}_{+}\in T_{\hat{z}_{+}} W^u_{\hat{y}_{+}}\subset\hat{C}^{u}_{1,\hat{z}_{+}}\subset\mathcal{Q}_{\hat{z}_{+}}(\mathrm{I,III})=\big\{(d\phi,d\theta)\in T_{\hat{z}_{+}}M\bigm|\frac{d\theta}{d\phi}\in[0,+\infty]\big\}$.
    
    And note that $\Mrin\ni\Fc^{-j_1(z')+1}(z'),\cdots,z'\in M_r\smallsetminus\Mrin$ are all collisions on $\Gamma_r$. Therefore \[D\Fc^{j_1(z')-2}_{\hat{z}_{+}}=\begin{pmatrix}
1 & 2(j_1(z')-2) \\
0 & 1 
\end{pmatrix}\] and $dz'=D\Fc^{j_1(z')}_{\hat{z}_{+}}(d\hat{z}_{+})\in\mathcal{Q}_{z'}(\mathrm{I,III})$.

If $U\subset\hat{\mathfrak{U}}_{+}$ and $x\in\hat{\mathfrak{U}}_+$ is the sufficient point with quadruple $(l,N,U,\mathcal{C})$ given in \cref{theorem:sufficentptsonMhatplus}, then let $k_{1}(z')\dfn\inf\big\{k'\ge0,\Fcal^{k'}(z')\in\hat{M}_1=\hat{M}_{+}\big\}.$ Note that if $k_1(z')=0$ then $W^u_{y'}$ is in the $\xi-$ neighborhood of line segment of $LL-LF$ in \cref{fig:hatUplusUminus} (boundary of the blue strip) so that $W^u_{y'}\subset\Fc(\Mrin\smallsetminus\Mrout)$. Then $\frac{\|D\Fc^{k_{1}(z')}_{z'}(dz')\|}{\|dz'\|}=1$.

If the orbit segment $z'$, $\Fc(z')$, $\cdots$, $\Fc^{k_1(z')}(z')\in\hat{M}_1$ are all collisions on $\Gamma_r$, since $D\Fc^{k_1(z')}_{z'}=\begin{pmatrix}
1 & 2k_1(z') \\
0 & 1 
\end{pmatrix}$ and $dz'\in\mathcal{Q}_{z'}(\mathrm{I,III})$, 
\begin{equation}\label{eq:subcase1-1L4conclusion0}
     \frac{\|D\Fc^{k_{1}(z')}_{z'}(dz')\|}{\|dz'\|}\ge1.
\end{equation}

Otherwise the orbit segment $z'$, $\Fc(z')$, $\cdots$, $\Fc^{k_1(z')}(z')\in\hat{M}_1$ contains $x_0\in\Mrout$, $\Fc(x_0)=x_1\in\MRin$, $\cdots$, $\Fc^{n_1+1}(x_1)=x_2\in\Mrin$ where $0\le n_1\dfn\max\big\{j\ge0:p(\Fc^i(x_1))\in\Gamma_R\text{ for }i=0,\cdots,j\big\}$ corresponding to \eqref{eqMhatOrbSeg} with cases described in \eqref{eqMainCases}. Here by requiring $\xi$ satisfying \eqref{eq:xicondition1}, the same analysis for \cref{TEcase_a0,TEcase_a1,TEcase_b,TEcase_c,thm:UniformExpansionOnMhat1,corollary:sinpmetricexpansionMhat1,Col:EuclideanExpansionMhat1} yield $\frac{\|D\Fc^{k_1(z')}_{z'}(dz')\|}{\|dz'\|}>C_1$ for some constant $C_1>0$.

Since $\Fc^{k_1(z')}(z')\in\hat{M}_1$ and $\Fc^{k}(z')\in\hat{M}_1$ with $k\ge k_1(z')$, $\Fc^k(z')=\hat{F}^{t}_1(z')$ for some $t\ge0$. Let $d\hat{z}'_1=D\Fc_{z'}^{k_1(z')}(dz')$, then $D\Fc^k_{z'}(dz')=D\hat{F}^{t}_{1,z'}(d\hat{z}'_1)$.

Invoking \cref{Col:EuclideanExpansionMhat1} again gives
\begin{equation}\label{eq:subcase1-1L4conclusion1}
\|D\Fc^{k}_{z'}\bigm|_{T_{z'}W_{y'}^u}\|=\frac{\|dz\|}{\|dz'\|}=\frac{\|D\Fc^{k}_{z'}(dz')\|}{\|dz'\|}=\underbracket{\frac{\|D\Fc^{k}_{z'}(dz')\|}{\|D\Fc^{k_1(z')}_{z'}(dz')\|}}_{\mathclap{=\frac{\|D\hat{F}^t_{1,z'}(d\hat{z}'_1)\|}{\|d\hat{z}'_1\|}\text{, \cref{Col:EuclideanExpansionMhat1}:}>c_2}}\cdot\frac{\|D\Fc^{k_{1}(z')}_{z'}(dz')\|}{\|dz'\|}>c_2\cdot\min\big\{C_1,1\big\}>0,
\end{equation} where constant $c_2$ is from \cref{Col:EuclideanExpansionMhat1}.

If $U\subset\hat{\mathfrak{U}}_{-}$ and $x\in\hat{\mathfrak{U}}_{-}$ is the sufficient point with quadruple $(l,N,U,\mathcal{C})$ given in \cref{theorem:sufficentptsonMhatminus}, then let $k_{2}(z')\dfn\inf\big\{k'\ge0,\Fcal^{k'}(z')\in\hat{M}=\hat{M}_{-}\big\}$.

If the orbit segment $z'$, $\Fc(z')$, $\cdots$, $\Fc^{k_2(z')}(z')\in\hat{M}$ are all collisions on $\Gamma_r$, since $D\Fc^{k_2(z')}_{z'}=\begin{pmatrix}
1 & 2k_2(z') \\
0 & 1 
\end{pmatrix}$ and $dz'\in\mathcal{Q}_{z'}(\mathrm{I,III})$, $\frac{\|D\Fc^{k_{2}(z')}_{z'}(dz')\|}{\|dz'\|}\ge1$.

Otherwise the orbit segment $z'$, $\Fc(z')$, $\cdots$, $\Fc^{k_2(z')}(z')\in\hat{M}$ contains $x_0\in\Mrout$, $\Fc(x_0)=x_1\in\MRin$, $\cdots$, $\Fc^{n_1+1}(x_1)=x_2\in\Mrin$ where $0\le n_1\dfn\max\big\{j\ge0:p(\Fc^i(x_1))\in\Gamma_R\text{ for }i=0,\cdots,j\big\}$ corresponding to \eqref{eqMhatOrbSeg} with cases described in \eqref{eqMainCases}. Here by requiring $\xi$ satisfying \eqref{eq:xicondition1}, the same analysis for \cref{TEcase_a0,TEcase_a1,TEcase_b,TEcase_c,thm:UniformExpansionOnMhat,corollary:sinpmetricexpansionMhat,Col:EuclideanExpansionMhat} yield $\frac{\|D\Fc^{k_2(z')}_{z'}(dz')\|}{\|dz'\|}>C_2$ for some constant $C_2>0$.

Since $\Fc^{k_2(z')}(z')\in\hat{M}$ and $\Fc^{k}(z')\in\hat{M}$ with $k\ge k_2(z')$, $\Fc^k(z')=\hat{F}^{t}(\Fc^{k_2(z')}(z'))$ for some $t\ge0$. 

Invoking \cref{Col:EuclideanExpansionMhat} again gives
\begin{equation}\label{eq:subcase1-1L4conclusion2}
\|D\Fc^{k}_{z'}\bigm|_{T_{z'}W_{y'}^u}\|=\frac{\|dz\|}{\|dz'\|}=\frac{\|D\Fc^{k}_{z'}(dz')\|}{\|dz'\|}=\underbracket{\frac{\|D\Fc^{k}_{z'}(dz')\|}{\|D\Fc^{k_2(z')}_{z'}(dz')\|}}_{\mathclap{=\frac{\|D\hat{F}^t(D\Fc^{k_2(z')}_{z'}(dz'))\|}{\|D\Fc^{k_2(z')}_{z'}(dz')\|}\text{, \cref{Col:EuclideanExpansionMhat}:}>c_0}}\cdot\frac{\|D\Fc^{k_{2}(z')}_{z'}(dz')\|}{\|dz'\|}>c_0\cdot\min\big\{C_2,1\big\}>0
\end{equation}
where constant $c_0$ is from \cref{Col:EuclideanExpansionMhat}.

\textbf{Singularity curve subcase (1.2):} $z'\in W^u_{y'}\subset\Mrin\cap \mathbb{S}(\xi)$, $dz'\in T_{z'}W^u_{y'}$.

We also make sure that $\xi$ is small enough to satisfy \begin{equation}\label{eq:xicondition2}
\begin{aligned}
    &\xi<0.5\Phistar
\end{aligned}
\end{equation} where $\Phistar=\sin^{-1}(r\sin\phistar/R)$ (see \cref{fig:1-petalExtensionFromCornerAndLemon}), $\phistar$, $r$, $R$ are given by the lemon billiard configuration $\mathbb{L}(r,R,\phistar)$ of \cref{notation:symplecticformmeasurenbhd}. $z'$ is the $x_2$ of some $\hat{M}$ return orbit segment in \eqref{eqMhatOrbSeg}. 

By the tangent vectors slopes range for Arch singularity curves $\mathcal{AS}^{\text{\upshape out}}_{n}$ in \cref{proposition:Archsingularitycurvesproperties}\eqref{item03Archsingularitycurvesproperties} (see \cref{fig:ArchSingularitiesComponent2,fig:ArchSingularitiesComponent3}), the $\xi$ neighborhood of $\Scal^{-}(\xi)\cap M_r$ cannot intersect with $\mathcal{D}^{out}_1$ (\cref{def:Curvilinearquadrilateral}). Hence condition \eqref{eq:xicondition2} ensures that this $\hat{M}$ return orbit segment must be case (a) or case (b) of \eqref{eqMainCases}. 

In case (a) of \eqref{eqMainCases} and in the context of \cref{lemma:ExpansionLowerBoundA0A1},
\begin{align*}
    dx_2=dz'=\left\{
    \begin{aligned}
        &D\Fc_{\Fc^{-3}(z')}^{3}(dx')\text{ with }dx'\in  T_{\Fc^{-3}(z')}W^u_{\Fc^{-3}(y')}\subset C_{\Fc^{-3}(z')},\text{ if }\Fc^{-2}(z')=x_0\in\Mrout\smallsetminus\Mrin,\\
        &D\Fc_{\Fc^{-2}(z')}^{2}(dx')\text{ with }dx'\in  T_{\Fc^{-2}(z')}W^u_{\Fc^{-2}(y')}\subset C_{\Fc^{-2}(z')},\text{ if }\Fc^{-2}(z')=x_0\in\Mrout\cap\Mrin.
    \end{aligned}\right.
\end{align*} And by \cref{caseA}, $dz'=dx_2=(d\phi_2,d\theta_2)\in\Int{\mathcal{Q}_{x_2}(\mathrm{I,III})}=\big\{(d\phi,d\theta)\bigm|\frac{d\theta}{d\phi}\in(0,+\infty)\big\}$.

In case (b) of \eqref{eqMainCases} and in the context of \cref{lemma:ExpansionLowerBoundB}, \cref{proposition:dx3dx2expansion}: $dz'=dx_2=(d\phi_2,d\theta_2)\in\big\{(d\phi,d\theta)\bigm|\frac{d\theta}{d\phi}\in(-\infty,-1)\cup(0,+\infty)\big\}$.

If the orbit segment $z',\Fc(z'),\cdots,\Fc^{k}(z')=z\in T_zW^u_y$ are all collisions on $\Gamma_r$, then $D\Fc^k_{z'}=\begin{pmatrix}
1 & 2k \\
0 & 1 
\end{pmatrix}$. With $dx_2=dz'=(d\phi_2, d\theta_2)$, $dz=D\Fc^k_{z'}(dz')=(d\phi_2+2kd\theta_2,d\theta_2)$. Therefore, for $dz'=dx_2=(d\phi_2,d\theta_2)$ in either case (a) or (b), we get the following.
\begin{equation}\label{eq:subcase1-2L4conclusion0}
    \begin{aligned}\|dz\|^2=&(\phi_2+2kd\theta_2)^2+(d\theta_2)^2=\overbracket{(d\phi_2)^2+(d\theta_2)^2}^{=\|dz'\|^2}+4k(k(d\theta_2)^2+d\theta_2d\phi_2)\\=&\|dz'\|^2+\underbracket{4k^2(d\phi_2)^2(\frac{d\theta_2}{d\phi_2}+\frac{1}{k})\frac{d\theta_2}{d\phi_2}}_{\mathclap{>0\text{ since }k\ge1,\frac{d\theta_2}{d\phi_2}\in(-\infty,-1)\cup(0,+\infty)}}>\|dz'\|^2\\
    \|D\Fc^{k}_{z'}\bigm|_{T_{z'}W_{y'}^u}\|=&\frac{\|dz\|}{\|dz'\|}>1
\end{aligned}
\end{equation}

Now, otherwise, the orbit segment $z',\Fc(z'),\cdots,\Fc^{k}(z')=z\in T_zW^u_y$ contains at least one collision on $\Gamma_R$. Then $W^u_{y'}$ must be contained in one connected component of $\Min_{r,i}$ (see \cref{fig:MrN}) with $0\le i<k$. Otherwise, $\Fc^{i-1}(W^u_{y'})$ or $\Fc^{i+1}(W^u_{y'})$, which are the images of $W^u_y$ under some $\Fc^{-1}$ iteration, has at least two components separately contained in $M_R$ and $M_r$. So it is not connected which contradicts \cref{lemma:UnstableManifoldInCone}\eqref{item04stablemanifoldthm} and the definition of the local unstable manifold. 

Therefore, by \cref{proposition:stablemanifoldthm}\eqref{item04stablemanifoldthm}, the following orbit approximation property for $z'$, $y'$:
\begin{equation}\label{eq:zprimeyprimeorbitclosing}
    z'\in W^u_{y'}\subset\Min_{r,i},\Fc(z')\in W^u_{\Fc(y')},\cdots,\Fc^{i}(z')\in W^u_{\Fc^i(y')} 
    \text{ for some }0\le i<k.
\end{equation}

If $U\subset\hat{\mathfrak{U}}_{+}$ and $x\in\hat{\mathfrak{U}}_+$ is the sufficient point with quadruple $(l,N,U,\mathcal{C})$ given in \cref{theorem:sufficentptsonMhatplus}, then let $k_{1}(z')\dfn\inf\big\{k'\ge0,\Fcal^{k'}(z')\in\hat{M}_1=\hat{M}_{+}\big\}$. Since $\hat{M}_1=\Fc(\Mrin\smallsetminus\Mrout)\sqcup(\Mrin\cap\Mrout)$, we immediately see the following. \begin{align*}
    k_1(z')=\left\{
    \begin{aligned}
        &0\text{, if } W^u_{y'}\subset \Min_{r,0},\\
        &1\text{, otherwise.}
    \end{aligned}
    \right.
\end{align*}

Since $z'$, $\Fc^{k_1 (z')}(z')$ are consecutive collisions on $\Gamma_r$, for the same reason as \eqref{eq:subcase1-2L4conclusion0} $\frac{\|D\Fc^{k_{1}(z')}_{z'}(dz')\|}{\|dz'\|}\ge1$, and by \eqref{eq:zprimeyprimeorbitclosing} and \cref{lemma:UnstableManifoldInCone}\eqref{item02UnstableManifoldInCone}: $D\Fc^{k_1(z')}_{z'}(dz')\in T_{\Fc^{k_1(z')}(z')}W_{\Fc^{k_1(z')}(y')}^u\in\hat{C}^{u}_{1,\Fc^{k_1(z')}(z')}$. Therefore, with $k\ge k_1(z')$, $\Fc^{k_1(z')}(z'),\Fc^{k}(z')\in\hat{M}_1$, $\Fc^{k}(z')=\hat{F}_1^t(\Fc^{k_1(z')}(z'))$ with some $t\ge0$.

Invoking \cref{Col:EuclideanExpansionMhat1} again gives
\begin{equation}\label{eq:subcase1-2L4conclusion1}
\|D\Fc^{k}_{z'}\bigm|_{T_{z'}W_{y'}^u}\|=\frac{\|dz\|}{\|dz'\|}=\frac{\|D\Fc^{k}_{z'}(dz')\|}{\|dz'\|}=\underbracket{\frac{\|D\Fc^{k}_{1}{z'}(dz')\|}{\|D\Fc^{k_1(z')}_{z'}(dz')\|}}_{\mathclap{\frac{\|D\hat{F}^t_1(D\Fc^{k_1(z')}_{z'}(dz'))\|}{\|D\Fc^{k_1(z')}_{z'}(dz')\|}\text{ ,\cref{Col:EuclideanExpansionMhat1}:}>c_2}}\cdot\underbracket{\frac{\|D\Fc^{k_{1}(z')}_{z'}(dz')\|}{\|dz'\|}}_{\ge1}>c_2>0,
\end{equation} where $c_2$ is from \cref{Col:EuclideanExpansionMhat1}.

If $U\subset\hat{\mathfrak{U}}_{-}$ and $x\in\hat{\mathfrak{U}}_{-}$ is the sufficient point with quadruple $(l,N,U,\mathcal{C})$ given in \cref{theorem:sufficentptsonMhatplus}, then let $k_{2}(z')\dfn\inf\big\{k'\ge0,\Fcal^{k'}(z')\in\hat{M}=\hat{M}_{-}\big\}$. Since $\hat{M}=\Fc^{-1}(\Mrout\smallsetminus\Mrin)\sqcup(\Mrin\cap\Mrout)$, we immediately see the following. \begin{align*}
    k_2(z')=\left\{
    \begin{aligned}
        &0\text{, if } W^u_{y'}\subset \Min_{r,0},\\
        &i-1 \text{, if } W^u_{y'}\subset \Min_{r,i}\text{ with }i\ge1.\\
    \end{aligned}
    \right.
\end{align*}

Since $z',\cdots,\Fc^{k_2 (z')}(z')$ are consecutive collisions in $\Gamma_r$, for the same reason as in \eqref{eq:subcase1-2L4conclusion0} $\frac{\|D\Fc^{k_{2}(z')}_{z'}(dz')\|}{\|dz'\|}\ge1$ and by \eqref{eq:zprimeyprimeorbitclosing} and \cref{lemma:UnstableManifoldInCone}\eqref{item03UnstableManifoldInCone} $D\Fc^{k_2(z')}_{z'}(dz')\in T_{\Fc^{k_2(z')}(z')}W_{\Fc^{k_2(z')}(y')}^u\in C_{\Fc^{k_2(z')}(z')}$.

Since $\Fc^{k_2(z')}(z')\in\hat{M}$ and $\Fc^{k}(z')\in\hat{M}$ with $k\ge k_2(z')$, $\Fc^k(z')=\hat{F}^{t}(\Fc^{k_2(z')}(z'))$ for some $t\ge0$. . Invoking \cref{Col:EuclideanExpansionMhat} gives
\begin{equation}\label{eq:subcase1-2L4conclusion2}
\|D\Fc^{k}_{z'}\bigm|_{T_{z'}W_{y'}^u}\|=\frac{\|dz\|}{\|dz'\|}=\frac{\|D\Fc^{k}_{z'}(dz')\|}{\|dz'\|}=\underbracket{\frac{\|D\Fc^{k}_{z'}(dz')\|}{\|D\Fc^{k_2(z')}_{z'}(dz')\|}}_{\mathclap{=\frac{\|D\hat{F}^t(D\Fc^{k_2(z')}_{z'}(dz'))\|}{\|D\Fc^{k_2(z')}_{z'}(dz')\|}\text{, \cref{Col:EuclideanExpansionMhat}:}>c_0}}\cdot\underbracket{\frac{\|D\Fc^{k_{2}(z')}_{z'}(dz')\|}{\|dz'\|}}_{\ge1}>c_0,
\end{equation}
where constant $c_0$ is from \cref{Col:EuclideanExpansionMhat}.

\textbf{Singularity curve case (2)}: Let $\mathbb{S}=\partial\Mrin\smallsetminus\partial M_R$, i.e., be the line segment as the boundary of $\MRin$ in \cref{fig:MR}. And by symmetry, we assume that $\mathbb{S}$ is the dashed line segment $SW-NE$ in \cref{fig:x1_n1equal1case}.

The local unstable manifold $W^u_{y'}$ is a connected line segment and cannot intersect with $(\mathcal{S}_{-1}\smallsetminus\partial M)\subset\partial\MRin$ at interior point of $W^u_{y'}$ by the definition of $W^u_{y'}$ in \cref{proposition:stablemanifoldthm}. Hence $z'\in W^u_{y'}\subset\MRin$ or $z'\in W^u_{y'}\subset(M_R\smallsetminus\MRin)$. We have the following two subcases of \textbf{case (2)}.

\textbf{Singularity curve subcase (2.1):} $z'\in W^u_{y'}\subset\MRin\cap \mathbb{S}(\xi)$, $dz'\in T_{z'}W^u_{y'}$. 

We also make sure that $\xi$ is small enough to satisfy \begin{equation}\label{eq:xicondition3}
\begin{aligned}
    &\xi<\sin^{-1}(2r/R)-\sin^{-1}{(\frac{r\sin\phi_*}{R})},
\end{aligned}
\end{equation} where $\phistar$, $r$, $R$ are given by the lemon billiard configuration $\mathbb{L}(r,R,\phistar)$.

Note that $W^u_{y'}$ must be contained in one connected component of $\Min_{R,i}$ (see \cref{fig:MrN,fig:x1_multi_collisions_caseC,fig:x1_n1cases_caseC}) with $0\le i<k$. Otherwise, $\Fc^{i-1}(W^u_{y'})$ or $\Fc^{i+1}(W^u_{y'})$, which are the images of $W^u_y$ under some $\Fc^{-1}$ iteration, has at least two components separately contained in $M_R$ and $M_r$. So, it is not connected, which contradicts \cref{proposition:stablemanifoldthm}\eqref{item04stablemanifoldthm} and the definition of the local unstable manifold. 

Each point on $W^u_{y'}$ is the $x_1$ in some $\hat{M}$ return orbit segment given by \eqref{eqPtsFromeqMhatOrbSeg}. By checking \cref{fig:A1contractionCase,fig:x1_n1cases_caseC,fig:A_compactRegion}, \eqref{eq:xicondition3} ensures that $z'=x_1=(\Phi_1,\theta_1)$ is in case (b) or (c) of \eqref{eqMainCases}. 

\textbf{In case (b)} of \eqref{eqMainCases}, $z'\in W^u_{y'}\cap\mathcal{L}$ region of \cref{fig:A1contractionCase}.

By \cref{contraction_region} and using the notations from \cref{lemma:ExpansionLowerBoundB}, we have the following: $z'=x_1$, $\Fc^{-1}(z')=\Fc^{-1}(x_1)=x_0\in\Nout$, $x\dfn\Fc^{-2}(z')=\Fc^{-2}(x_1)\in\Fc^{-1}(\Mrout\smallsetminus\Mrin)\subset\hat{M}$, $(\phi_2,\theta_2)=x_2=\Fc(x_1)\in\Nin$, $x_3=\Fc(x_2)=\Fc^2(x_1)=\Fc^3(x_0)=\Fc^4(x)\in\Fc(\Mrin\smallsetminus\Mrout)$.

By $z'\in W^u_{y'}$, \cref{proposition:stablemanifoldthm}\eqref{item04stablemanifoldthm}: the $\Fc^{-1}$ invariance of $W^u_{y'}$ will give $x\in W_{\Fc^{-2}(y')}^u$ and $D\Fc^{-2}_{z'}(dz')=dx\in T_{x}W^u_{\Fc^{-2}(y')}$. Therefore, by \cref{lemma:UnstableManifoldInCone}\eqref{item03UnstableManifoldInCone}, $dx\in C_x=\big\{(d\phi,d\theta)\bigm|\frac{d\theta}{d\phi}\in[0,1]\big\}$. In the context of \cref{lemma:ExpansionLowerBoundB}, $dx_0=D\Fc_{x}(dx)$, $dz'=dx_1=(d\Phi_1,d\theta_1)=D\Fc^2_{x}(dx)$, $(d\phi_3,d\theta_3)=dx_3=D\Fc_{x_0}^3(dx_0)\in\Int(\mathcal{HQ}_{x_3}(\mathrm{I,III}))=\Int(\big\{(d\phi,d\theta)\bigm|\frac{d\theta}{d\phi}\in[0,1]\big\})\subset \hat{C}^u_{1,x_3}$ by \cref{Proposition:DF4ConeUpperBound} and \eqref{eq:casebdtheta3dphi3ratio} with $\hat{C}^u_{1,x_3}$ defined in \cref{thm:UniformExpansionOnMhat1}. The same as in \cref{lemma:ExpansionLowerBoundA0A1}, we also use the notation $\Bcal^-_i,\Bcal^+_i$ as the before/after collision wave front curvatures at $x_i$, $i=0,1,2,3$. Therefore, by the definition of the $\p-$metric in \cref{def:SectionSetsDefs}, we have the following computations of the expansion with respect to the Euclidean metric, with $\mathcal{V}_1=\frac{d\theta_1}{d\Phi_1}, \mathcal{V}_3=\frac{d\theta_3}{d\phi_3}$:
\begin{equation}\label{eq:L4Case2-1EuclideanMetricRatio1}
\begin{aligned}
    \frac{\|dx_3\|}{\|dx_1\|}=&\frac{|d\phi_3|\sqrt{1+\mathcal{V}^2_3}}{|d\phi_1|\sqrt{1+\mathcal{V}^2_1}}\frac{R\sin\theta_1}{r\sin\theta_3}\frac{\sin\theta_3r\cdot|d\phi_3|\sqrt{1+\mathcal{V}^2_3}}{\sin\theta_1 R\cdot|d\phi_1|\sqrt{1+\mathcal{V}^2_1}}\\
    \overbracket{=}^{\mathllap{\substack{\p-\text{metric in}\\\text{\cref{def:SectionSetsDefs}}}}}&\frac{R\sin{\theta_1}}{r\sin\theta_3}\frac{\sqrt{1+\mathcal{V}^2_3}}{\sqrt{1+\mathcal{V}^2_1}}\frac{\|dx_3\|_{\p}}{\|dx_1\|_{\p}}\overbracket{=}^{\mathclap{\substack{d_1,d_2 \text{ in \cref{def:fpclf}}\\ p(x_2),p(x_3)\text{ are on }\Gamma_r:\,\theta_2=\theta_3\\\mid\\\mid}}}\frac{d_1\sqrt{1+\mathcal{V}^2_3}}{d_2\sqrt{1+\mathcal{V}^2_1}}\frac{\|dx_3\|_{\p}}{\|dx_1\|_{\p}}
    >\frac{d_1}{d_2\sqrt{1+\mathcal{V}^2_1}}\frac{\|dx_3\|_{\p}}{\|dx_1\|_{\p}}
\end{aligned}
\end{equation}
Note that for the arc length parameter $s$ on $\Gamma_R$, $ds_1=Rd\Phi_1$. Invoking \cite[Equation (3.31) and Mirror equation (3.39)]{cb} with the coordinate change: $\varphi=\pi/2-\theta$ and $\mathcal{K}=\frac{-1}{R}$ as the curvature at the boundary $\Gamma_R$ will give the following.
\begin{equation}\label{eq:L4Case2-1EuclideanMetricRatio2}
    \begin{aligned}
        \frac{-\mathcal{V}_1}{R}&=\frac{-1}{R}\frac{d\theta_1}{d\Phi_1}=\frac{d\varphi_1}{ds_1}\overbracket{=}^{\mathclap{\substack{\text{\cite[equation (3.39)]{cb} with}\mathcal{K}=-1/{R}\\\mid}}}\Bcal^{-}_1\cos\varphi_1-\frac{1}{R},\\
        \mathcal{V}_1&=(\Bcal^-_1\cos\varphi_1-1/R)(-R)=1-R\cos\varphi_1\Bcal^-_1\overbracket{=}^{\mathclap{R\cos\varphi_1=R\sin\theta_1=d_1\qquad\quad}}1-d_1\Bcal_1^-
        \overbracket{=}^{\mathclap{\qquad\text{\cite[Eq. (3.31)]{cb}}}}1-d_1\frac{\Bcal^+_0}{1+\tau_0\Bcal^+_0}.
    \end{aligned}
\end{equation}

\eqref{eq:L4Case2-1EuclideanMetricRatio1}, \eqref{eq:L4Case2-1EuclideanMetricRatio2}, \cite[equation (3.40)]{cb}, \cref{proposition:dx3dx2expansion,proposition:IktransionInCaseb} and \cref{lemma:ExpansionLowerBoundB}\eqref{itemExpansionLowerBoundB3} combined will give \begin{equation}\label{eq:L4Case2-1EuclideanMetricRatio3}\begin{aligned}
    \frac{\|dx_3\|}{\|dx_1\|}&\overbracket{>}^{\mathrlap{\text{\eqref{eq:L4Case2-1EuclideanMetricRatio1}, \eqref{eq:L4Case2-1EuclideanMetricRatio2}}}}\quad\frac{d_1}{d_2\sqrt{1+(1-\frac{\Bcal^+_0 d_1}{1+\tau_0\Bcal_0^+})^2}}\frac{\|dx_3\|_{\p}}{\|dx_1\|_{\p}}=\frac{d_1|1+\tau_0\Bcal^+_0|}{d_2\sqrt{(1+\tau_0\Bcal^+_0)^2+(1+\tau_0\Bcal^+_0-d_1\Bcal^+_0)^2}}\frac{\|dx_3\|_{\p}}{\|dx_1\|_{\p}}\\
    &\overbracket{=}^{\mathclap{\substack{\text{\cite[equation (3.40)]{cb}: }\\\frac{\|dx_1\|_{\p}}{\|dx_0\|_{\p}}=|1+\tau_0\Bcal^+_0|}}}\frac{d_1}{d_2\sqrt{(1+\tau_0\Bcal^+_0)^2+(1+\tau_0\Bcal^+_0-d_1\Bcal^+_0)^2}}\frac{\|dx_1\|_{\p}}{\|dx_0\|_{\p}}\frac{\|dx_3\|_{\p}}{\|dx_1\|_{\p}}\\
    &=\frac{d_1}{d_2\sqrt{(1+\tau_0\Bcal^+_0)^2+(1+\tau_0\Bcal^+_0-d_1\Bcal^+_0)^2}}\frac{\|dx_2\|_{\p}}{\|dx_0\|_{\p}}\cdot\underbracket{\frac{\|dx_3\|_{\p}}{\|dx_2\|_{\p}}}_{\mathrlap{\text{\cref{proposition:dx3dx2expansion}: }>1}}\\
    &\overbracket{>}^{\mathllap{\text{\cref{lemma:ExpansionLowerBoundB}\eqref{itemExpansionLowerBoundB3}}}}\frac{d_1}{d_2\sqrt{(1+\tau_0\Bcal^+_0)^2+(1+\tau_0\Bcal^+_0-d_1\Bcal^+_0)^2}}|\mathcal{I}(k)|,\text{ where } \mathcal{I}(k)\text{ is from \cref{lemma:ExpansionLowerBoundB}\eqref{itemExpansionLowerBoundB3} with }k\in[1,4/3].\\
    &\overbracket{>}^{\mathllap{\text{\cref{proposition:IktransionInCaseb}: }\mathcal{I}(k)<-0.05}}\frac{0.05d_1}{d_2\sqrt{(1+\tau_0\Bcal^+_0)^2+(1+\tau_0\Bcal^+_0-d_1\Bcal^+_0)^2}}=\frac{0.05}{\sqrt{\big(1+\tau_0\Bcal^+_0\big)^2\frac{d_2^2}{d_1^2}+\big(1+(\tau_0-d_1)\Bcal^+_0\big)^2\frac{d_2^2}{d_1^2}}}
\end{aligned}
\end{equation}
Note that in case (b), \cref{proposition:lengthsrestrictionincaseB}\eqref{eq:lengthrestrictions1}: $0.5r\sin\phi_*<d_1<2r$, \cref{lemma:ExpansionLowerBoundB}\eqref{itemExpansionLowerBoundB3}: $\Bcal^+_0\in[\frac{-4}{3d_0},\frac{-1}{d_0}]$. Hence, by \cref{remark:deffpclf}, $\tau_0<2d_0, \tau_0<2d_1$, thus $\tau_0<d_0+d_1$ and

We have the following estimate:
\begin{align*}
    1+\tau_0\Bcal^+_0&\overbracket{\in}^{\tau_0\in(0,2d_0),\Bcal^+_0<0}(1+2d_0\Bcal^+_0,1)\overbracket{\subset}^{\Bcal_0^+\ge\frac{-4}{3d_0}}(1+\frac{-4\times 2d_0}{3d_0},1)=(-\frac{5}{3},1)\\
    1+(\tau_0-d_1)\Bcal^+_0&\overbracket{\in}^{-d_1<\tau_0-d_1<d_0,\Bcal^+_0<0}(1+d_0\Bcal^+_0,1-d_1\Bcal^+_0)\overbracket{\subset}^{\frac{-4}{3d_0}<\Bcal^+_0\le\frac{-1}{d_0}}(1-\frac{4d_0}{3d_0},1+\frac{4d_1}{3d_0})=(-\frac{1}{3},1+\frac{4d_1}{3d_0})\\
    \text{Therefore, }|\frac{d_2}{d_1}(1+\tau_0\Bcal^+_0)|^2&<|\frac{5}{3}|^2(\frac{d_2}{d_1})^2\overbracket{<}^{\text{\cref{proposition:lengthsrestrictionincaseB}\eqref{eq:casebghd0lengthrestrictions5}},d_1>0.5r\sin\phi_*}\frac{25}{9}(\frac{1.0003r\sin\phi_*}{0.5r\sin\phi_*})^2<11.2\\
    \frac{d_2^2}{d_1^2}(1+(\tau_0-d_1)\Bcal^+_0)^2&<\max\Big\{\frac{d_2^2}{9d_1^2},(d_2+\frac{4d_2}{3d_0}d_1)^2\frac{1}{d_1^2}\Big\}\\
    &\overbracket{<}^{\mathllap{\substack{\text{\cref{proposition:lengthsrestrictionincaseB}\eqref{eq:casebghd0lengthrestrictions5},\eqref{eq:casebghd0lengthrestrictions2}},\\d_1>0.5r\sin\phi_*}}}\max\Big\{\frac{1.0003^2r^2\sin^2\phi_*}{9\cdot(0.5^2r^2\sin^2\phi_*)},(\frac{1.0003r\sin\phi_*}{0.5r\sin\phi_*}+\frac{4\cdot1.0003r\sin\phi_*}{3\cdot0.9997r\sin\phi_*})^2\Big\}<11.2
\end{align*}
Hence, \eqref{eq:L4Case2-1EuclideanMetricRatio3} implies that \[\frac{\|dx_3\|}{\|dx_1\|}>\frac{0.05}{\sqrt{11.2+11.2}}>0.01\]

\textbf{In case (c)} of $\hat{F}$ return orbit segment \cref{def:MhatReturnOrbitSegment} falling in case \eqref{eqMainCases}, $z'\in W^u_{y'}\subset\Min_{R,n_1}$ region in \cref{fig:x1_n1cases_caseC}, for some $n_1\ge1$.

By \cref{contraction_region} and using the notations from \cref{lemma:AformularForLowerBoundOfExpansionForn1largerthan0}, we have the following: $z'=x_1$, $\Fc^{-1}(z')=\Fc^{-1}(x_1)=x_0\in\Nout$, $x\dfn\Fc^{-2}(z')=\Fc^{-2}(x_1)\in\Fc^{-1}(\Mrout\smallsetminus\Mrin)\subset\hat{M}$, $(\phi_2,\theta_2)=x_2=\Fc^{n_1+1}(x_1)\in\Nin$, $x_3=\Fc(x_2)=\Fc^{n_1+2}(x_1)=\Fc^{n_1+3}(x_0)=\Fc^{n_1+4}(x)\in\Fc(\Mrin\smallsetminus\Mrout)\subset\hat{M}_1$.

By $z'\in W^u_{y'}$, \cref{proposition:stablemanifoldthm}\eqref{item04stablemanifoldthm}: the $\Fc^{-1}$ invariance of $W^u_{y'}$ will give $x\in W_{\Fc^{-2}(y')}^u$ and $D\Fc^{-2}_{z'}(dz')=dx\in T_{x}W^u_{\Fc^{-2}(y')}$. Therefore, by \cref{lemma:UnstableManifoldInCone}\eqref{item03UnstableManifoldInCone}, $dx\in C_x=\big\{(d\phi,d\theta)\bigm|\frac{d\theta}{d\phi}\in[0,1]\big\}$. In the context of \cref{lemma:ExpansionLowerBoundB}, $dx_0=D\Fc_{x}(dx)$, $dz'=dx_1=(d\Phi_1,d\theta_1)=D\Fc^2_{x}(dx)$, $(d\phi_3,d\theta_3)=dx_3=D\Fc_{x_0}^{n_1+3}(dx_0)$ satisfies $0<\frac{d\theta_3}{d\phi_2}<\lambda_2(r,R,\phistar)$ by \cref{theorem:B3minusdx3dphi3ratiobound} and thus $dx_3\in\hat{C}^u_{1,x_3}$ with $\hat{C}^u_{1,x_3}$ defined in \cref{thm:UniformExpansionOnMhat1}. The same as in \cref{lemma:AformularForLowerBoundOfExpansionForn1largerthan0}, we also use the notation $\Bcal^-_i,\Bcal^+_i$ as the before/after collision wave front curvatures at $x_i$, $i=0,1,2,3$. Therefore, by the $\p-$metric definition in \cref{def:SectionSetsDefs}, similar to \eqref{eq:L4Case2-1EuclideanMetricRatio1} with $\mathcal{V}_1=\frac{d\theta_1}{d\Phi_1}, \mathcal{V}_3=\frac{d\theta_3}{d\phi_3}$ we have the following computations for the $dx_3$, $dx_1$ Euclidean metric expansion.

\begin{equation}\label{eq:L4Case2-1EuclideanMetricRatio4}
\begin{aligned}
    \frac{\|dx_3\|}{\|dx_1\|}=&\frac{R\sin{\theta_1}}{r\sin\theta_3}\frac{\sqrt{1+\mathcal{V}^2_3}}{\sqrt{1+\mathcal{V}^2_1}}\frac{\|dx_3\|_{\p}}{\|dx_1\|_{\p}}=\frac{d_1\sqrt{1+\mathcal{V}^2_3}}{d_2\sqrt{1+\mathcal{V}^2_1}}\frac{\|dx_3\|_{\p}}{\|dx_1\|_{\p}}=\frac{d_1\sqrt{1+\mathcal{V}^2_3}}{d_2\sqrt{1+\mathcal{V}^2_1}}\frac{\|dx_3\|_{\p}}{\|dx_1\|_{\p}}\\
    \overbracket{=}^{\mathllap{\text{\eqref{eq:L4Case2-1EuclideanMetricRatio2}}}}&\frac{d_1\sqrt{1+\mathcal{V}^2_3}}{d_2\sqrt{1+(1-d_1\frac{\Bcal^+_0}{1+\tau_0\Bcal^+_0})^2}}\frac{\|dx_3\|_{\p}}{\|dx_1\|_{\p}}=\frac{d_1|1+\tau_0\Bcal^+_0|\sqrt{1+\mathcal{V}^2_3}}{d_2\sqrt{(1+\tau_0\Bcal^+_0)^2+\big(1+(\tau_0-d_1)\Bcal^+_0\big)^2}}\frac{\|dx_3\|_{\p}}{\|dx_1\|_{\p}}\\
    \overbracket{=}^{\mathclap{\substack{\text{\cite[equation (3.40)]{cb}: }\\\frac{\|dx_1\|_{\p}}{\|dx_0\|_{\p}}=|1+\tau_0\Bcal^+_0|}}}&\frac{d_1\sqrt{1+\mathcal{V}^2_3}}{d_2\sqrt{(1+\tau_0\Bcal^+_0)^2+\big(1+(\tau_0-d_1)\Bcal^+_0\big)^2}}\overbracket{\frac{\|dx_1\|_{\p}}{\|dx_0\|_{\p}}\frac{\|dx_3\|_{\p}}{\|dx_1\|_{\p}}}^{\mathclap{\text{\cref{lemma:AformularForLowerBoundOfExpansionForn1largerthan0}\eqref{eq:Expansionx3x0}: }=E(n_1,\Bcal_0^+)}}\\
    >&\frac{d_1}{d_2\sqrt{(1+\tau_0\Bcal^+_0)^2+\big(1+(\tau_0-d_1)\Bcal^+_0\big)^2}}E(n_1,\Bcal^+_0)\\
\end{aligned}
\end{equation}
By \cref{remark:deffpclf} $0<\tau_0<2d_0$, $-d_1<\tau_0-d_1<d_1$. In case (c) $d_1$ satisfies $d_1<r\sin\phi_*$ by $x_1$, $\Fc(x_1)$ are consecutive collisions on $\Gamma_R$: $2d_1<2R\sin\Phi_*=2r\sin\phi_*$. And by \cref{proposition:casecparametersestimate}\ref{itemcasecparametersestimate3}: $\frac{r\sin\phi_*}{n_1+3}<d_1$. By \cref{lemma:AformularForLowerBoundOfExpansionForn1largerthan0}\eqref{itemAformularForLowerBoundOfExpansionForn1largerthan01}, $\Bcal^+_0\in[-\frac{-4}{3d_0},\frac{-1}{d_0}]$. We get the following estimate.
\begin{align*}
1+\tau_0\Bcal^+_0&\overbracket{\in}^{\tau_0\in(0,2d_0),\Bcal^+_0<0}(1+2d_0\Bcal^+_0,1)\overbracket{\subset}^{\Bcal_0^+\ge\frac{-4}{3d_0}}(1+\frac{-4\times 2d_0}{3d_0},1)=(-\frac{5}{3},1)\\
1+(\tau_0-d_1)\Bcal^+_0&\overbracket{\in}^{\tau_0\in(0,2d_1),\Bcal^+_0<0}(1+d_1\Bcal^+_0,1-d_1\Bcal^+_0)\overbracket{\subset}^{\Bcal^+_0\in[-\frac{-4}{3d_0},\frac{-1}{d_0}]}(1+\frac{-4d_1}{3d_0},1+\frac{d_1}{d_0})\\
&\overbracket{\subset}^{d_1<0.5r\sin\phi_*}
(1-\frac{2r\sin\phi_*}{d_0},1+\frac{0.5r\sin\phi_*}{d_0})\overbracket{\subset}^{\mathrlap{\text{\cref{proposition:casecparametersestimate}\ref{itemcasecparametersestimate5}: }r\sin\phi_*/d_0<1.0026}}(-1.006,1.51)\\
\text{Therefore, }&\frac{d_1}{d_2\sqrt{(1+\tau_0\Bcal^+_0)^2+\big(1+(\tau_0-d_1)\Bcal^+_0\big)^2}}E(n_1,\Bcal^+_0)>\frac{d_1}{d_2\sqrt{1.51^2+25/9}}E(n_1,\Bcal^+_0)\\
&\overbracket{>}^{\mathllap{\substack{\text{\cref{proposition:casecparametersestimate}\ref{itemcasecparametersestimate3}\ref{itemcasecparametersestimate5}: }\\d_1>\frac{r\sin\phi_*}{n_1+3},r\sin\phi_*>0.9975d_2}}}\frac{0.9975E(n_1,\Bcal^+_0)}{\sqrt{1.51^2+25/9}(n_1+3)}>0.44\frac{E(n_1,\Bcal^+_0)}{(n_1+3)}
\end{align*}
Note that by \cref{proposition:casecexpansionn1eq1,proposition:casecexpansionn1eq2,proposition:casecexpansionn1eq3,proposition:casecexpansionn1ge4}, for each finite $n_1$, $E(n_1,\Bcal^+_0)>0.9$ and by \cref{proposition:casecexpansionn1ge4}, $\lim_{n_1\rightarrow\infty}\frac{E(n_1,\Bcal^+_0)}{n_1+3}=\infty$. We see that there is a uniform lower bound $C_4>0$ such that $\frac{\|dx_3\|}{\|dz'\|}=\frac{\|dx_3\|}{\|dx_1\|}>0.44\frac{E(n_1,\Bcal^+_0)}{(n_1+3)}>C_4$.

We can summarize in both case (b) and (c), with $z'\in W^u_{y'}\subset \Min_{R,n_1}$, $n_1\ge0$, $\hat{M}_1\ni x_3=\Fc^{n_1+1}(z')$, $dx_3=D\Fc_{z'}^{n_1+1}(dz')$ satisfying $\frac{\|dx_3\|}{\|dz'\|}>\min\big\{C_4,0.01\big\}$. 

If $U\subset\hat{\mathfrak{U}}_{+}$ and $x\in\hat{\mathfrak{U}}_+$ is the sufficient point with quadruple $(l,N,U,\mathcal{C})$ given in \cref{theorem:sufficentptsonMhatplus}, then the orbit segment $z',\Fc(z'),\cdots,\Fc^{n_1+1}(z')=x_3,\Fc^{n_1+2}(z')=x_3,\cdots,\Fc^{k}(z')=z\in T_zW_y^u\subset\hat{\mathfrak{U}}_{+}\subset\hat{M}_1=\hat{M}_{+}$ and $dx_3\in\hat{C}^u_{1,x_3}$. 
Since $x_3, \Fc^k(z')\in\hat{M}_1$, $\Fc^k(z')\in\hat{M}_1=\hat{F}_1^t(x_3)$ with some $t\ge0$. Hence, by \cref{thm:UniformExpansionOnMhat1,Col:EuclideanExpansionMhat1}, $\frac{\|D\Fc_{z'}^k(dz')\|}{\|dx_3\|}=\frac{\|D\hat{F}_1^t(dx_4)\|}{\|dx_4\|}>c_2$. Therefore, with $c_2>0$ from \cref{Col:EuclideanExpansionMhat1},
\begin{equation}\label{eq:subcase2-1L4conclusion1}
    \frac{\|D\Fc_{z'}^{k}(dz')\|}{\|dz'\|}=\frac{\|D\Fc^{k}_{z'}(dz')\|}{\|dx_3\|}\frac{\|dx_3\|}{\|dz'\|}>c_2\cdot\min\big\{C_4,0.01\big\}.
\end{equation}

If $U\subset\hat{\mathfrak{U}}_{-}$ and $x\in\hat{\mathfrak{U}}_{-}$ is the sufficient point with quadruple $(l,N,U,\mathcal{C})$ given in \cref{theorem:sufficentptsonMhatplus}, then by \cref{contraction_region} in erither case (b) or (c) we can pick $x\dfn\Fc^{-2}(z')\in\Fc^{-1}(\Mrout\smallsetminus\Mrin), dx=D\Fc^{-2}(dz')$ and let $x_4\dfn\hat{F}(x)$ and $dx_4=D\hat{F}_x(dx)$ so that the orbit segment $z',\Fc(z'),\cdots,\Fc^{n_1+1}(z')=x_2,\underbracket{\Fc^{n_1+2}(z')=x_3\in\hat{M}_1,\cdots,x_4\in\hat{M}}_{\text{consecutive collisions on }\Gamma_r},\cdots,\Fc^{k}(z')=z\in T_zW_y^u\subset\hat{\mathfrak{U}}_{-}\subset\hat{M}=\hat{M}_{-}$ and $dx_3\in\hat{C}^u_{1,x_3}$ satisfy $\frac{\|dx_4\|}{\|dx_3\|}>1$ and $\frac{\|dx_3\|}{\|dz'\|}>\min\big\{C_4,0.01\big\}$.
Since $x_4,\Fc^k(z')\in\hat{M}$, $\Fc^k(z')\in\hat{M}=\hat{F}^t(x_4)$ for some $t\ge0$. Then by \cref{thm:UniformExpansionOnMhat,Col:EuclideanExpansionMhat}, $\frac{\|D\Fc_{z'}^k(dz')\|}{\|dx_4\|}=\frac{\|D\hat{F}^t(dx_4)\|}{\|dx_4\|}>c_0$. Therefore, with $c_0>0$ from \cref{Col:EuclideanExpansionMhat},
\begin{equation}\label{eq:subcase2-1L4conclusion2}
    \frac{\|D\Fc_{z'}^{k}(dz')\|}{\|dz'\|}=\frac{\|D\Fc^{k}_{z'}(dz')\|}{\|dx_4\|}\frac{\|dx_4\|}{\|dx_3\|}\frac{\|dx_3\|}{\|dz'\|}>c_0\cdot\min\big\{C_4,0.01\big\}.
\end{equation}

\textbf{Singularity curve subcase (2.2):} $z'\in W^u_{y'}\subset(M_R\smallsetminus\MRin)\cap \mathbb{S}(\xi)$, $dz'\in T_{z'}W^u_{y'}$.

Note that $z'\in W^u_{y'}\subset\Fc^j(\Min_{R,n_1})$ with some $0<j\le n_1$, otherwise some $\Fc^{-1}$ iterations image of $W^u_y$ will be separate, i.e., it would have components in $M_R$ and $M_r$ which contradicts \cref{proposition:stablemanifoldthm}\eqref{item03stablemanifoldthm}. Hence $\Fc^{-j}(z')=x_1=(\Phi_1,\theta_1)$ in some $\hat{M}$ return orbit segment defined by \eqref{eqPtsFromeqMhatOrbSeg} in case (c) of \eqref{eqMainCases}. We also denote $z'=x_{1,j}=(\Phi_{1,j},\theta_{1,j})$ and $dx_{1,j}=dz'=(d\Phi_{1,j},d\theta_{1,j})\in T_{z'}W^u_{y'}$. Since $x_1,\cdots,\Fc^j(x_1)=x_{1,j}$ are consecutive collisions on arc $\Gamma_R$, $\theta_1=\theta_{1,j}$.

By \cref{contraction_region} and using the notations from \cref{lemma:AformularForLowerBoundOfExpansionForn1largerthan0}, we have the following: $\Fc^{-j}(z')=x_1$, $\Fc^{-j-1}(z')=\Fc^{-j-1}(x_1)=x_0\in\Nout$, $x\dfn\Fc^{-j-2}(z')=\Fc^{-j-2}(x_1)\in\Fc^{-1}(\Mrout\smallsetminus\Mrin)\subset\hat{M}$, $(\phi_2,\theta_2)=x_2=\Fc^{n_1+1}(x_1)\in\Nin$, $(\phi_3,\theta_3)=x_3=\Fc(x_2)=\Fc^{n_1+2}(x_1)=\Fc^{n_1+3}(x_0)=\Fc^{n_1+4}(x)\in\Fc(\Mrin\smallsetminus\Mrout)$.
Since $x_2$, $x_3$ are collisions on arc $\Gamma_r$, $\theta_3=\theta_2$.

By $z'\in W^u_{y'}$, \cref{proposition:stablemanifoldthm}\eqref{item04stablemanifoldthm}: the $\Fc^{-1}$ invariance of $W^u_{y'}$ will give $x\in W_{\Fc^{-j-2}(y')}^u$ and $D\Fc^{-j-2}_{z'}(dz')=dx\in T_{x}W^u_{\Fc^{-j-2}(y')}$. Therefore, by \cref{lemma:UnstableManifoldInCone}\eqref{item03UnstableManifoldInCone}, $dx\in C_x=\big\{(d\phi,d\theta)\bigm|\frac{d\theta}{d\phi}\in[0,1]\big\}$. In the context of \cref{lemma:ExpansionLowerBoundB}, $dx_0=D\Fc_{x}(dx)$, $dx_1=(d\Phi_1,d\theta_1)=D\Fc^2_{x}(dx)$, $D\Fc^j_{x_1}(dx_1)=dx_{1,j}=dz'$,
$(d\phi_3,d\theta_3)=dx_3=D\Fc_{x_0}^{n_1+2}(dx_0)$ satisfies $0<\frac{d\theta_3}{d\phi_3}=\mathcal{V}_3<\lambda_2(r,R,\phistar)$ by \cref{theorem:B3minusdx3dphi3ratiobound}, and thus $dx_3\in\hat{C}^u_{1,x_3}$ with $\hat{C}^u_{1,x_3}$ defined in \cref{thm:UniformExpansionOnMhat1}. The same as in \cref{lemma:AformularForLowerBoundOfExpansionForn1largerthan0}, we also use the notation $\Bcal^-_i,\Bcal^+_i$ as the before/after collision wave front curvatures at $x_i$, $i=0,1,2,3$ and $\Bcal^{-}_{1,j},\Bcal^{+}_{1,j}$ as the before/after collision wave front curvatures at $x_{1,j}$. Therefore, by the $\p-$metric definition in \cref{def:SectionSetsDefs}, similar to \eqref{eq:L4Case2-1EuclideanMetricRatio1} we have the following computations for the $dx_3$, $dx_{1,j}$ Euclidean metric expansion.
\begin{equation}\label{eq:L4Case2-2EuclideanMetricRatio1}
\begin{aligned}
    \frac{\|dx_3\|}{\|dz'\|}&=\frac{\|dx_3\|}{\|dx_{1,j}\|}=\frac{R\sin\theta_{1,j}}{r\sin\theta_3}\frac{\sqrt{1+\mathcal{V}^2_{3}}}{\sqrt{1+\mathcal{V}^2_{1,j}}}\frac{\|dx_3\|_{\p}}{\|dx_{1,j}\|_{\p}},\text{ where }\mathcal{V}_3=\frac{d\theta_3}{d\phi_3}, \mathcal{V}_{1,j}=\frac{d\theta_{1,j}}{d\Phi_{1,j}}\\
    &\overbracket{=}^{\mathllap{\theta_2=\theta_3,\theta_1=\theta_{1,j}}}\frac{d_1}{d_2}\frac{\sqrt{1+\mathcal{V}^2_{3}}}{\sqrt{1+\mathcal{V}^2_{1,j}}}\frac{\|dx_3\|_{\p}}{\|dx_{1,j}\|_{\p}},
\end{aligned}
\end{equation}

With the arc length parameter $s$ on the $\Gamma_R$, $ds_{1,j}=Rd\Phi_{1,j}$. Also note that $\theta_{1,j}=\theta_1$ since the collisions $x_1,\cdots, x_{1,j}$ are the consecutive collisions on the same arc $\Gamma_R$. Invoking \cite[Equation (3.31) and Mirror equation (3.39)]{cb} with the coordinate change: $\varphi=\pi/2-\theta$ and $\mathcal{K}=\frac{-1}{R}$ as the curvature at the boundary $\Gamma_R$ will give the following.
\begin{equation}\label{eq:L4Case2-2EuclideanMetricRatio2}
    \begin{aligned}
        \frac{-\mathcal{V}_{1,j}}{R}&=\frac{-1}{R}\frac{d\theta_{1,j}}{d\Phi_{1,j}}=\frac{d\varphi_{1,j}}{ds_{1,j}}\qquad\qquad\overbracket{=}^{\mathclap{\substack{\text{\cite[equation (3.39)]{cb} with}\\\mathcal{K}=\frac{-1}{R}}}}\qquad\qquad\Bcal^{+}_{1,j}\cos\varphi_{1,j}+\frac{1}{R},\\
        \mathcal{V}_{1,j}&=(\Bcal^+_{1,j}\cos\varphi_{1,j}+1/R)(-R)=-1-R\cos\varphi_{1,j}\Bcal^-_{1,j}\overbracket{=}^{\mathclap{\substack{R\cos\varphi_{1,j}=R\sin\theta_{1,j}=R\sin\theta_1=d_1\\\Bcal^{-}_{1,j}=\Bcal^{+}_{1,j}+\frac{2}{d_1}}}}\quad-3-d_1\Bcal_{1,j}^+,\\
        \Bcal^+_{1,j}&\overbracket{=}^{\mathllap{\text{\cite[exercise 8.28]{cb}}}}-\frac{1}{d_1}+\frac{1}{-2jd_1+\frac{1}{\frac{-1}{d_1}+\Bcal^-_1}}=\frac{-1}{d_1}+\frac{\Bcal^-_1-\frac{1}{d_1}}{(2j+1)-2jd_1\Bcal^-_1},\\
        \text{Therefore, } \mathcal{V}_{1,j}&=-3-d_1\big(\frac{-1}{d_1}+\frac{\Bcal^-_1-\frac{1}{d_1}}{(2j+1)-2jd_1\Bcal^-_1}\big)=-2+\frac{-\Bcal_1^-d_1+1}{(2j+1)-2jd_1\Bcal^-_1}
    \end{aligned}
\end{equation}
\eqref{eq:L4Case2-2EuclideanMetricRatio2}, \eqref{eq:L4Case2-2EuclideanMetricRatio1}, \cite[exercise 8.29]{cb} and \cite[equation (3.40)]{cb} give the following. 
\begin{equation}\label{eq:L4Case2-2EuclideanMetricRatio3}
\begin{aligned}
    &\frac{\|dx_3\|}{\|dx_{1,j}\|}\qquad\overbracket{=}^{\mathclap{\text{\eqref{eq:L4Case2-2EuclideanMetricRatio2}, \eqref{eq:L4Case2-2EuclideanMetricRatio1}}}}\qquad\frac{d_1\sqrt{1+\mathcal{V}^2_3}}{d_2\sqrt{1+\big(-2+\frac{1-\Bcal^-_1d_1}{1+2j-2jd_1\Bcal_1^-}\big)^2}}\frac{\|dx_3\|_{\p}}{\|dx_{1,j}\|_{\p}}\\
    \overbracket{=}^{\mathrlap{\substack{\text{\cite[exercise 8.29]{cb}: }\\\frac{\|dx_{1,j}\|_{\p}}{\|dx_1\|_{\p}}=|1-2jd_1(\Bcal_1^--\frac{1}{d_1})|}}}&\qquad\frac{d_1\sqrt{1+\mathcal{V}^2_3}}{d_2\sqrt{(1+2j-2jd_1\Bcal_1^-)^2+\big(-2\cdot(1+2j-2jd_1\Bcal_1^-)+(1-\Bcal^-_1d_1)\big)^2}}\frac{\|dx_{1,j}\|_{\p}}{\|dx_1\|_{\p}}\frac{\|dx_3\|_{\p}}{\|dx_{1,j}\|_{\p}}\\
    \overbracket{=}^{\mathrlap{\substack{\text{\cite[equation (3.31)]{cb}: }\\\Bcal^-_1=\frac{\Bcal^+_0}{1+\tau_0\Bcal^+_0}}}}&\qquad\frac{d_1\sqrt{1+\mathcal{V}^2_3}\bigm|1+\tau_0\Bcal^+_0\bigm|}{d_2\sqrt{\big((1+2j)(1+\tau_0\Bcal^+_0)-2jd_1\Bcal^+_0\big)^2+\big[-2\cdot\big((1+2j)(1+\tau_0\Bcal^+_0)-2jd_1\Bcal^+_0\big)+(\tau_0-d_1)\Bcal^+_0+1\big]^2}}\frac{\|dx_3\|_{\p}}{\|dx_1\|_{\p}}
    \end{aligned}
\end{equation}
Let $a=(1+2j)(1+\tau_0\Bcal^+_0)-2jd_1\Bcal^+_0$, $b=1+(\tau_0-d_1)\Bcal^+_0$, \cite[equation (3.40)]{cb}: $\frac{\|dx_1\|_{\p}}{\|dx_0\|_{\p}}=|1+\tau_0\Bcal^+_0|$, \eqref{eq:L4Case2-2EuclideanMetricRatio3} gives 
\begin{equation}\label{eq:L4Case2-2EuclideanMetricRatio4}
    \frac{\|dx_3\|}{\|dx_{1,j}\|}=\frac{d_1\sqrt{1+\mathcal{V}^2_3}}{d_2\sqrt{a^2+(b-2a)^2}}\frac{\|dx_3\|_{\p}}{\|dx_{0}\|_{\p}}
\end{equation}

By \cref{remark:deffpclf} $0<\tau_0<2d_0$, $-d_1<\tau_0-d_1<d_1$. In case (c) $d_1$ satisfies $d_1<r\sin\phi_*$ by $x_1$, $\Fc(x_1)$ are consecutive collisions on $\Gamma_R$: $2d_1<2R\sin\Phi_*=2r\sin\phi_*$. And by \cref{proposition:casecparametersestimate}\ref{itemcasecparametersestimate3}: $\frac{r\sin\phi_*}{n_1+3}<d_1$. By \cref{lemma:AformularForLowerBoundOfExpansionForn1largerthan0}\eqref{itemAformularForLowerBoundOfExpansionForn1largerthan01}, $\Bcal^+_0\in[-\frac{-4}{3d_0},\frac{-1}{d_0}]$. We get the following estimate.
\begin{align*}   1+\tau_0\Bcal^+_0\overbracket{\in}^{\mathllap{\tau_0\in(0,2d_0),\Bcal^+_0<0}}(1+2d_0\Bcal^+_0,1)\overbracket{\subset}^{\mathllap{\Bcal_0^+\ge\frac{-4}{3d_0}}}&(1+\frac{-4\times 2d_0}{3d_0},1)=(-\frac{5}{3},1),\\
\text{therefore }\beta=1+(\tau_0-d_1)\Bcal^+_0\overbracket{\in}^{\mathclap{\tau_0\in(0,2d_1),\Bcal^+_0<0}}&(1+d_1\Bcal^+_0,1-d_1\Bcal^+_0)\overbracket{\subset}^{\Bcal^+_0\in[-\frac{-4}{3d_0},\frac{-1}{d_0}]}(1+\frac{-4d_1}{3d_0},1+\frac{d_1}{d_0})\\
\overbracket{\subset}^{\mathllap{d_1<0.5r\sin\phi_*}}&(1-\frac{2r\sin\phi_*}{d_0},1+\frac{0.5r\sin\phi_*}{d_0})\qquad\overbracket{\subset}^{\mathrlap{\substack{\text{\cref{proposition:casecparametersestimate}\ref{itemcasecparametersestimate5}: }\\ r\sin\phi_*/d_0<1.0026}}}(-1.006,1.51).\\
\bigm|\alpha\bigm|=\bigm|(1+2j)(1+\tau_0\Bcal^+_0)-2jd_1\Bcal^+_0\bigm|\overbracket{\le}^{\mathllap{\substack{1>1+\tau_0\Bcal^+_0>-5/3\\\frac{-4}{3d_0}\le\Bcal^+_0\le\frac{-1}{d_0}}}}&\frac{5}{3}\bigm|1+2j\bigm|+\frac{8jd_1}{3d_0}\overbracket{\le}^{\mathclap{j\le n_1}}\quad\frac{5}{3}\bigm|1+2n_1\bigm|+\frac{8n_1d_1}{3d_0}\\
\overbracket{<}^{\mathllap{\text{\cref{proposition:casecparametersestimate}\ref{itemcasecparametersestimate4}\ref{itemcasecparametersestimate5}}}}&\frac{5}{3}(1+2n_1)+\frac{8(r\sin\phi_*+\frac{8.1r}{2R}r\sin\phi_*)}{3\cdot(r\sin\phistar)/1.0026}\overbracket{<}^{\mathrlap{R>1700r}}3.4n_1+4.4
\end{align*}

By the fact: $\forall a,b\in \mathbb{R}$, $\sqrt{a^2+(b-2a)^2}\le|a|+|b-2a|\le 3|a|+|b|$, the right-hand side of the last equation in \eqref{eq:L4Case2-2EuclideanMetricRatio4} satisfies
\begin{align*}
    \frac{d_1\sqrt{1+\mathcal{V}^2_3}}{d_2\sqrt{a^2+(b-2a)^2}}\frac{\|dx_3\|_{\p}}{\|dx_0\|_{\p}}\ge&\frac{d_1}{d_2(3|a|+|b|)}\frac{\|dx_3\|_{\p}}{\|dx_0\|_{\p}}\qquad\overbracket{=}^{\mathclap{\text{\cref{lemma:AformularForLowerBoundOfExpansionForn1largerthan0}\eqref{itemAformularForLowerBoundOfExpansionForn1largerthan04}}}}\frac{d_1}{d_2(3|\alpha|+|\beta|)}E(n_1,\Bcal^+_0)>\frac{d_1E(n_1,\Bcal^+_0)}{d_2\big(3\cdot(3.4n_1+4.4)+1.51\big)}\\
    \overbracket{>}^{\mathllap{\substack{\text{\cref{proposition:casecparametersestimate}\ref{itemcasecparametersestimate3}:}d_1>\frac{r\sin\phistar}{n_1+2}\\\text{and \ref{itemcasecparametersestimate5}}: d_2<\frac{r\sin\phistar}{0.9975}}}}&\frac{0.9975E(n_1,\Bcal_0^+)}{(n_1+2)(10.2n_1+14.8)}.
\end{align*}
By \cref{proposition:casecexpansionn1ge4,proposition:casecexpansionn1eq3,proposition:casecexpansionn1eq2,proposition:casecexpansionn1eq1}, $\frac{0.9975}{(n_1+2)(10.2n_1+14.8)}E(n_1,\Bcal_0^+)$ is positive and bounded away from $0$ for each $n_1\ge1$. Especially by \cref{proposition:casecexpansionn1ge4}, $\liminf_{n_1\rightarrow\infty}\frac{0.9975}{(n_1+2)(10.2n_1+14.8)}E(n_1,\Bcal_0^+)>constant_5$ for some constant $constant_5>0$. Therefore, there is some uniform constant $C_5>0$ such that $\frac{0.9975}{(n_1+2)(10.2n_1+14.8)}E(n_1,\Bcal_0^+)>C_5$ thus by \eqref{eq:L4Case2-2EuclideanMetricRatio1},\eqref{eq:L4Case2-2EuclideanMetricRatio2} and \eqref{eq:L4Case2-2EuclideanMetricRatio3} $\frac{\|dx_3\|}{\|dz'\|}=\frac{\|dx_3\|}{\|dx_{1,j}\|}>C_5$ for some uniform constant $C_5>0$.

If $U\subset\hat{\mathfrak{U}}_{+}$ and $x\in\hat{\mathfrak{U}}_+$ is the sufficient point with quadruple $(l,N,U,\mathcal{C})$ given in \cref{theorem:sufficentptsonMhatplus}, then the same reasoning as \eqref{eq:subcase2-1L4conclusion1} gives the following.

\begin{equation}\label{eq:subcase2-2L4conclusion1}
    \frac{\|D\Fc_{z'}^{k}(dz')\|}{\|dz'\|}=\frac{\|D\Fc^{k}_{z'}(dz')\|}{\|dx_3\|}\frac{\|dx_3\|}{\|dz'\|}>c_2\cdot C_5,
\end{equation} where $c_2>0$ from \cref{Col:EuclideanExpansionMhat1}.

If $U\subset\hat{\mathfrak{U}}_{-}$ and $x\in\hat{\mathfrak{U}}_{-}$ is the sufficient point with quadruple $(l,N,U,\mathcal{C})$ given in \cref{theorem:sufficentptsonMhatplus}, then the same reasoning as \eqref{eq:subcase2-1L4conclusion2} gives the following. 

\begin{equation}\label{eq:subcase2-2L4conclusion2}
    \frac{\|D\Fc_{z'}^{k}(dz')\|}{\|dz'\|}>c_0\cdot C_5,
\end{equation} where $c_0>0$ is from \cref{Col:EuclideanExpansionMhat}.

In all subcases (1.1), (1.2), (2.1) and (2.2), we choose $\xi$ satisfying \eqref{eq:xicondition1}, \eqref{eq:xicondition2}, \eqref{eq:xicondition3} and $\beta$ large enough so that $\frac{1}{\beta}$ is smaller than each constant on the right-hand side of \eqref{eq:subcase1-1L4conclusion0}, \eqref{eq:subcase1-1L4conclusion1}, \eqref{eq:subcase1-1L4conclusion2}, \eqref{eq:subcase1-2L4conclusion0}, \eqref{eq:subcase1-2L4conclusion1}, 
\eqref{eq:subcase1-2L4conclusion2}, 
\eqref{eq:subcase2-1L4conclusion1},
\eqref{eq:subcase2-1L4conclusion2},
\eqref{eq:subcase2-2L4conclusion1},
\eqref{eq:subcase2-2L4conclusion2}. Then claim \eqref{eq:L4AlternativeClaim} holds. Hence, the L4 conditions \eqref{item01L4condition}, \eqref{item02L4condition} hold.
\end{proof}
\subsection[Sinai-Chernov Ansatz (unbounded expansion near singularity curve)]{Sinai-Chernov Ansatz (unbounded expansion near singularity curve)}\hfill

Analogously to \cite[Definition 3.7]{MR3082535} of the u-\emph{essential} points, we make the following definition of the adapted u-\emph{essential} points.

\begin{definition}[adapted u(s)-\emph{essential}]\label{def:adaptedessential} For any sufficient point $x\in\hat{\mathcal{U}}_+$ with quadruple $(l,N,U,\mathcal{C})$ (\cref{def:sufficientPoint}) given by \cref{theorem:sufficentptsonMhatplus} and with $\Omega\subset\bigcup_{k\in\mathbb{Z}}\Fc^k(U)$ given in \cref{proposition:stablemanifoldthm}, a point $w\in M\smallsetminus\partial M$ is said to be adapted u-\emph{essential} w.r.t. sufficient point $x$ and $\Omega$, if for every $\alpha>0$, there exits:
\begin{enumerate}
    \item a neighborhood $V$ of $w$ and an interger $n>0$ such that $V\cap \mathcal{S}_n=\emptyset$;

    \item for all $y\in\Omega$, every $z\in \Int(W^u_y)\cap V$ and $0\ne dz\in T_z W^u_y$, the followings hold: $\Fc^n(z)\in\hat{\mathcal{U}}_+$ , $D\Fc^n_z(dz)\in\mathcal{C}(\Fc^n(z))$ and $\frac{\sqrt{Q_{\Fc^n(z)}(D\Fc^n_z(dz))}}{\|dz\|}>\alpha$.
\end{enumerate}

Analogously, we can also define the adapted u-\emph{essential} with respect to a sufficient point $x\in\hat{\mathcal{U}}_{-}$ by replacing $\hat{\mathcal{U}}_{+}$ given in \cref{theorem:sufficentptsonMhatplus} with $\hat{\mathcal{U}}_{-}$ given in \cref{theorem:sufficentptsonMhatminus}.

Analogously, we can also define the adapted s-\emph{essential} with respect to sufficient point $x$ by replacing the $\Fc$, $\mathcal{S}_n$, $\mathcal{C}$, $W^u_y$ and $Q$ norm in the above definition with the $\Fc^{-1}$, $\mathcal{S}_{-n}$, $\mathcal{C}'$, $W^s_y$ and $Q'$ norm.
\end{definition}

\begin{lemma}[L3' condition in \cref{thm:lemonLET}: Sinai-Chernov Ansatz in Lemon billiard: the unbounded expansion near singularity curve]\label{lemma:SCAnsatz}
        For the lemon billiard satisfying \cref{def:UniformHyperbolicLemonBilliard} and for a sufficient point $x\in\hat{\mathcal{U}}_+$ with the $(l,N,U,\mathcal{C})$ (\cref{def:sufficientPoint}) quadruple given in \cref{theorem:sufficentptsonMhatplus} and $\Omega\subset\bigcup_{k\in\mathbb{Z}}\Fc^k(U)$ given in \cref{proposition:stablemanifoldthm}, in ($\mathcal{L}_{-}$ measure \cref{notation:symplecticformmeasurenbhd}) Lebesgue measure almost every point on $\mathcal{S}^{-}_1$ is adapted u-\emph{essential} w.r.t. $x$ and $\Omega$.

        Similarly for a sufficient point $x\in\hat{\mathcal{U}}_-$ with the $(l,N,U,\mathcal{C})$ quadruple given in \cref{theorem:sufficentptsonMhatminus}, a.e. point of $\mathcal{S}^{-}_1$ is adapted u-\emph{essential} with respect to $x$ and $\Omega\subset\bigcup_{k\in\mathbb{Z}}\Fc^k(U)$ given in \cref{proposition:stablemanifoldthm}.

        Similarly in Lebesgue measure almost every point on $\mathcal{S}^{+}_1$ is adapted s-\emph{essential} w.r.t. suffcient point $x$ and $\Omega$ given $(l,N,U,\mathcal{C})$ in either \cref{theorem:sufficentptsonMhatplus} or \cref{theorem:sufficentptsonMhatminus}.
\end{lemma}
\begin{proof}
The proof of this lemma relies on the singularity curves having regularity and transversal intersection properties. The proof also uses the L4 condition proof reasoning \cref{theorem:L4condition} in \eqref{eq:subcase1-1L4conclusion1}, 
\eqref{eq:subcase1-2L4conclusion0},
\eqref{eq:subcase1-2L4conclusion1},
\eqref{eq:subcase2-1L4conclusion1} and \eqref{eq:subcase2-2L4conclusion1}. 

It suffice prove the conclusion for $x\in\hat{\mathcal{U}}_+$ with the $(l,N,U,\mathcal{C})$ (\cref{def:sufficientPoint}) quadruple given in \cref{theorem:sufficentptsonMhatplus} and $\Omega\subset\bigcup_{k\in\mathbb{Z}}\Fc^k(U)$ given in \cref{proposition:stablemanifoldthm}. The proof for $x\in\hat{\mathcal{U}}_-$ is the same by symmetry.

By the singularity curve regularity: \cref{corollary:RegSC}, for each $n>0$, $\Scal_n$ is a union of finitely many closed connected curve segments such that each connected curve segment only can have its endpoints be the intersection with another connected curve segment. By \cref{corollary:finiteintersectionwithS1}, for each $n>0$, $\mathcal{S}^{-}_{1}\cap\mathcal{S}_{n}$ is a finite set. Therefore, $\mathcal{S}^{-}_{1}\cap\mathcal{S}_{\infty}$ is a countable set and $\mathcal{S}^{-}_{1}\smallsetminus\mathcal{S}_{\infty}$ has full Lebesgue measure.

Hence, for each $w\in\mathcal{S}^{-}_{1}\smallsetminus\mathcal{S}_{\infty}$, in the billiard table its trajectory $w,\Fc(w),\cdots,\Fc^n(w),\cdots$ will experience arbitrarily many collisions on the interior of $\Gamma_r$ and $\Gamma_R$. 

Note that by \cref{remark:GaranteeLargeTheta}, in our configured lemon billiard, the sequence $w,\Fc(w),\cdots,\Fc^n(w),\cdots$ contains infinitely many elements in $\hat{\mathcal{U}}_{+}\subset \hat{M}_1$(\cref{def:UIntRegionofD}). Suppose that for each $n>0$, the finite trajectory segment $w,\Fc(w),\cdots,\Fc^n(w)$ has $k_n$ elements in $\hat{\mathcal{U}}_{+}$. As $n\rightarrow\infty$, $k_n\rightarrow\infty$.

For any given $\alpha>0$, we can pick some $n$ large enough so that $\hat{M}_1\ni w,\Fc(w),\cdots,\Fc^n(w)\in\hat{\mathcal{U}}_{+}\subset \hat{M}_1$ contains $k_n$ elements in $\hat{\mathcal{U}}_{+}$ and \begin{equation}\label{eq:EssentialConditionForPickingLargen}
C_6\hat{c}_1\Lambda^{k_n-1}_1\cdot c_2\cdot\min\big\{\frac{1}{c_2},C_1,C_4,0.01,C_5\big\}>\alpha,    
\end{equation}
where constant $c_2>0$ is from \cref{Col:EuclideanExpansionMhat1}, constant $C_6$ is defined in \eqref{eq:C6constantForCone}, constants $\hat{c}_1$, $\Lambda_1$ are from \cref{corollary:sinpmetricexpansionMhat1} and constants $C_1$, $C_4$, $0.01$, $C_5$ are from the L4 condition \cref{theorem:L4condition} proof arguments/conclusions in \eqref{eq:subcase1-1L4conclusion1}, \eqref{eq:subcase2-1L4conclusion1}, \eqref{eq:subcase2-2L4conclusion1}.

With the chosen $n$, $k_n$, we further pick a small enough open $V\ni w$ such that:
\begin{itemize}
    \item $\mathcal{S}^{-}_{1}(\xi)\supset V$ with $\xi$ given in \cref{theorem:L4condition};
    \item $V\cap\mathcal{S}_n=\emptyset$;
    \item For all $z\in V$, the trajectory segment $z, \Fc(z),\cdots,\Fc^n(z)$  at each step close to $w, \Fc(w),\cdots,\Fc^n(w)$ so that it also has exact $k_n$ times of visiting $\hat{\mathcal{U}}_{+}$.
\end{itemize}

Note that $W^u_y\cap V$ can only be located in one side of the singularity curve $\mathcal{S}^{-}_1$. We do the same analysis as in the proof \cref{theorem:L4condition} for 4 subcases of the locations of $W^u_y$.

Subcase (1.1): let $\mathbb{S}=\partial\Mrin\smallsetminus\partial M_r$, $z\in (\Int(W^u_y)\cap V)\subset(M_r\smallsetminus\Mrin)\cap\mathbb{S}(\xi)$;

Subcase (1.2): let $\mathbb{S}=\partial\Mrin\smallsetminus\partial M_r$, $z\in (\Int(W^u_y)\cap V)\subset\Mrin\cap\mathbb{S}(\xi)$;

Subcase (2.1): let $\mathbb{S}=\partial\MRin\smallsetminus\partial M_R$, $z\in (\Int(W^u_y)\cap V)\subset(M_R\smallsetminus\MRin)\cap\mathbb{S}(\xi)$;

Subcase (2.2): let $\mathbb{S}=\partial\MRin\smallsetminus\partial M_R$, $z\in (\Int(W^u_y)\cap V)\subset\MRin\cap\mathbb{S}(\xi)$;

Denote by 
\begin{equation}\label{eq:UEssentialExpansionNotation}
\begin{aligned}
    \iota(z)&\dfn\inf\big\{i>0\bigm|\Fc^i(z)\in\hat{\mathcal{U}}_+\big\}\\
    (d\hat{\phi},d\hat{\theta})&\nfd d\hat{z}=D\Fc^{\iota(z)}_z(dz) \text{ for }dz\in T_zW^u_y.\\ 
    d\hat{z}_n&\dfn\Fc^n(dz)\in T_{\Fc^n(z)}M.
\end{aligned}
\end{equation}
And since $\Fc^{\iota(z)}(z)\nfd\hat{z},\cdots,\Fc^n(z)\dfn\hat{z}_n=(\hat{\phi}_n,\hat{\theta}_n)$ contains $k_n$ elements in $\hat{\mathcal{U}}_{+}\subset\hat{M}_1$, thus $d\hat{z}_n=D\hat{F}^t_{1}(d\hat{z})$ for some $t\ge k_n-1$.

In subcase (1.1), the same analysis as for the conclusion in \eqref{eq:subcase1-1L4conclusion1} gives $\frac{\|d\hat{z}\|}{\|dz\|}>c_2\cdot\min\big\{C_1,1\big\}$.

In subcase (1.2), if $z,\Fc(z),\cdots,\Fc^{\iota(z)}(z)$ are all collisions on $\Gamma_r$, then the same analysis as for the conclusion in \eqref{eq:subcase1-2L4conclusion0} gives $\frac{\|d\hat{z}\|}{\|dz\|}>1$. If $z,\Fc(z),\cdots,\Fc^{\iota(z)}(z)$ contain some collision on $\Gamma_R$, then the same analysis as for the conclusion in \eqref{eq:subcase1-2L4conclusion1}, $\frac{\|d\hat{z}\|}{\|dz\|}>c_2$.

In subcase (2.1), the same analysis as for the conclusion in \eqref{eq:subcase2-1L4conclusion1} gives$\frac{\|d\hat{z}\|}{\|dz\|}>c_2\cdot\min\big\{C_4,0.01\big\}$.

In subcase (2.2), the same analysis as for the conclusion in \eqref{eq:subcase2-2L4conclusion1} gives $\frac{\|d\hat{z}\|}{\|dz\|}>c_2\cdot C_5$.
In all $4$ subcases, we get \begin{equation}\label{eq:UEssentialExpansion1}
    \frac{\|d\hat{z}\|}{\|dz\|}>c_2\cdot\min\big\{\frac{1}{c_2},C_1,C_4,0.01,C_5\big\}.
\end{equation}

Then denote by $\chi(y)\dfn\inf\big\{i>0\bigm|\Fc^{-i}(y)\in\hat{M}_{1}\big\}$, $\Fc^{-\chi}(y)\in W^u_{\Fc^{-\chi(y)}(y)}\subset\Int(\hat{M}_1)$ otherwise some $\Fc$-preimage of $W^u_{y}$ would separate ,i.e. had disconnected components on $M_r,M_R$ which contradicts to the $\Fc^{-1}$ invariance of unstable manifold (\cref{proposition:stablemanifoldthm}\eqref{item04stablemanifoldthm}) and by the definition of unstable manifold it should have one connected component.

For all $z\in\Int(W^u_y)$, denote by $z_0= \Fc^{-\chi(y)}(z)$,
$d\hat{z}_0\dfn D\Fc^{-\chi(y)}_z(dz)\in T_{\Fc^{-\chi(y)}(z)}W^u_{\Fc^{-\chi(y)}(y)}$. Therefore, $d\hat{z}_0\in \hat{C}^u_{1,\hat{z}_0}$ by \cref{lemma:UnstableManifoldInCone}\eqref{item01UnstableManifoldInCone} and $d\hat{z}=D\hat{F}_{1}(d\hat{z}_0)\in \hat{C}^u_{1,\hat{z}}$. $(d\hat{\phi}_n,d\hat{\theta}_n)\nfd d\hat{z}_n=D\hat{F}_1^{t}(d\hat{z})\in \hat{C}^u_{1,\hat{z}_n}\subset\mathcal{C}(\hat{z}_n)$ for some $t\ge k_n-1$. 

Moreover, since $\hat{z}\in\hat{\mathcal{U}}_{+}$, if $\hat{z}\in\Mrin\cap\Mrout$, then $d\hat{\theta}/d\hat{\phi}\in\big[\lambda_0(r,R,\phistar),\lambda_1(r,R,\phistar)\big]$. If $\hat{z}\in\Fc(\Mrin\smallsetminus\Mrout)$, then $d\hat{\theta}/d\hat{\phi}\in\big[\frac{1}{2+1/\lambda_0(r,R,\phistar)},\frac{1}{2+1/\lambda_1(r,R,\phistar)}\big]$.

Similarly, since $\hat{z}_n\in\hat{\mathcal{U}}_{+}$, if $\hat{z}_n\in\Mrin\cap\Mrout$, then $d\hat{\theta}_n/d\hat{\phi}_n\in\big[\lambda_0(r,R,\phistar),\lambda_1(r,R,\phistar)\big]$. If $\hat{z}\in\Fc(\Mrin\smallsetminus\Mrout)$, then $d\hat{\theta}_n/d\hat{\phi}_n\in\big[\frac{1}{2+1/\lambda_0(r,R,\phistar)},\frac{1}{2+1/\lambda_1(r,R,\phistar)}\big]$.

The $\lambda_0(r,R,\phistar)$, $\lambda_1(r,R,\phistar)$ are constants from \cref{corollary:cone-lowerbound} determined by the billiard configuration. Hence, $d\hat{\theta}_n/d\hat{\phi}_n$, $d\hat{\theta}/d\hat{\phi}$ are bounded away from $0$ and $\infty$. We denote by $\lambda_U=\max{\big\{\lambda_1(r,R,\phistar),\frac{1}{2+1/\lambda_1(r,R,\phistar)}\big\}}$, $\lambda_L=\min{\big\{\lambda_0(r,R,\phistar),\frac{1}{2+1/\lambda_0(r,R,\phistar)}\big\}}$ and a constant: 
\begin{equation}\label{eq:C6constantForCone}
    C_6=\frac{1}{\sqrt{1+\lambda^2_U}}\sqrt{\lambda_L}.
\end{equation} Then we can estimate

\begin{align*}
    \frac{\sqrt{Q_{\Fc^n(z)}(D\Fc^n_z(dz))}}{\|d\hat{z}\|}=&\frac{\|d\hat{z}_n\|_Q}{\|d\hat{z}\|}=\sqrt{\frac{\sin\hat{\theta}_n d\phi_n d\theta_n}{d\hat{\phi}^2+d\hat{\theta}^2}}=\sqrt{\frac{d\hat{\phi} d\hat{\theta}}{d\hat{\phi}^2+d\hat{\theta}^2}}\cdot\sqrt{\frac{\sin\hat{\theta}_n d\phi_n d\theta_n}{d\hat{\phi} d\hat{\theta}}}\\
    =&\sqrt{\frac{d\hat{\phi} d\hat{\theta}}{d\hat{\phi}^2+d\hat{\theta}^2}}\cdot\sqrt{\frac{\sin\hat{\theta}_n(d\phi_n)^2}{(d\hat{\phi})^2}}\cdot\sqrt{\frac{d\hat{\theta}_n}{d\hat{\phi}_n}\cdot\frac{d\hat{\phi}}{d\hat{\theta}}}
\end{align*}

Since $d\hat{\theta}_n/d\hat{\phi}_n$, $d\hat{\theta}/d\hat{\phi}$ are in $\big[\lambda_L,\lambda_U\big]$, $\sqrt{\frac{d\hat{\phi} d\hat{\theta}}{d\hat{\phi}^2+d\hat{\theta}^2}}\sqrt{\frac{d\hat{\theta}_n}{d\hat{\phi}_n}\cdot\frac{d\hat{\phi}}{d\hat{\theta}}}=\frac{1}{\sqrt{1+(d\hat{\theta}/d\hat{\phi})^2}}\sqrt{\frac{d\hat{\theta}_n}{d\hat{\phi}_n}}\ge C_6$. Hence, \begin{equation}\label{eq:UEssentialExpansion2}\begin{aligned}
    \frac{\|d\hat{z}_n\|_Q}{\|d\hat{z}\|}\ge&C_6\sqrt{\frac{\sin\hat{\theta}_n (d\hat{\theta}_n)^2}{(d\hat{\phi})^2}}=C_6\frac{\sin\hat{\theta}}{\sqrt{\sin{\hat{\theta}_n}}}\sqrt{\frac{\sin^2\hat{\theta}_n (d\hat{\theta}_n)^2}{\sin^2\hat{\theta}(d\hat{\phi})^2}}\ge C_6\frac{\sin\hat{\theta}\|d\hat{z}_n\|_{\p}}{\|d\hat{z}\|_{\p}}\qquad\overbracket{=}^{\mathclap{D\hat{F}^t_1(d\hat{z})=d\hat{z}_n}}\qquad C_6\frac{\sin\hat{\theta}\|D\hat{F}^t_1(d\hat{z})\|_{\p}}{\|d\hat{z}\|_{\p}}\\
    \overbracket{>}^{\mathllap{\text{\cref{corollary:sinpmetricexpansionMhat1}}}}&C_6\hat{c}_1\Lambda^{t}_1\overbracket{\ge}^{\mathrlap{t\ge k_n-1}} C_6\hat{c}_1\Lambda^{k_n-1}_1
\end{aligned}
\end{equation}
Therefore finally, \begin{equation}\label{eq:UEssentialExpansion3}
    \frac{\sqrt{Q_{\Fc^n(z)}(D\Fc^n_z(dz))}}{\|dz\|}=\frac{\sqrt{Q_{\Fc^n(z)}(D\Fc^n_z(dz))}}{\|d\hat{z}\|}\frac{\|d\hat{z}\|}{\|dz\|}\qquad\overbracket{>}^{\mathclap{\eqref{eq:UEssentialExpansion1},\eqref{eq:UEssentialExpansion2}}}C_6\hat{c}_1\Lambda^{k_n-1}_1\cdot c_2\cdot\min\big\{\frac{1}{c_2},C_1,C_4,0.01,C_5\big\}\overbracket{>}^{\mathclap{\text{condition \eqref{eq:EssentialConditionForPickingLargen}}}}\alpha
\end{equation}

Now we consider for $x\in\hat{\mathcal{U}}_-$ with the $(l,N,U,\mathcal{C})$ quadruple given in \cref{theorem:sufficentptsonMhatminus} and $\Omega$ given by \cref{proposition:stablemanifoldthm}, for all points $w\in\mathcal{S}^{-}_{1}\smallsetminus\mathcal{S}_{\infty}$, we pick some $n$ large enough so that $\hat{M}_1\ni w,\Fc(w),\cdots,\Fc^n(w)\in\hat{\mathcal{U}}_{-}\subset \hat{M}_1$ contains $k_n$ elements in $\hat{\mathcal{U}}_{+}$ and satisfies the condition \eqref{eq:EssentialConditionForPickingLargen}. 

If $\Fc^n(w)\in\hat{\mathcal{U}}_{-}\cap\hat{\mathcal{U}}_{+}$, then the previous analysis gives the desired conclusion. 

Otherwise, $\Fc^n(w)\in\hat{\mathcal{U}}_{-}\smallsetminus\hat{\mathcal{U}}_{+}$, we assume that $j<n$ is the largest integer such that $\Fc^{j}(w)\in\hat{\mathcal{U}}_{+}$ and $\Fc^{j+1}(w)$, $\cdots$, $\Fc^n(w)$ contain no element in $\hat{\mathcal{U}}_{+}$. Then $w, \Fc(w),\cdots,\Fc^j(w)$ has exact $k_n$ times of visiting $\hat{\mathcal{U}}_{+}$. Then with the chosen $n$, $j$, $k_n$, we further pick a small enough open $V\ni w$ such that:
\begin{itemize}
    \item $\mathcal{S}^{-}_{1}(\xi)\supset V$ with $\xi$ given in \cref{theorem:L4condition};
    \item $V\cap\mathcal{S}_n=\emptyset$;
    \item For all $z\in V$, the trajectory segment $z, \Fc(z),\cdots,\Fc^j(z)$  at each step close to $w, \Fc(w),\cdots,\Fc^j(w)$ so that it also has exact $k_n$ times of visiting $\hat{\mathcal{U}}_{+}$.
\end{itemize}
Using the notation from \eqref{eq:UEssentialExpansionNotation} and by the previous same analysis in \eqref{eq:UEssentialExpansion3} we get $\frac{\sqrt{Q_{\Fc^j(z)}(D\Fc^j_z(dz))}}{\|dz\|}>\alpha$. On the other hand, by \cref{caseA,caseB,caseC} the positive quadrant cone field $\mathcal{C}$ on $\hat{M}$ and positive cone field $\mathcal{C}$ on $\hat{M}_1$ are jointly invariant. Then by \cref{def:JointInvariant,remark:MonotoneCone} and \eqref{eq:ExpansionInCone}, $\sqrt{\frac{Q_{\Fc^n(z)}(D\Fc^n_z(dz))}{Q_{\Fc^j(z)}(D\Fc^j_z(dz))}}\ge1$. Hence, $\frac{\sqrt{Q_{\Fc^n(z)}(D\Fc^n_z(dz))}}{\|dz\|}>\alpha$.\qedhere
\end{proof}

\begin{corollary} [analogous to {\cite[lemma 5.23]{MR3082535}}]\label{corollary:essentialExpansionLemma} For a sufficient point $x\in\hat{\mathcal{U}}_+$, suppose $(l,N,U,\mathcal{C})$ (\cref{def:sufficientPoint}) is the quadruple given in \cref{theorem:sufficentptsonMhatplus} with the $\Omega$ given in \cref{proposition:stablemanifoldthm}. For every $t>0$ and $0<h<1$, there exists an integer $M_{t,h}>0$, two compact subsets $\mathcal{S}_{t,h}$ and $\mathcal{E}_{t,h}$ of $\mathcal{S}^{-}_{1}$ (defined in \cref{notation:symplecticformmeasurenbhd}), a real number $r_{t,h}>0$ such that \begin{enumerate}
        \item $\mathcal{S}^{-}_{1}=\mathcal{S}_{t,h}\cup\mathcal{E}_{t,h}$ and $\mathcal{L}_{-}(\mathcal{E}_{t,h})<h$;
        \item If $y\in\Omega$, $z\in \Int(W^u_y)\cap\mathcal{S}_{t,h}(r_{t,h})$ and $\Fc^j(z)\in U$ with $j\ge M_{t,h}$ and $z\notin\mathcal{S}_{j}$, then for $0\ne dz\in T_zW^u_y$ \[\frac{\sqrt{Q_{\Fc^jz}(D\Fc^j_z(dz))}}{\|dz\|}>t.
        \]\end{enumerate}
The same conclusion also holds for a sufficient point $x\in\hat{\mathcal{U}}_-$ with the quadruple $(l,N,U,\mathcal{C})$(\cref{def:sufficientPoint}) given in \cref{theorem:sufficentptsonMhatminus} and $\Omega$ given in \cref{proposition:stablemanifoldthm}. And note that the proof relies on \cref{lemma:SCAnsatz}.
\end{corollary}
\begin{proof}
    We first prove the case for sufficient point $x\in\hat{\mathcal{U}}_+$ with the quadruple $(l,N,U,\mathcal{C})$ given in \cref{theorem:sufficentptsonMhatplus}. The proof for a sufficient point $x\in\hat{\mathcal{U}}_-$ is the same. First, note that $\mathcal{S}^{-}_1$ for a lemon billiard are the $4$ line components of closure($\partial\Mrin\smallsetminus\partial M_r$) and closure($\partial\MRin\smallsetminus\partial M_R$). It suffices to prove the lemma with $\mathcal{S}^{-}_1$ replaced by $\Sigma$: arbitrary one of the line components of closure($\partial\Mrin\smallsetminus\partial M_r$) and closure($\partial\MRin\smallsetminus\partial M_R$). The regularity of the Lebesgue measure ($\mathcal{L}_{-}(\cdot)$) ensures that there is a compact subset $\Sigma_1$ of the $\Int(\Sigma)$ such that $\mathcal{L}_{-}(\Sigma\smallsetminus\Sigma_1)<h$ and every point of $\Sigma_1$ is adapted u-essential with respect to the sufficient point $x$ with the quadruple $(l,N,U,\mathcal{C})$ and $\Omega$ given in \cref{proposition:stablemanifoldthm}.
    
    For any $y\in\Omega$, by \cref{def:adaptedessential} of the adapted u-essential point, for every $w\in\Sigma_1$ there is a real number $p_w>0$ and an integer $n_w>0$ such that $B(w,p_w)\cap\mathcal{S}_{n_w}=\emptyset$, $D\Fc^{n_w}_z(dz)\in\mathcal{C}(\Fc^{n_w}(z))$ and $\frac{\sqrt{Q_{\Fc^{n_w}(z)}(D\Fc^{n_w}_z(dz))}}{\|dz\|}>t$ for all $z\in \Int(W^u_y)\cap B(w,p_w)$, where $B(w,p_w)$ is the open disk of phase space $M$ centered at $w$ with radius $p_w$. Since the cone $\mathcal{C}$ on $\hat{M}_1$ is invariant under the return map $\hat{F}_1$, if $z\notin\mathcal{S}_j$ and $\Fc^j(z)\in U$, then by \cref{remark:MonotoneCone} and $\Fc^j(z)\in U\subset\hat{M}_1$ we get $\frac{\sqrt{Q_{\Fc^j z}(D\Fc^j_z(dz))}}{\|dz\|}\ge\frac{\sqrt{Q_{\Fc^{n_w}(z)}(D\Fc^{n_w}_z(dz))}}{\|dz\|}>t$ for all $j\ge n_w$.

    Due to the compactness of $\Sigma_1$, there are finite $s>0$ points $w_1,\cdots,w_s$ of $\Sigma_1$, real numbers $p_{w_1},\cdots,p_{w_s}$ and positive integers $n_{w_1},\cdots,n_{w_s}$ such that $\mathcal{S}_{n_{w_i}}\cap B(w_i,p_{w_i})=\emptyset$, $\Sigma_1\subset\bigcup_{i=1}^s B(w_i,p_{w_i})$ and for all $z\in \Int(W^u_y)\cap B(w_i,p_{w_i})$, $j\ge n_i$ with $z\notin\mathcal{S}_j$, the following holds true: \[\frac{\sqrt{Q_{\Fc^j z}(D\Fc^j_z(dz))}}{\|dz\|}>t.\]  
    
    Using \cref{notation:symplecticformmeasurenbhd} for set's neighborhood, since the compact set $\Sigma_1$ is contained in the open set $\bigcup_{i=1}^s B(w_i,p_{w_i})$, there exist a compact set $K$ and two real numbers $q>0$ and $r_{t,h}>0$  such that \[\Sigma_1\subset\Sigma_1(q)\subset K\subset K(r_{t,h})\subset\bigcup_{i=1}^s B(w_i,p_{w_i}).\] Therefore, we set $M_{t,h}\dfn\max\big\{n_{w_1},\cdots,n_{w_s}\big\}$, $\mathcal{S}_{t,h}\dfn K\cap\Sigma$ and $\mathcal{E}_{t,h}\dfn\Sigma\smallsetminus\Sigma_1(q)$. For $j>M_{t,h}$, if $z\in(\Int(W^u_y)\cap \mathcal{S}_{t,h}(r_{t,h})\subset (\Int(W^u_y)\cap\bigcup_{i=1}^s B(w_i,p_{w_i}))$ with $z\notin\mathcal{S}_j$, then $\frac{\sqrt{Q_{\Fc^j z}(D\Fc^j_z(dz))}}{\|dz\|}>t$. And $\mathcal{L}_{-}(\mathcal{E}_{t,h})=\mathcal{L}_{-}(\Sigma\smallsetminus\Sigma_1(q))\le\mathcal{L}_{-}(\Sigma\smallsetminus\Sigma_1)<h$.
\end{proof}

\subsection[Adapted local ergodic theorem from {\cite[theorem 4.1]{MR3082535}} for lemon billiard] {Adapted local ergodic theorem from \cite[theorem 4.1]{MR3082535} for lemon billiard}\hfill
\begin{theorem}[Adapted from {\cite[theorem 4.1]{MR3082535}}]\label{thm:lemonLET}
    In the lemon billiard phase space, let $x\in (\hat{\mathfrak{U}}_+\cup\hat{\mathfrak{U}}_-)$ be a sufficient point (\cref{def:sufficientPoint}) defined in \cref{theorem:sufficentptsonMhatminus} or \cref{theorem:sufficentptsonMhatplus} for some quadruple $(l,N,U,\mathcal{C})$ and let $\Omega$ be a subset of $\bigcup_{k\in\mathbb{Z}}\Fc^kU$ defined in \cref{proposition:stablemanifoldthm}. If the following conditions L1, L2, L3', L4 hold: \begin{itemize}
        \item L1 (the regularity condition is the same as in \cite[theorem 4.1]{MR3082535}). The sets $\mathcal{S}_k$ and $\mathcal{S}_{-k}$ are regular (\cref{def:RegularSingularityCurve}) for every $k>0$,

        \item L2 (the alignment condition is the same as in \cite[theorem 4.1]{MR3082535}). For every $k\ge1$: 
        
        If $y\in \mathcal{S}\cap \Fc^{-N}(U)$ where $\mathcal{S}\subset\mathcal{S}_{-k}\smallsetminus \mathcal{S}_{-(k-1)}$, then the tangential vector of $\mathcal{S}$ at $y$ is contained in $\mathcal{C}$.
        
        If $y\in \mathcal{S}\cap U$ where $\mathcal{S}\subset\mathcal{S}_{k}\smallsetminus \mathcal{S}_{k-1}$, then the tangential vector of $\mathcal{S}$ at $y$ is contained in $\mathcal{C}'$ which is the complement of the cone $\mathcal{C}$,
        \item L3' (Sinai-Chernov ansatz condition adapted for lemon billiard). In the Lebesgue measure $\mathcal{L}_-$ from \cref{notation:symplecticformmeasurenbhd}, almost every point on $\mathcal{S}^{-}_1$ is adapted u-\emph{essential} (\cref{def:adaptedessential}) w.r.t. $x$ and $\Omega$. In the Lebesgue measure $\mathcal{L}_+$ from \cref{notation:symplecticformmeasurenbhd}, almost every point on $\mathcal{S}^{+}_1$ is adapted s-\emph{essential} (\cref{def:adaptedessential}) w.r.t. $x$ and $\Omega$,

        \item L4 (the contraction condition is the same as in \cite[theorem 4.1]{MR3082535}). There exist constants $\beta>0$ and $\xi>0$ such that:

        if $y\in\Omega\cap U$, $z\in W^u_y$ and $\Fc^{-k}z\in\mathcal{S}^{-}_{1}(\xi)$ with $k\ge1$, then $\|D\Fc^{-k}_z\mid_{T_{z}W^{u}_{y}}\|\le\beta$;

        if $y\in\Omega\cap U$, $z\in W^s_y$ and $\Fc^{k}z\in\mathcal{S}^{+}_{1}(\xi)$ with $k\ge1$, then $\|D\Fc^{k}_z\mid_{T_{z}W^{s}_{y}}\|\le\beta$, where $\mathcal{S}_1(\xi)$, $\mathcal{S}_{-1}(\xi)$ are the $\xi$ neighborhood of $\mathcal{S}^{+}_1$, $\mathcal{S}^{-}_{1}$ (see \cref{notation:symplecticformmeasurenbhd}), $W_y^u$ and $W_y^s$ are the local unstable/stable manifold at $y$. The existence of unstable/stable manifolds is guaranteed by \cite[Proposition 3.4]{MR3082535}, 
    \end{itemize} then there is an open neighborhood $\mathcal{O}$ of $x$ contained (mod 0) in a Bernoulli ergodic component of $\Fc$ with respect to the measure $\mu$ in \cref{notation:symplecticformmeasurenbhd}.
\end{theorem}
\begin{proof}
    We provide a proof adapted from \cite[section 5]{MR3082535} for the local ergodic theorem.

    \cref{lemma:SCAnsatz} implies the lemon billiard singularity curves $\Scal^-_1$ and $\Scal^+_1$ satisfy L3' thus  \cref{corollary:essentialExpansionLemma}.
        
        We see that in \cite[proof of Proposition 5.22]{MR3082535} the conditions L1, L2, L4 together with \cref{corollary:essentialExpansionLemma} replacing \cite[Lemma 5.23]{MR3082535} suffice to produce \cite[Proposition 5.22]{MR3082535} thus \cite[Proposition 5.20]{MR3082535}. 
        
        And L1, L2 conditions suffice to produce \cite[Proposition 5.19]{MR3082535}. Then \cite[Propositions 5.19 and 5.20]{MR3082535} imply the \emph{Sinai's Theorem \cite[Theorem 5.16]{MR3082535}} which means that there is a neighborhood $\mathcal{O}$ of the sufficient point $x$ (given in \cref{theorem:sufficentptsonMhatplus} or \cref{theorem:sufficentptsonMhatminus}) to be contained in one ergodic component of $\Fc$.
\end{proof}
\begin{remark}
    We note that the proof of \cref{lemma:SCAnsatz} relies on the regularity and transverse intersection properties for singularity curves (essential for L1, L2 properties). In next section, we proceed to prove that our lemon billiard satisfies the L1 and L2 conditions which are properties of singularity curves of lemon billiard.
\end{remark}

\section{Local ergodicity conditions, singularity curves, and global ergodicity}\label{sec:SCLB}
\subsection[Arch singularity curves of lemon billiard]{Arch singularity curves of lemon billiard}\label{sec:ArchSCLB}
\begin{figure}[h!]
\begin{minipage}{.5\textwidth}
\begin{center}
\begin{tikzpicture}[xscale=0.75, yscale=0.75]
    \tkzDefPoint(0,0){Or};
    \tkzDefPoint(0,-4.8){Y};
    \pgfmathsetmacro{\rradius}{4.5};
    \pgfmathsetmacro{\bdist}{6.3};
    \tkzDefPoint(0,\bdist){OR};
    \pgfmathsetmacro{\phistardeg}{40};
    \pgfmathsetmacro{\XValueArc}{\rradius*sin(\phistardeg)};
    \pgfmathsetmacro{\YValueArc}{\rradius*cos(\phistardeg)};
    \tkzDefPoint(\XValueArc,-1.0*\YValueArc){B};    \tkzDefPoint(-1.0*\XValueArc,-1.0*\YValueArc){A};
    \pgfmathsetmacro{\Rradius}{veclen(\XValueArc,\YValueArc+\bdist)};
    \tkzDrawArc[name path=Cr,thick](Or,B)(A);
    \tkzDrawArc[name path=CR,thick](OR,A)(B);
    \pgfmathsetmacro{\PXValue}{\Rradius*sin(7)};
    \pgfmathsetmacro{\PYValue}{\bdist-\Rradius*cos(7)};
    \tkzDefPoint(\PXValue,\PYValue){P};
    \tkzDefPoint(0,\bdist-\Rradius){C};
    \draw[ultra thin,dashed](OR)--(P);
    \draw[ultra thin,dashed](OR)--(A);
    \draw[ultra thin,dashed](OR)--(C);
    \tkzDefPointBy[projection=onto OR--P](A)\tkzGetPoint{M};
    \coordinate (Q_far) at ($(A)!2!(M)$);
    \path[name path=linePQ] (P)--(Q_far);
    \path[name intersections={of=Cr and linePQ,by={Q}}];
    \tkzDrawPoints(A,B,P,Or,OR,Q);
    \begin{scope}[ultra thin,decoration={
    markings,
    mark=at position 0.5 with {\arrow{>}}}
    ] 
    \draw[blue,postaction={decorate}] (A)--(P);
    \draw[blue,postaction={decorate}] (P)--(Q);
    \end{scope}
    \tkzLabelPoint[below](P){$P$};
    \tkzLabelPoint[below](A){$A$};
    \tkzLabelPoint[right](Q){$Q=p(x)$};
    \tkzLabelPoint[below](B){$B$};
    \tkzLabelPoint[below](Or){$O_r$};
    \tkzLabelPoint[above](OR){$O_R$};
    \tkzMarkAngle[arc=lll,ultra thin,size=0.7](A,OR,P);
    \node[] at ($(OR)+(3.2,0.1)$) {\small $\measuredangle{AO_RP}=2\theta_1=\Phi_*+\Phi$};
    \node[] at ($(OR)+(3.2,-0.6)$)  {\small $\measuredangle O_rO_RP=\Phi$};
    \node[] at ($(OR)+(3.2,-1.3)$)  {\small $\angle AO_RO_r=\Phi_*$};
\end{tikzpicture}
\caption{$x\in\textcolor{blue}{\mathcal{AS}^{in}_{-1}}$ with $2$-step past trajectory starting from corner $A$.}\label{fig:ArchSingularitiesneq1}
\end{center}
\end{minipage}\hfill
\begin{minipage}{.5\textwidth}
\begin{center}
\begin{tikzpicture}[xscale=0.75,yscale=0.75]
    \tkzDefPoint(0,0){Or};
    \tkzDefPoint(0,-4.8){Y};
    \pgfmathsetmacro{\rradius}{4.5};
    \pgfmathsetmacro{\phistardeg}{40};
    \pgfmathsetmacro{\XValueArc}{\rradius*sin(\phistardeg)};
    \pgfmathsetmacro{\YValueArc}{\rradius*cos(\phistardeg)};
    \pgfmathsetmacro{\Rradius}{
        \rradius*sin(\phistardeg)/sin(18)
    };
    \pgfmathsetmacro{\bdist}{(-\rradius*cos(\phistardeg))+cos(18)*\Rradius}
    \tkzDefPoint(0,\bdist){OR};
    \tkzDefPoint(\XValueArc,-1.0*\YValueArc){B};    \tkzDefPoint(-1.0*\XValueArc,-1.0*\YValueArc){A};
    \tkzDrawArc[name path=Cr, thick](Or,B)(A);
    \tkzDrawArc[name path=CR, thick](OR,A)(B);
    \pgfmathsetmacro{\PXValue}{-\Rradius*sin(3)};
    \pgfmathsetmacro{\PYValue}{\bdist-\Rradius*cos(3)};
    \tkzDefPoint(\PXValue,\PYValue){P1};
    
    \tkzDefPointBy[projection=onto OR--P1](A)\tkzGetPoint{M};
    \coordinate (P2) at ($(A)!2!(M)$);
    
    \tkzDefPointBy[projection=onto OR--P2](P1)\tkzGetPoint{M2};
    \coordinate (Q_far) at ($(P1)!2!(M2)$);
    \path[name path=Line_P2_Q_far](P2)-- (Q_far);
    \path[name intersections={of=Cr and Line_P2_Q_far,by={Q}}];
    \tkzDefPoint(0,\bdist-\Rradius){C};
    \draw[ultra thin,dashed](OR)--(P2);
    \draw[ultra thin,dashed](OR)--(A);
    \draw[ultra thin,dashed](OR)--(C);
    \begin{scope}[ultra thin,decoration={markings,mark=at position 0.5 with {\arrow{>}}}] 
    \draw[red,postaction={decorate}] (A)--(P1);
    \draw[red,postaction={decorate}] (P1)--(P2);
    \draw[red,postaction={decorate}] (P2)--(Q);
    \end{scope}
    \tkzLabelPoint[below](A){$A$};
    \tkzLabelPoint[below](P2){$P$};
    \tkzLabelPoint[above](Q){$Q=p(x)$};
    \tkzLabelPoint[below](B){$B$};
    \tkzLabelPoint[below](Or){$O_r$};
    \tkzLabelPoint[above](OR){$O_R$};
    \tkzDrawPoints(A,B,Or,P1,P2,Q,OR);
    \tkzMarkAngle[arc=lll,ultra thin,size=0.7](A,OR,P2);
    \node[] at ($(OR)+(3.0,0.8)$)  {$\measuredangle{AO_RP=4\theta_1=\Phi_*+\Phi}$};
    \node[] at ($(OR)+(4.0,0)$)  {$\measuredangle{O_rO_RP=\Phi}$};
    \node[] at ($(OR)+(4.0,-0.8)$)  {$\measuredangle{AO_RO_r=\Phi_*}$};
\end{tikzpicture}
\caption{$x\in\textcolor{red}{\mathcal{AS}^{in}_{-2}}$ with the $3$-step past trajectory starting from corner $A$}\label{fig:ArchSingularitiesneq2}
\end{center}
\end{minipage}\hfill
\end{figure}
\begin{figure}[h!]
\begin{minipage}{.38\textwidth}
\begin{center}
    \begin{tikzpicture}[xscale=0.45, yscale=0.45]
        \pgfmathsetmacro{\PHISTAR}{6};
        \pgfmathsetmacro{\WIDTH}{8};
        \tkzDefPoint(-\PHISTAR,0){LL};
        \tkzDefPoint(-\PHISTAR,\WIDTH){LU};
        \tkzDefPoint(\PHISTAR,0){RL};
        \tkzDefPoint(\PHISTAR,\WIDTH){RU};
        \tkzDefPoint(-\PHISTAR,\PHISTAR){X};
        \tkzDefPoint(\PHISTAR,\WIDTH-\PHISTAR){Y};
        \tkzDefPoint(-\PHISTAR,\WIDTH-\PHISTAR){IX};
        \tkzDefPoint(\PHISTAR,\PHISTAR){IY};
        \tkzDefPoint(\PHISTAR,0.5*\PHISTAR){S1};
        \tkzDefPoint(\PHISTAR,\PHISTAR/3){S2};
        \tkzDefPoint(\PHISTAR,\PHISTAR/4){S3};
        \tkzDefPoint(\PHISTAR,\PHISTAR/5){L4};
        \tkzDefPoint(-\PHISTAR,0.5*\PHISTAR){LS1};
        \tkzDefPoint(-\PHISTAR,\PHISTAR/3){LS2};
        \tkzDefPoint(-\PHISTAR,\PHISTAR/4){LS3};
        \draw [name path=lineLLIY,ultra thin,dashed] (RL) --(LS1);
        \draw [name path=lineLLIY,ultra thin,dashed] (RL) --(LS2);
        \draw [name path=lineLLIY,ultra thin,dashed] (RL) --(LS3);
        
        \draw [thin] (LL) --(LU);
        \draw [thin] (LU) --(RU);
        \draw [thin] (RL) --(RU);
        \draw [thin] (LL) --(RL);
        \draw [name path=lineLLIY, blue,ultra thin,dashed] (LL) --(IY);
        \draw [name path=lineXRL,red,ultra thin,dashed]  (X)--(RL);
        \draw [name path=lineS1LL,ultra thin,dashed]  (S1)--(LL);
        \draw [name path=lineS2LL,ultra thin,dashed]  (S2)--(LL);
        \draw [name path=lineS3LL,ultra thin,dashed]  (S3)--(LL);
        \path[name intersections= 
        {of=lineLLIY and lineXRL,by={M0}}];
        \draw[blue,thick](M0)--(IY);
        \path[name intersections= 
        {of=lineS1LL and lineXRL,by={M1}}];
        \draw[Red,thick](M1)--(S1);
        \path[name intersections= 
        {of=lineS2LL and lineXRL,by={M2}}];
        \draw[Green,thick](M2)--(S2);
        \path[name intersections= 
        {of=lineS3LL and lineXRL,by={M3}}];
        \draw[Purple,thick](M3)--(S3);
        \tkzDrawPoints(IY,S1,S2,S3,M0,M1,M2,M3,RL);
        \tkzLabelPoint[below left](RU){\Large $\MRout$};
        \node[] at (1.08*\PHISTAR,0.6*\PHISTAR){\Large $\xhookrightarrow{\Fc}$};
        \tkzLabelPoint[below](LL){$-\Phistar$};
        \tkzLabelPoint[below](RL){$+\Phistar$};
    \end{tikzpicture}
    \captionsetup{{width=0.7\linewidth}}
    \caption{The pre-image of Arch Singularity curves near $(\Phi_*,0)$: $\Fc^{-1}(\mathcal{AS}^{\text{\upshape in}}_{-n})$, $n=1,2,3,4,\cdots$, are straight line segments singularities $\mathcal{S}_{-n}\cap\MRout$}\label{fig:preimageArchSingularityCurves}
\end{center}
\end{minipage}\hfill
\begin{minipage}{.62\textwidth}
\begin{center}
    \scalebox{4.0}{\input{figs/Figure_1.pgf}}
    \caption{Arch Singularity curves: $\mathcal{AS}^{\text{\upshape in}}_{-n}$ with $n=1,2,3,4$ near $(\phi_*,\phi_*-\Phi_*)$ as a zoom in to $\Nin$ of \cref{fig:MrN}.}\label{fig:ArchSingularities}
\end{center}
\end{minipage}\hfill
\end{figure}
\begin{definition}\label{def:ASinCurves}
    For each $n\ge1$, we define \emph{Arch Singularity Curve} as $\mathcal{AS}^{\text{\upshape in}}_{-n}\dfn\big\{x\in\Mrin\bigm|\Fc^{-(n+1)}(x)\text{ is at corner and }\Fc^{-k}(x)\in\Int(M_R),k=1,2,\cdots,n\big\}$.
\end{definition}
\begin{remark}
    Each $x\in\mathcal{AS}^{\text{\upshape in}}_{-n}$ is the collision on $\Gamma_r$ having an in the past $n+1$ step trajectory starting from a corner followed by $n$ consecutive collisions on $\Gamma_R$ and lastly ending at $x$ (see \cref{fig:ArchSingularitiesneq1,fig:ArchSingularitiesneq2}). Each $\mathcal{AS}^{\text{\upshape in}}_{-n}$ has two components located in the neighborhoods of $(\phi_*,\phi_*-\Phi_*)$ and $(2\pi-\phi_*,\pi-\phi_*+\Phi_*)$ on the boundary of $M_r$. The $\Fc$ preimage of each component of every $\mathcal{AS}^{\text{\upshape in}}_{-n}$ is also located in the $(\Phi_*,0)$ and $(0,\pi)$ neighborhoods contained in $\MRout$. \cref{fig:preimageArchSingularityCurves,fig:ArchSingularities} show one component of $\mathcal{AS}^{\text{\upshape in}}_{-n}$ for $n=1,2,3,4$ near $(\phi_*,\phi_*-\Phi_*)$ and the $\Fc$ preimage in $\MRout$ near $(\Phi_*,0)$.
\end{remark}
\begin{notation}\label{def:lengthfunctiononArchSingularityCurve}
As shown in \cref{fig:ArchSingularitiesneq1,fig:ArchSingularitiesneq2}, suppose that $\mathcal{AS}^{\text{\upshape in}}_{-n}\ni x_2=(\phi_2,\theta_2)$, $(\Phi,\theta_1)\nfd x_1=\Fc^{-1}(x_2)$, $d_1=R\sin\theta_1$, $d_2=r\sin\theta_2$, $P=p(x_1)$, $Q=p(x_2)$ in the lemon billiard table and $\tau_1=|PQ|$.
\end{notation}
\begin{remark}\label{rmk:lengthfunctiononArchSingularityCurve}
    $\tau_1$, $d_2$, $d_1$ are defined similarly as in \cref{def:fpclf}. $\phi_2,\theta_2,\tau_1,d_1,d_2$ are smooth functions of $\Phi$ since the preimage $\mathcal{AS}^{\text{\upshape in}}_{-n}$ is smoothly parametrized in $\Phi$ and by \cite[equation (2.26)]{cb} the differential of the billiard map is smooth.
\end{remark}
\begin{proposition}[$\mathcal{AS}^{\text{\upshape in}}_{-n}$ tangential slopes]\label{proposition:Archsingularitycurvesproperties}
    Suppose that our lemon billiard satisfies \cref{def:UniformHyperbolicLemonBilliard} and has corners $A$ and $B$. Using the functions $\phi_2,\theta_2,d_1,d_2,\tau_1$ from \cref{def:lengthfunctiononArchSingularityCurve}, for each $x_2=(\phi_2,\theta_2)$ on one component of $\mathcal{AS}^{\text{\upshape in}}_{-n}$ (for each $n\ge1$) (in \cref{fig:ArchSingularities}) with the $n$ steps in the past on $\Gamma_R$ and the $n+1$-th step in the past starting from the corner $A$ (see \cref{fig:preimageArchSingularityCurves}), we assume that $x_2$ has preimage $\Fc^{-1}(x_2)=x_1=(\Phi,\theta_1)$. The curve component, the length $\tau_1$ and the parameters $\theta_2$, $\phi_2$ are smooth functions of $\Phi$ and satisfy the following properties.
    
    \begin{enumerate}
        \item\label{item01Archsingularitycurvesproperties} In phase space, this connected component of $\mathcal{AS}^{\text{\upshape in}}_{-n}$ is a smooth curve with $(\phistar,\phistar-\frac{n}{n+1}\Phistar)$ and $(\phistar,\phistar-\frac{n+1}{n}\Phistar)$ as its endpoints on the boundary of $M_r$ (see \cref{fig:ArchSingularities}).
        \item\label{item02Archsingularitycurvesproperties} The $\tau_1$, $\theta_2$, $\phi_2$ as smooth functions of $\Phi$ have $\Phi$ varying in $(\Phistar-\frac{2}{n+1}\Phistar,\Phistar)$, $\frac{d\tau_1}{d\Phi}<0$, i.e., $\tau_1$ is a monotonically decreasing function. And $\frac{d\phi_2}{d\Phi}=\frac{(1+\frac{1}{2n})\tau_1-d_1}{d_2}<0$, $\frac{d\theta_2}{d\Phi}=\frac{d\phi_2}{d\Phi}-(1+\frac{1}{2n})=\frac{(1+\frac{1}{2n})(\tau_1-d_2)-d_1}{d_2}$.
        \item\label{item03Archsingularitycurvesproperties} The tangential slope of the component of the curve at $x_2$: $\mathcal{V}(x_2)\dfn\frac{d\theta_2}{d\phi_2}$ is also a smooth function of $\Phi\in (\Phistar-\frac{2}{n+1}\Phistar,\Phistar)$ and $\mathcal{V}(x_2)$ has the range contained in $\big[-\infty,-0.495\big)\cup\big(1,+\infty\big]$. And especially if $n\ge2$, $\mathcal{V}(x_2)$ has the range contained in $\big[-\infty,-1.49\big)\cup\big(1,+\infty\big]$
    \end{enumerate}
\end{proposition}
\begin{proof}
    By \cref{rmk:lengthfunctiononArchSingularityCurve}, $\tau_1$, $\theta_2$, $\phi_2$, $d_1$, $d_2$ are smooth functions of $\Phi$. Since $x_1=(\Phi,\theta_1)=\Fc^{-1}(x_2)$ has the $n$-th previous step on corner $A$, it is on the straight line going through $(-\Phi_*,0)$ of slope $\frac{1}{2n}$ and also in the region $\MRout$ (see \cref{fig:preimageArchSingularityCurves}). For each line segment with slope $\frac{1}{2n}$ contained in $\Mrout$, in \cref{fig:preimageArchSingularityCurves} we observe that $\Phi$ varies in $(\Phistar-\frac{2}{n+1}\Phistar,\Phistar)$. We now proceed to prove each item.

    \textbf{Proof of \eqref{item01Archsingularitycurvesproperties}: } It suffices to check the limits of $\phi_2$, $\theta_2$ as $\Phi\longrightarrow\Phistar-\frac{2}{n+1}\Phistar$ and $\Phi\longrightarrow\Phistar$. Also note that on the line segment with slope $1/2n$ inside $\Mrout$ of \cref{fig:preimageArchSingularityCurves}. As $\Phi\rightarrow\Phi_*-\frac{2}{n+1}\Phistar$, $\theta_1\rightarrow\frac{1}{n+1}\Phistar$. As $\Phi\rightarrow\Phi_*$, $\theta_1\rightarrow\frac{1}{n}\Phistar$.

    In our billiard table as in \cref{fig:ArchSingularitiesneq1,fig:ArchSingularitiesneq2} and in the coordinate system given in \cref{def:StandardCoordinateTable}, $O_r=(0,0)$, $O_R=(0,b)$,
    \begin{equation}\label{eq:ArchSingularityPQcoordinates}
    \begin{aligned}
        (x_Q,y_Q)=Q&=p(x_2)=(r\sin\phi_2,-r\cos\phi_2),\\
        (x_P,y_P)=P&=p(x_1)=(R\sin\Phi,b-R\cos\Phi).
    \end{aligned}
    \end{equation}
    And note that given the counterclockwise orientation on $\Gamma_r$,
    \begin{equation}\label{eq:theta2phi2theta1PhiEquation}\theta_2=\measuredangle{(\overrightarrow{PQ},\quad\text{tangential direction of }\Gamma_r\text{ at } Q)}=\phi_2-(\Phi+\theta_1).
    \end{equation} 
    
    Therefore, the \textbf{left limit}: if $\Phi\rightarrow\Phistar-\frac{2}{n+1}\Phistar$, then $\theta_1\rightarrow\frac{1}{n+1}\Phistar$ in \cref{fig:preimageArchSingularityCurves} $\MRout\ni x_1$ is approaching $\mathcal{S}_{-1}\smallsetminus\partial M_R$, $x_2=\Fc(x_1)\rightarrow\partial M_r$ and in the billiard table $p(x_2)=Q\rightarrow B$, which means $\phi_2\longrightarrow\phi_*$. Thus, \eqref{eq:theta2phi2theta1PhiEquation} implies $\theta_2=\phi_2-(\Phi+\theta_1)\longrightarrow\phistar-(\Phistar-\frac{2}{n+1}\Phistar+\frac{1}{n+1}\Phistar)=\phistar-\frac{n}{n+1}\Phistar$.

     The \textbf{right limit}: if $\Phi\rightarrow\Phistar$, then $\theta_1\rightarrow\frac{1}{n}\Phistar$, $P=p(x_1)\rightarrow B$, $p(x_2)=Q\rightarrow B$, which means $\phi_2\rightarrow \phi_*$. Thus, \eqref{eq:theta2phi2theta1PhiEquation} implies $\theta_2=\phi_2-(\Phi+\theta_1)\longrightarrow\phistar-(\Phistar+\frac{1}{n}\Phistar)=\phistar-\frac{n+1}{n}\Phistar$.
     
     \textbf{Proof of \eqref{item02Archsingularitycurvesproperties}:} with the same notation as \cref{def:fpclf}, we denote $L_1=\tau_1=|PQ|$ (see \cref{fig:ArchSingularitiesneq1,fig:ArchSingularitiesneq2}). Since $x_1$ is on the line that runs through $(-\Phistar,0)$ with slope $1/2n$, \begin{equation}\label{eq:ArchSingularitytheta1derivative}\frac{\theta_1}{\Phi+\Phistar}=\frac{1}{2n}\Longrightarrow\theta_1=\frac{\Phi+\phistar}{2n}\Longrightarrow\frac{d\theta_1}{d\Phi}=\frac{1}{2n}.\end{equation}
     Since the vector $\overrightarrow{PQ}$ direction has an angle $\Phi+\theta_1$ with respect to the positive $x$ axis and Q is on the circle arc $\Gamma_r$, \begin{align*}
        &\left.
        \begin{aligned}
         &x_Q=L_1\cos{(\Phi+\theta_1)}+x_P\\
         &y_Q=L_1\sin{(\Phi+\theta_1)}+y_P\\
         &x^2_Q+y^2_Q=r^2
        \end{aligned}\right\}\Rightarrow L^2_1+2\big(x_P\cos{(\Phi+\theta_1)}+y_P\sin{(\Phi+\theta_1)}\big)L_1+x^2_P+y^2_P=r^2
     \end{align*}
Since $x_P=R\sin\Phi$, $y_P=b-R\cos\Phi$,
\begin{align*}
 x^2_P+y^2_P&=R^2+b^2-rbR\cos\Phi,\\ x_P\cos{(\Phi+\theta_1)}+y_P\sin{(\Phi+\theta_1)}&=b\sin{(\Phi+\theta_1)}-R\sin{(\Phi+\theta_1)}\cos\Phi+R\sin\Phi\cos{(\Phi+\theta_1)}\\
 &=b\sin{(\Phi+\theta_1)}-R\sin\theta_1.
\end{align*}
Therefore, \[L_1^2+2\big[b\sin{(\Phi+\theta_1)}-R\sin\theta_1\big]L_1+R^2-2bR\cos\Phi+b^2-r^2=0.\]
Then taking the derivative with respect to $\Phi$ in the equation above both sides and using the chain rule gives \[2L_1\frac{dL_1}{d\Phi}+2\big(b\sin{(\Phi+\theta_1)-R\sin\theta_1}\big)\frac{dL_1}{d\Phi}+2L_1\big[b\cos{(\Phi+\theta_1)}\frac{d(\Phi+\theta_1)}{d\Phi}-R\cos\theta_1\frac{d\theta_1}{d\Phi}\big]+2bR\sin\Phi=0.\]
Then with \eqref{eq:ArchSingularitytheta1derivative}: $\frac{d\theta_1}{d\Phi}=\frac{1}{2n}$, the above equation further gives \begin{align*}
    \frac{dL_1}{d\Phi}\big[L_1+b\sin{(\Phi+\theta_1)}-R\sin\theta_1\big]+L_1\big[b(1+\frac{1}{2n})\cos{(\Phi+\theta_1)}-\frac{1}{2n}R\cos\theta_1\big]+bR\sin\Phi=0.
\end{align*} Hence,
\begin{equation}\label{eq:dL1dPhiArchSingularityProperties}
\begin{aligned}
    \frac{dL_1}{d\Phi}\big[L_1+b\sin{(\Phi+\theta_1)}-R\sin\theta_1\big]&=\frac{1}{2n}L_1R\cos\theta_1-L_1b(1+\frac{1}{2n})\cos{(\theta_1+\Phi)}-bR\sin\Phi\\
    &=\frac{1}{2n}L_1\big[R\cos\theta_1-b\cos{(\theta_1+\Phi)}\big]-b\big[R\sin\Phi+L_1\cos{(\theta_1+\Phi)}\big]
\end{aligned}
\end{equation}
By the same analysis inside \eqref{dL1Dtheta} for $d_2=L_1+b\sin{(\Phi+\theta)}-R\sin\theta$ in \cref{fig:CaseAMonotone} (See the 6 cases analysis based on $\tau_1,d_2,d_1$ lengths), now in left hand side of \eqref{eq:dL1dPhiArchSingularityProperties}: 
\begin{equation}\label{eq:ArchSingularitydL1Dtheta}
L_1+b\sin{(\Phi+\theta_1)}-R\sin\theta_1=d_2,
\end{equation}
because it is implied by the same 6 cases analysis following \eqref{dL1Dtheta}.

In the current context, \cref{fig:CaseAMonotone} with $\theta=\theta_1$ also implies \[R\sin\Phi+L_1\cos{(\theta_1+\Phi)}=x_Q\] and \[R\cos\theta_1-b\cos{(\theta_1+\Phi)}=\text{dist}(O_r,\text{the line going through } P,Q),\], which means that the right-hand side of \eqref{eq:dL1dPhiArchSingularityProperties}: \[\frac{1}{2n}L_1\big[R\cos\theta_1-b\cos{(\theta_1+\Phi)}\big]-b\big[R\sin\Phi+L_1\cos{(\theta_1+\Phi)}\big]=\frac{1}{n}\text{Area}(\triangle O_rQP)-bx_Q.\] It is not hard to see it is negative by the following argument. 

Since $x_Q=r\sin\phi_2$ and by $x_2\in\Nin$ region containing $(\phistar,\phistar)$ in \cref{fig:MrN}, $\phi_2$ satisfies $0<\phi_2-\phistar<17\sqrt{\frac{r}{R}}\frac{1}{\sin{(\phistar/2)}}$. 

With $0<\phistar<\tan^{-1}(1/3)$, \eqref{eqJZHypCond} implies that
\[\phi_2<\phistar+\frac{17}{\sin{(\phistar/2)}}\sqrt{r/R}<\tan^{-1}(1/3)+\frac{17}{\sin{(\phistar/2)}}\sqrt{\frac{\sin^2\phistar}{30000}}<\tan^{-1}(1/3)+\frac{34}{\sqrt{30000}}<\pi/2.\] 

Hence $0<x_P<r\sin\phistar<x_Q<r$. Suppose in \cref{fig:CaseAMonotone} that the vertical line going through $P$ intersects the line segment $O_rQ$ at $P'$. It is clear that the line segment $PP'$ is contained within the disk surrounded by the circle $C_r$ since the disk is convex. 

Therefore, $|PP'|<2r$. We observe that Area($\triangle(O_rQP)$)=$\frac{1}{2}|PP'|\cdot x_Q$, then \[\frac{1}{n}\text{Area}(\triangle O_rQP)-bx_Q=(\frac{1}{2n}|PP'|-b)x_Q<0,\] because $n\ge1$ and \eqref{eqJZHypCond}: $b>R-r>1699r$. Hence $\frac{d\tau_1}{d\Phi}=\frac{dL_1}{d\Phi}<0$. We continue to compute $\frac{d\phi_2}{d\Phi}$ and $\frac{d\theta_2}{d\Phi}$.

Since the line $PQ$ has slope $\tan{(\Phi+\theta_1)}$, by $P$ and $Q$ coordinates in \eqref{eq:ArchSingularityPQcoordinates}, we get the followings.
\begin{equation}\label{eq:ArchSingulairtyTrajectorySlope}
\begin{aligned}
    &\frac{-r\cos\phi_2-b+R\cos\Phi}{r\sin\phi_2-R\sin\Phi}=\frac{\sin{(\Phi+\theta_1)}}{\cos{(\Phi+\theta_1)}}\\
    \Longrightarrow&-r\cos\phi_2\cos{(\Phi+\theta_1)}-b\cos{(\Phi+\theta_1)}+R\cos\Phi\cos{(\Phi+\theta_1)}=r\sin\phi_2\sin{(\Phi+\theta_1)}-R\sin\Phi\sin{(\Phi+\theta_1)}\\
    \Longrightarrow&R\cos\theta_1-b\cos{(\Phi+\theta_1)}=r\cos{(\phi_2-\Phi-\theta_1)}\underbracket{=}_{\mathclap{\eqref{eq:theta2phi2theta1PhiEquation}:\theta_2=\phi_2-\Phi-\theta_1}}r\cos\theta_2.
\end{aligned}
\end{equation}
By taking derivative with respect to $\Phi$ in both hand sides of last equation of \eqref{eq:ArchSingulairtyTrajectorySlope}, using chain rule with \eqref{eq:ArchSingularitytheta1derivative}: $\frac{d\theta_1}{d\Phi}=\frac{1}{2n}$, we further get the following:
\begin{align*}
    &-R\sin\theta_1\frac{d\theta_1}{d\Phi}+b(1+\frac{d\theta_1}{d\Phi})\sin{(\Phi+\theta_1)}=-r\sin\theta_2\frac{d\theta_2}{d\Phi}\\
    \overbracket{\Longrightarrow}^{\mathllap{\frac{d\theta_1}{d\Phi}=\frac{1}{2n}}}&-\frac{1}{2n}R\sin\theta_1+b(1+\frac{1}{2n})\sin{(\Phi+\theta_1)}=-r\sin\theta_2\frac{d\theta_2}{d\Phi}
\end{align*}
Then using \cref{def:lengthfunctiononArchSingularityCurve}: $d_1=R\sin\theta_1$, $d_2=r\sin\theta_2$, $L_1=\tau_1$ and by \eqref{eq:ArchSingularitydL1Dtheta}: $b\sin{(\Phi+\theta_1)}=d_2+R\sin\theta_1-L_1$, we get: 
\begin{equation}\label{eq:ArchSingularitytheta2derivative}
\begin{aligned}
    &d_2\frac{d\theta_2}{d\Phi}=\frac{1}{2n}d_1-(1+\frac{1}{2n})(d_2+d_1-L_1)\\
    \Longrightarrow&\frac{d\theta_2}{d\Phi}=\frac{(1+\frac{1}{2n})(\tau_1-d_2)-d_1}{d_2}=\frac{(1+\frac{1}{2n})\tau_1-d_1}{d_2}-(1+\frac{1}{2n}).
\end{aligned}\end{equation}
Note that by taking derivative with respect to $\Phi$ in both hand sides of \eqref{eq:theta2phi2theta1PhiEquation} and by using \eqref{eq:ArchSingularitytheta1derivative} we get the following: 
\begin{equation}\label{eq:ArchSingularityphi2derivative}
\frac{d\phi_2}{d\Phi}=\frac{d\theta_2}{d\Phi}+(1+\frac{d\theta_1}{d\Phi})\qquad\overbracket{=}^{\mathllap{\text{\eqref{eq:ArchSingularitytheta1derivative}}}}\frac{d\theta_2}{d\Phi}+(1+\frac{1}{2n})=\frac{(1+\frac{1}{2n})\tau_1-d_1}{d_2}.\end{equation}
Note that as functions of $\Phi\in\big(\Phistar-\frac{2}{n+1}\Phistar,\Phistar\big)$, $\tau_1$ monotone decreases and $d_1=R\sin\Phi$ monotone increases. Hence, $(1+\frac{1}{2n})\tau_1-d_1<\lim_{\Phi\rightarrow\Phistar-\frac{2\Phistar}{n+1}}\big[(1+\frac{1}{2n})\tau_1-d_1\big]$.

Note that as $\Phi\rightarrow\Phistar-\frac{2\Phistar}{n+1}$, $Q=p(x_2)\rightarrow B\in\Gamma_R$, since $P\in\Gamma_R$ and $|PQ|=\tau_1$,
\begin{equation}\label{eq:ArchSingularityLengthInEquality}
\begin{aligned}
\lim_{\Phi\rightarrow\Phistar-\frac{2\Phistar}{n+1}}\tau_1=\lim_{\Phi\rightarrow\Phistar-\frac{2\Phistar}{n+1}}(2d_1)\overbracket{=}^{\mathclap{\text{\cref{def:lengthfunctiononArchSingularityCurve}}}}&\lim_{\Phi\rightarrow\frac{(n-1)\Phistar}{n+1}}2R\sin\theta_1\overbracket{=}^{\mathclap{\text{\eqref{eq:ArchSingularitytheta1derivative}}}}\lim_{\Phi\rightarrow\frac{(n-1)\Phistar}{n+1}}2R\sin{(\frac{\Phi+\Phistar}{2n})}=2R\sin{(\frac{\Phistar}{n+1})}.\\
\lim_{\Phi\rightarrow\Phistar-\frac{2\Phistar}{n+1}}\big[(1+\frac{1}{2n})\tau_1-d_1\big]=&\lim_{\Phi\rightarrow\frac{(n-1)\Phistar}{n+1}}\big[(1+\frac{1}{2n})(2d_1)-d_1\big]=\lim_{\Phi\rightarrow\frac{(n-1)\Phistar}{n+1}}(1+\frac{1}{n})R\sin\theta_1\\
=&(1+\frac{1}{n})R\sin{(\frac{\Phistar}{n+1})}
\end{aligned}
\end{equation}
Therefore, if $n=1$, then 
\begin{equation}\label{eq:ArchSingularityLengthInEquality1}
\begin{aligned}
    (1+\frac{1}{2n})\tau_1-d_1<&\lim_{\Phi\rightarrow\Phistar-\frac{2\Phistar}{n+1}}\big[(1+\frac{1}{2n})\tau_1-d_1\big]\overbracket{=}^{\mathclap{\text{\eqref{eq:ArchSingularityLengthInEquality} with }n=1}}2R\sin{(\frac{\Phistar}{2})}=\frac{R\sin\Phistar}{\cos{(\frac{\Phistar}{2})}}=\frac{r\sin\phistar}{\cos{(\frac{\Phistar}{2})}}\\
    \overbracket{<}^{\mathllap{\substack{\eqref{eqJZHypCond}:R>1700r,\\0<\Phistar=\sin^{-1}{(r\sin\phistar/R)}<0.00059}}}&\frac{r\sin\phistar}{0.9999}.
\end{aligned}
\end{equation}
If $n\ge2$, then 
\begin{equation}\label{eq:ArchSingularityLengthInEquality2}
\begin{aligned}
    (1+\frac{1}{2n})\tau_1-d_1<&\lim_{\Phi\rightarrow\Phistar-\frac{2\Phistar}{n+1}}\big[(1+\frac{1}{2n})\tau_1-d_1\big]\quad\overbracket{\le}^{\mathllap{\text{\eqref{eq:ArchSingularityLengthInEquality} with }n\ge2}}\frac{3}{2}\cdot R\sin{(\frac{\Phistar}{2})}=\frac{3}{4}\cdot\frac{R\sin\Phistar}{\cos{(\frac{\Phistar}{2})}}=\frac{3}{4}\cdot\frac{r\sin\phistar}{\cos{(\frac{\Phistar}{2})}}\\
    \overbracket{<}^{\mathllap{\substack{\eqref{eqJZHypCond}:R>1700r,\\0<\Phistar=\sin^{-1}{(r\sin\phistar/R)}<0.00059}}}&\frac{3}{4}\cdot\frac{r\sin\phistar}{0.9999}.
\end{aligned}
\end{equation}

On the other hand, for $x_2\in\mathcal{AS}^{\text{\upshape in}}_{-n}$, if $n=1$, then $x_2$ has a neighborhood of points that are $x_2$ in \eqref{eqPtsFromeqMhatOrbSeg} points in case (b) or case (c) of \eqref{eqMainCases}. If $n\ge2$, then $x_2$ has a neighborhood of points which are all to be $x_2$ in \eqref{eqPtsFromeqMhatOrbSeg} points in case (b) or case (c) of \eqref{eqMainCases}. Since $\theta_2$ is continuous at $x_2$, by \eqref{eq:casebghd0lengthrestrictions5} for case (b) and by \cref{proposition:casecparametersestimate}\ref{itemcasecparametersestimate5} for case (c), we get 
\begin{equation}\label{eq:ArchSingularityd2range}
d_2=r\sin\theta_2\in \big[0.997r\sin\phistar,1.003r\sin\phistar\big].
\end{equation}
Therefore, if $n=1$, then $(1+\frac{1}{2n})d_2=\frac{3}{2}d_2\ge\frac{3\cdot0.997r\sin\phistar}{2}>\frac{r\sin\phistar}{0.9999}\overbracket{>}^{\eqref{eq:ArchSingularityLengthInEquality1}}(1+\frac{1}{2n})\tau_1-d_1$. 

If $n\ge2$, then $(1+\frac{1}{2n})d_2>d_2\ge0.997r\sin\phistar>\frac{3}{4}\cdot\frac{r\sin\phistar}{0.9999}\overbracket{>}^{\eqref{eq:ArchSingularityLengthInEquality2}}(1+\frac{1}{2n})\tau_1-d_1$. Hence, \[\frac{d\theta_2}{d\Phi}=\frac{(1+\frac{1}{2n})\tau_1-d_1}{d_2}-(1+\frac{1}{2n})=\frac{(1+\frac{1}{2n})\tau_1-d_1-(1+\frac{1}{2n})d_2}{d_2}<0.\] for all $n\ge1$.

\textbf{Proof of \eqref{item03Archsingularitycurvesproperties}:} from \eqref{eq:ArchSingularityphi2derivative},\eqref{eq:ArchSingularitytheta2derivative}, we get the tangential slope of $\mathcal{AS}^{\text{\upshape in}}_{-n}$ at $x_2$ is also a function of $\Phi$:
\begin{equation}\label{eq:ArchSingularityDtheta2Dphi2}
\begin{aligned}
    \mathcal{V}(x_2)&=\mathcal{V}(\Phi)=\frac{d\theta_2}{d\phi_2}=\frac{\frac{d\theta_2}{d\Phi}}{\frac{d\phi_2}{d\Phi}}=\frac{[(1+\frac{1}{2n})\tau_1-d_1]/d_2-(1+\frac{1}{2n})-(1+\frac{1}{2n})}{[(1+\frac{1}{2n})\tau_1-d_1]/d_2}\\
    &=1-\frac{(1+\frac{1}{2n})d_2}{(1+\frac{1}{2n})\tau_1-d_1}=1-\frac{(1+\frac{1}{2n})d_2}{(1+\frac{1}{2n})\tau_1-d_1}=1-\frac{d_2}{\tau_1-\frac{2n}{2n+1}d_1}.
\end{aligned}    
\end{equation}
In \eqref{eq:ArchSingularityLengthInEquality}, we get $\lim_{\Phi\rightarrow\Phistar-\frac{2\Phistar}{n+1}}\tau_1=2R\sin{(\frac{\Phistar}{n+1})}$, $\lim_{\Phi\rightarrow\Phistar-\frac{2\Phistar}{n+1}}d_1=R\sin{(\frac{\Phistar}{n+1})}$.

Since $\Phi\rightarrow\Phistar$, $P=p(x_1)\rightarrow B$ and $Q=p(x_2)\rightarrow B$, \[\lim_{\Phi\rightarrow\Phistar}\tau_1=\lim_{\Phi\rightarrow\Phistar}|PQ|=0.\] 

By \eqref{eq:ArchSingularitytheta1derivative}, \[\lim_{\Phi\rightarrow\Phistar-\frac{2\Phistar}{n+1}}d_1\overbracket{=}^{\mathclap{\text{\cref{def:lengthfunctiononArchSingularityCurve}}}}\lim_{\Phi\rightarrow\Phistar}R\sin\theta_1\overbracket{=}^{\eqref{eq:ArchSingularitytheta1derivative}}\lim_{\Phi\rightarrow\Phistar}R\sin{(\frac{\Phi+\Phistar}{2n})}=R\sin{(\frac{\Phistar}{n})}.\]

Therefore, $\lim_{\Phi\rightarrow\Phistar}(\tau_1-\frac{2n}{2n+1}d_1)=\frac{2n+2}{2n+1}R\sin{(\frac{\Phistar}{n})}>0$ and $\lim_{\Phi\rightarrow\Phistar}(\tau_1-\frac{2n}{2n+1}d_1)=-\frac{2n}{2n+1}R\sin{(\frac{\Phistar}{n})}<0$. Note that $\tau_1-\frac{2n}{2n+1}d_1$ is a strictly monotone decreasing smooth function of $\Phi\in(\Phistar-\frac{2\Phistar}{n+1},\Phistar)$. Therefore, by the intermediate value theorem, there exists a $\Phi_{c,n}\in(\Phistar-\frac{2\Phistar}{n+1},\Phistar)$ such that $(\tau_1-\frac{2n}{2n+1}d_1)\bigm|_{\Phi=\Phi_{c,n}}=0$, $\tau_1-\frac{2n}{2n+1}d_1>0$ for $\Phi\in(\Phistar-\frac{2\Phistar}{n+1},\Phi_{c,n})$, $\tau_1-\frac{2n}{2n+1}d_1<0$ for $\Phi\in(\Phi_{c,n},\Phistar)$. Then we get the following. 

If $\Phi\in(\Phistar-\frac{2\Phistar}{n+1},\Phi_{c,n}]$, then $0\le\tau_1-\frac{2n}{2n+1}d_1<\frac{2n+2}{2n+1}R\sin{(\frac{\Phistar}{n})}=(1+\frac{1}{2n+1})R\sin{(\frac{\Phistar}{n+1})}\overbracket{\le}^{\mathclap{n\ge1}}\frac{4}{3}R\sin{(\frac{\Phistar}{2})}$. 

Especially if $\Phi\in(\Phistar-\frac{2\Phistar}{n+1},\Phi_{c,n}]$ and $n\ge2$, then $0\le\tau_1-\frac{2n}{2n+1}d_1<\frac{4}{5}\sin{(\frac{\Phistar}{3})}<\frac{4}{5}\sin{(\frac{\Phistar}{2})}$.

If $\Phi\in(\Phi_{c,n},\Phistar)$, then $\tau_1-\frac{2n}{2n+1}d_1<0$.

Therefore, with the $d_2$ range in \eqref{eq:ArchSingularityd2range}, we also get the slope range:

If $\Phi\in(\Phistar-\frac{2\Phistar}{n+1},\Phi_{c,n}]$, then
\begin{equation}\label{eq:slopeForArchSingularityCruveInNeq1}
\begin{aligned}
    \mathcal{V}(x_2)=\mathcal{V}(\Phi)&\overbracket{=}^{\mathclap{\eqref{eq:ArchSingularityDtheta2Dphi2}}}1-\frac{d_2}{\tau_1-\frac{2n}{2n+1}d_1}<1-\frac{d_2}{\frac{4}{3}R\sin{(\frac{\Phistar}{2})}}=1-3\cdot\frac{d_2\cos{(\frac{\Phistar}{2})}}{4R\sin{(\frac{\Phistar}{2})}\cos{(\frac{\Phistar}{2})}}\\
    &\overbracket{<}^{\mathllap{\eqref{eq:ArchSingularityd2range}}}1-3\cdot\frac{0.997r\sin\phistar\cos{(\frac{\Phistar}{2})}}{2R\sin\Phistar}\qquad\quad\overbracket{=}^{\mathclap{R\sin\Phistar=r\sin\phistar}}\qquad 1-\frac{3\cdot0.997}{2}\cos{(\Phistar/2)}\qquad\overbracket{<}^{\mathclap{\substack{\eqref{eqJZHypCond}:R>1700r,\\0<\Phistar=\sin^{-1}{(r\sin\phistar/R)}<0.00059}}}\qquad -0.495
\end{aligned}    
\end{equation}

Especially if $\Phi\in(\Phistar-\frac{2\Phistar}{n+1},\Phi_{c,n}]$ and $n\ge2$, then\begin{align*}
    \mathcal{V}(x_2)=\mathcal{V}(\Phi)&\overbracket{=}^{\mathclap{\eqref{eq:ArchSingularityDtheta2Dphi2}}}1-\frac{d_2}{\tau_1-\frac{2n}{2n+1}d_1}<1-\frac{d_2}{\frac{4}{5}R\sin{(\frac{\Phistar}{2})}}<1-5\cdot\frac{d_2\cos{(\Phistar/2)}}{4R\sin{(\Phistar/2)}\cos{(\Phistar/2)}}\overbracket{<}^{\mathclap{\eqref{eq:ArchSingularityd2range}}}1-5\cdot\frac{0.997r\sin\phistar\cos{(\frac{\Phistar}{2})}}{2R\sin\Phistar}\\
    &\overbracket{=}^{\mathllap{R\sin\Phistar=r\sin\phistar}}1-\frac{5\cdot0.997}{2}\cos{(\Phistar/2)}\qquad\overbracket{<}^{\mathclap{\substack{\eqref{eqJZHypCond}:R>1700r,\\0<\Phistar=\sin^{-1}{(r\sin\phistar/R)}<0.00059}}}-1.49
\end{align*}

If $\Phi\in(\Phi_{c,n},\Phistar)$, then\begin{align*}
    \mathcal{V}(x_2)=\mathcal{V}(\Phi)&\overbracket{=}^{\mathclap{\eqref{eq:ArchSingularityDtheta2Dphi2}}}1+\underbracket{\frac{d_2}{-(\tau_1-\frac{2n}{2n+1}d_1)}}_{>0}>1.
\end{align*}
\end{proof}

\begin{remark}\label{remark:ArchSingualrityCurveAlmostBelowNegativeHalf}
    From the computation of \eqref{eq:slopeForArchSingularityCruveInNeq1}, except for a short curve piece  of $\mathcal{AS}^{\text{\upshape in}}_{-1}$ in the neighborhood of end points $(\phistar,\Phistar-\frac{\Phistar}{2})$ and $(2\pi-\phistar,\pi-(\Phistar-\frac{\Phistar}{2}))$, $\mathcal{AS}^{\text{\upshape in}}_{-1}$ tangential slopes on points of other points on $\mathcal{AS}^{\text{\upshape in}}_{-1}$ and all tangential slopes on $\mathcal{AS}^{\text{\upshape in}}_{-n}$ with $n\ge2$ are in $[-\infty,-0.5)\cup(1,\infty]$. The \emph{exceptional short curve pieces of $\mathcal{AS}^{\text{\upshape in}}_{-1}$} have tangential slopes in $[-0.5,-0.495)$.
\end{remark}

\begin{proposition}[The other component of $\mathcal{AS}^{\text{\upshape in}}_{-n}$]\label{proposition:Archsingularitycurvesproperties2}Since the two components of $\mathcal{AS}^{\text{\upshape in}}_{-n}$ are symmetric to each other by $I\circ J:(\phi,\theta)\rightarrow (2\pi-\phi,\pi-\theta)$, the connected component of $\mathcal{AS}^{\text{\upshape in}}_{-n}$ with the $n$ steps in the past on $\Gamma_R$ and the $n+1$-th step in the past starting from the corner $B$ satisfies the following. 
    \begin{enumerate}
        \item It is a smooth curve with endpoints $\big(2\pi-\phistar,\pi-(\phistar-\frac{n}{n+1}\phistar)\big)$ and $\big(2\pi-\phistar,\pi-(\phistar-\frac{n+1}{n}\phistar)\big)$.
        \item For $n\ge1$, the range of values of tangential slopes of $\mathcal{AS}^{\text{\upshape in}}_{-1}$ is contained in $\big[-\infty,-0.495\big)\cup\big(1,+\infty\big]$. For $n\ge2$, the range of values of tangential slopes is contained in $\big[-\infty,-1.49\big)\cup\big(1,+\infty\big]$. (See \cref{fig:ArchSingularitiesComponent2})
    \end{enumerate}
\end{proposition}
\begin{proof}
    This is just a symmetric conclusion of \cref{proposition:Archsingularitycurvesproperties2} by the symmetry $I\circ J$.
\end{proof}
\begin{figure}[ht!]
\begin{minipage}{.5\textwidth}
\begin{center}
    \scalebox{2.6}{\input{figs/Figure_2.pgf}}
    \caption{Arch Singularity curves: $\mathcal{AS}^{\text{\upshape in}}_{-n}$ with $n=1,2,3,4$ near $(2\pi-\phistar,\pi-(\phistar-\Phistar))$ as a zoom in to $\Nin$ of \cref{fig:MrN}.}\label{fig:ArchSingularitiesComponent2}
\end{center}
\end{minipage}\hfill
\begin{minipage}{.5\textwidth}
\begin{center}
    \scalebox{2.82}{\input{figs/Figure_3.pgf}}
    \caption{Arch Singularity curves: $\mathcal{AS}^{\text{\upshape out}}_{n}$ with $n=1,2,3,4$ near $(2\pi-\phistar,\phistar-\Phistar)$ as a zoom in to $\Nout$ of \cref{fig:MrN}.}\label{fig:ArchSingularitiesComponent3}
\end{center}
\end{minipage}\hfill
\end{figure}
\begin{definition}[Analogous to \cref{def:ASinCurves}]
    We define for each $n\ge1$ the \emph{Arch Singularity Curves}\[\mathcal{AS}^{\text{\upshape out}}_{n}\dfn\big\{x\in\Mrout\bigm|\Fc^{(n+1)}(x)\text{ is at corner and }\Fc^{k}(x)\in\Int(M_R),k=1,2,\cdots,n\big\}.\]
\end{definition}
\begin{remark}
    Each $x\in\mathcal{AS}^{\text{\upshape out}}_{n}$ is the collision on $\Gamma_r$ that starts a trajectory in future $n+1$ steps with $n$ consecutive collisions on $\Gamma_R$ ending at a corner. Each $\mathcal{AS}^{\text{\upshape out}}_{n}$ has two connected components located in the neighborhoods of $(\phi_*,\phi_*-\Phi_*)$ and $(2\pi-\phi_*,\pi-\phi_*+\Phi_*)$ on the boundary of $M_r$.
\end{remark}
\begin{proposition}[$\mathcal{AS}^{\text{\upshape out}}_n$ tangential slopes, analogous to \cref{proposition:Archsingularitycurvesproperties,proposition:Archsingularitycurvesproperties2}]\label{proposition:Archsingularitycurvesproperties3}\hfill
\begin{enumerate}
\item\label{item01Archsingularitycurvesproperties2} The $\mathcal{AS}^{\text{\upshape out}}_1$ tangential slopes at all points have the range of values contained in $[-\infty,-1)\cup(0.495,+\infty]$.
    \item\label{item02Archsingularitycurvesproperties2} For each $n\ge2$, the $\mathcal{AS}^{\text{\upshape out}}_n$ tangential slopes at all points have the range of values contained in $[-\infty,-1)\cup(1.49,+\infty]$.
\end{enumerate}    
\end{proposition}
\begin{proof}
    Conclusions are implied by \cref{proposition:Archsingularitycurvesproperties} and by the fact that in $M_r$ $\mathcal{AS}^{\text{\upshape out}}_n$ and $\mathcal{AS}^{\text{\upshape in}}_{-n}$ are curves symmetric to each other by $I$ or $J$ (given in \cref{DefSymmetries}).
\end{proof}
\begin{definition}[curvilinear quadrilateral]\label{def:Curvilinearquadrilateral}
    We define $\mathcal{D}^{\text{\upshape out}}_n$ as two components of the quadrilateral curvilinear region.

    A component near $(2\pi-\phistar,\phistar-\Phistar)$ (in \Nout) has the following four curves as boundaries: 
    \begin{itemize}
        \item curve $\mathcal{AS}^{\text{\upshape out}}_{n}$ (with endpoints $(2\pi-\phistar,\phistar-\frac{n}{n+1}\Phistar)$, $(2\pi-\phistar,\phistar-\frac{n+1}{n}\Phistar)$),
        \item straight line segment of $\partial M_r$ connecting $(2\pi-\phistar,\phistar-\frac{n}{n+1}\Phistar)$ and $(2\pi-\phistar,\phistar-\frac{n}{n+1}\Phistar)$,
        \item curve $\mathcal{AS}^{\text{\upshape out}}_{n+1}$ (with endpoints $(2\pi-\phistar,\phistar-\frac{n+1}{n+2}\Phistar)$, $(2\pi-\phistar,\phistar-\frac{n+2}{n+1}\Phistar)$),
        \item straight line segment of $\partial M_r$ connecting $(2\pi-\phistar,\phistar-\frac{n+2}{n+1}\Phistar)$ and $(2\pi-\phistar,\phistar-\frac{n+1}{n}\Phistar)$.
    \end{itemize}    
    
    Symmetrically, the other component of $\mathcal{AS}^{\text{\upshape out}}_{n+1}$ near $(\phistar,\pi-(\phistar-\Phistar))$ is defined as the curvilinear quadrilateral with boundaries $\mathcal{AS}^{\text{\upshape out}}_{n}$, $\mathcal{AS}^{\text{\upshape out}}_{n+1}$ and two straight line segments on $\partial M_r$.

    Similarly, we define $\mathcal{D}^{\text{\upshape in}}_{-n}$ having two components as the  curvilinear quadrilateral with boundaries $\mathcal{AS}^{\text{\upshape in}}_{-n}$, $\mathcal{AS}^{\text{\upshape in}}_{-(n+1)}$ and two straight line segments on $\partial M_r$.
\end{definition}
\begin{remark}[$\Fc$ action on curvilinear quadrilaterals]\label{remark:FcActionOnQuadrilaterals}\hfill

    $\mathcal{D}^{\text{out}}_{n}$ is diffeomorphic to $\mathcal{D}^{\text{in}}_{-n}$ by $\Fc^{n+2}$. The images in orbit $\mathcal{D}^{\text{out}}_{n},\Fc(\mathcal{D}^{\text{out}}_{n}),\cdots,\Fc^{n+2}(\mathcal{D}^{\text{out}}_{n})$ are $\mathcal{D}^{\text{out}}_{n},\Min_{R,n},\cdots,\Fc^{n}(\Min_{R,n}),\mathcal{D}^{\text{in}}_{-n}$ (shown in \cref{fig:D1celDynamics,fig:D2celDynamics} for $\mathcal{D}^{\text{out}}_{1}$ and $\mathcal{D}^{\text{out}}_{2}$). The $\Fc$ actions on curvilinear quadrilaterals shown in \cref{fig:D1celDynamics,fig:D2celDynamics} are analogous to the $\Fc$ actions rhombus regions in \cite[FIGURE 8.11]{cb} for the Stadium Billiard.
\end{remark}
\begin{remark}[$\Fc$ as actions on boundaries of $\MRout$ and $\MRin$ and of $\MRout$ and $\Mrin$]\hfill

In \cref{fig:F_actMroutMRin}, $\Fc$ as a diffeomorphism between $\Mrout$ and $\MRin$ has a continuous extension to $\partial\Mrout$ so that this extension acts as a homeomorphism from $\partial\Mrout$ to $\partial\MRout$.

In \cref{fig:F_actMRoutMrin}, $\Fc$ as a diffeomorphism between $\MRout$ and $\Mrin$ has a continuous extension to $\partial\MRout$ so that this extension acts as a homeomorphism from $\partial\MRout$ to $\partial\Mrin$. 
\end{remark}

\begin{figure}[h!]
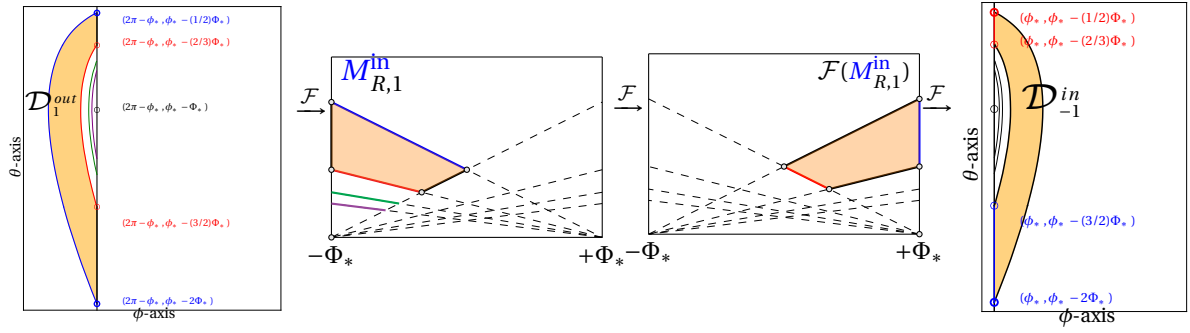

\begin{minipage}{.23\textwidth}
\begin{center}
\scalebox{2.05}{\input{figs/CellFigure_out.pgf}}
\end{center}
\end{minipage}
\begin{minipage}{.25\textwidth}
\begin{center}
    \begin{tikzpicture}[xscale=0.3, yscale=0.3]
        \pgfmathsetmacro{\PHISTAR}{6};
        \pgfmathsetmacro{\WIDTH}{8};
        \tkzDefPoint(-\PHISTAR,0){LL};
        \tkzDefPoint(-\PHISTAR,\WIDTH){LU};
        \tkzDefPoint(\PHISTAR,0){RL};
        \tkzDefPoint(\PHISTAR,\WIDTH){RU};
        \tkzDefPoint(-\PHISTAR,\PHISTAR){X};
        \tkzDefPoint(\PHISTAR,\WIDTH-\PHISTAR){Y};
        \tkzDefPoint(-\PHISTAR,\WIDTH-\PHISTAR){IX};
        \tkzDefPoint(\PHISTAR,\PHISTAR){IY};
        \tkzDefPoint(\PHISTAR,0.5*\PHISTAR){S1};
        \tkzDefPoint(\PHISTAR,\PHISTAR/3){S2};
        \tkzDefPoint(\PHISTAR,\PHISTAR/4){S3};
        \tkzDefPoint(\PHISTAR,\PHISTAR/5){L4};
        \tkzDefPoint(-\PHISTAR,0.5*\PHISTAR){LS1};
        \tkzDefPoint(-\PHISTAR,\PHISTAR/3){LS2};
        \tkzDefPoint(-\PHISTAR,\PHISTAR/4){LS3};
        \draw [name path=lineRLLS1,ultra thin,dashed] (RL) --(LS1);
        \draw [name path=lineRLLS2,ultra thin,dashed] (RL) --(LS2);
        \draw [name path=lineRLLS3,ultra thin,dashed] (RL) --(LS3);
        \draw [thin] (LL) --(LU);
        \draw [thin] (LU) --(RU);
        \draw [thin] (RL) --(RU);
        \draw [thin] (LL) --(RL);
        \draw [name path=lineLLIY, ultra thin,dashed] (LL) --(IY);
        \draw [name path=lineXRL,ultra thin,dashed]  (X)--(RL);
        \draw [name path=lineS1LL,ultra thin,dashed]  (S1)--(LL);
        \draw [name path=lineS2LL,ultra thin,dashed]  (S2)--(LL);
        \draw [name path=lineS3LL,ultra thin,dashed]  (S3)--(LL);
        \path[name intersections= 
        {of=lineLLIY and lineXRL,by={M0}}];
        \draw[blue,thick](M0)--(X);
        \path[name intersections= 
        {of=lineS1LL and lineXRL,by={M1}}];
        \path[name intersections= 
        {of=lineS2LL and lineXRL,by={M2}}];
        \path[name intersections= 
        {of=lineS3LL and lineXRL,by={M3}}];

        \path[name intersections= 
        {of=lineRLLS1 and lineLLIY,by={LM1}}];
        \draw[Red,thick](LM1)--(LS1);

        \path[name intersections= 
        {of=lineRLLS2 and lineLLIY,by={LM2}}];
        \draw[Green,thick](LM2)--(LS2);

        \path[name intersections= 
        {of=lineRLLS3 and lineLLIY,by={LM3}}];
        \draw[Purple,thick](LM3)--(LS3);
        \draw[thick](LM1)--(M0);
        \draw[thick](LS1)--(X);

         \fill [orange, opacity=10/30](X) -- (M0) -- (LM1) -- (LS1) -- cycle;
        \tkzDrawPoints(X,M0,LL,LM1,LS1);
        \tkzLabelPoint[above right](X){\Large $\Min_{R,1}$};

        \node[] at (-1.15*\PHISTAR,1.0*\PHISTAR){ $\xrightarrow{\Fc}$};
        \tkzLabelPoint[below](LL){$-\Phistar$};
        \tkzLabelPoint[below](RL){$+\Phistar$};
    \end{tikzpicture}
     \end{center}
\end{minipage}
\begin{minipage}{.25\textwidth}
\begin{center}
    \begin{tikzpicture}[xscale=0.3, yscale=0.3]
        \pgfmathsetmacro{\PHISTAR}{6};
        \pgfmathsetmacro{\WIDTH}{8};
        \tkzDefPoint(-\PHISTAR,0){LL};
        \tkzDefPoint(-\PHISTAR,\WIDTH){LU};
        \tkzDefPoint(\PHISTAR,0){RL};
        \tkzDefPoint(\PHISTAR,\WIDTH){RU};
        \tkzDefPoint(-\PHISTAR,\PHISTAR){X};
        \tkzDefPoint(\PHISTAR,\WIDTH-\PHISTAR){Y};
        \tkzDefPoint(-\PHISTAR,\WIDTH-\PHISTAR){IX};
        \tkzDefPoint(\PHISTAR,\PHISTAR){IY};
        \tkzDefPoint(\PHISTAR,0.5*\PHISTAR){S1};
        \tkzDefPoint(\PHISTAR,\PHISTAR/3){S2};
        \tkzDefPoint(\PHISTAR,\PHISTAR/4){S3};
        \tkzDefPoint(\PHISTAR,\PHISTAR/5){L4};
        \tkzDefPoint(-\PHISTAR,0.5*\PHISTAR){LS1};
        \tkzDefPoint(-\PHISTAR,\PHISTAR/3){LS2};
        \tkzDefPoint(-\PHISTAR,\PHISTAR/4){LS3};
        \draw [name path=lineRLLS1,ultra thin,dashed] (RL) --(LS1);
        \draw [name path=lineRLLS2,ultra thin,dashed] (RL) --(LS2);
        \draw [name path=lineRLLS3,ultra thin,dashed] (RL) --(LS3);

        \draw [thin] (LL) --(LU);
        \draw [thin] (LU) --(RU);
        \draw [thin] (RL) --(RU);
        \draw [thin] (LL) --(RL);
        \draw [name path=lineLLIY,ultra thin,dashed] (LL) --(IY);
        \draw [name path=lineXRL,ultra thin,dashed]  (X)--(RL);
        \draw [name path=lineS1LL,ultra thin,dashed]  (S1)--(LL);
        \draw [name path=lineS2LL,ultra thin,dashed]  (S2)--(LL);
        \draw [name path=lineS3LL,ultra thin,dashed]  (S3)--(LL);
        \path[name intersections= 
        {of=lineLLIY and lineXRL,by={M0}}];
        \draw[thick](M0)--(IY);
        \path[name intersections= 
        {of=lineS1LL and lineXRL,by={M1}}];
        \path[name intersections= 
        {of=lineS2LL and lineXRL,by={M2}}];
        \path[name intersections= 
        {of=lineS3LL and lineXRL,by={M3}}];

        \path[name intersections= 
        {of=lineRLLS1 and lineLLIY,by={LM1}}];
        \draw[thick](M1)--(S1);

        \path[name intersections= 
        {of=lineRLLS2 and lineLLIY,by={LM2}}];

        \path[name intersections= 
        {of=lineRLLS3 and lineLLIY,by={LM3}}];
        \draw[blue, thick](IY)--(S1);
        \draw[red, thick](M0)--(M1);

         \fill [orange, opacity=10/30](M1) -- (S1) -- (IY) -- (M0) -- cycle;
        \tkzDrawPoints(IY,M0,RL,M1,S1);
        \tkzLabelPoint[above left](IY){ $\Fc(\Min_{R,1})$};
        \node[] at (1.13*\PHISTAR,1.0*\PHISTAR){ $\xrightarrow{\Fc}$};

        \node[] at (-1.15*\PHISTAR,1.0*\PHISTAR){ $\xrightarrow{\Fc}$};
        \tkzLabelPoint[below](LL){$-\Phistar$};
        \tkzLabelPoint[below](RL){$+\Phistar$};
    \end{tikzpicture}
\end{center}
\end{minipage}\hfill
\begin{minipage}{.25\textwidth}
\begin{center}
\scalebox{2.75}{\input{figs/CellFigure_in.pgf}}
\end{center}
\end{minipage}
\caption{$\Fc$ action on $\mathcal{D}^{\text{\upshape out}}_{1}$ and its images. $\mathcal{D}^{\text{\upshape out}}_{1}$ is diffeomorphic to $\mathcal{D}^{\text{\upshape in}}_{-1}$ by $\Fc^3$. The same colored boundary are diffeomorphically mapped by $\Fc$ in each step.}\label{fig:D1celDynamics}
\end{figure}
\begin{figure}[h!]
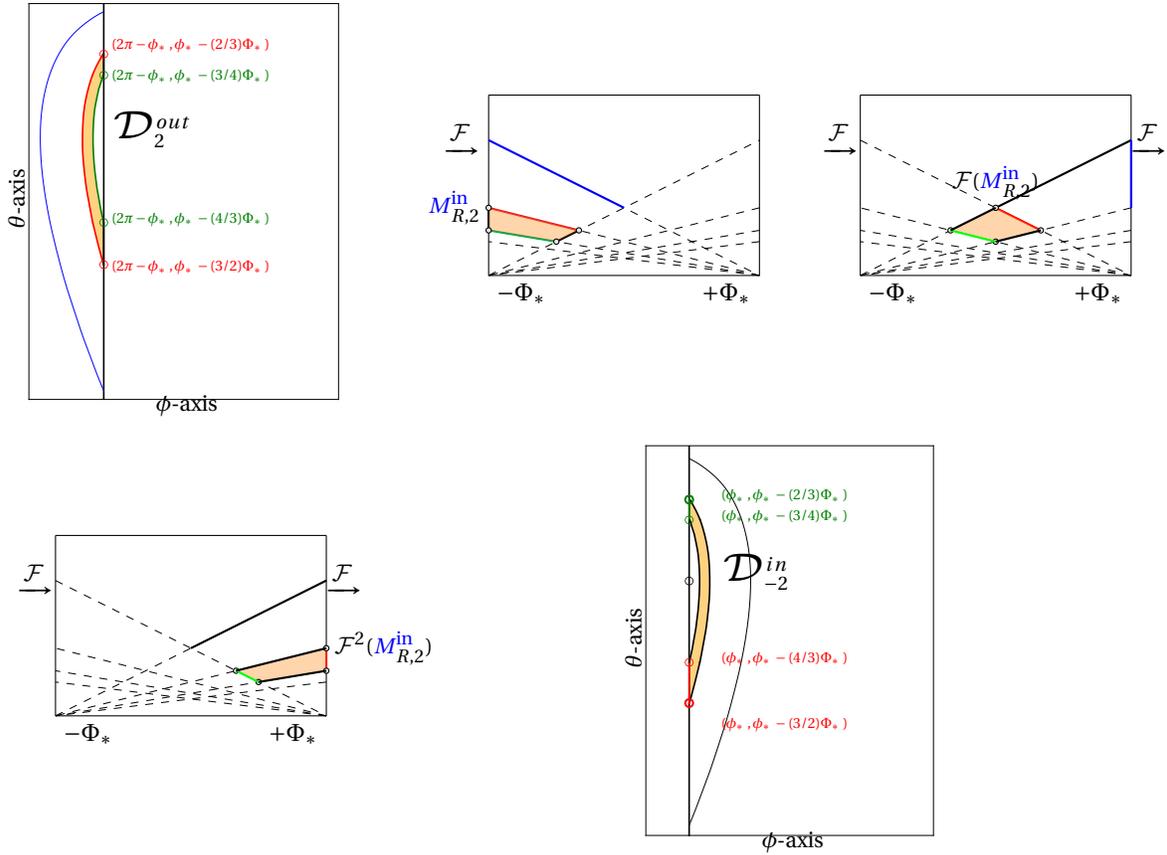

\begin{minipage}{.33\textwidth}
\begin{center}
\scalebox{3.0}{\input{figs/CellFigure_out2.pgf}}
\end{center}
\end{minipage}
\begin{minipage}{.32\textwidth}
\begin{center}
    \begin{tikzpicture}[xscale=0.3, yscale=0.3]
        \pgfmathsetmacro{\PHISTAR}{6};
        \pgfmathsetmacro{\WIDTH}{8};
        \tkzDefPoint(-\PHISTAR,0){LL};
        \tkzDefPoint(-\PHISTAR,\WIDTH){LU};
        \tkzDefPoint(\PHISTAR,0){RL};
        \tkzDefPoint(\PHISTAR,\WIDTH){RU};
        \tkzDefPoint(-\PHISTAR,\PHISTAR){X};
        \tkzDefPoint(\PHISTAR,\WIDTH-\PHISTAR){Y};
        \tkzDefPoint(-\PHISTAR,\WIDTH-\PHISTAR){IX};
        \tkzDefPoint(\PHISTAR,\PHISTAR){IY};
        \tkzDefPoint(\PHISTAR,0.5*\PHISTAR){S1};
        \tkzDefPoint(\PHISTAR,\PHISTAR/3){S2};
        \tkzDefPoint(\PHISTAR,\PHISTAR/4){S3};
        \tkzDefPoint(\PHISTAR,\PHISTAR/5){L4};
        \tkzDefPoint(-\PHISTAR,0.5*\PHISTAR){LS1};
        \tkzDefPoint(-\PHISTAR,\PHISTAR/3){LS2};
        \tkzDefPoint(-\PHISTAR,\PHISTAR/4){LS3};
        \draw [name path=lineRLLS1,ultra thin,dashed] (RL) --(LS1);
        \draw [name path=lineRLLS2,ultra thin,dashed] (RL) --(LS2);
        \draw [name path=lineRLLS3,ultra thin,dashed] (RL) --(LS3);
        \draw [thin] (LL) --(LU);
        \draw [thin] (LU) --(RU);
        \draw [thin] (RL) --(RU);
        \draw [thin] (LL) --(RL);
        \draw [name path=lineLLIY, ultra thin,dashed] (LL) --(IY);
        \draw [name path=lineXRL,ultra thin,dashed]  (X)--(RL);
        \draw [name path=lineS1LL,ultra thin,dashed]  (S1)--(LL);
        \draw [name path=lineS2LL,ultra thin,dashed]  (S2)--(LL);
        \draw [name path=lineS3LL,ultra thin,dashed]  (S3)--(LL);
        \path[name intersections= 
        {of=lineLLIY and lineXRL,by={M0}}];
        \draw[blue,thick](M0)--(X);
        \path[name intersections= 
        {of=lineS1LL and lineXRL,by={M1}}];
        \path[name intersections= 
        {of=lineS2LL and lineXRL,by={M2}}];
        \path[name intersections= 
        {of=lineS3LL and lineXRL,by={M3}}];

        \path[name intersections= 
        {of=lineRLLS1 and lineLLIY,by={LM1}}];
        \draw[Red,thick](LM1)--(LS1);

        \path[name intersections= 
        {of=lineRLLS2 and lineLLIY,by={LM2}}];
        \draw[Green,thick](LM2)--(LS2);

        \path[name intersections= 
        {of=lineRLLS3 and lineLLIY,by={LM3}}];
        \draw[thick](LM2)--(LM1);
        \draw[thick](LS1)--(LS2);

         \fill [orange, opacity=10/30](LS2) -- (LS1) -- (LM1) -- (LM2) -- cycle;
        \tkzDrawPoints(LS1,LS2,LM1,LM2);
        \tkzLabelPoint[left](LS1){\small $\Min_{R,2}$};

        \node[] at (-1.2*\PHISTAR,1.0*\PHISTAR){ \Large $\xrightarrow{\Fc}$};
        \tkzLabelPoint[below right](LL){$-\Phistar$};
        \tkzLabelPoint[below left](RL){$+\Phistar$};
    \end{tikzpicture}
     \end{center}
\end{minipage}
\begin{minipage}{.32\textwidth}
\begin{center}
    \begin{tikzpicture}[xscale=0.3, yscale=0.3]
        \pgfmathsetmacro{\PHISTAR}{6};
        \pgfmathsetmacro{\WIDTH}{8};
        \tkzDefPoint(-\PHISTAR,0){LL};
        \tkzDefPoint(-\PHISTAR,\WIDTH){LU};
        \tkzDefPoint(\PHISTAR,0){RL};
        \tkzDefPoint(\PHISTAR,\WIDTH){RU};
        \tkzDefPoint(-\PHISTAR,\PHISTAR){X};
        \tkzDefPoint(\PHISTAR,\WIDTH-\PHISTAR){Y};
        \tkzDefPoint(-\PHISTAR,\WIDTH-\PHISTAR){IX};
        \tkzDefPoint(\PHISTAR,\PHISTAR){IY};
        \tkzDefPoint(\PHISTAR,0.5*\PHISTAR){S1};
        \tkzDefPoint(\PHISTAR,\PHISTAR/3){S2};
        \tkzDefPoint(\PHISTAR,\PHISTAR/4){S3};
        \tkzDefPoint(\PHISTAR,\PHISTAR/5){L4};
        \tkzDefPoint(-\PHISTAR,0.5*\PHISTAR){LS1};
        \tkzDefPoint(-\PHISTAR,\PHISTAR/3){LS2};
        \tkzDefPoint(-\PHISTAR,\PHISTAR/4){LS3};
        \draw [name path=lineRLLS1,ultra thin,dashed] (RL) --(LS1);
        \draw [name path=lineRLLS2,ultra thin,dashed] (RL) --(LS2);
        \draw [name path=lineRLLS3,ultra thin,dashed] (RL) --(LS3);

        \draw [thin] (LL) --(LU);
        \draw [thin] (LU) --(RU);
        \draw [thin] (RL) --(RU);
        \draw [thin] (LL) --(RL);
        \draw [name path=lineLLIY,ultra thin,dashed] (LL) --(IY);
        \draw [name path=lineXRL,ultra thin,dashed]  (X)--(RL);
        \draw [name path=lineS1LL,ultra thin,dashed]  (S1)--(LL);
        \draw [name path=lineS2LL,ultra thin,dashed]  (S2)--(LL);
        \draw [name path=lineS3LL,ultra thin,dashed]  (S3)--(LL);
        \path[name intersections= 
        {of=lineLLIY and lineXRL,by={M0}}];
        \draw[thick](M0)--(IY);
        \path[name intersections= 
        {of=lineS1LL and lineXRL,by={M1}}];
        \path[name intersections= 
        {of=lineS2LL and lineXRL,by={M2}}];
        \path[name intersections= 
        {of=lineS3LL and lineXRL,by={M3}}];

        \path[name intersections= 
        {of=lineRLLS1 and lineLLIY,by={LM1}}];
        \path[name intersections= 
        {of=lineRLLS1 and lineS1LL,by={SM1}}];
        
        \path[name intersections= 
        {of=lineRLLS2 and lineLLIY,by={LM2}}];

        \path[name intersections= 
        {of=lineRLLS3 and lineLLIY,by={LM3}}];
        \draw[blue, thick](IY)--(S1);
        \draw[red, thick](M0)--(M1);

        \fill [orange, opacity=10/30](LM1) -- (SM1) -- (M1) -- (M0) -- cycle;
        \tkzDrawPoints(LM1,SM1,M1,M0);
        \draw[green, thick](LM1)--(SM1);
        \draw[thick](LM1)--(M0);
        \draw[thick](SM1)--(M1);
        \tkzLabelPoint[above](M0){\small $\Fc(\Min_{R,2})$};
        \node[] at (1.13*\PHISTAR,1.0*\PHISTAR){ \Large $\xrightarrow{\Fc}$};

        \node[] at (-1.15*\PHISTAR,1.0*\PHISTAR){ \Large $\xrightarrow{\Fc}$};
        \tkzLabelPoint[below right](LL){$-\Phistar$};
        \tkzLabelPoint[below left](RL){$+\Phistar$};
    \end{tikzpicture}
\end{center}
\end{minipage}\hfill
\begin{minipage}{.32\textwidth}
\begin{center}
    \begin{tikzpicture}[xscale=0.3, yscale=0.3]
        \pgfmathsetmacro{\PHISTAR}{6};
        \pgfmathsetmacro{\WIDTH}{8};
        \tkzDefPoint(-\PHISTAR,0){LL};
        \tkzDefPoint(-\PHISTAR,\WIDTH){LU};
        \tkzDefPoint(\PHISTAR,0){RL};
        \tkzDefPoint(\PHISTAR,\WIDTH){RU};
        \tkzDefPoint(-\PHISTAR,\PHISTAR){X};
        \tkzDefPoint(\PHISTAR,\WIDTH-\PHISTAR){Y};
        \tkzDefPoint(-\PHISTAR,\WIDTH-\PHISTAR){IX};
        \tkzDefPoint(\PHISTAR,\PHISTAR){IY};
        \tkzDefPoint(\PHISTAR,0.5*\PHISTAR){S1};
        \tkzDefPoint(\PHISTAR,\PHISTAR/3){S2};
        \tkzDefPoint(\PHISTAR,\PHISTAR/4){S3};
        \tkzDefPoint(\PHISTAR,\PHISTAR/5){L4};
        \tkzDefPoint(-\PHISTAR,0.5*\PHISTAR){LS1};
        \tkzDefPoint(-\PHISTAR,\PHISTAR/3){LS2};
        \tkzDefPoint(-\PHISTAR,\PHISTAR/4){LS3};
        \draw [name path=lineRLLS1,ultra thin,dashed] (RL) --(LS1);
        \draw [name path=lineRLLS2,ultra thin,dashed] (RL) --(LS2);
        \draw [name path=lineRLLS3,ultra thin,dashed] (RL) --(LS3);

        \draw [thin] (LL) --(LU);
        \draw [thin] (LU) --(RU);
        \draw [thin] (RL) --(RU);
        \draw [thin] (LL) --(RL);
        \draw [name path=lineLLIY,ultra thin,dashed] (LL) --(IY);
        \draw [name path=lineXRL,ultra thin,dashed]  (X)--(RL);
        \draw [name path=lineS1LL,ultra thin,dashed]  (S1)--(LL);
        \draw [name path=lineS2LL,ultra thin,dashed]  (S2)--(LL);
        \draw [name path=lineS3LL,ultra thin,dashed]  (S3)--(LL);
        \path[name intersections= 
        {of=lineLLIY and lineXRL,by={M0}}];
        \draw[thick](M0)--(IY);
        \path[name intersections= 
        {of=lineS1LL and lineXRL,by={M1}}];
        \path[name intersections= 
        {of=lineS2LL and lineXRL,by={M2}}];
        \path[name intersections= 
        {of=lineS3LL and lineXRL,by={M3}}];

        \path[name intersections= 
        {of=lineRLLS1 and lineLLIY,by={LM1}}];
        \path[name intersections= 
        {of=lineRLLS1 and lineS1LL,by={SM1}}];
        
        \path[name intersections= 
        {of=lineRLLS2 and lineLLIY,by={LM2}}];

        \path[name intersections= 
        {of=lineRLLS3 and lineLLIY,by={LM3}}];
         \draw[red, thick](S1)--(S2);

        \fill [orange, opacity=10/30](S1) -- (S2) -- (M2) -- (M1) -- cycle;
        \tkzDrawPoints(S1,S2,M2,M1);
        \draw[green, thick](M1)--(M2);
        \draw[thick](M1)--(S1);
        \draw[thick](M2)--(S2);
        \tkzLabelPoint[right](S1){\small $\Fc^2(\Min_{R,2})$};
        \node[] at (1.13*\PHISTAR,1.0*\PHISTAR){ \Large $\xrightarrow{\Fc}$};

        \node[] at (-1.15*\PHISTAR,1.0*\PHISTAR){ \Large $\xrightarrow{\Fc}$};
        \tkzLabelPoint[below right](LL){$-\Phistar$};
        \tkzLabelPoint[below left](RL){$+\Phistar$};
    \end{tikzpicture}
\end{center}
\end{minipage}
\begin{minipage}{0.6\textwidth}
\begin{center}
\scalebox{3.0}{\input{figs/CellFigure_in2.pgf}}
\end{center}
\end{minipage}
\caption{$\Fc$ action on $\mathcal{D}^{\text{\upshape out}}_{2}$ and its images. $\mathcal{D}^{\text{\upshape out}}_{2}$ is diffeomorphic to $\mathcal{D}^{\text{\upshape in}}_{-2}$ by $\Fc^4$. The same colored boundary are diffeomorphically mapped by $\Fc$ in each step.}\label{fig:D2celDynamics}
\end{figure}

\begin{figure}[h!]
\begin{minipage}{.75\textwidth}
\begin{center}
    \begin{tikzpicture}[xscale=.61,yscale=.61]
        \pgfmathsetmacro{\PHISTAR}{2.6}
        \pgfmathsetmacro{\LARGEPHISTAR}{0.7}
        \pgfmathsetmacro{\WIDTH}{11}
        \pgfmathsetmacro{\LENGTH}{11}
        \tkzDefPoint(-\LENGTH,0.5*\WIDTH){LLM};
        \tkzDefPoint(-\LENGTH,0){LLL};
        \tkzDefPoint(-\LENGTH+\PHISTAR,0){LL};
        \tkzDefPoint(-\LENGTH+\PHISTAR,\WIDTH){LU};
        \tkzDefPoint(-\LENGTH,\WIDTH){LLU};
        \tkzDefPoint(\LENGTH-\PHISTAR,0){RL};
        \tkzDefPoint(\LENGTH-\PHISTAR,\WIDTH){RU};
        \tkzDefPoint(\LENGTH,0){RRL};
        \tkzDefPoint(\LENGTH,\WIDTH){RRU};
        \tkzDefPoint(-\LENGTH+\PHISTAR,\PHISTAR){X};
        \tkzDefPoint(-\LENGTH+\PHISTAR,\WIDTH-\PHISTAR){IX};
        \tkzDefPoint(\LENGTH-\PHISTAR,\PHISTAR){IY};
        \tkzDefPoint(\LENGTH-\PHISTAR,\WIDTH-\PHISTAR){Y};
        \tkzLabelPoint[below right](LL){$(\phistar,0)$};
        \tkzLabelPoint[below left](RL){$(2\pi-\phistar,0)$};
        \tkzLabelPoint[left](IY){$(2\pi-\phistar,\phistar)$};
        \tkzLabelPoint[right](IX){$(\phistar,\pi-\phistar)$};

        \draw [ultra thin, dashed] (LL) --(IX);
        \draw [ultra thin, dashed] (LU) --(RU);
        \draw [ultra thin, dashed] (IY) --(RU);
        \draw [ultra thin, dashed] (LL) --(RL);
        \fill [red, opacity=6/30](LU) -- (IY) -- (RL) -- (IX) -- cycle;

        \tkzDefPoint(\LENGTH-\PHISTAR,\PHISTAR-\LARGEPHISTAR){R1};
        \tkzDefPoint(\LENGTH-\PHISTAR,\PHISTAR-2*\LARGEPHISTAR){R2};
        \tkzDefPoint(-\LENGTH+\PHISTAR,\WIDTH-\PHISTAR+\LARGEPHISTAR){L1};
        \tkzDefPoint(-\LENGTH+\PHISTAR,\WIDTH-\PHISTAR+2*\LARGEPHISTAR){L2};
        
        \tkzDrawPoints(L1,L2,LU,IX,R1,R2,IY,RL);
        \draw [red,thick,-triangle 45](IX)--(RL);
        \draw [yellow,thick,-triangle 45](RL)--(R2);
        \draw [orange,thick,-triangle 45](R2)--(R1);
        \draw [blue,thick,-triangle 45](R1)--(IY);
        \draw [Brown,thick,-triangle 45](IY)--(LU);
        \draw [green,thick,-triangle 45](LU)--(L2);
        \draw [teal,thick,-triangle 45](L2)--(L1);
        \draw [Purple,thick,-triangle 45](L1)--(IX);
        \tkzLabelPoint[left](R1){$(2\pi-\phistar,\phistar-\Phistar)$};
        \tkzLabelPoint[left](R2){$(2\pi-\phistar,\phistar-2\Phistar)$};
        \tkzLabelPoint[right](L1){$(\phistar,\pi-\phistar+\Phistar)$};
        \tkzLabelPoint[right](L2){$(\phistar,\pi-\phistar+2\Phistar)$};
        \tkzLabelPoint[above right](\LENGTH-\PHISTAR,0.5*\WIDTH){\Large $\xrightarrow{\Fc}$};
        \node[] at (0,0.5*\WIDTH)  {\Large $\Mrout$};
    \end{tikzpicture}
    \end{center}
\end{minipage}
\begin{minipage}{.2\textwidth}
\begin{center}
    \begin{tikzpicture}[scale=.24]
        \pgfmathsetmacro{\PHISTAR}{3};
        \pgfmathsetmacro{\SMALLPHISTAR}{12};
        \pgfmathsetmacro{\WIDTH}{36};
        \pgfmathsetmacro{\LENGTH}{18};
        \tkzDefPoint(-\LENGTH,0.5*\WIDTH){LLM};
        \tkzDefPoint(-\LENGTH,0){LLL};
        \tkzDefPoint(-\PHISTAR,0){LL};
        \tkzDefPoint(-\PHISTAR,\WIDTH){LU};
        \tkzDefPoint(-\LENGTH,\WIDTH){LLU};
        \tkzDefPoint(\PHISTAR,0){RL};
        \tkzDefPoint(\PHISTAR,\WIDTH){RU};
        \tkzDefPoint(\LENGTH,0){RRL};
        \tkzDefPoint(\LENGTH,\WIDTH){RRU};
        \tkzDefPoint(-\PHISTAR,\PHISTAR){X};
        \tkzDefPoint(-\PHISTAR,\WIDTH-\PHISTAR){IX};
        \tkzDefPoint(\PHISTAR,\WIDTH-\PHISTAR){Y};
        \tkzDefPoint(\PHISTAR,\PHISTAR){IY};
        \draw [ultra thin,dashed] (IX) --(LU);
        \draw [ultra thin,dashed] (RL) --(IY);
        \draw [ultra thin,dashed] (LL) --(RL);
        \draw [ultra thin,dashed] (LU) --(RU);
        \fill [blue, opacity=6/30](LL) -- (IY) -- (RU) -- (IX) -- cycle;
        \tkzDefPoint(-\PHISTAR,-\PHISTAR+\SMALLPHISTAR){L1};
        \tkzDefPoint(-\PHISTAR,\PHISTAR){L2};
        \tkzDefPoint(\PHISTAR,\WIDTH+\PHISTAR-\SMALLPHISTAR){R1};
        \tkzDefPoint(\PHISTAR,\WIDTH-\PHISTAR){R2};
        \tkzDrawPoints(L1,L2,LL,IY,R1,R2,RU,IX);
        
        \tkzLabelPoint[below](LL){$(-\Phistar,0)$};
        \tkzLabelPoint[below](RL){$(\Phistar,0)$};
        \tkzLabelPoint[above](IY){$(\Phistar,\Phistar)$};
        \draw [red,very thick,-triangle 45](IX)--(L1);
        \draw [yellow,very thick,-triangle 45](L1)--(L2);
        \draw [orange,very thick,-triangle 45](L2)--(LL);
        \draw [blue,very thick,-triangle 45](LL)--(IY);
        \draw [Brown,thick,-triangle 45](IY)--(R1);
        \draw [green,thick,-triangle 45](R1)--(R2);
        \draw [teal,thick,-triangle 45](R2)--(RU);
        \draw [Purple,thick,-triangle 45](RU)--(IX);
        \tkzLabelPoint[below right](L1){$(-\Phistar,\phistar-\Phistar)$};
        \tkzLabelPoint[above](L2){$(-\Phistar,\Phistar)$};
        \tkzLabelPoint[left](R1){$(\Phistar,\pi-\phistar+\Phistar)$};
        \tkzLabelPoint[below](R2){$(\Phistar,\pi-\Phistar)$};
        \tkzLabelPoint[above](RU){$(\Phistar,\pi)$};
        \tkzLabelPoint[above](IX){$(-\Phistar,\pi-\Phistar)$};
        \node[] at (0,0.5*\WIDTH)  {\Large $\MRin$};
    \end{tikzpicture}
\end{center}
\end{minipage}\caption{The $\Fc$ on $\Mrout$ has an continuous extension to $\partial\Mrout$ being a homeomephism from $\partial\Mrout$ to $\partial\MRin$ with the same color edge identification and orientation given in marked arrow.}\label{fig:F_actMroutMRin}
\end{figure}
\begin{figure}[h!]
\begin{minipage}{.2\textwidth}
\begin{center}
    \begin{tikzpicture}[scale=.24]
        \pgfmathsetmacro{\PHISTAR}{3};
        \pgfmathsetmacro{\SMALLPHISTAR}{12};
        \pgfmathsetmacro{\WIDTH}{36};
        \pgfmathsetmacro{\LENGTH}{18};
        \tkzDefPoint(-\LENGTH,0.5*\WIDTH){LLM};
        \tkzDefPoint(-\LENGTH,0){LLL};
        \tkzDefPoint(-\PHISTAR,0){LL};
        \tkzDefPoint(-\PHISTAR,\WIDTH){LU};
        \tkzDefPoint(-\LENGTH,\WIDTH){LLU};
        \tkzDefPoint(\PHISTAR,0){RL};
        \tkzDefPoint(\PHISTAR,\WIDTH){RU};
        \tkzDefPoint(\LENGTH,0){RRL};
        \tkzDefPoint(\LENGTH,\WIDTH){RRU};
        \tkzDefPoint(-\PHISTAR,\PHISTAR){X};
        \tkzDefPoint(-\PHISTAR,\WIDTH-\PHISTAR){IX};
        \tkzDefPoint(\PHISTAR,\WIDTH-\PHISTAR){Y};
        \tkzDefPoint(\PHISTAR,\PHISTAR){IY};
        \draw [ultra thin,dashed] (LL) --(X);
        \draw [ultra thin,dashed] (LU) --(RU);
        \draw [ultra thin,dashed] (Y) --(RU);
        \draw [ultra thin,dashed] (LL) --(RL);
        \fill [red, opacity=6/30](LU) -- (X) -- (RL) -- (Y) -- cycle;
        
        \tkzLabelPoint[below](LL){\small $(-\Phistar,0)$};
        \tkzLabelPoint[below](RL){\small $(\Phistar,0)$};
        \tkzDefPoint(\PHISTAR,\SMALLPHISTAR-\PHISTAR){R1};
        \tkzDefPoint(-\PHISTAR,\WIDTH-\SMALLPHISTAR+\PHISTAR){L1};
        \tkzDefPoint(\PHISTAR,\PHISTAR){R2};
        \tkzDefPoint(-\PHISTAR,\WIDTH-\PHISTAR){L2};

        \tkzDrawPoints(LU,X,RL,Y,L1,R1,L2,R2);
        \draw [blue,thick,-triangle 45](X)--(RL);
        \draw [orange,thick,-triangle 45](RL)--(R2);
        \draw [yellow,thick,-triangle 45](R2)--(R1);
        \draw [red,thick,-triangle 45](R1)--(Y);
        \draw [Purple,thick,-triangle 45](Y)--(LU);
        \draw [teal,thick,-triangle 45](LU)--(L2);
        \draw [green,thick,-triangle 45](L2)--(L1);
        \draw [Brown,thick,-triangle 45](L1)--(X);

        \tkzLabelPoint[right](\PHISTAR,0.5*\WIDTH){\Large $\xrightarrow{\Fc}$};
        \tkzLabelPoint[above left](R1){$(\Phistar,\phistar-\Phistar)$};
        \tkzLabelPoint[above](R2){$(\Phistar,\Phistar)$};
        \tkzLabelPoint[above](Y){\small $(\Phistar,\pi-\Phistar)$};
        \tkzLabelPoint[above](X){$(-\Phistar,\Phistar)$};
        \tkzLabelPoint[above](LU){\small $(-\Phistar,\pi)$};
        \tkzLabelPoint[below right](L2){\small $(-\Phistar,\pi-\Phistar)$};
        \tkzLabelPoint[below right](L1){\small $(-\Phistar,\pi-\phistar+\Phistar)$};
        \node[] at (0,0.5*\WIDTH)  {\Large $\MRout$};
    \end{tikzpicture}
\end{center} 
\end{minipage}\hfill
\begin{minipage}{.75\textwidth}
\begin{center}
    \begin{tikzpicture}[xscale=.63,yscale=.63]
        \pgfmathsetmacro{\PHISTAR}{2.6}
        \pgfmathsetmacro{\LARGEPHISTAR}{0.7}
        \pgfmathsetmacro{\WIDTH}{11}
        \pgfmathsetmacro{\LENGTH}{11}
        \tkzDefPoint(-\LENGTH,0.5*\WIDTH){LLM};
        \tkzDefPoint(-\LENGTH,0){LLL};
        \tkzDefPoint(-\LENGTH+\PHISTAR,0){LL};
        \tkzDefPoint(-\LENGTH+\PHISTAR,\WIDTH){LU};
        \tkzDefPoint(-\LENGTH,\WIDTH){LLU};
        \tkzDefPoint(\LENGTH-\PHISTAR,0){RL};
        \tkzDefPoint(\LENGTH-\PHISTAR,\WIDTH){RU};
        \tkzDefPoint(\LENGTH,0){RRL};
        \tkzDefPoint(\LENGTH,\WIDTH){RRU};
        \tkzDefPoint(-\LENGTH+\PHISTAR,\PHISTAR){X};
        \tkzDefPoint(-\LENGTH+\PHISTAR,\WIDTH-\PHISTAR){IX};
        \tkzDefPoint(\LENGTH-\PHISTAR,\PHISTAR){IY};
        \tkzDefPoint(\LENGTH-\PHISTAR,\WIDTH-\PHISTAR){Y};
        \tkzLabelPoint[below right](LL){$(\phistar,0)$};
        \tkzLabelPoint[below left](RL){$(2\pi-\phistar,0)$};
        \tkzLabelPoint[below left](Y){$(2\pi-\phistar,\pi-\phistar)$};
        \tkzLabelPoint[above right](X){$(\phistar,\phistar)$};

        \draw [ultra thin,dashed] (LU) --(X);
        \draw [ultra thin,dashed] (LU) --(RU);
        \draw [ultra thin,dashed] (RL) --(Y);
        \draw [ultra thin,dashed] (LL) --(RL);
        \fill [blue, opacity=6/30](LL) -- (Y) -- (RU) -- (X) -- cycle;
        \tkzDefPoint(-\LENGTH+\PHISTAR,\PHISTAR-\LARGEPHISTAR){L1};
        \tkzDefPoint(-\LENGTH+\PHISTAR,\PHISTAR-2*\LARGEPHISTAR){L2};
        \tkzDefPoint(\LENGTH-\PHISTAR,\WIDTH-\PHISTAR+\LARGEPHISTAR){R1};
        \tkzDefPoint(\LENGTH-\PHISTAR,\WIDTH-\PHISTAR+2*\LARGEPHISTAR){R2};
        \tkzDrawPoints(LL,Y,RU,X,L1,L2,R1,R2);
        
        \draw [blue,thick,-triangle 45](X)--(L1);
        \draw [orange,thick,-triangle 45](L1)--(L2);
        \draw [yellow,thick,-triangle 45](L2)--(LL);
        \draw [red,thick,-triangle 45](LL)--(Y);
        \draw [Purple,thick,-triangle 45](Y)--(R1);
        \draw [teal,thick,-triangle 45](R1)--(R2);
        \draw [green,thick,-triangle 45](R2)--(RU);
        \draw [Brown,thick,-triangle 45](RU)--(X);

        \tkzLabelPoint[right](L1){\small $(\phistar,\phistar-\Phistar)$};
        \tkzLabelPoint[right](L2){\small $(\phistar,\phistar-2\Phistar)$};
        \tkzLabelPoint[left](R2){\small $(2\pi-\phistar,\pi-\phistar+2\Phistar)$};
        \tkzLabelPoint[left](R1){\small $(2\pi-\phistar,\pi-\phistar+\Phistar)$};
        \tkzLabelPoint[above left](RU){\small $(2\pi-\phistar,\pi-\phistar)$};
        \node[] at (0,0.5*\WIDTH)  {\Large $\Mrin$};
    \end{tikzpicture}
    \end{center}
\end{minipage}
\caption{The $\Fc$ on $\MRout$ has an continuous extension to $\partial\MRout$ being a homeomephism from $\partial\MRout$ to $\partial\Mrin$ with the same color edge identification and orientation given in marked arrow.}\label{fig:F_actMRoutMrin}
\end{figure}

\NEWPAGE


\subsection[Alignment, transverse and regularity properties for singularity curves]{Alignment, transverse and regularity properties for singularity curves}

\begin{theorem}[Alignment of singularity curves on sections]\label{Theorem:SingularityAlignmentOnSection}
For a connected smooth curve segment $\mathcal{S}\subset(\mathcal{S}_{-t}\smallsetminus\mathcal{S}_{-(t-1)})$ with $t\ge1$, suppose $x$ is on $\mathcal{S}$ , then the following hold. 
\begin{enumerate}
    \item\label{item01:SingularityAlignmentOnSection} If $x\in(\mathcal{S}\cap\hat{M})$ with \eqref{def:Mhat}: $\hat{M}=\Fc^{-1}(\Mrout\smallsetminus\Mrin)\sqcup(\Mrout\cap\Mrin)$ given in and with the tangential vector $(d\phi,d\theta)\nfd dx\in T_x\mathcal{S}$, then $dx\in C_x$, $\frac{d\theta}{d\phi}\in(0,1)$, that is, the positive half-quadrant given by \cref{thm:UniformExpansionOnMhat}.
    \item\label{item02:SingularityAlignmentOnSection} If $x\in(\mathcal{S}\cap\Mrout)$ with and with the tangential vector $(d\phi,d\theta)\nfd dx\in T_x\mathcal{S}$, then $dx\in \tilde{C}^u_x$ with $\tilde{C}^u_x$ given in \cref{thm:MroutreturnUniformExpansion}, that is, $\frac{d\theta}{d\phi}$ satisfies 
    \begin{align*}
        \frac{d\theta}{d\phi}&\in\big(0,1\big)\text{, if } x\in \Mrout\smallsetminus\Nout,\\
        \frac{d\theta}{d\phi}&\in\big(0,\frac{1}{3}\big) \text{, if } x\in\Nout
    \end{align*}
    \item\label{item03:SingularityAlignmentOnSection} If $x\in(\mathcal{S}\cap\hat{M}_1)$ with \eqref{def:symetricsectionMhatPlusMinus}: $\hat{M}_1=\Fc(\Mrin\smallsetminus\Mrout)\sqcup(\Mrin\cap\Mrout)$ and with the tangential vector $(d\phi,d\theta)\nfd dx\in T_x\mathcal{S}$ and if $\Fc^{-1}(\mathcal{S})$ is not a part of the \emph{exceptional part of $\mathcal{AS}^{\text{\upshape in}}_{-1}$} described in \cref{remark:ArchSingualrityCurveAlmostBelowNegativeHalf} near the endpoints $(\phistar,\phistar-0.5\Phistar)$ or $(2\pi-\phistar,\pi-(\phistar-0.5\Phistar))$ , then $\frac{d\theta}{d\phi}\in(0,\infty)$, that is $dx$ is in the interior of the positive quadrant cone $\mathcal{C}(x)$ given in \eqref{eq:PositiveNegativeConeOnLemonSets} (or the $\mathcal{Q}_x(\mathrm{I,III})$ in \cref{corollary:InvQuadrant}).

\end{enumerate}
\end{theorem}
\begin{proof}
    \textbf{Proof of \eqref{item01:SingularityAlignmentOnSection}}: using \cref{notation:symplecticformmeasurenbhd} for singularity curves, we note that for any $x\in\mathcal{S}_{-t}\smallsetminus\mathcal{S}_{-(t-1)}=\mathcal{S}^{-}_{t}$ with $t\ge1$ has a trajectory with the $t$-th step in the past starting from a corner. 
    \[\underbracket{\Fc^{-t}(x)}_{\mathclap{\text{corner: }\in\partial M}},\Fc^{-t+1}(x),\cdots,x.\] In the following, we analyze the two cases for how this orbit starts from a corner.
    
    \textbf{Case 1}: If \begin{equation}\label{eq:startingAtMR}
        \underbracket{\Fc^{-t+1}(x),\cdots,-\Fc^{-(t+j)}(x)}_{\in M_R}
    \end{equation} $j\ge1$, are all collisions on $\Gamma_R$, then $\Fc^{-t+j}(\mathcal{S})$ is a curve segment of $\mathcal{AS}^{\text{\upshape in}}_{-j}$. For $x_{n,j}\dfn\Fc^{-t+j}(x)$, $dx_{t,j}\dfn D\Fc^{-t+j}(dx)$ is the tangential vector of $\mathcal{AS}^{\text{\upshape in}}_{-j}$ at $x_{t,j}$. By \cref{proposition:Archsingularitycurvesproperties,proposition:Archsingularitycurvesproperties2}, $dx_{t,j}=(d\phi_{t,j},d\theta_{t,j})$ satisfies: $\frac{d\theta_{t,j}}{\phi_{t,j}}\in [-\infty,-0.495)\cup(1,\infty]$

    Then for $x_{t,K}=\Fc^{-t+K}(x)\in\hat{M}$ with $K=\min{\Big\{k\ge j\bigm|\Fc^{-t+k}(x)=\Fc^{k-j}(x_{t,j})\in\hat{M}\Big\}}$. Note that since $\mathcal{AS}^{\text{\upshape in}}_{-j}\subset\Nin$, \cref{Prop:Udisjoint,contraction_region}: $K-j\ge2$, $x_{t,j},\Fc(x_{t,j}),\cdots,\Fc^{K-j}(x_{t,j})$ are collisions on the same boundary $\Gamma_r$, $(d\phi_{t,K},d\theta_{t,K})\nfd dx_{t,K}=D\Fc^{K-j}_{x_{n,j}}(dx_{t,j})$ with \[D\Fc^{K-j}_{x_{t,j}}=\begin{pmatrix}
        1 & 2(K-j)\\ 0 & 1
    \end{pmatrix}\] so that $\frac{d\theta_{t,K}}{d\phi_{t,K}}=\frac{d\theta_{t,j}}{\phi_{t,j}+2(K-j)d\theta_{t,j}}=\frac{1}{2(K-j)+\frac{1}{(d\theta_{t,j})/(d\phi_{t,j})}}$. 
    
    Then since $(d\theta_{t,j})/(d\phi_{t,j})\in [-\infty,-0.495)\cup(1,\infty]$, $(d\phi_{t,j})/(d\theta_{t,j})\in (-2.03,1)$, $2(K-j)+(d\phi_{t,j})/(d\theta_{t,j})>1.97$. Therefore, $\frac{d\theta_{t,K}}{d\phi_{t,K}}\in (0,\frac{1}{1.97})\subset (0,1)$. The image of the tangential vector at $x_{t,j}\in\mathcal{AS}^{\text{\upshape in}}_{-j}$ in its first arrival/return at $\hat{M}$ is already in the interior positive half-quadrant cone at $\hat{M}$. Then by \cref{thm:UniformExpansionOnMhat}: the $\hat{M}$ return map preserves the cone. Then, in all subsequent arrivals / returns on $\hat{M}$, the images of the tangent vector $dx$' at these points will also be in the positive half-quadrant interior cone at $\hat{M}$. Hence $x\in\mathcal{S}\cap\hat{M}$ has a tangential slope in the cone $C_x$.
    
    \textbf{Case 2}: If \begin{equation}\label{eq:startingAtMr}
        \underbracket{\Fc^{-t+1}(x),\cdots,\Fc^{-t+j}(x)}_{\in M_r}\in \Mrout 
    \end{equation} $j\ge1$, are all collisions on $\Gamma_r$, then $\Fc^{-t+1}(\mathcal{S})$ is a straight line segment of $\mathcal{S}^{-}_{1}$ that is part of the boundary of $\Mrin$: $\partial\Mrin$ with a constant slope $1/2$ (see \cref{fig:Mr}). 
    
    We suppose $K=\min{\Big\{k\ge 1\bigm|\Fc^{-t+k}(x)\in\hat{M}\Big\}}$. Let $dx\in T_x(\Fc^{-t+1}(\mathcal{S}))$, then $D\Fc_x^{K-1}(dx)$ is a tangential vector at a point of $\hat{M}$. By \cref{caseA,caseB,caseC,Prop:Udisjoint,contraction_region} (also see \cref{fig:MrN}), $D\Fc_x^{K-1}(dx)$ is an interior point of the positive half quadrant at a base point that is $x$'s first arrival on $\hat{M}$. By \cref{thm:UniformExpansionOnMhat1} in all subsequent arrivals / returns on $\hat{M}$, the images of the tangent vector $dx$' at these points will also be in the positive half-quadrant interior cone at $\hat{M}$. So $x\in\mathcal{S}\cap\hat{M}$ has a tangential vector inside the cone $C_x$.

    \textbf{Proof of \eqref{item02:SingularityAlignmentOnSection},\eqref{item03:SingularityAlignmentOnSection}}: the proofs are the same as for \eqref{item01:SingularityAlignmentOnSection} by analyzing the start of the trajectory from corner for two cases: \eqref{eq:startingAtMR}, \eqref{eq:startingAtMr} and invoking \cref{caseA,caseB,caseC,Prop:Udisjoint,contraction_region,thm:UniformExpansionOnMhat1,thm:MroutreturnUniformExpansion} (also see \cref{fig:MrN,remark:MoreUdisjoint})
\end{proof}
\begin{remark}[The locus of the exceptional curve segment's second-time arrival on $\hat{M}_1$]\label{remark:AllotherNonExceptionSingularityCurves} Note that by our chosen $R_{\text{\upshape HF}}(r,\phistar)$ as large as given in \cref{def:AcomputationForRHF}, (see \cref{remark:GaranteeLargeTheta}) orbits starting from $\Nin$ ( $\Nout$ ) cannot have 2 consecutive returns to $\Mrin$ ($\Mrout$) as the returns to $\Nin$ ( $\Nout$ ). Hence, the second arrival on $\hat{M}_1$ from the exceptional curve segment (\cref{remark:ArchSingualrityCurveAlmostBelowNegativeHalf}) of $\mathcal{AS}^{\text{\upshape in}}_{-1}$ near $(\phistar,\phistar-0.5\Phistar)$ or $(2\pi-\phistar,\pi-(\phistar-0.5\Phistar))$ will be a singularity curve with tangential direction vectors in the interior of the positive quadrant cone $\mathcal{C}(x)$.
\end{remark}
Then this immediately gives the L2 conditions in \cref{thm:lemonLET}.
\begin{corollary}[L2 condition in \cref{thm:lemonLET}]\label{corollary:L2AlignmentCondition} Let $k\ge1$ and $\Scal$ be a curve segment in $M$.\hfill

    \begin{enumerate}
        \item For the sufficient point $x\in\hat{\mathfrak{U}}_-$ with quadruple $(l,N,U,\mathcal{C})$ given in \cref{theorem:sufficentptsonMhatminus}, the following alignment property holds.

        If $y\in \mathcal{S}\cap \Fc^{-N}(U)$ where $\mathcal{S}\subset\mathcal{S}_{-k}\smallsetminus \mathcal{S}_{-(k-1)}$, then the tangential vector $\mathcal{S}$ at $y$ is contained in the cone $\mathcal{C}$.

        If $y\in \mathcal{S}\cap U$ where $\mathcal{S}\subset\mathcal{S}_{k}\smallsetminus \mathcal{S}_{k-1}$, then the tangential vector of $\mathcal{S}$ at $y$ is contained in $\mathcal{C}'$ which is the complement of the cone $\mathcal{C}$.
        \item For the sufficient point $x\in\hat{\mathfrak{U}}_+$ with quadruple $(l,N,U,\mathcal{C})$ given in \cref{theorem:sufficentptsonMhatplus}, the following alignment property holds.

        If $y\in \mathcal{S}\cap \Fc^{-N}(U)$ where $\mathcal{S}\subset\mathcal{S}_{-k}\smallsetminus \mathcal{S}_{-(k-1)}$, then the tangential vector of $\mathcal{S}$ at $y$ is contained in  $\mathcal{C}$.

        If $y\in \mathcal{S}\cap U$ where $\mathcal{S}\subset\mathcal{S}_{k}\smallsetminus \mathcal{S}_{k-1}$, then the tangential vector of $\mathcal{S}$ at $y$ is contained in $\mathcal{C}'$.
    \end{enumerate}
\end{corollary}
\begin{proof} We prove the $4$ cases in the above statement. And notice the last $2$ cases will be symmetric to the first $2$ cases since $DI$ and its inverse $DI^{-1}:(d\phi,d\theta)\rightarrow(d\phi,-d\theta)$ is an diffeomphism between the  positive cone $\mathcal{C}$ and its complement $\mathcal{C}'$ (see \eqref{eq:PositiveNegativeConeOnLemonSets}).

    \textbf{Case 1}: suppose $y\in\Scal\subset(\Scal_{-k}\smallsetminus\Scal_{-(k-1)})$, for the sufficient point $x\in\hat{\mathfrak{U}}_-$ with quadruple $(l,N,U,\mathcal{C})$ given in \cref{theorem:sufficentptsonMhatminus} and $y\in\Scal\cap\Fc^{-N}(U)$. \cref{Theorem:SingularityAlignmentOnSection}\eqref{item01:SingularityAlignmentOnSection} gives the conclusion that the tangential vector of $\Scal$ at $y$ is contained in $C_y$ the half positive quadrant thus in $\mathcal{C}(y)$ the positive quadrant.

    \textbf{Case 2}: alternatively, suppose $y\in\Scal\subset(\Scal_{-k}\smallsetminus\Scal_{-(k-1)})$, for the sufficient point $x\in\hat{\mathfrak{U}}_+$ with quadruple $(l,N,U,\mathcal{C})$ given in \cref{theorem:sufficentptsonMhatplus} and $y\in\Scal\cap\Fc^{-N}(U)$. Since \cref{theorem:sufficentptsonMhatplus}\eqref{item02sufficentptsonMhatplus} gives $\Fc^{-N}(U)\subset\hat{\mathfrak{U}}_+$, by \cref{def:UIntRegionofD}: the definition for $\hat{\mathfrak{U}}_+$, the exceptional curve segment per \cref{remark:ArchSingualrityCurveAlmostBelowNegativeHalf} is a part of $\mathcal{AS}^{in}_{-1}$. Since $\Fc(\mathcal{AS}^{in}_{-1})\cap\hat{\mathfrak{U}}_+=\emptyset$ and \cref{remark:AllotherNonExceptionSingularityCurves}, the tangential vector of $\mathcal{S}$ at $y$ is contained in $\mathcal{C}$. 

    \textbf{Case 3}: suppose $y\in\Scal\subset(\Scal_{k}\smallsetminus\Scal_{k-1})$, $y\in\Scal\cap U$,  for the sufficient $x\in\hat{\mathfrak{U}}_-$. For the tangential vector $v_y$ of $\Scal$ at $y$, $v_{I(y)}=DI(v_y)$ is a tangential vector of $I\Scal$ at $Iy$. Then by symmetry $Iy\in I\Scal\subset(\Scal_{-k}\smallsetminus\Scal_{-(k-1)})$ and $Iy\in IU\subset I\hat{\mathfrak{U}}_-=\hat{\mathfrak{U}}_+$. By the conclusion in the case $1$, $v_{I(y)}\in \mathcal{C}(I(y))$. Hence, $v_y=DI^{-1}(v_{I(y)})\in\mathcal{C}'(y)$.

    \textbf{Case 4}: suppose $y\in\Scal\subset(\Scal_{k}\smallsetminus\Scal_{k-1})$, $y\in\Scal\cap U$,  for the sufficient $x\in\hat{\mathfrak{U}}_+$. The reason same as for Case 3 using symmetry and the conclusion of Case 2 gives the tangential vector $v_y$ of $\Scal$ at $y$ is in $\mathcal{C}'(y)$.
\end{proof}

\begin{corollary}[Transverse intersection property]\label{Corollary:SingularityCurvesTransveralIntersection} For a connected smooth curve segment $\mathcal{S}\subset(\mathcal{S}_{-t}\smallsetminus\mathcal{S}_{-(t-1)})$ with $t\ge1$, suppose that $x$ is on $\mathcal{S}$ then the following hold.
\begin{enumerate}
    \item\label{item01SingularityCurvesTransveralIntersection} If $L$ is one of the two connected line components of $\partial\Mrout\smallsetminus\partial M_r=\Scal^+_1\cap M_r$, then $\mathcal{S}\cap L$ may contain at most one point.
    \item\label{item02SingularityCurvesTransveralIntersection} If $L$ is one of the two connected line components of $\partial\Mrout\cap\partial M_r$ (vertical line), then $\text{closure}(\mathcal{S})\cap L$ can contain at most one point which must be the endpoint of the curve segment $\mathcal{S}$.
    \item\label{item03SingularityCurvesTransveralIntersection} For $\forall n\ge1$, $\mathcal{S}\cap\mathcal{AS}^{\text{\upshape out}}_{n}$ can contain at most one point. 
    \textbf{$\text{closure}(\mathcal{S})\cap\text{closure}(\mathcal{AS}^{\text{\upshape out}}_{n})$} can contain at most one point.
    \item\label{item04SingularityCurvesTransveralIntersection} Suppose $t\ge3$, if $L$ is one of the connected line components of $\partial\MRout\smallsetminus\partial M_R=\Scal^+_1\cap M_R$, then $\mathcal{S}\cap L$ can contain at most one point.
    \item\label{item05SingularityCurvesTransveralIntersection} Suppose $t\ge3$, if $VL$ is one of the connected line components of $\partial\MRout\cap\partial M_R$, then \textbf{closure($(\mathcal{S})\cap VL$)} can contain only one or two points that must be one or two of the endpoints of the curve segment $\mathcal{S}$. (no more than $2$ intersections)
\end{enumerate}
\end{corollary}
\begin{figure}[h!]
\begin{minipage}{.45\textwidth}
\begin{center}
\begin{tikzpicture}[xscale=0.45, yscale=0.45]
        \pgfmathsetmacro{\PHISTAR}{6};
        \pgfmathsetmacro{\WIDTH}{8};
        \tkzDefPoint(-\PHISTAR,0){LL};
        \tkzDefPoint(-\PHISTAR,\WIDTH){LU};
        \tkzDefPoint(\PHISTAR,0){RL};
        \tkzDefPoint(\PHISTAR,\WIDTH){RU};
        \tkzDefPoint(-\PHISTAR,\PHISTAR){X};
        \tkzDefPoint(\PHISTAR,\WIDTH-\PHISTAR){Y};
        \tkzDefPoint(-\PHISTAR,\WIDTH-\PHISTAR){IX};
        \tkzDefPoint(\PHISTAR,\PHISTAR){IY};
        \tkzDefPoint(\PHISTAR,0.5*\PHISTAR){S1};
        \tkzDefPoint(\PHISTAR,\PHISTAR/3){S2};
        \tkzDefPoint(\PHISTAR,\PHISTAR/4){S3};
        \tkzDefPoint(\PHISTAR,\PHISTAR/5){L4};
        \tkzDefPoint(-\PHISTAR,0.5*\PHISTAR){LS1};
        \tkzDefPoint(-\PHISTAR,\PHISTAR/3){LS2};
        \tkzDefPoint(-\PHISTAR,\PHISTAR/4){LS3};
        \draw [name path=lineRLLS1,ultra thin,dashed] (RL) --(LS1);
        \draw [name path=lineRLLS2,ultra thin,dashed] (RL) --(LS2);
        \draw [name path=lineRLLS3,ultra thin,dashed] (RL) --(LS3);
        \draw [thin] (LL) --(LU);
        \draw [thin] (LU) --(RU);
        \draw [thin] (RL) --(RU);
        \draw [thin] (LL) --(RL);
        \draw [name path=lineLLIY, ultra thin,dashed] (LL) --(IY);
        \draw [name path=lineXRL,ultra thin,dashed]  (X)--(RL);
        \draw [name path=lineS1LL,ultra thin,dashed]  (S1)--(LL);
        \draw [name path=lineS2LL,ultra thin,dashed]  (S2)--(LL);
        \draw [name path=lineS3LL,ultra thin,dashed]  (S3)--(LL);
        \path[name intersections= 
        {of=lineLLIY and lineXRL,by={M0}}];
        \draw[blue,thick](M0)--(X);
        \path[name intersections= 
        {of=lineS1LL and lineXRL,by={M1}}];
        \path[name intersections= 
        {of=lineS2LL and lineXRL,by={M2}}];
        \path[name intersections= 
        {of=lineS3LL and lineXRL,by={M3}}];

        \path[name intersections= 
        {of=lineRLLS1 and lineLLIY,by={LM1}}];
        \draw[Red,thick](LM1)--(LS1);

        \path[name intersections= 
        {of=lineRLLS2 and lineLLIY,by={LM2}}];
        \draw[Green,thick](LM2)--(LS2);

        \path[name intersections= 
        {of=lineRLLS3 and lineLLIY,by={LM3}}];
        \draw[thick](LM2)--(LM1);
        \draw[thick](LS1)--(LS2);

        \fill [orange, opacity=10/30](LS2) -- (LS1) -- (LM1) -- (LM2) -- cycle;
        \tkzDrawPoints(LS1,LS2,LM1,LM2,M0,IY,X);
        \tkzLabelPoint[left](LS1){\small $\Min_{R,2}$};

        \node[] at (-1.2*\PHISTAR,1.0*\PHISTAR){ \Large $\xrightarrow{\Fc}$};
        \tkzLabelPoint[below](LL){$(-\Phistar,0)$};
        \tkzLabelPoint[above right](X){$(-\Phistar,\Phistar)$};
        \tkzLabelPoint[above left](IY){$(\Phistar,\Phistar)$};
        \tkzLabelPoint[below left](RL){$(\Phistar,0)$};

        \begin{scope}
            \clip (LS1) -- (LL) -- (RL) --cycle;
            \draw[name path=invalidSingCurve1,violet, thick] ellipse[rotate=165, at=(LM1),x radius=0.83, y radius =0.3];
            \path[name intersections= {of=invalidSingCurve1 and lineLLIY,by={SingPoint1}}];
            \tkzDrawPoints(SingPoint1);
        \end{scope}
    \end{tikzpicture}
     \end{center}
\end{minipage}
\begin{minipage}{.45\textwidth}
\begin{center}
    \begin{tikzpicture}[xscale=0.45, yscale=0.45]
        \pgfmathsetmacro{\PHISTAR}{6};
        \pgfmathsetmacro{\WIDTH}{8};
        \tkzDefPoint(-\PHISTAR,0){LL};
        \tkzDefPoint(-\PHISTAR,\WIDTH){LU};
        \tkzDefPoint(\PHISTAR,0){RL};
        \tkzDefPoint(\PHISTAR,\WIDTH){RU};
        \tkzDefPoint(-\PHISTAR,\PHISTAR){X};
        \tkzDefPoint(\PHISTAR,\WIDTH-\PHISTAR){Y};
        \tkzDefPoint(-\PHISTAR,\WIDTH-\PHISTAR){IX};
        \tkzDefPoint(\PHISTAR,\PHISTAR){IY};
        \tkzDefPoint(\PHISTAR,0.5*\PHISTAR){S1};
        \tkzDefPoint(\PHISTAR,\PHISTAR/3){S2};
        \tkzDefPoint(\PHISTAR,\PHISTAR/4){S3};
        \tkzDefPoint(\PHISTAR,\PHISTAR/5){L4};
        \tkzDefPoint(-\PHISTAR,0.5*\PHISTAR){LS1};
        \tkzDefPoint(-\PHISTAR,\PHISTAR/3){LS2};
        \tkzDefPoint(-\PHISTAR,\PHISTAR/4){LS3};
        \draw [name path=lineRLLS1,ultra thin,dashed] (RL) --(LS1);
        \draw [name path=lineRLLS2,ultra thin,dashed] (RL) --(LS2);
        \draw [name path=lineRLLS3,ultra thin,dashed] (RL) --(LS3);

        \draw [thin] (LL) --(LU);
        \draw [thin] (LU) --(RU);
        \draw [thin] (RL) --(RU);
        \draw [thin] (LL) --(RL);
        \draw [name path=lineLLIY,ultra thin,dashed] (LL) --(IY);
        \draw [name path=lineXRL,ultra thin,dashed]  (X)--(RL);
        \draw [name path=lineS1LL,ultra thin,dashed]  (S1)--(LL);
        \draw [name path=lineS2LL,ultra thin,dashed]  (S2)--(LL);
        \draw [name path=lineS3LL,ultra thin,dashed]  (S3)--(LL);
        \path[name intersections= 
        {of=lineLLIY and lineXRL,by={M0}}];
        \draw[thick](M0)--(IY);
        \path[name intersections= 
        {of=lineS1LL and lineXRL,by={M1}}];
        \path[name intersections= 
        {of=lineS2LL and lineXRL,by={M2}}];
        \path[name intersections= 
        {of=lineS3LL and lineXRL,by={M3}}];

        \path[name intersections= 
        {of=lineRLLS1 and lineLLIY,by={LM1}}];
        \path[name intersections= 
        {of=lineRLLS1 and lineS1LL,by={SM1}}];
        
        \path[name intersections= 
        {of=lineRLLS2 and lineLLIY,by={LM2}}];

        \path[name intersections= 
        {of=lineRLLS3 and lineLLIY,by={LM3}}];
        \draw[blue, thick](IY)--(S1);
        \draw[red, thick](M0)--(M1);

        \fill [orange, opacity=10/30](LM1) -- (SM1) -- (M1) -- (M0) -- cycle;
        \tkzDrawPoints(LM1,SM1,M1,M0);
        \draw[green, thick](LM1)--(SM1);
        \draw[thick](LM1)--(M0);
        \draw[thick](SM1)--(M1);
        \node[] at (1.13*\PHISTAR,1.0*\PHISTAR){ \Large $\xrightarrow{\Fc}$};
        \node[] at (1.2*\PHISTAR,0.47*\PHISTAR){ \Large $VL_0$};

        \node[] at (-1.15*\PHISTAR,1.0*\PHISTAR){ \Large $\xrightarrow{\Fc}$};
        \tkzLabelPoint[below right](LL){$(-\Phistar,0)$};
        \tkzLabelPoint[below left](RL){$(\Phistar,0)$};
        \tkzLabelPoint[above right](X){$(-\Phistar,\Phistar)$};
        \tkzLabelPoint[above left](IY){$(\Phistar,\Phistar)$};
        
        \begin{scope}
            \clip (X) -- (RL) -- (LL) --cycle;
            \draw[name path=invalidSingCurve1, violet, thick] ellipse[rotate=155, at=(M1),x radius=0.76, y radius =0.25];

            \draw[name path=invalidSingCurve2,teal, thick] ellipse[rotate=155, at=(M0),x radius=0.86, y radius =0.2];

            \path[name intersections= {of=invalidSingCurve1 and lineS1LL,by={SingPoint1}}];

            \path[name intersections= {of=invalidSingCurve2 and lineLLIY,by={SingPoint2}}];
            \tkzDrawPoints(SingPoint1,SingPoint2);
        \end{scope}
        \begin{scope}
            \clip (X) -- (RL) -- (IY) --cycle;
            \draw[name path=invalidSingCurve3,brown, thick] ellipse[rotate=155, at=(M0),x radius=1.36, y radius =0.23];
            
            \path[name intersections= {of=invalidSingCurve3 and lineLLIY,by={SingPoint3}}];
            \tkzDrawPoints(SingPoint3);
        \end{scope}
        \tkzLabelSegment[above=3pt](X,M0){$L(1)$};
        \tkzLabelPoint[right=2pt](M0){$L(2)$};
        
        \tkzLabelPoint[right](M1){$L(3)$};
        \tkzLabelPoint[right](M2){$L(4)$};
        \tkzDrawPoints(M2,M3,X,Y,RL,LL);
    \end{tikzpicture}
\end{center}
\end{minipage}\hfill
\caption{For a curve segment $\Scal\subset\mathcal{S}_{t}\smallsetminus\mathcal{S}_{t-1}$, $t\ge3$. the curves $\Scal$, $\text{closure}(\Scal)$ cannot intersect one connected component $\mathcal{S}^{+}_{1}\cap M_R$ at two points. Otherwise the following could happen:\newline (a): $\Scal$ or its some $\Fc$ iteration preimage (in left figure) could intersect with $\mathcal{S}^{-}_{1}$ (segment that connects $(-\Phistar,0)$ to $(\Phistar,\Phistar)$).\newline (b): The closure for some $\Fc$ iteration preimage of $\Scal$ could intersect 
$\text{closure} (\mathcal{AS}^{\text{out}}_{j})$ at two points (the \cref{fig:D2celDynamics} a boundary of the $\mathcal{D}^{\text{\upshape out}}_{n}$).}\label{fig:ImpossibleSingularityIntersection}
\end{figure}
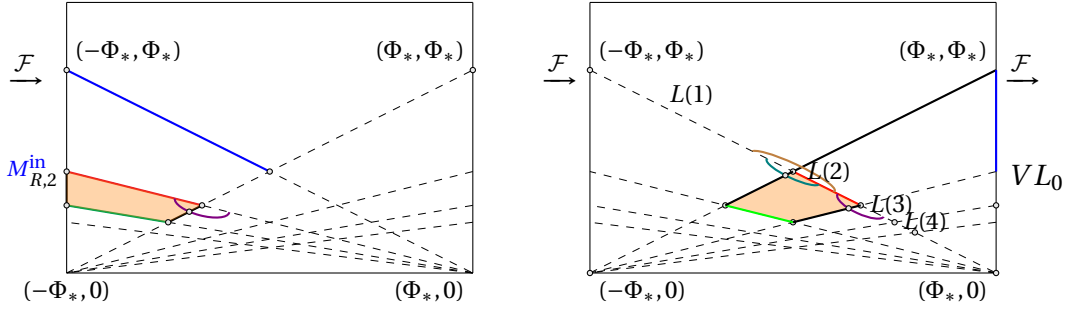
\begin{proof}
    \textbf{The proof of \eqref{item01SingularityCurvesTransveralIntersection}} follows \cref{Theorem:SingularityAlignmentOnSection}\eqref{item01:SingularityAlignmentOnSection}, \eqref{item02:SingularityAlignmentOnSection}, \eqref{item03:SingularityAlignmentOnSection}, \cref{caseA,caseB,caseC,contraction_region} that the singularity curve $\mathcal{S}\subset M_r$ with tangent vectors not being in the interior of positive quadrant $\mathcal{C}$ may only happen when $S\subset\Nin$ or $\mathcal{S}$ is the $\Fc$ image of the \emph{exceptional part of $\mathcal{AS}^{\text{\upshape in}}_{-1}$} described in \cref{remark:ArchSingualrityCurveAlmostBelowNegativeHalf}. Since both $\Nin$ and $\Fc(\Nin)$ have no intersection with a connected component $L\subset\partial\Mrout\smallsetminus\partial M_r$, $\mathcal{S}$ cannot have two intersections with $\partial\Mrout\smallsetminus\partial M_r$. Otherwise, by the Mean-Value Theorem, $\mathcal{S}$ could have a tangential vector with slope $-\frac{1}{2}$ that is not the interior of the cone $\mathcal{C}$. 

    $\mathcal{S}$ cannot have two intersections with a connected component $L\subset\partial\Mrout\cap\partial M_r$ (one vertical line segment). Otherwise, by the Mean-Value Theorem, $\mathcal{S}$ could have a tangential vector with slope $\infty$ (in $\phi\theta$ coordinate) that also is not the interior of the cone $\mathcal{C}$. 

    Example of components of singularity curves $\mathcal{S}_{-2}\smallsetminus\mathcal{S}_{-1}$: the two segments: $\Fc(\partial\Mrin\cap\Mrout)\subset\mathcal{S}_{-2}\smallsetminus\mathcal{S}_{-1}$ are two curves in $M_R$ that connect $(-\Phistar,\frac{\pi}{2}+\phistar-\Phistar)$ to $(\Phistar,\frac{\pi}{2}-\frac{\phistar}{2}+\Phistar)$ and $(-\Phistar,\frac{\pi}{2}+\frac{\phistar}{2}-\Phistar)$ to $(\Phistar,\frac{\pi}{2}-\phistar+\Phistar)$.
    
    \textbf{Proof for \eqref{item02SingularityCurvesTransveralIntersection}} also follows the same reason as \eqref{item01SingularityCurvesTransveralIntersection} that the tangential vector $(d\phi,d\theta)$ of the smooth curve $\mathcal{S}$ satisfies $0<\frac{d\theta}{d\phi}<1$. Hence, the closure of $\big(\mathcal{S}\big)$ can only intersect the vertical boundary $L=\partial\MRout\cap\partial M$ at one point.

    \textbf{Proof of \eqref{item03SingularityCurvesTransveralIntersection}}: suppose $\mathcal{S}$ and $\mathcal{AS}^{\text{\upshape out}}_{n}$ intersect at a point $v=(\phi_v,\theta_v)$.
    
    Since $\mathcal{AS}^{\text{\upshape out}}_{n}\subset \Nout$, the tangential vector $(d\phi,d\theta)$ at all points of $\mathcal{S}\subset\Nout$ has reciprocal slope $\frac{d\phi}{d\theta}>3$ by \cref{Theorem:SingularityAlignmentOnSection}\eqref{item02:SingularityAlignmentOnSection}. If we parameterize $\mathcal{S}$ as $(\phi_{\mathcal{S}}(\theta),\theta)$, then $\frac{d\phi_{\mathcal{S}}}{d\theta}>3$.
    
    On the other hand, if we parameterize the curve $\mathcal{AS}^{\text{\upshape out}}_n$ as $(\phi(\theta),\theta)$,  then \cref{proposition:Archsingularitycurvesproperties3}\eqref{item01Archsingularitycurvesproperties2}\eqref{item02Archsingularitycurvesproperties2} implies that the tangential vector $(d\phi,d\theta)$ at any point in $\mathcal{AS}^{\text{\upshape out}}_n$ has $\frac{d\phi}{d\theta}\in(-1,\frac{1}{0.495})$. 

    Then the difference in the tangential slopes mentioned above implies that the point $(\phi(\theta),\theta)$ on $\mathcal{AS}^{\text{\upshape out}}_n$ and $(\phi_{\mathcal{S}},\theta)$ on $\mathcal{S}$ with the same $\theta$  cannot have $\phi(\theta)=\phi_\mathcal{S}$ except $\theta=\theta_v$ and $\phi=\phi_\mathcal{S}=\phi_v$, that is, $\mathcal{S}$ cannot intersect $\mathcal{AS}^{\text{\upshape out}}_n$ twice. We can also suppose $\text{closure}(\Scal)\cap\text{closure}(\mathcal{AS}^{\text{\upshape out}}_n)$ at $v$ which is the endpoint, by the smoothness ($C^1$ is sufficient) of $\Scal$ and $\mathcal{AS}^{\text{\upshape out}}_n$, the difference in tangential slope will also imply that there cannot exist a second intersection.

    \textbf{Proof of }\eqref{item04SingularityCurvesTransveralIntersection}: In this case we can suppose $L$ is a connected component of the $\partial\MRout$ line with constant slope $-\frac{1}{2}$ that connects $(-\Phistar,\Phistar)$ and $(\Phistar,0)$ (see \cref{fig:D1celDynamics,fig:D2celDynamics}). And let $L'$ be the another connected component of the $\partial\MRout$ line with constant slope $-\frac{1}{2}$ which is the symmetric to $L$ by $I\circ J:(\Phi,\theta)\rightarrow(-\Phi,\pi-\theta)$ (see \cref{DefSymmetries}).
    
    Since $L'$, $L$ are contained the in the two different connected components of $M_R$ separated by a curve component of $\mathcal{S}_{-2}\smallsetminus\mathcal{S}_{-1}$ that connects $(-\Phistar,\frac{\pi}{2}+\phistar-\Phistar)$ to $(\Phistar,\frac{\pi}{2}-\frac{\phistar}{2}+\Phistar)$ mentioned in \eqref{item01SingularityCurvesTransveralIntersection}, if $\mathcal{S}\cap L=v$ and $\mathcal{S}'\cap L=v'$ then $\mathcal{S}$ could intersect with the component of $\mathcal{S}_{-2}\smallsetminus\mathcal{S}_{-1}$. This contradicts with  $\mathcal{S}\subset\mathcal{S}_{t}\smallsetminus\mathcal{S}_{t-1}$ of $t\ge3$. Therefore, $\mathcal{S}$ cannot have intersections with both $L$ and $L'$ by Jordan curve theorem.
    
    If $\mathcal{S}\subset(\mathcal{S}_{t}\smallsetminus\mathcal{S}_{t-1})$ is a curve segment with points having $1,\cdots,t-1$-th previous collisions on $\Int(\Gamma_R)$, i.e. a point on $\mathcal{S}$ has its previous trajectory starting at the corner as follows \[\underbracket{\Fc^{-t}(x)}_{\mathclap{\text{corner: }\in\partial M}},\underbracket{\Fc^{-t+1}(x),\cdots,x}_{\mathclap{\in M_R}},\]
    then $\mathcal{S}$ is a part of the straight line segment with slope $\frac{1}{2t}$ in \cref{fig:D1celDynamics,fig:D2celDynamics}. Of course, it cannot have more than $1$ intersection with the $-1/2$ slope line segment.

    Then, otherwise with some $1\le j<t$, $\mathcal{S}$ are points of collision with all its $1,\cdots,j-1$-th previous collisions on $\Int(\Gamma_R)$ and its $j$-th previous collision on $\Int(\Gamma_r)$ , i.e. $\Fc^{-j}(\mathcal{S})\subset\Mrout$. 

    Suppose $v\in\mathcal{S}\cap L$ with $\Fc^{-j}(v)\in\mathcal{AS}^{\text{\upshape out}}_{j}$. From \cref{fig:D1celDynamics}: the action of $\Fc$ on quadrilateral regions, we see that $v$ lies on the interior of the line segment connecting $\big(\Phistar-\frac{2}{j}\Phistar,\frac{1}{j}\Phistar\big)$ with $\big(\Phistar-\frac{2}{j+1}\Phistar,\frac{1}{j+1}\Phistar\big)$ in \cref{fig:F_actMroutMRin}.

    Now we denote the interior of the line segment that connects $\big(\Phistar-\frac{2}{j}\Phistar,\frac{1}{j}\Phistar\big)$ with $\big(\Phistar-\frac{2}{j+1}\Phistar,\frac{1}{j+1}\Phistar\big)$ by $L(j)$ for $j\ge1$. And denote the interior of the line segment that $\big(-\Phistar+\frac{2}{j}\Phistar,\pi-\frac{1}{j}\Phistar\big)$ with $\big(-\Phistar+\frac{2}{j+1}\Phistar,\pi-\frac{1}{j+1}\Phistar\big)$ by $L'(j)$ which is a segment on $L'$. $L(j)$ is symmetric to $L(j)$ by $I\circ J$.

    Hence, now in our assumption, $v\in L(j)$.
    
    Also note that $L\subset\bigcup_{j\ge1}\big(\text{closure}(L(j))\big)$, $L'\subset\bigcup_{j\ge1}(\text{closure}(L'(j)))$.
    
    Then suppose $v\ne v'\in\mathcal{S}\cap L$. Based on the position of $v'$ on $L$, there are the following three cases.

    If $v'\in\text{closure}(L(j'))$ with $j'<j$, then by the action of $\Fc$ on the quadrilaterals \cref{remark:FcActionOnQuadrilaterals,fig:D1celDynamics,fig:D2celDynamics,fig:ImpossibleSingularityIntersection} and by the Jordan curve theorem $(\mathcal{S}_{-(t-j'+1)}\smallsetminus\mathcal{S}_{-(t-j')})\supset\Fc^{-(j'-1)}(\mathcal{S})$ is a curve that intersects $\text{closure}(\mathcal{S}^-_1)\subset \mathcal{S}_{-1}$. It is a contradiction that some point is on $\mathcal{S}_{-1}$ but not on $\mathcal{S}_{-(t-j')}$ with $t-j'\ge 1$.
    
    If $v'\in\text{closure}(L(j'))$ with $j'>j$, then by the action of $\Fc$ on the quadrilaterals \cref{remark:FcActionOnQuadrilaterals,fig:D1celDynamics,fig:D2celDynamics,fig:ImpossibleSingularityIntersection} and by the Jordan curve theorem $(\mathcal{S}_{-(t-j+1)}\smallsetminus\mathcal{S}_{-(t-j)})\supset\Fc^{-(j-1)}(\mathcal{S})$ is a curve that intersects $\text{closure}(\mathcal{S}^-_1)\subset \mathcal{S}_{-1}$. It is a contradiction that some point is on $\mathcal{S}_{-1}$ but not on $\mathcal{S}_{-(t-j)}$ with $t-j\ge 1$.

    Otherwise, $v'\in L(j')=L(j)$, then from the actions of $\Fc$ on the quadrilaterals in \cref{remark:FcActionOnQuadrilaterals,fig:D1celDynamics,fig:D2celDynamics} we see that $\Fc^{-j}(L(j)\sqcup L'(j))=(\mathcal{AS}^{\text{out}}_{j})$ and that is, $\Fc^{-j}$ is a diffeomorphism from $L(j)\sqcup L'(j)$ to $\mathcal{AS}^{\text{out}}_{j}$ (two components). At two points $\Fc^{-j}(v)$ and $\Fc^{-j}(v')$, $(\mathcal{S}_{-(t-j)}\smallsetminus\mathcal{S}_{-(t-j-1)})\supset\Fc^{-j}(\mathcal{S})$ intersects $\Fc^{-j}(L(j))$, which is a connected component of $\mathcal{AS}^{\text{\upshape out}}_{j}$. This contradicts our conclusion in \eqref{item03SingularityCurvesTransveralIntersection}.

    Therefore, $\mathcal{S}\cap L$ contains at most $1$ point. 

    \textbf{Proof of }\eqref{item05SingularityCurvesTransveralIntersection}: Let $VL_0$, $VL_1$ be  the two line segments: connecting from $(-\Phistar,\pi)$ to $(-\Phistar,\pi-\Phistar)$ and connecting from $(\Phistar,0)$ to $(\Phistar,\Phistar)$ in \cref{fig:F_actMroutMRin}.
    
    If $\text{closure}(\mathcal{S})$ has intersection with $\partial M_R$, then $\mathcal{S}\subset M_R$.

    Since $\mathcal{S}\subset\mathcal{S}_{-t}\smallsetminus\mathcal{S}_{-(t-1)}$ with $t\ge3$, $\mathcal{S}$ has no intersection with the line segment connecting $(-\Phistar,0)$ to $(\Phistar,\Phistar)$ and the line segment connecting $(\Phistar,\pi)$ to $(\Phistar,\Phistar)$ (see $\MRin$ of \cref{fig:F_actMroutMRin}), which are curves of $\mathcal{S}^{-}_1$. Then according to Jordan curve theorem, either $\mathcal{S}$ is in $\MRin$ or in one connected component of $M_R\smallsetminus\MRin$ (also see $\MRin$ of \cref{fig:F_actMroutMRin}).

    case (i): $\mathcal{S}$ is in $\MRin$.
    If $\Scal\subset\MRin$ and \textbf{closure$(\Scal)$} has some intersection with $\partial M_R\cap\partial\MRout$, then $\Fc^{-1}(\Scal)$ is a curve segment in $\Mrout$, $\Scal_{-(t-1)}\smallsetminus\Scal_{-(t-2)}\supset\Fc^{-1}(\Scal)$. That is, $\Fc^{-1}(\Scal)$ is connected and is part of a smooth connected component of $(\Scal_{-(t-1)}\smallsetminus\Scal_{-(t-2)})\cap M_r$. Then According to \cref{fig:F_actMroutMRin}, $\text{closure}(\Fc^{-1}(\Scal))$ intersects with $\partial\Mrout$. 
    
    By \eqref{item02SingularityCurvesTransveralIntersection}, \textbf{closure$(\Fc^{-1}(\Scal))$} can have at most one intersection with $\partial M_r\cap\partial\Mrout$ (vertical boundary in \cref{fig:F_actMroutMRin}), which is one of the end point of $\Fc^{-1}(\Scal)$.

    By \eqref{item01SingularityCurvesTransveralIntersection}, \textbf{closure$(\Fc^{-1}(\Scal))$} can have at most one intersection with $\partial\Mrout\smallsetminus\partial M_r$ (the two lines with slope $-\frac{1}{2}$ as boundaries of $\Mrout$ in \cref{fig:F_actMroutMRin}), which is one of the end point of $\Fc^{-1}(\Scal)$.
    
    Therefore, in case (i), $\text{closure}(\mathcal{S})\cap VL$ can contain one or two points that must be the one or two of the endpoints of the curve segment $\mathcal{S}$.

    case (ii): $\mathcal{S}$ is in one connected component of $M_R\smallsetminus\MRin$ (also see $\MRin$ of \cref{fig:F_actMroutMRin}). Assume that $\Scal$ is in the "lower" component with $\theta<\pi/2$ as zoomed in / described in \cref{fig:ImpossibleSingularityIntersection}. Then \textbf{closure$(\Scal)$} can only possibly intersect a $VL_0$ component: the line segment from $(\Phistar,0)$ to $(\Phistar,\Phistar)$
    
    Since the $\Fc$ preimage of $VL_0$ is the line segment from $(-\Phistar,\Phistar)$ to $(\Phistar,0)$ which is a component of $\mathcal{S}^{+}_1$ and $\partial \MRout\smallsetminus M_R$ (the $L$ in the context of \eqref{item04SingularityCurvesTransveralIntersection}), \textbf{closure$(\Fc^{-1}(\Scal))$} intersects the line segment from $(-\Phistar,\Phistar)$ to $(\Phistar,0)$.

    Since $\Fc^{-1}(\Scal)$ is connected and part of a connected component of $\Scal_{-(t-1)}\smallsetminus\Scal_{-(t-2)}$, if \textbf{closure$(\Fc^{-1}(\Scal))$} could intersect \textbf{closure$(\Fc^{-1}(VL_0))$} at two points, then the same analysis performed for \eqref{item03SingularityCurvesTransveralIntersection} implies that the following could happen:

    (a): $\Fc^{-1}(\Scal)$ or its some $\Fc$ iteration preimage could intersect with $\mathcal{S}^{-}_{1}$ (segment that connects $(-\Phistar,0)$ to $(\Phistar,\Phistar)$).
    
    (b): The closure of some $\Fc$ iteration preimage of $\Fc^{-1}(\Scal)$ could intersect a component of \textbf{closure$(\mathcal{AS}^{\text{out}}_{j})$} at two points.
    
    Both are impossible by $\Scal\subset\big(\Scal_{-(t-1)}\smallsetminus\Scal_{-(t-2)}\big)$ and \eqref{item03SingularityCurvesTransveralIntersection}. So \textbf{closure$(\Fc^{-1}(\Scal))$} can intersect $\Fc^{-1}(LV_0)$ at most $1$ point, thus \textbf{closure$(\Scal)$} can intersect $LV_0$ at most $1$ point in case (ii) which must be the end point of \textbf{closure$(\Scal)$}.

    The conclusions for case (i) and (ii) together become \eqref{item05SingularityCurvesTransveralIntersection}.
\end{proof}
\begin{definition}{[\cite[Definition 2.1]{MR3082535} for $2$ dimension]}\label{def:RegularSingularityCurve} A compact subset of lemon billiard phase space $\mathcal{A}\subset M$ is said to be \emph{regular} if it is a union of finitely many compact sets $\mathcal{A}_1,\cdots,\mathcal{A}_p$ of $1$-dimensional $C^2$ submanifolds of $M$ such that 
\begin{enumerate}
    \item $\mathcal{A}_i=\text{closure}(\mathcal{A}_i)$
    \item $\mathcal{A}_i\cap\mathcal{A}_j\subset\partial\mathcal{A}_i$ for $i\ne j$.
    \item $\partial\mathcal{A}_i$ are two points.
\end{enumerate}
\end{definition}
\begin{corollary}{[L1 condition in \cref{thm:lemonLET}, regularity of singularity curves, \cite[Definition 2.1]{MR3082535} for dimension $2$]}\label{corollary:RegSC}
    For lemon billiard and the singularity curve $\mathcal{S}_{-k}$ in \cref{def:BasicNotations0} for each $k\ge0$ is a closed subset of phase space: $\mathcal{S}_{-k}\subset M$ such that \begin{enumerate}
        \item\label{item01RegSC} $\mathcal{S}_{-k}$ is a union of finitely many smooth closed connected curve segments ($C^2$ embedded closed interval). $\mathcal{A}_1,\cdots,\mathcal{A}_{n_k}$
        \item\label{item02RegSC} Each curve segment $\mathcal{A}_i$ intersects with some another curve segment $\mathcal{A}_j$ at boundary/endpoint i.e. $\mathcal{A}_i\cap\mathcal{A}_j\subset\partial\mathcal{A}_i$, $i\ne j$.
        \item\label{item03RegSC} $\partial\mathcal{A}_i$ are two points.
    \end{enumerate}
    And the singularity curve $\Scal_k$ also satisfies the 3 conditions above. 
\end{corollary}
\begin{figure}[h!]
\begin{minipage}{.45\textwidth}
\begin{center}
    \begin{tikzpicture}[xscale=.4,yscale=.4]
        \pgfmathsetmacro{\PHISTAR}{0.8}
        \pgfmathsetmacro{\WIDTH}{7}
        \pgfmathsetmacro{\LENGTH}{7}
        \tkzDefPoint(-\LENGTH,0.5*\WIDTH){LLM};
        \tkzDefPoint(-\LENGTH,0){LLL};
        \tkzDefPoint(-\LENGTH+\PHISTAR,0){LL};
        \tkzDefPoint(-\LENGTH+\PHISTAR,\WIDTH){LU};
        \tkzDefPoint(-\LENGTH,\WIDTH){LLU};
        \tkzDefPoint(\LENGTH-\PHISTAR,0){RL};
        \tkzDefPoint(\LENGTH-\PHISTAR,\WIDTH){RU};
        \tkzDefPoint(\LENGTH,0){RRL};
        \tkzDefPoint(\LENGTH,\WIDTH){RRU};
        \tkzDefPoint(-\LENGTH+\PHISTAR,\PHISTAR){X};
        \tkzDefPoint(-\LENGTH+\PHISTAR,\WIDTH-\PHISTAR){IX};
        \tkzDefPoint(\LENGTH-\PHISTAR,\PHISTAR){IY};
        \tkzDefPoint(\LENGTH-\PHISTAR,\WIDTH-\PHISTAR){Y};
        \tkzLabelPoint[below](LL){$(\phistar,0)$};
        \tkzLabelPoint[below](RL){$(2\pi-\phistar,0)\qquad$};
        \tkzLabelPoint[above](LU){$(\phistar,\pi)$};
        \tkzLabelPoint[above](RU){$(2\pi-\phistar,\pi)$};
        \tkzLabelPoint[right](IY){$x_*$};
        \tkzLabelPoint[left](IX){$y_*$};
        \node[] at (0.5*\LENGTH,0.8*\WIDTH)  {\Large $M_r\smallsetminus\Mrout$};
        \node[] at (-0.5*\LENGTH,0.2*\WIDTH)  {\Large $M_r\smallsetminus\Mrout$};

        \draw [thick,green, -triangle 45] (IX) --(LL);
        \draw [thick,blue, -triangle 45] (LL) --(RL);
        \draw [thick,red, -triangle 45] (RL) --(IX);
        \draw [thin,dashed] (IX) --(LU);
        \draw [thin,dashed] (IY) --(RL);
        \draw [thick,Purple, -triangle 45] (RU) --(LU);
        \draw [thick,teal, -triangle 45] (LU) --(IY);
        \draw [thick,yellow, -triangle 45] (IY) --(RU);

        \fill [orange, opacity=3/30](LL) -- (RL) -- (IX) -- cycle;
        \fill [orange, opacity=3/30](LU) -- (IY) -- (RU) -- cycle;
    \end{tikzpicture}
    \end{center}
\end{minipage}\hfill
\begin{minipage}{.45\textwidth}
\begin{center}
    \begin{tikzpicture}[xscale=.4,yscale=.4]
        \pgfmathsetmacro{\PHISTAR}{0.8}
        \pgfmathsetmacro{\WIDTH}{7}
        \pgfmathsetmacro{\LENGTH}{7}
        \tkzDefPoint(-\LENGTH,0.5*\WIDTH){LLM};
        \tkzDefPoint(-\LENGTH,0){LLL};
        \tkzDefPoint(-\LENGTH+\PHISTAR,0){LL};
        \tkzDefPoint(-\LENGTH+\PHISTAR,\WIDTH){LU};
        \tkzDefPoint(-\LENGTH,\WIDTH){LLU};
        \tkzDefPoint(\LENGTH-\PHISTAR,0){RL};
        \tkzDefPoint(\LENGTH-\PHISTAR,\WIDTH){RU};
        \tkzDefPoint(\LENGTH,0){RRL};
        \tkzDefPoint(\LENGTH,\WIDTH){RRU};
        \tkzDefPoint(-\LENGTH+\PHISTAR,\PHISTAR){X};
        \tkzDefPoint(-\LENGTH+\PHISTAR,\WIDTH-\PHISTAR){IX};
        \tkzDefPoint(\LENGTH-\PHISTAR,\PHISTAR){IY};
        \tkzDefPoint(\LENGTH-\PHISTAR,\WIDTH-\PHISTAR){Y};
        \tkzLabelPoint[below](LL){$(\phistar,0)$};
        \tkzLabelPoint[below left](RL){$(2\pi-\phistar,0)$};
        \tkzLabelPoint[above](LLM){$\xrightarrow{\Fc=\begin{pmatrix}1 &2\\0 &1\end{pmatrix}}\qquad\quad$};
        \tkzLabelPoint[above](LU){$(\phistar,\pi)$};
        \tkzLabelPoint[above left](RU){$(2\pi-\phistar,\pi)$};
        \tkzLabelPoint[right](Y){$I x_*$};
        \tkzLabelPoint[left](X){$I y_*$};

        \draw [thick,green, -triangle 45] (Y) --(LL);
        \draw [thick,blue, -triangle 45] (LL) --(RL);
        \draw [thick,red, -triangle 45] (RL) --(Y);
        \draw [thin,dashed] (Y) --(RU);
        \draw [thin,dashed] (X) --(LL);
        \draw [thick,Purple, -triangle 45] (RU) --(LU);
        \draw [thick,teal, -triangle 45] (LU) --(X);
        \draw [thick,yellow, -triangle 45] (X) --(RU);
        \fill [orange, opacity=3/30](X) -- (LU) -- (RU) -- cycle;
        \fill [orange, opacity=3/30](LL) -- (RL) -- (Y) -- cycle;
        
    \end{tikzpicture}
    \end{center}
\end{minipage}

\begin{minipage}{.45\textwidth}
\begin{center}
    \begin{tikzpicture}[scale=.35]
        \pgfmathsetmacro{\PHISTAR}{2.5};
        \pgfmathsetmacro{\WIDTH}{15};
        \pgfmathsetmacro{\LENGTH}{7.5};
        \tkzDefPoint(-\LENGTH,0.5*\WIDTH){LLM};
        \tkzDefPoint(-\LENGTH,0){LLL};
        \tkzDefPoint(-\PHISTAR,0){LL};
        \tkzDefPoint(-\PHISTAR,\WIDTH){LU};
        \tkzDefPoint(-\LENGTH,\WIDTH){LLU};
        \tkzDefPoint(\PHISTAR,0){RL};
        \tkzDefPoint(\PHISTAR,\WIDTH){RU};
        \tkzDefPoint(\LENGTH,0){RRL};
        \tkzDefPoint(\LENGTH,\WIDTH){RRU};
        \tkzDefPoint(-\PHISTAR,\PHISTAR){X};
        \tkzDefPoint(-\PHISTAR,\WIDTH-\PHISTAR){IX};
        \tkzDefPoint(\PHISTAR,\WIDTH-\PHISTAR){Y};
        \tkzDefPoint(\PHISTAR,\PHISTAR){IY};
        \draw [thick,blue, -triangle 45] (LL) --(X);
        \draw [thick,red, -triangle 45] (X) --(RL);
        \draw [thick,Green, -triangle 45] (RL) --(LL);
        \draw [thin,dashed] (X) --(LU);
        \draw [thin,dashed](RL)--(Y);
        \draw [thick,Purple, -triangle 45] (LU) --(RU);
        \draw [thick,teal, -triangle 45] (RU) --(Y);
        \draw [thick,yellow, -triangle 45] (Y) --(LU);
        \tkzLabelPoint[below left](LL){$(-\Phistar,0)$};
        \tkzLabelPoint[below right](RL){$(\Phistar,0)$};
        \tkzLabelPoint[right](RU){$(\Phistar,\pi)$};
        \tkzLabelPoint[left](LU){$(-\Phistar,\pi)$};
        \fill [green, opacity=3/30](LL) -- (X) -- (RL) -- cycle;
        \fill [green, opacity=3/30](LU) -- (Y) -- (RU) -- cycle;
        \node[] at (-0.8*\LENGTH,0.1*\WIDTH)  {\small $M_R\smallsetminus\MRout$};
        \node[] at (0.8*\LENGTH,0.9*\WIDTH)  {\small $M_R\smallsetminus\MRout$};
    \end{tikzpicture}
\end{center}
\end{minipage}
\begin{minipage}{.45\textwidth}
\begin{center}
    \begin{tikzpicture}[scale=.35]
        \pgfmathsetmacro{\PHISTAR}{2.5};
        \pgfmathsetmacro{\WIDTH}{15};
        \pgfmathsetmacro{\LENGTH}{7.5};
        \tkzDefPoint(-\LENGTH,0.5*\WIDTH){LLM};
        \tkzDefPoint(-\LENGTH,0){LLL};
        \tkzDefPoint(-\PHISTAR,0){LL};
        \tkzDefPoint(-\PHISTAR,\WIDTH){LU};
        \tkzDefPoint(-\LENGTH,\WIDTH){LLU};
        \tkzDefPoint(\PHISTAR,0){RL};
        \tkzDefPoint(\PHISTAR,\WIDTH){RU};
        \tkzDefPoint(\LENGTH,0){RRL};
        \tkzDefPoint(\LENGTH,\WIDTH){RRU};
        \tkzDefPoint(-\PHISTAR,\PHISTAR){X};
        \tkzDefPoint(-\PHISTAR,\WIDTH-\PHISTAR){IX};
        \tkzDefPoint(\PHISTAR,\WIDTH-\PHISTAR){Y};
        \tkzDefPoint(\PHISTAR,\PHISTAR){IY};
        \draw [thick,blue, -triangle 45] (LL) --(IY);
        \draw [thick,red, -triangle 45] (IY) --(RL);
        \draw [thick,Green, -triangle 45] (RL) --(LL);
        \draw [thin,dashed] (LL) --(IX);
        \draw [thin,dashed](RU)--(IY);
        \draw [thick,Purple, -triangle 45] (LU) --(RU);
        \draw [thick,teal, -triangle 45] (RU) --(IX);
        \draw [thick,yellow, -triangle 45] (IX) --(LU);
        \fill [green, opacity=3/30](LL) -- (IY) -- (RL) -- cycle;
        \fill [green, opacity=3/30](LU) -- (IX) -- (RU) -- cycle;
        \tkzLabelPoint[below left](LL){$(-\Phistar,0)$};
        \tkzLabelPoint[below right](RL){$(\Phistar,0)$};
        \tkzLabelPoint[right](RU){$(\Phistar,\pi)$};
        \tkzLabelPoint[left](LU){$(-\Phistar,\pi)$};
        \tkzLabelPoint[above](LLM){$\xrightarrow{\Fc=\begin{pmatrix}1 &2\\0 &1\end{pmatrix}}\qquad\quad$};
    \end{tikzpicture}
\end{center}
\end{minipage}\hfill
\caption{$\Fc$ restricted on the $4$ connected components of $M_R\smallsetminus\MRout$, $M_r\smallsetminus\Mrout$ and its continuous extension to the components boundary as marked by same color edges identifying map.}\label{fig:FrestrictedOnMrdeleteMrout}
\end{figure}
\begin{proof}
    It suffices to prove for $\Scal_{-k}$ since by symmetry $\Scal_k=I\Scal_{-k}$. In \cref{fig:Mr,fig:MR} it is easy to verify $\Scal_0$, $\Scal_{-1}$ and $\Scal_{-1}$ satisfy the conditions \eqref{item01RegSC}\eqref{item02RegSC}\eqref{item03RegSC} since they are the boundaries of the phase space, the line segments with slope $\pm1/2$, $\pm1/4$, the two curve segments in $M_R$ connecting $(-\Phistar,\pi/2+\phistar-\Phistar)$ with $(\Phistar,\pi/2-\phistar/2+\Phistar)$ and connecting $(-\Phistar,\pi/2-\phistar+\Phistar)$ with $(\Phistar,\pi/2+\phistar/2-\Phistar)$ and curve segments $\mathcal{AS}^{\text{\upshape in}}_{-1}$ with boundary(endpoints) on $\partial M$.
    
    Note that from \cref{def:BasicNotations0}, $\mathcal{S}_{-k}$ is recursively defined as $\Scal_{-k}=\Scal_{-k+1}\sqcup\Big\{x\bigm|x\in M\setminus \Scal_{-k+1},\mathcal{F}(x)\in \Scal_{-k+1}\Big\}=\Scal_{-(k-1)}\sqcup\big(\Scal_{-k}\smallsetminus\Scal_{-(k-1)}\big)$ and  \cref{notation:symplecticformmeasurenbhd}, $\Scal^{-}_k=\Scal_{-k}\smallsetminus\Scal_{-(k-1)}$. Hence $\Scal_{-k}=\Scal_{-(k-1)}\sqcup\Scal^{-}_k$. 

    $\Scal^{-}_1$ have $4$ connected components: $2$ straight line segments of $\partial\Mrin$ with slope $1/2$ and $2$ straight line segments of $\partial\MRin$ with slope $1/2$. It is not difficult to verify that $\Scal^{-}_2$ have $8$ connected components, $2$ straight line segments in $M_r$ with slope $1/4$, $2$ straight line segments $M_R$ with slope $1/4$, the two curve segments in $M_R$ connecting $(-\Phistar,\pi/2+\phistar-\Phistar)$ with $(\Phistar,\pi/2-\phistar/2+\Phistar)$ and connecting $(-\Phistar,\pi/2-\phistar+\Phistar)$ with $(\Phistar,\pi/2+\phistar/2-\Phistar)$ and the $2$ arched curve segment $\mathcal{AS}^{\text{\upshape in}}_{-1}$. Hence, $\Scal^{-}_2$ is a disjoint union of finitely many open curve segments with end points at $\mathcal{S}_{-1}$. Here, an open curve segment means that it is the interior of its closure.

    We suppose that $k\ge3$, $\Scal^-_{k-1}$ is a disjoint union of finitely many open curve segments with their end points (boundary) on $\mathcal{S}_{-(k-2)}$. Note that $k\ge3$, $\Scal^-_{k}$ is recursively defined as $\Scal^-_{k}=\Fc(\Scal^-_{k-1}\smallsetminus\Scal_1)$. $\Scal_1$ devided the phase space $M$ to $6$ connected components. $1$ component $\Mrout$, $1$ component $\MRout$ and $2$ components $M_r\smallsetminus\Mrout$ and $2$ components $M_R\smallsetminus\MRout$. $\Fc$ restricted on each component is a diffeomorphism to its image.
    
    $\Fc$ restricted on $\Mrout$ is a diffeomorpshism to $\MRin$ with continuous extension to the boundary shown in \cref{fig:F_actMroutMRin}, $\Fc$ restricted on $\MRout$ is a diffeomorphism to $\Mrin$ with continuous extension to the boundary shown in \cref{fig:F_actMRoutMrin}. \cref{Corollary:SingularityCurvesTransveralIntersection}\eqref{item01SingularityCurvesTransveralIntersection}\eqref{item02SingularityCurvesTransveralIntersection}\eqref{item04SingularityCurvesTransveralIntersection}\eqref{item05SingularityCurvesTransveralIntersection} imply that $\Scal^-_{k}\cap\Mrout$ and $\Scal^-_{k}\cap\MRout$ can have finitely many intersections with $\partial\MRout$ and $\partial \Mrout$, so images $\Fc(\Scal^-_{k}\cap\Mrout)$ and $\Fc(\Scal^-_{k}\cap\MRout)$ can have finitely many intersections with $\partial \Mrin$ and $\partial \MRin$ respectively. Therefore, images $\Fc(\Scal^-_{k}\cap\Mrout)$ and $\Fc(\Scal^-_{k}\cap\MRout)$ are disjoint unions of finitely many open curve segments with end points on $\mathcal{S}_{-k}$.

    As shown in \cref{fig:FrestrictedOnMrdeleteMrout}, $\Fc$ restricted on the connected components of $M_r\smallsetminus\Mrout$ and $M_R\smallsetminus\MRout$ are the linear maps $\begin{pmatrix}1 &2\\0 &1\end{pmatrix}$ in the coordinates $\phi\theta$ and $\Phi\theta$, respectively. \cref{Corollary:SingularityCurvesTransveralIntersection}\eqref{item01SingularityCurvesTransveralIntersection}\eqref{item02SingularityCurvesTransveralIntersection}\eqref{item04SingularityCurvesTransveralIntersection}\eqref{item05SingularityCurvesTransveralIntersection} imply that images $\Fc(\Scal^-_{k}\cap(M_r\smallsetminus\Mrout))$ and $\Fc(\Scal^-_{k}\cap(M_R\smallsetminus\MRout))$ are disjoint unions of finitely many open curve segments with end points on $\mathcal{S}_{-k}$.

    So, by induction, for every $k\ge1$, $\Scal_{-k}\smallsetminus\mathcal{S}_{-(k-1)}$ is a disjoint union of finitely many open curve segments so that by induction $\Scal_{-k}$ satisfies the regularity conditions \eqref{item01RegSC}\eqref{item02RegSC}\eqref{item03RegSC}.
\end{proof}
\begin{corollary}\label{corollary:finiteintersectionwithS1}
    For each $k\ge1$, $\Scal^+_1\cap\Scal_{-k}$ is a finite set. (By symmetry, so is $\Scal^-_1\cap\Scal_{k}$.) 
\end{corollary}
\begin{proof}
    $\Scal_{-k}=\Scal_0\sqcup\Scal^-_1\sqcup\cdots\sqcup\Scal^-_k$ with $\Scal^-_j=\Scal_{-j}\smallsetminus\Scal_{-(j-1)}$, $j=1,\cdots,k$. By \cref{Corollary:SingularityCurvesTransveralIntersection}\eqref{item01SingularityCurvesTransveralIntersection}\eqref{item04SingularityCurvesTransveralIntersection}, $\Scal^-_j\cap\Scal^+_1$ is a finite set, so $\Scal_{-k}\cap\Scal^+_1$ is a finite set.
\end{proof}

\subsection[Local ergodicity and global ergodicity]{Local ergodicity and global ergodicity}\label{subsec:LETGET}\hfill
\begin{corollary}[Local ergodicity on sections]\label{corollary:LETonSection}\hfill
    \begin{enumerate}
        \item\label{item01LETonSection} For each $x\in(\hat{\mathfrak{U}}_+\cup\hat{\mathfrak{U}}_-)\smallsetminus\Scal_{\infty}$, $\exists$ an open neighborhood $\mathcal{O}_x$ of $x$ such that $\mathcal{O}_x$ is contained in one ergodic component of $\Fc$ with respect to the measure $\mu$ in \cref{notation:symplecticformmeasurenbhd}.
        \item\label{item02LETonSection} For each $x\in(\hat{\mathfrak{U}}_+\cup\hat{\mathfrak{U}}_-)\smallsetminus\Scal_{-\infty}$, $\exists$ an open neighborhood $\mathcal{O}_x$ of $x$ such that $\mathcal{O}_x$ is contained in one ergodic component of $\Fc^{-1}$, thus of $\Fc$ with respect to the measure $\mu$ in \cref{notation:symplecticformmeasurenbhd}. 
    \end{enumerate}
\end{corollary}
\begin{proof}
\textbf{\eqref{item01LETonSection}:} For $x\in\hat{\mathfrak{U}}_+\smallsetminus\Scal_{\infty}$, it is a sufficient point with a quadruple $(l, N, U, \mathcal{C})$ by \cref{theorem:sufficentptsonMhatplus}. For $x\in\hat{\mathfrak{U}}_-\smallsetminus\Scal_{\infty}$, it is a sufficient point with a quadruple $(l, N, U, \mathcal{C})$ by \cref{theorem:sufficentptsonMhatminus}. Then \cref{corollary:RegSC} gives the L1 condition. \cref{corollary:L2AlignmentCondition} gives the L2 condition. \cref{lemma:SCAnsatz} gives the L3' condition. And \cref{theorem:L4condition} gives the L4 condition. Hence \cref{thm:lemonLET} gives the local ergodicity.

\textbf{\eqref{item02LETonSection}:} Since the symmetry $I:(\phi,\theta)\rightarrow(\phi,-\theta)$ is a diffeomorphism between $\hat{\mathfrak{U}}_+$ and $\hat{\mathfrak{U}}_-$ (see \cref{def:UIntRegionofD,fig:hatUplusUminus}), $I\Scal_{\infty}=\Scal_{-\infty}$, for every $x\in(\hat{\mathfrak{U}}_+\cup\hat{\mathfrak{U}}_-)\smallsetminus\Scal_{-\infty}$, $Ix\in(\hat{\mathfrak{U}}_+\cup\hat{\mathfrak{U}}_-)\smallsetminus\Scal_{\infty}$. By \eqref{item01LETonSection}, $Ix$ has an open neighborhood $\mathcal{O}_{Ix}$ to be contained in one ergodic component of $\Fc$. So $x$ has an open neighborhood $\mathcal{O}_x$ to be contained in one ergodic component of $\Fc^{-1}$, thus of $\Fc$.

\end{proof}
\begin{proposition}[analogous to {\cite[Proposition 6.19 and Theorem 6.20]{cb}}]\label{Proposition:globalErgodic} With respect to the measure $\mu$ in \cref{notation:symplecticformmeasurenbhd}, the following ergodicity conclusions hold.
\begin{enumerate} 
    \item\label{item01globalErgodic} Each of the three connected components of $\hat{\mathfrak{U}}_-$ (in \cref{fig:hatUplusUminus}) belongs (mod 0) to an ergodic component of $\Fc$.
    \item\label{item02globalErgodic} $\hat{\mathfrak{U}}_-$ belongs (mod 0) to one ergodic component of $\Fc$
    \item\label{item03globalErgodic} The entire phase space $M$ belongs (mod 0) to one ergodic component of $\Fc$.
\end{enumerate}
\end{proposition}
\begin{proof}
    \textbf{\eqref{item01globalErgodic}:} A point $x\in\hat{\mathfrak{U}}_-$ is said to be a double singularity if $x\in\Scal_{-m}\cap\Scal_{n}$ for some $n\ge0,m\ge0$; the set of these is $\bigcup_{m\ge0,n\ge0}(\Scal_{-m}\cap\Scal_{n})=\Scal_{-\infty}\cup\Scal_{\infty}$. By \cref{Theorem:SingularityAlignmentOnSection}\eqref{item01:SingularityAlignmentOnSection}, $\Scal_{-m}\cap\hat{\mathfrak{U}}_-$ is a finite union of curve segments with slopes of the tangent vectors $\frac{d\theta}{d\phi}\in(0,1)$. By \cref{Theorem:SingularityAlignmentOnSection}\eqref{item03:SingularityAlignmentOnSection} $\hat{\mathfrak{U}}_{+}\cap\Scal_{-n}$ has tangent vectors in the positive quadrant $\mathcal{C}$. Hence, by symmetry $\hat{\mathfrak{U}}_{-}\cap\Scal_{n}= I(\hat{\mathfrak{U}}_{+}\cap\Scal_{-n})$ has tangent vectors in the negative quadrant $I\mathcal{C}=\mathcal{C}'$. That means $\Scal_{-m}\cap\Scal_{n}\cap\hat{\mathfrak{U}}_-$ is finite. So, the set of double singularities: $\Scal_{-\infty}\cup\Scal_{\infty}$ is countable. 
    
    By \cref{corollary:LETonSection}\eqref{item01LETonSection}\eqref{item02LETonSection}, every $x\in\hat{\mathfrak{U}}_-\smallsetminus(\Scal_{-\infty}\cup\Scal_{\infty})$ has an open neighborhood contained in one ergodic component of $\Fc$. 
    
    Consider $x,y\in\hat{\mathfrak{U}}_-\smallsetminus(\Scal_{-\infty}\cup\Scal_{\infty})$ belonging to a connected component of $\hat{\mathfrak{U}}_-$. Since $\Scal_{-\infty}\cup\Scal_{\infty}$ is countable, there is an compact arc $l$ belonging to this connected component and connecting $x$ to $y$ such that no $z\in l$ is a double singularity point (see the proof of \cite[Proposition 6.19]{cb}). 
    
    This means that for every $z\in l$, there is an open neighborhood $\mathcal{O}_z$ of $z$ such that $\mathcal{O}_z$ belongs to one ergodic component. We use the same argument as in \cite[Proposition 6.19 and Figure 6.15]{cb}.
    
    Being compact, $l$ is covered by finitely many open (and overlapping) $\mathcal{O}_z$ with $z\in l$. Then $x$ and $y$ are in the same ergodic component because these are a.e.\ disjoint. 
Thus, almost all points of a connected component of $\hat{\mathfrak{U}}_-$ are in the same ergodic component of $\Fc$.

    \textbf{\eqref{item02globalErgodic}:} The preimage $\Fc^{-3}(\Mrin\cap\Mrout)$ of the connected component $\Mrin\cap\Mrout$ (the $K1$-$K2$-$K3$-$K4$ rhombus in \cref{fig:hatUplusUminus}) intersects each of the two connected components of $\Fc^{-1}(\Mrout\smallsetminus\Mrin)$ (two pink strips in \cref{fig:hatUplusUminus}) in a set of positive measure.
    
    By disjointness of ergodic components, 
    all three connected components of $\hat{\mathfrak{U}}_-$ (mod 0) belong to one ergodic component (also see the proof of \cite[Theorem 6.20]{cb}).
    
    \textbf{\eqref{item03globalErgodic}:} By \cref{remark:UIntRegionofD},  $M\aeq\bigcup_{k\in\mathbb{Z}}\Fc^k(\hat{\mathfrak{U}}_{-})$, so $M$ is in one ergodic component of $\Fc$.
\end{proof}
\begin{remark}[Next steps]\label{REMFutureWork}\strut\COMMENT{The first two are not worded like they are in the future}\hfill

    \begin{enumerate}
        \item We will prove the K-property in a way similar to the proof for Wojtkowski billiards in \cite{OKPPHB}.
        \item Using uniform hyperbolicity of the return map, we will prove that the singularity curves on those sections have uniformly bounded curvatures, even though this condition is not required in our local ergodic theorem. With uniformly bounded curvatures and if we can prove the so-called $1$-step or $m$-step expansion condition in \cite{MR2358286,MR2540156,MR2810510}, we will obtain mixing rate estimates.
        \item One might next consider whether the same techniques succeed when the flatter side is curved inward rather than outward; these are called \emph{moon billiards}  \cite{MR3456006}. On one hand, such inward-bent side is dispersing, which tends to help establish hyperbolicity. On the other hand, such a dispersing boundary arc produces more serious singularities (grazing collisions), and this adds significant difficulty.\COMMENT{Edit when done writing.\wentao{I have add that our next step is the proving the LET \Cref{thm:lemonLET} in \cite{MR3082535} for lemon billiard.}}
    \end{enumerate}
\end{remark}
\bibliographystyle{alpha}
\bibliography{lemonExpansion}  
\end{document}